\documentclass[envcountchap, envcountsame]{svmono}
\usepackage{etex}

\usepackage{helvet}
\usepackage{courier}
\usepackage{type1cm}         

\usepackage{makeidx}         

\usepackage{hyperref} 
\usepackage{xcolor}
\usepackage{amsmath}
\usepackage{amscd}
\usepackage{amssymb}
\usepackage{euscript}
\usepackage[frame,ps,matrix,arrow,curve,rotate,all,2cell,tips]{xy}
\usepackage{epic,eepic}
\usepackage{graphicx}

\usepackage{multicol}        

\usepackage{enumerate}

\usepackage{tikz}
\usetikzlibrary{arrows,decorations.pathmorphing}
\usetikzlibrary{backgrounds,positioning}
\usetikzlibrary{fit,petri,shapes.misc}

\usepackage{tikz-cd}

\tikzset{auto}

\tikzset{empty/.style={circle,inner sep=0pt,minimum size=6mm}}
\tikzset{emptyvt/.style={circle,inner sep=0pt,minimum size=0mm}}

\tikzset{plain/.style={circle,draw,very thick,
inner sep=0pt,minimum size=6mm}}

\tikzset{fatplain/.style={rounded rectangle,draw,very thick,minimum size=6mm}}

\tikzset{bigplain/.style={rounded rectangle,draw,very thick,minimum size=.8cm}}

\tikzset{yellowvt/.style={circle,draw,fill=yellow,very thick,inner sep=0pt,minimum size=6mm}}

\tikzset{bluevt/.style={circle,draw,fill=blue!20,very thick,inner sep=0pt,minimum size=6mm}}

\tikzset{greenvt/.style={circle,draw,fill=green!30,very thick,inner sep=0pt,minimum size=6mm}}

\tikzset{redvt/.style={circle,draw,fill=red!30,very thick,inner sep=0pt,minimum size=6mm}}

\tikzset{arrow/.style={->,thick}}
\tikzset{dashedarrow/.style={->,dashed,thick}}
\tikzset{dottedarrow/.style={->,dotted,thick}}
\tikzset{mapto/.style={|->,thick}}

\tikzset{implies/.style={thick,double,double equal sign distance,-implies}}

\tikzset{line/.style={thick}}
\tikzset{dottedline/.style={dotted,thick}}
\tikzset{dashedline/.style={dashed,thick}}

\tikzset{inputleg/.style={<-,thick}}
\tikzset{outputleg/.style={->,thick}}
\tikzset{dottedinput/.style={<-,dotted,thick}}

\spnewtheorem{sublemma}{Sublemma}{\bfseries}{\itshape}
\spnewtheorem{convention}{Convention}{\bfseries}{\rmfamily}
\spnewtheorem{notation}{Notation}{\bfseries}{\rmfamily}

\newcommand{\minus}{{-}}

\newcommand{\be}{\mathbf{e}}

\newcommand{\bone}{\mathbf{1}}

\newcommand{\bonew}{\mathbf{1}^{\wheel}}

\newcommand{\bi}{\mathbf{i}}
\newcommand{\bN}{\mathbf{N}}

\newcommand{\frakC}{\mathfrak{C}}
\newcommand{\fC}{\mathfrak{C}}

\newcommand{\fD}{\mathfrak{D}}

\newcommand{\fE}{\mathfrak{E}}

\newcommand{\fF}{\mathfrak{F}}

\newcommand{\sK}{\mathsf{K}}

\newcommand{\sO}{\mathsf{O}}

\newcommand{\sP}{\mathsf{P}}

\newcommand{\pk}{\sP_{\sK}}

\newcommand{\sQ}{\mathsf{Q}}
\newcommand{\qk}{\sQ_{\sK}}
\newcommand{\sR}{\mathsf{R}}
\newcommand{\sS}{\mathsf{S}}

\newcommand{\sW}{\mathsf{W}}

\newcommand{\ua}{\underline{a}}
\newcommand{\ub}{\underline{b}}
\newcommand{\uc}{\underline{c}}

\newcommand{\ud}{\underline{d}}

\newcommand{\ue}{\underline{e}}
\newcommand{\uf}{\underline{f}}
\newcommand{\ug}{\underline{g}}

\newcommand{\uv}{\underline{v}}
\newcommand{\uw}{\underline{w}}
\newcommand{\ux}{\underline{x}}
\newcommand{\uy}{\underline{y}}

\newcommand{\ghat}{\widehat{G}}
\newcommand{\hhat}{\widehat{H}}

\newcommand{\khat}{\widehat{K}}
\newcommand{\lhat}{\widehat{L}}
\newcommand{\phat}{\widehat{P}}
\newcommand{\qhat}{\widehat{Q}}

\newcommand{\adjoint}{
\nicexy{ \ar@<2pt>[r] & \ar@<2pt>[l]}}

\newcommand{\defn}{~\buildrel \text{def} \over=~}

\renewcommand{\hookrightarrow}{\nicexy{\ar@{^{(}->}[r] &}}

\newcommand{\nicearrow}{\SelectTips{cm}{10}}
\newcommand{\nicexy}{\nicearrow\xymatrix}

\newcommand{\wheel}{{\circlearrowright}}

\newcommand{\wheelc}{{\wheel_{c}}}

\renewcommand{\to}{\longrightarrow}

\newcommand{\xiij}{\xi^i_j}

\newcommand{\fbar}{\overline{f}}

\newcommand{\catc}{\mathcal{C}}
\newcommand{\catd}{\mathcal{D}}

\newcommand{\catp}{\mathcal{P}}

\newcommand{\catcs}{\catc^{\dis(S)}}

\newcommand{\category}{\mathtt{Category}}

\newcommand{\alg}{\mathbf{Alg}}

\newcommand{\Set}{\mathtt{Set}}

\newcommand{\Gammai}{\varGamma_i}
\newcommand{\Gammao}{\varGamma_o}
\newcommand{\Gammaw}{\varGamma_\wheel}
\newcommand{\Lambdaw}{\Lambda_\wheel}

\newcommand{\gupcset}{\Set^{\varGamma^{\op}}}
\newcommand{\gupciset}{\Set^{\Gammai^{\op}}}
\newcommand{\gupcoset}{\Set^{\Gammao^{\op}}}

\newcommand{\gwheelcset}{\Set^{\Gammaw^{\op}}}

\newcommand{\gupdset}{\Set^{\varTheta^{\op}}}
\newcommand{\omegaset}{\Set^{\varOmega^{\op}}}
\newcommand{\sset}{\Set^{\varDelta^{\op}}}

\newcommand{\cg}{\mathcal G}
\newcommand{\Gr}{\mathtt{Gr}}

\newcommand{\gup}{\Gr^{\uparrow}}
\newcommand{\gupc}{\gup_{\text{c}}}
\newcommand{\gupci}{\gup_{\text{ci}}}
\newcommand{\gupco}{\gup_{\text{co}}}
\newcommand{\gupcs}{\gup_{\text{cs}}}
\newcommand{\gupd}{\gup_{\text{di}}}
\newcommand{\gwheelc}{\Gr^\wheel_{\text{c}}}

\newcommand{\ULin}{\mathtt{ULin}}

\newcommand{\uoperad}{\mathtt{UTree}}
\newcommand{\utree}{\mathtt{UTree}}

\newcommand{\properad}{\mathtt{Properad}}
\newcommand{\properadi}{\mathtt{Properad}_{\mathtt{i}}}
\newcommand{\properado}{\mathtt{Properad}_{\mathtt{o}}}
\newcommand{\properads}{\mathtt{Properad}_{\mathtt{s}}}

\newcommand{\wproperad}{\mathtt{Properad}^{\wheel}}

\newcommand{\compi}{\circ_i}
\newcommand{\compj}{\circ_j}

\newcommand{\jcompi}{{\,{}_j\circ_i\,}}
\newcommand{\icompj}{{{}_i\circ_j}}
\newcommand{\onecompone}{{\,{}_1\circ_1\,}}
\newcommand{\onecompi}{{\,{}_1\circ_i\,}}
\newcommand{\jcompone}{{\,{}_j\circ_1\,}}
\newcommand{\kcompone}{{\,{}_k\circ_1\,}}
\newcommand{\scompr}{{{}_s\circ_r}}

\newcommand{\pofc}{\catp(\frakC)}
\newcommand{\pofcop}{\pofc^{\mathrm{op}}}
\newcommand{\ptwoc}{\pofcop \times \pofc}

\newcommand{\SC}{\mathsf{S}(\fC)}
\newcommand{\SD}{\mathsf{S}(\fD)}
\newcommand{\SE}{\mathsf{S}(\fE)}
\newcommand{\SCD}{\mathsf{S}(\fC \times \fD)}

\newcommand{\ba}{\binom{\ub}{\ua}}

\newcommand{\ccsingle}{\binom{c}{c}}

\newcommand{\ucuc}{\binom{\uc^2}{\uc^1}}

\newcommand{\dc}{\binom{\ud}{\uc}}

\newcommand{\dcsigma}{\binom{\sigma\ud}{\uc\tau}}

\newcommand{\udud}{\binom{\ud^2}{\ud^1}}

\newcommand{\nm}{\binom{n}{m}}

\newcommand{\yx}{\binom{\uy}{\ux}}
\newcommand{\yxlambda}{\binom{\lambda\uy}{\ux\pi}}

\newcommand{\emptyprofh}{(\varnothing;\varnothing)}
\newcommand{\emptyprof}{\binom{\varnothing}{\varnothing}}

\newcommand{\profileg}{\binom{\out(G)}{\inp(G)}}

\newcommand{\profileu}{\binom{\out(u)}{\inp(u)}}
\newcommand{\profilev}{\binom{\out(v)}{\inp(v)}}

\newcommand{\bah}{(\ua;\ub)}

\newcommand{\dch}{(\uc;\ud)}

\newcommand{\feh}{(\ue;\uf)}

\newcommand{\wvh}{(\uv;\uw)}
\newcommand{\yxh}{(\ux;\uy)}

\newcommand{\andspace}{\quad\text{and}\quad}

\newcommand{\orspace}{\quad\text{or}\quad}

\newcommand{\ving}{v \in \vertex(G)}

\DeclareMathOperator{\dis}{dis}
\DeclareMathOperator{\End}{End}
\DeclareMathOperator{\Hom}{Hom}

\DeclareMathOperator*{\colim}{colim}

\DeclareMathOperator{\edge}{Edge}
\DeclareMathOperator{\edgei}{Edge_i}

\DeclareMathOperator{\Flag}{Flag}

\DeclareMathOperator{\Leg}{Leg}

\DeclareMathOperator{\vertex}{Vt}

\DeclareMathOperator{\id}{id}
\DeclareMathOperator{\Id}{Id}
\DeclareMathOperator{\image}{Im}

\DeclareMathOperator{\inp}{in}

\newcommand{\Iso}{\operatorname{Iso}}

\DeclareMathOperator{\out}{out}

\DeclareMathOperator{\Ob}{Ob}
\DeclareMathOperator{\op}{op}

\DeclareMathOperator{\Sc}{Sc}

\makeindex

\begin{document}

\author{Philip Hackney, Marcy Robertson, and Donald Yau}
\title{Infinity Properads and Infinity Wheeled Properads}


\maketitle

\frontmatter







%
%
%

\begin{dedication}
To Chlo\"e and Elly. \\
\noindent To Rosa.\\
\noindent To Eun Soo and Jacqueline.
\end{dedication}


\preface

This monograph fits in the intersection of two long and intertwined stories. The first part of our story starts in the mid-twentieth century, when it became clear that a new conceptual framework was necessary for the study of higher homotopical structures arising in algebraic topology.
Some better known examples of these higher homotopical structures appear in work of J. F. Adams and S. Mac Lane \cite{maclane_catalg} on the coproduct in the bar construction and work of J. Stasheff \cite{stasheff1}, J. M. Boardman and R. M. Vogt \cite{bvpaper, bv}, and J. P. May \cite{may} on recognition principles for (ordinary, $n$-fold, or infinite) loop spaces.
The notions of `operad' and `prop' were precisely formulated for the purpose of this work; the former is suitable for modeling algebraic or coalgebraic structures, while the latter is also capable of modeling bialgebraic structures, such as Hopf algebras.
Operads came to prominence in other areas of mathematics beginning in the 1990s (but see, e.g. \cite{kad, smi} for earlier examples) through the work of V. Ginzburg and M. Kapranov on Koszul duality \cite{gk}, E. Getzler and J. Jones on two-dimensional topological field theories \cite{getzj,get94}, and M. Kontsevich on deformation quantization \cite{kontsevich}. 

This renaissance in the world of operads \cite{lod96,lsv97} and the popularity of quantum groups \cite{Dr1,Dr2} in the 1990s lead to a resurgence of interest in props, which had long been in the shadow of their little single-output nephew. 
Properads were invented around this time, during an effort of B. Vallette to formulate a Koszul duality for props \cite{vallette}. Properads and props both model algebraic structures with several inputs and outputs, but properads govern a smaller class of such structures, those whose generating operations and relations among operations can be taken to be \emph{connected}. 
This class includes most types of bialgebras that arise in nature, such as biassociative bialgebras, (co)module bialgebras, Lie bialgebras, and Hopf algebras.

Wheeled variants of operads, properads, and props where introdued by M. Markl, S. Merkulov and S. Shadrin \cite{mms} to model algebraic structures with traces.
For instance, one of the simplest examples of a wheeled properad controls \emph{finite-dimensional} associative algebras.
There are numerous applications of wheeled properads in geometry, deformation theory, and mathematical physics \cite{merkulov3}.

\noindent\includegraphics[width=\textwidth]{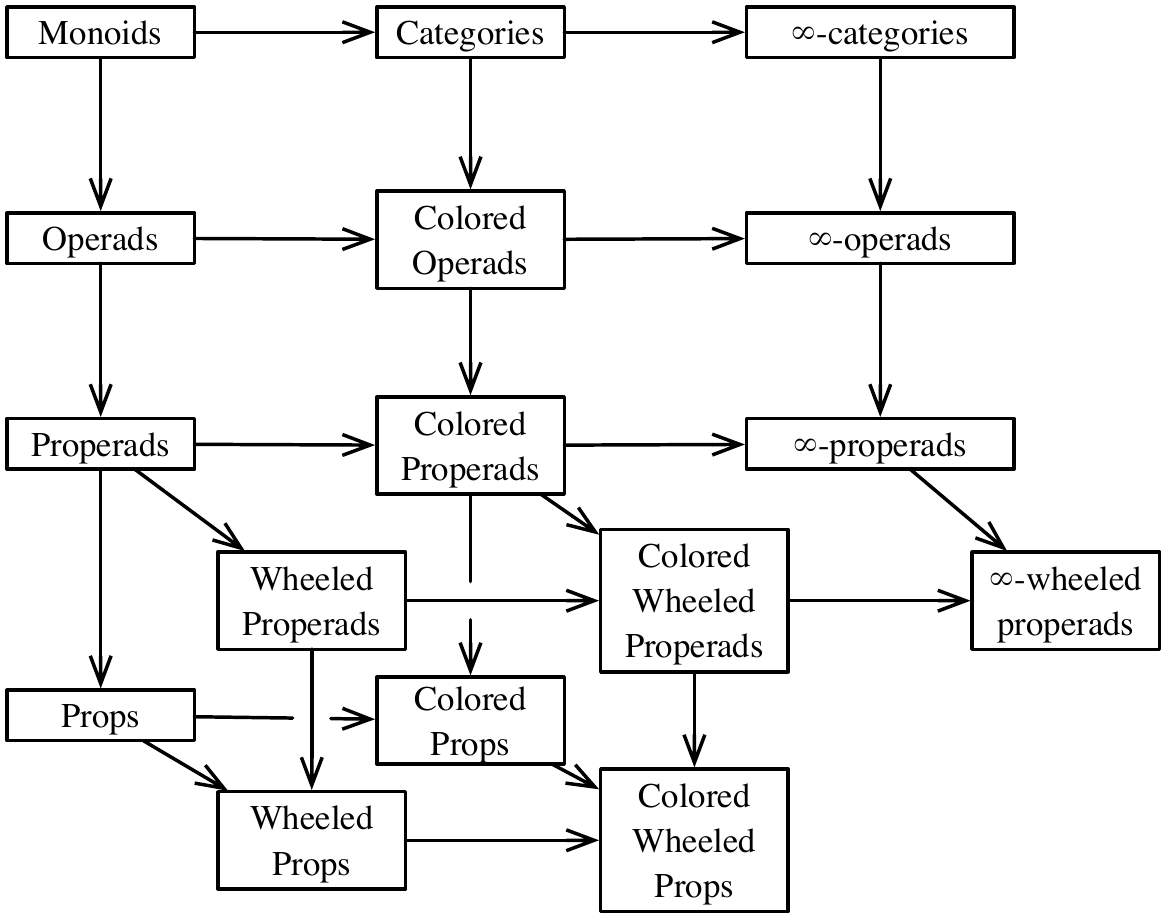}

The use of `colored' or `multisorted' variants of operads or props \cite{brinkmeier, bv}, where composition is only partially defined, allows one to address many other situations of interest.
It allows one, for instance, to model morphisms of algebras associated to a given operad. 
There is a two-colored operad which encodes the data of two associative algebras as well as a map from one to the other; a resolution of this operad precisely gives the correct notion of morphism of $A_\infty$-algebras \cite{markl04,markl02}.
It also provides a unified way to treat operads, cyclic operads, modular operads, properads, and so on: for each there is a colored operad which controls the structure in question.

The second part of our story is an extension of the theory of categories. 
Categories are pervasive in pure mathematics and, for our purposes, can be loosely described as tools for studying collections of objects up to isomorphism and comparisons between collections of objects up to isomorphism. 
When the objects we want to study have a homotopy theory, we need to generalize category theory to identify two objects which are, while possibly not isomorphic in the categorical sense, equivalent up to homotopy. 
For example, when discussing topological spaces we might replace `homeomorphism' with `homotopy equivalence'.
This leads to the theory of $\infty$-categories (or restricted Kan complexes \cite{bv}, quasi-categories, and so on); 
A. Joyal \cite{joyal_notes, joyal_theory} and J. Lurie \cite{lurie} have extended many tools from traditional category theory to $\infty$-category theory.
This extension of category theory has led to new applications in various subjects ranging from  a convenient framework for the study of $A_{\infty}$ and  $E_{\infty}$-ring spectra in stable homotopy theory \cite{abghr}, derived algebraic geometry \cite{dag,toen1,toen_vezzosi1,toen_vezzosi2} and geometric representation theory \cite{toen2, bzn1, bzn2}. 

I. Moerdijk and I. Weiss \cite{mw1} realized that the connection between colored operads and weakly composable maps was worth further exploration, and introduced a way to think about the notion of $\infty$-operad.
They introduced a category of trees $\varOmega$, the `dendroidal category', which contains the simplicial category $\varDelta$.
Dendroidal sets, or presheaves on $\varOmega$, are an extension of simplicial sets, and this extension allows us to consider `quasi-operads' in the category of dendroidal sets which are analogs of quasi-categories.
Moerdijk and Weiss proposed a model for weak $n$-categories based on this formalism, which has been partially validated \cite{lukacs}. There are also relations to connective spectra \cite{bn}, $E_\infty$-spaces \cite{heuts2}, algebraic K-theory \cite{nik}, and group actions on operads \cite{bh}.

This book is a thorough initial investigation into the theory of $\infty$-properads and $\infty$-wheeled properads.  
We here lay the foundation for our aim of exploring the homotopy theory of properads in depth.
This work also serves as a complete guide to the generalized graphs which are pervasive in the study of operads and properads.  
In the final chapter we include a preliminary list of potential applications ranging from string topology to category theory.

This monograph is written for mathematicians in the fields of topology, algebra, category theory, and related areas.  It is written roughly at the second year graduate level.  We assume some very basic knowledge of category theory, as discussed in the standard references \cite{bor1,bor2,maclane98}.  Topics such as monads, simplicial objects, generalized PROPs, and so forth, will be reviewed in the text.


The first two authors would like to thank the third author for his considerable patience and role as a mentor in the preparation of this monograph.
They would also like to thank Tom Fiore for planting the seeds of this project during his talk at the Graduate Student Geometry and Topology Conference in 2010.

The authors would like to thank 
Alexander Berglund, 
Julie Bergner, 
Benoit Fresse, 
David Gepner, 
Mark W. Johnson, 
Andr\'e Joyal, 
Sergei Merkulov, 
Ieke Moerdijk,  
Bruno Vallette, 
Rainer Vogt, and 
Ben Ward 
for their interest and encouragement while completing this project. 
We would particularly like to thank Joachim Kock for sharing his point of view on the definition of generalized graph \cite{kock}. 
Finally, much credit goes to all five anonymous referees who each provided insightful comments and suggestions.

\vspace{\baselineskip}
\begin{flushright}\noindent
London, Los Angeles, Stockholm, Newark\hfill {\it Philip Hackney}\\
January 2014 \hfill {\it Marcy Robertson}\\
February 2015 \hfill {\it Donald Yau} \\
\end{flushright}

\tableofcontents

\mainmatter

\chapter{Introduction}


\abstract{
A theory of $\infty$-properads is developed, extending both the Joyal-Lurie $\infty$-categories and the Cisinski-Moerdijk-Weiss $\infty$-operads.  Every connected wheel-free graph generates a properad, giving rise to the graphical category $\varGamma$ of properads.  Using graphical analogs of coface maps and the properadic nerve functor, an $\infty$-properad is defined as an object in the graphical set category $\gupcset$ that satisfies some inner horn extension property.  Symmetric monoidal closed structures are constructed in the categories of properads and of graphical sets. Strict $\infty$-properads, in which inner horns have unique fillers, are given two alternative characterizations, one in terms of graphical analogs of the Segal maps, and the other as images of the properadic nerve.  The fundamental properad of an $\infty$-properad is characterized in terms of homotopy classes of $1$-dimensional elements. Using all connected graphs instead of connected wheel-free graphs, a parallel theory of $\infty$-wheeled properads is also developed.
}


Let us first recall the notions of a properad and of a wheeled properad.

\section*{(Wheeled) Properads as Generalized Categories}

In an ordinary category, a morphism $\nicexy{x \ar[r]^-{f} & y}$ has one input and one output.  If $\nicexy{y \ar[r]^-{g} & z}$ is another morphism, then the composition $gf$ is uniquely defined.  Moreover, the identity and associativity axioms hold in the strict sense.  There are two natural ways in which the notion of a category can be extended.

The first natural way to extend the notion of a category is to allow morphisms with finite lists of objects as inputs and outputs, together with appropriately chosen axioms that hold in the strict sense.  For example, an operad \cite{may} is a generalization of a category in which a morphism has one output and finitely many inputs, say 
\[
\nicexy{
(x_1, \ldots, x_n) \ar[r]^-{f} & y
}\]
with $n \geq 0$.  We often call such a morphism an operation and denote it by the following decorated graph.
\begin{center}
\begin{tikzpicture}
\matrix[row sep=1cm, column sep=1cm]{
\node [plain, label=below:$...$] (f) {$f$};\\
};
\draw [outputleg] (f) to node[above=.1cm]{$y$} +(0,.7cm);
\draw [inputleg] (f) to node[below left=.1cm]{$x_1$} +(-.7cm,-.5cm);
\draw [inputleg] (f) to node[below right=.1cm]{$x_n$} +(.7cm,-.5cm);
\end{tikzpicture}
\end{center}
Composition in a category is represented linearly.  With multiple inputs, composition in an operad takes on the shape of a tree.  Explicitly, if there are operations
\[
\nicexy{
\left(w^i_1, \ldots , w^i_{k_i}\right) \ar[r]^-{g_i} & x_i
}\]
for each $i$, then the operadic composition $\gamma(f;g_1,\ldots,g_n)$ is represented by the following decorated $2$-level tree.
\begin{center}
\begin{tikzpicture}
\matrix[row sep=.4cm, column sep=1.2cm]{
& \node [plain, label=below:$...$] (f) {$f$}; &\\
\node [plain, label=below:$...$] (g1) {$g_1$}; &&
\node [plain, label=below:$...$] (gn) {$g_n$};\\
};
\draw [outputleg] (f) to node[above=.1cm]{$y$} +(0,.8cm);
\foreach \x in {g1,gn}
{
\draw [arrow] (\x) to (f);
}
\draw [inputleg] (g1) to node[below left=.1cm]{$w^1_1$} +(-.8cm,-.6cm);
\draw [inputleg] (g1) to node[below right=.1cm]{$w^1_{k_1}$} +(.8cm,-.6cm);
\draw [inputleg] (gn) to node[below left=.1cm]{$w^n_1$} +(-.8cm,-.6cm);
\draw [inputleg] (gn) to node[below right=.1cm]{$w^n_{k_n}$} +(.8cm,-.6cm);
\end{tikzpicture}
\end{center}
In particular, its inputs are the concatenation of the lists $\left(w^i_1, \ldots , w^i_{k_i}\right)$ as $i$ runs from $1$ to $n$.  Associativity of the operadic composition takes the form of a $3$-level tree.  There are also unity and equivariance axioms, which come from permutations of the inputs.  We should point out that what we call an operad here is sometimes called a colored operad or a multicategory in the literature.

A properad \cite{vallette} allows even more general operations, where both inputs and outputs are finite lists of objects, say
\[
\nicexy{
(x_1, \ldots , x_m) \ar[r]^-{f} & (y_1, \ldots , y_n)
}\]
with $m, n \geq 0$.  Such an operation is visualized as the following decorated corolla.  
\begin{center}
\begin{tikzpicture}
\matrix[row sep=1cm, column sep=1cm]{
\node [plain, label=above:$...$, label=below:$...$] (f) {$f$};\\
};
\draw [outputleg] (f) to node[above left=.1cm]{$y_1$} +(-.7cm,.4cm);
\draw [outputleg] (f) to node[above right=.1cm]{$y_n$} +(.7cm,.4cm);
\draw [inputleg] (f) to node[below left=.1cm]{$x_1$} +(-.7cm,-.4cm);
\draw [inputleg] (f) to node[below right=.1cm]{$x_m$} +(.7cm,-.4cm);
\end{tikzpicture}
\end{center}
Since both the inputs and the outputs can be permuted, there are bi-equivariance axioms.  The properadic composition is represented by the following graph, called a partially grafted corollas.
\begin{center}
\begin{tikzpicture}
\matrix[row sep=.7cm,column sep=1.5cm] {
\node [plain,label=below:$...$,label=above:$...$] (p1) {$g$}; \\
\node [plain,label=below:$...$,label=above:$...$] (p2) {$f$}; \\
};
\foreach \x in {1,2}
{
\draw [outputleg] (p\x) to +(-.6cm,.4cm);
\draw [outputleg] (p\x) to +(.6cm,.4cm);
\draw [inputleg] (p\x) to +(-.6cm,-.4cm);
\draw [inputleg] (p\x) to +(.6cm,-.4cm);
}
\draw [arrow,bend left=40] (p2) to (p1);
\draw [arrow,bend right=40] (p2) to (p1);
\end{tikzpicture}
\end{center}
For this properadic composition to be defined, a non-empty  sub-list of the outputs of $f$ must match a non-empty sub-list of the inputs of $g$.  Associativity says that if a $3$-vertex connected wheel-free graph has its vertices decorated by operations, then there is a well-defined operation.

One particular instance of the properadic composition is when the partially grafted corollas has only one edge connecting the two vertices, as in the following simply connected graph.
\begin{center}
\begin{tikzpicture}
\matrix[row sep=1.2cm,column sep=1.5cm] {
\node [plain,label=above:$...$] (p1) {$g$}; \\
\node [plain,label=below:$...$] (p2) {$f$}; \\
};
\foreach \x in {1,2}
{
\draw [outputleg] (p\x) to +(-.6cm,.4cm);
\draw [outputleg] (p\x) to +(.6cm,.4cm);
\draw [inputleg] (p\x) to +(-.6cm,-.4cm);
\draw [inputleg] (p\x) to +(.6cm,-.4cm);
}
\draw [arrow] (p2) to (p1);
\end{tikzpicture}
\end{center}
It is called a dioperadic graph and represents the dioperadic composition in a dioperad \cite{gan}, which is sometimes called a polycategory in the literature.

Thus, a category is a special case of an operad, which in turn is a special case of a dioperad.  Moreover, properads contain dioperads.  These generalizations of categories are very powerful tools that encode operations.  For example, the little $n$-cube operads introduced by May \cite{may} provide a recognition principle for connected $n$-fold loop spaces.  Dioperads can model algebraic structures with multiple inputs and multiple outputs, whose axioms are represented by simply connected graphs, such as Lie bialgebras.  Going even further, properads can model algebraic structures with multiple inputs and multiple outputs, whose axioms are represented by connected wheel-free graphs, such as biassociative bialgebras and (co)module bialgebras.  There are numerous other applications of these generalized categories in homotopy theory, string topology, deformation theory, and mathematical physics, among many other subjects.  See, for example, \cite{markl08}, \cite{mss}, and the references therein.

A wheeled properad takes this line of thought one step further by allowing a contraction operation.  To motivate this structure, recall from above that a dioperadic composition is a special case of a properadic composition.  It is natural to ask if properadic compositions can be generated by dioperadic compositions in some way, since the latter are much simpler than the former.  A dioperadic graph has only one internal edge, while a partially grafted corollas can have finitely many internal edges.  Starting with a dioperadic graph, to create a general partially grafted corollas with the same vertices and containing the given internal edge, one needs to connect some output legs of the bottom vertex with some input legs of the top vertex.  So we need a \emph{contraction} operation that connects an output $d_i$ of an operation $h$ with an input $c_j$ of the same color (i.e., $d_i=c_j$) of the same operation.  The following picture, called a contracted corolla, represents such a contraction.
\begin{center}
\begin{tikzpicture}
\matrix[row sep=.8cm,column sep=1cm] {
\node [plain] (p) {$h$}; \\
};
\draw [outputleg] (p) to +(-.5cm,.3cm);
\draw [outputleg] (p) to +(.5cm,.3cm);
\draw [inputleg] (p) to +(-.5cm,-.3cm);
\draw [inputleg] (p) to +(.5cm,-.3cm);
\draw [arrow, looseness=1, in=-45, out=45, loop] (p) to node[near start]{\footnotesize{$d_i$}} node[near end]{\footnotesize{$c_j$}} ();
\end{tikzpicture}
\end{center}
A general properadic composition is a composition of a dioperadic composition and finitely many contractions.

A wheeled properad  is an object that has a bi-equivariant structure, units, a dioperadic composition, and a contraction, satisfying suitable axioms.  Wheeled properads in the linear setting are heavily used in applications \cite{kwz,mms,merkulov,merkulov2,merkulov3}.  Foundational discussion of wheeled properads can be found in \cite{jy2,jy3}.  As the picture above indicates, when working with wheeled properads, one must allow graphs to have loops and directed cycles in general.

\section*{Infinity Categories and Infinity Operads}

Another natural way to extend the notion of a category comes from relaxing the axioms, so they do not need to hold in the strict sense, resulting in what is called a weak category.  There are quite a few variations of this concept; see, for example,  \cite{leinster,simpson}.  We concentrate only on $\infty$-categories in the sense of Joyal \cite{joyal1} and Lurie \cite{lurie}.  These objects were actually first defined by Boardman and Vogt \cite{bv} as simplicial sets satisfying the restricted Kan condition.

The basic idea of an $\infty$-category is very similar to that of the path space of a topological space $X$.  Given two composable paths $f$ and $g$ in $X$, there are many ways to form their composition, but any two such compositions are homotopic.  Moreover, any two such homotopies are themselves homotopic, and so forth.  Composition of paths is not associative in the strict sense, but it holds up to homotopy again.

Similarly, in an $\infty$-category, one asks that there be a composition of two given morphisms.  
\begin{center}
\begin{tikzpicture}
\matrix[row sep=.7cm, column sep=1cm]{
& \node [empty] (one) {$1$}; &\\
\node [empty] (zero) {$0$}; &&
\node [empty] (two) {$2$};\\
};
\draw [arrow] (zero) to node{$f$} (one);
\draw [arrow] (one) to node{$g$} (two);
\draw [dottedarrow] (zero) to node{$\exists$} (two);
\end{tikzpicture}
\end{center}
There may be many such compositions, but any two compositions are homotopic.  Associativity of such composition holds up to homotopy, and there are also higher coherence conditions.  To make these ideas precise, observe that the above triangle can be phrased in a simplicial set $X$.  Namely, $f$ and $g$ are two $1$-simplices in $X$ that determine a unique inner horn $\Lambda^1[2] \to X$, with $g$ as the $0$-face and $f$ as the $2$-face.  To say that a composition exists, one can say that this inner horn has an extension to $\varDelta[2] \to X$, so its $1$-face is such a composition.  In fact, an $\infty$-category is defined as a simplicial set in which every inner horn
\[
\nicexy@C+12pt{
\Lambda^k[n] \ar[d] \ar[r]^-{\forall} & X\\
\varDelta[n] \ar@{.>}[ur]_{\exists} &
}\]
with $0 < k < n$ has a filler.  There are several other models of $\infty$-categories, which are discussed in \cite{bergner07,bergner,jt,rezk}.


Similarly, an $\infty$-operad captures the idea of an up-to-homotopy operad.  It was first developed by Moerdijk and Weiss \cite{mw1,mw2}, and more recently in \cite{bn,cma,cmb,cmc,heuts1,heuts2,hhm,lukacs,nik}.  To define $\infty$-operads, recall that $\infty$-categories are simplicial sets satisfying an inner horn extension property.  The finite ordinal category $\varDelta$ can be represented using linear graphs.  In fact, the object
\[
[n] = \{0 < 1 < \cdots < n\} \in \varDelta
\]
is the category generated by the following linear graph with $n$ vertices.
\begin{center}
\begin{tikzpicture}
\matrix[row sep=1cm,column sep=1.2cm] {
\node [plain] (v1) {$v_1$}; &
\node [plain] (v2) {$v_2$}; &
\node [empty] (vdot) {$\cdots$}; &
\node [plain] (vn) {$v_n$}; &\\
};
\draw [inputleg] (v1) to node[swap]{\footnotesize{$0$}} +(-1.4cm,0cm);
\draw [arrow] (v1) to node{\footnotesize{$1$}} (v2);
\draw [arrow] (v2) to node{\footnotesize{$2$}} (vdot);
\draw [arrow] (vdot) to node{\footnotesize{$n-1$}} (vn);
\draw [outputleg] (vn) to node{\footnotesize{$n$}} +(1.4cm,0cm);
\end{tikzpicture}
\end{center}
Here each vertex $v_i$ is the generating morphism $i-1 \to i$.

Likewise, each unital tree freely generates an operad.  The resulting full subcategory $\varOmega$ of operads generated by unital trees is called the dendroidal category.  Using the dendroidal category $\varOmega$ instead of the finite ordinal category $\varDelta$, one can define analogs of coface and  codegeneracy maps.  Objects in the presheaf category $\omegaset$ are called dendroidal sets, which are tree-like analogs of simplicial sets.  An $\infty$-operad is then defined as a dendroidal set in which every dendroidal inner horn has a filler.

The objects discussed above are summarized in the following table.
\begin{center}
\setlength{\tabcolsep}{4mm}
\begin{tabular}{|c|c|c|}\hline
categories 
& $\sset$
& $\infty$-categories \\
operads
& $\omegaset$
& $\infty$-operads \\
properads
& ? & ?\\
wheeled properads
& ? & ?\\ \hline
\end{tabular}
\end{center}
Starting in the upper left corner, moving downward represents the first type of generalized categories discussed above, while moving to the right yields $\infty$-categories.  In view of this table, a natural question is this:
\begin{quotation}
What are $\infty$-properads and $\infty$-wheeled properads?
\end{quotation}
The answer should simultaneously capture the notion of an up-to-homotopy (wheeled) properad and also extend $\infty$-category and $\infty$-operad.

\section*{Infinity (Wheeled) Properads}

The purpose of this monograph is to initiate the study of $\infty$-properads and $\infty$-wheeled properads.  Let us very briefly describe how $\infty$-properads are defined.

Both $\infty$-categories and $\infty$-operads are objects in some presheaf categories satisfying some inner horn extension property.  In each case, the presheaf category is induced by graphs that parametrize composition of operations and their axioms.  For categories (resp., operads), one uses linear graphs (resp., unital trees), which freely generate categories (resp., operads) that form the finite ordinal category $\varDelta$ (resp., dendroidal category $\varOmega$).

Properadic compositions and their axioms are parametrized by connected graphs without directed cycles, which we call connected wheel-free graphs.  Each connected wheel-free graph freely generates a properad, called a graphical properad.  With carefully defined morphisms called graphical maps, such graphical properads form a non-full subcategory $\varGamma$, called the graphical category, of the category of properads.  There are graphical analogs of coface and codegeneracy maps in the graphical category.  Objects in the presheaf category $\gupcset$ are called graphical sets.  Infinity properads are defined as graphical sets that satisfy an inner horn extension property.

The graphical category $\varGamma$ contains a full subcategory $\varTheta$ corresponding to simply connected graphs, which as mentioned above parametrize the dioperadic composition.  Likewise, by restricting to linear graphs and unital trees, the finite ordinal category $\varDelta$ and the dendroidal category $\varOmega$ may be regarded as full subcategories of the graphical category.  The full subcategory inclusions
\[
\nicexy{
\varDelta \ar[r]^-{i} & \varOmega \ar[r]^{i} & \varTheta \ar[r]^-{i} & \varGamma
}\]
induce restriction functors $i^*$ and fully faithful left adjoints $i_!$,
\[
\nicexy{
\sset \ar@<.4ex>[r]^-{i_!} & 
\omegaset \ar@<.4ex>[l]^-{i^*} \ar@<.4ex>[r]^-{i_!} &
\gupdset \ar@<.4ex>[l]^-{i^*} \ar@<.4ex>[r]^-{i_!} &
\gupcset \ar@<.4ex>[l]^-{i^*},
}\]
on presheaf categories.  Moreover, we have $i^*i_! = \Id$ in each case.  In particular, the graphical set generated by an $\infty$-category or an $\infty$-operad is an $\infty$-properad.

One difficulty in studying $\infty$-properads is that working with connected wheel-free graphs requires a lot more care than working with simply connected graphs, such as linear graphs and unital trees.  For example, even the graphical analogs of the cosimplicial identities are not entirely straightfoward to prove.  Likewise, whereas every object in the finite ordinal category $\varDelta$ and the dendroidal category $\varOmega$ has a finite set of elements, most objects in the graphical category $\varGamma$ have infinite sets of elements.  In fact, a graphical properad is finite if and only if it is generated by a simply connected graph.  Furthermore, since graphical properads are often very large, general properad maps betweem them may exhibit bad behavior that would never happen in $\varDelta$ and $\varOmega$.  To obtain graphical analogs of, say, the epi-mono factorization, we will need to impose suitable restrictions on the maps in the graphical category $\varGamma$.

To explain such differences in graph theoretical terms, note that in a simply connected graph, any two vertices are connected by a unique internal path.  In particular, any finite subset of edges can be shrunk away in any order.  This is far from the case in a connected wheel-free graph, where there may be finitely many edges adjacent to two given vertices.  Thus, one cannot in general shrink away an edge in a connected wheel-free graph, or else one may end up with a directed loop.  Therefore, in the development of the graphical category $\varGamma$, one must be extremely careful about the various graph operations.

Infinity wheeled properads are defined similarly, using all connected graphs, where loops and directed cycles are allowed, instead of connected wheel-free graphs.  One major difference between $\infty$-wheeled properads and $\infty$-properads is that, when working with all connected graphs, there are more types of coface maps.  This is due to the fact that a properad has only one generating operation, namely the properadic composition, besides the bi-equivariant structure and the units.  In a wheeled properad, there are two generating operations, namely the dioperadic composition and the contraction, each with its own corresponding coface maps.

\section*{Chapter Summaries}

This monograph is divided into two parts.  In part 1 we set up the graph theoretic foundation and discuss $\infty$-properads.  Part 2 has the parallel theory of $\infty$-wheeled properads.

In chapter \ref{ch:graph} we recall the definition of a graph and the operation of graph substitution developed in detail in the monograph \cite{jy2}.  This discussion is needed because graphical properads (resp., graphical wheeled properads) are generated by connected wheel-free graphs (resp., all connected graphs).  Also, many properties of the graphical category are proved using the associativity and unity of graph substitution. The first four sections of this chapter are adapted from that monograph.  In the remaining sections, we develop graph theoretical concepts that will be needed later to define coface and codegeneracy maps in the graphical categories for connected (wheeled-free) graphs.  We give graph substitution characterization of each of these concepts.

A major advantage of using graph substitution to define coface and codegeneracy maps is that the resulting definition is very formal.  From the graph substitution point of view, regardless of the graphs one uses (linear graphs, unital trees, simply connected graphs, or connected (wheel-free) graphs), all the coface maps have almost the same definition.  What changes in these definitions of coface maps, from one type of graphs to another, is the set of minimal generating graphs.  For linear graphs, a minimal generating graph is the linear graph with two vertices.  For unital trees (resp., simply connected graphs), one takes as minimal generating graphs the unital trees with one ordinary edge (resp., dioperadic graphs).  For connected wheel-free graphs, the partially grafted corollas form the set of minimal generating graphs.  For all connected graphs, the set of minimal generating graphs consists of the dioperadic graphs and the contracted corollas.

In chapter \ref{ch:properads} we recall both the biased and the unbiased definitions of a properad.  The former describes a properad in terms of generating operations, namely, units, $\Sigma$-bimodule structure, and properadic composition.  The unbiased definition of a properad describes it as an algebra over a monad induced by connected wheel-free graphs.  The equivalence of the two definitions of a properad are proved in detail in \cite{jy2} as an example of a general theory of generating sets for graphs.

In chapter  \ref{ch:tensor} we equip the category of properads with a symmetric monoidal closed structure.   For topological operads, a symmetric monoidal product was already defined by Boardman and Vogt \cite{bv}.  One main result of this chapter gives a simple description of the tensor product of two free properads in terms of the two generating sets.  In particular, when the free properads are finitely generated, their tensor product is finitely presented.  This is not immediately obvious from the definition because free properads are often infinite sets. 
In future work we will compare our symmetric monoidal tensor product with the Dwyer-Hess box product of operadic bimodules \cite{dh13} and show they agree in the special case when our properads come from operads.

In chapter \ref{ch:grproperad} we define graphical properads as the free properads generated by connected wheel-free graphs.  We observe that a graphical properad has an infinite set of elements precisely when the generating graph is not simply connected.  The discussion of the tensor product of free properads in chapter  \ref{ch:tensor}  applies in particular to graphical properads.  Then we illustrate with several examples that a general properad map between graphical properads may exhibit bad behavior that would never happen when working with simply connected graphs.  These examples serve as the motivation of the restriction on the morphisms in the graphical category.

In chapter \ref{ch:mapgrproperad} we define the properadic graphical category $\varGamma$, whose objects are graphical properads.  Its morphisms are called properadic graphical maps.  To define such graphical maps, we first discuss coface and codegeneracy maps between graphical properads.  We establish graphical analogs of the cosimplicial identities.  The most interesting case is the graphical analog of the cosimplicial identity
\[
d^j d^i = d^i d^{j-1}
\]
for $i<j$ because it involves iterating the operations of deleting an almost isolated vertex and of smashing two closest neighbors together.  Graphical maps do not have the bad behavior discussed in the examples in chapter \ref{ch:grproperad}.  In particular, it is observed that each graphical map has a factorization into codegeneracy maps followed by coface maps.  Such factorizations do not exist for general properad maps between graphical properads.
Finally, we show that the properadic graphical category admits the structure of a (dualizable) generalized Reedy category, in the sense of \cite{reedyextension}.


In chapter \ref{ch:graphicalset} we first define the category $\gupcset$ of graphical sets.  There is an adjoint pair
\[
\nicexy@C+10pt{
\gupcset \ar@<.4ex>[r]^-{L} & \properad \ar@<.4ex>[l]^-{N},
}\]
in which the right adjoint $N$ is called the properadic nerve.  The symmetric monoidal product of properads in chapter  \ref{ch:tensor} induces, via the properadic nerve, a symmetric monoidal closed structure on $\gupcset$.  Then we define an  $\infty$-properad as a graphical set in which every inner horn has a filler.  If, furthermore, every inner horn filler is unique, then it is called a strict $\infty$-properad.  The rest of this chapter contains two alternative descriptions of a strict $\infty$-properad.  One description is in terms of the graphical analogs of the Segal maps, and the other is in terms of the properadic nerve.

In chapter \ref{ch:fundprop} we give an explicit description of the fundamental properad $L\sK$ of an $\infty$-properad $\sK$.  The fundamental properad of an $\infty$-properad consists of homotopy classes of $1$-dimensional elements.  It takes a bit of work to prove that there is a well-defined homotopy relation among $1$-dimensional elements and that a properad structure can be defined on homotopy classes.  This finishes part 1 on $\infty$-properads.

Part 2 begins with chapter \ref{ch:wproperads}.  We first recall from \cite{jy2} the biased and the unbiased definitions of a wheeled properad.  There is a symmetric monoidal structure on the category of wheeled properads.  Then we define graphical wheeled properads as free wheeled properads generated by connected graphs, possibly with loops and directed cycles.  With the exception of the exceptional wheel, a graphical wheeled properad has a finite set of elements precisely when the generating graph is simply connected.  So most graphical wheeled properads are infinite.  In the rest of this chapter, we discuss wheeled versions of coface maps, codegeneracy maps, and graphical maps, which are used to define the wheeled properadic graphical category $\Gammaw$.  Every wheeled properadic graphical map has a decomposition into codegeneracy maps followed by coface maps. 

In chapter \ref{ch:infwproperad} we define the adjunction
\[
\nicexy@C+10pt{
\gwheelcset \ar@<2pt>[r]^-{L} & \wproperad \ar@<2pt>[l]^-{N}}
\]
between wheeled properads and wheeled properadic graphical sets.  Then we define $\infty$-wheeled properads as wheeled properadic graphical sets that satisfy an inner horn extension property.  Next we give two alternative characterizations of \emph{strict} $\infty$-wheeled properads, one in terms of the wheeled properadic Segal maps, and the other in terms of the wheeled properadic nerve.  In the last section, we give an explicit description of the fundamental wheeled properad $L\sK$ of an $\infty$-wheeled properad $\sK$ in terms of homotopy classes of $1$-dimensional elements.

In the final chapter, we mention some potential future applications and extensions of this work.
These include applications to string topology and deformation theory.
We also discuss other models for $\infty$-properads and the categorical machinery of M. Weber which produces categories like $\varDelta$ and $\varOmega$.


\part{Infinity Properads}

\chapter{Graphs}
\label{ch:graph}

\abstract*{We recall the definition of a graph and the operation of graph substitution developed in detail in the monograph \cite{jy2}.  This discussion is needed because graphical properads (resp., graphical wheeled properads) are generated by connected wheel-free graphs (resp., all connected graphs).  Also, many properties of the graphical category are proved using the associativity and unity of graph substitution. The first four sections of this chapter are adapted from that monograph.  In the remaining sections, we develop graph theoretical concepts that will be needed later to define coface and codegeneracy maps in the graphical categories for connected (wheeled-free) graphs.  We give graph substitution characterization of each of these concepts.}

The purpose of this chapter is to discuss graphs and graph operations that are suitable for the study of $\infty$-properads and $\infty$-wheeled properads.  We will need connected (wheel-free) graphs for two reasons.  First, the free (wheeled) properad monad is a coproduct parametrized by suitable subsets of connected wheel-free (connected) graphs.  Second, each connected wheel-free (resp., connected) graph generates a free properad (resp., wheeled properad), called a \emph{graphical (wheeled) properad}.  The graphical (wheeled) properads form the graphical category, which in turn yields graphical sets and $\infty$-(wheeled) properads.

 In section \ref{sec:wheeledgraph} we discuss graphs in general.  In section \ref{sec:exgraph} we provide some examples of graphs, most of which are important later.  In section \ref{sec:cwfree} we discuss connected (wheel-free) graphs as well as some important subsets of graphs, including simply connected graphs, unital trees, and linear graphs.  In section \ref{sec:graphsub} we discuss the operation of graph substitution, which induces the multiplication of the free (wheeled) properad monad.  Graph substitution will also be used in the second half of this chapter to characterize some graph theoretic concepts that will be used to define coface maps in the graphical category in section \ref{sec:coface}.  The reference for sections \ref{sec:wheeledgraph}--\ref{sec:graphsub} on graphs and graph substitution is \cite{jy2}.

In sections \ref{sec:closestnbd} and \ref{sec:almostisolated} we discuss \emph{closest neighbors} and \emph{almost isolated vertices} in connected wheel-free graphs.  These concepts are needed when we define inner and outer coface maps in the graphical category for connected wheel-free graphs.  In sections \ref{sec:deletable} and \ref{sec:disconnectable} we discuss analogous graph theoretic concepts that will be used in Part 2 to define inner and outer coface maps in the graphical category for connected graphs.

\section{Wheeled Graphs}
\label{sec:wheeledgraph}

The intuitive idea of a graph is quite simple, but the precise definition is slightly abstract.  We first give the more general definition of a wheeled graph before restricting to the connected (wheel-free) ones.

\subsection{Profiles of Colors}

We begin by describing profiles of colors, which provide a way to parametrize lists of inputs and outputs.  This discussion of profiles of colors can be found in \cite{jy1,jy2}.

Fix a set $\fC$.\label{note:colors}  

\begin{definition}
\label{def:colors}
Let $\Sigma_n$ denote the symmetric group on $n$ letters.
\begin{enumerate}
\item
An element in $\fC$\index{color} will be called a \textbf{color}.
\item
A \textbf{$\fC$-profile of length $n$}\index{profile} is a finite sequence
\[
\uc = (c_1, \ldots , c_n) = c_{[1,n]}\label{note:profile}
\]
of colors.  We write $|\uc| = n$ for the length. The \textbf{empty profile},\index{profile!empty} with $n=0$, is denoted by $\varnothing$.
\item
Given a $\fC$-profile $\uc = c_{[1,n]}$ and $0 < k \leq n$, a \textbf{$k$-segment}\index{profile!segment} of $\uc$ is a sub-$\fC$-profile
\[
\uc' = (c_i,\ldots , c_{i+k-1})
\]
of length $k$ for some $1 \leq i \leq n+1-k$.
\item
Given two $\fC$-profiles $\uc$ and $\ud = d_{[1,m]}$ and a $k$-segment $\uc' \subseteq \uc$ as above, define the $\fC$-profile
\[
\uc \circ_{\uc'} \ud = \left(c_1, \ldots , c_{i-1}, d_1, \ldots , d_m , c_{i+k}, \ldots , c_n\right).\label{note:subprofile}
\]
If $\uc'$ happens to be the $1$-segment $(c_i)$, we also write $\uc \compi \ud$ for $\uc \circ_{\uc'} \ud$.
\item
For a $\fC$-profile $\uc$ of length $n$ and $\sigma \in \Sigma_n$, define the left and right actions
\[
\sigma\uc = \left(c_{\sigma(1)}, \ldots , c_{\sigma(n)}\right) \andspace 
\uc\sigma = \left(c_{\sigma^{-1}(1)}, \ldots , c_{\sigma^{-1}(n)}\right).
\]
\item
The groupoid of all $\fC$-profiles\index{profile!groupoid of} with left (resp., right) symmetric group actions as morphisms is denoted by $\pofc$\label{note:pc} (resp., $\pofcop$).
\item
Define the product category
\[
\SC = \ptwoc,\label{note:sc}
\]
which will be abbreviated to $\sS$ if $\fC$ is clear from the context.  Its elements are pairs of $\fC$-profiles and are written either horizontally as $\dch$ or vertically as $\dc$.\label{note:dc}
\end{enumerate}
\end{definition}

\subsection{Generalized Graphs}

Fix an infinite set $\fF$ once and for all.

\begin{definition}
\label{def:graph}
A \textbf{generalized graph}\index{graph} $G$ is a finite set $\Flag(G) \subset \fF$ with
\begin{itemize}
\item
a partition $\Flag(G)\label{note:flag} = \coprod_{\alpha \in A} F_\alpha$ with $A$ finite,
\item
a distinguished partition subset $F_\epsilon$ called the \textbf{exceptional cell},\index{exceptional cell}
\item
an involution $\iota$ satisfying $\iota F_\epsilon \subseteq F_\epsilon$, and
\item
a free involution $\pi$ on the set of $\iota$-fixed points in $F_\epsilon$.
\end{itemize}
\end{definition}

Next we introduce some intuitive terminology associated to a generalized graph.

\begin{definition}
\label{def:edges}
Suppose $G$ is a generalized graph.
\begin{enumerate}
\item
The elements in $\Flag(G)$ are called \textbf{flags}.\index{flag}  Flags in a non-exceptional cell are called \textbf{ordinary flags}.\index{flag!ordinary}  Flags in the exceptional cell $F_\epsilon$ are called \textbf{exceptional flags}.\index{flag!exception}
\item
Call $G$ an \textbf{ordinary graph}\index{graph!ordinary} if its exceptional cell is empty.
\item
Each non-exceptional partition subset $F_\alpha \not= F_\epsilon$ is a \textbf{vertex}.\index{vertex}  The set of vertices is denoted by $\vertex(G)$.\label{note:vt}  An empty vertex is an \textbf{isolated vertex},\index{vertex!isolated} which is often written as $C_{\emptyprofh}$.  A flag in a vertex is said to be \textbf{adjacent to}\index{adjacent} or \textbf{attached to}\index{attached} that vertex.
\item
An $\iota$-fixed point is a \textbf{leg}\index{leg} of $G$.  The set of legs of $G$ is denoted by $\Leg(G)$.\label{note:leg}  An \textbf{ordinary leg} \index{leg!ordinary} (resp., \textbf{exceptional leg})\index{leg!exceptional} is an ordinary (resp., exceptional) flag that is also a leg. For an $\iota$-fixed point $x \in F_\epsilon$, the pair $\{x,\pi x\}$ is an \textbf{exceptional edge}.\index{exceptional edge} \index{edge!exceptional}
\item
A $2$-cycle of $\iota$ consisting of ordinary flags is an \textbf{ordinary edge}.\index{edge!ordinary} \index{ordinary edge}   A $2$-cycle of $\iota$ contained in a vertex $v$ is a \textbf{loop}\index{loop} at $v$.  A vertex that does not contain any loop is \textbf{loop-free}.\index{loop-free}   A $2$-cycle of $\iota$ in the exceptional cell is an \textbf{exceptional loop}.\index{exceptional loop} \index{loop!exceptional}  
\item
An \textbf{internal edge}\index{internal edge} \index{edge!internal} is a $2$-cycle of $\iota$, i.e., either an ordinary edge or an exceptional loop.  An \textbf{edge}\index{edge} means an internal edge, an exceptional edge, or an ordinary leg.  The set of edges (resp., internal edges) in $G$ is denoted by $\edge(G)$ (resp., $\edgei(G)$).\label{note:edge}
\item
An ordinary edge $e = \{e_{-1},e_1\}$ is said to be \textbf{adjacent to} or \textbf{attached to} a vertex $v$ if either (or both) $e_i \in v$.
\end{enumerate}
\end{definition}

\begin{remark}
In plain language, flags are half-edges. The elements in a vertex are the flags attached to it.  On the other hand, exceptional flags are \emph{not} attached to any vertex.  So an exceptional edge (resp., exceptional loop) is an edge (resp., a loop) not attached to any vertex.  A leg has at most one end attached to a vertex.  Also note that an exceptional edge is \emph{not} an internal edge.
\end{remark}

\subsection{Structures on Generalized Graphs}

To describe free (wheeled) properads as well as certain maps in the graphical category later, we need some extra structures on a generalized graph, which we now discuss.  Intuitively, we need an orientation for each edge, a color for each edge, and a labeling of the incoming/outgoing flags of each vertex as well as the generalized graph. Fix a set of colors $\fC$.

\begin{definition}
Suppose $G$ is a generalized graph.
\begin{enumerate}
\item
A \textbf{coloring}\index{coloring} \index{graph!coloring} for $G$ is a function
\[\nicexy{
\Flag(G) \ar[r]^-{\kappa} & \fC
}\]
that is constant on orbits of both involutions $\iota$ and $\pi$.
\item
A \textbf{direction}\index{direction} \index{graph!direction} for $G$ is a function
\[\nicexy{
\Flag(G) \ar[r]^-{\delta} & \{-1,1\}
}\]
such that
\begin{itemize}
\item
if $\iota x \not= x$, then $\delta(\iota x) = -\delta(x)$, and
\item
if $x \in F_\epsilon$, then $\delta(\pi x) = -\delta(x)$.
\end{itemize}
\item
For $G$ with direction, an \textbf{input}\index{input} \index{vertex!input} (resp., \textbf{output})\index{output} \index{vertex!output} of a vertex $v$ is a flag $x \in v$ such that $\delta(x) = 1$ (resp., $\delta(x) = -1$).  An \textbf{input}\index{graph!input} (resp., \textbf{output})\index{graph!output} of the graph $G$ is a leg $x$ such that $\delta(x) = 1$ (resp., $\delta(x) = -1$).  For $u \in \{G\} \cup \vertex(G)$, the set of inputs (resp., outputs) of $u$ is written as $\inp(u)$ (resp., $\out(u)$). \label{note:in}
\item
A \textbf{listing}\index{listing} \index{graph!listing} for $G$ with direction is a choice for each $u \in \{G\} \cup \vertex(G)$ of a bijection of pairs of sets
\[\nicexy{
\left(\inp(u), \out(u)\right)
\ar[r]^-{\ell_u} &
\left(\{1,\ldots,|\inp(u)|\}, \{1, \ldots , |\out(u)|\}\right),
}\]
where for a finite set $T$ the symbol $|T|$ denotes its cardinality.
\item
A \textbf{$\fC$-colored wheeled graph}, or just a \textbf{wheeled graph},\index{graph!wheeled} is a generalized graph together with a choice of a coloring, a direction, and a listing.
\end{enumerate}
\end{definition}

\begin{definition}
Suppose $G$ is a wheeled graph.
\begin{enumerate}
\item
Its \textbf{input profile} (resp., \textbf{output profile}) is the $\fC$-profile $\kappa(\inp(u))$ (resp., $\kappa(\out(u))$), where $\inp(u)$ and $\out(u)$ are regarded as ordered sets using the listing.  A \textbf{$\dch$-wheeled graph} is a wheeled graph with input profile $\uc$ and output profile $\ud$.
\item
Suppose $e = \{e_1,e_{-1}\}$ is an ordinary edge attached to a vertex $v$ with $\delta(e_i) = i$.  If $e_{-1} \in v$ (resp., $e_1 \in v$), then $v$ is called the \textbf{initial vertex} \index{initial vertex} \index{vertex!initial} (resp., \textbf{terminal vertex}) \index{terminal vertex} \index{vertex!terminal} of $e$.  If $e$ has initial vertex $u$ and terminal vertex $v$, then it is also denoted by $\nicexy{u \ar[r]^-{e} & v}$.
\end{enumerate}
\end{definition}

\begin{remark}
It \emph{is} possible that a vertex $v$ is both the initial vertex and the terminal vertex of an ordinary edge $e$, which is then a loop at the vertex $v$, as in the following picture.
\begin{center}
\begin{tikzpicture}
\matrix[row sep=.5cm, column sep=1cm]{
\node [plain] (v) {$v$};\\
};
\draw [arrow, out=45, in=-45, loop] (v) to node{$e$} ();
\end{tikzpicture}
\end{center}
\end{remark}

\begin{definition}
A \textbf{strict isomorphism}\index{strict isomorphism!graph} between two wheeled graphs is a bijection of partitioned sets preserving the exceptional cells, both involutions, the colorings, the directions, and the listings.
\end{definition}

\begin{convention}
\label{graphconvention}
In what follows, we will mostly be talking about strict isomorphism classes of wheeled graphs, and we will lazily call them \textbf{graphs}.  Furthermore, we will sometimes ignore the listings, since they can always be dealt with using input/output relabeling permutations, but writing all of them down explicitly tends to obscure the simplicity of the ideas and constructions involved.
\end{convention}

\section{Examples of Graphs}
\label{sec:exgraph}

\begin{example}
\label{ex:emptygraph}
The \textbf{empty graph} \index{graph!empty} \index{empty graph} $\varnothing$ has
\[
\Flag(\varnothing) = \varnothing = \coprod \varnothing,
\]
whose exceptional cell is $\varnothing$, and it has no non-exceptional partition subsets.  In particular, it has no vertices and no flags.
\end{example}

\begin{example}
\label{ex:isolatedvt}
Suppose $n$ is a positive integer.  The \textbf{union of $n$ isolated vertices} \index{isolated vertex} is the graph $V_n$ with
\[
\Flag(V_n) = \varnothing = \coprod_{i=1}^{n+1} \varnothing.
\]
It has an empty set of flags, an empty exceptional cell, and $n$ empty non-exceptional partition subsets, each of which is an isolated vertex.  For example, we can represent $V_3$ pictorially as
\[
\nicexy{
\bullet & \bullet & \bullet
}\]
with each $\bullet$\label{note:isovt} representing an isolated vertex.
\end{example}

\begin{example}
\label{ex:exceptionaledge}
Pick a color $c \in \fC$.  The \textbf{$c$-colored exceptional edge} \index{exceptional edge} is the graph $G$ whose only partition subset is the exceptional cell
\[
\Flag(G) = F_\epsilon = \{f_1, f_{-1}\},
\]
with
\[
\iota(f_i) = f_i, \quad \kappa(f_i) = c, \andspace \delta(f_i) = i.
\]
It can be represented pictorially as
\begin{center}\label{note:uparrow}
\begin{tikzpicture}
\matrix[row sep=1cm,column sep=1cm] {
\node [empty] (t) {}; \\
\node [empty] (b) {};\\
};
\draw [arrow] (b) to node[swap,near start]{$c$} (t);
\end{tikzpicture}
\end{center}
in which the top (resp., bottom) half is $f_{-1}$ (resp., $f_1$).  Note that this graph has no vertices and has one exceptional edge.  The $c$-colored exceptional edge will be referred to in Remark \ref{rk:visualizeunits}.
\end{example}

\begin{example}
\label{ex:exceptionalloop}
The \textbf{$c$-colored exceptional loop} \index{exceptional loop} is defined exactly like the exceptional edge $\uparrow_c$, except for
\[
\iota(f_i) = f_{-i}.
\]
It can be represented pictorially as
\[\label{note:wheel}
\wheelc
\]
in which the left (resp., right) half is $f_{-1}$ (resp., $f_1$).  This graph has no vertices and has one exceptional loop.
\end{example}

\begin{example}
\label{ex:corolla}
Suppose $\uc = c_{[1,m]}$ and $\ud = d_{[1,n]}$ are $\fC$-profiles.  The \textbf{$\dch$-corolla} \index{corolla} $C_{\dch}$\label{note:corolla} is the $\dch$-wheeled graph with
\[
\Flag\left(C_{\dch}\right) = \left\{i_1, \ldots , i_m, o_1, \ldots , o_n\right\}.
\]
Its only non-exceptional partition subset is $v = \Flag(G)$, which is its only vertex, and its exceptional cell is empty.  The structure maps are defined as:
\begin{itemize}
\item
$\iota(i_k) = i_k$ and $\iota(o_j) = o_j$ for all $k$ and $j$.
\item
$\kappa(i_k) = c_k$ and $\kappa(o_j) = d_j$.
\item
$\delta(i_k) = 1$ and $\delta(o_j) = -1$.
\item
$\ell_u(i_k) = k$ and $\ell_u(o_j) = o_j$ for $u \in \{C_{\dch}, v\}$.
\end{itemize}
The $\dch$-corolla can be represented pictorially as the following graph. 
\begin{center}
\begin{tikzpicture}
\matrix[row sep=1cm,column sep=1cm] {
\node [plain,label=above:$...$,label=below:$...$] (p) {$v$}; \\
};
\draw [outputleg] (p) to node[above left=.1cm]{$d_1$} +(-.6cm,.5cm);
\draw [outputleg] (p) to node[above right=.1cm]{$d_n$} +(.6cm,.5cm);
\draw [inputleg] (p) to node[below left=.1cm]{$c_1$} +(-.6cm,-.5cm);
\draw [inputleg] (p) to node[below right=.1cm]{$c_m$} +(.6cm,-.5cm);
\end{tikzpicture}
\end{center}
This corolla will be referred to in Remark \ref{rk:coloredob}.  In particular, when $\uc = \ud = \varnothing$, the corolla $C_{\emptyprofh}$ consists of a single isolated vertex.  Notice that we omit some of the structure of a wheeled graph in its pictorial representation in order not to overload it with too many symbols.

In the $1$-colored case with $\fC = \{*\}$, we also write $C_{\dch}$ as $C_{(m;n)}$.  In other words, $C_{(m;n)}$ is the $1$-colored corolla with $m$ inputs and $n$ outputs.
\end{example}

\begin{example}
\label{ex:permutedcor}
Suppose $\uc = c_{[1,m]}$ and $\ud = d_{[1,n]}$ are $\fC$-profiles, $\sigma \in \Sigma_n$, and $\tau \in \Sigma_m$.  Define the \textbf{permuted corolla} $\sigma C_{\dch} \tau$,\label{note:permcor} \index{permuted corolla} \index{corolla!permuted} which is a $(\uc\tau;\sigma\ud)$-wheeled graph, with
\[
\Flag\left(\sigma C_{\dch} \tau\right) = \Flag\left(C_{\dch}\right).
\]
All the structure maps are the same as for the corolla $C_{\dch}$, except for the listing of the full graph, which in this case is
\[
\ell_G(i_k) = \tau(k) \andspace \ell_G(o_j) = \sigma^{-1}(j).
\]
The $\dch$-corolla is also a permuted corolla with $\sigma$ and $\tau$ the identity permutations.  In the unbiased formulation of a properad, permuted corollas give rise to the bi-equivariant structure.  Moreover, permuted corollas can be used to change the listing of a graph via graph substitution.  Substituting a graph into a suitable permuted corolla changes the listing of the full graph.  On the other hand, substituting a suitable permuted corolla into a vertex changes the listing of that vertex.  The permuted corollas generate the $\Sigma$-bimodule structures in all variants of generalized PROPs \cite{jy2}.
\end{example}

\begin{example}
\label{ex:contractedcor}
Suppose $\uc = c_{[1,m]}$, $\ud = d_{[1,n]}$, $1 \leq j \leq m$, and $1 \leq i \leq n$.  The \textbf{contracted corolla} $\xiij C_{\dch}$\label{note:contcor} \index{contracted corolla} \index{corolla!contracted} is defined exactly like the $\dch$-corolla (Example \ref{ex:corolla}) with the following two exceptions.
\begin{enumerate}
\item
The involution $\iota$ is given by
\[
\iota(i_j) = o_i, \quad 
\iota(o_i) = i_j,
\]
and is the identity on other flags.
\item
The listing for the full graph is given by
\[
\ell_G(i_k) =
\begin{cases}
k & \text{ if $k<j$},\\
k-1 & \text{ if $k>j$}
\end{cases}
\]
and
\[
\ell_G(o_k) =
\begin{cases}
k & \text{ if $k<i$},\\
k-1 & \text{ if $k>i$}.
\end{cases}
\]
\end{enumerate}
This contracted corolla can be represented pictorially as follows.
\begin{center}
\begin{tikzpicture}
\matrix[row sep=1.5cm,column sep=.8cm] {
& \node [emptyvt] (a) {}; &\\
\node [plain] (v) {$v$}; &&
\node [emptyvt] (b) {}; \\
& \node [emptyvt] (c) {}; &\\
};
\draw [outputleg] (v) to node[above left=.1cm]{$d_1$} +(-.6cm,.5cm);
\draw [outputleg] (v) to node[above right=.1cm]{$d_n$} +(.6cm,.5cm);
\draw [inputleg] (v) to node[below left=.1cm]{$c_1$} +(-.6cm,-.5cm);
\draw [inputleg] (v) to node[below right=.1cm]{$c_m$} +(.6cm,-.5cm);
\draw [->,thick] (v) 
to [out=90,in=180] node{$d_i$} (a)
to [out=0,in=90] (b) 
to [out=-90,in=0] (c)
to [out=180,in=-90] node{$c_j$} (v);
\end{tikzpicture}
\end{center}
Its input/output profiles are
\[
\left(\uc \setminus \{c_j\}; \ud \setminus \{d_i\}\right)
\]
and has $1$ internal edge.  On the other hand, the incoming/outgoing profiles of the vertex $v$ are still $\dch$.  If the internal edge is named $e$, then we also write
\[
\xi_e C_{\dch} = \xiij C_{\dch}.
\]
The contracted corollas generate the contraction in wheeled properads (Definition \ref{def:wproperad}).
\end{example}

\begin{example}
\label{ex:pgcor}
Suppose $\ua = a_{[1,k]}$, $\ub = b_{[1,l]}$, $\uc = c_{[1,m]}$, $\ud = d_{[1,n]}$ are $\fC$-profiles, $\uc \supseteq \uc' = \ub' \subseteq \ub$ are equal $\alpha$-segments for some $\alpha > 0$ with
\[
\ub' = (b_i, \ldots , b_{i+\alpha-1}) \andspace \uc' = (c_j, \ldots , c_{j+\alpha-1}).
\]
Then the \textbf{partially grafted corollas}\label{note:pgcor} \index{partially grafted corollas} \index{corolla!partially grafted}
\[
G = C_{\dch} \boxtimes^{\uc'}_{\ub'} C_{\bah}
\]
is defined as follows.
\begin{itemize}
\item
$\Flag(G) = \left\{f_{a_p}, f_{b_q}, f_{c_r}, f_{d_s}\right\}$ with $1 \leq p \leq k$, $1 \leq q \leq l$, $1 \leq r \leq m$, and $1 \leq s \leq n$.
\item
The exceptional cell is empty.
\item
There are two vertices,
\begin{itemize}
\item
$v = \left\{f_{c_r}, f_{d_s}\right\}$ with  $1 \leq r \leq m$ and $1 \leq s \leq n$, and
\item
$u = \left\{f_{a_p}, f_{b_q}\right\}$ with  $1 \leq p \leq k$ and $1 \leq q \leq l$.
\end{itemize}
\item
$\iota$ fixes the flags with subscripts in $\ua$, $\ud$, $\ub \setminus \ub' = (b_1, \ldots , b_{i-1}, b_{i+\alpha}, \ldots , b_l)$, and $\uc \setminus \uc' = (c_1, \ldots , c_{j-1}, c_{j+\alpha}, \ldots , c_m)$.
\item
$\iota\left(b_{i+t}\right) = c_{j+t}$ and $\iota\left(c_{j+t}\right) = b_{i+t}$ for $0 \leq t \leq \alpha-1$.
\item
The coloring is defined as
\[
\kappa\left(f_{e_x}\right) = e_x
\]
for each possible subscript $e_x$.
\item
Flags with subscripts in $\ua$ and $\uc$ have $\delta = 1$, while flags with subscripts in $\ub$ and $\ud$ have $\delta = -1$.
\item
At each vertex $w \in \{u,v\}$, the listing is given by
\[
\ell_w\left(f_{e_x}\right) = x
\]
for $e \in \{a,b,c,d\}$.
\item
The listing for the full graph $G$ at its inputs is given by
\[
\ell_G\left(f_{e_x}\right) = 
\begin{cases}
y & \text{ if $e_x = c_y$ with $1 \leq y < j$},\\
y+j-1 & \text{ if $e_x = a_y$ for $1 \leq y \leq k$},\\
y-\alpha+k & \text{ if $e_x = c_y$ with $j+\alpha \leq y \leq m$}.
\end{cases}
\]
\item
The listing for the full graph $G$ at its outputs is given by
\[
\ell_G\left(f_{e_x}\right) = 
\begin{cases}
y & \text{ if $e_x = b_y$ with $1 \leq y < i$},\\
y + i-1 & \text{ if $e_x = d_y$ with $1 \leq y \leq n$},\\
y-\alpha+n & \text{ if $e_x = b_y$ with $i+\alpha \leq y \leq l$}.
\end{cases}
\]
\end{itemize}
This partially grafted corollas can be represented pictorially as the following graph. 
\begin{center}
\begin{tikzpicture}
\matrix[row sep=.2cm,column sep=1.5cm] {
\node [plain,label=below:$\uc'$,label=above:$\ud$] (p1) {$v$}; \\
\node [empty]  {...}; \\
\node [plain,label=below:$\ua$,label=above:$\ub'$] (p2) {$u$}; \\
};
\foreach \x in {1,2}
{
\draw [outputleg] (p\x) to +(-.8cm,.6cm);
\draw [outputleg] (p\x) to +(.8cm,.6cm);
\draw [inputleg] (p\x) to +(-.8cm,-.6cm);
\draw [inputleg] (p\x) to +(.8cm,-.6cm);
}
\draw [arrow,bend left=40] (p2) to (p1);
\draw [arrow,bend right=40] (p2) to (p1);
\end{tikzpicture}
\end{center}
This is the same graph as the one in Remark \ref{rk:visualizeunits}, which explains why we use the same symbol $\boxtimes^{\uc'}_{\ub'}$ for both the properadic composition and the partially grafted corollas.  Note that the input and output profiles are
\[
\inp(G) = \uc \circ_{\uc'} \ua \andspace \out(G) = \ub \circ_{\ub'} \ud,
\]
which form the input/output profiles of the target of the properadic composition.  There are $\alpha$ internal edges.  If these internal edges are names $\be = (e_1,\ldots,e_{\alpha})$, then we also write
\[
C_{\dch} \boxtimes_{\be} C_{\bah} = 
C_{\dch} \boxtimes^{\uc'}_{\ub'} C_{\bah}.
\]
The partially grafted corollas generate the properadic composition (Definition \ref{def:biasedproperad}).
\end{example}

\begin{example}
\label{ex:dioperadic}
Using the same symbols as in Example \ref{ex:pgcor}, suppose that $\alpha=1$, i.e., $\ub' = (b_i) = (c_j) = \uc'$ are equal $1$-segments.  Define the \textbf{dioperadic graph}\label{note:dioperadic} \index{dioperadic graph} \index{graph!dioperadic}
\[
C_{\dch} \icompj C_{\bah} 
= C_{\dch} \boxtimes^{(c_j)}_{(b_i)} C_{\bah}
\]
as a special case of a partially grafted corollas.  This dioperadic graph can be represented pictorially as follows.
\begin{center}
\begin{tikzpicture}
\matrix[row sep=2cm,column sep=1cm] {
\node [plain,label=above:$\ud$] (p1) {$v$}; \\
\node [plain,label=below:$\ua$] (p2) {$u$}; \\
};
\foreach \x in {1,2}
{
\draw [outputleg] (p\x) to +(-.7cm,.5cm);
\draw [outputleg] (p\x) to +(.7cm,.5cm);
\draw [inputleg] (p\x) to +(-.7cm,-.5cm);
\draw [inputleg] (p\x) to +(.7cm,-.5cm);
}
\draw [arrow] (p2) to node[swap,near end]{\footnotesize{$c_j$}} node[swap,near start]{\footnotesize{$b_i$}} (p1);
\end{tikzpicture}
\end{center}
It has input/output profiles
\[
\left(\uc \circ_j \ua; \ub \circ_i \ud\right)
\]
and  $1$ internal edge.  If this internal edge is named $e$, then we also write
\[
C_{\dch} \circ_e C_{\bah} = C_{\dch} \icompj C_{\bah}
\]
The dioperadic graphs generate the dioperadic composition in a dioperad \cite{gan} and in a wheeled properad (Definition \ref{def:wproperad}).
\end{example}

\section{Connected Graphs}
\label{sec:cwfree}

In this section, we discuss connected graphs and connected wheel-free graphs.  They will be used to define free (wheeled) properads and the graphical category later.  We will also define simply connected graphs, unital trees, and linear graphs. To discuss connectivity and wheels, we need the concept of a path.

\subsection{Paths}

\begin{definition}
\label{def:path}
Suppose $G$ is a graph.
\begin{enumerate}
\item
A \textbf{path} \index{path} in $G$ is a pair
\[
P = \left(\left(e^j\right)_{j=1}^r, \left(v_i\right)_{i=0}^r\right)
\]
with $r \geq 0$, in which
\begin{itemize}
\item
the $v_i$ are distinct vertices except possibly for $v_0 = v_r$,
\item
the $e^j$ are distinct ordinary edges, and
\item
each $e^j$ is adjacent to both $v_{j-1}$ and $v_j$.
\end{itemize}
Such a path is said to have \textbf{length} \index{path!length} $r$.
\item
A path of length $0$ is called a \textbf{trivial path}.\index{trivial path} \index{path!trivial}  A path of length $\geq 1$ is called an \textbf{internal path}.\index{internal path} \index{path!internal}
\item
Given a path $P$ as above, call $v_0$ (resp., $v_r$) its \textbf{initial vertex} \index{initial vertex} (resp., \textbf{terminal vertex}).\index{terminal vertex}  An \textbf{end vertex} \index{end vertex} means either an initial vertex or a terminal vertex.
\item
An internal path whose initial vertex is equal to its terminal vertex is called a  \textbf{cycle}.\index{cycle}  Otherwise, it is called a \textbf{trail}.\index{trail}
\item
A \textbf{directed path} \index{directed path} \index{path!directed} in $G$ is an internal path $P$ as above such that each $e^j$ has initial vertex $v_{j-1}$ and terminal vertex $v_j$.
\item
A \textbf{wheel} \index{wheel} in $G$ is a directed path that is also a cycle.
\end{enumerate}
\end{definition}

\begin{remark}
If $P$ is a wheel, then a cyclic permutation of its edges and vertices is also a wheel.  In what follows, we will not distinguish between a wheel and its cyclic permutations.
\end{remark}

\begin{remark}
When we specify a path, we will sometimes just specify the vertices $v_i$ or just the edges $e^j$.
\end{remark}

\begin{example}
\label{ex:oppositepath}
Given an internal path $P$ as above, by reversing the labels of the $e^j$ and the $v_i$, we obtain an internal path $P^{\op}$, called the \textbf{opposite internal path}.\index{internal path!opposite}  It contains the same sets of ordinary edges and vertices as $P$.  Its initial (resp., terminal) vertex is the terminal (resp., initial) vertex of $P$.
\end{example}

\begin{example}
Suppose $v_i$ for $0\leq i \leq 4$ are distinct vertices in a graph $G$.  Suppose $e^j$ is an ordinary edge adjacent to both $v_{j-1}$ and $v_j$ as in the following picture.
\begin{center}
\begin{tikzpicture}
\matrix[row sep=1cm,column sep=1cm] {
\node [plain] (v0) {$v_0$}; 
& \node [plain] (v1) {$v_1$};
& \node [plain] (v2) {$v_2$};
& \node [plain] (v3) {$v_3$};
& \node [plain] (v4) {$v_4$};\\
};
\draw [arrow] (v0) to node{$e^1$} (v1);
\draw [arrow] (v2) to node[swap]{$e^2$} (v1);
\draw [arrow] (v2) to node{$e^3$} (v3);
\draw [arrow] (v3) to node{$e^4$} (v4);
\end{tikzpicture}
\end{center}
Then this is an internal path of length $4$ that is also a trail.  Likewise, the picture
\begin{center}
\begin{tikzpicture}
\matrix[row sep=1cm,column sep=1cm] {
\node [plain] (v2) {$v_2$};
& \node [plain] (v3) {$v_3$};
& \node [plain] (v4) {$v_4$};\\
};
\draw [arrow] (v2) to node{$e^3$} (v3);
\draw [arrow] (v3) to node{$e^4$} (v4);
\end{tikzpicture}
\end{center}
depicts a directed path of length $2$, but it is not a wheel.
\end{example}

\begin{example}
\label{ex:cyclenotwheel}
Suppose $v_i$ for $0\leq i \leq 3$ are distinct vertices in a graph $G$.  Suppose $e^j$ is an ordinary edge adjacent to both $v_{j-1}$ and $v_j$ as in the following picture.
\begin{center}
\begin{tikzpicture}
\matrix[row sep=.6cm,column sep=1cm] {
& \node [plain] (v2) {$v_2$}; &\\
\node [plain] (v1) {$v_1$};
&& \node [plain] (v3) {$v_3$};\\
& \node [plain] (v0) {$v_0$}; &\\
};
\draw [arrow] (v0) to node{$e^1$} (v1);
\draw [arrow] (v2) to node[swap]{$e^2$} (v1);
\draw [arrow] (v2) to node{$e^3$} (v3);
\draw [arrow] (v3) to node{$e^4$} (v0);
\end{tikzpicture}
\end{center}
Then this is a cycle of length $4$, but it is not a wheel.  
\end{example}

\begin{example}
If the orientation of $e_2$ is reversed in Example \ref{ex:cyclenotwheel}, then we have
\begin{center}
\begin{tikzpicture}
\matrix[row sep=.6cm,column sep=1cm] {
& \node [plain] (v2) {$v_2$}; &\\
\node [plain] (v1) {$v_1$};
&& \node [plain] (v3) {$v_3$};\\
& \node [plain] (v0) {$v_0$}; &\\
};
\draw [arrow] (v0) to node{$e^1$} (v1);
\draw [arrow] (v1) to node{$e^2$} (v2);
\draw [arrow] (v2) to node{$e^3$} (v3);
\draw [arrow] (v3) to node{$e^4$} (v0);
\end{tikzpicture}
\end{center}
which is a wheel of length $4$.
\end{example}

\subsection{Connected Graphs}

First we define the concept of connected graphs.

\begin{definition}
Suppose $G$ is a graph.  Then $G$ is called a \textbf{connected graph} \index{connected graph} \index{graph!connected} if one of the following three statements is true.
\begin{enumerate}
\item
$G$ is a single exceptional edge (Example \ref{ex:exceptionaledge}).
\item
$G$ is a single exceptional loop (Example  \ref{ex:exceptionalloop}).
\item
$G$ satisfies all of the following conditions.
\begin{itemize}
\item
$G$ is ordinary (i.e., has no exceptional flags).
\item
$G$ is \emph{not} the empty graph (Example \ref{ex:emptygraph}).
\item
For any two distinct vertices $u$ and $v$ in $G$, there exists an internal path in $G$ with $u$ as its initial vertex and $v$ as its terminal vertex.
\end{itemize}
\end{enumerate}
\end{definition}

\begin{example}
The following are examples of connected graphs.
\begin{itemize}
\item
A single isolated vertex $V_1$ (Example \ref{ex:isolatedvt}).
\item
The $c$-colored exceptional edge $\uparrow_c$ (Example \ref{ex:exceptionaledge}).
\item
The $c$-colored exceptional loop $\wheelc$ (Example \ref{ex:exceptionalloop}).
\item
The $\dch$-corolla $C_{\dch}$ (Example \ref{ex:corolla}).
\item
The permuted corolla $\sigma C_{\dch} \tau$ (Example \ref{ex:permutedcor}).
\item
The contracted corolla $\xiij C_{\dch}$ (Example \ref{ex:contractedcor}).
\item
The partially grafted corollas $C_{\dch} \boxtimes^{\uc'}_{\ub'} C_{\bah}$ (Example \ref{ex:pgcor}).
\item
The dioperadic graph $C_{\dch} \icompj C_{\bah}$ (Example \ref{ex:dioperadic}).
\end{itemize}  
\end{example}

\begin{example}
On the other hand, the following graphs are \emph{not} connected.
\begin{itemize}
\item
The empty graph $\varnothing$ (Example \ref{ex:emptygraph}).
\item
The union of $n$ isolated vertices $V_n$ with $n \geq 2$ (Example \ref{ex:isolatedvt}).
\end{itemize}
\end{example}

\subsection{Wheel-Free Graphs}

Next we define some related classes of graphs.

\begin{definition}
\label{def:wheelfreegraphs}
Suppose $G$ is a graph.
\begin{enumerate}
\item
We say $G$ is \textbf{wheel-free} \index{wheel-free} \index{graph!wheel-free} if it contains neither exceptional loops nor wheels.
\item
We say $G$ is \textbf{simply connected} \index{simply connected} \index{graph!simply connected} if it is connected, is not an exceptional loop, and contains no cycles.
\item
We call $G$ a \textbf{unital tree} \index{unital tree} if it is simply connected in which each vertex has exactly one output flag.
\item
We call $G$ a \textbf{linear graph} \index{linear graph} \index{graph!linear} if it is a unital tree in which each vertex has exactly one input flag.
\item
We say $G$ has \textbf{non-empty inputs} \index{graph!non-empty inputs} (resp., \textbf{non-empty outputs}) \index{graph!non-empty outputs} if $\inp(v)$ (resp., $\out(v)$) is non-empty for each vertex $v$ in $G$.
\item
We say $G$ is \textbf{special} \index{graph!special} if it has non-empty inputs and non-empty outputs.
\end{enumerate}
\end{definition}

\begin{definition}
\label{def:setsofgraphs}
Define the following sets of graphs.
\begin{enumerate}
\item
$\gwheelc$\label{note:gwheelc} is the set of connected graphs.
\item
$\gupc$\label{note:gupc} is the set of connected wheel-free graphs.
\item
$\gupci$\label{note:gupci} (resp., $\gupco$) is the set of connected wheel-free graphs with non-empty inputs (resp., outputs).
\item
$\gupcs$\label{note:gupcs} is the set of special connected wheel-free graphs.
\item
$\gupd$\label{note:gupd} is the set of simply connected graphs.
\item
$\uoperad$\label{note:uoperad} is the set of unital trees.
\item
$\ULin$\label{note:ulin} is the set of linear graphs.
\end{enumerate}
\end{definition}

\begin{remark}
\label{rk:colorsforgraphs}
In view of Convention \ref{graphconvention}, the sets in Definition \ref{def:setsofgraphs} are really sets of strict isomorphism classes of $\fC$-colored wheeled graphs with the indicated properties for a fixed color set $\fC$.  In particular, if we wish to emphasize the color set $\fC$, we will write $\gupc(\fC)$ and so forth.
\end{remark}

\begin{remark}
In each of $\ULin$, $\uoperad$, $\gupd$, $\gupci$, $\gupco$, $\gupcs$, and $\gupc$, the only graphs without any vertex are the exceptional edges $\uparrow_c$ of a single color $c$.  On the other hand, the set $\gwheelc$ also contains the exceptional loops $\wheel_c$ of a single color.  
\end{remark}

\begin{remark}
We have the following strict inclusions:
\[\nicexy{
& \gupcs \ar[r] \ar[dr] & \gupci \ar[dr] &&\\
\ULin \ar[dr] \ar[ur] & &  \gupco \ar[r] & \gupc \ar[r] & \gwheelc\\
& \uoperad \ar[r] \ar[ur] & \gupd \ar[ur] &&
}\]
Moreover, $\gupcs$ is the intersection of $\gupci$ and $\gupco$ within $\gupc$.
\end{remark}

\begin{remark}
\label{rk:nonemptyinout}
If $G \in \gupci$ (resp., $G \in \gupco$), then $\inp(G)$ (resp., $\out(G)$) is non-empty.  Indeed, an exceptional edge has both an input leg and an output leg.  If $G \in \gupci$ (resp., $\gupco$) is ordinary, then the initial (resp., terminal) vertex of any maximal directed path in $G$ has an incoming (resp., outgoing) flag that must be an input (resp., output) leg of $G$.  On the other hand, even if $G \in \gupc$ has $\inp(G)$ (resp., $\out(G)$) non-empty, it does \emph{not} follow that each vertex in $G$ has non-empty inputs (resp., outputs).  For example, the connected graph
\begin{center}
\begin{tikzpicture}
\matrix[row sep=.6cm,column sep=1cm] {
& \node [plain] (u) {$u$}; \\
\node [plain] (v) {$v$}; &\\
};
\draw [inputleg] (u) to +(0,-1cm);
\draw [arrow] (v) to (u);
\end{tikzpicture}
\end{center}
has $\inp(G)$ non-empty, but $\inp(v) = \varnothing$.  So $G$ does not have non-empty inputs.  
\end{remark}

\begin{remark}
\label{rk:lineargraphln}
A linear graph is precisely a simply connected graph in which each vertex has exactly one input flag and one output flag.  There is a canonical bijection between the set of linear graphs and the set $\coprod_{n \geq 1} \fC^n$.  Indeed, the linear graphs with $0$ vertex are the $c$-colored exceptional edges $\uparrow_c$ with $c \in \fC$.  For $k \geq 1$ the linear graphs with $k$ vertices all have the form
\begin{center}
\begin{tikzpicture}
\matrix[row sep=.6cm,column sep=1cm] {
\node [plain] (vk) {$v_k$}; \\
\node [empty] (vdot) {$\vdots$}; \\
\node [plain] (v2) {$v_2$};\\
\node [plain] (v1) {$v_1$};\\
};
\draw [inputleg] (v1) to +(0,-.7cm);
\draw [arrow] (v1) to (v2);
\draw [arrow] (v2) to (vdot);
\draw [arrow] (vdot) to (vk);
\draw [outputleg] (vk) to +(0,.7cm);
\end{tikzpicture}
\end{center}
The $k-1$ internal edges, the single input leg, and the single output leg can be arbitrary colors.
\end{remark}

\begin{remark}
In a simply connected graph, given any two distinct vertices $u$ and $v$, there is a \emph{unique} internal path $P$ with initial vertex $u$ and terminal vertex $v$.  The opposite internal path $P^{\op}$ (Example \ref{ex:oppositepath}) has initial vertex $v$ and terminal vertex $u$.  The internal paths $P$ and $P^{\op}$ are the only ones that have $u$ and $v$ as end vertices.  Therefore, with a slight abuse of terminology, we say that in a simply connected graph, there is a unique path between any two distinct vertices.
\end{remark}

\begin{remark}
A simply connected graph is both connected and wheel-free.  On the other hand, a wheel-free graph may contain cycles, such as the one in Example \ref{ex:cyclenotwheel}, although it cannot contain exceptional loops or wheels.
\end{remark}

\begin{remark}
A graph $G$ is connected wheel-free if and only if one of the following two statements is true.
\begin{enumerate}
\item
$G$ is a single exceptional edge (Example \ref{ex:exceptionaledge}).
\item
$G$ satisfies all of the following conditions.
\begin{itemize}
\item
$G$ is ordinary (i.e., has no exceptional flags).
\item
$G$ is \emph{not} the empty graph (Example \ref{ex:emptygraph}).
\item
$G$ has no wheels.
\item
For any two distinct vertices $u$ and $v$ in $G$, there exists an internal path in $G$ with $u$ as its initial vertex and $v$ as its terminal vertex.
\end{itemize}
\end{enumerate}
\end{remark}

\begin{example}
The exceptional edge $\uparrow_c$ and the permuted corolla $\sigma C_{\dch} \tau$ with $|\ud| = 1$ are unital trees.
\end{example}

\begin{example}
Examples of simply connected graphs that are \emph{not} unital trees include:
\begin{itemize}
\item
a single isolated vertex $V_1$ and
\item
a permuted corolla $\sigma C_{\dch} \tau$ with $|\ud| \not= 1$.
\end{itemize}
\end{example}

\begin{example}
\label{ex:dioperadicgraph}
The partially grafted corollas $C_{\dch} \boxtimes^{\uc'}_{\ub'} C_{\bah}$ (Example \ref{ex:pgcor}) is connected wheel-free. 
\begin{enumerate}
\item
Moreover, such a partially grafted corollas is simply connected if and only if it is a dioperadic graph (Example \ref{ex:dioperadic}).
\item
A dioperadic graph is a unital tree if and only if $|\ud| = 1 = |\ub|$.  In this case, it looks like
\begin{center}
\begin{tikzpicture}
\matrix[row sep=2cm,column sep=1.5cm] {
\node [plain] (p1) {$v$}; \\
\node [plain,label=below:$\ua$] (p2) {$u$}; \\
};
\foreach \x in {1,2}
{
\draw [inputleg] (p\x) to +(-.6cm,-.5cm);
\draw [inputleg] (p\x) to +(.6cm,-.5cm);
}
\draw [outputleg] (p1) to node[above=.1cm]{$d_1$} +(0,.7cm);
\draw [arrow] (p2) to node{$b_1=c_j$} (p1);
\end{tikzpicture}
\end{center}
and is called a \textbf{simple tree}.  It generates the operadic $\circ_j$ operation \cite{ger,may97}.
\end{enumerate}
\end{example}

\section{Graph Substitution}
\label{sec:graphsub}

The free (wheeled) properad monad is induced by the operation of graph substitution.  The reader is referred to \cite{jy2} for the detailed construction of graph substitution and proof of its associativity and unity properties.  Intuitively, at each vertex $v$ in a given graph $G$, we drill a small hole containing $v$ and replace $v$ with a scaled down version of another graph $H_v$ whose profiles are the same as those of $v$.

\subsection{Properties of Graph Substitution}

\begin{definition}
Suppose $G$ is a graph with profiles $\dc$, and $H_v$ is a graph with profiles $\profilev$ for each $v \in \vertex(G)$.  Define the \textbf{graph substitution}\label{note:graphsub} \index{graph substitution}
\[
G\left(\{H_v\}_{v\in \vertex(G)}\right),
\]
or just $G(\{H_v\})$, as the graph obtained from $G$ by
\begin{enumerate}
\item
replacing each vertex $v \in \vertex(G)$ with the graph $H_v$, and
\item
identifying the legs of $H_v$ with the incoming/outgoing flags of $v$.
\end{enumerate}
We say that $H_v$ is \textbf{substituted into} $v$.
\end{definition}

\begin{remark}
Let us make a few observations.
\begin{enumerate}
\item
The graph substitution $G(\{H_v\})$ has the same input/output profiles as $G$. Moreover, there is a canonical identification
\[
\vertex\left(G(\{H_v\})\right) = \coprod_{\ving} \vertex(H_v).
\]
All the internal edges in the $H_v$ become internal edges in $G(\{H_v\})$.
\item
Corollas are units for graph substitution, in the sense that
\[
C(G) = G = G(\{C_v\}),
\]
where on the left $C$ denotes the corolla whose unique vertex has the same profiles as $G$.  On the right, $C_v$ is the corolla with the same profiles as $v$.
\item
Graph substitution is associative in the following sense.  Suppose $I_u$ is a graph with the same profiles as $u$ for each $u \in \vertex(G\{H_v\})$.  If $u$ is a vertex in $H_v$, we write $I_u$ as $I^v_u$ as well.  Then
\[
\left[G(\{H_v\})\right] \left(\{I_u\}\right) = G\left(\left\{ H_v(\{I^v_u\})\right\}\right).
\]
\item
Each of the sets of graphs in Definition \ref{def:setsofgraphs} is closed under graph substitution.  For example, if $G$ and all the $H_v$ are connected (wheel-free) graphs, then so is the graph substitution $G(\{H_v\})$.
\end{enumerate}
\end{remark}

\begin{notation}
\label{notation:gh}
If $H_w$ is a graph with the profiles of a vertex $w \in \vertex(G)$, then we use the abbreviation
\[\label{note:ghw}
G(H_w) = G(\{H_v\}),
\]
where for vertices $u \not= w$,  $H_u$ is the corolla with the profiles of $u$.  There are a canonical bijection
\[
\vertex\left(G(H_w)\right) 
= \vertex(H_w) \coprod \left[\vertex(G) \setminus \{w\}\right],
\]
a canonical injection
\[
\edgei(H_w) \hookrightarrow
\edgei\left(G(H_w)\right),
\]
and a canonical map
\[
\Leg(H_w) \to \edge\left(G(H_w)\right)
\]
that identifies each leg of $H_w$ with an element in $\inp(w) \cup \out(w)$.  Note that for a vertex $v \in \vertex(G)$, the sets $\inp(v)$ and $\out(v)$, regarded as subsets of $\edge(G)$, may have non-empty intersection.  In fact, $e \in \edge(G)$ lies in both $\inp(v)$ and $\out(v)$ precisely when $e$ is a loop at $v$.  So if $G$ does not have loops, then there is an injection
\[
\edge(H_w) \hookrightarrow \edge\left(G(H_w)\right).
\]
\end{notation}

\subsection{Examples}

\begin{example}
The exceptional loop can be obtained by substituting an exceptional edge into a contracted corolla, i.e.,
\[
\left(\xi^1_1 C_{(c;c)}\right) \left(\uparrow_c\right) = \wheel_c
\]
for any color $c$.
\end{example}

\begin{example}
Substituting an exceptional wheel into an isolated vertex yields the same exceptional wheel, i.e.,
\[
\left(\bullet\right)(\wheel_c) = \wheel_c
\]
for any color $c$.
\end{example}

\begin{example}
\label{ex:pgcorgeneration}
A partially grafted corollas can be obtained from a dioperadic graph by repeated contractions.  To be precise, we use the notations in Examples \ref{ex:contractedcor}, \ref{ex:pgcor}, and \ref{ex:dioperadic}.  Then we have a graph substitution decomposition
\[
C_{\dch} \boxtimes_{\be} C_{\bah} 
= (\xi_{e_\alpha}C_\alpha) \cdots (\xi_{e_2} C_2) \left(C_{\dch} \circ_{e_1} C_{\bah}\right).
\]
Here the $C_i$ are corollas with appropriate profiles such that the graph substitutions make sense.  There are many other such decompositions of the partially grafted corollas.  For example, we can start with a dioperadic graph whose internal edge is some $e_r$.  Then we use $\alpha-1$ contractions to create the other internal edges in any order.
\end{example}

\begin{example}
Suppose $G$ is the graph
\begin{center}
\begin{tikzpicture}
\matrix[row sep=1cm,column sep=1cm] {
\node [plain] (v) {$v$}; \\
\node [plain] (u) {$u$}; \\
\node [plain] (t) {$t$}; \\
};
\draw [outputleg] (v) to +(-.5cm,.4cm);
\draw [outputleg] (v) to +(.5cm,.4cm);
\draw [inputleg] (t) to +(0,-.6cm);
\draw [arrow, bend left=60] (t) to node{\footnotesize{$a$}} (v);
\draw [arrow] (t) to node[swap]{\footnotesize{$b$}} (u);
\draw [arrow] (u) to node[swap]{\footnotesize{$b$}} (v);
\end{tikzpicture}
\end{center}
with three vertices, three internal edges, one input, and two outputs, where $a$ and $b$ are colors.  The legs can be given arbitrary colors.  Suppose $H_u =~ \uparrow_b$, the $b$-colored exceptional edge.  Suppose $H_v$ is the contracted corolla
\begin{center}
\begin{tikzpicture}
\matrix[row sep=1cm,column sep=1cm] {
\node [plain] (z) {$z$}; \\
};
\draw [outputleg] (z) to +(-.5cm,.4cm);
\draw [outputleg] (z) to +(.5cm,.4cm);
\draw [inputleg] (z) to node[below left=.1cm]{\footnotesize{$a$}} +(-.5cm,-.4cm);
\draw [inputleg] (z) to node[below right=.1cm]{\footnotesize{$b$}} +(.5cm,-.4cm);
\draw [arrow, out=20, in=-20, loop] (z) to ();
\end{tikzpicture}
\end{center}
with two inputs with colors $a$ and $b$, respectively, two outputs, and a loop of arbitrary color.  Suppose $H_t$ is the graph
\begin{center}
\begin{tikzpicture}
\matrix[row sep=1cm,column sep=1cm] {
\node [plain] (y) {$y$}; \\
\node [plain] (w) {$w$}; \\
\node [plain] (x) {$x$}; \\
};
\draw [outputleg] (y) to node[above=.1cm]{\footnotesize{$a$}}+(0,.6cm);
\draw [inputleg] (y) to +(-.6cm,-.3cm);
\draw [outputleg] (x) to node[right=.1cm]{\footnotesize{$b$}} +(.5cm,.4cm);
\draw [arrow] (y) to (w);
\draw [arrow, bend left=45] (x) to (y);
\draw [arrow] (w) to (x);
\end{tikzpicture}
\end{center}
with three vertices, three internal edges, one input, and two outputs of colors $a$ and $b$, respectively.  The input and the internal edges can have arbitrary colors.  Then the graph substitution $G(\{H_t,H_u,H_v\})$ is the graph
\begin{center}
\begin{tikzpicture}
\matrix[row sep=1cm,column sep=1cm] {
\node [plain] (z) {$z$}; \\
\node [plain] (y) {$y$}; \\
\node [plain] (w) {$w$}; \\
\node [plain] (x) {$x$}; \\
};
\draw [outputleg] (z) to +(-.5cm,.4cm);
\draw [outputleg] (z) to +(.5cm,.4cm);
\draw [inputleg] (y) to +(-.6cm,-.3cm);
\draw [arrow] (y) to (w);
\draw [arrow, bend left=45] (x) to (y);
\draw [arrow] (w) to (x);
\draw [arrow, bend left=35] (y) to node{\footnotesize{$a$}} (z);
\draw [arrow, bend right=45] (x) to node[swap]{\footnotesize{$b$}} (z);
\draw [arrow, out=20, in=-20, loop] (z) to ();
\end{tikzpicture}
\end{center}
with four vertices, six internal edges (one of which is a loop at $z$), one input, and two outputs.
\end{example}

\section{Closest Neighbors}
\label{sec:closestnbd}

In this section, we discuss a connected wheel-free analog of two vertices in a simply connected graph connected by an ordinary edge.  This concept will be useful for various purposes later, for example, in the definition of inner coface maps in the graphical category for connected wheel-free graphs and in studying the tensor product of two free properads.

All the assertions in this section concerning connected wheel-free graphs have obvious analogs for connected wheel-free graphs with non-empty inputs or non-empty outputs.  Since the proofs are the same in all three cases, we will not state the non-empty input/output cases separately.

\subsection{Motivating Examples}

In a simply connected graph (e.g., a unital tree or a linear graph), if two vertices $u$ and $v$ are connected by an ordinary edge $e$, then there is only one such ordinary edge.  Moreover, the two vertices can be combined into a single vertex with $e$ deleted, and the resulting graph is still simply connected.  Such combination of two vertices into one is how inner coface maps are defined in the finite ordinal category $\varDelta$ and the Moerdijk-Weiss dendroidal category $\varOmega$.

The obvious analog is \emph{not} true for connected wheel-free graphs in general.

\begin{example}
For example, in a partially grafted corollas (Example \ref{ex:pgcor}) with at least two ordinary edges, if the two vertices are combined  with some but not all of the ordinary edges deleted, then the resulting graph has a directed loop at the combined vertex.  Here is an example with two internal edges, in which $e$ is deleted when the two vertices are combined.
\begin{center}
\begin{tikzpicture}
\matrix[row sep=.3cm, column sep=1cm]{
\node [plain] (v) {$v$};&&&\\
& \node [empty] (s) {}; & \node [empty] (t) {}; & \node [plain] (w) {$w$};\\
\node [plain] (u) {$u$}; &&&\\
};
\draw [mapto] (s) to (t);
\draw [arrow,bend left=40] (u) to node{$e$} (v);
\draw [arrow,bend right=40] (u) to node[swap]{$f$} (v);
\draw [arrow, looseness=15, in=-45, out=45, loop] (w) to node{$f$} ();
\end{tikzpicture}
\end{center}
In particular, it would not be wheel-free any more.
\end{example}

\begin{example}
Even if all the ordinary edges between such vertices are deleted when the vertices are combined, it still does not guarantee that the result is wheel-free.  For example, when $u$ and $v$ are combined with $e$ deleted in the graph on the left
\begin{center}
\begin{tikzpicture}
\matrix[row sep=.3cm, column sep=1cm]{
\node [plain] (v) {$v$};&&&&&\\
& \node [plain] (x1) {$x$}; & \node [empty] (s) {}; & \node [empty] (t) {}; & \node [plain] (y) {$y$}; & \node [plain] (x2) {$x$};\\
\node [plain] (u) {$u$}; &&&&&\\
};
\draw [mapto] (s) to (t);
\draw [arrow] (u) to node{$e$} (v);
\draw [arrow] (u) to node[swap]{$f$} (x1);
\draw [arrow] (x1) to node[swap]{$g$} (v);
\draw [arrow,bend right=40] (y) to node[swap]{$f$} (x2);
\draw [arrow,bend right=40] (x2) to node[swap]{$g$} (y);
\end{tikzpicture}
\end{center}
the resulting graph on the right has a wheel.
\end{example}

\subsection{Closest Neighbors}

In order to smash two vertices together and keep the resulting graph connected wheel-free, we need the following analog of two vertices connected by an ordinary edge.

\begin{definition}
\label{def:closestnbd}
Suppose $G$ is a connected wheel-free graph, and $u$ and $v$ are two distinct vertices in $G$.  Call $u$ and $v$ \textbf{closest neighbors} \index{closest neighbors} \index{vertex!closest neighbors} if:
\begin{enumerate}
\item
there is at least one ordinary edge adjacent to both of them, and
\item
there are no directed paths with initial vertex $u$ and terminal vertex $v$  that involve a third vertex.
\end{enumerate}
In this case, we also say $v$ is a \textbf{closest neighbor of $u$}.
\end{definition}

\begin{remark}
In \cite{vallette} (section 4) the author defined a very similar and probably equivalent concept called \emph{adjacent vertices}.
\end{remark}

\begin{example}
\label{ex:cneighbor}
Consider the following connected wheel-free graph $K$.
\begin{center}
\begin{tikzpicture}
\matrix[row sep=.3cm, column sep = 1.5cm]{
\node [plain] (w) {$w$};&\\
& \node [plain] (v) {$v$};\\
\node [plain] (u) {$u$}; & \\
};
\draw [arrow] (u) to (v);
\draw [arrow] (u) to (w);
\draw [arrow,bend left=30] (v) to (w);
\draw [arrow,bend right=30] (v) to (w);
\end{tikzpicture}
\end{center}
Then $u$ and $v$ are closest neighbors, as are $v$ and $w$.  The vertices $u$ and $w$ are \emph{not} closest neighbors.
\end{example}

\begin{example}
The two vertices in a partially grafted corollas (Example \ref{ex:pgcor}) are closest neighbors.
\end{example}

\begin{example}
\label{ex:cnbdexist}
For $G \in \gupc$ with at least two vertices, every vertex $u$ has a closest neighbor.  Indeed, either $\inp(u)$ or $\out(u)$ is not empty, so let us assume $\out(u)$ is not empty.  Among all the \emph{directed} paths with initial vertex $u$, pick one, say $P$, with maximal length.  Then the vertex $v$ in $P$ right after $u$ must be a closer neighbor of $u$, since otherwise $P$ would not be maximal.
\end{example}

\subsection{Inner Properadic Factorization}

Closest neighbors are intimately related to graph substitution in $\gupc$ involving partially grafted corollas.  To make this precise, we need the following definition, where we use the notation in \ref{notation:gh}.

\begin{definition}
\label{def:ipropfact}
Suppose $K$ is a connected wheel-free graph.  An \textbf{inner properadic factorization of $K$} \index{inner properadic factorization} \index{properadic factorization!inner} is a graph substitution decomposition
\[
K = G(\{H_v\}) = G(H_w)
\]
in which
\begin{itemize}
\item
$G$ is connected wheel-free,
\item
a chosen $H_w$ is a partially grafted corollas, and 
\item
all other $H_u$ are corollas.
\end{itemize}
In this case, $H_w$ is called the \textbf{distinguished subgraph}.\index{distinguished subgraph}
\end{definition}

\begin{remark}
The ``inner" in Definition \ref{def:ipropfact} refers to the assumption that the distinguished partially grafted corollas $H_w$ is an inner graph in the graph substitution $G(\{H_v\})$.
\end{remark}

\begin{example}
\label{ex0:ipropfact}
Suppose $K$ is a partially grafted corollas with profiles $\dch$.  Then there is an inner properadic factorization
\[
K = C_{\dch}(K),
\] 
where $K$ itself is the distinguished subgraph.
\end{example}

\begin{example}
\label{ex:properadicfact}
The graph $K$ in Example \ref{ex:cneighbor} admits an inner properadic factorization in which $G$ is
\begin{center}
\begin{tikzpicture}
\matrix[row sep=1cm, column sep = 1cm]{
\node [plain] (y) {$y$};\\
\node [plain] (u) {$u$}; \\
};
\draw [arrow,bend left=40] (u) to (y);
\draw [arrow,bend right=40] (u) to (y);
\end{tikzpicture}
\end{center}
and $H_u$ is the corolla $C_u$ with the same profiles as $u$.  The distinguished subgraph $H_y$ is
\begin{center}
\begin{tikzpicture}
\matrix[row sep=1cm, column sep = 1cm]{
\node [plain] (w) {$w$};\\
\node [plain] (v) {$v$}; \\
};
\draw [arrow,bend left=40] (v) to (w);
\draw [arrow,bend right=40] (v) to (w);
\draw [inputleg] (w) to +(-.7cm,-.5cm);
\draw [inputleg] (v) to +(0,-.7cm);
\end{tikzpicture}
\end{center}
which has two inputs and empty output.  This inner properadic factorization corresponds to the closest neighbors $v$ and $w$, in the sense that $H_y$ is defined by the flags in $v$ and $w$ as well as the ordinary edges adjacent to both of them.
\end{example}

\begin{example}
\label{ex2:properadicfact}
The only other  inner properadic factorization of $K$ in Example \ref{ex:cneighbor} is the one in which $G'$ is
\begin{center}
\begin{tikzpicture}
\matrix[row sep=1cm, column sep = 1cm]{
\node [plain] (w) {$w$};\\
\node [plain] (z) {$z$}; \\
};
\draw [arrow,bend left=40] (z) to (w);
\draw [arrow] (z) to (w);
\draw [arrow,bend right=40] (z) to (w);
\end{tikzpicture}
\end{center}
and $H_w$ is the corolla $C_w$ with the same profiles as $w$.  The distinguished subgraph $H_z$ is
\begin{center}
\begin{tikzpicture}
\matrix[row sep=1cm, column sep = 1cm]{
\node [plain] (v) {$v$};\\
\node [plain] (u) {$u$}; \\
};
\draw [arrow] (u) to (v);
\draw [outputleg] (u) to +(-.7cm,.5cm);
\draw [outputleg] (v) to +(-.7cm,.5cm);
\draw [outputleg] (v) to +(.7cm,.5cm);
\end{tikzpicture}
\end{center}
which has three outputs and empty input.  This inner properadic factorization corresponds to the closest neighbors $u$ and $v$.
\end{example}

Example \ref{ex:properadicfact} suggests a close relationship between inner properadic factorization and closest neighbors.  In fact, the two notions are equivalent.

\begin{theorem}
\label{thm:cnbdfact}
Suppose $K$ is a connected wheel-free graph, and $x$ and $y$ are distinct vertices in $K$.  Then the following statements are equivalent.
\begin{enumerate}
\item
The vertices $x$ and $y$ are closest neighbors in $K$. 
\item
$K$ admits an inner properadic factorization $G(H_w)$ in which the two vertices in the distinguished subgraph $H_w$ are $x$ and $y$.\index{closest neighbors!characterization} \index{inner properadic factorization!characterization}
\end{enumerate}
\end{theorem}

\begin{proof}
First suppose $x$ and $y$ are closest neighbors.  Define $G$ as the graph obtained from $K$ by:
\begin{itemize}
\item
combining the closest neighbors $x$ and $y$ into one vertex $w$, and
\item
deleting all the ordinary edges adjacent to both $x$ and $y$ in $K$.
\end{itemize}
The distinguished subgraph $H_w$ is defined using all the flags in the vertices $x$ and $y$ in $K$, with ordinary edges the ones adjacent to both $x$ and $y$ in $K$.  Then we have $K = G(\{H_v\})$, $H_w$ is a partially grafted corollas, and $G$ is still connected.  

It remains to see that $G$ is wheel-free.  Suppose to the contrary that $G$ has a wheel $Q$.   Then $Q$ must contain $w$, since all other $H_u$ are just corollas and $K$ is wheel-free.  Moreover, $Q$ must contain some other vertex $z$ because all the ordinary edges between $x$ and $y$ in $K$ have already been deleted during the passage to $G$.  As a vertex in $K$, $z$ is different from $x$ and $y$. Suppose $e$ is an ordinary edge in $K$ adjacent to both $x$ and $y$, and $Q'$ is the internal path in $K$ consisting of $Q$ and $e$.  Then the directed path $Q_0$ within $Q'$ consisting of the ordinary edges from $Q$ has $x$ and $y$ as end vertices.  Graphically, we have either one of the following two scenarios in $K$:
\begin{center}
\begin{tikzpicture}
\matrix[row sep=1.5cm, column sep=4cm]{
\node [plain] (y1) {$y$}; &
\node [plain] (y2) {$y$};\\
\node [plain] (x1) {$x$}; &
\node [plain] (x2) {$x$};\\
};
\draw [arrow] (x1) to node{$e$} (y1);
\draw [arrow] (x2) to node{$e$} (y2);
\draw [dashedarrow, bend right=80] (x1) to node[swap]{$Q_0$} (y1);
\draw [dashedarrow, bend left=80] (y2) to node{$Q_0$} (x2);
\end{tikzpicture}
\end{center}
The scenario on the right cannot happen because $K$ is wheel-free.  The scenario on the left cannot happen either because $Q_0$ contains $z$, and $x$ and $y$ are closest neighbors. So $G$ is connected wheel-free.

For the converse, note that any ordinary edge in $H_w$ is still one in $K$, and hence adjacent to both $x$ and $y$.  If $x$ and $y$ are not closest neighbors in $K$, then there is a directed path in $K$ that has $x$ and $y$ as end vertices and that contains a third vertex $z$.  But this implies that $G$ has a wheel containing $w$ and $z$, contradicting the wheel-free assumption on $G$.
\end{proof}

\section{Almost Isolated Vertices}
\label{sec:almostisolated}

In a linear graph with at least two vertices, the top or bottom vertex can be deleted, and the resulting graph is still a linear graph.  Such deletions correspond to top and bottom coface maps in the finite ordinal category $\varDelta$. Likewise, in a simply connected graph (resp., unital tree), if a vertex that is adjacent to only one ordinary edge is deleted, then the resulting graph is still simply connected (resp., a unital tree).  This is how outer coface maps are defined in the Moerdijk-Weiss dendroidal category $\varOmega$.

Here we develop an analogous concept of a vertex in a connected wheel-free graph that is extreme in some sense and whose deletion yields a connected wheel-free graph.  We will need this later to define outer coface maps in the graphical category for connected wheel-free graphs and also to study the tensor product of two free properads.

As in section \ref{sec:closestnbd}, all the assertions in this section concerning connected wheel-free graphs also hold for connected wheel-free graphs with non-empty inputs or non-empty outputs.  The arguments for the three cases are again the same.

\subsection{Definition and Examples}

\begin{definition}
\label{def:almostisolated}
Let $v$ be a vertex in a connected wheel-free graph $G$.  
\begin{enumerate}
\item
Call $v$ \textbf{weakly initial} \index{weakly initial vertex} \index{vertex!weakly initial} (resp., \textbf{weakly terminal}) \index{weakly terminal vertex} \index{vertex!weakly terminal} if all the ordinary edges adjacent to $v$ have $v$ as the initial (resp., terminal) vertex.  Call $v$ \textbf{extremal} \index{extremal} \index{vertex!extremal} if it is either weakly initial or weakly terminal.
\item
Call $v$ an \textbf{almost isolated} \index{almost isolated vertex} \index{vertex!almost isolated} vertex if either:
\begin{enumerate}
\item
$|\vertex(G)| = 1$ (i.e., $G$ is a permuted corolla), or
\item
\begin{itemize}
\item
$|\vertex(G)| \geq 2$,
\item
$v$ is extremal, and
\item
deleting $v$ from $G$ yields a connected wheel-free graph $G_v$.
\end{itemize}
\end{enumerate}
\end{enumerate}
\end{definition}

\begin{remark}
To be precise, in the previous definition, the graph $G_v$ is obtained from $G$ by deleting the non-exceptional cell $v$.  If $\{e_{-1},e_1\}$ is an ordinary edge in $G$ (i.e., $2$-cycle of $\iota$ within the non-exceptional cells) with one $e_i \in v$ and $e_{-i} \not\in v$, then we redefine $\iota(e_{-i}) = e_{-i}$ in $G_v$.  In other words, the flag $e_{-i}$ is a leg in $G_v$.
\end{remark}

\begin{remark}
If $|\vertex(G)| \geq 2$ and $v \in \vertex(G)$ is almost isolated, then we can visualize $G$ as follows.
\begin{center}
\begin{tikzpicture}
\matrix[row sep=1.5cm, column sep=2.5cm]{
\node [plain,label=above:$...$] (v1) {$v$}; &
\node [fatplain,label=above:$...$] (gv2) {$G_v$}; \\
\node [fatplain,label=below:$...$] (gv1) {$G_v$};& 
\node [plain,label=below:$...$] (v2) {$v$};\\
};
\draw [arrow,bend left=30] (gv1) to node[swap]{$.....$} (v1);
\draw [arrow,bend right=30] (gv1) to (v1);
\draw [arrow,bend left=30] (v2) to node[swap]{$.....$} (gv2);
\draw [arrow,bend right=30] (v2) to (gv2);
\foreach \x in {1,2}
{
\draw [inputleg] (gv\x) to +(-.6cm,-.5cm);
\draw [inputleg] (gv\x) to +(.6cm,-.5cm);
\draw [outputleg] (gv\x) to +(-.6cm,.5cm);
\draw [outputleg] (gv\x) to +(.6cm,.5cm);
}
\foreach \x in {1,2}
{
\draw [inputleg] (v\x) to +(-.6cm,-.5cm);
\draw [inputleg] (v\x) to +(.6cm,-.5cm);
\draw [outputleg] (v\x) to +(-.6cm,.5cm);
\draw [outputleg] (v\x) to +(.6cm,.5cm);
}
\end{tikzpicture}
\end{center}
On the left (resp., right), $v$ is weakly terminal (resp., weakly initial).  There are canonical bijections
\[
\begin{split}
\vertex(G) 
&=\vertex(G_v) \coprod \{v\},\\
\edge(G) 
&= \edge(G_v) \coprod \left(\edge(C_v) \setminus \{\be\}\right),\\
\edgei(G) 
&= \edgei(G_v) \coprod \{\be\}.
\end{split}
\]
Here $C_v$ is the corolla with the same profiles as the vertex $v$, and $\be$ is the non-empty set of internal edges between $G_v$ and $v$.  
\end{remark}

\begin{example}
In a partially grafted corollas (Example \ref{ex:pgcor}), both vertices are almost isolated, with the top (resp., bottom) vertex weakly terminal (resp., weakly initial).
\end{example}

\begin{example}
For the graph $K$ in Example \ref{ex:cneighbor}, the vertices $u$ and $w$ are almost isolated, but $v$ is not because it is not extremal.
\end{example}

\begin{example}
\label{ex:vgraph}
The vertex $u$ in the graph
\begin{center}
\begin{tikzpicture}
\matrix[row sep=.5cm, column sep=1cm]{
\node [plain] (v) {$v$}; && \node [plain] (w) {$w$};\\
& \node [plain] (u) {$u$}; &\\
};
\draw [arrow] (u) to (v);
\draw [arrow] (u) to (w);
\end{tikzpicture}
\end{center}
is weakly initial, and hence extremal.  However, it is \emph{not} almost isolated because deleting it would yield a graph that is not connected.
\end{example}

The following observation gives another characterization of a vertex, extremal or not, that can be deleted from a connected wheel-free graph.

\begin{lemma}
\label{deletev}
Let $v$ be a vertex in a connected wheel-free graph $G$ with at least two vertices.  Then the following statements are equivalent.
\begin{enumerate}
\item
Deleting $v$ from $G$ yields a connected wheel-free graph.
\item
Given any pair of vertices $x$ and $y$ in $G$ different from $v$, there is an internal path with end vertices $x$ and $y$ that does not contain $v$.
\end{enumerate}
\end{lemma}

\begin{proof}
Since $G$ is connected wheel-free and contains a vertex, it is also ordinary.  An ordinary wheeled graph is connected if and only if:
\begin{enumerate}
\item
it is not empty, and
\item
any two distinct vertices are connected by an internal path.
\end{enumerate}
The lemma follows from this characterization of connectivity in an ordinary wheeled graph.
\end{proof}

The following observation provides a simple suffiicient, but not necessary, condition that guarantees that a vertex is almost isolated.

\begin{lemma}
\label{oneordedge}
Suppose $G \in \gupc$, and $v \in \vertex(G)$ is adjacent to only one ordinary edge.  Then $v$ is almost isolated.
\end{lemma}

\begin{proof}
By definition $v$ is extremal.  It remains to show that deleting $v$ from $G$ yields a connected wheel-free graph.  So pick two vertices $x$ and $y$ different from $v$.  By connectivity there is an internal path $P$ with end vertices $x$ and $y$.  The internal path $P$ does \emph{not} contain $v$, since otherwise $P$ contains two different ordinary edges adjacent to $v$.  Therefore, by Lemma \ref{deletev}, deleting $v$ from $G$ yields a connected wheel-free graph.
\end{proof}

\begin{remark}
The converse of Lemma \ref{oneordedge} is \emph{not} true.  For example, in a partially grafted corollas with at least two ordinary edges, both vertices are almost isolated and adjacent to multiple ordinary edges.  However, the converse holds for simply connected graphs, as we now observe.
\end{remark}

\begin{proposition}
\label{prop:aisconn}
Suppose $v$ is a vertex in a simply connected graph $G$.  Then the following statements are equivalent.
\begin{enumerate}
\item
$v$ is almost isolated.
\item
There is only one ordinary edge adjacent to $v$.
\end{enumerate}
\end{proposition}

\begin{proof}
Lemma \ref{oneordedge} says that the second statement implies the first.  For the converse, suppose $v$ is almost isolated.  Then $G$ has at least two vertices.  Suppose to the contrary that there are at least two ordinary edges adjacent to $v$, which we depict as follows.
\begin{center}
\begin{tikzpicture}
\matrix[row sep=1cm, column sep=2cm]{
\node [plain] (u) {$u$};
& \node [plain] (v) {$v$};
& \node [plain] (w) {$w$};\\
};
\draw [thick] (u) to node{$e$} (v);
\draw [thick] (v) to node{$f$} (w);
\end{tikzpicture}
\end{center}
Since $G$ is simply connected, the concatenation of $e$ and $f$ is the \emph{unique} internal path with end vertices $u$ and $w$.  Therefore, once $v$ is deleted, the resulting graph \emph{cannot} be connected.  This contradicts the assumption on $v$, so there is only one ordinary edge adjacent to $v$.
\end{proof}

\subsection{Extremal Paths}

Later we will need to know that a connected wheel-free graph with at least two vertices always has at least two almost isolated vertices.  For a partially grafted corollas (Example \ref{ex:pgcor}), both vertices are almost isolated.  When there are more than two vertices in a connected wheel-free graph, we need the following concept to show the existence of almost isolated vertices.  To simplify the notation, when we write down an internal path below, we sometimes omit the ordinary edges and only exhibit the vertices involved.

\begin{definition}
Suppose $G$ is a connected wheel-free graph.
\begin{enumerate}
\item
An \textbf{extremal path} \index{extremal path} \index{path!extremal} in $G$ is an internal path
\[
P = \left(v_0, v_1, \ldots , v_r\right)
\]
such that:
\begin{itemize}
\item
$v_0 \not= v_r$, and
\item
both $v_0$ and $v_r$ are extremal vertices.
\end{itemize}
\item
An extremal path is \textbf{maximal} if there are no extremal paths that properly contain it.
\end{enumerate}
\end{definition}

\begin{remark}
\begin{enumerate}
\item
There can be many extremal paths involving exactly the same vertices because there may be multiple ordinary edges between $v_i$ and $v_{i+1}$.
\item
An extremal path is a trail, i.e., the vertices in it are \emph{all} different from each other.
\item
An end vertex of an extremal path need not be almost isolated.  For example, in the $V$-shape graph in Example \ref{ex:vgraph}, $(u,v)$ is an extremal path, but $u$ is not almost isolated.  On the other hand, we will show that the end vertices of a \emph{maximal} extremal path are both almost isolated.
\end{enumerate}
\end{remark}

First we want to establish the existence of a maximal extremal path, for which we need the following preliminary observation.

\begin{lemma}
\label{weakvertexexist}
Suppose $G$ is a connected wheel-free graph with at least two vertices.  Then there exist:
\begin{enumerate}
\item
a weakly initial vertex $u$ and a weakly terminal vertex $v \not= u$, and
\item
an extremal directed path $P$ with end vertices $u$ and $v$.\index{weakly initial vertex!existence} \index{weakly terminal vertex!existence}
\end{enumerate}
\end{lemma}

\begin{proof}
Let $P$ be a \emph{maximal} directed path in $G$.  Call its initial vertex $u$ and terminal vertex $v \not= u$.  Such a maximal directed path exists because there must be an ordinary edge in $G$.  Simply take $P$ as the longest directed path containing such an ordinary edge. Then $u$ is weakly initial, and $v$ is weakly terminal.  Indeed, if $u$ is not weakly initial, then there exists a vertex $w$ and an ordinary edge $e$ from $w$ to $u$.  Since $G$ is wheel-free, $w$ is different from all the vertices in $P$.  The current situation is depicted in the following picture.
\begin{center}
\begin{tikzpicture}
\matrix[row sep=1cm, column sep=1.5cm]{
\node [plain] (w) {$w$}; &
\node [plain] (u) {$u$}; &
\node [plain] (v) {$v$};\\
};
\draw [arrow] (w) to node{$e$} (u);
\draw [dashedarrow] (u) to node{$P$} (v);
\end{tikzpicture}
\end{center}
The concatenation of $e$ and $P$ is a directed path in $G$ that properly contains $P$, contradicting the maximality of $P$.  Therefore, $u$ must be weakly initial.  A similar argument shows that $v$ must be weakly terminal.  Since $u$ and $v$ are both extremal, $P$ is an extremal directed path.
\end{proof}

\begin{example}
Even an extremal directed path $P$ as in Lemma \ref{weakvertexexist} is \emph{not} necessarily maximal.  For example, in the connected wheel-free graph
\begin{center}
\begin{tikzpicture}
\matrix[row sep=.3cm, column sep=1cm]{
& \node [plain] (v) {$v$}; &\\
& &\\
\node [plain] (u) {$u$}; & \node [plain] (w) {$w$}; & \node [plain] (x) {$x$};\\
};
\draw [arrow] (u) to (v);
\draw [arrow] (u) to (w);
\draw [arrow] (w) to (v);
\draw [arrow] (w) to (x);
\end{tikzpicture}
\end{center}
$u$ is weakly initial, and $v$ is weakly terminal.  So
\[
P = (u,v)
\]
is an extremal directed path.  However, it is \emph{not} maximal because it is properly contained in the maximal extremal path
\[
Q = (u,v,w,x),
\]
where $x$ is weakly terminal.
\end{example}

\begin{proposition}
\label{lem1:gupcenoughface}
Suppose $G$ is a connected wheel-free graph with at least two vertices.  Then $G$ has a maximal extremal path.
\end{proposition}

\begin{proof}
Take an extremal directed path $P$ as in Lemma \ref{weakvertexexist}.  If it is a maximal extremal path, then we are done.  If $P$ is not maximal, then take the longest extremal path $Q$ containing $P$.  By construction the extremal path $Q$ must be maximal.
\end{proof}

\subsection{Existence of Almost Isolated Vertices}

Next we observe that, in a maximal extremal path, both end vertices are almost isolated.

\begin{theorem}
\label{lem2:gupcenoughface}
Suppose $G$ is a connected wheel-free graph with at least two vertices, and
\[
P = \left(v_0,\ldots,v_r\right)
\]
is a maximal extremal path in $G$.  Then both end vertices $v_0$ and $v_r$ are almost isolated.\index{almost isolated vertex!existence}
\end{theorem}

\begin{proof}
By symmetry it suffices to show that $v_0$ is almost isolated.  We assume that the extremal vertex $v_0$ is weakly initial.  There is a similar argument if $v_0$ is weakly terminal.  We must show that, if $v_0$ is deleted from $G$, then the resulting graph is still connected.  We will use the characterization in Lemma \ref{deletev}.  First, any two vertices in $P$ different from $v_0$ can be connected by an internal sub-path of $P$ that does not contain $v_0$.

Next suppose $w \in \vertex(G) \setminus P$.  It suffices to show that there exists an internal path from $w$ to $v_r$ that does not contain $v_0$.  Indeed, if this is true, then $w$ can be connected to any $v_i$ with $i>0$ via an internal path that does not contain $v_0$.  Similarly, if $z$ is another vertex not in $P$, then there is an internal path connecting  each of $w$ and $z$ to $v_r$, which does not contain $v_0$.  Splicing these internal paths at $v_r$ and taking only a subset of the ordinary edges if necessarily, we obtain an internal path from $w$ to $z$ not containing $v_0$.

To show the existence of the desired internal path from $w$ to $v_r$, we argue by contradiction.  So suppose every internal path from $w$ to $v_r$ must contain $v_0$.  Then every internal path from $w$ to any $v_i$ with $i>0$ must also contain $v_0$.  In fact, an internal path from $w$ to $v_i$ not containing $v_0$ together with the internal sub-path $(v_i,\ldots,v_r)$ of $P$ would give an internal path from $w$ to $v_r$ not containing $v_0$.  

Among all the internal paths from $v_0$ to $w$, pick the longest one,  and call it $Q$.  By maximality of $Q$, the first vertex $x$ in $Q$ after $v_0$ must be a closest neighbor of $v_0$, and $x \not= v_i$ for $i>0$.  Extend any ordinary edge $v_0 \to x$ to a maximal \emph{directed} path $R$.  The maximal directed path $R$ begins at the weakly initial vertex $v_0$ and ends at some weakly terminal vertex $y$, which may be equal to $x$.  Here is a diagram of the constructions so far, where a dashed line (resp., dashed arrow) represents an internal path (resp., directed path):
\begin{center}
\begin{tikzpicture}
\matrix[row sep=1cm, column sep=2cm]{
\node [plain] (w) {$w$}; &
\node [plain] (y) {$y$}; &
\node [plain] (vr) {$v_r$};\\
& \node [plain] (x) {$x$}; &\\
& \node [plain] (v0) {$v_0$}; &\\
};
\draw [arrow] (v0) to (x);
\draw [dashedline,bend right=40] (v0) to node[swap]{$P$} (vr);
\draw [dashedline] (x) to node{in $Q$} (w);
\draw [dashedarrow] (x) to node[swap]{in $R$} (y);
\end{tikzpicture}
\end{center}
The internal paths $P$ and $R$ are disjoint except at the vertex $v_0$, since otherwise there would be an internal path from $w$ to some $v_i$ with $i>0$ not containing $v_0$.  Splicing $P$ and $R$ together at $v_0$, we obtain an extremal path that properly contains $P$, contradicting the maximality of $P$.  Therefore, there must exist an internal path from $w$ to $v_r$ not containing $v_0$.
\end{proof}

\begin{corollary}
\label{aiexist}
Suppose $G$ is a connected wheel-free graph with at least two vertices.  Then $G$ has at least two almost isolated vertices.
\end{corollary}

\begin{proof}
By Proposition \ref{lem1:gupcenoughface} $G$ has a maximal extremal path $P$.  By Theorem \ref{lem2:gupcenoughface} the end vertices of $P$ are both almost isolated.
\end{proof}

\subsection{Outer Properadic Factorization}

Just like closest neighbors, the existence of almost isolated vertices is closely related to graph substitution involving partially grafted corollas (Example \ref{ex:pgcor}).  To make this relationship precise, we need the following definition, where we use the notation in \ref{notation:gh}.

\begin{definition}
\label{def:opropfact}
Suppose $K$ is an ordinary connected wheel-free graph.  An \textbf{outer properadic factorization of $K$} \index{outer properadic factorization} \index{properadic factorization!outer} is a graph substitution decomposition
\[
K = G(\{C_u,H_w\}) = G(H_w)
\]
in which
\begin{itemize}
\item
$G$ is a partially grafted corollas,
\item
a chosen $H_w$ is a connected wheel-free graph, and
\item
$C_u$ is a corolla.
\end{itemize}
In this case, $H_w$ is called the \textbf{distinguished subgraph}.\index{distinguished subgraph}
\end{definition}

\begin{remark}
\begin{enumerate}
\item
The ``outer" in Definition \ref{def:opropfact} refers to the assumption that the outer graph $G$ in the graph substitution is a partially grafted corollas.
\item
Outer properadic factorization is in some sense dual to inner properadic factorization.  In the latter, a partially grafted corollas along with a finite set of corollas are substituted into a connected wheel-free graph.  In an outer properadic factorization, a connected wheel-free graph and a corolla are substituted into the two vertices of a partially grafted corollas.
\end{enumerate}
\end{remark}

\begin{remark}
\label{rk:hwexceptional}
Recall the notion of a dioperadic graph from Example \ref{ex:dioperadicgraph}.
\begin{enumerate}
\item
A connected wheel-free graph is ordinary if and only if it has at least one vertex.  In particular, if $K$ has at least two vertices, then in any outer properadic factorization of $K$, the distinguished subgraph $H_w$ is ordinary, hence \emph{not} an exceptional edge.
\item
In Definition \ref{def:opropfact}, the distinguished subgraph $H_w$ is an exceptional edge $\uparrow$ if and only if
\begin{itemize}
\item
$K = H_u$ is a corolla, and
\item
$G$ is a dioperadic graph in which one vertex, corresponding to $w$, has one input and one output.
\end{itemize}
When $K$ is a corolla, there is one outer properadic factorization $G(\uparrow)$ for each leg of $K$, to which the vertex $w$ in $G$ is attached.
\end{enumerate}
\end{remark}

\begin{example}
This example refers to the graph $K$ in Example \ref{ex:cneighbor}
\begin{enumerate}
\item
The inner properadic factorization
\[
K = G(\{C_u,H_y\})
\]
in Example \ref{ex:properadicfact} is also an outer properadic factorization because $G$ is a partially grafted corollas, and $C_u$ is a corolla.
\item
Likewise, the inner properadic factorization
\[
K = G'(\{C_w, H_z\})
\]
in Example \ref{ex2:properadicfact} is also an outer properadic factorization because $G'$ is a partially grafted corollas, and $C_w$ is a corolla.
\end{enumerate}
\end{example}

\begin{example}
If
\[
K = G(\{C_u,H_w\})
\]
is an outer properadic factorization with $|\vertex(H_w)| \not= 2$, then it is \emph{not} an inner properadic factorization.  Likewise, if $K = G(H_w)$ is an inner properadic factorization with $|\vertex(G)| \not=2$, then it is \emph{not} an outer properadic factorization.
\end{example}

We now observe that outer properadic factorizations are equivalent to almost isolated vertices.

\begin{theorem}
\label{thm:opropfact}
Suppose $K$ is a connected wheel-free graph, and $x \in \vertex(K)$.  Then the following two statements are equivalent.
\begin{enumerate}
\item
$x$ is an almost isolated vertex.
\item
There is an outer properadic factorization $K = G(\{C_x,H_w\})$.\index{almost isolated vertex!characterization} \index{outer properadic factorization!characterization}
\end{enumerate}
\end{theorem}

\begin{proof}
We may assume $|\vertex(K)| \geq 2$.  Indeed, if $|\vertex(K)| = 1$, then $x$ is by definition almost isolated.  Moreover, in this case, $K$ is a permuted corolla, so there is an outer properadic factorization
\[
K = G(\{K, \uparrow\})
\]
for some dioperadic graph $G$ (Remark \ref{rk:hwexceptional}).

With $|\vertex(K)| \geq 2$, first suppose the vertex $x$ is almost isolated.  Define
\begin{itemize}
\item
$K_x$ as the connected wheel-free graph obtained from $K$ by deleting $x$, and 
\item
$C_x$ as the corolla with the same profiles as $x$.
\end{itemize}
Then there is an outer properadic factorization
\[
K = G(\{C_x,K_x\}),
\]
in which $G$ is the partially grafted corollas whose two vertices have the profiles of $K_x$ and $C_x$, and whose ordinary edges are those in $K$ between $K_x$ and $x$.  This graph substitution decomposition of $K$ is the required outer properadic factorization.

Conversely, suppose $K$ has an outer properadic factorization as stated.  Since $|\vertex(K)| \geq 2$, the distinguished subgraph $H_w$ has at least one vertex.  Since $G$ is a partially grafted corollas with $C_x$ a corolla, the connected wheel-free graph $H_w$ is obtained from $K$ by deleting $x$.  It remains to see that $x$ is extremal in $K$.  However, since $C_x$ is a corolla and since both vertices in the partially grafted corollas $G$ are extremal, it follows that $x$ is extremal in $K$.
\end{proof}

\section{Deletable Vertices and Internal Edges}
\label{sec:deletable}

In this section we describe a concept about vertices that we will use in part 2 to define one of two types of outer coface maps in the graphical category for connected graphs.  We also discuss a closely related concept that describes an internal edge connecting two distinct vertices.  This concept will be used to define one of two types of inner coface maps in the graphical category for connected graphs. Throughout this section we work with the collection $\gwheelc$ of all connected graphs.

\subsection{Deletable Vertices}

Intuitively, we want to describe a vertex that can be deleted from a connected graph such that the resulting graph is still connected.

\begin{definition}
\label{def:deletable}
Suppose $G$ is a connected graph, and $v$ is a vertex in $G$.  Then we say $v$ is \textbf{deletable} \index{deletable vertex} \index{vertex!deletable} if one of the following two conditions holds.
\begin{enumerate}
\item
$G$ is a permuted corolla with either $\inp(G)$ or $\out(G)$ (or both) non-empty.
\item
The vertex $v$ is loop-free, and there is exactly one internal edge in $G$ that is adjacent to $v$.
\end{enumerate}
In the second case, write $G_v$ for the graph obtained from $G$ by deleting $v$.
\end{definition}

\begin{remark}
The two conditions in the above definition are mutually exclusive because a permuted corolla has no internal edges.  In the second case, $G$ has at least two vertices, so by connectivity $G$ is ordinary.  To simplify the presentation, as in previous sections, in the discussion that follows we will mostly ignore input/output relabelings.  So we will treat the first condition in the above definition as saying $v$ is the unique vertex in a corolla.
\end{remark}

The following observation will justify our terminology.

\begin{lemma}
\label{lem:deletable}
Suppose $G$ is a connected graph that is not a corolla, and $v$ is a deletable vertex in $G$.  Then $G_v$ is connected.
\end{lemma}

\begin{proof}
First note that $G$ has at least two vertices and is ordinary.  So $G_v$ has at least one vertex and is also ordinary.  If $G_v$ has only one vertex, then it is connected.  So suppose $G_v$ has at least two vertices.  Suppose $e$ is the unique internal edge in $G$ that is adjacent to $v$.  Pick two distinct vertices $w$ and $x$ in $G_v$.  Since they are already vertices in $G$, there is an internal path $P$ in $G$ with $w$ and $x$ as end vertices.  Moreover, since there are no other internal edges in $G$ that are adjacent to $v$ besides $e$, the path $P$ cannot contain $e$.  So $P$ is a path in $G_v$, which shows that $G_v$ is connected.
\end{proof}

\begin{remark}
\label{rk:deletable}
In the setting of Lemma \ref{lem:deletable}, we may visualize $G$ in one of two ways.
\begin{center}
\begin{tikzpicture}
\matrix[row sep=1.5cm, column sep=2.5cm]{
\node [plain,label=above:$...$] (v1) {$v$}; &
\node [fatplain,label=above:$...$] (gv2) {$G_v$}; \\
\node [fatplain,label=below:$...$] (gv1) {$G_v$};& 
\node [plain,label=below:$...$] (v2) {$v$};\\
};
\draw [arrow] (gv1) to node{$e$} (v1);
\draw [arrow] (v2) to node{$e$} (gv2);
\foreach \x in {1,2}
{
\draw [inputleg] (gv\x) to +(-.6cm,-.5cm);
\draw [inputleg] (gv\x) to +(.6cm,-.5cm);
\draw [outputleg] (gv\x) to +(-.6cm,.5cm);
\draw [outputleg] (gv\x) to +(.6cm,.5cm);
}
\foreach \x in {1,2}
{
\draw [inputleg] (v\x) to +(-.6cm,-.5cm);
\draw [inputleg] (v\x) to +(.6cm,-.5cm);
\draw [outputleg] (v\x) to +(-.6cm,.5cm);
\draw [outputleg] (v\x) to +(.6cm,.5cm);
}
\end{tikzpicture}
\end{center}
On the left (resp., right), the unique internal edge adjacent to $v$ is incoming (resp., outgoing).  The flag in $G_v$ corresponding to $e$ is a leg.  There are canonical bijections
\[
\begin{split}
\vertex(G) 
&=\vertex(G_v) \coprod \{v\},\\
\edge(G) 
&= \edge(G_v) \coprod \left(\edge(C_v) \setminus \{e\}\right),\\
\edgei(G) 
&= \edgei(G_v) \coprod \{e\},
\end{split}
\]
where as usual $C_v$ is the corolla with the same profiles as the vertex $v$.
\end{remark}

The following observation guarantees the existence of deletable vertices in ordinary simply connected graphs.

\begin{lemma}
\label{deletableexists}
Suppose $G$ is a simply connected graph with at least two vertices.  Then it has at least two deletable vertices.
\end{lemma}

\begin{proof}
Suppose $P = (v_i)_{i=0}^r$ is a maximal path in $G$.  In other words, it is a path (Definition \ref{def:path}) that is not properly contained in any other path.  Since $G$ has at least two vertices, we have $r \geq 1$.  Moreover, since $G$ has no cycles, we have $v_0 \not= v_r$.   Then both end vertices $v_0$ and $v_r$ are deletable.  Indeed, if $v_0$ is not deletable, then since it is loop-free, there is an internal edge $e \not\in P$ connecting $v_0$ and some vertex $u$.  By simple connectivity $u \not= v_i$ for any $i$.  But then the concatenation of $P$ and $e$, with the new initial vertex $u$ replacing $v_0$, is a path that properly contains $P$.  This cannot happen by the maximality assumption on $P$.  Therefore, $v_0$ is deletable, and similarly the terminal vertex $v_r$ is also deletable.
\end{proof}

\begin{remark}
The conclusion of Lemma \ref{deletableexists} does \emph{not} hold in general if $G$ is not assumed to be simply connected.  For example, the non-simply connected graph
\begin{center}
\begin{tikzpicture}
\matrix[row sep=1cm, column sep=2.5cm]{
\node [plain] (v) {$v$};\\
\node [plain] (u) {$u$};\\
};
\draw [arrow, bend left=45] (u) to (v);
\draw [arrow, bend right=45] (u) to (v);
\end{tikzpicture}
\end{center}
has two vertices, neither of which is deletable.
\end{remark}

As in the previous two sections, we want to describe deletable vertices more systematically in terms of graph substitution.  For this purpose, we will use dioperadic graphs (Example \ref{ex:dioperadic}).

\begin{definition}
\label{def:outerdiopfact}
Suppose $G$ is a connected graph.  An \textbf{outer dioperadic factorization} \index{outer dioperadic factorization} \index{dioperadic factorization!outer} of $G$ is a graph substitution decomposition
\[
G = D(\{C_v, H\})
\]
in which 
\begin{itemize}
\item
$D$ is a dioperadic graph with vertices $u$ and $v$, 
\item
a chosen $H$ substituted into $u$ is connected, and
\item
$C_v$ is the corolla with the same profiles as $v$.
\end{itemize}
Call $H$ the \textbf{distinguished subgraph}.\index{distinguished subgraph}
\end{definition}

\begin{remark}
In an outer dioperadic factorization, the distinguished subgraph $H$ has the same profiles as a vertex in a dioperadic graph.  It follows that $H$ must have a leg.  In particular, $H$ cannot be an exceptional loop or the empty graph.
\end{remark}

\begin{theorem}
\label{thm:outerdiopfact}
Suppose $G$ is a connected graph, and $v$ is a vertex in $G$.  Then the following two statements are equivalent.
\begin{enumerate}
\item
$v$ is deletable.
\item
There is an outer dioperadic factorization $G = D(\{C_v,H\})$.\index{deletable vertex!characterization} \index{outer dioperadic factorization!characterization}
\end{enumerate}
\end{theorem}

\begin{proof}
First suppose $v$ is deletable.  If $G$ is a corolla with $\inp(G)$ non-empty, then define $D_1$ as the dioperadic graph in which
\begin{itemize}
\item
the top vertex has the same profiles as $v$,
\item
the bottom vertex $u$ has exactly one incoming flag and one outgoing flag, both with the same color $c$ as the first incoming flag $f$ of $v$, and
\item
the unique internal edge $e$ is adjacent to $v$ via $f$.
\end{itemize}
We may visualize $D_1$ as follows.
\begin{center}
\begin{tikzpicture}
\matrix[row sep=1cm, column sep=2cm]{
\node [plain,label=above:$...$] (v) {$v$}; \\
\node [plain] (u) {$u$};\\
};
\draw [arrow] (u) to node{$e$} node[swap,near end]{$..$} (v);
\draw [inputleg] (u) to +(0,-.7cm);
\draw [inputleg] (v) to +(.6cm,-.5cm);
\draw [outputleg] (v) to +(-.6cm,.5cm);
\draw [outputleg] (v) to +(.6cm,.5cm);
\end{tikzpicture}
\end{center}
Then there is an outer dioperadic factorization
\[
G = D_1(\{C_v, \uparrow_c\})
\]
with distinguished subgraph $\uparrow_c$.  There is a similar argument for the case $\out(G) \not= \varnothing$.

Next suppose $v$ is loop-free, and there is exactly one internal edge $e$ in $G$ that is adjacent to $v$.  We assume that $e$ is an incoming flag of $v$.  By Lemma \ref{lem:deletable} the graph $G_v$ obtained from $G$ by deleting $v$ is connected.  Define $D_e$ as the dioperadic graph in which
\begin{itemize}
\item
the top vertex has the same profiles as $v$,
\item
the bottom vertex $u$ as the same profiles as $G_v$, and
\item
the unique internal edge connects the flags in $u$ and $v$ corresponding to $e$ in $G$.
\end{itemize}
Then there is an outer dioperadic factorization
\[
G = D_e(\{C_v,G_v\})
\]
with distinguished subgraph $G_v$.  There is a similar argument if $e$ is an outgoing flag of $v$.  This proves that $(1) \Longrightarrow (2)$.

For the other direction, suppose there is an outer dioperadic factorization
\[
G = D(\{C_v,H\}).
\]
If $H$ is an exceptional edge, then $G = C_v$, so $v$ is deletable.  Next assume that $H$ is ordinary.  Since $D$ is a dioperadic graph, we may visualize $G$ as in Remark \ref{rk:deletable}.  As a vertex in $G$, $v$ is loop-free, and the unique internal edge in $D$ is the only internal edge in $G$ that is adjacent to $v$.  This shows that $v$ is deletable.
\end{proof}

\begin{remark}
In the proof of $(1) \Longrightarrow (2)$ in Theorem \ref{thm:outerdiopfact}, in the construction of the dioperadic graph $D_1$, instead of the first incoming flag of $v$, we could also have used any other incoming or outgoing flag of $v$.  In other words, if $G$ is a corolla with vertex $v$ and either $\inp(G)$ or $\out(G)$ non-empty, then there is an outer dioperadic factorization
\[
G = D(\{C_v,\uparrow_e\})
\]
for \emph{each} leg $e$ of $G$, where $\uparrow_e$ is the exceptional edge with the same color as the chosen leg $e$.
\end{remark}

\subsection{Internal Edges}

Now we discuss an inner analog of outer dioperadic factorization.  We will use this concept to define one of two types of inner coface maps in the graphical category for connected graphs.

\begin{definition}
\label{def:innerdiopfact}
Suppose $K$ is a connected graph.  An \textbf{inner dioperadic factorization} \index{inner dioperadic factorization} \index{dioperadic factorization!inner} of $K$ is a graph substitution decomposition
\[
K = G(H_w)
\]
in which
\begin{itemize}
\item
$G$ is connected,
\item
$H_w$ is a dioperadic graph with the same profiles as a chosen vertex $w$ in $G$, and
\item
for each vertex $u \not= v$ in $G$, a corolla $C_u$ is substituted into $u$.
\end{itemize}
Call $H_w$ the \textbf{distinguished subgraph}.\index{distinguished subgraph}
\end{definition}

\begin{remark}
If $K=G(H_w)$ is an inner dioperadic factorization, then both $K$ and $G$ are ordinary.  Indeed, $G$ is connected and has at least one vertex $w$, while $K$ is connected and has at least two vertices (namely, the two vertices in $H_w$).
\end{remark}

\begin{remark}
Suppose $K=G(H_w)$ is an inner dioperadic factorization of $K$.  Then there are canonical bijections
\[
\begin{split}
\vertex(K)
&= \vertex(H_w) \coprod \left(\vertex(G) \setminus \{w\}\right),\\
\edgei(K)
&= \edgei(G) \coprod \{e\},\\
\edge(K)
&= \edge(G) \coprod \{e\},
\end{split}
\]
where $e$ is the unique internal edge in the dioperadic graph $H_w$.
\end{remark}

We now observe that inner dioperadic factorizations correspond to internal edges connecting distinct vertices.

\begin{theorem}
\label{thm:innerdiopfact}
Suppose $K$ is a connected graph.  Then there is a canonical bijection between the following two sets.
\begin{enumerate}
\item
The set of internal edges in $K$ that connect two distinct vertices.
\item
The set of inner dioperadic factorizations of $K$.\index{internal edge!characterization} \index{inner dioperadic factorization!characterization}
\end{enumerate}
Moreover, such an internal edge in $K$ corresponds to the unique internal edge in the distinguished subgraph in an inner dioperadic factorization of $K$.
\end{theorem}

\begin{proof}
Let us first describe the desired maps between the two sets.  Suppose $\nicexy{u \ar[r]^{e} & v}$ is an internal edge in $K$ connecting distinct vertices $u$ and $v$.  Define $G$ as the graph obtained from $K$ by shrinking away $e$ and combining the vertices $u$ and $v$ into a single vertex $w$.  To be more precise, $e$ is a $2$-cycle $\{e_{-1},e_1\}$ of the involution $\iota$ of the graph $K$, with $e_{-1} \in u$ and $e_1 \in v$.  From the non-exceptional cells $u$ and $v$ in $K$, we form the new non-exceptional cell
\[
w = \left(u \coprod v\right) \setminus \{e_{-1},e_1\} 
\]
in $G$.  All other non-exceptional cells and structure maps in $G$ are the same as those in $K$, except that the vertex $w$ is given the dioperadic listing.  

Define $H_w$ as the dioperadic graph whose bottom and top vertices have the same profiles as $u$ and $v$, respectively, and whose unique internal edge corresponds to $e$ in $K$.  Then there is an inner dioperadic factorization $K=G(H_w)$.

Conversely, suppose $K=G(H_w)$ is an inner dioperadic factorization, and $\nicexy{u \ar[r]^{e} & v}$ is the internal edge in $H_w$.  Under the graph substitution, $e$ is identified with an internal edge in $K$ connecting the two distinct vertices $u$ and $v$.

By inspection the two maps described above are mutual inverses.
\end{proof}

\section{Disconnectable Edges and Loops}
\label{sec:disconnectable}

In this section, we discuss two concepts about internal edges that we will later use to define the other types of inner and outer coface maps in the graphical category for connected graphs.  As in the previous section, here we work with the collection $\gwheelc$ of all connected graphs.

\subsection{Disconnectable Edges}

First we discuss the concept that corresponds to a type of outer coface maps in the graphical category for connected graphs.  Intuitively, we want to describe an internal edge in a connected graph that can be disconnected (\emph{not} deleted) such that the resulting graph is still connected.

\begin{definition}
\label{def:disconnectable}
Suppose $e$ is an internal edge in a connected graph $G$.  We say $e$ is \textbf{disconnectable} \index{disconnectable edge} \index{edge!disconnectable} if one of the following statements holds.
\begin{enumerate}
\item
$G$ is an exceptional loop.
\item
$\nicexy{u \ar[r]^{e} & v}$ is an ordinary edge, in which $v$ may be equal to $u$, such that there exists a path $P$ in $G$
\begin{itemize}
\item
with end vertices $u$ and $v$, and
\item
that does not contain $e$.
\end{itemize}
\end{enumerate}
In such cases, write $G_e$ for the graph obtained from $G$ by disconnecting $e$.
\end{definition}

\begin{remark}
\label{rk:disconnedge}
In the previous definition, by disconnecting $e$, we mean if $e$ is the $2$-cycle $\{e_{-1},e\}$ of $\iota$ in $G$, then we redefine $\iota$ such that both $e_i$ are $\iota$-fixed points in $G_e$.  In particular, in $G_e$ both flags $e_i$ are legs.  Moreover, there are canonical bijections
\[
\begin{split}
\vertex(G) &= \vertex(G_e),\\
\edge(G) \setminus \{e\} 
&= \edge(G_e) \setminus \{e_{-1},e_1\},\\
\edgei(G) &= \edgei(G_e) \coprod \{e\}.
\end{split}
\]
in which the middle bijection requires $G \not= \wheel$.  Indeed, if $G = \wheel$, then $G_e = \uparrow$, which means that the two flags $\{e_{-1},e\}$ in $G_e$ form a single exceptional edge.  The middle bijection can also be rephrased as
\[
\edge(G) = \frac{\edge(G_e)}{(e_{-1} \sim e_1)},
\]
where on the right-hand side the quotient identifies the legs $e_{\pm 1}$ in $G_e$ to form the internal edge $e$ in $G$. 
\end{remark}

\begin{remark}
A loop at a vertex $v$ is always disconnectable because we can use the trivial path containing only $v$ in Definition \ref{def:disconnectable}.  If $e$ is an ordinary disconnectable edge that is not a loop, then we have a picture like
\begin{center}
\begin{tikzpicture}
\matrix[row sep=1cm, column sep=2cm]{
\node [plain] (u) {$u$}; &
\node [plain] (v) {$v$};\\
};
\draw [arrow] (u) to node{$e$} (v);
\draw [dashedline,bend right=60] (u) to node[swap]{$P$} (v);
\end{tikzpicture}
\end{center}
inside $G$, in which $P$ is an internal path that does not contain $e$.
\end{remark}

The following observation justifies our terminology.

\begin{lemma}
\label{lem:geconnected}
Suppose $e$ is a disconnectable edge in a connected graph $G$.  Then $G_e$ is connected.
\end{lemma}

\begin{proof}
If $G = \wheel$, then $G_e =~ \uparrow$, which is connected.  So suppose  $G$ is ordinary, which implies that $G_e$ is also ordinary.  If $e$ is a loop at $v$, then $G_e$ is connected.

So suppose $P$ is a path in $G$ not containing $\nicexy{u \ar[r]^e & v}$ with distinct end vertices $u$ and $v$.  To see that $G_e$ is connected, it is enough to observe that for any two distinct vertices $x$ and $y$ in $G$, there is a path $Q$ in $G$ with end vertices $x$ and $y$ and that does not contain $e$.  Since $G$ is connected, there must be a path $Q$ in $G$ with end vertices $x$ and $y$.  If $Q$ contains $e$, then we may replace $e$ in $Q$ with $P$, removing some redundant edges and vertices if necessary, to obtain a path $Q'$ with end vertices $x$ and $y$ and that does not contain $e$.
\end{proof}

\begin{remark}
\label{rk:gepicture}
In the context of Lemma \ref{lem:geconnected}, we may visualize $G$ as follows.
\begin{center}
\begin{tikzpicture}
\matrix[row sep=.8cm,column sep=.8cm] {
& \node [emptyvt] (a) {}; &\\
\node [bigplain] (v) {$G_e$}; &&
\node [emptyvt] (b) {}; \\
& \node [emptyvt] (c) {}; &\\
};
\draw [outputleg] (v) to +(-.6cm,.5cm);
\draw [outputleg] (v) to +(.6cm,.5cm);
\draw [inputleg] (v) to +(-.6cm,-.5cm);
\draw [inputleg] (v) to +(.6cm,-.5cm);
\draw [arrow, looseness=25, in=-60, out=60, loop] (v) to node[near start]{$e_{-1}$} node[near end]{$e_1$} ();
\end{tikzpicture}
\end{center}
Note that $e$ is an internal edge in $G$ but not in $G_e$, in which the flags $e_{\pm 1}$ are both legs.
\end{remark}

The following observation gives an alternative characterization of an ordinary disconnectable edge as an edge in a cycle.

\begin{lemma}
\label{disconnectablecycle}
Suppose $\nicexy{u \ar[r]^e & v}$ is an ordinary edge in a connected graph $G$, in which $v$ may be equal to $u$.  Then the following statements are equivalent.
\begin{enumerate}
\item
$e$ is disconnectable.
\item
There exists a cycle in $G$ that contains $e$.\index{disconnectable edge!characterization}
\end{enumerate}
\end{lemma}

\begin{proof}
First suppose $e$ is disconnectable.  If $e$ is a loop at $v$, then $e$ is the desired cycle $P$.  So suppose $e$ is not a loop.  Since $G_e$ is connected (Lemma \ref{lem:geconnected}), there is a path $Q$ in $G_e$ with $u$ and $v$ as end vertices.  Since $e$ is \emph{not} an ordinary edge in $G_e$, it is not contained in $Q$.  Therefore, the concatenation of $Q$ and $e$ is a cycle in $G$ that contains $e$.

Conversely, suppose there exists a cycle $P$ in $G$ that contains $e$.  If $P$ contains only $e$, then $e$ is a loop at $v$, which means it is disconnectable.  On the other hand, if $P$ contains at least two ordinary edges, then removing $e$ from $P$ and cyclically relabeling the other edges in $P$ if necessary, the resulting internal path $P'$ has end vertices $u$ and $v$ and does not contain $e$. Therefore, $e$ is disconnectable.
\end{proof}

As usual, we want to describe disconnectable edges in terms of graph substitution.

\begin{definition}
\label{def:outercontfact}
Suppose $G$ is a connected graph.  An \textbf{outer contracting factorization} \index{outer contracting factorization} \index{contracting factorization!outer} of $G$ is a graph substitution decomposition
\[
G = (\xi_e C)(H)
\]
in which
\begin{itemize}
\item
$\xi_eC$ is a contracted corolla with internal edge $e$ (Example \ref{ex:contractedcor}), and
\item
$H$ is connected.
\end{itemize}
Call $H$ the \textbf{distinguished subgraph}.\index{distinguished subgraph}
\end{definition}

\begin{theorem}
\label{thm:outercontfact}
Suppose $e$ is an internal edge in a connected graph $G$.  Then the following two statements are equivalent.
\begin{enumerate}
\item
$e$ is disconnectable.
\item
There exists an outer contracting factorization $G = (\xi_eC)(H)$.\index{disconnectable edge!characterization} \index{outer contracting factorization!characterization}
\end{enumerate}
\end{theorem}

\begin{proof}
Suppose $e$ is disconnectable.  If $G = \wheel_c$, then there is an outer contracting factorization
\[
\wheel_c = \left(\xi_e C_{(c;c)}\right)(\uparrow_c)
\]
with distinguished subgraph $\uparrow_c$.  On the other hand, suppose $\nicexy{u \ar[r]^e & v}$ is ordinary as in Definition \ref{def:disconnectable}.  If $C$ is the corolla with the same profiles as $G_e$, then there is an outer contracting factorization
\[
G = (\xi_eC)(G_e)
\]
with distinguished subgraph $G_e$.

Conversely, suppose there exists an outer contracting factorization 
\[
G = (\xi_eC)(H).
\]
The vertex $v$ in the corolla $C$ must have non-empty inputs and non-empty outputs.  Since $H$ has the same profiles as $v$, $H$ cannot be an exceptional loop or the empty graph.  If $H$ is an exceptional edge, then $G$ is an exceptional loop, which means $e$ is disconnectable.  Next suppose $H$ is  ordinary.  By connectivity of $H$, there is an internal path $P$ in $H$ with end vertices $u$ and $v$.  Since $e$ is \emph{not} an internal edge in $H$, $e$ is not in $P$.  Thus, $P$ is an internal path in $G$ with end vertices $u$ and $v$ that does not contain $e$.  Therefore, $e$ is disconnectable.
\end{proof}

\subsection{Loops}

Now we discuss an inner analog of outer contracting factorization.  We will use this concept to define the other type of inner coface maps in the graphical category for connected graphs.

\begin{definition}
\label{def:innercontfact}
Suppose $K$ is a connected graph.  An \textbf{inner contracting factorization} \index{inner contracting factorization} \index{contracting factorization!inner} of $K$ is a graph substitution decomposition
\[
K = G(H_w)
\]
in which
\begin{itemize}
\item
$G$ is connected,
\item
$H_w = \xi_e C_v$ is a contracted corolla with the same profiles as a chosen vertex $w$ in $G$, and
\item
for each vertex $u \not= v$ in $G$, a corolla $C_u$ is substituted into $u$.
\end{itemize}
Call $H_w$ the \textbf{distinguished subgraph}.\index{distinguished subgraph}
\end{definition}

\begin{remark}
If $K=G(H_w)$ is an inner contracting factorization, then both $K$ and $G$ are ordinary because each of them has at least one vertex.
\end{remark}

\begin{remark}
Suppose $K=G(H_w)$ with $H_w = \xi_e C_v$ is an inner contracting factorization of $K$.  Then there are canonical bijections
\[
\begin{split}
\vertex(K)
&= \{v\} \coprod \left(\vertex(G) \setminus \{w\}\right),\\
\edgei(K)
&= \edgei(G) \coprod \{e\},\\
\edge(K)
&= \edge(G) \coprod \{e\}.
\end{split}
\]
Here $e$ is the unique internal edge in the contracted corolla $H_w = \xi_eC_v$, and $v$ is the unique vertex in $H_w$.
\end{remark}

We now observe that inner contracting factorizations correspond to loops.

\begin{theorem}
\label{thm:innercontfact}
Suppose $K$ is a connected graph.  Then there is a canonical bijection between the following two sets.
\begin{enumerate}
\item
The set of loops in $K$.
\item
The set of inner contracting factorizations of $K$.\index{loop!characterization} \index{inner contracting factorization!characterization}
\end{enumerate}
Moreover, such a loop in $K$ corresponds to the unique internal edge in the distinguished subgraph in an inner contracting factorization of $K$.
\end{theorem}

\begin{proof}
Let us first describe the desired maps between the two sets.  First suppose $e=\{e_{-1},e_1\}$ is a loop at $v$ in $K$.  Define $G$ as the graph obtained from $K$ by \emph{deleting} the loop $e$.  In other words, from the ordinary cell $v$ in $K$, we form a new ordinary cell
\[
w = v \setminus \{e_{-1},e_1\}
\]
in $G$, which is still connected.  Define $H$ as the contracted corolla $\xi_eC_v$ whose unique vertex has the same profiles as $v$ in $K$ and whose unique internal edge corresponds to $e$ in $K$.  Then there is an inner contracting factorization
\[
K=G(H)
\]
with distinguished subgraph $H$.

Conversely, suppose given an inner contracting factorization
\[
K=G(\xi_eC_v).
\]
Then the internal edge $e$ in the distinguished subgraph $\xi_eC_v$ becomes a loop at $v$ in $K$.

Finally, observe that the two maps defined above are mutual inverses.
\end{proof}

The following table provides a summary of the four types of graph substitution factorizations discussed in this and the previous sections.

\begin{center}
\begin{tabular}{|c|c|c|}\hline
& Graph substitution & $G$ is 
\\ \hline \hline
\shortstack{outer\\dioperadic\\factorization} 
& \shortstack{$K=D(\{C_v,G\})$\\$D$ dioperadic\\Theorem \ref{thm:outerdiopfact}} 
& \shortstack{$K$ with a deletable\\vertex $v$ deleted}
\\ \hline
\shortstack{inner\\dioperadic\\factorization} 
& \shortstack{$K=G(D)$\\$D$ dioperadic\\Theorem \ref{thm:innerdiopfact}}
& \shortstack{$K$ with an internal edge\\connecting
two distinct\\vertices shrunk away}
\\ \hline
\shortstack{outer\\contracting\\factorization} 
& \shortstack{$K=(\xi_eC)(G)$\\$\xi_eC$ contracted corolla\\Theorem \ref{thm:outercontfact}} 
& \shortstack{$K$ with a disconnectable\\edge $e$ disconnected}
\\ \hline
\shortstack{inner\\contracting\\factorization} 
& \shortstack{$K=G(\xi_eC)$\\$\xi_eC$ contracted corollas\\Theorem \ref{thm:innercontfact}}
& \shortstack{$K$ with a loop $e$ at a\\vertex deleted}
\\ \hline
\end{tabular}
\end{center}


\chapter{Properads}
\label{ch:properads}

\abstract*{We recall both the biased and the unbiased definitions of a properad.  The former describes a properad in terms of generating operations, namely, units, $\Sigma$-bimodule structure, and properadic composition.  The unbiased definition of a properad describes it as an algebra over a monad induced by connected wheel-free graphs.  The equivalence of the two definitions of a properad are proved in detail in \cite{jy2} as an example of a general theory of generating sets for graphs.}

This chapter is a brief introduction to properads.  We recall both the biased (section \ref{sec:biasedproperad}) and the unbiased (section \ref{sec:unbiasedproperad}) descriptions of a properad.  These are two equivalent ways to define a properad.  We emphasize that what we call a properad here is sometimes called a \emph{colored} properad in the literature.

Properads are objects that effectively parametrize operations with multiple inputs, multiple outputs, symmetric group actions, units, and associativity axioms along connected wheel-free graphs.   One-colored properads were first introduced by Vallette \cite{vallette} in the linear setting.  These objects are general enough to describe, for example,  biassociative bialgebras, Lie bialgebras, and (co)module bialgebras as algebras over suitable properads.  Properads are more general than operads \cite{may} in the sense that the latter are properads whose operations have only one output.  Also, properads are more general than dioperads \cite{gan} in the sense that the dioperadic composition is also a properadic operation, but the converse is not true.

The biased version of a properad describes it as a suitably parametrized set of objects with some extra structures, namely, units, symmetric group actions, and a properadic composition, satisfying suitable axioms.  This is similar in spirit to the original definition of an operad given by May \cite{may}.  In even more familiar terms, a biased properad is similar to the usual definition of a category, where the properadic composition generalizes the categorical composition of two morphisms.

The unbiased version of a properad is more formal and describes it as an algebra over a monad associated to connected wheel-free graphs $\gupc$.  This is the free properad monad.  Free properads are needed not only to define the graphical category of connected wheel-free graphs, but also to define the symmetric monoidal closed structure on the category of properads.

One main difference between our free properads and those in, for example, \cite{vallette} is that we take as our underlying object a suitably parametrized set of objects \emph{without} symmetric group actions.  The symmetric group actions on a properad are generated by certain structures on connected wheel-free graphs called a \emph{listing}.  This approach to the free properads and other variants is developed in \cite{jy2}.  Also, the proof of the equivalence between the biased and the unbiased descriptions of properads is not trivial.  The full detail is given in \cite{jy2}, along with many other variants of operads and properads.

Eventually we will work over the category of sets.  However, throughout this chapter, we work more generally over a symmetric monoidal category $(\catc, \otimes, I)$ with all small colimits and initial object $\varnothing$ such that $\otimes$ commutes with colimits on both sides.  Unless otherwise specified, the reference for this chapter is \cite{jy2}, where all the details can be found.

Everything in this chapter about properads has obvious analogs for properads with non-empty inputs or non-empty outputs.  Instead of restating everything for these close variants, we point out the simple modifications in Remarks \ref{rk:properadio} and \ref{rk:unbiasedtable}.

\section{Biased Properads}
\label{sec:biasedproperad}

In this section, we recall the definition of a properad in biased form, algebras over a properad, and a few examples.  A properad is a mechanism for organizing operations with multiple inputs and multiple outputs.  We can think of each input/output as a color.   The reader may wish to review Definition \ref{def:colors} on colors and profiles, which we will use below.

\subsection{\texorpdfstring{$\Sigma$}{Σ}-Bimodules and Colored Objects}

\begin{definition}
\label{def:coloredobject}
Fix a set $\fC$ of colors.
\begin{enumerate}
\item
The category of \textbf{$\Sigma_{\SC}$-bimodules} \index{sigma bimodule} is the diagram category $\catc^{\SC}$.
\item
The discrete category associated to $\SC$ is written as $\dis(\SC)$ or $\dis(\sS)$\label{note:dissc} if $\fC$ is clear from the context.
\item
An \textbf{$\SC$-colored object}, or simply a \index{colored object}\textbf{colored object}, is an object in the diagram category $\catc^{\dis(\SC)} = \prod_{\SC} \catc$.
\item
A colored object $\sP$ is said to be \textbf{special} \index{special colored object} if the components $\sP\binom{\varnothing}{\uc}$ and $\sP\binom{\ud}{\varnothing}$ are all equal to the initial object $\varnothing$ for all $\fC$-profiles $\uc$ and $\ud$.
\item
Suppose $\sP$ is an $\SC$-colored object and $\sQ$ is an $\SD$-colored object.  Then a \textbf{map} $\nicexy{\sP \ar[r]^-{f} & \sQ}$ of colored objects consists of:
\begin{itemize}
\item
a function $\nicexy{\fC \ar[r]^-{f_0} & \fD}$ on color sets, and
\item
a map
\[
\nicexy{
\sP\dc \ar[r]^-{f_1} & \sQ\binom{f_0\ud}{f_0\uc},
}\]
where $f_0\uc = \left(f_0(c_1), \ldots , f_0(c_m)\right)$.
\end{itemize} 
\end{enumerate}
\end{definition}

\begin{remark}
\label{rk:coloredob}
A colored object $\sP \in \catc^{\dis(\SC)}$ consists of a set of objects $\sP\dc \in \catc$, one for each pair of $\fC$-profiles $\dch \in \SC$.  We call $\uc$ (resp., $\ud$) the \textbf{input profile} \index{input profile} (resp., \textbf{output profile}) \index{output profile} of the component $\sP\dc$.  We think of $\sP\dc$ as consisting of operations with inputs $\uc = (c_1,\ldots,c_m)$ and outputs $\ud = (d_1,\ldots,d_n)$.  When $\sP\dc$ has an underlying set, we can visualize an element $p$ in it as a decorated corolla
\begin{center}
\begin{tikzpicture}
\matrix[row sep=1cm,column sep=1cm] {
\node [plain,label=above:$...$,label=below:$...$] (p) {$p$}; \\
};
\draw [outputleg] (p) to node[above left=.1cm]{$d_1$} +(-.6cm,.5cm);
\draw [outputleg] (p) to node[above right=.1cm]{$d_n$} +(.6cm,.5cm);
\draw [inputleg] (p) to node[below left=.1cm]{$c_1$} +(-.6cm,-.5cm);
\draw [inputleg] (p) to node[below right=.1cm]{$c_m$} +(.6cm,-.5cm);
\end{tikzpicture}
\end{center}
with one vertex decorated by $p$, $|\uc|$ input legs colored by the $c_i$, and $n$ output legs colored by the $d_j$.
\end{remark}

\begin{remark}
\label{rk:sigmabimodule}
On the other hand, a $\Sigma_{\SC}$-bimodule is a colored object $\sP$ together with isomorphisms
\[
\nicexy{
\sP\dc \ar[r]^-{(\tau;\sigma)}_-{\cong} & \sP\dcsigma
}\]
with $\sigma \in \Sigma_{|\ud|}$ and $\tau \in \Sigma_{|\uc|}$ such that
\begin{itemize}
\item
$(\id;\id)$ is the identity map, and
\item
$(\tau';\sigma') \circ (\tau;\sigma) = (\tau\tau';\sigma'\sigma)$.
\end{itemize}
A map $f \colon \sP \to \sQ$ of $\Sigma_{\SC}$-bimodules is a map of the underlying $\SC$-colored objects with $f_0 = \Id$ such that all the squares
\[
\nicexy{
\sP\dc \ar[r]^-{f} \ar[d]_{(\tau;\sigma)} & \sQ\dc \ar[d]^{(\tau;\sigma)}\\
\sP\dcsigma \ar[r]^-{f} & \sQ\dcsigma
}\]
are commutative.
\end{remark}

\begin{example}
Consider the $1$-colored case with $\fC = \{*\}$.  Then $\pofc$ is the groupoid $\bN$ whose objects are non-negative integers $\{0, 1, 2, \ldots \}$ and whose only morphisms are the symmetric groups $\bN(n,n) =  \Sigma_n$ for $n \geq 0$.  A colored object is a double sequence of objects $\sP = \{\sP(m;n)\}_{m,n \geq 0}$.  A $\Sigma_{\sS}$-bimodule is a colored object $\sP$ together with isomorphisms
\[
\nicexy{
\sP\nm \ar[r]^-{(\tau;\sigma)}_{\cong} & \sP\nm
}\]
satisfying the two conditions stated in Remark \ref{rk:sigmabimodule}.
\end{example}

\subsection{Biased Definition of a Properad}

The following definition of a properad will be used in most of the later chapters in Part I.

\begin{definition}
\label{def:biasedproperad}
Let $\fC$ be a set of colors.
\begin{enumerate}
\item
A \textbf{$\fC$-colored properad} \index{properad} $(\sP, \bone, \boxtimes)$ consists of:
\begin{enumerate}
\item
a $\Sigma_{\SC}$-bimodule $\sP$,
\item
a \textbf{$c$-colored unit} \index{colored unit}
\[\label{note:coloredunit}
\nicexy{
I \ar[r]^-{\bone_c} & \sP\ccsingle
}\]
for each $c \in \fC$, and
\item
a \textbf{properadic composition} \index{properadic composition}
\[\label{note:propcomp}
\nicexy{
\sP\dc \otimes \sP\ba \ar[r]^-{\boxtimes^{\uc'}_{\ub'}} 
& \sP\binom{\ub \circ_{\ub'} \ud}{\uc \circ_{\uc'} \ua}
}\]
whenever $\uc' \subseteq \uc$ and $\ub' \subseteq \ub$ are equal $k$-segments for some $k>0$.
\end{enumerate}
These structures are required to satisfy suitable bi-equivariant, unity, and associativity axioms.
\item
A \textbf{morphism} $\nicexy{\sP \ar[r]^-{f} & \sQ}$ from a $\fC$-colored properad $\sP$ to a $\fD$-colored properad $\sQ$ is a map of the underlying colored objects that respects the bi-equivariant structure, colored units, and properadic compositions.
\item
A \textbf{properad} is a $\fC$-colored properad for some color set $\fC$.
\item
Denote by $\properad$\label{note:properad} the category of all properads and morphisms.
\end{enumerate}
\end{definition}

\begin{definition}
\label{rk:properadio}
\textbf{Properads with non-empty inputs} (resp., \textbf{properads with non-empty outputs} \index{properad!with non-empty inputs} \index{properad!with non-empty outputs} and \textbf{special properads}) \index{properad!special} \index{special properad} are defined just like properads, except that in a $\Sigma_{\SC}$-bimodule we replace $\SC$ by the full subcategory $\SC_i$ (resp., $\SC_o$ and $\SC_{s}$) consisting of pairs of profiles $\dch$ with $\uc \not=\varnothing$ (resp., $\ud \not= \varnothing$ and $\uc \not= \varnothing \not= \ud$).  The category of properads with non-empty inputs (resp., properads with non-empty outputs and special properads) is denoted by $\properadi$ (resp., $\properado$ and $\properads$).\label{note:properadi}
\end{definition}

\begin{remark}
The category $\properadi$ (resp., $\properado$ and $\properads$) is canonically isomorphic to the full subcategory of $\properad$ consisting of properads $\sP$ with $\sP\dc = \varnothing$ whenever $\uc = \varnothing$ (resp., $\ud = \varnothing$, and either $\uc = \varnothing$ or $\ud = \varnothing$).  
\end{remark}

\begin{remark}
When the color set $\fC$ is clear from the context, we will omit mentioning it.  In what follows, we usually abbreviate the properadic composition $\boxtimes^{\uc'}_{\ub'}$ to just $\boxtimes$. 
\end{remark}

\begin{remark}
In Definition \ref{def:biasedproperad} if we insist that $\ub' = \uc'$ are equal $1$-segments, then the resulting object is exactly a dioperad.  In the linear setting, a $1$-colored dioperad was introduced in \cite{gan}.  Moreover, if we further insist that $\sP\dc = \varnothing$ unless $|\ud| = 1$, then the resulting structure is equivalent to an operad.  A $1$-colored operad in the topological setting was introduced in \cite{may}. 
\end{remark}

\begin{remark}
\label{rk:visualizeunits}
Following Remark \ref{rk:coloredob} we visualize the $c$-colored unit in a properad as the $c$-colored exceptional edge
\begin{center}
\begin{tikzpicture}
\matrix[row sep=1cm,column sep=1cm] {
\node [empty] (t) {}; \\
\node [empty] (b) {};\\
};
\draw [arrow] (b) to node[swap,near start]{$c$} (t);
\end{tikzpicture}
\end{center}
that contains no vertices. The properadic composition $\boxtimes^{\uc'}_{\ub'}$ is visualized as the assignment
\begin{center}
\begin{tikzpicture}
\matrix[row sep=.2cm,column sep=1.5cm] {
\node [plain,label=below:$\uc'$,label=above:$\ud$] (p1) {$p$}; &&& \\
\node [empty]  {...}; & 
\node [empty] (s) {}; & \node [empty] (t) {}; &
\node [fatplain,label=below:$\uc \circ_{\uc'} \ua$,label=above:$\ub \circ_{\ub'} \ud$] (p3) {$p \boxtimes q$}; \\
\node [plain,label=below:$\ua$,label=above:$\ub'$] (p2) {$q$}; &&& \\
};
\draw [mapto] (s) to (t);
\foreach \x in {1,2,3}
{
\draw [outputleg] (p\x) to +(-.6cm,.4cm);
\draw [outputleg] (p\x) to +(.6cm,.4cm);
\draw [inputleg] (p\x) to +(-.6cm,-.4cm);
\draw [inputleg] (p\x) to +(.6cm,-.4cm);
}
\draw [arrow,bend left=40] (p2) to (p1);
\draw [arrow,bend right=40] (p2) to (p1);
\end{tikzpicture}
\end{center}
for $p \in \sP\dc$ and $q \in \sP\ba$ with $\uc \supseteq \uc' = \ub' \subseteq \ub$, although the components of $\sP$ need not have underlying sets.
\end{remark}

\begin{remark}
We refer to the above definition of a properad as \emph{biased}, in the sense that only some of the generating structure maps and generating axioms are involved.  This biased definition is similar to the usual definition of an associative algebra, where one states that there is a binary operation $A \otimes A \to A$, satisfying one associativity axiom
\[
\left(a_1 a_2\right) a_3= a_1\left(a_2 a_3\right).
\]
However, in an associativie algebra there are actually many more operations, such as
\[
\left(a_1,a_2,a_3,a_4\right) \longmapsto 
a_1 \left[ (a_2 a_3) a_4\right].
\]
An \emph{unbiased} definition of an associative algebra should involve all such $n$-ary operations, together with all the relations among them.
\end{remark}

\subsection{Algebras over a Properad}

If a properad is regarded as an object that parametrizes operations with multiple inputs and multiple outputs, then the objects on which these operations act are called algebras.  For the following definition, we also require the symmetric monoidal category $\catc$ to be closed with internal hom $\Hom(-,-)$.  The endomorphism properad below is the canonical example of a properad.

\begin{definition}
\label{def:endobject}
Fix a set $\fC$ of colors.
\begin{enumerate}
\item
Regard $\fC$ as a discrete category with only identity morphisms.  Objects in the diagram category $\catc^{\fC} = \prod_{\fC} \catc$ are called \textbf{$\fC$-colored objects}.
\item
For a $\fC$-colored object $X = \{X_c\}_{c\in \fC}$ with each $X_c \in \catc$ and a $\fC$-profile $\uc = c_{[1,m]}$, define the object
\[
X_{\uc} = X_{c_1} \otimes \cdots \otimes X_{c_m},
\]
with $X_{\varnothing} = I$.
\item
The \textbf{endomorphism properad} \index{endomorphism properad} \index{properad!endomorphism} $\End(X)$ of a $\fC$-colored object $X$ is defined as the $\fC$-colored properad with components
\[
\End(X)\dc = \Hom\left(X_{\uc}, X_{\ud}\right).
\]
The bi-equivariant structure is induced by permutations of tensor factors and the symmetric monoidal structure, while the colored units are adjoints to the identity maps on the $X_c$.  The properadic composition is induced by categorical composition, the colored units, the bi-equivariant structure, and tensor products of maps.
\item
Suppose $\sP$ is a $\fC$-colored properad, and $X$ is a $\fC$-colored object.  Then a \textbf{$\sP$-algebra} \index{algebra} \index{properad!algebra} structure on $X$ is a morphism
\[
\nicexy{\sP \ar[r]^-{\lambda} & \End(X)}
\]
of properads, called the \textbf{structure map}, such that $\lambda_0 = \Id$.
\end{enumerate}
\end{definition}

\begin{remark}
The adjoints of the components of the structure map take the form
\[
\nicexy{
\sP\dc \otimes X_{\uc} \ar[r]^-{\lambda} & X_{\ud}.
}\]
This is the reason why $\sP\dc$ is regarded as parametrizing operations with inputs $\uc$ and outputs $\ud$. 
\end{remark}

\section{Unbiased Properads}
\label{sec:unbiasedproperad}

In this section, we discuss the unbiased definition of a properad.  This requires a discussion of the free properad monad, which in turn uses the notion of connected wheel-free graphs from section \ref{sec:cwfree}.

\subsection{Monads and their Algebras}

Before we discuss the free properad monad, let us first recall the definitions of a monad and of an algebra over a monad \cite{maclane98}.

\begin{definition}
A \textbf{monad}\index{monad} $(T,\mu,\nu)$ on a category $\catc$ consists of
\begin{itemize}
\item
a functor $\nicexy{\catc \ar[r]^-{T} & \catc}$\label{note:T} and
\item
natural transformations $\nicexy{
T^2 \ar[r]^-{\mu} & T}$ and $\nicexy{\Id \ar[r]^-{\nu} & T}$, called the \textbf{multiplication} and the \textbf{unit}, respectively, 
\end{itemize}
such that the following associativity and unity diagrams are commutative.
\begin{equation}
\label{monadaxioms}
\nicearrow
\xymatrix@C+12pt{
T^3 \ar[r]^{T\mu} \ar[d]_{\mu T} & T^2 \ar[d]^{\mu} & & T \ar[r]^{T \nu} \ar[dr]_{\Id} &
	T^2 \ar[d]^{\mu} & T \ar[l]_{\nu T} \ar[dl]^{\Id} \\
T^2 \ar[r]_{\mu} & T & & & T &
}
\end{equation}
\end{definition}

\begin{definition}
Suppose $(T,\mu,\nu)$ is a monad on $\catc$.
\begin{enumerate}
\item
A \textbf{$T$-algebra}\index{monadic algebra} 
$(X,\gamma)$\label{note:algebrat} consists of
\begin{itemize}
\item
an object $X \in \catc$ and
\item
a morphism $\nicexy{TX \ar[r]^-{\gamma} & X \in \catc}$, called the \textbf{structure map}, 
\end{itemize}
such that the following associativity and unity diagrams are commutative.
\[
\nicearrow
\xymatrix@C+12pt{
T^2X \ar[r]^{T\gamma} \ar[d]_{\mu_X} & TX \ar[d]^{\gamma} & & X \ar[r]^{\nu_X} \ar[dr]_{\Id} &
	TX \ar[d]^{\gamma}\\
TX \ar[r]_{\gamma} & X & & & X
}
\]
\item
If $(Y,\gamma^Y)$ is another $T$-algebra, then a 
\textbf{morphism}\index{monadic algebra!morphism of} of $T$-algebras
\[
\nicexy{(X,\gamma^X) \ar[r]^-{f} & (Y,\gamma^Y)}
\]
is a morphism $\nicexy{X \ar[r]^-{f} & Y}$ in $\catc$ such that the square
\[
\nicearrow
\xymatrix@C+10pt{
TX \ar[r]^{T(f)} \ar[d]_{\gamma^X} & TY \ar[d]^{\gamma^Y}\\
X \ar[r]_{f} & Y
}
\]
commutes.  The category of algebras over $T$ will be denoted 
$\alg(T)$\label{note:talg}.
\end{enumerate}
\end{definition}

\begin{example}
An adjoint pair $L \colon \catc \adjoint \catd \colon R$, with $L$ the left adjoint, defines a monad on $\catc$, whose underlying functor is $RL$.
\end{example}

\subsection{Decorated Graphs}

The free properad monad is made up of components that appear in the following definition.  Intuitively, we want to decorate the vertices of a graph with appropriate components of a colored object.  

Fix a set $\fC$ of colors, so $\sS = \SC$ and $\catcs = \catc^{\dis(\SC)}$.

\begin{definition}
\label{def:decoratedgraph}
Given $\sP \in \catcs$ and a graph $G$, define the \textbf{$\sP$-decorated graph}\index{decorated graph} 
as the unordered tensor product
\begin{equation}
\label{wheeldecoration}
\sP[G] = \bigotimes_{v \in \vertex(G)} \sP\profilev.
\end{equation}
If $\vertex(G)$ is empty, then $\sP[G]=I$, the unit of the tensor
product.
\end{definition}

\begin{remark}
Note that the $\sP$-decorated graph \eqref{wheeldecoration} is an \emph{unordered} tensor product.  This makes sense because the set of vertices in a graph is not ordered.
\end{remark}

\begin{example}
\label{ex:pofg}
For the $c$-colored exceptional edge $\uparrow_c$, we have
\[
\sP[\uparrow_c] = I.
\]
For the permuted corolla $\sigma C_{\dch}\tau$, we have
\[
\sP[\sigma C_{\dch} \tau] = \sP\dc.
\] 
For the partially grafted corollas $C_{\dch} \boxtimes^{\uc'}_{\ub'} C_{\bah}$, we have
\[
\sP\left[C_{\dch} \boxtimes^{\uc'}_{\ub'} C_{\bah}\right] \cong \sP\dc \otimes \sP\ba.
\]
\end{example}

\subsection{Free Properad Monad}

Recall from Definition \ref{def:setsofgraphs} that $\gupc$ is the set of connected wheel-free graphs.

\begin{definition}
\label{def:gupdc}
For a pair of $\fC$-profiles $\dc$, let $\gupc\dc$ be the subset of $\gupc$ consisting of those connected wheel-free graphs $G$ that satisfy $\profileg = \dc$.
\end{definition}

\begin{remark}
There is a disjoint union decomposition
\[
\gupc = \coprod_{\dc \in \sS} \gupc\dc.
\]
Moreover, the subset $\gupc\dc$ contains at least the permuted corollas
\[
\sigma C_{(\uc\tau^{-1}; \sigma^{-1}\ud)} \tau,
\]
and in particular the $\dch$-corolla $C_{\dch}$.
\end{remark}

Now we define the free properad monad.

\begin{definition}
\label{def:fgupc}
Suppose $\sP \in \catcs$.
\begin{enumerate}
\item
Define the functor \index{properad!free} \index{free properad}
\[
F = F_{\gupc}\index{monad!associated to a pasting scheme} \colon \catcs \longrightarrow \catcs
\]
by
\begin{equation}
\label{fgx}
F\sP\dc = \coprod_{G \in \gupc\dc} \sP[G] 
= \coprod_{G \in \gupc\dc} \bigotimes_{v \in \vertex(G)} \sP\profilev.
\end{equation}
for $\dch \in \sS$.
\item
Define the natural transformation $\nicexy{F^2 \ar[r]^-{\mu} & F}$ as the one induced by graph substitution.
\item
Define the natural transformation $\nicexy{\Id \ar[r]^-{\nu} & F}$ as the one induced by the $\dch$-corollas as $\dch$ runs through $\sS$.
\end{enumerate}
\end{definition}

The associativity and unity properties of graph substitution imply the following observation.

\begin{theorem}
\label{thm:freepropmonad}
For each non-empty set $\fC$, there is a monad $\left(F_{\gupc},\mu,\nu\right)$ on $\catcs$.
\end{theorem}

\begin{remark}
\label{rk:falgstructure}
Recall that each monad, such as $F=F_{\gupc}$, has associated algebras.  An $F$-algebra $\sP$ has a structure map $\nicexy{F\sP \ar[r]^-{\gamma} & \sP \in \catcs}$, which consists of component maps
\[
\nicexy{
\sP[G] = \bigotimes_{\ving} \sP\profilev \ar[r]^-{\gamma_G} & \sP\profileg
}\]
for $G \in \gupc$.  For example:
\begin{enumerate}
\item
For the permuted corolla, we have the component structure map
\[
\nicexy{
\sP\left[\sigma C_{\dch} \tau\right] = \sP\dc \ar[r]^-{\gamma_{\sigma C\tau}} & \sP\dcsigma.
}\]
\item
For the $c$-colored exceptional edge, we have the component structure map
\[
\nicexy{
\sP[\uparrow_c] = I \ar[r]^-{\gamma_{\uparrow}} & \sP\ccsingle.
}\]
\item
For the partially grafted corollas, we have the component structure map
\[
\nicexy{
\sP\left[C_{\dch} \boxtimes^{\uc'}_{\ub'} C_{\bah}\right] 
\cong 
\sP\dc \otimes \sP\ba \ar[r]^-{\gamma} 
& \sP\binom{\ub \circ_{\ub'} \ud}{\uc \circ_{\uc'} \ua}.
}\]
\end{enumerate}
These component structure maps correspond to the $\Sigma_{\sS}$-bimodule structure, the $c$-colored unit, and the properadic composition, respectively, of a $\fC$-colored properad.  The required bi-equivariant, unity, and associativity axioms are implied by the $F$-algebra axioms.  The converse is also true.  Namely, a $\fC$-colored properad uniquely determines an $F$-algebra structure on the underlying $\sS$-colored object.  So we have the following observation.
\end{remark}

\begin{corollary}
\label{cor:unbiasedequiv}
There is a natural bijection between $\fC$-colored properads and $F_{\gupc}$-algebras.\index{properad!unbiased}
\end{corollary}

By general category theory, for a colored object $\sP$, the free $F$-algebra $F\sP$ is precisely the free properad of $\sP$.

\begin{remark}
Following Markl \cite{markl08}, we refer to the $F$-algebra description of a properad as \emph{unbiased} in the sense that all connected wheel-free graphs are used in the monad $F$.  
\end{remark}

\begin{remark}
\label{rk:unbiasedtable}
With similar definitions, Theorem \ref{thm:freepropmonad} and Corollary \ref{cor:unbiasedequiv} also hold for the other sets of graphs in Definition \ref{def:setsofgraphs}.  In fact, Theorem \ref{thm:freepropmonad} holds for any \emph{pasting scheme} \index{pasting scheme} as defined in \cite{jy2}, while Corollary \ref{cor:unbiasedequiv} holds for most variants of operads and PROPs.  For a fixed set $\fC$ of colors and suitable subcategories of $\SC$, the following table provides the objects associated to each set of graphs in Definition \ref{def:setsofgraphs}.
\begin{center}
\begin{tabular}{|c|c|c|}\hline
$\fC$-colored & $F_{?}$-algebras & Types of Graphs \\ \hline\hline
categories & $\ULin$ & linear graphs \\ \hline
operads & $\uoperad$ & unital trees \\ \hline
dioperads & $\gupd$ & simply connected graphs \\ \hline
properads & $\gupc$ & connected wheel-free graphs \\ \hline
special properads & $\gupcs$ & \shortstack{special connected\\ wheel-free graphs} \\ \hline
\shortstack{properads with\\non-empty inputs} & $\gupci$ &
\shortstack{connected wheel-free graphs\\with non-empty inputs}\\  \hline
\shortstack{properads with\\non-empty outputs} & $\gupco$ &
\shortstack{connected wheel-free graphs\\with non-empty outputs}\\  \hline
wheeled properads & $\gwheelc$ & connected graphs\\ \hline
\end{tabular}
\end{center}
Wheeled properads will be discussed in part 2.  Note that an $F_{\ULin}$-algebra is exactly a small $\catc$-enriched category with object set $\fC$.
\end{remark}

\subsection{Maps from Free Properads}

We will be dealing with free properads and maps between them very often.  Here we record what constitutes a map out of a free properad.

\begin{lemma}
\label{lem:mapfromfree}
Suppose $\sP$ is an $\SC$-colored object, and $\sQ$ is a $\fD$-colored properad.  Then a map of properads
\[
\nicexy{F\sP \ar[r]^-{f} & \sQ}
\]
is equivalent to a pair of functions:\index{free properad!map from}
\begin{enumerate}
\item
A function $\nicexy{\fC \ar[r]^-{f_0} & \fD}$ between color sets.
\item
A function $\nicexy{\sP\dc \ar[r]^-{f_1} & \sQ\binom{f_0\ud}{f_0\uc}}$ for each pair of $\fC$-profiles $\dc$.
\end{enumerate}
\end{lemma}

\begin{proof}
Since $F\sP$ is the free properad generated by $\sP$, a properad map $f \colon F\sP \to \sQ$ is equivalent to a map $\sP \to \sQ$ of colored objects (Definition \ref{def:coloredobject}), which consists of a pair of functions as stated.
\end{proof}


\chapter{Symmetric Monoidal Closed Structure on Properads}
\label{ch:tensor}

\abstract*{We equip the category of properads with a symmetric monoidal closed structure.   For topological operads, a symmetric monoidal product was already defined by Boardman and Vogt \cite{bv}.  One main result of this chapter gives a simple description of the tensor product of two free properads in terms of the two generating sets.  In particular, when the free properads are finitely generated, their tensor product is finitely presented.  This is not immediately obvious from the definition because free properads are often infinite sets. }

In this chapter, we construct a symmetric monoidal structure on the category of properads.  For special properads, there is also an internal hom, giving the category of special properads a symmetric monoidal closed structure.  In section \ref{sec:smgrsets} we will show that this symmetric monoidal structure induces a symmetric monoidal closed structure on the category of graphical sets for connected wheel-free graphs.  Our symmetric monoidal product of properads extends the one on operads defined by Boardman and Vogt \cite{bv} in the topological setting.

The symmetric monoidal product of properads is constructed in section \ref{sec:properadtensor}.  One slight complication of the properadic context comes from the fact that a vertex in a connected wheel-free graph may have multiple inputs/outputs or zero input/output.  Therefore, we need to be careful when we define the so-called distributivity relation, which tells us how two elements from the two component properads commute in the tensor product.

The symmetric monoidal product of two \emph{free} properads with special generating sets is given a much simpler description in sections \ref{sec:productfreeproperad} and \ref{sec:prdistsimplified}.  The main observation here is Theorem \ref{distsimplified}.  It says that the tensor product of two free properads can be directly generated from the original special generating sets with a much shorter and simplified lists of relations, called generating distributivity.  This is important because, as we will see later, a graphical properad, which is by definition freely generated, is an infinite set \emph{precisely} when the graph is not simply connected.  So for these infinite graphical properads generated by special connected wheel-free graphs (Definition \ref{def:wheelfreegraphs}), Theorem \ref{distsimplified} will tell us that the tensor product is finitely presented, which is not obvious from the definition of the tensor product.

The internal hom on special properads is constructed in section \ref{sec:properadhom}.  The generating distributivity mentioned above is very closely related to the naturality condition on internal hom.  In fact, as we will see in Lemma \ref{extendnaturality}, the way generating distributivity generates all distributivity is the same as the way natural transformations extend along connected wheel-free graphs.  This is, of course, no accident.  Indeed, in a symmetric monoidal closed category there is a natural isomorphism with the tensor product on one side and the internal hom on the other side \eqref{gpropclosed}.

Throughout this chapter, we work over the symmetric monoidal closed category $\Set$ of sets, so $\properad$ is the category of all colored properads in sets.  As in the previous chapter, with suitable restrictions everything in this chapter has obvious analogs for properads with non-empty inputs or non-empty outputs.

\section{Symmetric Monoidal Product}
\label{sec:properadtensor}

In this section, we define a symmetric monoidal product on properads that extends the Boardman-Vogt tensor product of operads \cite{bv}.  The symmetric monoidal product of two properads will be defined as a quotient of some free properad generated by the following colored object.

\subsection{Smash Product}

\begin{definition}
\label{def:smashproduct}
Suppose $\sP$ is an $\SC$-colored object and $\sQ$ is an $\SD$-colored object.  Define their \textbf{smash product} \index{smash product} $\sP \wedge \sQ$\label{note:smash} as the $\SCD$-colored object with two types of elements:
\begin{enumerate}
\item
For each color $c \in \fC$ and element $q \in \sQ\udud$, there is an element
\[
c \otimes q \in (\sP \wedge \sQ)\binom{c \times \ud^2}{c \times \ud^1},
\]
where
\[
c \times \ud^1 = \left((c,d^1_1),\ldots,(c,d^1_q)\right)
\]
if $\ud^1 = (d^1_1,\ldots,d^1_q)$.
\item
For each element $p \in \sP\ucuc$ and color $d \in \fD$, there is an element
\[
p \otimes d \in (\sP \wedge \sQ)\binom{\uc^2 \times d}{\uc^1 \times d},
\]
where
\[
\uc^1 \times d = \left((c^1_1,d),\ldots,(c^1_p,d)\right)
\]
if $\uc^1 = (c^1_1,\ldots,c^1_p)$.
\end{enumerate}
\end{definition}

\begin{remark}
As in Remark \ref{rk:coloredob}, we visualize the elements in the smash product $\sP \wedge \sQ$ as the following decorated corollas.
\begin{center}
\begin{tikzpicture}
\matrix[row sep=1cm,column sep=3cm] {
\node [fatplain,label=above:$c \times \ud^2$,label=below:$c \times \ud^1$] (p1) {$c \otimes q$}; & 
\node [fatplain,label=above:$\uc^2 \times d$,label=below:$\uc^1 \times d$] (p2) {$p \otimes d$}; \\
};
\foreach \x in {1,2}
{
\draw [outputleg] (p\x) to +(-.8cm,.6cm);
\draw [outputleg] (p\x) to +(.8cm,.6cm);
\draw [inputleg] (p\x) to +(-.8cm,-.6cm);
\draw [inputleg] (p\x) to +(.8cm,-.6cm);
}
\end{tikzpicture}
\end{center}
Also, we may write the element set of the smash product as
\[
\sP \wedge \sQ = \left[\sP \times \fD\right] \coprod \left[\fC \times \sQ\right],
\]
with the understanding that the profiles of an element are those of an element in $\sP$ (resp., $\sQ$) paired with a color in $\fD$ (resp., $\fC$).
\end{remark}

\begin{remark}
Smash product is associative in the sense that, if $\sR$ is an $\SE$-colored object, then there is a canonical isomorphism
\[
(\sP \wedge \sQ) \wedge \sR \cong \sP \wedge (\sQ \wedge \sR)
\]
of $\sS(\fC \times \fD \times \fE)$-colored objects.  The elements of this iterated smash product are of the forms
\[
p \otimes d \otimes e, \quad c \otimes q \otimes e, \quad c \otimes d \otimes r
\]
with $e \in \fE$ and $r \in \sR$.
\end{remark}

\subsection{Tensor Product}

\begin{definition}
\label{def:gpropmonoidalproduct}
Suppose $\sP$ is a $\fC$-colored properad and $\sQ$ is a $\fD$-colored properad.  Define the quotient $\fC \times \fD$-colored properad\index{tensor product!of properads} \index{properad!tensor product}
\begin{equation}
\label{gproptensor}
\sP \otimes \sQ \defn \frac{F(\sP \wedge \sQ)}{\text{3 types of relations}},
\end{equation}
where $F(\sP \wedge \sQ)$ is the free properad of the colored object $\sP\wedge\sQ$ (Definition \ref{def:fgupc}).  The relations are of the following three types.
\begin{enumerate}
\item
For each color $d \in \fD$, the functions
\[
\begin{split}
\fC \ni c & \longmapsto (c,d),\\
\sP \ni p & \longmapsto p \otimes d
\end{split}
\]
are required to define a map of properads $\sP \longrightarrow \sP \otimes \sQ$.
\item
For each color $c \in \fC$, the functions
\[
\begin{split}
\fD \ni d & \longmapsto (c,d),\\
\sQ \ni q & \longmapsto c \otimes q
\end{split}
\]
are required to define a map of properads $\sQ \longrightarrow \sP \otimes \sQ$.
\item
Suppose $p \in \sP\left(\ua;\ub\right)$ and $q \in \sQ\left(\uc;\ud\right)$ with $|\ua| = k$, $|\ub| = l$, $|\uc| = m$, $|\ud| = n$, $(k,l) \not= (0,0)$, and $(m,n) \not= (0,0)$.  The relation is then the equality
\begin{equation}
\label{distributivity}
\begin{split}
& \{p \otimes d_j\}_{j=1}^n \times \{a_i \otimes q\}_{i=1}^k \\
&=
\sigma^n_l\left[
 \{b_j \otimes q\}_{j=1}^l \times \{p \otimes c_i\}_{i=1}^m
\right]\sigma^m_k
\end{split}
\end{equation}
in $\sP \otimes \sQ$, called \index{distributivity} \textbf{distributivity}, which we will explain in detail below.
\end{enumerate}
\end{definition}

\subsection{Distributivity}
\label{subsec:distributivity}

Explicitly, the distributivity relation \eqref{distributivity} can be graphically represented as follows.  The left side of distributivity is the following decorated graph in $F(\sP \wedge \sQ)$:
\begin{center}
\begin{tikzpicture}
\matrix[row sep=2cm,column sep=2cm] {
\node [fatplain,label=above:$...$,label=below:$...$] (pd1) {$p \otimes d_1$}; &
\node [empty] (pd) {$\cdots$}; &
\node [fatplain,label=above:$...$,label=below:$...$] (pdn) {$p \otimes d_n$}; \\
\node [fatplain,label=above:$...$,label=below:$\cdots$] (a1q) {$a_1 \otimes q$}; & \node [empty] (aq) {$\cdots$}; &
\node [fatplain,label=above:$...$,label=below:$...$] (akq) {$a_k \otimes q$};
\\
};
\draw [arrow,bend left=45] (a1q) to node{$(a_1,d_1)$} (pd1);
\draw [arrow] (a1q) to node[swap,near end]{$(a_1,d_n)$} (pdn);
\draw [arrow] (akq) to node[near end]{$(a_k,d_1)$} (pd1);
\draw [arrow,bend right=45] (akq) to node[swap]{$(a_k,d_n)$} (pdn);
\draw [outputleg] (pd1) to node[above left]{$(b_1,d_1)$} +(-.8cm,.5cm);
\draw [outputleg] (pd1) to node[above right]{$(b_l,d_1)$} +(.8cm,.5cm);
\draw [outputleg] (pdn) to node[above left]{$(b_1,d_n)$} +(-.8cm,.5cm);
\draw [outputleg] (pdn) to node[above right]{$(b_l,d_n)$} +(.8cm,.5cm);
\draw [inputleg] (a1q) to node[below left]{$(a_1,c_1)$} +(-.8cm,-.5cm);
\draw [inputleg] (a1q) to node[below right]{$(a_1,c_m)$} +(.8cm,-.5cm);
\draw [inputleg] (akq) to node[below left]{$(a_k,c_1)$} +(-.8cm,-.5cm);
\draw [inputleg] (akq) to node[below right]{$(a_k,c_m)$} +(.8cm,-.5cm);
\end{tikzpicture}
\end{center}
In the top (resp., bottom) row, the $i$th vertex from the left is decorated by $p\otimes d_i$ (resp., $a_i \otimes q$).  The $r$th output of the vertex $a_i \otimes q$ is connected to the $i$th input of $p \otimes d_r$.  This internal edge has color $(a_i,d_r)$.

The right side of the distributivity relation, \emph{before} applying the input and output relabeling permutations $\sigma^m_k$ and $\sigma^n_l$, is the following decorated graph in $F(\sP \wedge \sQ)$:
\begin{center}
\begin{tikzpicture}
\matrix[row sep=2cm,column sep=2cm] {
\node [fatplain,label=above:$...$,label=below:$...$] (b1q) {$b_1 \otimes q$}; &
\node [empty] (bq) {$\cdots$}; &
\node [fatplain,label=above:$...$,label=below:$...$] (blq) {$b_l \otimes q$}; \\
\node [fatplain,label=above:$...$,label=below:$\cdots$] (pc1) {$p \otimes c_1$}; & \node [empty] (pc) {$\cdots$}; &
\node [fatplain,label=above:$...$,label=below:$...$] (pcm) {$p \otimes c_m$};
\\
};
\draw [arrow,bend left=45] (pc1) to node{$(b_1,c_1)$} (b1q);
\draw [arrow] (pc1) to node[swap,near end]{$(b_l,c_1)$} (blq);
\draw [arrow] (pcm) to node[near end]{$(b_1,c_m)$} (b1q);
\draw [arrow,bend right=45] (pcm) to node[swap]{$(b_l,c_m)$} (blq);
\draw [outputleg] (b1q) to node[above left]{$(b_1,d_1)$} +(-.8cm,.5cm);
\draw [outputleg] (b1q) to node[above right]{$(b_1,d_n)$} +(.8cm,.5cm);
\draw [outputleg] (blq) to node[above left]{$(b_l,d_1)$} +(-.8cm,.5cm);
\draw [outputleg] (blq) to node[above right]{$(b_l,d_n)$} +(.8cm,.5cm);
\draw [inputleg] (pc1) to node[below left]{$(a_1,c_1)$} +(-.8cm,-.5cm);
\draw [inputleg] (pc1) to node[below right]{$(a_k,c_1)$} +(.8cm,-.5cm);
\draw [inputleg] (pcm) to node[below left]{$(a_1,c_m)$} +(-.8cm,-.5cm);
\draw [inputleg] (pcm) to node[below right]{$(a_k,c_m)$} +(.8cm,-.5cm);
\end{tikzpicture}
\end{center}
In the top (resp., bottom) row, the $j$th vertex from the left is decorated by $b_j \otimes q$ (resp., $p \otimes c_j$).  The $s$th output of the vertex $p \otimes c_j$ is connected to the $j$th input of $b_s \otimes q$.  This internal edge has color $(b_s,c_j)$.

The above two $(\sP \wedge \sQ)$-decorated graphs have the same input and output color sets, but the orders are \emph{not} the same.  For example, in the first decorated graph above the output profile is
\[
\left(\ub \times d_1, \ldots , \ub \times d_n\right),
\]
while in the second decorated graph above the output profile is
\[
\left(b_1 \times \ud, \ldots , b_l \times \ud\right).
\]
The output relabeling $\sigma^n_l$ is the permutation that yields the first output profile when it is applied to the second output profile from the left, i.e.,
\[
\sigma^n_l \left(b_1 \times \ud, \ldots , b_l \times \ud\right) 
= \left(\ub \times d_1, \ldots , \ub \times d_n\right).
\] 
Likewise, the input relabeling $\sigma^m_l$ is the permutation that yields the input profile of the first decorated graph when it is applied to the input profile of the second decorated graph from the right, i.e.,
\[
\left(a_1 \times \uc, \ldots , a_k \times \uc\right)
= \left(\ua \times c_1, \ldots , \ua \times c_m\right) \sigma^m_l.
\]

If some of the parameters $k,l,m,n$ are $0$, then the distributivity relation requires special interpretation to keep the decorated graphs connected.  Explicitly, the additional restrictions for distributivity are as follows:
\begin{enumerate}
\item
If $|\ua| = k = 0$, then it is required that $|\ud| = n = 1$.  Conversely, if $n = 0$, then it is required that $k=1$.
\item
If $|\ub| = l = 0$, then it is required that $|\uc| = m = 1$.  Conversely, if $m=0$, then $l=1$.
\end{enumerate}
Graphically, these exceptional distributivity relations can be represented as follows.
\begin{enumerate}
\item
Suppose $|\ua| = k = 0$, so $ l > 0$ and $\ud = (d)$.
\begin{enumerate}
\item
Suppose $|\uc| = m > 0$.  Then distributivity identifies the decorated corolla
\bigskip
\begin{center}
\begin{tikzpicture}
\node [fatplain,label=above:$...$] (pd) {$p \otimes d$};
\draw [outputleg] (pd) to node[above left]{$(b_1,d)$} +(-.8cm,.5cm);
\draw [outputleg] (pd) to node[above right]{$(b_l,d)$} +(.8cm,.5cm);
\end{tikzpicture}
\end{center}
\bigskip
with the decorated graph:
\bigskip
\begin{center}
\begin{tikzpicture}
\matrix[row sep=2cm,column sep=2cm] {
\node [fatplain,label=below:$...$] (b1q) {$b_1 \otimes q$}; &
\node [empty] (bq) {$\cdots$}; &
\node [fatplain,label=below:$...$] (blq) {$b_l \otimes q$}; \\
\node [fatplain,label=above:$...$] (pc1) {$p \otimes c_1$}; & \node [empty] (pc) {$\cdots$}; &
\node [fatplain,label=above:$...$] (pcm) {$p \otimes c_m$};
\\
};
\draw [arrow,bend left=45] (pc1) to node{$(b_1,c_1)$} (b1q);
\draw [arrow] (pc1) to node[swap,near end]{$(b_l,c_1)$} (blq);
\draw [arrow] (pcm) to node[near end]{$(b_1,c_m)$} (b1q);
\draw [arrow,bend right=45] (pcm) to node[swap]{$(b_l,c_m)$} (blq);
\draw [outputleg] (b1q) to node[above=.1cm]{$(b_1,d)$} +(0,.7cm);
\draw [outputleg] (blq) to node[above=.1cm]{$(b_l,d)$} +(0,.7cm);
\end{tikzpicture}
\end{center}
\bigskip
The decorated corolla has no inputs.  In this case, input and output relabeling are not needed.
\item
Suppose further that $|\uc| = m = 0$, so $\ub = (b)$. Then distributivity identifies the following two decorated corollas:
\bigskip
\begin{center}
\begin{tikzpicture}
\matrix[row sep=2cm,column sep=1cm] {
\node [fatplain] (pd) {$p \otimes d$}; &
\node [empty] {$\sim$}; &
\node [fatplain] (bq) {$b \otimes q$}; \\
};
\draw [outputleg] (pd) to node[above]{$(b,d)$} +(0,.7cm);
\draw [outputleg] (bq) to node[above]{$(b,d)$} +(0,.7cm);
\end{tikzpicture}
\end{center}
\bigskip
Each decorated corolla has no inputs.
\end{enumerate}
\item
Suppose $|\ud| = n = 0$, so $m > 0$ and $\ua = (a)$.
\begin{enumerate}
\item
Suppose $|\ub| = l > 0$.  Then distributivity identifies the decorated corolla
\bigskip
\begin{center}
\begin{tikzpicture}
\node [fatplain,label=below:$...$] (aq) {$a \otimes q$};
\draw [inputleg] (aq) to node[below left]{$(a,c_1)$} +(-.8cm,-.5cm);
\draw [inputleg] (aq) to node[below right]{$(a,c_m)$} +(.8cm,-.5cm);
\end{tikzpicture}
\end{center}
\bigskip
with the decorated graph:
\bigskip
\begin{center}
\begin{tikzpicture}
\matrix[row sep=2cm,column sep=2cm] {
\node [fatplain,label=below:$...$] (b1q) {$b_1 \otimes q$}; &
\node [empty] (bq) {$\cdots$}; &
\node [fatplain,label=below:$...$] (blq) {$b_l \otimes q$}; \\
\node [fatplain,label=above:$...$] (pc1) {$p \otimes c_1$}; & \node [empty] (pc) {$\cdots$}; &
\node [fatplain,label=above:$...$] (pcm) {$p \otimes c_m$};
\\
};
\draw [arrow,bend left=45] (pc1) to node{$(b_1,c_1)$} (b1q);
\draw [arrow] (pc1) to node[swap,near end]{$(b_l,c_1)$} (blq);
\draw [arrow] (pcm) to node[near end]{$(b_1,c_m)$} (b1q);
\draw [arrow,bend right=45] (pcm) to node[swap]{$(b_l,c_m)$} (blq);
\draw [inputleg] (pc1) to node[below]{$(a,c_1)$} +(0,-.7cm);
\draw [inputleg] (pcm) to node[below]{$(a,c_m)$} +(0,-.7cm);
\end{tikzpicture}
\end{center}
\bigskip
The decorated corolla has no outputs.  Input and output relabeling are not needed.
\item
Suppose further that $|\ub| = l = 0$, so $\uc = (c)$.  Then distributivity identifies the following two decorated corollas:
\bigskip
\begin{center}
\begin{tikzpicture}
\matrix[row sep=2cm,column sep=1cm] {
\node [fatplain] (pc) {$p \otimes c$}; &
\node [empty] {$\sim$}; &
\node [fatplain] (aq) {$a \otimes q$}; \\
};
\draw [inputleg] (pc) to node[below]{$(a,c)$} +(0,-.7cm);
\draw [inputleg] (aq) to node[below]{$(a,c)$} +(0,-.7cm);
\end{tikzpicture}
\end{center}
\bigskip
\end{enumerate}
\end{enumerate}

\begin{remark}
The tensor product $\sP \otimes \sQ$ above can be similarly defined if $\sP$ and $\sQ$ are both in $\properadi$, $\properado$, or $\properads$ (Definition \ref{rk:properadio}).  Simply insert the phrase \emph{with non-empty inputs}, \emph{with non-empty outputs}, or \emph{special} in suitable places throughout the definition.  The functor $F$ in $F(\sP \wedge \sQ)$ means $F_{\gupci}$, $F_{\gupco}$, or $F_{\gupcs}$.  In particular, in defining the tensor product in $\properads$, distributivity requires no extra interpretation because the parameters $k,l,m,n$ are all positive.
\end{remark}

\subsection{Symmetric Monoidal Structure}

\begin{theorem}
\label{thm:propgmonoidal}
With respect to $\otimes$, the categories $\properad$, $\properadi$, $\properado$, and $\properads$ are symmetric monoidal.
\end{theorem}

\begin{proof}
We verify the various axioms for a symmetric monoidal structure.
\begin{enumerate}
\item
The symmetry is induced by the isomorphism of smash products
\[
\sP \wedge \sQ \cong \sQ \wedge \sP
\]
that switches the two entries in both colors and elements.
\item
The unit element is the $1$-colored properad $\bone$ whose only element is the identity $1 \in \bone(*;*)$.
\item
The associativity of $\otimes$ follows from the fact that both $(\sP \otimes \sQ) \otimes \sR$ and $\sP \otimes (\sQ \otimes \sR)$ are isomorphic to the quotient $(\fC \times \fD \times \fE)$-colored properad
\[
\sO = \frac{F(\sP \wedge \sQ \wedge \sR)}{\text{relations}}.
\]
There are two types of relations in this quotient, similar to the ones in Definition \ref{def:gpropmonoidalproduct}.  The first type of relations ensures that there are properad maps from each of $\sP$, $\sQ$, and $\sR$ to $\sO$ if we fix two colors in two different color sets $\fC$, $\fD$, and $\fE$.  The second type of relations is the three-component version of distributivity, which equates the element
\[
\{p \otimes d \otimes e\} \times \{c \otimes q \otimes e\} \times \{c \otimes d \otimes r\},
\]
with suitable subscripts and superscripts in the colors, 
with the element
\[
\{c \otimes q \otimes e\} \times \{p \otimes d \otimes e\} \times \{c \otimes d \otimes r\},
\]
and so forth.
\end{enumerate}
\end{proof}

\begin{remark}
The unit element $\bone$ above is the free $1$-colored properad generated by the $1$-colored object with all empty components.
\end{remark}

\begin{remark}
For colored PROPs, a similar symmetric monoidal product is constructed in \cite{hr1}.
\end{remark}


\section{Symmetric Monoidal Product of Free Properads}
\label{sec:productfreeproperad}

In this and the next sections, we give a simple description of the tensor product of two \emph{free} properads generated by special colored objects (Definition \ref{def:coloredobject}).  Later we will use this description on graphical properads, which are freely generated by construction.  This in turn yields a description of the tensor product of two graphical sets.

\subsection{Motivation}

Let us provide some motivation for the theorems we are about to prove.  Suppose now that $\sP = F(\phat)$ and $\sQ = F(\qhat)$ are both free properads with colored generating sets $\phat$ and $\qhat$, respectively.  Then their tensor product is defined as the quotient
\begin{equation}
\label{ptensorq}
\sP \otimes \sQ = \frac{F\left(F(\phat) \wedge F(\qhat)\right)}{\text{3 types of relations}}.
\end{equation}
Note that this construction of the tensor product seems somewhat redundant for two reasons.
\begin{enumerate}
\item
There are two layers of the  free properad construction.  Indeed, $F(\phat)$ consists of $\phat$-decorated graphs, and similarly for the free properad $F(\qhat)$.  So the free properad $F\left(F(\phat) \wedge F(\qhat)\right)$ consists of connected wheel-free graphs in which each vertex is itself a decorated graph.  It is therefore desirable to have a more direct description of the tensor product in terms of the generating sets $\phat$ and $\qhat$ that has a single free construction.
\item
Moreover, even if $\phat$ has a finite set of colors and a finite set of elements, the free properad $F(\phat)$ is usually an infinite set.  This is the case, for example, for the graphical properad of a graph that is \emph{not} simply connected (Lemma \ref{lem:omegaginfinite}).  Therefore, each of the three relations imposed on the free properad $F\left(F(\phat) \wedge F(\qhat)\right)$ is usually an infinite list of relations.  We would like to simplify the relations, and hence the construction of the tensor product, using the generating sets.
\end{enumerate}

Our plan to understanding the tensor product of two free properads \eqref{ptensorq} then consists of two stages.
\begin{enumerate}
\item
The first part is a Poincar\'{e}-Birkhoff-Witt type theorem.  We consider the quotient of the free properad $F\left(F(\phat) \wedge F(\qhat)\right)$ by the first two relations.  We will provide a very direct description of this quotient using the generating sets $\phat$ and $\qhat$.  Indeed, we will show that this quotient is canonically isomorphic to the free properad $F(\phat \wedge \qhat)$.  
\item
Once this is done, we consider the last relation, namely, distributivity.  We will show that imposing the distributivity relations on the quotient of the previous step is equivalent to imposing a very small subset of distributivity relations, called generating distributivity, on $F(\phat \wedge \qhat)$.  Each generating distributivity corresponds to an element in the Cartesian product $\phat \times \qhat$.  In particular, not only are we greatly reducing the number of relations, but the remaining relations are all relatively simple.
\end{enumerate}

\subsection{The First Two Relations}

We first give a simpler description of the part of the tensor product \eqref{ptensorq} about the first two relations, corresponding to the first stage discussed above.  First we need a preliminary observation.

\begin{lemma}
\label{pqgenerates}
Suppose $\sP = F(\phat)$ and $\sQ = F(\qhat)$ are $\fC$-colored and $\fD$-colored free properads with colored generating sets $\phat$ and $\qhat$, respectively.  Then the natural map 
\[
\nicearrow
\xymatrix{
\phat \wedge \qhat \ar[r]^-{\iota} & F(\phat) \wedge F(\qhat)
}\]
of colored sets induces a surjection
\begin{equation}
\label{sthequotient}
\nicearrow
\xymatrix{
F\left(\phat \wedge \qhat\right) \ar[r]^-{F(\iota)} & 
F\left(F(\phat) \wedge F(\qhat)\right) \ar[r]^-{\pi} &
\frac{F\left(F(\phat) \wedge F(\qhat)\right)}{\text{relations (1) and (2)}} \defn \sW
}
\end{equation}
of properads, where $\pi$ is the quotient map. 
\end{lemma}

To be precise, relations (1) and (2) in \eqref{sthequotient} mean that
\begin{enumerate}
\item
$\nicearrow\xymatrix{F(\phat) \ar[r]^-{(-)\otimes d} & \sW}$ is a map of properads for each $d \in \fD$, and
\item
$\nicearrow\xymatrix{F(\qhat) \ar[r]^-{c \otimes (-)} & \sW}$ is a map of properads for each $c \in \fC$.
\end{enumerate}

\begin{proof}
This is mostly about unraveling the definitions of the objects involved.  A general element in $F(\phat)$ is a $\phat$-decorated graph
\[
\left\{p_v\right\}_{\ving} \in \phat[G] = \prod_{v \in \vertex(G)} \phat\profilev
\]
for some $\fC$-colored connected wheel-free graph $G$.  An entry $p_v$ is called a \textbf{vertex decoration}.  Similar descriptions hold for the free properads $F(\qhat)$, $F(\phat \wedge \qhat)$, and $F\left(F(\phat) \wedge F(\qhat)\right)$. In particular, an element
\[
R = \{R_v\}_{\ving} \in F\left(F(\phat) \wedge F(\qhat)\right)
\]
is a $(\fC \times \fD)$-colored decorated connected wheel-free graph with each vertex decoration
\[
R_v \in \left[F(\phat) \times \fD\right] \coprod \left[\fC \times F(\qhat)\right].
\]
In other words, $R_v$ is itself either
\begin{itemize}
\item
a $\fC$-colored $\phat$-decorated graph paired with a color in $\fD$, i.e.,
\[
R_v = \{p^v_u\}_{u \in H_v} \otimes d_v,
\]
or 
\item
a color in $\fC$ paired with a $\fD$-colored $\qhat$-decorated graph, i.e.,
\[
R_v = c_v \otimes \{q^v_u\}_{u \in H_v}.
\]
\end{itemize}
Here we abbreviated $u \in \vertex(H_v)$ to $u \in H_v$. Each generator $p \in \phat$ (resp., $q \in \qhat$) is regarded as a decorated corolla in $F(\phat)$ (resp., $F(\qhat)$), which defines the injection $\iota$.  The map $F(\iota)$ sends a $(\phat \wedge \qhat)$-decorated graph to the same decorated graph, but with each vertex decoration regarded as a decorated corolla.

Recall that the free properad structure comes from graph substitution.  Therefore, the unity property of corollas with respect to graph substitution and relation (1) imply that for each $\phat$-decorated graph $\{p_u\}_{u \in H} \in F(\phat)$ and color $d \in \fD$, there is an equality
\[
\{p_u\}_{u\in H} \otimes d = \{p_u \otimes d\}_{u \in H}
\]
in $\sW$, where on the right-hand side $p_u \otimes d$ is a decorated corolla.  Likewise, we have
\[
c \otimes \{q_u\}_{u\in H} = \{c \otimes q_u\}_{u\in H}
\]
in $\sW$ for each $\qhat$-decorated graph $\{q_u\}_{u \in H} \in F(\qhat)$ and color $c \in \fC$.

For an element $R = \{R_v\}$ as above, apply these equalities to its image in $\sW$.  Using the associativity of the free properad structure map (i.e., associativity of graph substitution), it follows that the image of $R$ in $\sW$ takes the form
\[
R = \left\{r_w\right\}_{w \in G(\{H_v\})},
\]
in which each $r_w$ is a decorated corolla of the form
\[
p^v_u \otimes d_v \orspace c_v \otimes q^v_u
\]
with $u$ a vertex in some $H_v$.  Here $G(\{H_v\})$ is the graph substitution of the $H_v$ into $G$.  Now the decorated graph $\{r_w\}$ is in the image of $F(\iota)$, which shows that $\pi \circ F(\iota)$ is surjective.
\end{proof}

\begin{remark}
Lemma \ref{pqgenerates} says that the image of $\phat \wedge \qhat$ generates the quotient $\sW$.  However, $\phat \wedge \qhat$ does \emph{not} generate the free properad $F\left(F(\phat) \wedge F(\qhat)\right)$, which is much bigger than $\sW$.  For example, for a $\phat$-decorated graph $\{p\}_v \in F(\phat)$ with at least two vertices and a color $d \in \fD$, the element
\[
\{p\}_v \otimes d \in F\left(F(\phat) \wedge F(\qhat)\right)
\]
does not belong to the properad generated by the image of $\phat \wedge \qhat$.  In other words, taking the quotient $\sW$ is necessary to make the map from $F(\phat \wedge \qhat)$ surjective.
\end{remark}

\begin{theorem}
\label{relationsonetwo}
Suppose $\sP = F(\phat)$ and $\sQ = F(\qhat)$ are free properads with colored generating sets $\phat$ and $\qhat$, respectively.  Then the properad map
\begin{equation}
\label{simpleonetwo}
\nicearrow
\xymatrix@C+10pt{
F\left(\phat \wedge \qhat\right) \ar[r]^-{\pi \circ F(\iota)} & 
\frac{F\left(F(\phat) \wedge F(\qhat)\right)}{\text{relations (1) and (2)}} = \sW
}
\end{equation}
in Lemma \ref{pqgenerates} is an isomorphism.
\end{theorem}

\begin{proof}
We will use the notations in the proof of Lemma \ref{pqgenerates}.
 
It suffices to observe that $\sW$ has the same universal property that characterizes the free properad $F\left(\phat \wedge \qhat\right)$.  So suppose $\sR$ is an $\fE$-colored properad, and $f \colon \phat \wedge \qhat \to \sR$ is a map of colored objects.  We must show that $f$ extends uniquely to a map $\fbar \colon \sW \to \sR$ of properads.  First note that, since the composition $\pi \circ F(\iota)$ is surjective by Lemma \ref{pqgenerates}, if such an extension $\fbar$ exists, then it must be unique.

To see that an extension $\fbar$ exists, first note that $f$ has a properad map extension to the intermediate free properad $F\left(F(\phat) \wedge F(\qhat)\right)$.  Indeed, since this properad is freely generated, such a properad map is equivalent to a map
\[
F(\phat) \wedge F(\qhat) \to \sR
\]
of colored objects.  Each element in $F(\phat) \wedge F(\qhat)$ is an element in either
\begin{enumerate}
\item
$\phat[G] \times d$ for some $\fC$-colored graph $G \in \gupc$ and color $d \in \fD$, or
\item
$c \times \qhat[H]$ for some $\fD$-colored graph $H \in \gupc$ and color $c \in \fC$.
\end{enumerate}
In the first case with $G \in \gupc\yxh$, we can send it to $\sR$ using the following composition:
\[
\nicearrow
\xymatrix@C+10pt{
\phat[G] \times d \ar@{=}[d] \ar[r]^-{\fbar} & \sR\binom{f(\uy \times d)}{f(\ux \times d)}\\
\left[\prod_{v \in \vertex(G)} \phat\profilev\right] \times d \ar[d]_{\cong} & \\
\prod_{v\in \vertex(G)} \left[\phat\profilev \times d\right] \ar[r]^-{\prod f} &  \sR[f(G \otimes d)]. \ar[uu]_{\gamma}
}\]
Here $G \otimes d$ is the $(\fC\times\fD)$-colored graph obtained from $G$ by pairing each flag with $d$.  Likewise, $f(G \otimes d)$ is the $\fE$-colored graph obtained from $G \otimes d$ by applying $f$ to each flag color.  The map $\gamma$ is a properad structure map of $\sR$, using the unbiased description of a properad (Remark \ref{rk:falgstructure}).  There is a similar composition for the second case.

This composition $\fbar$ extends $f$ because:
\begin{itemize}
\item
an element in $\phat$ is regarded as a decorated corolla in $F(\phat)$, and
\item
the properad structure map $\gamma$ corresponding to a corolla is the identity map.
\end{itemize}
Next observe that this map $\fbar$ factors through the quotient $\sW$ by construction.  This is true because:
\begin{itemize}
\item
the free properad structure is given by graph substitution,
\item
$\gamma$ is associative with respect to graph substitution, and
\item
$\fbar$ of a decorated graph is defined at each vertex.
\end{itemize}
Therefore, an extension $\fbar$ of $f$ to $\sW$ exists.  The map $\fbar \colon \sW \to \sR$ is a map of properads (i.e., of $F_{\gupc}$-algebras) because the properad structure on $\sW$ is given by graph substitution, and the properad structure map $\gamma$ on $\sR$ is associative with respect to graph substitution.  As noted earlier, this suffices to finish the proof.
\end{proof}

\subsection{Generating Distributivity}

Next we give a simpler description of the tensor product of two free properads regarding distributivity, for which we need the following special case.

\begin{definition}
\label{def:gendist}
Suppose $\sP = F(\phat)$ and $\sQ = F(\qhat)$ are free properads with colored generating sets $\phat$ and $\qhat$, respectively.  By \textbf{generating distributivity} \index{generating distributivity} we mean the distributivity relation on $F(\phat \wedge \qhat)$, as in section \ref{subsec:distributivity}, with $p \in \phat$ and $q \in \qhat$.
\end{definition}

\begin{remark}
\begin{enumerate}
\item
In other words, generating distributivity is the distributivity relation for the generators only.  In particular, if $\phat$ and $\qhat$ both have finite sets of elements, then there are only finitely many generating distributivity relations.
\item
By Theorem \ref{relationsonetwo}, generating distributivity may be interpreted in either the free properad $F(\phat \wedge \qhat)$ or the isomorphic properad $\sW$ \eqref{sthequotient}.
\item
In Theorem \ref{relationsonetwo} we already described the first two relations involved in constructing the properadic tensor product $F(\phat) \otimes F(\qhat)$.  To describe the tensor product itself, it remains to describe the distributivity relation, which is the next observation.
\end{enumerate}
\end{remark}

Recall that a colored object $X$ is \emph{special} if $X\dc = \varnothing$ whenever either $\uc = \varnothing$ or $\ud=\varnothing$.

\begin{theorem}
\label{distsimplified}
Suppose $\sP = F(\phat)$ and $\sQ = F(\qhat)$ are free properads with special colored generating sets $\phat$ and $\qhat$, respectively.  Then there is an isomorphism
\[
\sP \otimes \sQ \cong \frac{F(\phat \wedge \qhat)}{\text{generating distributivity}}
\]
of properads.\index{tensor product!of free properads}
\end{theorem}

\begin{remark}
The point of this theorem is that both the colored generating set and the relations are much simpler than in the definition of the tensor product \eqref{ptensorq}.  More explicitly, this description uses the colored generating set $\phat \wedge \qhat$, as opposed to $F(\phat) \wedge F(\qhat)$ as in the original definition.  Moreover, instead of three types of relations, there is only one here.  Furthermore, here we only need a very special kind of distributivity relations, namely, those involving the generators (i.e., elements in $\phat$ and $\qhat$), as opposed to elements in $F(\phat)$ and $F(\qhat)$.
\end{remark}

\subsection{Consequences of Theorem \ref{distsimplified}}

The proof of Theorem \ref{distsimplified} will be given in the next section.  Before that let us first discuss a few consequences.  The first consequence of Theorem \ref{distsimplified} is that the properadic tensor product of two finitely generated free properads is finitely presented.

\begin{definition}
Suppose $\sP$ is a properad.
\begin{enumerate}
\item
We say that $\sP$ is \textbf{finitely generated} \index{finitely generated} if there is an isomorphism
\[
\sP \cong \frac{F(\phat)}{\text{some relations}}
\]
for some finite colored object $\phat$ on finitely many colors and for some relations.  We call such an isomorphism a \textbf{finite generation presentation} \index{finite generation presentaion} of $\sP$.
\item
We say that $\sP$ is \textbf{finitely presented} \index{finitely presented} if it admits a finite generation presentation in which there are only finitely many relations in the quotient.  In this case, we call the quotient a \textbf{finite presentation} \index{finite presentation} of $\sP$.
\end{enumerate}
\end{definition}

\begin{corollary}
\label{tensorfpresented}
Suppose $\phat$ and $\qhat$ are finite special colored objects on finitely many colors.  Then the properad $F(\phat) \otimes F(\qhat)$ is finitely presented. 
\end{corollary}

\begin{proof}
Theorem \ref{distsimplified} exhibits the generating set
\[
\phat \wedge \qhat = \left[\phat \times \fD\right] \coprod \left[\fC \times \qhat\right],
\]
in which $\fC$ and $\fD$ are the sets of colors of $\phat$ and $\qhat$, respectively.  The hypotheses imply that this coproduct is a finite set.  Moreover, there are only finitely many generating distributivity relations because both $\phat$ and $\qhat$ are finite.
\end{proof}

Another consequence of the theorem is that it yields a simple description of maps out of $F(\phat) \otimes F(\qhat)$.

\begin{corollary}
\label{mapfromtensor}
Suppose $\sP = F(\phat)$ and $\sQ = F(\qhat)$ are free properads with special colored generating sets $\phat$ and $\qhat$, respectively.  Suppose $\sR$ is a properad.  Then a map $\sP \otimes \sQ \to \sR$ of properads is equivalent to a map 
\[
\nicearrow
\xymatrix{
\phat \wedge \qhat \ar[r]^-{\theta} & \sR
}\]
of colored objects such that the following condition holds:  For any $p \in \phat\bah$ and $q \in \qhat\dch$, the equality
\[
\begin{split}
&\gamma_T \left(\left\{\theta(p \otimes d_j)\right\}_{j=1}^n \times \left\{\theta(a_i \otimes q)\right\}_{i=1}^k\right) \\
&= 
\sigma^n_l\left[
\gamma_{T'} \left( \left\{\theta(b_j \otimes q)\right\}_{j=1}^l \times \left\{\theta(p \otimes c_i)\right\}_{i=1}^m\right) 
\right]\sigma^m_k
\end{split}
\]
holds in $\sR$, where $\gamma$ is the properad structure map of $\sR$, and the elements are as in the distributivity relation  \eqref{distributivity}.
\end{corollary}

\begin{proof}
By Theorem \ref{distsimplified} a properad map $\sP \otimes \sQ \to \sR$ is equivalent to a properad map $F(\phat \wedge \qhat) \to \sR$ such that, for each generating distributivity relation, the images of the two decorated graphs are equal in $\sR$.  A properad map out of the free properad $F(\phat \wedge \qhat)$ is equivalent to a colored object map out of $\phat \wedge \qhat$.  The condition about respecting the generating distributivity relations is the stated equality.
\end{proof}


\section{Proof of Theorem \ref{distsimplified}}
\label{sec:prdistsimplified}

This section contains the proof of Theorem \ref{distsimplified}.  One point to keep in mind is that, although it looks like the proof is quite long, there is really only one main point in the proof.  It consists of several pictures in the proof of Lemma \ref{gendisttypical}.

\subsection{Reduction of Proof}

By definition the tensor product $F(\phat) \otimes F(\qhat)$ is obtained from the quotient $\sW$ \eqref{sthequotient} by imposing the distributivity relations \eqref{distributivity}. Moreover, by Theorem \ref{relationsonetwo} there is a canonical isomorphism between the free properad $F(\phat \wedge \qhat)$ and $\sW$.  Thus, it suffices to show that in $\sW$ every distributivity relation (section \ref{subsec:distributivity}) decomposes into a finite string of \emph{generating} distributivity relations (Definition \ref{def:gendist}).

Recall that an element in $F(\phat)$ is a $\phat$-decorated connected wheel-free graph, i.e., an element in the Cartesian product
\[
\phat[G] = \prod_{v\in \vertex(G)} \phat\profilev
\]
for some $\fC$-colored $G \in \gupc$, where $\fC$ is the color set of $\phat$.  Likewise, an element in $F(\qhat)$ is an element in
\[
\qhat[H] = \prod_{u \in \vertex(H)} \qhat\profileu
\]
for some $\fD$-colored $H \in \gupc$, where $\fD$ is the color set of $\qhat$.

Pick a $\phat$-decorated graph
\[
p = \{p_v\}_{v\in \vertex(G)} \in \phat[G]
\]
with each $p_v \in \phat\profilev$ and a $\qhat$-decorated graph
\[
q = \{q_u\}_{u \in \vertex(H)} \in \qhat[H]
\]
with each $q_u \in \qhat\profileu$. Suppose $|\vertex(G)| = m$ and $|\vertex(H)| = n$.   We will show by induction on $m+n$ that the distributivity relation involving $p$ and $q$ decomposes into a finite string of generating distributivity relations. We will safely suppress the permutations in what follows.  This is possible because if the distributivity involving $p$ and $q$ decomposes into a finite string of generating distributivity, then the same is true for permuted versions of $p$ and $q$.

Note that if either one of these decorated graphs has $0$ vertex (i.e., it is an exceptional edge), then distributivity is trivially true because it is a consequence of the unity axiom in a properad.  Therefore, in what follows we may assume without loss of generality that $m,n > 0$.  Moreover, if $m=n=1$, then the distributivity relation is a generating distributivity relation.  This is true because the only connected wheel-free graphs with exactly one vertex are the permuted corollas.  Therefore, we may further assume that $m+n \geq 3$.  To improve readability, we will split the rest of the proof into several steps.

\subsection{Induction Step}

First we consider a special case.  Its proof will be reused several times below for more general cases.

\begin{lemma}
\label{gendisttypical}
Suppose $n  = |\vertex(H)| = 1$ and $m = |\vertex(G)| = 2$.  Then the distributivity relation involving $p = \{p_v\}$ and $q = \{q_u\}$ decomposes into a finite string of generating distributivity relations.
\end{lemma}

\begin{proof}
Since $n=1$, $H$ is a corolla, and there is only one $q_u$.  So we will simply write $q_u$ as $q$ and draw it as a decorated corolla:
\begin{center}
\begin{tikzpicture}
\node [plain,label=above:$...$,label=below:$...$] (q) {$q$};
\draw [inputleg] (q) to +(-.6cm,-.4cm);
\draw [inputleg] (q) to +(.6cm,-.4cm);
\draw [outputleg] (q) to +(-.6cm,.4cm);
\draw [outputleg] (q) to +(.6cm,.4cm);
\end{tikzpicture}
\end{center}
Since $G \in \gupc$ has $2$ vertices, it is a partially grafted corollas.  Write $p_1$ (resp., $p_2$) for the entry in $p$ corresponding to the lower (resp., upper) vertex in $G$.  We will draw $p$ as a decorated partially grafted corollas:
\begin{center}
\begin{tikzpicture}
\matrix[row sep=.5cm,column sep=1cm] {
\node [plain,label=above:$...$] (p2) {$p_2$}; \\
\node [empty] (dots) {$\cdots$}; \\ 
\node [plain,label=below:$...$] (p1) {$p_1$}; \\
};
\draw [arrow,bend left=40] (p1) to (p2);
\draw [arrow,bend right=40] (p1) to (p2);
\draw [inputleg] (p2) to +(-.6cm,-.4cm);
\draw [inputleg] (p2) to +(.6cm,-.4cm);
\draw [outputleg] (p2) to +(-.6cm,.4cm);
\draw [outputleg] (p2) to +(.6cm,.4cm);
\draw [outputleg] (p1) to +(-.6cm,.4cm);
\draw [outputleg] (p1) to +(.6cm,.4cm);
\draw [inputleg] (p1) to +(-.6cm,-.4cm);
\draw [inputleg] (p1) to +(.6cm,-.4cm);
\end{tikzpicture}
\end{center}
The assumption that $\phat$ and $\qhat$ are special ensures that, in these decorated graphs, each vertex has both incoming flags and outgoing flags.  This point will be important in the proof below.

For ease of drawing the decorated graphs below, we assume that $q$ has two inputs and two outputs.  Likewise, we assume that $p_2$ has two outputs and four inputs, two of which are ordinary edges connected to $p_1$, and that $p_1$ has two inputs and four outputs, two of which are ordinary edges connected to $p_2$.  In other words, remove the $\cdots$ in the two decorated graphs above.  Note that in this case $G$ has four inputs (resp., outputs), two of which are connected to each $p_i$. The general case will be obtained as a minor modification of this special case.

We now show that the distributivity  involving $p$ and $q$ decomposes into \emph{two} generating distributivity.  One side of distributivity is the decorated graph:
\begin{center}
\begin{tikzpicture}
\matrix[row sep=1cm,column sep=1.5cm] {
& \node [plain] (p21) {$p_2$}; &
& \node [plain] (p22) {$p_2$}; & \\
& \node [plain] (p11) {$p_1$}; &
& \node [plain] (p12) {$p_1$}; & \\
&&&& \\
\node [plain] (q1) {$q$}; 
& \node [plain] (q2) {$q$}; 
&& \node [plain] (q3) {$q$}; 
& \node [plain] (q4) {$q$};\\
};
\draw [arrow,bend left=20] (p11) to (p21);
\draw [arrow,bend right=20] (p11) to (p21);
\draw [arrow,bend left=20] (p12) to (p22);
\draw [arrow,bend right=20] (p12) to (p22);
\draw [arrow,bend left=20] (q1) to (p21);
\draw [arrow] (q1) to (p22);
\draw [arrow] (q2) to (p11);
\draw [arrow] (q2) to (p12);
\draw [arrow] (q3) to (p11);
\draw [arrow] (q3) to (p12);
\draw [arrow] (q4) to (p21);
\draw [arrow,bend right=20] (q4) to (p22);
\draw [outputleg] (p21) to +(-.6cm,.4cm);
\draw [outputleg] (p21) to +(.6cm,.4cm);
\draw [outputleg] (p22) to +(-.6cm,.4cm);
\draw [outputleg] (p22) to +(.6cm,.4cm);
\draw [outputleg] (p11) to +(-.6cm,.4cm);
\draw [outputleg] (p11) to +(.6cm,.4cm);
\draw [outputleg] (p12) to +(-.6cm,.4cm);
\draw [outputleg] (p12) to +(.6cm,.4cm);
\draw [inputleg] (q1) to +(-.6cm,-.4cm);
\draw [inputleg] (q1) to +(.6cm,-.4cm);
\draw [inputleg] (q2) to +(-.6cm,-.4cm);
\draw [inputleg] (q2) to +(.6cm,-.4cm);
\draw [inputleg] (q3) to +(-.6cm,-.4cm);
\draw [inputleg] (q3) to +(.6cm,-.4cm);
\draw [inputleg] (q4) to +(-.6cm,-.4cm);
\draw [inputleg] (q4) to +(.6cm,-.4cm);
\end{tikzpicture}
\end{center}
The top left (resp., right) copy of $p$ is paired with the first (resp., second) output color of $q$.  Likewise, along the bottom row from left to right, the $i$th copy of $q$ is paired with the $i$th input color of $p$.  We will use such conventions below to simplify the notations in such decorated graphs.  Furthermore, we will not draw permutations for input and output relabeling.  

In the above decorated graph, consider the two copies of $p_1$ and the two copies of $q$ directly connected to them, i.e., the middle two copies of $q$.  Using the generating distributivity involving $p_1$ and $q$, the above decorated graph is identified with the following decorated graph:
\begin{center}
\begin{tikzpicture}
\matrix[row sep=2cm,column sep=1.5cm] {
& \node [plain] (p21) {$p_2$}; &&&  \node [plain] (p22) {$p_2$}; &\\
\node [plain] (q1) {$q$};
& \node [plain] (q2) {$q$};
& \node [plain] (q3) {$q$};
& \node [plain] (q4) {$q$};
& \node [plain] (q5) {$q$};
& \node [plain] (q6) {$q$};\\
& \node [plain] (p11) {$p_1$}; &&&  \node [plain] (p12) {$p_1$}; &\\
};
\draw [arrow,bend left=20] (q1) to (p21);
\draw [arrow] (q1) to (p22);
\draw [arrow] (q3) to (p21);
\draw [arrow] (q3) to (p22);
\draw [arrow] (q4) to (p21);
\draw [arrow] (q4) to (p22);
\draw [arrow] (q6) to (p21);
\draw [arrow,bend right=20] (q6) to (p22);
\foreach \x in {1,2}
\foreach \y in {2,3,4,5}
{
 \draw [arrow] (p1\x) to (q\y);
}
\draw [outputleg] (p21) to +(-.6cm,.4cm);
\draw [outputleg] (p21) to +(.6cm,.4cm);
\draw [outputleg] (p22) to +(-.6cm,.4cm);
\draw [outputleg] (p22) to +(.6cm,.4cm);
\draw [outputleg] (q2) to +(-.5cm,.3cm);
\draw [outputleg] (q2) to +(.5cm,.3cm);
\draw [outputleg] (q5) to +(-.5cm,.3cm);
\draw [outputleg] (q5) to +(.5cm,.3cm);
\draw [inputleg] (q1) to +(-.6cm,-.4cm);
\draw [inputleg] (q1) to +(.6cm,-.4cm);
\draw [inputleg] (q6) to +(-.6cm,-.4cm);
\draw [inputleg] (q6) to +(.6cm,-.4cm);
\draw [inputleg] (p11) to +(-.6cm,-.4cm);
\draw [inputleg] (p11) to +(.6cm,-.4cm);
\draw [inputleg] (p12) to +(-.6cm,-.4cm);
\draw [inputleg] (p12) to +(.6cm,-.4cm);
\end{tikzpicture}
\end{center}
There is another generating distributivity that can be applied to this latest decorated graph.  To clarify how this is applied, we will redraw the above decorated graph by swapping the positions of the first two copies of $q$ and also the last two copies of $q$:
\begin{center}
\begin{tikzpicture}
\matrix[row sep=2cm,column sep=1.5cm] {
& \node [plain] (p21) {$p_2$}; &&&  \node [plain] (p22) {$p_2$}; &\\
\node [plain] (q1) {$q$};
& \node [plain] (q2) {$q$};
& \node [plain] (q3) {$q$};
& \node [plain] (q4) {$q$};
& \node [plain] (q5) {$q$};
& \node [plain] (q6) {$q$};\\
& \node [plain] (p11) {$p_1$}; &&&  \node [plain] (p12) {$p_1$}; &\\
};
\foreach \x in {2,3,4,5}
\foreach \y in {1,2}
{
 \draw [arrow] (q\x) to (p2\y);
}
\draw [arrow,bend left=25] (p11) to (q1);
\draw [arrow] (p11) to (q3);
\draw [arrow] (p11) to (q4);
\draw [arrow] (p11) to (q6);
\draw [arrow] (p12) to (q1);
\draw [arrow] (p12) to (q3);
\draw [arrow] (p12) to (q4);
\draw [arrow,bend right=25] (p12) to (q6);
\draw [outputleg] (p21) to +(-.6cm,.4cm);
\draw [outputleg] (p21) to +(.6cm,.4cm);
\draw [outputleg] (p22) to +(-.6cm,.4cm);
\draw [outputleg] (p22) to +(.6cm,.4cm);
\draw [outputleg] (q1) to +(-.5cm,.4cm);
\draw [outputleg] (q1) to +(.5cm,.4cm);
\draw [outputleg] (q6) to +(-.5cm,.4cm);
\draw [outputleg] (q6) to +(.5cm,.4cm);
\draw [inputleg] (q2) to +(-.5cm,-.3cm);
\draw [inputleg] (q2) to +(.5cm,-.3cm);
\draw [inputleg] (q5) to +(-.5cm,-.3cm);
\draw [inputleg] (q5) to +(.5cm,-.3cm);
\draw [inputleg] (p11) to +(-.6cm,-.4cm);
\draw [inputleg] (p11) to +(.6cm,-.4cm);
\draw [inputleg] (p12) to +(-.6cm,-.4cm);
\draw [inputleg] (p12) to +(.6cm,-.4cm);
\end{tikzpicture}
\end{center}

In the above decorated graph, consider the two copies of $p_2$ and the four copies of $q$ that are directly connected to them, i.e., the middle four copies of $q$.  The generating distributivity involving $p_2$ and $q$ now implies that the above decorated graph is identified with the following decorated graph:
\begin{center}
\begin{tikzpicture}
\matrix[row sep=1cm,column sep=1.5cm] {
\node [plain] (q1) {$q$};
& \node [plain] (q2) {$q$}; &
& \node [plain] (q3) {$q$};
& \node [plain] (q4) {$q$};\\
&&&& \\
& \node [plain] (p21) {$p_2$}; &&
 \node [plain] (p22) {$p_2$}; &\\
& \node [plain] (p11) {$p_1$}; &&
 \node [plain] (p12) {$p_1$}; &\\
};
\draw [arrow] (p21) to (q2);
\draw [arrow] (p21) to (q3);
\draw [arrow] (p22) to (q2);
\draw [arrow] (p22) to (q3);
\draw [arrow,bend left=20] (p11) to (q1);
\draw [arrow] (p11) to (q4);
\draw [arrow,bend left=20] (p11) to (p21);
\draw [arrow,bend right=20] (p11) to (p21);
\draw [arrow] (p12) to (q1);
\draw [arrow,bend right=20] (p12) to (q4);
\draw [arrow,bend left=20] (p12) to (p22);
\draw [arrow,bend right=20] (p12) to (p22);
\draw [outputleg] (q1) to +(-.5cm,.4cm);
\draw [outputleg] (q1) to +(.5cm,.4cm);
\draw [outputleg] (q2) to +(-.5cm,.4cm);
\draw [outputleg] (q2) to +(.5cm,.4cm);
\draw [outputleg] (q3) to +(-.5cm,.4cm);
\draw [outputleg] (q3) to +(.5cm,.4cm);
\draw [outputleg] (q4) to +(-.5cm,.4cm);
\draw [outputleg] (q4) to +(.5cm,.4cm);
\draw [inputleg] (p21) to +(-.5cm,-.4cm);
\draw [inputleg] (p21) to +(.5cm,-.4cm);
\draw [inputleg] (p22) to +(-.5cm,-.4cm);
\draw [inputleg] (p22) to +(.5cm,-.4cm);
\draw [inputleg] (p11) to +(-.5cm,-.4cm);
\draw [inputleg] (p11) to +(.5cm,-.4cm);
\draw [inputleg] (p12) to +(-.5cm,-.4cm);
\draw [inputleg] (p12) to +(.5cm,-.4cm);
\end{tikzpicture}
\end{center}
Up to input and output relabeling, this is the other side of distributivity involving $p$ and $q$.  Therefore, we have shown that the distributivity relation involving $p$ and $q$ decomposes into two generating distributivity relations.

Finally, the general case of $(n,m) = (1,2)$ can be proved by a slight modification of the above four decorated graphs.  For example, if $q$ has $3$ outputs instead of $2$, then in the first decorated graph for distributivity we begin with $3$ copies of $p$:
\begin{center}
\begin{tikzpicture}
\matrix[row sep=1.5cm,column sep=1.5cm] {
\node [plain] (p21) {$p_2$}; 
&& \node [plain] (p2i) {$p_2$}; 
&& \node [plain] (p22) {$p_2$}; \\
\node [plain] (p11) {$p_1$}; 
&& \node [plain] (p1i) {$p_1$}; 
&& \node [plain] (p12) {$p_1$}; \\
&&&& \\
&&&& \\
\node [plain] (q1) {$q$}; 
& \node [plain] (q2) {$q$}; 
&& \node [plain] (q3) {$q$}; 
& \node [plain] (q4) {$q$};\\
};
\draw [arrow,bend left=20] (p11) to (p21);
\draw [arrow,bend right=20] (p11) to (p21);
\draw [arrow,bend left=20] (p1i) to (p2i);
\draw [arrow,bend right=20] (p1i) to (p2i);
\draw [arrow,bend left=20] (p12) to (p22);
\draw [arrow,bend right=20] (p12) to (p22);
\draw [arrow,bend left=45] (q1) to (p21);
\draw [arrow,bend left=20] (q1) to (p2i);
\draw [arrow] (q1) to (p22);
\draw [arrow,bend left=20] (q2) to (p11);
\draw [arrow,bend left=20] (q2) to (p1i);
\draw [arrow] (q2) to (p12);
\draw [arrow] (q3) to (p11);
\draw [arrow,bend right=20] (q3) to (p1i);
\draw [arrow,bend right=20] (q3) to (p12);
\draw [arrow] (q4) to (p21);
\draw [arrow,bend right=20] (q4) to (p2i);
\draw [arrow,bend right=45] (q4) to (p22);
\draw [outputleg] (p21) to +(-.6cm,.4cm);
\draw [outputleg] (p21) to +(.6cm,.4cm);
\draw [outputleg] (p2i) to +(-.6cm,.4cm);
\draw [outputleg] (p2i) to +(.6cm,.4cm);
\draw [outputleg] (p22) to +(-.6cm,.4cm);
\draw [outputleg] (p22) to +(.6cm,.4cm);
\draw [outputleg] (p11) to +(-.6cm,.4cm);
\draw [outputleg] (p11) to +(.6cm,.4cm);
\draw [outputleg] (p1i) to +(-.5cm,.4cm);
\draw [outputleg] (p1i) to +(.5cm,.4cm);
\draw [outputleg] (p12) to +(-.6cm,.4cm);
\draw [outputleg] (p12) to +(.6cm,.4cm);
\draw [inputleg] (q1) to +(-.6cm,-.4cm);
\draw [inputleg] (q1) to +(.6cm,-.4cm);
\draw [inputleg] (q2) to +(-.6cm,-.4cm);
\draw [inputleg] (q2) to +(.6cm,-.4cm);
\draw [inputleg] (q3) to +(-.6cm,-.4cm);
\draw [inputleg] (q3) to +(.6cm,-.4cm);
\draw [inputleg] (q4) to +(-.6cm,-.4cm);
\draw [inputleg] (q4) to +(.6cm,-.4cm);
\end{tikzpicture}
\end{center}
The two generating distributivity are very similar to the ones presented above.  For instance, the first generating distributivity involves the three copies of $p_1$ and the two middle copies of $q$.  In other words, in the general case the two generating distributivity steps are essentially identical to the ones above, but the horizontal sizes of the decorated graphs may change.
\end{proof}

We next extend the previous lemma by allowing $G$ to have more than two vertices.

\begin{lemma}
\label{distnone}
Suppose $n  = |\vertex(H)| = 1$ and $m = |\vertex(G)| \geq 2$.  Then the distributivity relation involving $p = \{p_v\}$ and $q = \{q_u\}$ decomposes into a finite string of generating distributivity relations.
\end{lemma}

\begin{proof}
The proof is an induction on $m \geq 2$.  The initial case $m=2$ is Lemma \ref{gendisttypical}.  For the induction step, suppose $|\vertex(G)| > 2$.  Suppose $v$ is an almost isolated vertex in $G$, which exists by Corollary \ref{aiexist}.  We assume that $v$ is weakly terminal.  The case where $v$ is weakly initial has a similar proof.

Since $v$ is weakly terminal, it has only incoming but not outgoing ordinary edges.  Moreover, deleting it from $G$ yields a connected wheel-free graph $K$.  Therefore, $G$ has the following form:
\begin{center}
\begin{tikzpicture}
\matrix[row sep=.5cm,column sep=1cm] {
\node [plain,label=above:$...$] (v) {$v$}; \\
\node [empty] (dots) {$\cdots$}; \\ 
\node [fatplain,label=below:$...$] (K) {$K$}; \\
};
\draw [arrow,bend left=40] (p1) to (p2);
\draw [arrow,bend right=40] (p1) to (p2);
\draw [inputleg] (p2) to +(-.6cm,-.4cm);
\draw [inputleg] (p2) to +(.6cm,-.4cm);
\draw [outputleg] (p2) to +(-.6cm,.4cm);
\draw [outputleg] (p2) to +(.6cm,.4cm);
\draw [outputleg] (p1) to +(-.6cm,.4cm);
\draw [outputleg] (p1) to +(.6cm,.4cm);
\draw [inputleg] (p1) to +(-.6cm,-.4cm);
\draw [inputleg] (p1) to +(.6cm,-.4cm);
\end{tikzpicture}
\end{center}
As in the proof of Lemma \ref{gendisttypical}, for ease of drawing the decorated graphs, we assume that $v$ has two outputs and four inputs, two of which are ordinary edges connected to $K$, and that $K$ has two inputs, two outputs, and two ordinary edges connected to $v$.  We will therefore draw $p$ as the following decorated graph:
\begin{center}
\begin{tikzpicture}
\matrix[row sep=1.5cm,column sep=1cm] {
\node [plain] (pv) {$p_v$}; \\
\node [fatplain] (pK) {$p_K$}; \\
};
\draw [arrow,bend left=40] (pK) to (pv);
\draw [arrow,bend right=40] (pK) to (pv);
\draw [inputleg] (pv) to +(-.6cm,-.4cm);
\draw [inputleg] (pv) to +(.6cm,-.4cm);
\draw [outputleg] (pv) to +(-.6cm,.4cm);
\draw [outputleg] (pv) to +(.6cm,.4cm);
\draw [outputleg] (pK) to +(-.6cm,.4cm);
\draw [outputleg] (pK) to +(.6cm,.4cm);
\draw [inputleg] (pK) to +(-.6cm,-.4cm);
\draw [inputleg] (pK) to +(.6cm,-.4cm);
\end{tikzpicture}
\end{center}
Here the decorated graph
\[
p_K = \{p_w\}_{w \in \vertex(K)} \in \phat[K]
\]
is the restriction of $p$ to $K$.

One side of distributivity involving $p$ and $q$ is the decorated graph:
\begin{center}
\begin{tikzpicture}
\matrix[row sep=1cm,column sep=1.5cm] {
& \node [plain] (p21) {$p_v$}; &
& \node [plain] (p22) {$p_v$}; & \\
& \node [fatplain] (p11) {$p_K$}; &
& \node [fatplain] (p12) {$p_K$}; & \\
&&&& \\
\node [plain] (q1) {$q$}; 
& \node [plain] (q2) {$q$}; 
&& \node [plain] (q3) {$q$}; 
& \node [plain] (q4) {$q$};\\
};
\draw [arrow,bend left=20] (p11) to (p21);
\draw [arrow,bend right=20] (p11) to (p21);
\draw [arrow,bend left=20] (p12) to (p22);
\draw [arrow,bend right=20] (p12) to (p22);
\draw [arrow,bend left=20] (q1) to (p21);
\draw [arrow] (q1) to (p22);
\draw [arrow] (q2) to (p11);
\draw [arrow] (q2) to (p12);
\draw [arrow] (q3) to (p11);
\draw [arrow] (q3) to (p12);
\draw [arrow] (q4) to (p21);
\draw [arrow,bend right=20] (q4) to (p22);
\draw [outputleg] (p21) to +(-.6cm,.4cm);
\draw [outputleg] (p21) to +(.6cm,.4cm);
\draw [outputleg] (p22) to +(-.6cm,.4cm);
\draw [outputleg] (p22) to +(.6cm,.4cm);
\draw [outputleg] (p11) to +(-.6cm,.4cm);
\draw [outputleg] (p11) to +(.6cm,.4cm);
\draw [outputleg] (p12) to +(-.6cm,.4cm);
\draw [outputleg] (p12) to +(.6cm,.4cm);
\draw [inputleg] (q1) to +(-.6cm,-.4cm);
\draw [inputleg] (q1) to +(.6cm,-.4cm);
\draw [inputleg] (q2) to +(-.6cm,-.4cm);
\draw [inputleg] (q2) to +(.6cm,-.4cm);
\draw [inputleg] (q3) to +(-.6cm,-.4cm);
\draw [inputleg] (q3) to +(.6cm,-.4cm);
\draw [inputleg] (q4) to +(-.6cm,-.4cm);
\draw [inputleg] (q4) to +(.6cm,-.4cm);
\end{tikzpicture}
\end{center}
The rest of the proof is similar to that of Lemma \ref{gendisttypical}, with $p_v$ and $p_K$ playing the roles of $p_2$ and $p_1$, respectively.

Since $|\vertex(K)| = |\vertex(G)|-1$, the induction hypothesis applies to $p_K$ and $q$.  In other words, the distributivity relation involving $p_K$ and $q$ decomposes into a finite string of generating distributivity relations.  Using these generating distributivity, the above decorated graph is identified with the following decorated graph:
\begin{center}
\begin{tikzpicture}
\matrix[row sep=2cm,column sep=1.5cm] {
& \node [plain] (p21) {$p_v$}; &&&  \node [plain] (p22) {$p_v$}; &\\
\node [plain] (q1) {$q$};
& \node [plain] (q2) {$q$};
& \node [plain] (q3) {$q$};
& \node [plain] (q4) {$q$};
& \node [plain] (q5) {$q$};
& \node [plain] (q6) {$q$};\\
& \node [fatplain] (p11) {$p_K$}; &&&  \node [fatplain] (p12) {$p_K$}; &\\
};
\draw [arrow,bend left=20] (q1) to (p21);
\draw [arrow] (q1) to (p22);
\draw [arrow] (q3) to (p21);
\draw [arrow] (q3) to (p22);
\draw [arrow] (q4) to (p21);
\draw [arrow] (q4) to (p22);
\draw [arrow] (q6) to (p21);
\draw [arrow,bend right=20] (q6) to (p22);
\draw [arrow] (p11) to (q2);
\draw [arrow] (p11) to (q3);
\draw [arrow] (p11) to (q4);
\draw [arrow] (p11) to (q5);
\draw [arrow] (p12) to (q2);
\draw [arrow] (p12) to (q3);
\draw [arrow] (p12) to (q4);
\draw [arrow] (p12) to (q5);
\draw [outputleg] (p21) to +(-.6cm,.4cm);
\draw [outputleg] (p21) to +(.6cm,.4cm);
\draw [outputleg] (p22) to +(-.6cm,.4cm);
\draw [outputleg] (p22) to +(.6cm,.4cm);
\draw [outputleg] (q2) to +(-.5cm,.3cm);
\draw [outputleg] (q2) to +(.5cm,.3cm);
\draw [outputleg] (q5) to +(-.5cm,.3cm);
\draw [outputleg] (q5) to +(.5cm,.3cm);
\draw [inputleg] (q1) to +(-.6cm,-.4cm);
\draw [inputleg] (q1) to +(.6cm,-.4cm);
\draw [inputleg] (q6) to +(-.6cm,-.4cm);
\draw [inputleg] (q6) to +(.6cm,-.4cm);
\draw [inputleg] (p11) to +(-.6cm,-.4cm);
\draw [inputleg] (p11) to +(.6cm,-.4cm);
\draw [inputleg] (p12) to +(-.6cm,-.4cm);
\draw [inputleg] (p12) to +(.6cm,-.4cm);
\end{tikzpicture}
\end{center}

Next we apply the generating distributivity relation involving $p_v$ and $q$.  In the previous decorated graph, this is applied to the two copies of $p_v$ and the four copies of $q$ directly connected to them.  It is then identified with the following decorated graph:
\begin{center}
\begin{tikzpicture}
\matrix[row sep=1cm,column sep=1.5cm] {
\node [plain] (q1) {$q$};
& \node [plain] (q2) {$q$}; &
& \node [plain] (q3) {$q$};
& \node [plain] (q4) {$q$};\\
&&&& \\
& \node [plain] (p21) {$p_v$}; &&
 \node [plain] (p22) {$p_v$}; &\\
& \node [fatplain] (p11) {$p_K$}; &&
 \node [fatplain] (p12) {$p_K$}; &\\
};
\draw [arrow] (p21) to (q2);
\draw [arrow] (p21) to (q3);
\draw [arrow] (p22) to (q2);
\draw [arrow] (p22) to (q3);
\draw [arrow,bend left=20] (p11) to (q1);
\draw [arrow] (p11) to (q4);
\draw [arrow,bend left=20] (p11) to (p21);
\draw [arrow,bend right=20] (p11) to (p21);
\draw [arrow] (p12) to (q1);
\draw [arrow,bend right=20] (p12) to (q4);
\draw [arrow,bend left=20] (p12) to (p22);
\draw [arrow,bend right=20] (p12) to (p22);
\draw [outputleg] (q1) to +(-.5cm,.4cm);
\draw [outputleg] (q1) to +(.5cm,.4cm);
\draw [outputleg] (q2) to +(-.5cm,.4cm);
\draw [outputleg] (q2) to +(.5cm,.4cm);
\draw [outputleg] (q3) to +(-.5cm,.4cm);
\draw [outputleg] (q3) to +(.5cm,.4cm);
\draw [outputleg] (q4) to +(-.5cm,.4cm);
\draw [outputleg] (q4) to +(.5cm,.4cm);
\draw [inputleg] (p21) to +(-.5cm,-.4cm);
\draw [inputleg] (p21) to +(.5cm,-.4cm);
\draw [inputleg] (p22) to +(-.5cm,-.4cm);
\draw [inputleg] (p22) to +(.5cm,-.4cm);
\draw [inputleg] (p11) to +(-.5cm,-.4cm);
\draw [inputleg] (p11) to +(.5cm,-.4cm);
\draw [inputleg] (p12) to +(-.5cm,-.4cm);
\draw [inputleg] (p12) to +(.5cm,-.4cm);
\end{tikzpicture}
\end{center}
This last decorated graph is the other side of the distributivity relation involving $p$ and $q$, except for the input and output relabeling.

If the almost isolated vertex $v$ is weakly initial, then the roles of $p_K$ and $p_v$ are reversed.  In this case, starting at one side of the distributivity relation for $p$ and $q$, we first apply the generating distributivity for $p_v$ and $q$.  Next we apply the induction hypothesis to $p_K$ and $q$.

Finally, for the general case of $n=1$ and $m \geq 2$, we appropriately adjust the numbers of copies of $p$ and $q$ in the first decorated graph for distributivity.  The next two steps, using the induction hypothesis and a generating distributivity, are essentially identical to the ones given above with minor adjustment in the horizontal sizes of the decorated graphs.
\end{proof}

Next we consider the dual case.

\begin{lemma}
\label{distmone}
Suppose $n  = |\vertex(H)| \geq 2$ and $m = |\vertex(G)| = 1$.  Then the distributivity relation involving $p = \{p_v\}$ and $q = \{q_u\}$ decomposes into a finite string of generating distributivity relations.
\end{lemma}

\begin{proof}
The proof is a minor modification of those of Lemmas \ref{gendisttypical} and \ref{distnone}.  More precisely, the first step is to establish the special case $(n,m) = (2,1)$, which is similar to Lemma \ref{gendisttypical}.  This time $p$ is a decorated corolla
\begin{center}
\begin{tikzpicture}
\node [plain,label=above:$...$,label=below:$...$] (q) {$p$};
\draw [inputleg] (q) to +(-.6cm,-.4cm);
\draw [inputleg] (q) to +(.6cm,-.4cm);
\draw [outputleg] (q) to +(-.6cm,.4cm);
\draw [outputleg] (q) to +(.6cm,.4cm);
\end{tikzpicture}
\end{center} 
while $q$ is a decorated partially grafted corollas:
\begin{center}
\begin{tikzpicture}
\matrix[row sep=.5cm,column sep=1cm] {
\node [plain,label=above:$...$] (p2) {$q_2$}; \\
\node [empty] (dots) {$\cdots$}; \\ 
\node [plain,label=below:$...$] (p1) {$q_1$}; \\
};
\draw [arrow,bend left=40] (p1) to (p2);
\draw [arrow,bend right=40] (p1) to (p2);
\draw [inputleg] (p2) to +(-.6cm,-.4cm);
\draw [inputleg] (p2) to +(.6cm,-.4cm);
\draw [outputleg] (p2) to +(-.6cm,.4cm);
\draw [outputleg] (p2) to +(.6cm,.4cm);
\draw [outputleg] (p1) to +(-.6cm,.4cm);
\draw [outputleg] (p1) to +(.6cm,.4cm);
\draw [inputleg] (p1) to +(-.6cm,-.4cm);
\draw [inputleg] (p1) to +(.6cm,-.4cm);
\end{tikzpicture}
\end{center}
The distributivity relation involving $p$ and $q$ then decomposes into two generating distributivity relations.  In fact, one may simply reuse the proof of Lemma \ref{gendisttypical} with the roles of $p$ and $q$ reversed, and with the orientations of all the edges reversed.

The general case of $n \geq 2$ and $m=1$ is proved by induction on $n$, with the initial case $n=2$ as in the previous paragraph.  For the induction step, one first decomposes $H$ into a partial grafting of an almost isolated vertex $u$ and a connected wheel-free graph $J$, which must have one fewer vertex than $H$.  For example, if the almost isolated vertex $u$ is weakly terminal, then $H$ decomposes as follows:
\begin{center}
\begin{tikzpicture}
\matrix[row sep=.5cm,column sep=1cm] {
\node [plain,label=above:$...$] (v) {$u$}; \\
\node [empty] (dots) {$\cdots$}; \\ 
\node [fatplain,label=below:$...$] (K) {$J$}; \\
};
\draw [arrow,bend left=40] (p1) to (p2);
\draw [arrow,bend right=40] (p1) to (p2);
\draw [inputleg] (p2) to +(-.6cm,-.4cm);
\draw [inputleg] (p2) to +(.6cm,-.4cm);
\draw [outputleg] (p2) to +(-.6cm,.4cm);
\draw [outputleg] (p2) to +(.6cm,.4cm);
\draw [outputleg] (p1) to +(-.6cm,.4cm);
\draw [outputleg] (p1) to +(.6cm,.4cm);
\draw [inputleg] (p1) to +(-.6cm,-.4cm);
\draw [inputleg] (p1) to +(.6cm,-.4cm);
\end{tikzpicture}
\end{center}
Next, as in the proof of Lemma \ref{distnone}, one reuses essentially the same argument as in the special case $(n,m) = (2,1)$, using the induction hypothesis once and a generating distributivity.
\end{proof}

We now consider the general case.

\begin{lemma}
\label{distmngeneral}
The distributivity relation involving $p = \{p_v\}$ and $q = \{q_u\}$ decomposes into a finite string of generating distributivity relations.
\end{lemma}

\begin{proof}
We consider the induction step for the induction on $m+n$ started just before Lemma \ref{gendisttypical}.  Because of Lemmas \ref{distnone} and \ref{distmone}, we may further assume that both $n \geq 2$ and $m \geq 2$.  By Corollary \ref{aiexist} and Theorem \ref{thm:opropfact}, $G$ has an outer properadic factorization
\[
G = P(\{C_v,K\}),
\]
so $P$ is a partially grafted corollas, $C_v$ is a corolla, and $K$ is an ordinary connected wheel-free graph.  In particular, $|\vertex(K)| = |\vertex(G)| - 1$, and $v$ is an almost isolated vertex.  Then we recycle the proof of Lemma \ref{distnone} essentially verbatim.  More precisely, the distributivity relation involving $p_K$ and $q$ decomposes into finitely many generating distributivity relations by the induction hypothesis.  Next, the distributivity relation involving $p_v$ and $q$ decomposes into finitely many generating distributivity relations by the induction hypothesis again. (Alternatively, one may also use Lemmas \ref{distnone} and \ref{distmone} for this last step because $v$ is a single vertex.)  This finishes the induction step.
\end{proof}

The proof of Theorem \ref{distsimplified} is now complete.

\begin{remark}
The above proof of Theorem \ref{distsimplified} reveals a bit more than the assertion.  Namely, the distributivity relation involving $p$ and $q$ decomposes into $mn$ generating distributivity relations.
\end{remark}

\section{Internal Hom of Special Properads}
\label{sec:properadhom}

In this section, we observe that the category $\properads$ of special properads (Definition \ref{rk:properadio}) is symmetric monoidal \emph{closed} by constructing an internal hom of special properads.

\subsection{Natural Transformation}

Recall that a properad is a generalization of a category in which morphisms have finite lists of objects as source and target.  For ordinary categories, functors and natural transformations provide the internal hom.  For properads, the analogs of functors are maps between properads.  So to define the internal hom of special  properads, we still need a properadic analog of natural transformations, which we now define.

\begin{definition}
\label{def:naturaltransf}
Suppose $\sP$ is a $\fC$-colored special properad, and $\sQ$ is a $\fD$-colored special properad.
\begin{enumerate}
\item
Suppose $f_i, g_j \in \properads(\sP,\sQ)$ with $1 \leq i \leq k$, $1 \leq j \leq l$, and $k,l \geq 1$.  A \textbf{transformation} \index{transformation} $\psi$ with profiles $(\uf;\ug)$ is a function that assigns to each color $c \in \fC$ an element
\[
\psi_c \in \sQ\binom{\ug(c)}{\uf(c)} = \sQ\binom{g_1(c),\ldots ,g_l(c)}{f_1(c), \ldots ,f_k(c)}.
\]
\item
Suppose $p \in \sP(\uc_{[1,m]};\ud_{[1,n]})$, and $\psi$ is a transformation with profiles $(\uf;\ug)$.  Then we say $\psi$ is \textbf{natural} with respect to $p$ if the equality
\begin{equation}
\label{naturaltransf}
\begin{split}
& \gamma_T\left(\{g_j(p)\}_{j=1}^l \times 
\{\psi_{c_i}\}_{i=1}^m\right)\\
&= 
\gamma_{T'}\left(\sigma^l_n \left[
\{\psi_{d_j}\}_{j=1}^n \times \{f_i(p)\}_{i=1}^k\right] 
\sigma^k_m\right)
\end{split}
\end{equation}
holds in
\[
\sQ\binom{g_1(\ud), \ldots , g_l(\ud)}{\uf(c_1), \ldots , \uf(c_m)}.
\]
Here $\gamma_T$ and $\gamma_{T'}$ are properad structure maps of $\sQ$, and $T$ is the following $\sQ$-decorated graph.
\begin{center}
\begin{tikzpicture}
\matrix[row sep=2cm,column sep=2cm] {
\node [fatplain,label=above:$g_1(\ud)$,label=below:$g_1(\uc)$] (pd1) {$g_1(p)$}; &
\node [empty] (pd) {$\cdots$}; &
\node [fatplain,label=above:$g_l(\ud)$,label=below:$g_l(\uc)$] (pdn) {$g_l(p)$}; \\
\node [fatplain,label=below:$\uf(c_1)$,label=above:$\ug(c_1)$] (a1q) {$\psi_{c_1}$}; & \node [empty] (aq) {$\cdots$}; &
\node [fatplain,label=below:$\uf(c_m)$,label=above:$\ug(c_m)$] (akq) {$\psi_{c_m}$};
\\
};
\draw [arrow,bend left=45] (a1q) to (pd1);
\draw [arrow] (a1q) to (pdn);
\draw [arrow] (akq) to (pd1);
\draw [arrow,bend right=45] (akq) to (pdn);
\draw [outputleg] (pd1) to +(-.8cm,.5cm);
\draw [outputleg] (pd1) to +(.8cm,.5cm);
\draw [outputleg] (pdn) to +(-.8cm,.5cm);
\draw [outputleg] (pdn) to +(.8cm,.5cm);
\draw [inputleg] (a1q) to +(-.8cm,-.5cm);
\draw [inputleg] (a1q) to +(.8cm,-.5cm);
\draw [inputleg] (akq) to +(-.8cm,-.5cm);
\draw [inputleg] (akq) to +(.8cm,-.5cm);
\end{tikzpicture}
\end{center}
In the top (resp., bottom) row, the $i$th vertex from the left is decorated by $g_i(p)$ (resp., $\psi_{c_i}$).  The $r$th output of the vertex $\psi_{c_i}$ is connected to the $i$th input of $g_r(p)$.  This internal edge has color $g_r(c_i)$.

Likewise, \emph{before} applying the input and output relabeling permutations $\sigma^l_n$ and $\sigma^k_m$, $T'$ is the following $\sQ$-decorated graph.
\begin{center}
\begin{tikzpicture}
\matrix[row sep=2cm,column sep=2cm] {
\node [fatplain,label=above:$\ug(d_1)$,label=below:$\uf(d_1)$] (b1q) {$\psi_{d_1}$}; &
\node [empty] (bq) {$\cdots$}; &
\node [fatplain,label=above:$\ug(d_n)$,label=below:$\uf(d_n)$] (blq) {$\psi_{d_n}$}; \\
\node [fatplain,label=above:$f_1(\ud)$,label=below:$f_1(\uc)$] (pc1) {$f_1(p)$}; & \node [empty] (pc) {$\cdots$}; &
\node [fatplain,label=above:$f_k(\ud)$,label=below:$f_k(\uc)$] (pcm) {$f_k(p)$};
\\
};
\draw [arrow,bend left=45] (pc1) to  (b1q);
\draw [arrow] (pc1) to (blq);
\draw [arrow] (pcm) to (b1q);
\draw [arrow,bend right=45] (pcm) to (blq);
\draw [outputleg] (b1q) to +(-.8cm,.5cm);
\draw [outputleg] (b1q) to +(.8cm,.5cm);
\draw [outputleg] (blq) to +(-.8cm,.5cm);
\draw [outputleg] (blq) to +(.8cm,.5cm);
\draw [inputleg] (pc1) to +(-.8cm,-.5cm);
\draw [inputleg] (pc1) to +(.8cm,-.5cm);
\draw [inputleg] (pcm) to +(-.8cm,-.5cm);
\draw [inputleg] (pcm) to +(.8cm,-.5cm);
\end{tikzpicture}
\end{center}
In the top (resp., bottom) row, the $j$th vertex from the left is decorated by $\psi_{d_j}$ (resp., $f_j(p)$).  The $s$th output of the vertex $f_j(p)$ is connected to the $j$th input of $\psi_{d_s}$.  This internal edge has color $f_j(d_s)$.
\item
A \textbf{natural transformation}  \index{natural transformation} is a transformation that is natural with respect to all elements in $\sP$ for which \eqref{naturaltransf} is defined.
\end{enumerate}
\end{definition}

\begin{remark}
\begin{enumerate}
\item
If $k=l=m=n=1$, then the naturality condition \eqref{naturaltransf} is the usual commutative square that defines a natural transformation in a category.
\item
If $l = n = 1$,  then the naturality condition \eqref{naturaltransf} becomes the one for operads  \cite{weiss}.
\item
Natural transformations of PROPs are defined in  \cite{hr1}.
\item
The naturality condition \eqref{naturaltransf} is very similar to the distributivity relation (section \ref{subsec:distributivity}).  In fact, the two graphs $T$ and $T'$ (but of course not their decorations) are the same as in the distributivity relation.  Likewise, the restrictions are introduced to make sure that the graphs $T$ and $T'$ are connected.  A consequence of this close resemblance between distributivity and naturality is that we will be able to extend naturality, just like the way we extended generating distributivity to distributivity in section \ref{sec:prdistsimplified}.
\end{enumerate}
\end{remark}

The following observation about extending natural transformations along connected wheel-free graphs is needed to define the properad structure on the internal hom that we will define shortly. Roughly speaking, it says that any coherent composition of natural transformations is again a natural transformation. In the very special case of natural transformations of categories, the following observation is the simple fact that natural transformations are closed under vertical composition.

To extend naturality along connected wheel-free graphs, we will need the following definition.

\begin{definition}
Suppose $\sP$ is a $\fC$-colored special properad, $\sQ$ is a $\fD$-colored special properad, $G$ is a $\properads(\sP,\sQ)$-colored connected wheel-free graph, and $c \in \fC$.  Define the $\fD$-colored connected wheel-free graph $G_c$ as the one obtained from $G$ by applying its coloring to $c$.
\end{definition}

\begin{lemma}
\label{extendnaturality}
Suppose
\begin{itemize}
\item
$\sP$ is a $\fC$-colored special properad,
\item
$\sQ$ is a $\fD$-colored special properad,
\item
$G$ is a $\properads(\sP,\sQ)$-colored connected wheel-free graph with profiles $(\uf;\ug)$ in which each vertex $v$ has non-empty $\inp(v)$ and $\out(v)$, and
\item
$\psi^v$ is a natural transformation with the same profiles as the vertex $v$ for each $v \in \vertex(G)$.
\end{itemize}
Then the transformation $\psi^G$ defined by
\begin{equation}
\label{gammainternalhom}
\psi^G_c 
\defn 
\gamma_{G_c}\left(\left\{\psi^v_c\right\}\right)
\quad
\text{for $c \in \fC$}
\end{equation}
is a natural transformation, where $\gamma_{G_c}$ is a properad structure map in $\sQ$.
\end{lemma}

\begin{proof}
This is proved by induction on $|\vertex(G)|$.  First note that if $|\vertex(G)| \leq 1$ (i.e., $G$ is an exceptional edge or a corolla), then the assertion is trivially true.

So we may assume $|\vertex(G)| \geq 2$.  Pick an element $p \in \sP$, which is as usual depicted as a decorated corolla.  
\begin{center}
\begin{tikzpicture}
\node [plain,label=above:$...$,label=below:$...$] (q) {$p$};
\draw [inputleg] (q) to +(-.6cm,-.4cm);
\draw [inputleg] (q) to +(.6cm,-.4cm);
\draw [outputleg] (q) to +(-.6cm,.4cm);
\draw [outputleg] (q) to +(.6cm,.4cm);
\end{tikzpicture}
\end{center} 
To show that $\psi^G$ is natural, we must show that the two decorated graphs in Definition \ref{def:naturaltransf}, with $\psi^G$ in place of $\psi$, become equal  after applying the properad structure map $\gamma$ in $\sQ$.  To prove this, take an outer properadic factorization
\[
G = P\left(\{C_v, K\}\right),
\]
in which $P$ is a partially grafted corolla, $K$ is ordinary connected wheel-free with one fewer vertex than $G$, $C_v$ is a corolla, and $v$ is almost isolated in $G$.  Such a factorization exists by Corollary \ref{aiexist} and Theorem \ref{thm:opropfact}.  Thus, assuming $v$ is weakly terminal, we may graphically represent the transformation $\psi^G$ as the following decorated graph.
\begin{center}
\begin{tikzpicture}
\matrix[row sep=1cm,column sep=1cm] {
\node [fatplain, label=above:$...$, label=below:$...$] (pv) {$\psi^v$}; \\
\node [fatplain, label=above:$...$, label=below:$...$] (pK) {$\psi^K$}; \\
};
\draw [arrow,bend left=40] (pK) to (pv);
\draw [arrow,bend right=40] (pK) to (pv);
\draw [inputleg] (pv) to +(-.7cm,-.4cm);
\draw [inputleg] (pv) to +(.7cm,-.4cm);
\draw [outputleg] (pv) to +(-.7cm,.4cm);
\draw [outputleg] (pv) to +(.7cm,.4cm);
\draw [outputleg] (pK) to +(-.7cm,.4cm);
\draw [outputleg] (pK) to +(.7cm,.4cm);
\draw [inputleg] (pK) to +(-.7cm,-.4cm);
\draw [inputleg] (pK) to +(.7cm,-.4cm);
\end{tikzpicture}
\end{center}
If $v$ is weakly initial, then $\psi^v$ decorates the bottom vertex, and $\psi^K$ decorates the top vertex. We know $\psi^v$ is natural, and $\psi^K$ is natural by the induction hypothesis.

Recall that the graphs (but not the decorations) in the distributivity relation \eqref{distributivity} and naturality \eqref{naturaltransf} are the same.  Therefore, using the associativity of the properad structure map $\gamma$ in $\sQ$, we may recycle the proof of Lemma \ref{distnone} to show that $\psi^G$ is natural.  Indeed, all we really have to do is to replace $(q, p_v, p_K)$ in the proof of Lemma  \ref{distnone} by $\left(p, \psi^v, \psi^K\right)$ here, using the associativity of $\gamma$ every time naturality is used to replace one decorated graph by another.
\end{proof}

\subsection{Internal Hom}

\begin{definition}
\label{def:gpropinthom}
Suppose $\sP$ and $\sQ$ are special properads.  Define the \textbf{internal hom} \index{internal hom} $\Hom(\sP,\sQ)$\label{note:hom} as the $\sS(\properads(\sP,\sQ))$-colored object whose elements are natural transformations in the sense of Definition \ref{def:naturaltransf}.
\end{definition}

First we observe that the internal hom has a canonical properad structure.

\begin{lemma}
\label{lem:propginternalhom}
Suppose $\sP$ and $\sQ$ are special properads.  Then the internal hom $\Hom(\sP,\sQ)$ is a $\properads(\sP,\sQ)$-colored special properad with structure map
\[
\gamma_G\left(\left\{\psi^v\right\}_{\ving}\right)
\defn 
\psi^G
\]
as in \eqref{gammainternalhom}.
\end{lemma}

\begin{proof}
Lemma \ref{extendnaturality} guarantees that, if each $\psi^v$ for $v \in \vertex(G)$ is a natural transformation, then so is $\psi^G$.  Associativity and unity of the structure maps $\gamma$ on $\Hom(\sP,\sQ)$ are implied by those of the properad $\sQ$ because
\[
\left[\gamma_G\left(\left\{\psi^v\right\}_{\ving}\right)\right]_c
= \psi^G_c
=
\gamma_{G_c}\left( \left\{ \psi^v_c\right\} \right),
\]
in which the right-hand side is an operation in $\sQ$.
\end{proof}

\subsection{Symmetric Monoidal Closed Structure}

\begin{theorem}
\label{thm:gpropmonoidalclosed}
The category $\properads$ is symmetric monoidal closed with monoidal product $\otimes$ and internal hom $\Hom$.
\end{theorem}

\begin{proof}
Use Theorem \ref{thm:propgmonoidal} and Lemma \ref{lem:propginternalhom}.  The required natural bijection
\begin{equation}
\label{gpropclosed}
\properads\left(\sP \otimes \sQ, \sR\right) 
\cong 
\properads\left(\sP, \Hom(\sQ,\sR)\right)
\end{equation}
follows by unwrapping the constructions of $\otimes$ \eqref{gproptensor} and $\Hom$ (Definition \ref{def:naturaltransf} and \eqref{gammainternalhom}).  To be more explicit, the correspondence goes as follows.  

Suppose $\sP$, $\sQ$, and $\sR$ have color sets $\fC$, $\fD$, and $\fE$, respectively.  Given a properad map $\nicexy{\sP \otimes \sQ \ar[r]^-{\varphi} & \sR}$, the corresponding map $\nicexy{\sP \ar[r]^-{\varphi^\sharp} & \Hom(\sQ,\sR)}$ of properads is defined as follows:
\begin{itemize}
\item
For a color $c \in \fC$, the color $\varphi^\sharp_0(c)$ of the internal hom is the properad map composition
\[
\nicexy{
\sQ \ar[r]^-{c \otimes (-)} & \sP \otimes \sQ \ar[r]^-{\varphi} & \sR.
}\]
\item
For an element $p \in \sP$, the element $\varphi^\sharp(p)$ of the internal hom is the transformation
\[
\varphi(p \otimes -) = \left\{ \fD \ni d \longmapsto \varphi(p \otimes d) \in \sR \right\}.
\]
The naturality of this transformation is a consequence of the distributivity relation in $\sP \otimes \sQ$.
\end{itemize}

Conversely, suppose given a properad map $\nicexy{\sP \ar[r]^-{\phi} & \Hom(\sQ,\sR)}$.  The corresponding properad map $\nicexy{\sP \otimes \sQ \ar[r]^-{\phi^\flat} & \sR}$ is defined as follows:
\begin{itemize}
\item
For a pair of colors $(c,d) \in \fC \times \fD$, define the color
\[
\phi^\flat_0(c,d) = \left[\phi_0(c)\right]_0 (d) \in \fE.
\]
Here $\phi_0(c)$ is a color of $\Hom(\sQ,\sR)$, i.e., a properad map $\sQ \to \sR$.  So it has a map $ \left[\phi_0(c)\right]_0 \colon \fD \to \fE$ on color sets.
\item
For a generator $p \otimes d \in \sP \wedge \sQ$, define
\[
\phi^\flat(p \otimes d)
=
\phi(p)_d \in \sR,
\]
where $\phi(p) \in \Hom(\sQ,\sR)$, i.e., a natural transformation.
\item
For a generator $c \otimes q \in \sP \wedge \sQ$, define
\[
\phi^\flat(c \otimes q)
=
\left[ \phi_0(c)\right] (q) \in \sR.
\]
\end{itemize}
It is now a simple exercise to check that the above constructions define mutually inverse maps.
\end{proof}

In fact, the bijection \eqref{gpropclosed} is part of an  enriched isomorphism.

\begin{corollary}
\label{cor:enrichediso}
Suppose $\sP$, $\sQ$, and $\sR$ are special properads.  Then there is a natural isomorphism
\[
\Hom\left(\sP \otimes \sQ, \sR\right) 
\cong 
\Hom\left(\sP, \Hom(\sQ,\sR)\right)
\]
of special properads.
\end{corollary}

\begin{proof}
We need to check several things.
\begin{enumerate}
\item
The bijection between color sets, i.e.,
\[
\properads(\sP \otimes \sQ, \sR) 
\cong 
\properads(\sP, \Hom(\sQ,\sR)),
\]
is \eqref{gpropclosed}.
\item
That the underlying sets of the two internal hom objects are in bijection with each other follows from a simple adjoint exercise.  Indeed, suppose given an element $\varphi \in \Hom(\sP \otimes \sQ, \sR)$, i.e., a natural transformation
\[
\left\{\fC \times \fD \ni (c,d) \longmapsto \varphi_{(c,d)} \in \sR\right\}.
\]
The corresponding element $\varphi^\sharp \in \Hom(\sP, \Hom(\sQ,\sR))$ is the natural transformation that sends $c \in \fC$ to the natural transformation
\[
\left\{\fD \ni d \longmapsto \varphi_{(c,d)} \in \sR\right\}.
\]
The map going backward is similar.
\item
Finally, the two internal hom properad structures are also isomorphic. Indeed, as defined in \eqref{gammainternalhom} and Lemma \ref{lem:propginternalhom}, the properad structure of the internal hom object is induced by that of the target.  Therefore, the two internal hom objects under consideration both have their properad structures induced by that of $\sR$.
\end{enumerate}
\end{proof}

\begin{remark}
\begin{enumerate}
\item
The symmetric monoidal closed structure on operads is \cite{weiss} (Theorem 2.22).
\item
The symmetric monoidal closed structure on PROPs is \cite{hr1} (Theorem 35).  In this setting, the enriched isomorphism in Corollary \ref{cor:enrichediso} is \cite{hr1} (Proposition 34).
\end{enumerate}
\end{remark}


\chapter{Graphical Properads}
\label{ch:grproperad}

\abstract*{We define graphical properads as the free properads generated by connected wheel-free graphs.  We observe that a graphical properad has an infinite set of elements precisely when the generating graph is not simply connected.  The discussion of the tensor product of free properads in chapter  \ref{ch:tensor}  applies in particular to graphical properads.  Then we illustrate with several examples that a general properad map between graphical properads may exhibit bad behavior that would never happen when working with simply connected graphs.  These examples serve as the motivation of the restriction on the morphisms in the graphical category.}

In this chapter, we define and study the objects, called graphical properads, in the graphical category for connected wheel-free graphs.  The graphical category is defined in chapter \ref{ch:mapgrproperad}.  Graphical properads are freely generated by connected wheel-free graphs.  The construction of these graphical properads is analogous to the way linear graphs generate the finite ordinal category $\varDelta$.  In fact, the graphical category $\varGamma$ for connected wheel-free graphs contains a full subcategory isomorphic to $\varDelta$, and another full subcategory equivalent to the Moerdijk-Weiss dendroidal category $\varOmega$.

Although graphical properads have similar construction as the objects in finite ordinal category and the dendroidal category, they have very different properties.  In section \ref{sec:graphicalproperad}, after defining a graphical properad, it is observed that the graphical properad $\varGamma(G)$ is a finite set \emph{precisely} when $G$ is simply connected.  Thus, most graphical properads are infinite sets.  This is very different from the finite ordinal category and the dendroidal category, in which every object is finite.

In section \ref{sec:tensorgrproperad} we discuss the tensor product of two graphical properads, using the results in section \ref{sec:productfreeproperad}.  We then consider an explicit example of the tensor product of two graphical properads.

In section \ref{sec:mapsfromgprop} we discuss some important properties that does \emph{not} extend from the finite ordinal category and the dendroidal category to general properad maps between graphical properads.  A general graphical properad is quite large, so a general properad map between graphical properads can have strange behavior that would never happen in the dendroidal category or the finite ordinal category.  Therefore, unlike the finite ordinal category or the dendroidal category, when we work with all  connected wheel-free graphs, we should not take the full subcategory generated by the graphical properads.  As we will discuss in chapter \ref{ch:mapgrproperad}, to have a well-behaved graphical category for connected wheel-free graphs, we will have to impose certain restrictions on the maps.

In section \ref{sec:sctarget} it is observed that the bad behavior of general properad maps between graphical properads in section \ref{sec:mapsfromgprop} does not happen when the target is simply connected.  In other words, for graphical properads generated by simply connected graphs, general properad maps among them have good properties, as in the finite ordinal category and the dendroidal category.

As in the previous chapter, the underlying category here is $\Set$.  With suitable restrictions (i.e., using $\gupci$ or $\gupco$ instead of $\gupc$), everything in sections \ref{sec:graphicalproperad} and \ref{sec:tensorgrproperad}  has obvious analogs for properads with non-empty inputs or non-empty outputs.

\section{Properads Generated by Connected Wheel-Free Graphs}
\label{sec:graphicalproperad}

In this section, we define free properads generated by connected wheel-free graphs.  They will form the objects in the graphical category for connected wheel-free graphs.

\subsection{\texorpdfstring{$1$}{1}-Colored Graphs}

\begin{convention}
Unless otherwise specified, we use the one-colored set $\{*\}$ to define the set $\gupc$ of connected wheel-free graphs, i.e., $\gupc = \gupc(\{*\})$.  We will call them \textbf{$1$-colored connected wheel-free graphs}, or just \textbf{$1$-colored graphs}.\index{graph!$1$-colored}
\end{convention}

One could think of a $1$-colored graph as one in which every flag has the same color.  However, eventually we will give every edge a separate color, so it is psychologically easier to think of the flags as not having any colors to begin with.  The current objective is to define the free properad generated by a $1$-colored graph $G \in \gupc$.  The set of colors used in the free properad construction will be the set of edges of $G$ defined below.

The intuitive idea is that the edges of  $G \in \gupc$ can all be given different colors, so we can color each edge using its own name.  Each vertex is then a generating operation with inputs and outputs parametrized by the incoming and outgoing flags.  This defines a colored object, where the set of colors is the set of edges.  Then we take the free properads generated by this colored object.

\begin{remark}
Suppose $G \in \gupc$ is a $1$-colored graph.
\begin{enumerate}
\item
The set $\edge(G)$ of edges is finite, since a graph can only have finitely many flags to begin with.  It is empty precisely when $G$ consists of a single isolated vertex.
\item
An edge that is actually an ordinary leg consists of one flag.  On the other hand, an edge that is actually an ordinary edge consists of two flags.
\item
Using the set $\edge(G)$ as colors, each vertex $v \in \vertex(G)$ determines a corresponding pair of $\edge(G)$-profiles
\[\label{note:profv}
\profilev \in \sS(\edge(G)) 
= 
\catp(\edge(G))^{op} \times \catp(\edge(G))
\]
In other words, the incoming flags of $v$ determines a profile (i.e., finite ordered sequence) of edges, each of which is either an input of $G$ or an ordinary edge of $G$, and similarly for the outgoing flags of $v$.  Moreover, for two distinct vertices in $G$, their corresponding pairs of $\edge(G)$-profiles are different because $G$ cannot have two isolated vertices.  For an isolated vertex, the pair of $\edge(G)$-profiles is $\emptyprof$.
\end{enumerate}
\end{remark}

\begin{remark}
If $G$ is a $\fC$-colored connected graph, then we may also regard it as a $1$-colored graph by forgetting about the coloring, i.e., by reassigning to every edge the same color $*$.  This is true even if $\fC = \varnothing$ because this implies $\edge(G) = \varnothing$.  The only connected graph with an empty set of edges is the isolated vertex $\bullet$.  But then we can also regard $\bullet$ as a $1$-colored graph, where there is no edge to color.  In what follows, wherever necessary this switch from a $\fC$-colored graph to its associated $1$-colored graph will be done automatically without further comment.
\end{remark}

\subsection{Generating Set}

Next we define the colored object associated to a $1$-colored graph, consisting of essentially its vertices with the input/output $\edge(G)$-profiles. It will serve as the generating set for the associated graphical properad.

\begin{definition}
\label{def:clessobject}
Suppose $G \in \gupc$.  Define an $\sS(\edge(G))$-colored object $\ghat$ \index{graphical properad!generating set} as follows:
\[\label{note:ghat}
\ghat\dc = 
\begin{cases}
\{v\} & \text{ if $\dc = \profilev$ for $\ving$},\\
\varnothing & \text{ otherwise}.
\end{cases}
\]
\end{definition}

\begin{remark}
Each component of the $\sS(\edge(G))$-colored object $\ghat$ is either empty or consists of a single element.  In particular, $\ghat$ has only finitely many elements.
\end{remark}

\begin{example}
Consider the $1$-colored graph $G \in \gupc$
\begin{center}
\begin{tikzpicture}
\matrix[row sep=1cm,column sep=1.5cm] {
\node [plain] (p1) {$v$}; \\
\node [plain] (p2) {$u$}; \\
};
\draw [arrow,bend left=50] (p2) to node{$e$} (p1);
\draw [arrow,bend right=50] (p2) to node[swap]{$f$} (p1);
\end{tikzpicture}
\end{center}
with profiles $\emptyprofh$, vertices $\{u,v\}$, and $\edge(G) = \{e,f\}$.  Then 
\[
\ghat\binom{e,f}{\varnothing} = \{u\}, \quad
\ghat\binom{\varnothing}{e,f} = \{v\},
\]
and all other components of $\ghat$ are empty.
\end{example}

\begin{example}
Consider the $1$-colored graph $K$
\begin{center}
\begin{tikzpicture}
\matrix[row sep=.7cm, column sep = 2cm]{
\node [plain] (w) {$w$};&\\
& \node [plain] (v) {$v$};\\
\node [plain] (u) {$u$}; & \\
};
\draw [arrow] (u) to node[swap]{$f$} (v);
\draw [arrow] (u) to node{$e$} (w);
\draw [arrow,bend left=30] (v) to node{$g$} (w);
\draw [arrow,bend right=30] (v) to node[swap]{$h$} (w);
\draw [inputleg] (w) to node[below left=.1cm]{$i_1$} +(-.6cm,-.5cm);
\draw [inputleg] (u) to node[below left=.1cm]{$i_2$} +(-.6cm,-.5cm);
\draw [inputleg] (u) to node[below right=.1cm]{$i_3$} +(.6cm,-.5cm);
\draw [inputleg] (v) to node[below right=.1cm]{$i_4$} +(.6cm,-.5cm);
\draw [outputleg] (v) to node[above right=.1cm]{$o$} +(.6cm,.5cm);
\end{tikzpicture}
\end{center}
with profiles $\binom{o}{i_1,i_2,i_3,i_3}$, three vertices $\{u,v,w\}$, and 
\[
\edge(K) = \{i_1, i_2, i_3, i_4, o, e, f, g, h\}.
\]
Then 
\[
\khat\binom{e,f}{i_2,i_3} = \{u\}, \quad
\khat\binom{g,h,o} {f,i_4}= \{v\}, \quad
\khat\binom{\varnothing}{i_1,e,g,h} = \{w\},
\]
and all other components of $\khat$ are empty.
\end{example}

Recall the underlying functor $F = F_{\gupc}$ of the free properad monad \eqref{fgx}.

\begin{definition}
\label{def:freegprop}
Suppose $G \in \gupc$.  Define the free $\edge(G)$-colored properad
\[\label{note:gammag}
\varGamma(G) = F(\ghat),
\]
called the \textbf{graphical properad} \index{graphical properad} generated by $G$.
\end{definition}

\begin{remark}
To explain what $\varGamma(G)$ means intuitively, recall from \eqref{fgx} that the free properad $F\sP$ of an $\SC$-colored object $\sP$ consists of $\sP$-decorated $\fC$-colored connected wheel-free graphs, whose properad structure maps $\gamma_G$ are given by graph substitution.  For a vertex $v$ in $G$, the corresponding element in $\ghat$ can only decorate a vertex with equal profiles in an $\edge(G)$-colored connected wheel-free graph.  So an element in $\varGamma(G)$ is an $\edge(G)$-colored $\ghat$-decored connected wheel-free graph.
\end{remark}

\begin{notation}
If there is no danger of confusion, we sometimes abbreviate:
\begin{itemize}
\item
a graphical properad $\varGamma(G)$ as $G$, and
\item
a morphism $\varGamma(G) \to \varGamma(H) \in \properad$ as $G \to H$.
\end{itemize}
\end{notation}

\subsection{Size of Graphical Properads}
\label{sec:nonfiniteness}

Next we characterize the finiteness of a graphical properad $\varGamma(G)$ using a graph theoretic property of $G$.  The upshot is that there is a big jump in the size of $\varGamma(G)$ once $G$ is \emph{not} simply connected.

\begin{theorem}
\label{omegainfinite}
Suppose $G \in \gupc$.  Then the graphical properad $\varGamma(G)$ is a finite set if and only if $G$ is simply connected.\index{graphical properad!finiteness}
\end{theorem}

Theorem \ref{omegainfinite} is a combination of two observations, the first of which requires the following preliminary fact.

\begin{lemma}
\label{vertexappearsonce}
Suppose $G \in \gupc$ is simply connected.  Then each vertex $v$ in $G$ can appear in each element in $\varGamma(G)$ at most once.
\end{lemma}

\begin{proof}
Suppose to the contrary that there exists a vertex $v$ in $G$ that appears in some $K \in \varGamma(G)$ at least twice.  Since $K$ is connected, there is an internal path $P$ in $K$ from one copy of $v$, say $v_1$, to another copy of $v$, say $v_2$.  Since $G$ has no loops, $P$ must contain another vertex $u \not= v$ of $G$.  We may visualize $P$ as follows.
\begin{center}
\begin{tikzpicture}
\matrix[row sep=1cm, column sep=1.5cm]{
\node [plain] (v1) {$v_1$}; &
\node [plain] (u) {$u$}; &
\node [plain] (v2) {$v_2$};\\
};
\draw [dashedline] (v1) to (u);
\draw [dashedline] (u) to (v2);
\end{tikzpicture}
\end{center}
The two ordinary edges of $P$ adjacent to $u$ must be \emph{distinct} ordinary edges of $G$.
Taking an internal sub-path of $P$ if necessary, it follows that $G$ has a cycle containing $u$, contradicting the simple connectivity assumption on $G$.  Therefore, each $v \in \vertex(G)$ can appear at most once in each $\edge(G)$-colored $\ghat$-decorated graph $K \in \varGamma(G)$.
\end{proof}

The next observation should be of some interest on its own.

\begin{lemma}
\label{lem:gammasconn}
Suppose $G \in \gupc$ is simply connected.  Then each element in the graphical properad $\varGamma(G)$, regarded as an $\edge(G)$-colored $\ghat$-decorated connected wheel-free graph, is also simply connected.
\end{lemma}

\begin{proof}
Suppose $K \in \varGamma(G)$ and $u,v \in \vertex(K)$ are distinct vertices.  We must show that there exists a unique internal path in $K$ with end vertices $u$ and $v$.  Since $u$ and $v$ are distinct vertices in $K$, they are also distinct vertices in $G$ by Lemma \ref{vertexappearsonce}.  Moreover, since $K$ is connected, there is an internal path $P$ in $K$ with end vertices $u$ and $v$.   To see that $P$ is unique, note that the vertices in $P$ are distinct vertices in $K$ by definition, so they are also distinct vertices in $G$ by Lemma \ref{vertexappearsonce} again.  It follows that $P$ is also an internal path in $G$ with end vertices $u$ and $v$.  Since $G$ is simply connected, such an internal path is uniquely determined by its end vertices.
\end{proof}

The next observation gives half of Theorem \ref{omegainfinite}.

\begin{lemma}
\label{lem:omegagfinite}
Suppose $G \in \gupc$ is simply connected.  Then the graphical properad $\varGamma(G)$ is a finite set.
\end{lemma}

\begin{proof}
If $G$ is an exceptional edge $\uparrow$, then $\varGamma(G)$ consists of only the unit colored by $\uparrow$.  So suppose $G$ has at least one vertex.  Each $K \in \varGamma(G)$ is either an exceptional edge colored by some edge in $G$, or else it is, up to listing, uniquely determined by its set of vertices $\vertex(K) \subseteq \vertex(G)$ by Lemma \ref{lem:gammasconn}.  Since $\edge(G)$ and $\vertex(G)$ are both finite, it follows that $\varGamma(G)$ is finite as well.
\end{proof}

The next observation is the other half of Theorem \ref{omegainfinite}.

\begin{lemma}
\label{lem:omegaginfinite}
Suppose $G \in \gupc$ is not simply connected.  Then the graphical properad $\varGamma(G)$ is an infinite set.
\end{lemma}

\begin{proof}
We will exhibit an infinite sequence of distinct elements in $\varGamma(G)$.  Since $G$ is not simply connected, it contains a cycle. So there exist
\begin{itemize}
\item
two distinct vertices $u$ and $v$ in $G$, and 
\item
two internal paths $P$ and $Q$ with end vertices $u$ and $v$ that do \emph{not} have any common ordinary edges or common vertices besides $u$ and $v$.
\end{itemize}
We represent this situation as in the following picture.
\begin{center}
\begin{tikzpicture}
\matrix[row sep=1cm, column sep=2cm]{
\node [plain] (u) {$u$}; &
\node [plain] (v) {$v$};\\
};
\draw [dashedline,bend left=50] (u) to node{$P$} (v);
\draw [dashedline,bend right=50] (u) to node[swap]{$Q$} (v);
\end{tikzpicture}
\end{center}
Each of $u$ and $v$ may have flags that are not part of $P$ and $Q$, which for simplicity we are not drawing in the above picture.  For each integer $n \geq 1$, consider the following $\edge(G)$-colored $\ghat$-decorated graph $K_n$.
\begin{center}
\begin{tikzpicture}
\matrix[row sep=1cm, column sep=2.5cm]{
\node [plain] (v1) {$v_1$}; &
\node [plain] (v2) {$v_2$}; &
\node [empty] (v3) {$\cdots$}; &
\node [plain] (v4) {$v_n$}; \\
\node [plain] (u1) {$u_1$}; &
\node [plain] (u2) {$u_2$}; &
\node [empty] (u3) {$\cdots$}; &
\node [plain] (u4) {$u_n$}; \\
};
\foreach \x in {1,2,4}
\draw [dashedline] (u\x) to node[swap]{$P$} (v\x);
\draw [dashedline] (v1) to node[near start]{$Q$} (u2);
\draw [dashedline] (v2) to node[near start]{$Q$} (u3);
\draw [dashedline] (v3) to node[near start]{$Q$} (u4);
\draw [thick,dashed] (v4) to [out=-115, in=-30] node[swap]{$Q$} (u1);
\end{tikzpicture}
\end{center}
Each $v_i$ (resp., $u_i$) is a copy of $v$ (resp., $u$).  There are $n$ copies of $P$, each represented by a vertical dashed line.  There are $n$ copies of $Q$, including the slanted dashed line and the smooth dashed curve from $v_n$ to $u_1$.  Since the $K_n$ are all different elements in the free properad $\varGamma(G)$, the latter is infinite.
\end{proof}

Lemmas \ref{lem:omegagfinite} and \ref{lem:omegaginfinite} now combine to yield Theorem \ref{omegainfinite}.

\begin{example}
Consider the connected wheel-free graph $G$
\begin{center}
\begin{tikzpicture}
\matrix[row sep=.3cm,column sep=1.5cm]{
\node [plain] (v) {$v$}; & \\
& \node [plain] (w) {$w$}; \\
\node [plain] (u) {$u$}; & \\
};
\draw [arrow] (u) to node{$e$} (v);
\draw [arrow] (w) to node[swap]{$g$} (v);
\draw [arrow] (u) to node[swap]{$f$} (w);
\end{tikzpicture}
\end{center}
with vertices $\{u,v,w\}$, ordinary edges $\{e,f,g\}$, and profiles $\emptyprofh$.  Using the construction in the proof of Lemma \ref{lem:omegaginfinite}, for each $n \geq 1$ there is an $\edge(G)$-colored $\ghat$-decorated graph
\begin{center}
\begin{tikzpicture}
\matrix[row sep=.6cm,column sep=1.5cm]{
\node [plain] (v1) {$v_1$}; && 
\node [plain] (v2) {$v_2$}; && 
\node [plain] (v3) {$v_n$};  \\
& \node [plain] (w1) {$w_1$}; && 
\node [empty] (w2) {$\cdots$}; & \\
\node [plain] (u1) {$u_1$}; && 
\node [plain] (u2) {$u_2$}; && 
\node [plain] (u3) {$u_n$}; \\
&& \node [plain] (wn) {$w_n$}; && \\
};
\draw [arrow] (u1) to node{$e$} (v1);
\draw [arrow] (w1) to node[swap]{$g$} (v1);
\draw [arrow] (u2) to node[swap,near end]{$f$} (w1);
\draw [arrow] (u2) to node{$e$} (v2);
\draw [arrow] (w2) to node[swap]{$g$} (v2);
\draw [arrow] (u3) to node[swap,near end]{$f$} (w2);
\draw [arrow] (u3) to node[swap]{$e$} (v3);
\draw [arrow, bend right=15] (u1) to node[swap]{$f$} (wn);
\draw [->,thick] (wn) to [out=0, in=-110] node[swap]{$g$} (v3);
\end{tikzpicture}
\end{center}
with profiles $\emptyprofh$.  There are $n$ copies of each of $u$, $v$, and $w$.  This gives an infinite list of distinct elements in $\varGamma(G)$.

We can also twist the above construction as follows.  There is an $\edge(G)$-colored $\ghat$-decorated graph
\begin{center}
\begin{tikzpicture}
\matrix[row sep=.6cm,column sep=1.5cm]{
\node [plain] (v1) {$v_1$}; && 
\node [plain] (v2) {$v_2$}; && 
\node [plain] (v3) {$v_n$};  \\
& \node [plain] (w1) {$w_1$}; && 
\node [empty] (w2) {$\cdots$}; & \\
\node [plain] (u1) {$u_1$}; && 
\node [plain] (u2) {$u_2$}; && 
\node [plain] (u3) {$u_n$}; \\
};
\draw [outputleg] (u1) to node[swap]{$f$} +(.6cm,.4cm);
\draw [arrow] (u1) to node{$e$} (v1);
\draw [arrow] (w1) to node[swap]{$g$} (v1);
\draw [arrow] (u2) to node[swap,near end]{$f$} (w1);
\draw [arrow] (u2) to node{$e$} (v2);
\draw [arrow] (w2) to node[swap]{$g$} (v2);
\draw [arrow] (u3) to node[swap,near end]{$f$} (w2);
\draw [arrow] (u3) to node{$e$} (v3);
\draw [inputleg] (v3) to node{$g$} +(.6cm,-.4cm);
\end{tikzpicture}
\end{center}
with profiles $(g;f)$, $n$ copies of $v$, $n$ copies of $u$, and $n-1$ copies of $w$. This gives another infinite family of distinct elements in $\varGamma(G)$.
\end{example}

\section{Symmetric Monoidal Product}
\label{sec:tensorgrproperad}

Recall that $\properad$ is a symmetric monoidal category (Theorem \ref{thm:propgmonoidal}).  In this section, we consider the tensor product of two graphical properads along with one elaborate example.

\subsection{Symmetric Monoidal Product of Graphical Properads}

In Theorem \ref{distsimplified} it was shown that the tensor product of two free properads, $F(\phat)$ and $F(\qhat)$ with $\phat$ and $\qhat$ special, is the quotient of the free properad $F(\phat \wedge \qhat)$ by the generating distributivity relations.  Moreover, in Corollary \ref{tensorfpresented} it was observed that, when the special generating sets $\phat$ and $\qhat$ are finite on finitely many colors, the tensor product $F(\phat) \otimes F(\qhat)$ is finitely presented.

A graphical properad $\varGamma(G)$ with $G \in \gupc$ is the free properad $F(\ghat)$.  The generating set $\ghat$ is finite, whose elements are the vertices in $G$.  Its color set is the edge set $\edge(G)$.  Therefore, Theorem \ref{distsimplified} and  Corollary \ref{tensorfpresented} both apply to graphical properads generated by special graphs (Definition \ref{def:wheelfreegraphs}).

\begin{corollary}
\label{graphicalproptensor}
Suppose $G,H \in \gupcs$.  Then there is an isomorphism
\[
\varGamma(G) \otimes \varGamma(H) 
\cong 
\frac{F(\ghat \wedge \hhat)}{\text{generating distributivity}}
\]
of properads.\index{tensor product!of graphical properads}  In particular, this tensor product is finitely presented.
\end{corollary}

\begin{remark}
In the context of Corollary \ref{graphicalproptensor}, using the original definition of the tensor product, the generating set is $\varGamma(G) \wedge \varGamma(H)$.  Suppose either $G$ or $H$ is \emph{not} simply connected.  Then this generating set is infinite by Theorem \ref{omegainfinite}.  There are infinitely many distributivity relations, one for each element in $\varGamma(G) \times \varGamma(H)$, with a few exceptions.  The previous corollary says that the original description of the tensor product for graphical properads is redundant.  The tensor product of two graphical properads is always finitely generated with finitely many relations.
\end{remark}

\begin{example}
Corollary \ref{mapfromtensor} also applies to any two graphical properads $\varGamma(G) = F(\ghat)$ and $\varGamma(H) = F(\hhat)$ with $G,H \in \gupcs$.  In this case, the colored object $\ghat \wedge \hhat$ has
\begin{itemize}
\item
color set $\edge(G) \times \edge(H)$, and
\item
element set the coproduct
\[
\left[\vertex(G) \times \edge(H)\right] 
\coprod 
\left[\edge(G) \times \vertex(H)\right]. 
\]
\end{itemize}
The equality in Corollary \ref{mapfromtensor} corresponds to generating distributivity.  In particular, there are only finitely many such equalities to check because both $\vertex(G)$ and $\vertex(H)$ are finite.
\end{example}

The rest of this section contains an explicit example of the tensor product of two graphical properads.

\subsection{The Connected Wheel-Free Graphs}

In the rest of this section, we consider the following two connected wheel-free graphs.  
\begin{enumerate}
\item 
Suppose $G$ is the connected wheel-free graph:
\bigskip
\begin{center}
\begin{tikzpicture}
\matrix[row sep=.7cm,column sep=2cm] {
\node [plain] (w) {$w$}; & \\
& \node [plain] (v) {$v$}; \\ 
\node [plain] (u) {$u$}; & \\
};
\draw [arrow] (u) to node{$e_6$} (w);
\draw [arrow] (u) to node[swap]{$e_7$} (v);
\draw [arrow,bend left=15] (v) to node{$e_8$} (w);
\draw [arrow,bend right=15] (v) to node[swap]{$e_9$} (w);
\draw [inputleg] (u) to node[below=.1cm]{$e_1$} +(0cm,-.7cm);
\draw [inputleg] (v) to node[below right]{$e_2$} +(.6cm,-.4cm);
\draw [outputleg] (v) to node[above right]{$e_5$} +(.6cm,.4cm);
\draw [outputleg] (w) to node[above left]{$e_3$} +(-.6cm,.4cm);
\draw [outputleg] (w) to node[above right]{$e_4$} +(.6cm,.4cm);
\end{tikzpicture}
\end{center}
\bigskip
There are three vertices $\{u,v,w\}$, two inputs $\{e_1,e_2\}$, three outputs $\{e_3,e_4,e_5\}$, and four internal edges $\{e_6,e_7,e_8,e_9\}$.  The set of edges is
\[
\edge(G) = \{e_1,\ldots,e_9\},
\]
and $\ghat$ is the $\edge(G)$-colored set with three elements
\[
u \in \ghat\binom{e_6,e_7}{e_1},\quad
v \in \ghat\binom{e_8,e_9,e_5}{e_7,e_2}, \andspace
w \in \ghat\binom{e_3,e_4}{e_6,e_8,e_9}.
\]
As before, we will draw these elements as decorated corollas:
\bigskip
\begin{center}
\begin{tikzpicture}
\matrix[row sep=2cm,column sep=3cm] {
\node [plain] (u) {$u$}; &
\node [plain] (v) {$v$}; & 
\node [plain] (w) {$w$}; \\ 
};
\draw [inputleg] (u) to node[below=.1cm]{$e_1$} +(0,-.6cm);
\draw [outputleg] (u) to node[above left]{$e_6$} +(-.6cm,.4cm);
\draw [outputleg] (u) to node[above right]{$e_7$} +(.6cm,.4cm);
\draw [inputleg] (v) to node[below left]{$e_7$} +(-.6cm,-.4cm);
\draw [inputleg] (v) to node[below right]{$e_2$} +(.6cm,-.4cm);
\draw [outputleg] (v) to node[above left]{$e_8$} +(-.6cm,.4cm);
\draw [outputleg] (v) to node[above=.1cm]{$e_9$} +(0cm,.6cm);
\draw [outputleg] (v) to node[above right]{$e_5$} +(.6cm,.4cm);
\draw [inputleg] (w) to node[below left]{$e_6$} +(-.6cm,-.4cm);
\draw [inputleg] (w) to node[below=.1cm]{$e_8$} +(0cm,-.6cm);
\draw [inputleg] (w) to node[below right]{$e_9$} +(.6cm,-.4cm);
\draw [outputleg] (w) to node[above left]{$e_3$} +(-.6cm,.4cm);
\draw [outputleg] (w) to node[above right]{$e_4$} +(.6cm,.4cm);
\end{tikzpicture}
\end{center}
\bigskip
Note that $G$ is \emph{not} simply connected.  In particular, the graphical properad $\varGamma(G)$ is an infinite set by Theorem \ref{omegainfinite}.
\item
Suppose $H$ is the connected wheel-free graph:
\bigskip
\begin{center}
\begin{tikzpicture}
\matrix[row sep=2cm,column sep=2cm] {
\node [plain] (y) {$y$};  \\
\node [plain] (x) {$x$}; \\ 
};
\draw [arrow,bend left=35] (x) to node{$f_6$} (y);
\draw [arrow,bend right=35] (x) to node[swap]{$f_7$} (y);
\draw [inputleg] (x) to node[below left]{$f_2$} +(-.6cm,-.4cm);
\draw [inputleg] (x) to node[below right]{$f_3$} +(.6cm,-.4cm);
\draw [outputleg] (x) to node[above right]{$f_5$} +(.6cm,.4cm);
\draw [inputleg] (y) to node[above left]{$f_1$} +(-.6cm,-.4cm);
\draw [outputleg] (y) to node[above=.1cm]{$f_4$} +(0cm,.7cm);
\end{tikzpicture}
\end{center}
\bigskip
There are two vertices $\{x,y\}$, three inputs $\{f_1,f_2,f_3\}$, two outputs $\{f_4,f_5\}$, and two internal edges $\{f_6,f_7\}$.  The set of edges is
\[
\edge(H) = \{f_1,\ldots,f_7\},
\]
and $\hhat$ is the $\edge(H)$-colored set with two elements
\[
x \in \hhat\binom{f_6,f_7,f_5}{f_2,f_3} \andspace
y \in \hhat\binom{f_4}{f_1,f_6,f_7}.
\]
They will be drawn as the decorated corollas:
\bigskip
\begin{center}
\begin{tikzpicture}
\matrix[row sep=2cm,column sep=4cm] {
\node [plain] (x) {$x$}; &
\node [plain] (y) {$y$}; \\ 
};
\draw [inputleg] (x) to node[below left]{$f_2$} +(-.6cm,-.4cm);
\draw [inputleg] (x) to node[below right]{$f_3$} +(.6cm,-.4cm);
\draw [outputleg] (x) to node[above left]{$f_6$} +(-.6cm,.4cm);
\draw [outputleg] (x) to node[above=.1cm]{$f_7$} +(0cm,.6cm);
\draw [outputleg] (x) to node[above right]{$f_5$} +(.6cm,.4cm);
\draw [inputleg] (y) to node[below left]{$f_1$} +(-.6cm,-.4cm);
\draw [inputleg] (y) to node[below=.1cm]{$f_6$} +(0cm,-.6cm);
\draw [inputleg] (y) to node[below right]{$f_7$} +(.6cm,-.4cm);
\draw [outputleg] (y) to node[above=.1cm]{$f_4$} +(0cm,.6cm);
\end{tikzpicture}
\end{center}
\bigskip
Note that $H$ is also \emph{not} simply connected, so the graphical properad $\varGamma(H)$ is infinite by Theorem \ref{omegainfinite}.
\end{enumerate}

\subsection{The Smash Product}

Recall from Corollary \ref{graphicalproptensor}
that there is an isomorphism
\[
\varGamma(G) \otimes \varGamma(H) 
\cong 
\frac{F(\ghat \wedge \hhat)}{\text{generating distributivity}}
\]
of properads. To provide a description of the tensor product $\varGamma(G) \otimes \varGamma(H)$, let us first describe the smash product of the generating sets $\ghat$ and $\hhat$.
\begin{enumerate}
\item
The set of colors of $\ghat \wedge \hhat$ is the Cartesian product
\[
\edge(G) \times \edge(H) = 
\left\{(e_i,f_j) \colon 1 \leq i \leq 9,~ 1 \leq j \leq 7\right\}.
\]
\item
The element set of $\ghat \wedge \hhat$ is the coproduct
\[
\begin{split}
\ghat \wedge \hhat 
&= \left[\ghat \times \edge(H)\right] 
\coprod 
\left[\edge(G) \times \hhat\right]\\
&= \left\{u \otimes f_j, v \otimes f_j, w \otimes f_j\right\}_{1\leq j \leq 7} \coprod 
\left\{e_i \otimes x, e_i \otimes y\right\}_{1 \leq i \leq 9}.
\end{split}
\]
\end{enumerate}
In particular, the colored object $\ghat \wedge \hhat$ has $16$ elements on $63$ colors.

\subsection{Generating Distributivity}

The generating distributivity relations correspond to the elements in the Cartesian product $\ghat \times \hhat$, so there are $6$ generating distributivity relations.  We will list them explicitly.  For this purpose, we will use the same conventions for drawing decorated graphs as in the proof of Theorem \ref{distsimplified}.  In particular, input and output relabeling permutations will not be drawn.
\begin{enumerate}
\item
The generating distributivity relation for $u \in \ghat$ and $x \in \hhat$ identifies the following two decorated graphs, i.e., elements in $F(\ghat \wedge \hhat)$:
\bigskip
\begin{center}
\begin{tikzpicture}
\matrix[row sep=1.5cm,column sep=1.5cm] {
\node [plain] (u1) {$u$}; 
& \node [plain] (u2) {$u$}; 
& \node [plain] (u3) {$u$};
& \node [plain] (x2) {$x$}; 
& \node [plain] (x3) {$x$};\\
& \node [plain] (x1) {$x$}; 
&& \node [plain] (u4) {$u$};
& \node [plain] (u5) {$u$};
\\ 
};
\foreach \x in {1,2,3}
{
 \draw [outputleg] (u\x) to +(-.6cm,.4cm);
 \draw [outputleg] (u\x) to +(.6cm,.4cm);
}
\draw [inputleg] (x1) to +(-.6cm,-.4cm);
\draw [inputleg] (x1) to +(.6cm,-.4cm);
\foreach \x in {1,2,3}
 \draw [arrow] (x1) to (u\x);
\foreach \x in {2,3}
{
 \draw [outputleg] (x\x) to +(-.6cm,.4cm);
 \draw [outputleg] (x\x) to +(0,.6cm);
 \draw [outputleg] (x\x) to +(.6cm,.4cm);
}
\foreach \x in {4,5}
 \draw [inputleg] (u\x) to +(0,-.6cm);
\foreach \x in {4,5}
 \foreach \y in {2,3}
 \draw [arrow] (u\x) to (x\y);
\end{tikzpicture}
\end{center}
\bigskip
For the decorated graph on the left, the three copies of $u$ from left to right are $u \otimes f_6$, $u \otimes f_7$, and $u \otimes f_5$.  The bottom copy of $x$ is $e_1 \otimes x$. The input and output profiles are
\[
\binom{(e_6,f_6), (e_7,f_6), (e_6,f_7), (e_7,f_7), (e_6,f_5), (e_7,f_5)}{(e_1,f_2), (e_1,f_3)}.
\]
The three internal edges from left to right have colors $(e_1,f_6)$, $(e_1,f_7)$, and $(e_1,f_5)$.

For the decorated graph on the right, the two copies of $x$ from left to right are $e_6 \otimes x$ and $e_7 \otimes x$.  The two copies of $u$ from left to right are $u \otimes f_2$ and $u \otimes f_3$.  The input and output profiles are
\[
\binom{(e_6,f_6), (e_6,f_7), (e_6,f_5), (e_7,f_6), (e_7,f_7), (e_7,f_5)}{(e_1,f_2), (e_1,f_3)}
\]
The two internal edges from the left copy of $u$ have colors $(e_6,f_2)$ and $(e_7,f_2)$.  The two internal edges from the right copy of $u$ have colors $(e_6,f_3)$ and $(e_7,f_3)$. The decorated graph on the right requires an output relabeling.  Similar descriptions will apply to the decorated graphs in the next five generating distributivity relations, and we will therefore omit them.
\item
The generating distributivity relation for $u \in \ghat$ and $y \in \hhat$ identifies the following two decorated graphs:
\bigskip
\begin{center}
\begin{tikzpicture}
\matrix[row sep=1.5cm,column sep=1.5cm] {
\node [plain] (u1) {$u$}; 
&& \node [plain] (y2) {$y$}; 
&& \node [plain] (y3) {$y$};\\
\node [plain] (y1) {$y$}; 
&& \node [plain] (u2) {$u$};
& \node [plain] (u3) {$u$}; 
& \node [plain] (u4) {$u$};
\\ 
};
\draw [outputleg] (u1) to +(-.6cm,.4cm);
\draw [outputleg] (u1) to +(.6cm,.4cm);
\draw [inputleg] (y1) to +(-.6cm,-.4cm);
\draw [inputleg] (y1) to +(0,-.6cm);
\draw [inputleg] (y1) to +(.6cm,-.4cm);
\draw [arrow] (y1) to (u1);
\foreach \x in {2,3}
 \draw [outputleg] (y\x) to +(0,.6cm);
\foreach \x in {2,3,4}
 \draw [inputleg] (u\x) to +(0,-.6cm);
\foreach \x in {2,3,4}
\foreach \y in {2,3}
 \draw [arrow] (u\x) to (y\y);
\end{tikzpicture}
\end{center}
\bigskip
Neither input nor output relabeling is needed in this case.
\item
The generating distributivity relation for $v \in \ghat$ and $x \in \hhat$ identifies the following two decorated graphs:
\bigskip
\begin{center}
\begin{tikzpicture}
\matrix[row sep=1.5cm,column sep=1cm] {
\node [plain] (v1) {$v$}; 
& \node [plain] (v2) {$v$}; 
& \node [plain] (v3) {$v$};
& \node [plain] (x3) {$x$}; 
& \node [plain] (x4) {$x$}; 
& \node [plain] (x5) {$x$};
\\
\node [plain] (x1) {$x$}; 
&& \node [plain] (x2) {$x$};
& \node [plain] (v4) {$v$}; 
&& \node [plain] (v5) {$v$};
\\ 
};
\foreach \x in {1,2,3}
{
 \draw [outputleg] (v\x) to +(-.5cm,.4cm);
 \draw [outputleg] (v\x) to +(0,.6cm);
 \draw [outputleg] (v\x) to +(.5cm,.4cm);
}
\foreach \y in {3,4,5}
{
 \draw [outputleg] (x\y) to +(-.5cm,.4cm);
 \draw [outputleg] (x\y) to +(0,.6cm);
 \draw [outputleg] (x\y) to +(.5cm,.4cm);
}
\foreach \x in {1,2}
{
 \draw [inputleg] (x\x) to +(-.5cm,-.4cm);
 \draw [inputleg] (x\x) to +(.5cm,-.4cm);
}
\foreach \y in {4,5}
{
 \draw [inputleg] (v\y) to +(-.5cm,-.4cm);
 \draw [inputleg] (v\y) to +(.5cm,-.4cm);
}
\foreach \x in {1,2}
 \foreach \y in {1,2,3}
 \draw [arrow] (x\x) to (v\y);
\foreach \x in {4,5}
 \foreach \y in {3,4,5}
 \draw [arrow] (v\x) to (x\y);
\end{tikzpicture}
\end{center}
\bigskip
Both input and output relabeling are needed on the right-hand side.
\item
The generating distributivity relation for $v \in \ghat$ and $y \in \hhat$ identifies the following two decorated graphs:
\bigskip
\begin{center}
\begin{tikzpicture}
\matrix[row sep=1.5cm,column sep=1cm] {
& \node [plain] (v1) {$v$}; 
&& \node [plain] (y3) {$y$}; 
& \node [plain] (y4) {$y$}; 
& \node [plain] (y5) {$y$};
\\
\node [plain] (y1) {$y$}; 
&& \node [plain] (y2) {$y$};
& \node [plain] (v2) {$v$}; 
& \node [plain] (v3) {$v$};
& \node [plain] (v4) {$v$};
\\ 
};
\draw [outputleg] (v1) to +(-.5cm,.4cm);
\draw [outputleg] (v1) to +(0,.6cm);
\draw [outputleg] (v1) to +(.5cm,.4cm);
\foreach \x in {1,2}
{
 \draw [inputleg] (y\x) to +(-.5cm,-.4cm);
 \draw [inputleg] (y\x) to +(0,-.6cm);
 \draw [inputleg] (y\x) to +(.5cm,-.4cm);
}
\foreach \x in {1,2}
 \draw [arrow] (y\x) to (v1);
\foreach \x in {3,4,5}
 \draw [outputleg] (y\x) to +(0,.6cm);
\foreach \x in {2,3,4}
{
 \draw [inputleg] (v\x) to +(-.5cm,-.4cm);
 \draw [inputleg] (v\x) to +(.5cm,-.4cm);
}
\foreach \x in {2,3,4}
 \foreach \y in {3,4,5}
 \draw [arrow] (v\x) to (y\y);
\end{tikzpicture}
\end{center}
\bigskip
An input relabeling is needed on the right-hand side.
\item
The generating distributivity relation for $w \in \ghat$ and $x \in \hhat$ identifies the following two decorated graphs:
\bigskip
\begin{center}
\begin{tikzpicture}
\matrix[row sep=1.5cm,column sep=1.5cm] {
\node [plain] (w1) {$w$}; 
& \node [plain] (w2) {$w$};
& \node [plain] (w3) {$w$};
& \node [plain] (x4) {$x$}; 
& \node [plain] (x5) {$x$};
\\
\node [plain] (x1) {$x$}; 
& \node [plain] (x2) {$x$};
& \node [plain] (x3) {$x$};
& \node [plain] (w4) {$w$}; 
& \node [plain] (w5) {$w$};
\\ 
};
\foreach \x in {1,2,3}
{
 \draw [outputleg] (w\x) to +(-.5cm,.4cm);
 \draw [outputleg] (w\x) to +(.5cm,.4cm);
}
\foreach \x in {1,2,3}
{
 \draw [inputleg] (x\x) to +(-.5cm,-.4cm);
 \draw [inputleg] (x\x) to +(.5cm,-.4cm);
}
\foreach \x in {1,2,3}
 \foreach \y in {1,2,3}
 \draw [arrow] (x\x) to (w\y);
\foreach \x in {4,5}
{
 \draw [outputleg] (x\x) to +(-.5cm,.4cm);
 \draw [outputleg] (x\x) to +(0cm,.6cm);
 \draw [outputleg] (x\x) to +(.5cm,.4cm);
}
\foreach \x in {4,5}
{
 \draw [inputleg] (w\x) to +(-.5cm,-.4cm);
 \draw [inputleg] (w\x) to +(0,-.6cm);
 \draw [inputleg] (w\x) to +(.5cm,-.4cm);
}
\foreach \x in {4,5}
 \foreach \y in {4,5}
 \draw [arrow] (w\x) to (x\y);
\end{tikzpicture}
\end{center}
\bigskip
Both input and output relabeling are needed on the right-hand side.
\item
The generating distributivity relation for $w \in \ghat$ and $y \in \hhat$ identifies the following two decorated graphs:
\bigskip
\begin{center}
\begin{tikzpicture}
\matrix[row sep=1.5cm,column sep=1cm] {
& \node [plain] (w1) {$w$}; 
&& \node [plain] (y4) {$y$}; 
&& \node [plain] (y5) {$y$};
\\
\node [plain] (y1) {$y$}; 
& \node [plain] (y2) {$y$};
& \node [plain] (y3) {$y$};
& \node [plain] (w2) {$w$};
& \node [plain] (w3) {$w$};
& \node [plain] (w4) {$w$}; 
\\ 
};
\draw [outputleg] (w1) to +(-.5cm,.4cm);
\draw [outputleg] (w1) to +(.5cm,.4cm);
\foreach \x in {1,2,3}
{
 \draw [inputleg] (y\x) to +(-.5cm,-.4cm);
 \draw [inputleg] (y\x) to +(0cm,-.6cm);
 \draw [inputleg] (y\x) to +(.5cm,-.4cm);
}
\foreach \x in {1,2,3}
 \draw [arrow] (y\x) to (w1);
\foreach \x in {4,5}
 \draw [outputleg] (y\x) to +(0cm,.6cm);
\foreach \x in {2,3,4}
{
 \draw [inputleg] (w\x) to +(-.5cm,-.4cm);
 \draw [inputleg] (w\x) to +(0,-.6cm);
 \draw [inputleg] (w\x) to +(.5cm,-.4cm);
}
\foreach \x in {2,3,4}
 \foreach \y in {4,5}
 \draw [arrow] (w\x) to (y\y);
\end{tikzpicture}
\end{center}
\bigskip
An input relabeling is needed on the right-hand side.
\end{enumerate}

\subsection{Generating Distributivity in Action}

Let us provide an explicit example of how generating distributivity applies to decorated graphs.  Consider the following element $K$ in $F(\ghat \wedge \hhat)$:
\bigskip
\begin{center}
\begin{tikzpicture}
\matrix[row sep=1.4cm,column sep=1.5cm] {
&& \node [plain] (w) {$w$}; &\\
& \node [plain] (v1) {$v$}; &&
\node [plain] (v2) {$v$};\\
\node [plain] (y1) {$y$};
& \node [plain] (y2) {$y$};
&& \node [plain] (u) {$u$};\\
&&& \node [plain] (y3) {$y$};\\
&&& \node [plain] (x) {$x$};\\
};
\foreach \x in {1,2}
{ 
 \draw [arrow] (v\x) to (w);
 \draw [arrow] (y\x) to (v1);
}
\draw [arrow] (u) to (v2);
\draw [arrow] (y3) to (u);
\draw [arrow,bend left=20] (x) to (y3);
\draw [arrow,bend right=20] (x) to (y3);
\foreach \x in {1,2}
{
 \draw [inputleg] (y\x) to +(-.6cm,-.4cm);
 \draw [inputleg] (y\x) to +(0,-.6cm);
 \draw [inputleg] (y\x) to +(.6cm,-.4cm);
}
\draw [inputleg] (w) to +(-.8cm,-.2cm);
\draw [inputleg] (y3)to +(-.6cm,-.4cm);
\draw [inputleg] (x) to +(-.6cm,-.4cm);
\draw [inputleg] (x) to +(.6cm,-.4cm);
\draw [inputleg] (v2) to +(.6cm,-.4cm);
\begin{scope}[outputleg]
\draw (v1) to +(.6cm,.4cm);
\draw (v1) to +(.7cm,.2cm);
\draw (w) to +(-.6cm,.4cm);
\draw (w) to +(.6cm,.4cm);
\draw (u) to +(-.6cm,.4cm);
\draw (v2) to +(-.6cm,.4cm);
\draw (v2) to +(.6cm,.4cm);
\draw (x) to +(.6cm,.4cm);
\end{scope}
\end{tikzpicture}
\end{center}
Let us describe the vertices in each row from top to bottom. The top $w$ is $w \otimes f_4$. On the next row, both copies of $v$ are $v \otimes f_4$. On the third row, the two copies of $y$ are $e_7 \otimes y$ and $e_2 \otimes y$, while the right-most $u$ is $u \otimes f_4$.  The remaining bottom right vertices are $e_1 \otimes y$ and $e_1 \otimes x$.  Observe that there are two copies of $v \otimes f_4$ and three copies of $y$, each paired with a different $e_i$.

Next we list the colors of the internal edges from top to bottom and from left to right. The two internal edges from the two copies of $v$ to the top $w$ have colors $(e_8,f_4)$ and $(e_9,f_4)$. On the next row, the two internal edges from $e_7 \otimes y$ and $e_2 \otimes y$ to the left copy of $v$ have colors $(e_7,f_4)$ and $(e_2,f_4)$.  The internal edge from $u$ to the right copy of $v$ has color $(e_7,f_4)$.  The internal edge from the right copy of $y$ to $u$ has color $(e_1,f_4)$.  The two internal edges out of the bottom right $x$ have colors $(e_1,f_6)$ and $(e_1,f_7)$.

The decorated graph $K$ has $11$ inputs and $8$ outputs.  More precisely, its input and output profiles are
\begin{footnotesize}
\[
\binom{(e_9, f_4), (e_5, f_4), (e_3, f_4), (e_4, f_4), (e_6, f_4), (e_8, f_4), (e_5, f_4), (e_1, f_5)}{(e_7, f_1), (e_7, f_6), (e_7, f_7), (e_2, f_1), (e_2, f_6), (e_2, f_7), (e_6, f_4), (e_1, f_1), (e_1, f_2), (e_1, f_3), (e_2, f_4)}.
\]
\end{footnotesize}
Similar descriptions apply to the next decorated graph, so we will once again omit them.

The generating distributivity relation involving $v$ and $y$ can be applied on the left side of $K$, while the generating distributivity relation involving $u$ and $y$ can be applied on the right side of $K$.  Therefore, as an element in the tensor product $\varGamma(G) \otimes \varGamma(H)$, the image of $K$ is equal to the following decorated graph:
\bigskip
\begin{center}
\begin{tikzpicture}
\begin{scope}
\matrix[row sep=1.5cm,column sep=1cm] {
&&& 
\node[plain] (w) {$w$}; &&&\\
\node[plain] (y1) {$y$}; &
\node[plain] (y2) {$y$}; & 
\node[plain] (y3) {$y$}; &&&& 
\node[plain] (v4) {$v$};\\
\node[plain] (v1) {$v$}; &
\node[plain] (v2) {$v$}; & 
\node[plain] (v3) {$v$}; && 
\node[plain] (y4) {$y$}; &&
\node[plain] (y5) {$y$};\\
&&&&
\node[plain] (u1) {$u$}; &
\node[plain] (u2) {$u$}; & 
\node[plain] (u3) {$u$};\\
&&&&&
\node[plain] (x) {$x$}; &\\
};
\end{scope}
\begin{scope}[arrow]
\draw (y1) to (w);
\draw (v4) to (w);
\foreach \x in {1,2,3}
 \foreach \y in {1,2,3}
 \draw (v\x) to (y\y);
\draw (y5) to (v4);
\foreach \x in {1,2,3}
 \foreach \y in {4,5}
 \draw (u\x) to (y\y);
\foreach \x in {2,3}
 \draw (x) to (u\x);
\end{scope}
\begin{scope}[inputleg]
\foreach \x in {1,2,3}
{ 
 \draw (v\x) to +(-.6cm,-.4cm);
 \draw (v\x) to +(.6cm,-.4cm);
}
\draw (w) to +(-.8cm,-.2cm);
\draw (u1) to +(0,-.6cm);
\draw (x) to +(-.6cm,-.4cm);
\draw (x) to +(.6cm,-.4cm);
\draw (v4) to +(.6cm,-.4cm);
\end{scope}
\begin{scope}[outputleg]
\foreach \x in {2,3}
 \draw (y\x) to +(0,.6cm);
\draw (w) to +(-.6cm,.4cm);
\draw (w) to +(.6cm,.4cm);
\draw (y4) to +(0,.6cm);
\draw (x) to +(.7cm,.3cm);
\draw (v4) to +(-.8cm,.2cm);
\draw (v4) to +(.8cm,.2cm);
\end{scope}
\end{tikzpicture}
\end{center}
Both input and output relabeling are needed for this decorated graph.  Note that here there are five copies of $y$, four copies of $v$, and three copies of $u$.

\section{Maps of Graphical Properads}
\label{sec:mapsfromgprop}

The graphical category $\varGamma$ for connected wheel-free graphs will be defined as a certain non-full subcategory of $\properad$ with objects the graphical properads.  In this section, we give some examples of maps between graphical properads to illustrate that, without suitable restrictions, such maps can have very wild behavior compared to the finite ordinal category and the dendroidal category.  These examples will serve as motivation for the restrictions on the maps in $\varGamma$ that will be defined in Definition \ref{def:graphicalmap}.

\subsection{Maps Between Graphical Properads}

First let us make it explicit what a properad map out of a graphical properad constitutes.

\begin{lemma}
\label{lem:mapfromgprop}
Suppose $G \in \gupc$, and $\sQ$ is a $\fD$-colored properad.  Then a map
\[
\nicexy{\varGamma(G) \ar[r]^-{f} & \sQ}
\]
of properads is equivalent to a pair of functions:
\begin{enumerate}
\item
A function $\nicexy{\edge(G) \ar[r]^-{f_0} & \fD.}$
\item
A function $f_1$ that assigns to each vertex $v \in \vertex(G)$ an element $f_1(v) \in \sQ\binom{f_0 \out (v)}{f_0 \inp(v)}$.\index{graphical properad!map from}
\end{enumerate}
\end{lemma}

\begin{proof}
This is a special case of Lemma \ref{lem:mapfromfree}, since $\varGamma(G) = F(\ghat)$ has color set $\edge(G)$, and its elements are the vertices in  $G$ with their $\sS(\edge(G))$-profiles.
\end{proof}

If we apply Lemma \ref{lem:mapfromgprop} when $\sQ$ is a graphical properad as well, then we obtain the following observation.

\begin{lemma}
\label{lem:mapbtwgprop}
Suppose $G,H \in \gupc$.  Then a map
\[
\nicexy{\varGamma(G) \ar[r]^-{f} & \varGamma(H)}
\]
of properads is equivalent to a pair of functions:
\begin{enumerate}
\item
A function $\nicexy{\edge(G) \ar[r]^-{f_0} & \edge(H).}$
\item
A function $f_1$ that assigns to each vertex $v \in \vertex(G)$ an $\edge(H)$-colored $\hhat$-decorated connected wheel-free graph
\[
f_1(v) \in \varGamma(H)\binom{f_0 \out (v)}{f_0 \inp(v)}.
\]
\end{enumerate} 
\end{lemma}

\begin{example}
\label{ex:mapfromarrow}
Suppose $G =~ \uparrow$, which has no vertices, and $\sQ$ is a $\fD$-colored properad.  Then a map 
\[
\nicexy{\varGamma(\uparrow) \ar[r]^-{f} & \sQ}
\]
of properads is equivalent to a choice of a color $d \in \fD$, since $\edge(G) = \{\uparrow\}$.
\end{example}

\begin{example}
Suppose there is a graph substitution decomposition $K=H(G)$ in $\gupc$, in which:
\begin{itemize}
\item
$H \not=~ \uparrow$, 
\item
$G$ is substituted into a vertex $w$ in $H$, and 
\item
a corolla is substituted into every other vertex in $H$. 
\end{itemize}
As discussed in Notation \ref{notation:gh}, there are canonical injections
\[
\edge(G) \hookrightarrow \edge(K) \andspace
\vertex(G) \hookrightarrow \vertex(K).
\]
These injections induce a map $G \to K$, which simply sends each edge/vertex in $G$ to its corresponding image in $K$.
\end{example}

\begin{convention}
In what follows, we will often write down a properad map $\nicexy{\varGamma(G) \ar[r]^-{\varphi} & \varGamma(H)}$ for specific graphs $G$ and $H$ by drawing the graphs and specifying the functions $\varphi_0$ and $\varphi_1$.
\end{convention}

\subsection{Non-Determination by Edge Sets}

An important fact in the finite ordinal category $\varDelta$ and the Moerdijk-Weiss dendroidal category $\varOmega$ is the  following statement.
\begin{quotation}
Every map $\nicexy{S \ar[r]^-{f} & T}$ in $\varDelta$ and $\varOmega$ is uniquely determined by what it does on edges, i.e., the function between color sets $\nicexy{\edge(S) \ar[r]^-{f_0} & \edge(T).}$
\end{quotation}
The reason this statement is true is that, for a simply connected graph $T$, each vertex in $T$ can appear in each element in $\varGamma(T)$ at most once by Lemma \ref{vertexappearsonce}, and that in a simply connected graph any two vertices are connected by a unique internal path.  This fact is used, for example, in the proof of \cite{mt} Lemma 2.3.2 ($=$ epi-mono factorization in $\varOmega$).  The following example illustrates that a general properad map between graphical properads is not determined by its action on edge sets.

\begin{example}
\label{ex:phizero}
Consider the properad map $\nicexy{\varGamma(G) \ar[r]^-{\varphi} & \varGamma(H),}$
\begin{center}
\begin{tikzpicture}
\matrix[row sep=.5cm, column sep=1.7cm]{
&&& \node [plain] (w) {$w$};\\
\node [plain] (u) {$u$}; 
& \node [empty] (s) {}; 
& \node [empty] (t) {}; &\\
&&& \node [plain] (v) {$v$};\\
};
\draw [inputleg] (u) to node[below left=.1cm]{$i_1$} +(-.7cm,-.5cm);
\draw [inputleg] (u) to node[below right=.1cm]{$i_2$} +(.7cm,-.5cm);
\draw [outputleg] (u) to node[above left=.1cm]{$o_1$} +(-.7cm,.5cm);
\draw [outputleg] (u) to node[above right=.1cm]{$o_2$} +(.7cm,.5cm);
\draw [arrow, bend left=50] (v) to node{$e$} (w);
\draw [arrow] (v) to node{$f$} (w);
\draw [arrow, bend right=50] (v) to node[swap]{$g$} (w);
\draw [arrow] (s) to node{$\varphi$} (t);
\end{tikzpicture}
\end{center}
defined as follows:
\begin{itemize}
\item
$\varphi_0(i_1) = \varphi_0(o_1) = e$, and $\varphi_0(i_2) = \varphi_0(o_2) = g$.
\item
$\varphi_1(u)$ is the following $\edge(H)$-colored $\hhat$-decorated connected wheel-free graph.
\begin{center}
\begin{tikzpicture}
\matrix[row sep=2cm, column sep=1cm]{
\node [plain] (w) {$w$};\\
\node [plain] (v) {$v$};\\
};
\draw [arrow] (v) to node{$f$} (w);
\draw [inputleg] (w) to node[below left=.1cm]{$e$} +(-.7cm,-.5cm);
\draw [inputleg] (w) to node[below right=.1cm]{$g$} + (.7cm,-.5cm);
\draw [outputleg] (v) to node[above left=.1cm]{$e$} +(-.7cm,.5cm);
\draw [outputleg] (v) to node[above right=.1cm]{$g$} +(.7cm,.5cm);
\end{tikzpicture}
\end{center}
This is well-defined because $\varphi_1(u)$ should be an element in $\varGamma(H)$ with profiles
\[
\binom{\varphi_0(o_1), \varphi_0(o_2)}{\varphi_0(i_1), \varphi_0(i_2)} 
=
\binom{e,g}{e,g}.
\]
\end{itemize}
There is another properad map $\nicexy{\varGamma(G) \ar[r]^-{\phi} & \varGamma(H)}$ defined as follows:
\begin{itemize}
\item
$\phi_0 = \varphi_0$.
\item
$\phi_1(u)$ is the following graph.
\begin{center}
\begin{tikzpicture}
\matrix[row sep=1cm, column sep=2cm]{
\node [plain] (w1) {$w$}; &\\
& \node [plain] (v1) {$v$};\\
& \node [plain] (w2) {$w$};\\
\node [plain] (v2) {$v$}; &\\
};
\draw [arrow] (v1) to node[near end]{$f$} (w1);
\draw [arrow] (v2) to node{$e$} (w1);
\draw [arrow] (v2) to node[near start]{$f$} (w2);
\draw [arrow, bend right=45] (v2) to node[swap]{$g$} (w2);
\draw [inputleg] (w1) to node[right=.2cm]{$g$} +(.8cm,0cm);
\draw [outputleg] (v1) to node[left=.2cm]{$e$} +(-.8cm,.0cm);
\draw [outputleg] (v1) to node[above right=.1cm]{$g$} +(.7cm,.4cm);
\draw [inputleg] (w2) to node[left=.2cm]{$e$} +(-.8cm,-0cm);
\end{tikzpicture}
\end{center}
\end{itemize}
This example shows that a general properad map between two graphical properads is \emph{not} determined by its function on edge sets.
\end{example}

\begin{example}
\label{ex:idonedge}
Here is an even more dramatic example that illustrates that a properad map between two graphical properads is \emph{not} determined by its function on edge sets.  Suppose $G$ is the following connected wheel-free graph.
\begin{center}
\begin{tikzpicture}
\matrix[row sep=1cm,column sep=1cm] {
\node [plain] (v) {$v$}; \\
\node [plain] (u) {$u$}; \\
};
\draw [arrow,bend left=40] (u) to node{$e$} (v);
\draw [arrow,bend right=40] (u) to node[swap]{$f$} (v);
\end{tikzpicture}
\end{center}
Define the properad map $\nicexy{\varGamma(G) \ar[r]^-{\varphi} & \varGamma(G)}$ as follows:
\begin{itemize}
\item
$\varphi_0 = \Id$ on $\edge(G)$.
\item
$\varphi_1(v) = C_v$, the corolla at $v$, i.e., the $\edge(G)$-colored $\ghat$-decorated corolla
\begin{center}
\begin{tikzpicture}
\matrix[row sep=1.5cm,column sep=1cm] {
\node [plain] (v) {$v$}; \\
};
\begin{scope}[inputleg]
\draw (v) to node[below left=.1cm]{$e$} +(-.6cm,-.5cm);
\draw (v) to node[below right]{$f$}+(.6cm,-.5cm);
\end{scope}
\end{tikzpicture}
\end{center}
with two inputs.
\item
$\varphi_1(u)$ is the following $\edge(G)$-colored $\ghat$-decorated connected wheel-free graph.
\begin{center}
\begin{tikzpicture}
\matrix[row sep=1.5cm,column sep=1.5cm] {
& \node [plain] (v) {$v$}; & \\
\node [plain] (u1) {$u$}; && \node [plain] (u2) {$u$}; \\
};
\draw [arrow,bend left=40] (u1) to node{$e$} (v);
\draw [arrow,bend right=40] (u2) to node[swap]{$f$} (v);
\draw [outputleg] (u1) to node[above right=.1cm]{$f$} +(.6cm,.5cm);
\draw [outputleg] (u2) to node[above left=.1cm]{$e$} +(-.6cm,.5cm);
\end{tikzpicture}
\end{center}
This is well-defined because the vertex $u$ in $G$ has two outputs with colors $e$ and $f$, and no inputs.
\end{itemize}
Notice that $\varphi_0 = \Id$, but $\varphi$ is not the identity map on $\varGamma(G)$.  Moreover, $G$ is simplest connected wheel-free graph that is not simply connected.  So this sort of non-determination behavior can happen even for very simple connected wheel-free graphs.
\end{example}

\subsection{Non-Injections}

We need the following definition of a linear branch to discuss the next cosimplicial/dendroidal fact that does not extend to a general properad map between graphical properads.

\begin{definition}
\label{def:linearbranch}
Suppose $G$ is a connected wheel-free graph.
\begin{enumerate}
\item
A \textbf{lane} \index{lane} in $G$ is a pair
\[
P = \left(\left(e^j\right)_{j=1}^r, \left(v_i\right)_{i=1}^{r-1}\right)
\]
with $r \geq 2$ such that all the conditions of an internal path (Definition \ref{def:path}) not involving the end vertices are satisfied, except that $e^1$ and $e^r$ may be ordinary legs.
\item
A \textbf{branch} \index{branch} in $G$ is a lane such that every vertex $v_i$ has one incoming flag and one outgoing flag.
\item
A \textbf{linear branch} \index{linear branch} in $G$ is a branch in which
\begin{itemize}
\item
$e^1$ has terminal vertex $v_1$,
\item
$e^j$ has initial vertex $v_{j-1}$ and terminal vertex $v_j$ for $1 < j < r$, and
\item
$e^r$ has initial vertex $v_{r-1}$.
\end{itemize}
\end{enumerate}
\end{definition}

Another important fact in the finite ordinal category $\varDelta$ and the  dendroidal category $\varOmega$ is the following statement.
\begin{quotation}
Suppose a map $\nicexy{S \ar[r]^-{f} & T}$ in $\varDelta$ or $\varOmega$ satisfies $f_0(a) = f_0(b)$ for two distinct edges $a$ and $b$ in $S$.  Then
\begin{itemize}
\item
$a$ and $b$ belong to the same linear branch of $S$, and
\item
$f$ sends all the vertices between $a$ and $b$ to the $f_0(a)$-colored exceptional edge.
\end{itemize}
\end{quotation}

The above fact in the dendroidal category is used in the proof of \cite{mt} Lemma 2.3.2 ($=$ epi-mono factorization in $\varOmega$).  It   is a consequence of the fact that every map in $\varOmega$ is uniquely determined by the function on edges, which was discussed above.   This fact does \emph{not} extend to general properad maps between graphical properads.  In fact, it fails in several different ways, as we illustrate in the following three examples.

\begin{example}
\label{ex1:noninjection}
In this example, we exhibit a properad map between graphical properads such that there exist two edges with the same image that do \emph{not} belong to the same linear branch.  Consider the following properad map $\nicexy{H \ar[r]^-{\psi} & G}$,
\begin{center}
\begin{tikzpicture}
\matrix[row sep=.6cm,column sep=1.5cm] {
\node [plain] (v') {$v'$}; &&& \node [plain] (v) {$v$}; \\
& \node [empty] (s) {}; & \node [empty] (t) {}; &\\
\node [plain] (u') {$u'$}; &&& \node [plain] (u) {$u$}; \\
};
\draw [arrow,bend right=40] (u') to node[swap]{$f'$} (v');
\draw [inputleg] (v') to node[below left=.1cm]{$e_{-1}$} +(-.6cm,-.5cm);
\draw [outputleg] (u') to node[above left=.1cm]{$e_1$} +(-.6cm,.5cm);
\draw [arrow,bend left=40] (u) to node{$e$} (v);
\draw [arrow,bend right=40] (u) to node[swap]{$f$} (v);
\draw [arrow] (s) to node{$\psi$} (t);
\end{tikzpicture}
\end{center}
that is defined as follows.
\begin{itemize}
\item
$H$ has two vertices $\{u',v'\}$, one internal edge $f'$ from $u'$ to $v'$, one input leg $e_{-1}$ attached to $v'$, and one output leg $e_1$ attached to $u'$.
\item
$G$ is as in Example \ref{ex:idonedge}.
\item
$\psi_0(e_i) = e$ and $\psi_0(f') = f$.
\item
$\psi_1(u') = C_u$ (the corolla at $u$) and $\psi_1(v') = C_v$ (the corolla at $v$).
\end{itemize}
Now observe that $\psi_0(e_1) = e = \psi_0(e_{-1})$, but the edges $e_1$ and $e_{-1}$ do \emph{not} belong to the same linear branch in $H$.
\end{example}

\begin{example}
\label{ex2:noninjection}
From Example \ref{ex1:noninjection}, one might think that for the dendroidal fact under discussion to fail for graphical properads, the two offending edges should not lie on the same linear branch.  This is not the case.  In this example, we exhibit a properad map between graphical properads such that there exist two edges on the same linear branch with the same image such that a vertex in between is \emph{not} sent to the identity. 

Consider the following properad map $\nicexy{J \ar[r]^-{\phi} & G}$ in $\varGamma$,
\begin{center}
\begin{tikzpicture}
\matrix[row sep=.3cm,column sep=1.5cm] {
\node [empty] (o) {}; &&& \node [plain] (v) {$v$}; \\
\node [plain] (w) {$w$}; & \node [empty] (s) {}; & \node [empty] (t) {}; &\\
\node [empty] (i) {}; &&& \node [plain] (u) {$u$}; \\
};
\draw [arrow] (i) to node[near start]{$e_1$} (w);
\draw [arrow] (w) to node[near end]{$e_{-1}$} (o);
\draw [arrow,bend left=40] (u) to node{$e$} (v);
\draw [arrow,bend right=40] (u) to node[swap]{$f$} (v);
\draw [arrow] (s) to node{$\phi$} (t);
\end{tikzpicture}
\end{center}
that is defined as follows.
\begin{itemize}
\item
$J$ is the corolla with one input $e_1$ and one output $e_{-1}$.
\item
$G$ is as in Example \ref{ex:idonedge}.
\item
$\phi_0(e_i) = e$.
\item
$\phi_1(w)$ is the $\edge(G)$-colored $\ghat$-decorated graph
\begin{center}
\begin{tikzpicture}
\matrix[row sep=2cm,column sep=1cm] {
\node [plain] (v) {$v$}; \\
\node [plain] (u) {$u$}; \\
};
\draw [arrow,bend right=40] (u) to node[swap]{$f$} (v);
\draw [inputleg] (v) to node[below left=.1cm]{$e$} +(-.6cm,-.5cm);
\draw [outputleg] (u) to node[above left=.1cm]{$e$} +(-.6cm,.5cm);
\end{tikzpicture}
\end{center}
with two vertices $\{u,v\}$, one internal edge $f$ from $u$ to $v$, one $e$-colored input leg attached to $v$, and one $e$-colored output leg attached to $u$.
\end{itemize}
Now observe that $\phi_0(e_1) = e = \phi_0(e_{-1})$, that $e_1$ and $e_{-1}$ belong to the same linear branch in $J$, and that the vertex $w$ in $J$ is \emph{not} sent to the $e$-colored identity.
\end{example}

\begin{example}
\label{ex3:noninjection}
Examples \ref{ex1:noninjection} and \ref{ex2:noninjection} may lead the reader to think that for the dendroidal fact under discussion to fail for graphical properads, any two offending edges should not be both input (or both output) flags of the same vertex.  This is once again not the case.  In this example, we exhibit a properad map between graphical properads that sends two input flags of a vertex to the same image.  

Consider the following properad map $\nicexy{K \ar[r]^-{\theta} & L}$ in $\varGamma$,
\begin{center}
\begin{tikzpicture}
\matrix[row sep=.5cm,column sep=1.5cm] {
&&& \node [plain] (w) {$w$}; \\
\node [plain] (u) {$u$}; & \node [empty] (s) {}; & \node [empty] (t) {}; &\\
&&& \node [plain] (v) {$v$}; \\
&&& \node [empty] (i) {}; \\
};
\draw [inputleg] (u) to node[below left=.1cm]{$e_1$} +(-.6cm,-.5cm);
\draw [inputleg] (u) to node[below right=.1cm]{$e_2$} +(.6cm,-.5cm);
\draw [arrow,bend left=40] (v) to node{$f$} (w);
\draw [arrow,bend right=40] (v) to node[swap]{$g$} (w);
\draw [arrow] (i) to node[near start]{$e$} (v);
\draw [arrow] (s) to node{$\theta$} (t);
\end{tikzpicture}
\end{center}
that is defined as follows.
\begin{itemize}
\item
$K$ is the corolla with two inputs $e_1$ and $e_2$, and no outputs.
\item
$L$ has two vertices $\{v,w\}$, two internal edges $\{f,g\}$ from $v$ to $w$, one input leg $e$ attached to $v$, and no outputs.
\item
$\theta_0(e_i) = e$.
\item
$\theta_1(u)$ is the $\edge(L)$-colored $\lhat$-decorated graph
\begin{center}
\begin{tikzpicture}
\matrix[row sep=.6cm,column sep=1.5cm] {
& \node [plain] (w1) {$w$}; & \\
& \node [plain] (w2) {$w$}; & \\
\node [plain] (v1) {$v$}; & & \node [plain] (v2) {$v$}; \\
\node [empty] (i1) {}; && \node [empty] (i2) {}; \\
};
\draw [arrow,bend left=35] (v1) to node{$f$} (w1);
\draw [arrow] (v1) to node[swap]{$g$} (w2);
\draw [arrow,bend right=35] (v2) to node[swap]{$g$} (w1);
\draw [arrow] (v2) to node{$f$} (w2);
\foreach \x in {1,2}
\draw [arrow] (i\x) to node[near start]{$e$} (v\x);
\end{tikzpicture}
\end{center}
with two copies of $v$, two copies of $w$, two $e$-colored input legs, and no outputs.
\end{itemize}
Now observe that both $e_1$ and $e_2$ are incoming flags of the vertex $u$, and they both have image $e$.
\end{example}

\section{Maps with Simply Connected Targets}
\label{sec:sctarget}

In this section we make several observations concerning maps between graphical properads with simply connected targets.  The upshot is that the bad behavior of general properad maps between graphical properads, as exhibited in section \ref{sec:mapsfromgprop}, does \emph{not} occur when the graphs are simply connected.

\subsection{Simply Connected Targets}

The following observation says that a simply connected target guarantees a simply connected source.

\begin{proposition}
\label{sconnsource}
Suppose $\nicexy{G \ar[r]^-{\varphi} & H}$ is a properad map between graphical properads with $H$ simply connected.  Then $G$ is also simply connected.
\end{proposition}

\begin{proof}
Suppose to the contrary that $G$ is \emph{not} simply connected.  Since $G$ is connected wheel-free to begin with, this means that there exists a cycle $P$ (Definition \ref{def:path}) in $G$.  We may regard $P$ as an element in $\varGamma(G)$ by using the vertices and ordinary edges in $P$ as well as all the flags attached to these vertices.  By definition, the $\edge(H)$-colored $\hhat$-decorated graph $\varphi(P) \in \varGamma(H)$ is obtained as the following graph substitution:
\[
\varphi(P) = \left[\varphi_0(P)\right] \left(\{\varphi_1 (u)\}_{u \in \vertex(P)}\right).
\]
Since each $\varphi_1(u) \in \varGamma(H)$ is connected and since $P$ has a cycle, it follows that $\varphi(P) \in \varGamma(H)$ also has a cycle.  But this cannot happen by Lemma \ref{lem:gammasconn}.  Therefore, $G$ must be simply connected.
\end{proof}

\subsection{Injectivity on Inputs and Outputs}

In Example \ref{ex3:noninjection} we exhibited a properad map between graphical properads in which two input edges of a vertex are sent to the same image.  The next observation says that for this to happen, the target cannot be simply connected.

\begin{proposition}
\label{inoutinject}
Suppose $\nicexy{G \ar[r]^-{\varphi} & H}$ is a properad map between graphical properads with $H$ simply connected, and $u \in \vertex(G)$.  Then the restrictions of the function
\[
\nicexy{\edge(G) \ar[r]^-{\varphi_0} & \edge(H)}
\]
to $\inp(u)$ and to $\out(u)$ are both injective.
\end{proposition}

\begin{proof}
Suppose to the contrary that there are two distinct input edges $a,b \in \inp(u)$ with
\[
\varphi_0(a) = e = \varphi_0(b) \in \edge(H).
\]
Recall that $\varphi_1(u) \in \varGamma(H)$ has input/output profiles the $\varphi_0$-images of those of $u$.  Since $a$ and $b$ are input edges of $u$, the color $e$ appears at least twice in the input profile of $\varphi_1(u)$. In particular, $e$ is an input edge of a unique vertex $w \in \vertex(H)$, which appears in $\varphi_1(u)$ at least twice.  But this cannot happen by Lemma \ref{vertexappearsonce}.  The proof for the case $a,b \in \out(u)$ is nearly identical.
\end{proof}

\subsection{Non-Injections}

As discussed above, in the dendroidal category, when two edges are sent to the same image, the two edges must lie in the same linear branch (Definition \ref{def:linearbranch}) of the source.  Moreover, the map must send all the vertices between these two edges to identities.  The following observation is an extension of this dendroidal fact to graphical properads generated by simply connected graphs.

\begin{proposition}
\label{graphicalnoninject}
Suppose $\nicexy{G \ar[r]^-{\varphi} & H}$ is a properad map between graphical properads with $H$ simply connected.  Suppose $a$ and $b$ are distinct edges in $G$ with $\varphi_0(a) = \varphi_0(b)$. Then:
\begin{enumerate}
\item
There is a unique linear branch
\[
P = \left(\left(e^j\right)_{j=1}^r, \left(v_i\right)_{i=1}^{r-1}\right)
\]
in $G$ with $\{e^1,e^r\} = \{a,b\}$ as sets.
\item
The map $\varphi_1$ sends all the vertices in $P$ to the $\varphi_0(a)$-colored identity.
\end{enumerate}
\end{proposition}

\begin{proof}
By Proposition \ref{sconnsource} $G$ is simply connected.  If $H =~ \uparrow$, then $\varGamma(H)$ contains only the unit element.  This forces every vertex in $G$ to have one incoming flag and one outgoing flag, i.e., $G$ is a linear graph.  In this case, the two assertions are immediate.  So now we assume that $H \not=~ \uparrow$. In particular, every edge in $H$ is adjacent to at least one vertex.

By Proposition \ref{inoutinject}, there does \emph{not} exist $u\in\vertex(G)$ such that $a,b \in \inp(u)$ or $a,b \in \out(u)$.  So there exists a unique lane $P$ in $G$ (Definition \ref{def:linearbranch}) whose end edges are $a$ and $b$.  Regard $P$ as an element in $\varGamma(G)$ by using its ordinary edges and vertices as well as all the flags attached to these vertices in $G$.  Suppose
\[
\varphi_0(a) = \varphi_0(b) = e \in \edge(H).
\]
By definition $\varphi(P) \in \varGamma(H)$ is the graph substitution
\[
\varphi(P) = \left[\varphi_0 (P)\right] \left(\{\varphi_1 u\}_{u \in \vertex(P)}\right).
\]
The edges $a$ and $b$ \emph{cannot} be both input legs of $P$ because otherwise the terminal vertex of $e$ would appear in $\varphi(P) \in \varGamma(H)$ at least twice.  By the simple connectivity of $H$ and Lemma \ref{vertexappearsonce}, this cannot happen.  Likewise, $a$ and $b$ cannot be both output legs of $P$.  Switching names if necessary, we may assume without loss of generality that $a$ is an input leg of $P$, and $b$ is an output leg of $P$.

To prove the two stated assertions, it suffices to show that $\varphi(P)$ is the $e$-colored exceptional edge $\uparrow_e ~ \in \varGamma(H)$.  Indeed, the only way to get such an exceptional edge as a graph substitution of connected wheel-free graphs is to substitute exceptional edges into a linear graph.  To show that $\varphi(P)$ is the $e$-colored exceptional edge, it in turn suffices to show that $\varphi(P)$ contains no vertices.  Indeed, $\varphi(P) \in \varGamma(H)$ is simply connected by Lemma \ref{lem:gammasconn} and has the color $e$ among its input legs and also among its output legs.  The only simply connected graph with no vertices is the exceptional edge $\uparrow$.  So if $\varphi(P)$ has no vertices, then it must be the $e$-colored exceptional edge.

Since $\varphi(P)$ has $e \in \edge(H)$ among its input legs and among its output legs, $e$ must have both an initial vertex $v \in \vertex(H)$ and a terminal vertex $w \in \vertex(H)$, as in the following picture.
\begin{center}
\begin{tikzpicture}
\matrix[row sep=1.4cm,column sep=1.5cm] {
\node [plain,label=above:$...$] (v2) {$w$};\\
\node [plain,label=below:$...$] (v1) {$v$};\\
};
\draw [arrow] (v1) to node{$e$} (v2);
\foreach \x in {1,2}
{
\draw [inputleg] (v\x) to +(-.5cm,-.4cm);
\draw [inputleg] (v\x) to +(.5cm,-.4cm);
\draw [outputleg] (v\x) to +(-.5cm,.4cm);
\draw [outputleg] (v\x) to +(.5cm,.4cm);
}
\end{tikzpicture}
\end{center}
Suppose $\varphi(P)$ has at least one vertex.  Then the vertex $w$ appears in $\varphi(P)$ and contributes an $e$-colored input leg.  Likewise, the vertex $v$ appears in $\varphi(P)$ and contributes an $e$-colored output leg.  This is, however, impossible because in the simply connected graph $\varphi(P)$ (Lemma \ref{lem:gammasconn}), the only way to connect $w$ with $v$ is through the $e$-colored ordinary edge.  This means that $\varphi(P)$ contains no vertices.
\end{proof}

\subsection{Unique Determination by Edge Sets}

In Examples \ref{ex:phizero} and \ref{ex:idonedge}, we showed that a general properad map between two graphical properads is \emph{not} necessarily determined by its action on edges.  We now observe that such non-determination by edge sets can only happen when the target is not simply connected.  The following preliminary observation is needed.

\begin{lemma}
\label{onlyh}
Suppose $H \in \gupc$ is simply connected such that the component $\varGamma(H)\emptyprof$ is non-empty.  Then the following statements hold.
\begin{enumerate}
\item
$H$ has input/output profiles $\emptyprofh$.
\item
$\varGamma(H)\emptyprof$ is the singleton consisting of $H$.
\end{enumerate}
\end{lemma}

\begin{proof}
Suppose $K \in \varGamma(H)\emptyprof$.  We first claim that $K$ contains all the vertices in $H$.  If this is not true, then we can pick $w \in \vertex(H)$ that does \emph{not} appear in $K$ and $u \in \vertex(H)$ that \emph{does} appear in $K$.  Since $H$ is simply connected, there exists a unique internal path in $H$ with end vertices $u$ and $w$.  Choosing closer vertices if necessary, we may assume without loss of generality that there is an internal edge $e$ in $H$ that is adjacent to both $w$ and $u$.  On the other hand, since $K$ has empty input and output profiles, every vertex adjacent to $u$ in $H$ must also be in $K$.  This contradicts the existence of $w$.  Therefore, $K$ must contain all the vertices in $H$. By Lemma \ref{vertexappearsonce} each vertex in $H$ can appear in $K$ at most once, and hence exactly once by the previous sentence.  By the simple connectivity of $H$ and $K$ (Lemma \ref{lem:gammasconn}), we must have $K = H$.
\end{proof}

\begin{proposition}
\label{mapdetermined}
Suppose $\nicexy{G \ar[r]^-{\varphi} & H}$ is a properad map between graphical properads with $H$ simply connected.  Then $\varphi$ is uniquely determined by the function $\nicexy{\edge(G) \ar[r]^-{\varphi_0} & \edge(H)}$.
\end{proposition}

\begin{proof}
Suppose $u \in \vertex(G)$.  Then $\varphi_1(u)$ is an $\edge(H)$-colored $\hhat$-decorated graph with input/output profiles the $\varphi_0$-images of those of $u$, i.e.,
\[
\varphi_1(u) \in \varGamma(H)\binom{\varphi_0 \out(u)}{\varphi_0 \inp(u)}.
\]
The input (resp., output) profile $\varphi_0 \inp(u)$ (resp., $\varphi_0 \out(u)$) is a subset of $\edge(H)$.  If it is non-empty, then it corresponds to a \emph{unique} subset of $\vertex(H)$ by Lemma \ref{vertexappearsonce}.  In other words, the vertices in $\varphi_1(u)$ that contribute to the profiles $\varphi_0 \inp(u)$ and $\varphi_0 \out(u)$ are uniquely determined by the function $\varphi_0$.  But then the other vertices in $\varphi_1(u)$--namely, those that do \emph{not} contribute to either the input or the output profiles of $\varphi_1(u)$--are also determined because, in a simply connected graph, any two vertices are connected by a unique internal path.

The only remaining case is when both $\inp(u)$ and $\out(u)$ are emtpy, which can only happen if $G$ itself is a single isolated vertex $u$.  In this case, we have $\edge(G) = \varnothing$, and $\varphi_0$ is the trivial function. Since $\varphi_1(u) \in \varGamma(H)\emptyprof$, by Lemma \ref{onlyh} we have $\varphi_1(u) = H$.
\end{proof}



\chapter{Properadic Graphical Category}
\label{ch:mapgrproperad}

\abstract*{We define the properadic graphical category $\varGamma$, whose objects are graphical properads.  Its morphisms are called properadic graphical maps.  To define such graphical maps, we first discuss coface and codegeneracy maps between graphical properads.  We establish graphical analogs of the cosimplicial identities.  The most interesting case is the graphical analog of the cosimplicial identity
\[
d^j d^i = d^i d^{j-1}
\]
for $i<j$ because it involves iterating the operations of deleting an almost isolated vertex and of smashing two closest neighbors together.  Graphical maps do not have the bad behavior discussed in the examples in chapter \ref{ch:grproperad}.  In particular, it is observed that each graphical map has a factorization into codegeneracy maps followed by coface maps.  Such factorizations do not exist for general properad maps between graphical properads.
Finally, we show that the properadic graphical category admits the structure of a (dualizable) generalized Reedy category, in the sense of \cite{reedyextension}.}

The main purpose of this chapter is to define the maps in the graphical category $\varGamma$ for connected wheel-free graphs, whose objects are graphical properads (Definition \ref{def:freegprop}).

As the examples in section \ref{sec:mapsfromgprop} illustrate, in order to avoid bad behavior for maps between graphical properads, we should \emph{not} take the full subcategory of $\properad$ generated by graphical properads.  Some restrictions on the maps are necessary.  In particular, in each of the offending examples in section \ref{sec:mapsfromgprop}, the source is \emph{not} sent to a subgraph, which we will define precisely in this chapter, of the target.  In contrast, as a consequence of simple connectivity, every category (resp., operad) map between linear graphs (resp., unital trees) automatically sends the source to a subgraph of the target.  It turns out that this condition is sufficient to yield graphical analog of the epi-mono factorization and other good properties.  So we will use this condition--that the image of the source is a subgraph of the target--to define maps in the graphical category.

There are several equivalent ways to define the concept of a \emph{subgraph} of a connected wheel-free graph.  We will define it using outer coface maps and relabelings.  In section \ref{sec:coface} we define inner coface maps, outer coface maps, and codegeneracy maps between graphical properads.  In addition to the definition of a subgraph, coface maps will be used in section \ref{sec:infinityproperad} to define $\infty$-properads.  When restricted to linear graphs and unital trees, our coface and codegeneracy maps restrict to those in the finite ordinal category $\varDelta$ and the Moerdijk-Weiss dendroidal category $\varOmega$ \cite{mw1}.

In section \ref{sec:graphicalid} we show that the cosimplicial identities in $\varDelta$ have precise analogs for connected wheel-free graphs.  The graphical analog of the cosimplicial identity
\[
d^jd^i = d^id^{j-1} \quad \text{for $i<j$}
\]
is more involved than in either $\varDelta$ or $\varOmega$.  The reason is that graphical coface maps, just like in $\varDelta$ and $\varOmega$, are defined by either deleting an almost isolated vertex or by smashing two closest neighbors together.  These concepts were discussed in chapter \ref{ch:graph}.  For general connected wheel-free graphs, mixing these operations requires a lot of care.  This is why the graphical analog of the above cosimplicial identity needs more work.

In section \ref{sec:gupcgraphcat} we define relabelings, subgraphs, and maps in the graphical category $\varGamma$ for connected wheel-free graphs, called \emph{properadic graphical maps}.  All inner/outer coface maps, codegeneracy maps, subgraphs, and changes of listings are in $\varGamma$.  Subgraphs are characterized in terms of graph substitution decompositions.  A properadic graphical map is defined as a map between graphical properads such that the image is a subgraph of the target.  This condition is shown to be equivalent to the requirement that \emph{every} subgraph of the source is sent to a subgraph of the target.  Then we establish the graphical analog of the epi-mono factorization in $\varDelta$ and $\varOmega$ for properadic graphical maps.  We show that each properadic graphical map has a unique decomposition, up to isomorphism, into codegeneracy maps followed by coface maps.  Furthermore, a properadic graphical map is uniquely determined by its action on edge sets.

In section \ref{sec:reedygamma} we recall the definition of generalized Reedy category due to Berger and Moerdijk. This generalization is needed to deal with categories which have non-identity isomorphisms. We use results from previous sections to show that the graphical category $\varGamma$ and its opposite $\varGamma^{op}$ admit such generalized Reedy structures.

As in previous chapters, all the definitions and results about properads in this chapter have obvious analogs for properads with non-empty inputs or non-empty outputs.

\section{Coface and Codegeneracy Maps}
\label{sec:coface}

In this section, we define coface and codegeneracy maps between graphical properads, which will be used to define maps in the graphical category and the notion of an $\infty$-properad.  These maps extend the coface and codegeneracy maps in the finite ordinal category $\varDelta$ and the dendroidal category $\varOmega$.

\subsection{Inner Coface Maps}

To motivate our construction below, consider the finite ordinal category $\varDelta$.  In $\varDelta$ a coface map is \emph{inner} if it is neither the top coface map nor the bottom coface map.  In terms of linear graphs, an inner coface map $d^i \colon [n] \to [n+1]$ in $\varDelta$ ($0 < i < n$) corresponds to substituting the $1$-colored linear graph $L_2$, depicted as
\begin{center}
\begin{tikzpicture}
\matrix[row sep=.8cm, column sep=1cm]{
\node [empty] (v3) {};\\
\node [plain] (v2) {$v_{i+1}$};\\
\node [plain] (v1) {$v_i$};\\
\node [empty] (v0) {};\\
};
\draw [arrow] (v0) to node{\footnotesize{$i-1$}} (v1);
\draw [arrow] (v1) to node{\footnotesize{$i$}} (v2);
\draw [arrow] (v2) to node{\footnotesize{$i+1$}} (v3);
\end{tikzpicture}
\end{center}
into the vertex
\begin{center}
\begin{tikzpicture}
\matrix[row sep=.7cm, column sep=1cm]{
\node [empty] (v2) {$\vdots$};\\
\node [plain] (v1) {$v_i$};\\
\node [empty] (v0) {$\vdots$};\\
};
\draw [arrow] (v0) to node{\footnotesize{$i-1$}} (v1);
\draw [arrow] (v1) to node{\footnotesize{$i$}} (v2);
\end{tikzpicture}
\end{center}
of the $1$-colored linear graph $L_n$ with $n$ vertices.  In other words, the inner coface map $d^i \colon [n] \to [n+1]$ corresponds to the inner properadic factorization (Definition \ref{def:ipropfact})
\[
L_{n+1} = L_n\left(L_2\right)
\]
of $L_{n+1}$, with distinguished subgraph $L_2$ substituted into the vertex $v_i$.

Generalizing the setting of linear graphs, for connected wheel-free graphs an inner coface map corresponds to an inner properadic factorization, with the outer graph as the source.

\begin{definition}
\label{def:innercoface}
Suppose $G,K \in \gupc$.  An \textbf{inner coface map} \index{inner coface map!for graphical properad}
\[
\nicexy{
\varGamma(G) \ar[r]^-{d^w} & \varGamma(K)
}\]
is a properad map, corresponding to an inner properadic factorization $K = G(H_w)$ of $K$ (Definition \ref{def:ipropfact}),  defined by
\[
d^w_1(v) = 
\begin{cases}
H_w & \text{ if $w = v$},\\
C_v & \text{ otherwise},
\end{cases}
\]
where $C_v$ is the corolla with the same profiles as $v \in \vertex(G)$.  The map
\[
\nicexy{\edge(G) \ar[r]^-{d^w_0} & \edge(K)}
\]
sends
\[
\left(\inp(w); \out(w)\right) \longmapsto 
\left(\inp(H_w); \out(H_w)\right)
\]
and all other edges in $G$ to the corresponding ones in $K=G(H_w)$.
\end{definition}

\begin{remark}
Recall from Theorem \ref{thm:cnbdfact} that $K$ admits an inner properadic factorization $G(H_w)$ precisely when the vertices $u$ and $v$ in $H_w$ are closest neighbors in $K$.  In this case, $G$ is obtained from $K$ by smashing these closest neighbors together and deleting all the ordinary edges adjacent to both of them.  Therefore, we say that an inner coface map $\varGamma(G) \to \varGamma(K)$ as in Definition \ref{def:innercoface}  \textbf{corresponds to the closest neighbors} $u$ and $v$.
\end{remark}

\begin{example}
When restricted to linear graphs, an inner coface map as above is the same as one in $\varDelta$.  Indeed, in an inner properadic factorization
\[
K=G(H_w),
\]
the distinguished subgraph $H_w$ is assumed to be a partially grafted corollas.  The only partially grafted corollas in $\ULin$ is the linear graph $L_2$ with two vertices.  So Definition \ref{def:innercoface} restricts to the usual definition of an inner coface map in $\varDelta$.
\end{example}

\begin{example}
Similarly, when restricted to unital trees, Definition \ref{def:innercoface} restricts to the Moerdijk-Weiss inner coface map in the dendroidal category \cite{mw1}.  Note that those authors actually call these inner \emph{face} maps.
\end{example}

\begin{example}
\label{ex0:innercoface}
Suppose $G$ is a partially grafted corollas with profiles $\dch$.  Then there is an inner coface map
\[
C_{\dch} \to G
\]
corresponding to the inner properadic factorization $G=C_{\dch}(G)$.
\end{example}

\begin{example}
For the  graph $K$ in Example \ref{ex:cneighbor}, its only two inner properadic factorizations are in Examples \ref{ex:properadicfact} and \ref{ex2:properadicfact}.  Using the notation in those examples, the diagram
\[
\nicexy{
C_{\emptyprofh} \ar[r] \ar[d]
& G \ar[d]^-{d^y} \\
G' \ar[r]^-{d^z} & K
}\]
contains all the iterated inner coface maps ending at $K$.  Here $C_{\emptyprofh}$ is an isolated vertex, and the two unnamed inner coface maps are as in Example \ref{ex0:innercoface}.
\end{example}

As we will see later, inner coface maps can be used to factor maps in the properadic graphical category.  The following observation is a preliminary version of such a factorization.  We will use Notation \ref{notation:gh}.

\begin{lemma}
\label{innercofacefactor}
Suppose $K=G(H_w)$ in $\gupc$ with $H_w \not=~ \uparrow$.  Then the map $\nicexy{G \ar[r]^-{\varphi} & K}$ determined by
\[
\varphi_1(u) = 
\begin{cases}
C_u & \text{ if $u\not= w$},\\
H_w & \text{ if $u=w$},
\end{cases}
\]
has a decomposition into inner coface maps and isomorphisms induced by changes of vertex listings.
\end{lemma}

\begin{proof}
Write $H$ for $H_w$.  We prove the assertion by induction on $n = |\vertex(H)| > 0$.  If $n=1$, then $H$ is a permuted corolla, and $K$ is obtained from $G$ by changing the listing at the vertex $w$.  There is a canonical properad isomorphism
\[
\nicexy{G \ar[r]^-{\cong} & K}
\]
because there are canonical bijections
\[
\edge(G) \cong \edge(K) \andspace
\vertex(G) \cong \vertex(K)
\]  
induced by the change of listing at $w$.

Suppose $n>1$.  Then $H$ must have a pair of closest neighbors $u$ and $v$ (Example \ref{ex:cnbdexist}).  So there is a graph substitution decomposition
\[
H = J(P)
\]
in which $P$ is a partially grafted corollas with vertices $u$ and $v$ (Theorem \ref{thm:cnbdfact}).  We have $|\vertex(J)| = n-1$ because $P$ has two vertices.  Thus, by induction hypothesis the map
\[
\nicexy{G \to G(J)}
\]
decomposes into inner coface maps and isomorphisms induced by changes of vertex listings.  Now we have
\[
K = G[J(P)] = [G(J)](P),
\]
so the map $\varphi$ factors as follows.
\[
\nicexy{
G \ar[r]^-{\varphi} \ar[d] & K = [G(J)](P)\\
G(J) \ar[ur] & 
}\]  
The slanted map is by definition an inner coface map because $P$ is a partially grafted corollas.  So $\varphi$ is also a composition of inner coface maps and isomorphisms induced by changes of vertex listings.
\end{proof}

\subsection{Outer Coface Maps}

To motivate our definition below, we once again consider the finite ordinal category $\varDelta$.  A \emph{bottom} coface map $d^0 \colon [n] \to [n+1]$ in $\varDelta$ corresponds to deleting the bottom corolla $L_1$, depicted as
\begin{center}
\begin{tikzpicture}
\matrix[row sep=.7cm, column sep=1cm]{
\node [empty] (v2) {};\\
\node [plain] (v1) {$v_1$};\\
\node [empty] (v0) {};\\
};
\draw [arrow] (v0) to node{\footnotesize{$0$}} (v1);
\draw [arrow] (v1) to node{\footnotesize{$1$}} (v2);
\end{tikzpicture}
\end{center}
from the linear graph $L_{n+1}$ with $n+1$ vertices.  In other words, it corresponds to the outer properadic factorization  (Definition \ref{def:opropfact})
\[
L_{n+1} = L_2(L_n)
\]
of $L_{n+1}$, with distinguished subgraph $L_n$ substituted into the top vertex of $L_2$.

Likewise, a top coface map $d^{n+1} \colon [n] \to [n+1]$ in $\varDelta$ corresponds to deleting the top corolla $L_1$ in $L_{n+1}$.  This in turn corresponds to the outer properadic factorization
\[
L_{n+1} = L_2(L_n),
\]
in which the distinguished subgraph $L_n$ is substituted into the bottom vertex of $L_2$.  So the essence of a top or bottom coface map is about deleting a top or bottom vertex in a linear graph.

The analog of such a deletable vertex in a connected wheel-free graph is an almost isolated vertex (Definition \ref{def:almostisolated}).  Moreover, by Theorem \ref{thm:opropfact}, an almost isolated vertex corresponds to the vertex of the non-distinguished subgraph in an outer properadic factorization.  Generalizing the setting of linear graphs, for connected wheel-free graphs an outer coface map corresponds to an outer properadic factorization, with the distinguished subgraph as the source.

\begin{definition}
\label{def:outercoface}
Suppose $H_w, K \in \gupc$.  An \textbf{outer coface map} \index{outer coface map!for graphical properad}
\[
\nicexy{
\varGamma(H_w) \ar[r]^-{d^u} & \varGamma(K)
}\]
is a properad map, corresponding to an outer properadic factorization
\[
K = G(\{C_u,H_w\}) = G(H_w)
\]
of $K$ (Definition \ref{def:opropfact}), defined by
\[
d^u_1(v) = C_v
\]
for all $v \in \vertex(H_w)$, where $C_v$ denotes the corolla with the same profiles as $v$.  The map
\[
\nicexy{\edge(H_w) \ar[r]^-{d^u_0} & \edge(K)}
\]
sends every edge in $H_w$ to the corresponding one in $K=G(H_w)$.
\end{definition}

\begin{definition}
A \textbf{coface map} means either an inner coface map or an outer coface map.\index{coface map!for graphical properad}
\end{definition}

\begin{notation}
In later chapters, we sometimes denote a coface map with a subscript, e.g., $\nicexy{G \ar[r]^-{d_u} & H}$, instead of a superscript to improve typography.  In any case, the context should make it clear that we are talking about coface maps.
\end{notation}

\begin{remark}
Recall from Theorem \ref{thm:opropfact} that a vertex $u$ in $K$ is almost isolated precisely when there is an outer properadic factorization
\[
K=G(\{C_u,H_w\})
\]
in which $u$ is the non-distinguished vertex in the partially grafted corollas $G$. There is a canonical bijection
\[
\vertex(K) = \vertex(H_w) \coprod \{u\}.
\]  
Therefore, we say that an outer coface map $\varGamma(H_w) \to \varGamma(K)$ as in Definition \ref{def:outercoface} \textbf{corresponds to the almost isolated vertex $u$} in $K$.
\end{remark}

\begin{example}
As before, when restricted to linear graphs and unital trees, we recover an outer coface map in the finite ordinal category $\varDelta$ and the dendroidal category $\varOmega$.  Note that Moerdijk and Weiss \cite{mw1} called it an outer \emph{face} map.
\end{example}

\begin{example}
\label{ex:legcoface}
Suppose $C$ is a corolla but not an isolated vertex, and $e$ is one of its legs.  There is an outer properadic factorization
\[
C = G\left(\{C, \uparrow\}\right),
\]
in which $G$ is the dioperadic graph (Example \ref{ex:dioperadicgraph}) with one vertex having the same profiles as $C$ and the other vertex having one input and one output.  The latter is adjacent to the unique ordinary edge corresponding to $e$.  Therefore, there is an outer coface map
\[
\uparrow~ \to C,
\]
one for each leg of the corolla $C$.  This map simply sends the unique element on the left to the unit on the right whose color is the chosen leg of $C$.
\end{example}

\begin{example}
\label{ex:pgcorcoface}
Suppose $G$ is a partially grafted corollas (Example \ref{ex:pgcor}) with top vertex $v$ and bottom vertex $u$.  Using the corollas $C_u$ and $C_v$ with the profiles of $u$ and $v$, there is an outer properadic factorization
\[
G = G(\{C_v,C_u\}).
\]
Therefore, there are two outer coface maps
\[
C_u \to G 
\andspace
C_v \to G,
\]
corresponding to the almost isolated vertices $v$ and $u$, respectively.
\end{example}

\begin{example}
For the graph $K$ in Example \ref{ex:cneighbor}, its only almost isolated vertices are $u$ and $w$.  Using the notations in Examples \ref{ex:properadicfact} and \ref{ex2:properadicfact}, the diagram
\[
\nicexy{
\uparrow \ar@<.5ex>[r] \ar@<-.5ex>[r] \ar[d] \ar@<1ex>[d] \ar@<-1ex>[d] & 
H_y \ar[d]^{d^u} \\
H_z \ar[r]^-{d^w} & K
}\]
contains all the iterated outer coface maps ending at $K$.  The outer coface maps from the exceptional edge $\uparrow$ correspond to the legs of $H_y$ and $H_z$.
\end{example}

\subsection{Codegeneracy Maps}

Before we define a codegeneracy map between graphical properads, let us motivate the definition using the familiar finite ordinal category $\varDelta$.   A codegeneracy map $s^i \colon [n] \to [n-1]$ in $\varDelta$ corresponds to substituting the exceptional edge $\uparrow$ into the corolla
\begin{center}
\begin{tikzpicture}
\matrix[row sep=.7cm, column sep=1cm]{
\node [empty] (v2) {};\\
\node [plain] (v1) {$v_{i+1}$};\\
\node [empty] (v0) {};\\
};
\draw [arrow] (v0) to node{\footnotesize{$i$}} (v1);
\draw [arrow] (v1) to node{\footnotesize{$i+1$}} (v2);
\end{tikzpicture}
\end{center}
of the linear graph $L_n$ with $n$ vertices. In other words, a codegeneracy map in $\varDelta$ corresponds to the graph substitution decomposition
\[
L_{n-1} = L_n(\uparrow),
\]
in which the exceptional edge is substituted into a chosen vertex of $L_n$, and a corolla is substituted into every other vertex of $L_n$.

For graphical properads, the intuitive idea of a codegeneracy map $G \longrightarrow K$ is that $K$ is obtained from $G$ by substituting an exceptional edge into a chosen vertex with one input and one output.  To formalize this idea, we make the following definition.

\begin{definition}
\label{def:codegdecomp}
Suppose $G \in \gupc$, and $v \in \vertex(G)$ has exactly one incoming flag and one outgoing flag.  The \textbf{degenerate reduction} \index{degenerate reduction} of $G$ at $v$ is the graph
\[
G_v = G(\uparrow)
\]
obtained from $G$ by substituting the exceptional edge into $v$ and a corolla into every other vertex of $G$.  The edge in $G_v$ corresponding to the two edges adjacent to $v$ is denoted by $e_v$.
\end{definition}

\begin{remark}
Given a degenerate reduction $G_v = G(\uparrow)$, there are canonical bijections
\[
\begin{split}
\vertex(G)
&= \vertex(G_v) \coprod \{v\},\\
\edge(G) \setminus \{v_{\inp}, v_{\out}\}
&= \edge(G_v) \setminus \{e_v\},
\end{split}
\]
where $v_{\inp}$ and $v_{\out}$ are the incoming and outgoing flags of $v$. 
\end{remark}

\begin{remark}
The notation $G_v$ was also used in Definition \ref{def:almostisolated} when $v$ is almost isolated and in Definition \ref{def:deletable} when $v$ is deletable.  The context should make it clear what $G_v$ denotes.
\end{remark}

\begin{definition}
\label{def:codegenmap}
Suppose $G \in \gupc$.  A \textbf{codegeneracy map} \index{codegeneracy map}
\[
\nicexy{
\varGamma(G) \ar[r]^-{s^v} & \varGamma(G_v)}
\]
is a properad map, corresponding to a degenerate reduction $G_v=G(\uparrow)$ of $G$ at some vertex $v$ with one incoming flag and one outgoing flag, defined as follows:
\begin{itemize}
\item
For an edge $e \in \edge(G)$, set
\[
s^v_0(e) = 
\begin{cases}
e_v & \text{ if $e$ is adjacent to $v$},\\
e & \text{ otherwise}.
\end{cases}
\]
\item
For a vertex $u \in \vertex(G)$, set
\[
s^v_1(u) =
\begin{cases}
\uparrow_{e_v} & \text{ if $u=v$},\\
C_u & \text{ otherwise}.
\end{cases}
\]
\end{itemize}
\end{definition}

\begin{example}
As before, when restricted to linear graphs and unital trees, we recover a codegeneracy map in the finite ordinal category and the dendroidal category.  Note that Moerdijk and Weiss \cite{mw1} called it a \emph{degeneracy} map.
\end{example}

\begin{remark}
A codegeneracy map is much simpler than inner and outer coface maps because we do not need to worry about closest neighbors or almost isolated vertices.  All we need to define a codegeneracy map is a vertex with exactly one incoming flag and one outgoing flag.
\end{remark}

\section{Graphical Identities}
\label{sec:graphicalid}

The purpose of this section is to prove analogs of the cosimplicial identities for coface and codegeneracy maps between graphical properads, as defined in section \ref{sec:coface}.  The dendroidal analogs of the cosimplicial identities are discussed in \cite{mt} (section 2.2.3).  As usual, we recover both the cosimplicial and the dendroidal cases when we restrict to linear graphs and unital trees.

\subsection{Identities of Codegeneracy Maps}

First we consider a graphical analog of the cosimplicial identity
\[
s^j s^i = s^i s^{j+1} \quad \text{for $i \leq j$}.
\]
Recall from Definition \ref{def:codegdecomp} the notation $G_v = G(\uparrow)$ for a degenerate reduction of $G$ at a vertex $v$ with one incoming flag and one outgoing flag.

\begin{lemma}
\label{graphidss}
Suppose $u$ and $v$ are vertices in $G \in \gupc$, each having one incoming flag and one outgoing flag.  Suppose $G_{u,v}$ is obtained from $G$ by substituting an exceptional edge into each of $u$ and $v$, and a corolla into every other vertex.  Then there is a commutative diagram
\[
\nicexy{
G \ar[d]_{s^v} \ar[r]^-{s^u} &
G_u \ar[d]^-{s^v}\\
G_v \ar[r]^-{s^u} &
G_{u,v}
}\]
of codegeneracy maps.
\end{lemma}

\begin{proof}
This follows from the fact that
\[
\begin{split}
G_{u,v} 
&= G_u(\uparrow) = (G_u)_v\\
&= G_v(\uparrow) = (G_v)_u.
\end{split}
\]
In other words, $v$ is still a vertex in the degenerate reduction $G_u$ of $G$ at $u$, and substituting an exceptional edge into this $v$ yields
\[
(G_u)_v = G_{u,v}.
\]
A similar comment applies to $(G_v)_u$.
\end{proof}

\subsection{Identities Involving Both Codegeneracy and Coface Maps}

Next we consider graphical analogs of the cosimplicial identities
\[
s^j d^i = 
\begin{cases}
d^i s^{j-1} & \text{ if $i<j$},\\
\Id = s^j d^{j+1} & \text{ if $i=j$},\\
d^{i-1} s^j & \text{ if $i > j+1$}.
\end{cases}
\]
There are two cases.  First we consider an \emph{inner} coface map followed by a codegeneracy map.

\begin{lemma}
\label{graphidds}
Suppose $K=G(H_w)$ is an inner properadic factorization, and $v \in \vertex(K)$ has one incoming flag and one outgoing flag.
\begin{enumerate}
\item
If $v \not\in \vertex(H_w)$, then the diagram
\[
\nicexy{
G \ar[d]_{s^v} \ar[r]^-{d^w} 
& K \ar[d]^{s^v}\\
G_v \ar[r]^-{d^w} 
& K_v
}\]
is commutative.
\item
If $v \in \vertex(H_w)$, then the diagram
\[
\nicexy{
G \ar[dr]_{\Id} \ar[r]^-{d^w} 
& K \ar[d]^{s^v}\\
& G = K_v
}\]
is commutative.
\end{enumerate}
\end{lemma}

\begin{proof}
Here we use the decomposition
\[
\vertex(K) = \vertex(H_w) \coprod \left[\vertex(G) \setminus \{w\}\right].
\]
The first commutative diagram follows from the computation:
\[
\begin{split}
K_v 
&= K(\uparrow) = \left[ G(H_w)\right] (\uparrow)\\
&= G\left(\{H_w, \uparrow\}\right)\\
&= \left[G(\uparrow)\right](H_w)\\
&= G_v(H_w).
\end{split}
\]
The second commutative diagram follows from the computation:
\[
\begin{split}
K_v 
&= \left[G(H_w)\right](\uparrow)\\
&= G\left[H_w(\uparrow)\right]\\
&= G(C) = G.
\end{split}
\]
Here we used the fact that $H_w$ is a partially grafted corollas and $v \in \vertex(H_w)$ to infer that $H_w$ is a dioperadic graph (Example \ref{ex:dioperadicgraph}).  This implies that
\[
H_w(\uparrow) = C,
\]
the corolla with the same profiles as the other vertex in $H_w$.
\end{proof}

Next we consider an \emph{outer} coface map followed by a codegeneracy map.

\begin{lemma}
\label{graphiddstwo}
Suppose $K=G\left(\{C_u,H_w\}\right)$ is an outer properadic factorization, and $v \in \vertex(K)$ has one incoming flag and one outgoing flag.
\begin{enumerate}
\item
If $v \in \vertex(H_w)$, then the diagram
\[
\nicexy{
H_w \ar[d]_{s^v} \ar[r]^-{d^u} 
& K \ar[d]^{s^v}\\
(H_w)_v \ar[r]^-{d^u} 
& K_v
}\]
is commutative.
\item
If $v=u$, then the diagram
\[
\nicexy{
H_w \ar[dr]_{\Id} \ar[r]^-{d^u} 
& K \ar[d]^{s^u}\\
& H_w = K_u
}\]
is commutative.
\end{enumerate}
\end{lemma}

\begin{proof}
Here we use the decomposition
\[
\vertex(K) = \vertex(H_w) \coprod \{u\}.
\]
The first commutative diagram follows from the computation:
\[
\begin{split}
K_v
&= \left[G(H_w)\right] (\uparrow)\\
&= G\left[H_w(\uparrow)\right]\\
&= G\left[(H_w)_v\right].
\end{split}
\]
The second commutative diagram follows from the computation:
\[
\begin{split}
K_u 
&= \left[G(H_w)\right](\uparrow)\\
&= G\left(\{H_w, \uparrow\}\right)\\
&= \left[G(\uparrow)\right](H_w)\\
&= C(H_w) = H_w.
\end{split}
\]
Here we used the fact that $G$ is a partially grafted corollas and $v=u$ to infer that $G$ is a dioperadic graph.  This implies
\[
G(\uparrow) = C,
\]
the corolla with the same profiles as $w$.
\end{proof}

\subsection{Identities of Coface Maps}

Here we consider a graphical analog of the cosimplicial identity
\[
d^jd^i = d^i d^{j-1} \quad \text{for $i < j$}.
\]
This is the most interesting case because each of the coface maps can be either inner or outer.  In other words, each coface map may be associated with either a pair of closest neighbors (Theorem \ref{thm:cnbdfact}) or an almost isolated vertex (Theorem \ref{thm:opropfact}).  The result of mixing closest neighbors with almost isolated vertices is not immediately obvious. Therefore, the interaction of such coface maps takes on a much more complicated form than in either the finite ordinal category $\varDelta$ or the dendroidal category $\varOmega$.

The graphical analog of the above cosimplicial identity will take the following form.

\begin{definition}
\label{def:codimentwo}
We say that $\gupc$ has the  \textbf{codimension $2$ property} \index{codimension $2$ property!for graphical properads} if the following statement holds:  Given any two composable coface maps
\[
\nicearrow
\xymatrix{
\varGamma(K) \ar[r]^-{d^v} 
& \varGamma(H) \ar[r]^-{d^u} 
& \varGamma(G)
}
\]
in $\properad$, there exists a commutative square
\[
\nicearrow
\xymatrix{
\varGamma(K) \ar[r]^-{d^v} \ar[d]_{d^y} 
& \varGamma(H) \ar[d]^{d^u}\\
\varGamma(J) \ar[r]^-{d^x} 
&\varGamma(G)
}
\]
of coface maps such that $d^x$ is not obtainable from $d^u$ by changing the listing.
\end{definition}

\begin{remark}
The codimension $2$ property says that for each coface of codimension $2$, there are at least two different ways, even up to listing, to obtain $K$ from $G$ using iterated inner or outer properadic factorizations.  In what follows, when we say that $d^ud^v$ has another decomposition into two coface maps $d^xd^y$, we automatically mean that $d^x$ is not obtainable from $d^u$ by changing the listing. 
\end{remark}

\begin{theorem}
\label{gupccodimtwo}
$\gupc$ has the codimension $2$ property.
\end{theorem}

\begin{proof}
Suppose
\begin{equation}
\label{codimtwogupc}
\nicearrow
\xymatrix{
\varGamma(K) \ar[r]^-{d^v} 
& \varGamma(H) \ar[r]^-{d^u} 
& \varGamma(G)
}
\end{equation}
are two coface maps.  We must show that the composition $d^ud^v$ has another decomposition into two coface maps.  Each of $d^u$ and $d^v$ can be either an inner coface map or an outer coface map.  To improve readability, we will check these four cases in the next four lemmas.
\end{proof}

First we consider the case when $d^u$ and $d^v$ are both \emph{inner} coface map.

\begin{lemma}
\label{lem1:gupccodimtwo}
If both $d^u$ and $d^v$ are inner coface maps, then the composition $d^ud^v$ in \eqref{codimtwogupc} has another decomposition into two coface maps.
\end{lemma}

\begin{proof}
Here $H$ is obtained from $G$ by smashing together two closest neighbors $u_1$ and $u_2$ to form $u$, or equivalently there is an inner properadic factorization
\[
G = H(H_u),
\]
where $H_u$ is a partially grafted corollas with vertices the $u_i$.  Likewise, $K$ is obtained from $H$ by smashing together two closest neighbors $v_1$ and $v_2$ to form $v$, or equivalently there is an inner properadic factorization
\[
H = K(K_v),
\]
where $K_v$ is a partially grafted corollas with vertices the $v_i$.  In particular, we have
\[
G = H(H_u) = \left[K(K_v)\right](H_u).
\]
There are two cases.
\begin{enumerate}
\item
First suppose both $v_i$ are different from $u$.  We will show that $d^v$ and $d^u$ commute in a suitable sense.  We have
\[
G = K(\{K_v,H_u\}) = [K(H_u)](K_v)
\]
by associativity and unity of graph substitution.  Since both $H_u$ and $K_v$ are partially grafted corollas, there is a commutative diagram
\begin{equation}
\label{dudvcommute}
\nicexy{
K \ar[d]_{d^u} \ar[r]^-{d^v} 
& H \ar[d]^{d^u}\\
K(H_u) \ar[r]^-{d^v} & G
}
\end{equation}
of inner coface maps.
\item
Next suppose one $v_i$ is $u$.  Switching indices if necessary, we may assume that $v_1 = u$ and that there is an ordinary edge from $u_2$ to $u_1$ in $G$.  There are six possibilities regarding these vertices in $G$, depicted locally as follows.
\bigskip
\begin{center}
\begin{tikzpicture}
\matrix[row sep=1cm, column sep=1cm]{
\node [plain] (u11) {$u_1$};
&& \node [plain] (u12) {$u_1$};
&& \node [plain] (u13) {$u_1$}; &\\
\node [plain] (u21) {$u_2$}; 
& \node[plain] (v21) {$v_2$}; 
& \node [plain] (u22) {$u_2$};
& \node [plain] (v22) {$v_2$};
& \node [plain] (u23) {$u_2$};
& \node [plain] (v23) {$v_2$};\\
};
\draw [implies] (v21) to (u21);
\draw [implies] (u21) to (u11);
\draw [implies] (u22) to (u12);
\draw [implies] (v22) to (u12);
\draw [implies] (v23) to (u23);
\draw [implies] (v23) to (u13);
\draw [implies] (u23) to (u13);
\end{tikzpicture}
\end{center}

\bigskip
\begin{center}
\begin{tikzpicture}
\matrix[row sep=1cm, column sep=1cm]{
\node [plain] (u11) {$u_1$}; 
& \node [plain] (v21) {$v_2$};
& \node [plain] (u12) {$u_1$};
& \node [plain] (v22) {$v_2$}; 
& \node [plain] (u13) {$u_1$};
& \node [plain] (v23) {$v_2$};\\
\node[plain] (u21) {$u_2$};
&& \node [plain] (u22) {$u_2$};
&& \node [plain] (u23) {$u_2$};
&\\
};
\draw [implies] (u21) to (u11);
\draw [implies] (u11) to (v21);
\draw [implies] (u22) to (u12);
\draw [implies] (u22) to (v22);
\draw [implies] (u23) to (u13);
\draw [implies] (u23) to (v23);
\draw [implies] (u13) to (v23);
\end{tikzpicture}
\end{center}
Here the double arrow $\Longrightarrow$ means there is at least one ordinary edge in that direction between the indicated vertices.  The top left case says there is an ordinary edge from $v_2$ to $u_2$ but not to $u_1$.  The top middle case says there is an ordinary edge from $v_2$ to $u_1$ but not to $u_2$.  The top right case says there are ordinary edges from $v_2$ to both $u_1$ and $u_2$.  The three cases in the bottom row are interpreted similarly.  Going left to right along the top row and then along the bottom row, we will refer to them as case 1 through case 6 below.

In all six cases, we have
\begin{equation}
\label{gkkh}
G = [K(K_v)](H_u) = K[K_v(H_u)]
\end{equation}
by associativity of graph substitution, where $H_u$ is substituted into $v_1 \in \vertex(K_v)$.  The three cases in the top (resp., bottom) row correspond to the case where $v_1$ is the top (resp., bottom) vertex in $K_v$.  The graph $K_v(H_u)$ is a connected wheel-free graph with three vertices.  It is always possible to rewrite it using a different graph substitution involving partially grafted corollas as
\begin{equation}
\label{kyjx}
K_v(H_u) = K_y(J_x).
\end{equation}
Here $J_x$ is the partially grafted corollas defined by $v_2$ and $u_1$ (resp., $u_2$) in cases 2, 4, and 6 (resp., 1, 3, and 5), while $K_y$ is the partially grafted corollas containing $u_2$ (resp., $u_1$) and the combined vertex from $J_x$.  Using \eqref{kyjx} in \eqref{gkkh} and associativity of graph substitution, we obtain
\[
G = K[K_y(J_x)] = [K(K_y)](J_x).
\]
This in turn yields the commutative diagram
\[
\nicexy{
K \ar[d]_{d^y} \ar[r]^-{d^v} 
& H \ar[d]^{d^u}\\
K(K_y) \ar[r]^-{d^x} & G
}\]
of inner coface maps.
\end{enumerate}
\end{proof}

Next we consider the case where we begin with two \emph{outer} coface maps.  The graph-drawing convention in the proof of the previous lemma will be used below.

\begin{lemma}
\label{lem2:gupccodimtwo}
If both $d^u$ and $d^v$ are outer coface maps, then the composition $d^ud^v$ in \eqref{codimtwogupc} has another decomposition into two coface maps.
\end{lemma}

\begin{proof}
In this case, $H$ is obtained from $G$ by deleting an almost isolated vertex $v_1$.  Equivalently, there is an outer properadic factorization
\[
G = P_1(\{C_{v_1}, H\}) = P_1(H)
\]
for some partially grafted corollas $P_1$ with vertices $\{v_1,u_1\}$, where $C_{v_1}$ is a corolla, and $H$ is substituted into $u_1$.  Likewise, $K$ is obtained from $H$ by deleting an almost isolated vertex $v_2$.   Equivalently, there is an outer properadic factorization
\[
H = P_2(\{C_{v_2}, K\}) = P_2(K)
\]
for some partially grafted corollas $P_2$ with vertices $\{v_2,u_2\}$, where $C_{v_2}$ is a corolla, and $K$ is substituted into $u_2$.  In particular, we have
\[
G = P_1(H) = P_1[P_2(K)] = [P_1(P_2)](K)
\]
by associativity of graph substitution.

There are four cases, depending on whether each $u_i$ is the top or the bottom vertex in $P_i$.  Since the four cases have very similar proofs, we only consider in detail the case where each $u_i$ is the \emph{bottom} vertex in $P_i$.  So we may visualize $G$ and $H$ as follows.
\begin{center}
\begin{tikzpicture}
\matrix[row sep=1cm, column sep=3cm]{
\node [plain] (v1) {$v_1$};
& \node [plain] (v2) {$v_2$};\\
\node [fatplain] (H) {$H$};
& \node [fatplain] (K) {$K$};\\
};
\draw [implies] (H) to (v1);
\draw [implies] (K) to (v2);
\end{tikzpicture}
\end{center}
The composition $d^ud^v$ corresponds to the construction of $K$ from $G$ by first deleting $v_1$ and then deleting $v_2$.  The desired alternative factorization will come from a different construction of $K$ from $G$.

To obtain this alternative construction, note that the graph substitution $P_1(P_2)$ is a connected wheel-free graph with three vertices with each $P_i$ having two vertices.  There are three possible shapes for $P_1(P_2)$:
\begin{center}
\begin{tikzpicture}
\matrix[row sep=1cm, column sep=1cm]{
\node [plain] (v11) {$v_1$};
&& \node [plain] (v12) {$v_1$};
&& \node [plain] (v13) {$v_1$}; &\\
\node [plain] (v21) {$v_2$};
& \node [plain] (u21) {$u_2$};
& \node [plain] (v22) {$v_2$};
& \node [plain] (u22) {$u_2$};
& \node [plain] (v23) {$v_2$};
& \node [plain] (u23) {$u_2$};\\
};
\draw [implies] (u21) to (v21);
\draw [implies] (u21) to (v11);
\draw [implies] (v21) to (v11);
\draw [implies] (u22) to (v22);
\draw [implies] (v22) to (v12);
\draw [implies] (u23) to (v23);
\draw [implies] (u23) to (v13);
\end{tikzpicture}
\end{center}
From left to right, we call them case 1 through case 3.  In case 1 and case 2, the vertices $v_i$ are closest neighbors in $G$.  When they are smashed together, the combined vertex is almost isolated, whose deletion yields $K$.  In case 3, $v_2$ is already almost isolated in $G$, so it can be deleted first, followed by the deletion of $v_1$.

More formally, in each of the three cases we can rewrite $P_1(P_2)$ using a different graph substitution involving partially grafted corollas as
\[
P_1(P_2) = P_1'(P_2'),
\]
just like we did in \eqref{kyjx}.  In cases 1 and 2 (resp., case 3), $P_2'$ is the partially grafted corollas defined by the $v_i$ (resp., $v_1$ and $u_2$), and $P_1'$ is the partially grafted corollas defined by $u_2$ (resp., $v_2$) and the combined vertex from $P_2'$.  Recall that in all cases $K$ is supposed to be substituted into $u_2$.
\begin{enumerate}
\item
In cases 1 and 2, we have
\[
G = [P_1'(P_2')](K) = P_1'(\{P_2', K\}) = [P_1'(K)](P_2')
\]
by associativity and unity of graph substitution.  This yields the alternative factorization
\begin{equation}
\label{kpkouter}
\nicexy@C+12pt{
K \ar[r]^-{\text{outer}} & P_1'(K) \ar[r]^-{\text{inner}} & G
}
\end{equation}
of $d^ud^v$ into an outer coface map followed by an inner coface map.
\item
Similarly, in case 3 we have
\[
G = [P_1'(P_2')](K) = P_1'[P_2'(K)].
\]
This yields the alternative factorization
\[
\nicexy@C+12pt{
K \ar[r]^-{\text{outer}} & P_2'(K) \ar[r]^-{\text{outer}} & G
}\]
of $d^ud^v$ into two outer coface maps.
\end{enumerate}
\end{proof}

We now consider the third case of the proof of Theorem \ref{gupccodimtwo}.

\begin{lemma}
\label{lem3:gupccodimtwo}
Suppose $d^u$ is an inner coface map, and $d^v$ is an outer coface map.  Then the composition $d^ud^v$ in \eqref{codimtwogupc} has another decomposition into two coface maps.
\end{lemma}

\begin{proof}
In this case, $H$ is obtained from $G$ by smashing together two closest neighbors $u_1$ and $u_2$ to form $u$.  Equivalently, there is an inner properadic factorization
\[
G = H(H_u),
\]
in which $H_u$ is a partially grafted corollas with vertices the $u_i$, and $H_u$ is substituted into $u \in \vertex(H)$.  Likewise, $K$ is obtained from $H$ by deleting an almost isolated vertex $w$, or equivalently there is an outer properadic factorization
\[
H = P(\{C_w,K\}) = P(K)
\]
in which $P$ is a partially grafted corollas, and $C_w$ is a corolla.  In particular, we have the decomposition
\[
\vertex(H) = \vertex(K) \coprod \{w\}
\]
and the graph substitution decomposition
\[
G = [P(K)](H_u).
\]
There are two sub-cases.
\begin{enumerate}
\item
If $u \in \vertex(H)$ is actually in $K$, then we have
\[
G = [P(K)](H_u) = P[K(H_u)]
\]
by associativity of graph substitution.  This yields the alternative decomposition
\[
\nicexy@C+12pt{
K \ar[r]^-{\text{inner}} & K(H_u) \ar[r]^-{\text{outer}} & G
}\]
of $d^ud^v$ into an inner coface map followed by an outer coface map.
\item
If $u = w$, then we have
\[
G = P(\{H_u,K\}) = [P(H_u)](K).
\]
Note that $P(H_u)$ is a connected wheel-free graph with three vertices.  Thus, as before \eqref{kyjx} there is another graph substitution decomposition
\[
P(H_u) = P'(H')
\]
involving two other partially grafted corollas.  Then we have
\[
G = [P(H_u)](K) = [P'(H')](K) = P'[H'(K)].
\]
Since both $H'$ and $P'$ are partially grafted corollas, there is an alternative decomposition
\[
\nicexy@C+12pt{
K \ar[r]^-{\text{outer}} & H'(K) \ar[r]^-{\text{outer}} & G
}\]
of $d^ud^v$ into two outer coface maps.
\end{enumerate}
\end{proof}

The final case of the proof of Theorem \ref{gupccodimtwo} is considered next.

\begin{lemma}
\label{lem4:gupccodimtwo}
Suppose $d^u$ is an outer coface map, and $d^v$ is an inner coface map.  Then the composition $d^ud^v$ in \eqref{codimtwogupc} has another decomposition into two coface maps.
\end{lemma}

\begin{proof}
As before, we will obtain an alternative construction relating $K$ and $G$.  Here $H$ is obtained from $G$ by deleting an almost isolated vertex $u$.  Equivalently, we have
\[
G = P(\{C_u,H\}) = P(H),
\]
in which $P$ is a partially grafted corollas, and $C_u$ is a corolla.  Likewise, $K$ is obtained from $H$ by smashing together two closest neighbors $v_1$ and $v_2$ to form $v$, or equivalently
\[
H = K(K_v),
\]
in which $K_v$ is a partially grafted corollas with vertices the $v_i$.  In particular, we have
\[
G = P[K(K_v)] = [P(K)](K_v)
\]
by associativity of graph substitution.  Since $P$ and $K_v$ are partially grafted corollas, this yields the alternative decomposition
\[
\nicexy@C+12pt{
K \ar[r]^-{\text{outer}} & P(K) \ar[r]^-{\text{inner}} & G
}\]
of $d^ud^v$ into an outer coface map followed by an inner coface map.
\end{proof}

With Lemmas \ref{lem1:gupccodimtwo}, \ref{lem2:gupccodimtwo}, \ref{lem3:gupccodimtwo}, and \ref{lem4:gupccodimtwo}, the proof of Theorem \ref{gupccodimtwo} is complete.

\begin{remark}
\label{rk:uniquecodimtwo}
If one goes through the proof of Theorem \ref{gupccodimtwo} carefully, one can see that the given composition $d^ud^v$ of coface maps has a \emph{unique} alternative decomposition $d^xd^y$ into coface maps.  Here uniqueness is understood to be up to listing.
\end{remark}

\section{Graphical Category}
\label{sec:gupcgraphcat}

The main purpose of this section is to define the maps in the graphical category $\varGamma$ for connected wheel-free graphs via the concept of subgraphs.  The point of using subgraphs is to avoid some bad behavior for general properad maps between graphical properads, as exhibited in section \ref{sec:mapsfromgprop}.  At the end of this section, we establish a graphical analog of the epi-mono factorization.

\subsection{Input and Output Relabeling}
\label{inoutputrelabeling}

The concept of a subgraph is defined via outer coface maps, which were discussed in section \ref{sec:coface}, and input/output relabeling, which we now discuss.

Suppose $G \in \gupc$ has input/output profiles $\dch$.  Recall from chapter \ref{ch:graph} that $G$ has a listing, which consists of a labeling of the inputs and outputs of each vertex as well as of the whole graph $G$.  Suppose $\sigma \in \Sigma_{|\ud|}$ and $\tau \in \Sigma_{|\uc|}$.  Then there is a graph
\[\label{note:relabeling}
\sigma G \tau \in \gupc\dcsigma
\]
with profiles $(\uc\tau;\sigma\ud)$ obtained from $G$ by permuting the input/output profiles $\dch$ using $\sigma$ and $\tau$.  Since only the listing of the whole graph has changed, there are canonical bijections
\[
\begin{split}
\vertex(G) &= \vertex(\sigma G\tau),\\
\edge(G) &= \edge(\sigma G\tau),\\
\edgei(G) &= \edgei(\sigma G \tau).
\end{split}
\]
The first two bijections induce a canonical isomorphism
\[
\nicexy{
G \ar[r]^-{(\tau;\sigma)}_-{\cong} & \sigma G \tau \in \properad
}\]
of graphical properads, called \textbf{input/output relabeling}, or simply \textbf{relabeling}.\index{relabeling}  Note that the identity map is a relabeling, with both $\sigma$ and $\tau$ the respective identity permutations.

\subsection{Subgraphs}

In the following definition, we will make use of the concepts of relabeling just discussed and of outer coface maps in Definition \ref{def:outercoface}.

\begin{definition}
\label{def:subgraph}
Suppose $G,K \in \gupc$.  A map
\[
\nicexy{G \ar[r]^-{f} & K}
\]
of graphical properads is called a \textbf{subgraph} \index{subgraph!of graphical properads} if $f$ admits a decomposition into outer coface maps and relabelings.  In this case, we also call $G$ a \textbf{subgraph of} $K$.
\end{definition}

\begin{remark}
Subgraphs are closed under compositions.  The identity map is a relabeling, hence a subgraph.
\end{remark}

\begin{remark}
Since a relabeling is a fairly simple canonical isomorphism of graphical properads induced by input/output relabeling permutations, we will often not mention it in discussing subgraphs.
\end{remark}

\begin{remark}
To say that $f$ is a subgraph is \emph{not} the same thing as saying that $f$ cannot be factored by inner coface maps.  In fact, by the proof of Lemma \ref{lem2:gupccodimtwo}, in particular  \eqref{kpkouter}, a compostion of outer coface maps may have other factorizations that involve inner coface maps.  The requirement to being a subgraph is that there exists one factorization of the map into outer coface maps up to relabeling.
\end{remark}

We will use Notation \ref{notation:gh} below.

\begin{remark}
Let us explain the geometric meaning of a subgraph.  Recall that an outer coface map $H \to K$ corresponds to an outer properadic factorization $K=P(H)$ (Definition \ref{def:opropfact}), where $P$ is a partially grafted corollas with vertices $\{u,w\}$, and $H$ is substituted into one of its two vertices, say $w$.  By Theorem \ref{thm:opropfact} this is equivalent to saying that $u$ is an almost isolated vertex in $K$ and that $H$ is obtained from $K$ by deleting $u$.  If $\nicexy{G \ar[r]^-{f} & K}$ is a subgraph, then it has a decomposition into outer coface maps.  Therefore, $G$ is obtained from $K$ by repeatedly deleting almost isolated vertices.

Furthermore, given an outer properadic factorization $K=P(H)$ as above, there are canonical injections
\[
\vertex(H) \hookrightarrow \vertex(K) \andspace
\edge(H) \hookrightarrow \edge(K).
\]
These injections induce the outer coface map $H \to K$.  If $f$ as above is a subgraph, then there are injections
\[
\vertex(G) \hookrightarrow \vertex(K) \andspace
\edge(G) \hookrightarrow \edge(K).
\]
So we can think of the graph $G$ as sitting inside the graph $K$.  The map $f$ records the part of $K$ outside of $G$, or equivalently the part of $K$ that gets deleted to obtain $G$.  We formalize this idea in the following observation.
\end{remark}

In the next result, we characterize subgraphs in terms of graph substitution.  This result will be used many times in what follows.

\begin{theorem}
\label{thm:gupcsubgraph}
Suppose $\nicexy{G \ar[r]^-{f} & K}$ is a map of graphical properads.  Then the following statements are equivalent.
\begin{enumerate}
\item
$f$ is a subgraph.
\item
There exists a graph substitution decomposition
\[
K = H(G)
\]
in $\gupc$ such that $f$ sends the edges and vertices in $G$ to their corresponding images in $H(G)$. \index{subgraph!characterization}
\end{enumerate}
\end{theorem}

\begin{proof}
First suppose $f$ is a subgraph.  If $f$ consists of a single relabeling, then $K = \sigma G \tau$ for some suitable permutations $\sigma$ and $\tau$.  So there is a graph substitution decomposition
\[
K = (\sigma C\tau)(G),
\]
where $C$ is the corolla with the same profiles as $G$.  In the rest of this proof, we will omit mentioning relabelings.

Next suppose
\[
f=d_n \cdots d_1
\]
in which each $d_i$ is an outer coface map.  We show by induction on $n$ that $K$ admits a graph substitution decomposition as stated.  If $n=1$, then $f = d_1$ is an outer coface map, and $K=P(G)$ for some partially grafted corollas $P$.

Suppose $n > 1$.  Then there is a factorization
\[
\nicexy{
G \ar[r]^-{g} & H_{n-1} \ar[r]^-{d_n} & K
}\]
of $f$ in which
\[
g = d_{n-1} \cdots d_1.
\]
Since $d_n$ is an outer coface map, there is a graph substitution decomposition
\[
K = P(H_{n-1})
\]
for some partially grafted corollas $P$.  Moreover, $g$ is by definition a subgraph.  By induction hypothesis, there is a graph substitution decomposition
\[
H_{n-1} = H'(G)
\]
in $\gupc$.  Therefore, we have the desired decomposition
\[
K = P[H'(G)] = [P(H')](G)
\]
by associativity of graph substitution.

Conversely, suppose $K=H(G)$ as stated.  We show that $f$ is a subgraph by induction on $m=|\vertex(H)|$.  If $m=1$, then $H$ is a corolla, and $f$ is the identity map, which is a subgraph.

Suppose $m > 1$, and $G$ is substituted into $w \in \vertex(H)$.  There must exist an almost isolated vertex $v \not= w$ in $H$ (Theorem \ref{lem2:gupcenoughface}).  So there is an outer properadic factorization
\[
H = P(H_v),
\]
in which $P$ is a partially grafted corollas and $H_v$ is obtained from $H$ by deleting the almost isolated vertex $v$.  There is a corresponding outer coface map
\[
H_v \to H.
\]
Since $w \in \vertex(H_v)$, the graph substitution $H_v(G)$ makes sense.  So we have
\[
K = [P(H_v)](G) = P[H_v(G)],
\]
and $f$ factors into
\[
\nicexy@C+10pt{
G \ar[r] \ar[dr]_-{f} & H_v(G) \ar[d]\\
& K = P[(H_v(G)].
}\]
We have $|\vertex(H_v)| = m-1$, so the first map $G \to H_v(G)$ is a subgraph by induction hypothesis .  The second map $H_v(G) \to K$ is an outer coface map because $P$ is a partially grafted corollas.  Therefore, $f$ is the composition of a subgraph and an outer coface map, so it is a subgraph.
\end{proof}

The following observation provides some small examples of subgraphs.

\begin{corollary}
\label{cor:gupcsubgraph}
Suppose $K$ is a connected wheel-free graph.
\begin{enumerate}
\item
Suppose $x$ and $y$ are closest neighbors in $K$, and $G$ is the partially grafted corolla defined by them.  Then the map $G \to K$ is a subgraph.
\item
For each vertex $v$ in $K$, the corolla inclusion $C_v \to K$ is a subgraph.
\item
For each edge $e$ in $K$, the edge inclusion $\uparrow_e ~\to K$ is a subgraph.
\end{enumerate}
\end{corollary}

\begin{proof}
For the first assertion, by Theorem \ref{thm:cnbdfact} there is a graph substitution decomposition $K=H(G)$.  So by Theorem \ref{thm:gupcsubgraph} the map $G \to K$ is a subgraph.  For the second assertion, use the graph substitution decomposition
\[
K = K(C_v),
\]
in which a corolla is substituted into each vertex, and Theorem \ref{thm:gupcsubgraph}.  

For the last assertion, if $K$ itself is an exceptional edge, then $\uparrow_e \to K$ is the identity map.  If $K$ is ordinary, then $e$ is adjacent to some vertex $v$ in $K$.  The map $\uparrow_e ~\to K$ factors into
\[
\nicexy{
\uparrow_e \ar[r] & C_v \ar[r] & K.
}\]
The first map is an outer coface map that identifies the edge $e$ as a leg of the corolla $C_v$ with the same profiles as $v$ (Example \ref{ex:legcoface}).  The second map is a corolla inclusion, which is a subgraph by the previous part.  Therefore, their composition is also a subgraph.
\end{proof}

The following observation says that a subgraph is uniquely determined by its input/output profiles.

\begin{lemma}
\label{subgraphunique}
Suppose $G \in \gupc$, and $(\ue;\uf)$ is a pair of $\edge(G)$-profiles.  Then there exists at most one subgraph
\[
H \to G
\]
in which $H$, when regarded as an element in $\varGamma(G)$, has profiles $\feh$.\index{subgraph!uniqueness}
\end{lemma}

\begin{proof}
A subgraph $J \to G$ is by definition a composition of outer coface maps, so $J$ is obtained from $G$ by repeatedly deleting almost isolated vertices.  If two subgraphs have the same $\edge(G)$-profiles, then they are obtained from $G$ by deleting the same subset of vertices.  So they have the same sets of vertices.  Finally, note that if $x$ and $y$ are vertices in a subgraph $J$, then all the internal edges between $x$ and $y$ in $G$ are also in $J$.  Therefore, two subgraphs with the same $\edge(G)$-profiles have the same sets of internal edges as well.
\end{proof}

\subsection{Images}

Before we can define a map in the properadic graphical category, we need to define the image of a properad map between graphical properads.

Recall from Lemma \ref{lem:mapbtwgprop} that a properad map $f$ between two graphical properads is determined by a pair of functions $(f_0,f_1)$.   The function $f_0$ goes between the sets of edges (i.e., colors), and $f_1$ sends each vertex in the source to an element in the target with the correct profiles.

\begin{definition}
Suppose $\nicexy{G \ar[r]^{f} & K}$ is a properad map of graphical properads.  The \textbf{image of $G$} \index{image!of graphical properad} is defined as the graph substitution
\[\label{note:image}
f(G) = \left[f_0 G\right] \left(\{f_1(u)\}_{u \in \vertex(G)}\right) \in \varGamma(K),
\]
in which $f_0G$ is the graph obtained from $G$ by applying $f_0$ to its edges.
\end{definition}

\begin{example}
\label{ex:gupcimage}
Here we describe the images of codegeneracy maps, coface maps, subgraphs, and changes of vertex listings.
\begin{enumerate}
\item
For a codegeneracy map
\[
\nicexy{G \ar[r]^-{s} & G_v},
\]
by definition we have $G_v = G(\uparrow)$.  So the image $s(G)$ is $G_v$.
\item
For an inner coface map
\[
\nicexy{G \ar[r]^-{d_{\inp}} & K},
\]
by definition we have $K=G(P)$ for some partially grafted corollas $P$.  So once again the image $d_{\inp}(G)$ is $K$.
\item
For an outer coface map
\[
\nicexy{G \ar[r]^-{d_{\out}} & K},
\]
by definition we have $K=P(G)$ for some partially grafted corollas $P$.  So the image $d_{\out}(G)$ is $G \in \varGamma(K)$, where $G$ is regarded as an $\edge(K)$-colored $\khat$-decorated graph via the inclusions
\[
\edge(G) \hookrightarrow \edge(K) \andspace
\vertex(G) \hookrightarrow \vertex(K).
\]
\item
For a subgraph 
\[
\nicexy{G \ar[r]^-{f} & K},
\]
by Theorem \ref{thm:cnbdfact} we have $K = H(G)$ for some $H \in \gupc$.  So the image $f(G)$ is $G \in \varGamma(K)$.
\item
Suppose $K$ is obtained from $G$ by changing the listings at a subset of vertices, and
\[
\nicexy{G \ar[r]^-{f} & K}
\]
is the corresponding properad isomorphism.  Then the image $f(G)$ is $K$.
\end{enumerate}
\end{example}

Just as one would expect from an image, the original map factors through it.

\begin{lemma}
\label{lem1:gupcimage}
Suppose $\nicexy{G \ar[r]^{f} & K}$ is a properad map of graphical properads.  Then there is a canonical commutative diagram
\[
\nicexy@C+10pt{
G \ar[r]^-{g} \ar[dr]_-{f} & f(G) \ar[d]^-{h}\\
& K}
\]
of properad maps between graphical properads.
\end{lemma}

\begin{proof}
First we define the maps $g$ and $h$.  Note that the image $f(G) \in \varGamma(K)$ is by definition an $\edge(K)$-colored $\khat$-decorated graph.  So each of its edges/vertices has a canonical image in $K$.  This defines the properad map $h$, which sends every vertex $w$ in $f(G)$ to the corolla $C_w \in \varGamma(K)$.

The map $g$ is defined by sending each edge $e$ in $G$ to the edge $f_0(e)$ in $f_0(G)$, which in turn yields an edge in $f(G)$.  For a vertex $u$ in $G$, $g_1(u)$ is $f_1(u)$ in $\varGamma(f(G))$.

To see that $f = hg$, observe that for an edge $e$ in $G$, $h$ sends $f_0(e)$ in $f(G)$ to the edge in $K$ with the same name.  For a vertex $u$ in $G$, $h$ sends $f_1(u)$ in $\varGamma(f(G))$ to $f_1(u)$ in $\varGamma(K)$ because $h$ sends each vertex to the corresponding corolla.
\end{proof}

\begin{example}
\label{ex:cofacesubgraph}
If $\nicexy{G \ar[r]^-{f} & K}$ is a codegeneracy map, a coface map, a subgraph, or a change of vertex listings, then the map $f(G) \to K$ is a subgraph.  Indeed, if $f$ is a codegeneracy map, an inner coface map, or a change of vertex listings, then $f(G)=K$.  If $f$ is an outer coface map or, more generally, a subgraph, then $f(G)$ is $G$.
\end{example}

\begin{example}
For the map $\nicexy{G \ar[r]^-{\varphi} & H}$ in Example \ref{ex:phizero}, the image $\varphi(G)$ is
\begin{center}
\begin{tikzpicture}
\matrix[row sep=2cm, column sep=1cm]{
\node [plain] (w) {$w$};\\
\node [plain] (v) {$v$};\\
};
\draw [arrow] (v) to node{$f$} (w);
\draw [inputleg] (w) to node[below left=.1cm]{$e_i$} +(-.7cm,-.5cm);
\draw [inputleg] (w) to node[below right=.1cm]{$g_i$} + (.7cm,-.5cm);
\draw [outputleg] (v) to node[above left=.1cm]{$e_o$} +(-.7cm,.5cm);
\draw [outputleg] (v) to node[above right=.1cm]{$g_o$} +(.7cm,.5cm);
\end{tikzpicture}
\end{center}
such that the map $G \to \varphi(G)$ satisfies
\[
(i_1,i_2) \longmapsto (e_i,g_i) \andspace
(o_1,o_2) \longmapsto (e_o,g_o).
\]
The map $\varphi(G) \to H$ sends the edges $e_*$ (resp., $g_*$ and $f$) to $e$ (resp., $g$ and $f$) in $H$, and sends the vertices $v$ and $w$ to the corresponding corollas.  Note that $\varphi(G) \to H$ is \emph{not} a subgraph because each outer coface map (and inner coface map as well) must increase the number of vertices by $1$.
\end{example}

\begin{example}
For the map $\nicexy{G \ar[r]^-{\varphi} & G}$ in Example \ref{ex:idonedge}, the image $\varphi(G)$ is
\begin{center}
\begin{tikzpicture}
\matrix[row sep=1cm,column sep=2cm] {
& \node [plain] (v1) {$v_1$}; & \\
& \node [plain] (v2) {$v_2$}; & \\
\node [plain] (u1) {$u_1$}; && \node [plain] (u2) {$u_2$}; \\
};
\draw [arrow,bend left=40] (u1) to node{$e_1$} (v1);
\draw [arrow,bend right=40] (u2) to node[swap]{$f_1$} (v1);
\draw [arrow] (u1) to node[swap]{$f_2$} (v2);
\draw [arrow] (u2) to node{$e_2$} (v2);
\end{tikzpicture}
\end{center}
The map $G \to \varphi(G)$ sends the edges $e$ and $f$ in $G$ to the edges $e_2$ and $f_2$ in $\varphi(G)$.  The map $\varphi(G) \to G$ sends the edges $e_*$ (resp., $f_*$) in $\varphi(G)$ to the edge $e$ (resp., $f$) in $G$, and sends the vertices $v_*$ (resp., $u_*$) to the corolla $C_v$ (resp., $C_u$).  Note that, once again, $\varphi(G) \to G$ is \emph{not} a subgraph.
\end{example}

One can similarly check that for each of the maps in Examples \ref{ex1:noninjection}, \ref{ex2:noninjection}, and \ref{ex3:noninjection}, the map from the image to the target is \emph{not} a subgraph.

\subsection{Graphical Maps}

We will use the factorization in Lemma \ref{lem1:gupcimage} in the next definition.

\begin{definition}
\label{def:graphicalmap}
A \index{properadic graphical map} \index{graphical map!properadic} \textbf{properadic graphical map}, or simply a \textbf{graphical map}, is defined as a properad map
\[
\nicexy{G \ar[r]^-{f} & K}
\]
between graphical properads such that the map $f(G) \to K$ is a subgraph.
\end{definition}

\begin{example}
Codegeneracy maps, coface maps, subgraphs, and changes of vertex listings are all graphical maps by Example \ref{ex:cofacesubgraph}.  On the other hand, the properad maps in Examples \ref{ex:phizero}, \ref{ex:idonedge}, \ref{ex1:noninjection}, \ref{ex2:noninjection}, and \ref{ex3:noninjection} are \emph{not} graphical maps.
\end{example}

\begin{example}
Suppose $\nicexy{G \ar[r]^-{f} & K}$ is a properad map between graphical properads.  If both $G$ and $K$ are linear graphs (resp., unital trees or simply connected graphs), then the map $f(G) \to K$ is a subgraph, so $f$ is a graphical map.  This is a consequence of Lemma \ref{lem:gammasconn}.
\end{example}

We will need the following observation, which says that graphical maps are closed under compositions.

\begin{lemma}
\label{lem1:graphicalmapclosed}
Suppose $\nicexy{G \ar[r]^-{f} & K}$ and $\nicexy{K \ar[r]^-{g} & M}$ are graphical maps.  Then the composition $\nicexy{G \ar[r]^-{gf} & M}$ is also a graphical map.
\end{lemma}

\begin{proof}
Both $f(G) \to K$ and $g(K) \to M$ are subgraphs by definition.  By Theorem \ref{thm:gupcsubgraph} there are graph substitution decompositions
\[
K = H_1(f(G)) \andspace 
M = H_2(g(K)).
\]
Therefore, by associativity of graph substitution we have
\[
\begin{split}
M &= H_2\left[g(H_1(f(G)))\right]\\
&= H_2\left[(g_0H_1)(gf(G))\right]\\
&= \left[H_2(g_0H_1)\right] (gf(G)).
\end{split}
\]
Therefore, by Theorem \ref{thm:gupcsubgraph} again, we conclude that $gf(G) \to M$ is a subgraph.
\end{proof}

Next we give another characterization of a graphical map.  By definition a graphical map sends the source to a subgraph of the target.  The following observation says that a graphical map also sends every subgraph of the source to a subgraph of the target.

\begin{theorem}
\label{thm:graphicalmapchar}
Suppose $\nicexy{G \ar[r]^-{f} & K}$ is a properad map between graphical properads.  Then the following statements are equivalent.
\begin{enumerate}
\item
$f$ is a graphical map.
\item
For each subgraph $\nicexy{H \ar[r]^-{\varphi} & G}$, the map $\nicexy{f(H) \ar[r] & K}$ is a subgraph.\index{graphical map!characterization}
\end{enumerate}
\end{theorem}

\begin{proof}
The direction $(2) \Longrightarrow (1)$ holds because the identity map on $G$ is a subgraph.

To prove $(1) \Longrightarrow (2)$, suppose $f$ is a graphical map, so $f(G) \to K$ is a subgraph.  If $G =~ \uparrow$, then its only subgraph is itself.  So we may assume that $G \not=~ \uparrow$.  

Pick a subgraph $\nicexy{H \ar[r]^-{\varphi} & G}$.  We show by downward induction on $m = |\vertex(H)|$ that $f(H)$ is a subgraph of $K$.  First suppose $m=|\vertex(G)|$.  In this case $H$ has all the vertices of $G$.  But it is also true that for vertices $x$ and $y$ in $H$, all the internal edges between them in $G$ are also in $H$ because $H$ is obtained from $G$ by repeatedly deleting almost isolated vertices.  Therefore, we conclude that $H=G$, and by assumption $f(G) \to K$ is a subgraph.

Next suppose $m < |\vertex(G)|$.  Since $H$ is a subgraph of $G$, by Theorem \ref{thm:gupcsubgraph} there is a graph substitution decomposition
\[
G = J(H)
\]
for some $J \in \gupc$, which has at least two vertices because $m < |\vertex(G)|$.  The vertex $v \in \vertex(J)$ into which $H$ is substituted must have a closest neighbor $u$ (Example \ref{ex:cnbdexist}).  By Theorem \ref{thm:cnbdfact}, there is an inner properadic factorization
\[
J = N(P),
\]
in which $P$ is a partially grafted corollas with vertices $u$ and $v$. Therefore, we have
\[
G = \left[N(P)\right](H) = N\left[P(H)\right],
\]
so by Theorem \ref{thm:gupcsubgraph} $P(H)$ is a subgraph of $G$.  Since $P$ has two vertices, $P(H)$ has $m+1$ vertices. By induction hypothesis $f(P(H))$ is a subgraph of $K$.  By Theorem \ref{thm:gupcsubgraph} we have a graph substitution decomposition
\[
\begin{split}
K 
&= M\left[f(P(H))\right]\\
&= \left[M(f_0P)\right](f(H))
\end{split}
\]
for some $M \in \gupc$.  Therefore, by Theorem \ref{thm:gupcsubgraph} once again, $f(H)$ is a subgraph of $K$.
\end{proof}

\begin{corollary}
\label{cor1:gupcgraphicalmap}
Suppose $\nicexy{G \ar[r]^-{f} & K}$ is a graphical map.  Then for each vertex $v$ in $G$, $f_1(v)$ is a subgraph of $K$.
\end{corollary}

\begin{proof}
Use Theorem \ref{thm:graphicalmapchar} and the fact that the corolla inclusion $C_v \to G$ is a subgraph (Corollary \ref{cor:gupcsubgraph}).
\end{proof}

\subsection{Graphical Category}

Lemma \ref{lem1:graphicalmapclosed} is used in the next definition to ensure that $\varGamma$ is really a category.

\begin{definition}
\label{def:graphicalcat}
The \index{properadic graphical category} \index{graphical category!properadic} \textbf{properadic graphical category}, or simply the \textbf{graphical category}, is the category $\varGamma$\label{note:gamma} with
\begin{itemize}
\item
objects the graphical properads $\varGamma(G)$ for $G \in \gupc$, and
\item
morphisms $\varGamma(G) \to \varGamma(H) \in \properad$ the properadic graphical maps.
\end{itemize}
Denote by
\[
\nicexy{\varGamma \ar[r]^-{\iota} & \properad}
\]
the non-full subcategory inclusion.
\end{definition}

\begin{remark}
The graphical category $\varGamma$ is small because there is only a set of $1$-colored graphs.
\end{remark}

\begin{remark}
\label{rk:deltaingamma}
The graphical category contains a full subcategory that is isomorphic to the finite ordinal category.  Indeed, there is a full subcategory $\varGamma(\ULin)$ of $\varGamma$ consisting of the objects $\varGamma(L_n)$ for $n \geq 0$, where $L_n$ is the $1$-colored linear graph with $n$ vertices (Definition \ref{def:wheelfreegraphs} and Remark \ref{rk:lineargraphln}).  There is an isomorphism of categories
\[\label{note:delta}
\varGamma(\ULin) \cong \varDelta,
\]
sending $\varGamma(L_n)$ to the finite ordered set $[n] = \{0<1<\cdots<n\} \in \varDelta$.
\end{remark}

\begin{remark}
\label{rk:omegaingamma}
The graphical category also contains a full subcategory that is equivalent to the Moerdijk-Weiss dendroidal category $\varOmega$ \cite{mw1}.  Indeed, there is a full subcategory $\varGamma(\uoperad)$ of $\varGamma$ consisting of the objects $\varGamma(T)$, where $T$ is any unital tree (Definition \ref{def:wheelfreegraphs}).  There is an equivalence of categories
\[\label{note:omega}
\varGamma(\uoperad) \simeq \varOmega.
\]
This is true because the free properad $\varGamma(T)$ is obtained from its underlying operad, which is the free operad generated by $T$, by adding empty components.  

There are similar full subcategories of $\varGamma$ defined by the indicated subsets of graphs:
\begin{itemize}
\item
$\varTheta = \varGamma(\gupd)$,\label{note:theta} where $\gupd$ is the set of simply connected graphs,
\item
$\varGamma_i = \varGamma(\gupci)$,\label{note:gammai} where $\gupci$ is the set of connected wheel-free graphs with non-empty inputs, and
\item
$\varGamma_o = \varGamma(\gupco)$,\label{note:gammao} where $\gupco$ is the set of connected wheel-free graphs with non-empty outputs.
\end{itemize}
Using the above isomorphism $\varGamma(\ULin) \cong \varDelta$ and equivalence $\varGamma(\uoperad) \simeq \varOmega$ of categories,  there is a commutative diagram of subcategory inclusions:
\begin{equation}
\label{differerntgammas}
\nicexy{
&& \varGamma_i \ar[dr] &\\
\varDelta \ar[urr] \ar[r] & \varOmega \ar[dr] \ar[r] & \varTheta \ar[r] &  \varGamma\\
&& \varGamma_o \ar[ur] &}
\end{equation}
\end{remark}

\begin{remark}
By Theorem \ref{omegainfinite}, every element in $\varTheta$ is a finite set.  In particular, every element in the Moerdijk-Weiss dendroidal category $\varOmega$ is a finite set.  On the other hand, every graphical properad \emph{not} in $\varTheta$ is an infinite set.
\end{remark}

\subsection{Factorization of Graphical Maps}

Now we observe that properadic graphical maps have an analog of the epi-mono factorization in the finite ordinal category $\varDelta$ and the dendroidal category $\varOmega$ \cite{mt} (Lemma 2.3.2).  We will use the factorization and notations in Lemma \ref{lem1:gupcimage}. We will often omit mentioning isomorphisms induced by changes of listings.

\begin{theorem}
\label{thm:gupcepimono}
Suppose $\nicexy{G \ar[r]^-{f} & K}$ is a graphical map.  Then there is a factorization
\[
\nicexy@C+10pt{
G \ar[rr]^-{f}  \ar[d]_{\sigma} \ar[drr]^-{g} && K\\
G_1 \ar[r]^-{\cong}_-{i} & G_2 \ar[r]_-{\delta} & f(G) \ar[u]_{h}
}\]
in which:
\begin{itemize}
\item
$\sigma$ is a composition of codegeneracy maps, 
\item
$i$ is an isomorphism,
\item
$\delta$ is a composition of inner coface maps, and 
\item
$h$ is a composition of outer coface maps.\index{graphical map!factorization}
\end{itemize}
\end{theorem}

\begin{proof}
By assumption $\nicexy{f(G) \ar[r]^-{h} & K}$ is a subgraph, i.e., a composition of outer coface maps.  So it suffices to show that $g$ decomposes as $\delta \sigma$ up to isomorphism as stated.  If $G =~\uparrow$, then $g$ is an isomorphism.  So we may assume that $|\vertex(G)| \geq 1$.

Suppose $T \subseteq \vertex(G)$ is the subset of vertices $w$ with
\begin{itemize}
\item
precisely one incoming flag and one outgoing flag, and
\item
$f_1(w) =~ \uparrow_e$ for some $e \in \edge(K)$.
\end{itemize}
Define
\[
G_1 = G\left(\{\uparrow\}_{w \in T}\right),
\]
in which an exceptional edge is substituted into each $w \in T$, and a corolla is substituted into each vertex not in $T$.  There is a composition of codegeneracy maps
\[
\nicexy{G \ar[r]^-{\sigma} & G_1.
}\]
The number of codegeneracy maps in $\sigma$ is equal to $|T|$.

Next define
\[
G_2 = f_0(G_1),
\]
which is obtained from $G_1$ by applying $f_0$ to its edges.  This is well-defined because for $w \in T$, its incoming flag and outgoing flag have the same $f_0$-image.  There is an isomorphism
\[
\nicexy{G_1 \ar[r]^-{\cong} & G_2
}\]
given by changing the names of edges using $f_0$.  There is a canonical bijection
\[
\vertex(G_2) \cong \vertex(G) \setminus T.
\]
The image of $G$ can now be constructed from $G_2$ as the graph substitution
\[
f(G) = G_2\left(\{f_1(u)\}_{u \in \vertex(G)\setminus T}\right),
\]
where each $f_1(u) \not=~ \uparrow$ (i.e., $f_1(u)$ has at least one vertex) by the construction of $T$.  The map
\[
\nicexy{G_2 \ar[r]^-{\delta} & f(G)}
\]
is determined by sending each $u \in \vertex(G) \setminus T$ to $f_1(u)$.  It remains to show that $\delta$ is a composition of inner coface maps.

Order $\vertex(G) \setminus T$ arbitrarily as $\{u_1,\ldots,u_k\}$.  Using the associativity and unity of graph substitution, we can construct $f(G)$ from $G_2$ in $k$ steps as follows:
\[
\begin{split}
H_0 &= G_2,\\
H_{i+1} &= H_i\left(f_1(u_{i+1})\right)
\end{split}
\]
for $0 \leq i \leq k-1$.  Then we have
\[
H_k = H_{k-1}\left(f_1(u_k)\right) = f(G). 
\]
The map $\delta$ now factors as follows.
\begin{equation}
\label{deltafactor}
\nicexy{
G_2 = H_0 \ar[d] \ar[rr]^-{\delta} && H_k = f(G)\\
H_1 \ar[r] & \cdots \ar[r] & H_{k-1} \ar[u]
}
\end{equation}
Each map $H_i \to H_{i+1}$ in this diagram is a composition of inner coface maps by Lemma \ref{innercofacefactor}.  Therefore, $\delta$ itself is a composition of inner coface maps.
\end{proof}

\begin{remark}
In Theorem \ref{thm:gupcepimono} we could also have factored $f$ as $h\delta\sigma i$, so the isomorphism $i$ is the first map.  To do this, we instead define
\[
G_1 = f_0(G)
\]
with $i$ the isomorphism induced by applying $f_0$.  Then we define
\[
G_2 = G_1\left(\{\uparrow\}_{w \in T}\right)
\]
with $\nicexy{G_1 \ar[r]^-{\sigma} & G_2}$ the corresponding composition of codegeneracy maps.
\end{remark}

\begin{remark}
In what follows, to simplify the presentation we will often omit mentioning the isomorphism $i$ in a factorization $f = h\delta i \sigma$ as above.
\end{remark}

Next we observe that the codegeneracies-cofaces factorization above is unique.

\begin{lemma}
\label{gupcfactorunique}
Suppose $\nicexy{G \ar[r]^-{f} & K}$ is a properad map between graphical properads.   Up to isomorphism, there exists at most one factorization
\[
\nicexy@C+10pt{
G \ar[d]_{\sigma} \ar[r]^-{f} & K\\
G_1 \ar[r]^-{i} & G_2 \ar[u]_{\partial}
}\]
of $f$, in which
\begin{itemize}
\item
$\sigma$ a composition of codegeneracy maps,
\item
$i$ is an isomorphism, and
\item
$\partial$ is a composition of coface maps.
\end{itemize}
\end{lemma}

\begin{proof}
Suppose $\partial i \sigma = \partial' i' \sigma'$ are two such factorizations of $f$.  On the sets of edges (i.e., colors), $\sigma_0$ is surjective, while $\partial_0 i_0$ is injective.  Therefore, the decomposition
\[
f_0 = (\partial_0 i_0) \sigma_0
\]
is the unique factorization of the function
\[
\nicexy{\edge(G) \ar[r]^-{f_0} & \edge(K)}
\]
into a surjection followed by an injection.  By uniqueness there is a bijection $i'' \colon \edge(G_1) \to \edge(G_1')$ such that the diagram
\[
\nicexy@R+15pt{
& \edge(G_1) \ar[r]^-{i_0}_{\cong} \ar[dd]^-{i''}_-{\cong} 
& \edge(G_2) \ar[dr]^-{\partial_0} &\\
\edge(G) \ar[ur]^-{\sigma_0} \ar[dr]_-{\sigma'_0} 
&&& \edge(K)\\
& \edge(G_1') \ar[r]^-{i'_0}_-{\cong}
& \edge(G_2') \ar[ur]_-{\partial'_0} &
}\]
is commutative. Since both $\sigma$ and $\sigma'$ are compositions of codegeneracy maps, the left commutative triangle  implies that $\sigma'=\sigma$.  In particular, the maps $\partial i$ and $\partial' i'$ have the same source $G_1$, and $\partial_0 i_0 = \partial'_0 i'_0$.

For each vertex $u \in \vertex(G_1)$, we have
\[
\partial_1 i_1(u) = f_1(u) = \partial'_1 i'_1(u).
\]
So the maps $\partial i$ and $\partial' i'$ coincide on colors (i.e., edges) and generators (i.e., vertices).  Therefore, we have $\partial i = \partial' i'$ as well.
\end{proof}

\begin{corollary}
\label{epimonounique}
Each graphical map $f$ admits a decomposition $\partial i \sigma$ in which:
\begin{itemize}
\item
$\sigma$ is a composition of codegeneracy maps, 
\item
$i$ is an isomorphism, and
\item
$\partial$ is a composition of coface maps.
\end{itemize}
Moreover, this decomposition is  unique up to isomorphism.
\end{corollary}

\begin{proof}
Use Theorem \ref{thm:gupcepimono} for the existence of a factorization, and use Lemma \ref{gupcfactorunique} for the uniqueness.
\end{proof}

The next observation says that a graphical map is uniquely determined by its action on edge sets.

\begin{corollary}
\label{edgesetdet}
Suppose $f,f' \colon G \to K$ are graphical maps such that $f_0 = f'_0$ as functions $\edge(G) \to \edge(K)$.  Then $f=f'$.
\end{corollary}

\begin{proof}
Suppose $f=\partial i \sigma$ and $f'=\partial' i' \sigma'$ are codegeneracies-cofaces decompositions, which exist by Theorem \ref{thm:gupcepimono}.  Then we have
\[
(\partial_0 i_0) \sigma_0 = f_0 = f'_0 
= (\partial'_0 i'_0) \sigma'_0
\]
on $\edge(G)$.  As in the proof of Lemma \ref{gupcfactorunique}, we have $\sigma = \sigma'$, that the maps $\partial i$ and $\partial' i'$ have the same source $G_1$, and $\partial_0 i_0 = \partial'_0 i'_0$. Pick a vertex $u \in \vertex(G_1)$.  We have
\[
\partial_1 i_1(u) = f_1(u) \andspace
\partial'_1 i'_1(u) = f'_1(u).
\]
These are subgraphs of $K$ (Corollary \ref{cor1:gupcgraphicalmap}), and they have the same input/output $\edge(K)$-profiles
\[
\binom{f_0 \out(u)}{f_0\inp(u)} =
\binom{f'_0 \out(u)}{f'_0\inp(u)}.
\]
By Lemma \ref{subgraphunique} these subgraphs are equal.  Therefore, we have $\partial i= \partial' i'$ as well.
\end{proof}

\section{A Generalized Reedy Structure on \texorpdfstring{$\varGamma$}{Γ}}
\label{sec:reedygamma}

The simplicial category $\varDelta$ and its opposite $\varDelta^{op}$ are both \emph{Reedy categories}, which allows us to work inductively on diagrams of these shapes. If $\mathcal{M}$ is any Quillen model category, there is an associated model structure on the categories of diagrams $\mathcal{M}^{\varDelta}$ and $\mathcal{M}^{\varDelta^{op}}$ \cite{reedy}. A generalization of this structure was introduced by Berger and Moerdijk in \cite{reedyextension} in order to describe the inductive structure of the dendroidal category $\varOmega$. In this section, we observe that the graphical category $\varGamma$ also possesses a generalized Reedy structure, implying the existence of a Berger-Moerdijk-Reedy model structure on $\mathcal{M}^{\varGamma}$ and $\mathcal{M}^{\varGamma^{op}}$.

Recall that a \textbf{wide subcategory} \index{wide subcategory} of a category $\mathcal{C}$ is a subcategory which contains all objects of $\mathcal{C}$. We will write $\Iso(\mathcal{C})$ \label{note:isoc} for the maximal sub-groupoid of a category $\mathcal{C}$. The following definition appears in \cite{reedyextension}.

\begin{definition}
\label{def:greedy}
A \textbf{generalized Reedy structure} \index{generalized Reedy structure} on a small category $\mathcal{R}$ consists of
\begin{itemize}
\item
wide subcategories $\mathcal{R}^+$ and $\mathcal{R}^-$, and 
\item
a degree function $\deg: \Ob(\mathcal{R}) \to \mathbb{N}$
\end{itemize}
satisfying the following four axioms.
\begin{enumerate}[(i)]
\item 
Non-invertible morphisms in $\mathcal{R}^+$ (resp., $\mathcal{R}^-$) raise (resp., lower) the degree.  Isomorphisms in $\mathcal{R}$ preserve the degree.
\item 
$\mathcal{R}^+ \cap \mathcal{R}^- = \Iso(\mathcal{R})$.
\item
Every morphism $f$ of $\mathcal{R}$ factors as $f = gh$ with $g  \in \mathcal{R}^+$ and $h \in \mathcal{R}^-$, and this
factorization is unique up to isomorphism.
\item 
If $\theta f=f$ for $\theta \in \Iso(\mathcal{R})$ and $f\in \mathcal{R}^-$, then $\theta$ is an identity.
\end{enumerate} If, morever, the condition
\begin{enumerate}[(iv')]
\item If $f \theta=f$ for $\theta \in \Iso(\mathcal{R})$ and $f\in \mathcal{R}^+$, then $\theta$ is an identity
\end{enumerate}
holds, then we call this a generalized \textbf{dualizable} \index{dualizable} Reedy structure.
\end{definition}

We now turn to the generalized Reedy structure on $\varGamma$.

\begin{definition}
Define the \textbf{degree} \index{degree} \label{note:deg} of a graph $G \in \gupc$ to be 
\[ \deg (G) = |\vertex(G)|.\]
Define two wide subcategories of $\varGamma$ as follows, using the fact that graphical maps send vertices to subgraphs (Corollary~\ref{cor1:gupcgraphicalmap}):
\begin{itemize}
\item $\varGamma^+$ \label{note:gammaplus} consists of those maps $f:H\to G$ which are injective on edge sets.
\item $\varGamma^-$ \label{note:gammaminus} consists of those maps $f:H\to G$ which are surjective on edge sets and, for every vertex $v\in \vertex(G)$, there is a vertex $\tilde{v} \in \vertex(H)$ so that $f_1(\tilde{v})$ is a corolla containing $v$.
\end{itemize}
\end{definition}

The following observation gives a characterization of $\varGamma^+$ and $\varGamma^-$.

\begin{lemma}\label{L:pmaltchar}
Suppose $f \colon G \to K \in \varGamma$.  Then:
\begin{enumerate}
\item
$f$ is in $\varGamma^+$ if and only if $f$ is a composition of isomorphisms and coface maps.
\item
$f$ is in $\varGamma^-$ if and only if $f$ is a composition of isomorphisms and codegeneracy maps.
\end{enumerate}
\end{lemma}

\begin{proof}
For the first statement, the ``if" direction is true because each coface map is injective on edge sets.  For the other direction, suppose $f \in \varGamma^+$.  Consider the decomposition
\[
\nicexy@C+10pt{
G \ar[r]^-{\sigma} 
& G_1 \ar[r]^-{h\delta i} 
& K
}\]
of $f$ from Theorem \ref{thm:gupcepimono}, where $\sigma$ is a composition of codegeneracy maps, $h\delta$ is a composition of coface maps, and $i$ is an isomorphism. Since $\sigma$ is a composition of codegeneracy maps, it is surjective on edges. But $\sigma$ is also injective on edges since $f$ is, so $\sigma$ is bijective on edges and is then an identity. Hence $f=h\delta i$.

For the second statement, the ``if" direction is true because each codegeneracy map is in $\varGamma^-$.  For the other direction, suppose $f \in \varGamma^-$.  Consider the same decomposition of $f$ from Theorem \ref{thm:gupcepimono},
\[
\nicexy@C+10pt{
G \ar[r]^-{i\sigma}
& G_2 \ar[r]^-{h\delta}
& K,
}\]
where $h$ is a composition of outer face maps and $\delta$ is a composition of inner face maps. We wish to show that both $h$ and $\delta$ are identities.  If $h$ is not an identity, then there exists an almost isolated vertex $v$ of $K$ such that $f$ factors as
\[
\nicexy{
G \ar[r] & G' \ar[r]^-{d^v} & K,
}\]
where $d^v$ is an outer coface map. Then there is no $\tilde v$ such that $f_1(\tilde v)$ is $C_v$, contradicting the fact that $f$ is in $\varGamma^-$. Hence $h$ is an identity. 

If $\delta$ is not an identity, then $f$ factors as
\[
\nicexy{
G \ar[r] & G' \ar[r]^-{d^v} & K,
}\]
where $d^v$ is an inner coface map. But then $d^v$ is not surjective on edges, so neither is $f$, contradicting the fact that $f\in \varGamma^-$. Thus $\delta$ is an identity as well, and $f=i\sigma$ is a composition of codegeneracy maps followed by an isomorphism.
\end{proof}

\begin{proposition}
The graphical category $\varGamma$ satisfies condition (i) in Definition \ref{def:greedy}.
\end{proposition}

\begin{proof}
It is clear that isomorphisms preserve the degree.  Suppose $f \in \varGamma^+$ is not an isomorphism. By Lemma \ref{L:pmaltchar} $f$ is a composition of isomorphisms and at least one coface map, which must raise the degree.  Therefore, $f$ raises the degree as well.  Likewise, suppose that $f\in \varGamma^-$ is not an isomorphism. By Lemma \ref{L:pmaltchar} $f$ is a composition of isomorphisms and at least one codegeneracy map, which must lower the degree.  Therefore, $f$ also lowers the degree.
\end{proof}

\begin{proposition}
\label{P:cond2}  
The graphical category $\varGamma$ satisfies condition (ii) in Definition \ref{def:greedy}, namely
\[
\varGamma^+ \cap \varGamma^- = \Iso(\varGamma).
\]
\end{proposition}

\begin{proof}
        Inclusion from right to left is obvious. For the reverse, suppose that $f: G\to K$ is in $\varGamma^+ \cap \varGamma^-$. In Lemma \ref{L:pmaltchar} we showed that, since $f\in \varGamma^+$, $f$ admits a factorization
\[
f = \partial i,
\]
where $\partial$ is a composition of coface maps and $i$ is an isomorphism.  In particular, we have
\[
\deg(G) \leq \deg(K).
\]
Since $f\in \varGamma^-$, this same lemma also gives a factorization
\[
f= i' \sigma,
\]
where $\sigma$ is a composition of codegeneracies and $i'$ is an isomorphism.  So we have
\[
\deg(G) \geq \deg(K).
\]
Together with the previous inequality, we have
\[
\deg(G) = \deg(K).
\]
This equality implies that $\partial$ is the identity, so $f=i$ is an isomorphism.
\end{proof}

\begin{proposition}
The graphical category $\varGamma$ satisfies condition (iii) in Definition \ref{def:greedy}.  In other words, every map in $f\in \varGamma$ factors as
\[
f=gh,
\]
where $h\in \varGamma^-$ and $g\in \varGamma^+$, and this factorization is unique up to isomorphism.
\end{proposition}

\begin{proof}
This follows from Lemma \ref{L:pmaltchar} and Corollary \ref{epimonounique}.
\end{proof}

\begin{proposition}
The graphical category $\varGamma$ satisfies conditions (iv) and (iv')  in Definition \ref{def:greedy}.  In other words:
\begin{enumerate}
\item If $f\in \varGamma^-$, $\theta \in \Iso(\varGamma)$, and $\theta f =f$, then $\theta = \Id$. 
\item If $f\in \varGamma^+$, $\theta \in \Iso(\varGamma)$, and $f\theta =f$, then $\theta = \Id$. 
\end{enumerate}
\end{proposition}

\begin{proof}
        Since graphical maps are determined by their actions on edge sets (Corollary \ref{edgesetdet}), 
        it is enough to show that $\theta_0$ is an identity. But now this comes down to the same fact in $\mathbf{Set}$: If $f_0$ is surjective and $\theta_0f_0 = f_0$, then $\theta_0 = \Id$. Similarly, if $f_0$ is injective and $f_0 \theta_0 = f_0$, then $\theta_0 = \Id$.
\end{proof}

\begin{theorem}
The graphical category $\varGamma$ is a dualizable generalized Reedy category.
\end{theorem}

\begin{proof}
Combine the previous four propositions.
\end{proof}


\chapter{Properadic Graphical Sets and Infinity Properads}
\label{ch:graphicalset}

\abstract*{We define the category $\gupcset$ of graphical sets.  There is an adjoint pair
\[
\nicexy@C+10pt{
\gupcset \ar@<.4ex>[r]^-{L} & \properad \ar@<.4ex>[l]^-{N},
}\]
in which the right adjoint $N$ is called the properadic nerve.  The symmetric monoidal product of properads in chapter  \ref{ch:tensor} induces, via the properadic nerve, a symmetric monoidal closed structure on $\gupcset$.  Then we define an  $\infty$-properad as a graphical set in which every inner horn has a filler.  If, furthermore, every inner horn filler is unique, then it is called a strict $\infty$-properad.  The rest of this chapter contains two alternative descriptions of a strict $\infty$-properad.  One description is in terms of the graphical analogs of the Segal maps, and the other is in terms of the properadic nerve.}

The main purposes of this chapter are to define $\infty$-properads and to characterize strict $\infty$-properads.

Infinity properads and strict $\infty$-properads are defined in section \ref{sec:infinityproperad}.  Joyal-Lurie $\infty$-categories (resp., Moerdijk-Weiss $\infty$-operads)  are simplicial sets (resp., dendroidal sets) that satisfy a weaker version of the Kan extension property that only involves inner horns. To define $\infty$-properads, first we define the category $\gupcset$ of properadic graphical sets, which is the presheaf category of the properadic graphical category $\varGamma$ (Definition \ref{def:graphicalcat}).  This makes sense because simplicial sets and dendroidal sets are presheaf categories corresponding to the graphical subcategories $\varDelta \cong \varGamma(\ULin)$ and $\varOmega \simeq \varGamma(\uoperad)$ (Remarks \ref{rk:deltaingamma} and \ref{rk:omegaingamma}).

Similar to the nerve of a small category, there is a properadic nerve $N$, which is part of an adjunction
\[
\nicexy@C+10pt{
\gupcset \ar@<.5ex>[r]^-{L} 
& \properad \ar@<.5ex>[l]^-{N}.
}\]
The symmetric monoidal product in $\properad$ (Theorem \ref{thm:propgmonoidal}) induces, via the properadic nerve, a symmetric monoidal closed structure on $\gupcset$.  Once the graphical analogs of horns are defined, we define an $\infty$-properad as a properadic graphical set in which every inner horn has a filler.  A strict $\infty$-properad has the further property that all the inner horn fillers are unique.

In sections \ref{sec:properadsegal} and  \ref{sec:fundpropsiprop}, we provide two alternative characterizations of strict $\infty$-properads.  Within the category of simplicial sets, it is known that a strict $\infty$-category is isomorphic to the nerve of a category, which in turn is equivalent to a simplicial set whose Segal maps are bijections.  The dendroidal analogs of these statements are also known to be true.  We will prove the graphical analogs of these statements.

In section  \ref{sec:properadsegal} we define the properadic Segal maps of a properadic graphical set.  Then we observe that for the properadic nerve of a properad, the properadic Segal maps are all bijections.  Furthermore, a properadic graphical set whose properadic Segal maps are bijections must be a strict $\infty$-properad.  In section \ref{sec:fundpropsiprop}, we show that every strict $\infty$-properad is, up to isomorphism, a properadic nerve.  Thererfore, a strict $\infty$-properad is precisely a properadic  graphical set whose properadic Segal maps are bijections, which is, up to isomorphism, equivalent to the properadic nerve of a properad.

As in previous chapters, all the properadic concepts in this chapter have obvious analogs for properads with non-empty inputs or non-empty outputs.  Instead of the properadic graphical category $\varGamma$, one uses the full subcategories $\varGamma_i$ or $\varGamma_o$ generated by $\gupci$ or $\gupco$ (Remark \ref{rk:omegaingamma}).  To form the adjunction involving the nerve, instead of the category $\properad$, one uses the categories of properads with non-empty inputs $\properadi$ or non-empty outputs $\properado$.

\section{\texorpdfstring{$\infty$}{∞}-Properads}
\label{sec:infinityproperad}

In this section, we define properadic graphical sets, the properadic nerve functor from properads to properadic graphical sets, the symmetric monoidal closed structure on properadic graphical sets, and (strict) $\infty$-properads.

\subsection{Graphical Sets}
\label{subsec:graphicalset}

Recall from chapter \ref{ch:mapgrproperad} that $\varGamma$ is the properadic graphical category.  It is the non-full subcategory of $\properad$ with objects the graphical properads $\varGamma(G) = F(\ghat)$ for connected wheel-free graphs $G$.  Its morphisms are properadic graphical maps.  Also recall that we often abbreviate $\varGamma(G)$ to just $G$ and drop the adjective \emph{properadic}.

The usual category of simplicial sets is the presheaf category of the finite ordinal category $\varDelta$.  Likewise, the Moerdijk-Weiss category of dendroidal sets is the presheaf category of the dendroidal category $\varOmega$.  Now if we use the graphical category $\varGamma$ instead, then the presheaf category is the desired category of graphical sets.

\begin{definition}
\label{def:gupcgraphicalset}
The diagram category $\gupcset$\label{note:gupcset} is called the category of \textbf{properadic graphical sets},\index{properadic graphical set} \index{graphical set!properadic} or simply \textbf{graphical sets}.
\begin{enumerate}
\item
An object $X \in \gupcset$ is called a \textbf{properadic graphical set}, or simply a \textbf{graphical set}.
\item
For $G \in \gupc$, an element in the set $X(G)$ is called a \textbf{graphex with shape $G$}.\index{graphex}  The plural form of \emph{graphex} is \emph{graphices}.
\item
A graphical set $X$ is \textbf{reduced} \index{graphical set!reduced} if the set $X(C_{\emptyprofh})$ is a singleton, where $C_{\emptyprofh}$ is the single isolated vertex.
\end{enumerate}
\end{definition}

\begin{remark}
If we use the full subcategories $\varGamma(\ULin) \cong \varDelta$ or $\varGamma(\uoperad) \simeq \varOmega$ instead of $\varGamma$ itself, then the definition above gives the categories of simplicial sets or of dendroidal sets \cite{mw1}.  Indeed, the full subcategory inclusions
\[
\nicexy{
\varDelta \ar[r]^-{i} & \varOmega \ar[r]^-{i} & \varGamma
}\]
induce adjoint pairs
\[
\nicexy{
\sset \ar@<.4ex>[r]^-{i_!} & 
\omegaset \ar@<.4ex>[l]^-{i_*} \ar@<.4ex>[r]^-{i_!} &
\gupcset \ar@<.4ex>[l]^-{i_*}.
}\]
The right adjoints $i_*$ are given by restrictions, while the left adjoints $i_!$ are left Kan extensions.  In particular, for $W \in \omegaset$, we have
\[
(i_!W)(G) = 
\begin{cases}
W(G) & \text{ if $G \in \varOmega$},\\
\varnothing & \text{ otherwise},
\end{cases}
\]
and similarly for the first left adjoint.  Since the left adjoints $i_!$ are full and faithful, we may regard $\sset$ as a full subcategory of $\omegaset$, which in turn is regarded as a full subcategory of $\gupcset$.  There are similar adjoint pairs
\begin{equation}
\label{gupcioleftadjoint}
\nicexy{
\gupciset \ar@<.4ex>[r]^-{i_!} & 
\gupcset \ar@<.4ex>[l]^-{i_*} & \text{and} &
\gupcoset \ar@<.4ex>[r]^-{i_!} & 
\gupcset \ar@<.4ex>[l]^-{i_*}
}
\end{equation}
for the full subcategories $\Gammai = \varGamma(\gupci)$ and $\Gammao = \varGamma(\gupco)$.
\end{remark}

\begin{remark}
To be more explicit, a graphical set $X \in \gupcset$ consists of
\begin{itemize}
\item
a set $X(G)$ of graphices of shape $G$ for each connected wheel-free graph $G$, and
\item
a function
\[
\nicexy{
X(H)  \ar[r]^-{f^*} & X(G)
}\]
for each map $\nicexy{G \ar[r]^-{f} & H \in \varGamma}$
\end{itemize}
such that
\[
\Id^* = \Id \andspace (gf)^* = f^*g^*.
\]
A map $\nicexy{X \ar[r]^-{\varphi} & Y \in \gupcset}$ consists of a function $\nicexy{X(G) \ar[r]^-{\varphi_G} & Y(G)}$ for each $G \in \gupc$ such that the diagram
\[
\nicexy{
X(H) \ar[d]_{f^*} \ar[r]^-{\varphi_H} & Y(H) \ar[d]^{f^*}\\
X(G) \ar[r]^-{\varphi_G} & Y(G)
}\]
is commutative for each map $\nicexy{G \ar[r]^-{f} & H \in \varGamma}$.
\end{remark}

\begin{definition}
Suppose $X \in \gupcset$ and $\nicexy{G \ar[r]^-{f} & H \in \varGamma}$.
\begin{enumerate}
\item
If $f$ is an inner/outer coface (resp., codegeneracy) map, then $f^*$ is called an \textbf{inner/outer face} \index{inner face!of properadic graphical set} \index{outer face!of properadic graphical set} (resp., \textbf{degeneracy}) \index{degeneracy!of properadic graphical set} \textbf{map}.
\item
A \textbf{face map} \index{face map!of graphical properad} means either an inner face map or an outer face map.
\end{enumerate}
\end{definition}

\subsection{Why a Graphical Set Resembles a Properad}
\label{rk:gsetproperad}

Here explain why a graphical set $X \in \gupcset$ looks somewhat like a properad, an idea that will eventually lead to our notion of an $\infty$-properad.  The point of the discussion below is to convince the reader that a graphical set is a reasonable candidate for the underlying set of an $\infty$-properad.

We will refer to the biased definition of a properad (Definition \ref{def:biasedproperad}) below.   Let us consider connected wheel-free graphs with few vertices.
\begin{enumerate}
\item
The only $1$-colored connected wheel-free graph with $0$ vertex is the exceptional edge $\uparrow$, so there is a set $X(\uparrow)$, whose elements will be called \textbf{colors} of $X$.
\item
The only $1$-colored connected wheel-free graphs with $1$ vertex are the corollas and their input/output relabeling (Example \ref{ex:permutedcor}), which must induce an isomorphism of component sets of $X$.

Define a \textbf{$1$-dimensional element} \index{$1$-dimensional element} in $X$ as an element in $X(\sigma C \tau)$ for some corolla $C$ and permutations $\sigma$ and $\tau$.  A $1$-dimensional element $f \in X(\sigma C \tau)$, in which $C = C_{(m;n)}$ is the corolla with $m$ inputs and $n$ outputs, is said to have \textbf{$m$ inputs} and \textbf{$n$ outputs}.

From Example \ref{ex:legcoface}, we know that each leg of a permuted corolla $\sigma C\tau$ has a corresponding outer coface map $\uparrow ~ \to \sigma C\tau$.  This yields an outer face map
\[
X(\sigma C\tau) \to X(\uparrow)
\]
corresponding to this leg of $\sigma C\tau$.  Assembling these maps, we then obtain a map
\begin{equation}
\label{gsetprofile}
\nicexy{
X(\sigma C\tau) \ar[d]^-{\rho} \\
\left[\prod_{\inp(\sigma C\tau)} X(\uparrow)\right] 
\times 
\left[\prod_{\out(\sigma C\tau)} X(\uparrow)\right]
}
\end{equation}
For a $1$-dimensional element $x$, we call the image $\rho(x)$ the \textbf{profiles of $x$}. \index{$1$-dimensional element!profiles}

We may reorganize the $1$-dimensional elements in $X$ according to their profiles as follows.  Suppose $\dch$ is an element in $\sS(X(\uparrow))$, i.e., a pair of $X(\uparrow)$-profiles.  Define the component set
\begin{equation}
\label{gsetcobject}
X\dc = \left\{
x \in \coprod_{\sigma, \tau} X(\sigma C\tau) \colon \rho(x) = \dch
\right\}
\end{equation}
where $C$ is the $1$-colored corolla with $|\uc|$ inputs and $|\ud|$ outputs, $\sigma \in \Sigma_{|\ud|}$, and $\tau \in \Sigma_{|\uc|}$.  In other words, $X\dc$ is the set of $1$-dimensional elements in $X$ with profiles $\dc$.  Input and output relabeling among corollas yields isomorphisms
\[
\nicexy@C+10pt{
X\dc \ar[r]^-{(\tau;\sigma)}_-{\cong} & X\dcsigma.
}\]
In particular, the collection of sets
\[
\left\{ X\dc \colon \dc \in \sS(X(\uparrow))\right\}
\]
is a $\Sigma_{\sS(X(\uparrow))}$-bimodule.
\item
For the $1$-colored corolla $C_{(1;1)}$ with $1$ input and $1$ output, there is a codegeneracy map $\nicexy{C_{(1;1)} \ar[r]^-{s} & \uparrow}$.  This yields a degeneracy map
\begin{equation}
\label{gsetunit}
\nicexy{
X(\uparrow) \ar[r]^-{s^*} & X\left(C_{(1;1)}\right).
}
\end{equation}
For a color $c \in X(\uparrow)$, the image
\[
\bone_c \defn s^*(c) \in X(C_{(1;1)})
\]
is called the \textbf{$c$-colored unit} \index{colored unit} of $X$.  Note that $\bone_c$ has profiles $(c;c)$ because the compositions
\[
\nicexy@C+12pt{
\uparrow \ar[r]^-{\text{either}}_-{\text{coface}} 
& C_{(1;1)} \ar[r]^-{s} & \uparrow
}\]
are both equal to the identity.
\item
The only connected wheel-free graphs with $2$ vertices are, up to listing, the partially grafted corollas (Example \ref{ex:pgcor}).  Consider a partially grafted corollas $P$.   There is an inner coface map $C \to P$ (Example \ref{ex0:innercoface}) corresponding to smashing together the two vertices in $P$, which are closest neighbors.  This yields an inner face map
\begin{equation}
\label{gsetproduct}
\nicexy{
X(P) \ar[r]^-{\mu} & X(C),}
\end{equation}
which we would like to think of as a properadic composition.  But of course this is not in general a properadic composition.  

Indeed, $\mu$ need not be a binary product.  To obtain the two underlying elements of an element in $X(P)$, the most natural thing to do is to use the two outer coface maps
\[
C_u \to P
\andspace
C_v \to P
\]
for the two vertices in $P$ (Example \ref{ex:pgcorcoface}), which are both almost isolated.  They yield the map
\begin{equation}
\label{fakeproduct}
\nicexy{
X(P) \ar[r]^-{d} & X(C_u) \times X(C_v).
}
\end{equation}
If $p \in X(P)$ and $d(p) = (x,y)$, then we want to think of $\mu(p) \in X(C)$ as the properadic composition of $x$ and $y$.
\[
\nicexy{
X(P) \ar[dr]^{\mu} \ar[d]_{d} &\\
X(C_u) \times X(C_v) \ar@{.>}[r] & X(C)
}\]
This picture fails to give a binary product because $d$ goes in the wrong direction.
\end{enumerate}

For a properad, the properadic composition is uniquely defined, and the axioms hold strictly.  In particular, we would need the map $d$ to be an isomorphism.  

On the other hand, for $\infty$-properads, in which the properadic composition and the axioms hold only up to homotopy, we do \emph{not} need to insist that $d$ be an isomorphism.  Just like different choices of compositions of paths in a topological space, we only ask that a composition exists without insisting on its uniqueness. In other words, we ask that $d$ be surjective.  So any $d$-preimage followed by $\mu$, the inner face map, will yield such a \emph{weak} properadic composition.  The surjectivity of $d$ will be phrased as part of an extension property that resembles the inner Kan extension property (Example \ref{ex:pgcorhornfiller} below).  In particular, for a partially grafted corollas, this extension property corresponds to the surjectivity of $d$.  For connected wheel-free graphs with more than two vertices, this extension property corresponds to the \emph{weak} properad axioms.

\subsection{Properadic Nerve}

The properadic nerve is a way to associate a graphical set to a properad.  Before we give the precise definition, let us provide some motivation for the properadic nerve functor.

Recall that for a category $\catc$, its \emph{nerve} \index{nerve!categorical} $N\catc$ is a simplicial set that encodes the categorical information of $\catc$.  More precisely, the set of $n$-simplices is
\[
(N\catc)_n = \category([n],\catc).
\]
So a $0$-simplex is just an object of $\catc$, while an $n$-simplex for $n \geq 1$ is a sequence 
\[
\nicexy{
x_0 \ar[r] & x_1 \ar[r] & \cdots \ar[r] & x_{n-1} \ar[r] & x_n
}\]
of $n$ composable maps in $\catc$.  The top and bottom face maps are given by deletion of the last and the first maps, while the other face maps are given by categorical composition.  The degeneracy maps are given by insertion of identity maps.  The nerve construction gives a functor
\[
\nicexy{
\category \ar[r]^-{N} & \sset,
}\]
which by general non-sense has a left adjoint.

Now the category $[n]$ can be represented as the linear graph $L_n \in \varGamma(\ULin)$ with $n$ vertices.  So to define the nerve of a properad, we only need to change $\category(-,\catc)$ to $\properad(-,\sP)$.

\begin{definition}
\label{graphicalnerve}
The \textbf{properadic nerve} \index{properadic nerve} \index{nerve!properadic} is the functor
\[\label{note:nerve}
\nicexy{
\properad \ar[r]^-{N} & \gupcset
}\]
defined by
\[
(N\sP)(G) = \properad(G,\sP)
\]
for $\sP \in \properad$ and $G \in \varGamma$.  In the context of the unique factorization of a graphical map $G \to K \in \varGamma$ (Corollary \ref{epimonounique}), the map
\[
\nicexy{
(N\sP)(K) = \properad(K,\sP) \ar[r] &
\properad(G,\sP) = (N\sP)(G)}
\]
is a composition of maps of the following form:
\begin{itemize}
\item
deletion of an entry (for an outer coface map);
\item
properad structure map of $\sP$ (for an inner coface map);
\item
isomorphism;
\item
colored units of $\sP$ (for a codegeneracy map).
\end{itemize}
\end{definition}

\begin{remark}
To see that $N\sP$ is indeed in $\gupcset$, recall that maps in $\varGamma$ are compositions of coface maps, codegeneracy maps, and isomorphisms.  Coface and codegeneracy maps correspond to graph substitutions.  The properad structure maps of a properad $\sP$ are associative and unital with respect to graph substitutions.
\end{remark}

\begin{lemma}
\label{graphicalnerveadjoint}
The properadic nerve \index{properadic nerve!left adjoint} admits a left adjoint
\[
\nicexy{
\gupcset \ar[r]^-{L} & \properad}
\]
such that the diagram
\[
\nicexy{
\varGamma \ar[r]^-{\iota} \ar[d]_{\mathrm{Yoneda}} & \properad\\
\gupcset \ar[ur]_{L} &}
\]
is commutative up to natural isomorphism.
\end{lemma}

\begin{proof}
By standard category theory, the left adjoint of the properadic nerve is the left Kan extension of the subcategory inclusion $\iota$ along the Yoneda embedding \cite{maclane98} (Ch.X).  It exists because the graphical category $\varGamma$ is small, and $\properad$ has all small colimits.  The left Kan extension makes the diagram commutative up to natural isomorphism because the Yoneda embedding is full and faithful.
\end{proof}

The following observation says that a graphex in $(N\sP)(G)$ is really a $\sP$-decoration of $G$, which consists of a coloring of the edges in $G$ by the colors of $\sP$ and a decoration of each vertex in $G$ by an element in $\sP$ with the corresponding profiles.

\begin{lemma}
\label{npgelement}
Suppose $\sP$ is a $\fC$-colored properad, and $G \in \varGamma$.  Then an element in $(N\sP)(G)$ consists of:
\begin{itemize}
\item
a function $\nicexy{\edge(G) \ar[r]^-{\varphi_0} & \fC}$, and
\item
a function $\varphi_1$ that assigns to each $v \in \vertex(G)$ an element $\varphi_1(v) \in \sP\binom{\varphi_0\out (v)}{\varphi_0 \inp(v)}$.
\end{itemize}
\end{lemma}

\begin{proof}
This is a special case of Lemma \ref{lem:mapfromfree}, since $\varGamma(G)$ is the free properad $F(\ghat)$, where $\ghat$ has color set $\edge(G)$ and element set $\vertex(G)$.
\end{proof}

Next we define the graphical analog of the representable simplicial set $\varDelta[n]$.

\begin{definition}
Suppose $G \in \varGamma$.  The \index{representable graphical set} \textbf{representable graphical set}  $\varGamma[G] \in \gupcset$\label{note:representable} is defined, for $H \in \varGamma$, by
\[
\varGamma[G](H) = \varGamma(H,G),
\]
i.e., the set of graphical maps $H \to G$ (Definition \ref{def:graphicalmap}).
\end{definition}

\begin{remark}
By Yoneda's Lemma a map $\varGamma[G] \to X$ of graphical sets is equivalent to a graphex in the set $X(G)$.
\end{remark}

\subsection{Symmetric Monoidal Closed Structure on Graphical Sets}
\label{sec:smgrsets}

Here we observe that the symmetric monoidal product of properads and the properadic nerve induce a symmetric monoidal closed structure on the category of graphical sets.

First note that each graphical set $X \in \gupcset$ can be expressed up to isomorphism as a colimit of representable graphical sets,
\[
X \cong \colim_{\varGamma[G] \to X} \varGamma[G],
\]
where the colimit is indexed by maps of graphical sets $\varGamma[G] \to X$, i.e., graphices in $X$.

\begin{definition}
Suppose $X$ and $Y$ are graphical sets.
\begin{enumerate}
\item
Define the graphical set
\[
X \otimes Y \defn \colim_{\varGamma[G]\to X, \varGamma[G']\to Y} \left(\varGamma[G] \otimes \varGamma[G']\right),
\]
where
\[
\varGamma[G] \otimes \varGamma[G'] 
\defn
N\left(\varGamma(G) \otimes \varGamma(G')\right)
\]
with $N$ the properadic nerve and $\otimes$ the symmetric monoidal product in $\properad$ (Definition \ref{def:gpropmonoidalproduct}).
\item
Define the graphical set $\Hom(X,Y)$ by
\[
\Hom(X,Y)(G) 
\defn 
\gupcset\left(X \otimes \varGamma[G], Y\right)
\]
for $G \in \varGamma$.
\end{enumerate}
\end{definition}

\begin{remark}
Tensor products of graphical properads $\varGamma(G) \otimes \varGamma(G')$ are described in Corollary \ref{graphicalproptensor}.
\end{remark}

\begin{theorem}
\label{dendroidalclosed}
The category $\gupcset$ of graphical sets is symmetric monoidal closed with monoidal product $\otimes$ and internal hom $\Hom$.\index{graphical set!symmetric monoidal closed structure}
\end{theorem}

\begin{proof}
That $\otimes$ gives a symmetric monoidal product follows from the fact that $\properad$ is symmetric monoidal.  Suppose $X$, $Y$, and $Z$ are graphical sets.  To prove the required adjunction
\[
\gupcset\left(X \otimes Y,Z\right) 
\cong 
\gupcset\left(X, \Hom(Y,Z)\right),
\]
it suffices to consider the case when $X$ is a representable graphical set $\varGamma[G]$.  Then we have natural isomorphisms:
\[
\begin{split}
\gupcset\left(\varGamma[G] \otimes Y,Z\right) 
&\cong \gupcset\left(Y \otimes \varGamma[G],Z\right)\\
&= \Hom(Y,Z)(G)\\
&\cong \gupcset\left(\varGamma[G], \Hom(Y,Z)\right).
\end{split}
\]
\end{proof}

\begin{remark}
The symmetric monoidal closed structure of dendroidal sets is in \cite{mw1,weiss}.
\end{remark}

\subsection{Faces and Horns}

Here we define the graphical analogs of the horns $\Lambda^k[n] \subset \varDelta[n]$, which we will soon use to define an $\infty$-properad.  The formalism here is almost exactly the same as in the category of simplicial sets.

\begin{definition}
Suppose $G \in \varGamma$.  A \textbf{face of $G$} is a coface map in $\varGamma$ whose target is $G$.  An \textbf{inner/outer face of $G$} \index{inner face!of graphical properad} \index{outer face!of graphical properad} is an inner/outer coface map whose target is $G$.
\end{definition}

\begin{definition}
Suppose $X$ is a graphical set.  A \textbf{graphical subset} \index{graphical subset} of $X$ is a graphical set $W$ that is equipped with a map $W \to X$ of graphical sets such that each component map
\[
W(G) \to X(G)
\]
for $G \in \varGamma$ is a subset inclusion.
\end{definition}

\begin{definition}
Suppose $G \in \varGamma$, and $\nicexy{K \ar[r]^{d} & G}$ is a face of $G$.
\begin{enumerate}
\item
The \textbf{$d$-face} of $\varGamma[G] \in \gupcset$ is the graphical subset $\varGamma^d[G]$ defined by
\[\label{note:dface}
\begin{split}
& \varGamma^d[G](J) \\
&= 
\left\{\text{composition of }
\nicexy{J \ar[r]^-{f} & K \ar[r]^-{d} & G}
\text{ with $f \in \varGamma[K](J)$}\right\}
\end{split}
\]
for $J \in \varGamma$.
\item
The \textbf{$d$-horn} \index{horn!of properadic graphical set} of $\varGamma[G]$ is the graphical subset $\Lambda^d[G]$ defined by
\[\label{note:dhorn}
\Lambda^d[G](J) 
= 
\bigcup_{\text{faces $d'\not=d$}} \varGamma^{d'}[G] (J),
\]
where the union is indexed by all the faces of $G$ \emph{except} $d$.  Write
\[
\nicexy{\Lambda^d[G] \ar[r]^-{i} & \varGamma[G]}
\]
for the graphical subset inclusion.
\item
A \textbf{horn}  of $\varGamma[G]$ is a $d$-horn for some coface map $d$.  An \textbf{inner horn} \index{inner horn!of properadic graphical set} is a $d$-horn in which $d$ is an inner coface map.
\item
Given  $X \in \gupcset$, a \textbf{horn of $X$} is a map
\[
\nicexy{\Lambda^d[G] \ar[r] & X}
\]
of graphical sets.  It is an \textbf{inner horn of $X$} if $\Lambda^d[G]$ is an inner horn.
\end{enumerate}
\end{definition}

\begin{example}
Restricted to the full subcategories $\varGamma(\ULin) \cong \varDelta$  and $\varGamma(\utree) \simeq \varOmega$ (Remarks \ref{rk:deltaingamma} and \ref{rk:omegaingamma}), the above definitions become the ones in simplicial sets and in dendroidal sets.
\end{example}

\begin{notation}
If there is a face $\nicexy{K \ar[r] & G}$ of $G$ denoted by $d^u$ or $d_u$, we will sometimes abbreviate the face $\varGamma^{d^u}[G]$ or $\varGamma^{d_u}[G]$ to $\varGamma^u[G]$.  The same goes for horns.
\end{notation}

\begin{remark}
\begin{enumerate}
\item
If $\nicexy{K \ar[r]^{d} & G}$ is a face of $G$, then the $d$-face $\varGamma^d[G]$ is the image of the induced map $\nicexy{\varGamma[K] \ar[r]^{d \circ (-)} & \varGamma[G]}$ of graphical sets.
\item
In the definition of the $d$-horn $\Lambda^d[G]$, we used the same convention as before about ignoring listings (Convention \ref{graphconvention}).  In other words, when we say $d$ is excluded, we mean $d$ and every coface map $K' \to G$ obtainable from $d$ by changing the listing are all excluded.
\item
If there is an \emph{inner} horn $\Lambda^d[G]$ for some inner coface map $\nicexy{K \ar[r]^-{d} & G}$, then $G$ has at least two vertices because there is an inner properadic factorization (Definition \ref{def:ipropfact}) of $G$,
\[
G = K(H_w),
\]
in which the distinguished subgraph $H_w$ is a partially grafted corollas (Example \ref{ex:pgcor}).
\end{enumerate}
\end{remark}

The following observation gives a more explicit description of a horn of a graphical set.

\begin{lemma}
\label{horndescription}
Suppose $X \in \gupcset$, and $\nicexy{K \ar[r]^{d} & G}$ is a face of $G \in \varGamma$.  Then a horn
\[
\nicexy{\Lambda^d[G] \ar[r]^-{f} & X}
\]
of $X$ \index{inner horn!characterization} is equivalent to a collection of maps
\[
\left\{\nicexy{\varGamma[H] \ar[r]^-{f_H} & X} \colon 
\nicexy{H \ar[r]^{d'} & G} \text{ face $\not= d$}
\right\}
\]
such that, if
\[
\nicexy{
J \ar[r]^-{a^1} \ar[d]_-{a^2} 
& H^1 \ar[d]^{d^1}\\
H^2 \ar[r]^-{d^2} & G
}\]
is a commutative diagram of coface maps with each $d^i \not= d$, then the diagram
\[
\nicexy{
\varGamma[J] \ar[r]^-{a^1} \ar[d]_-{a^2} 
& \varGamma[H^1] \ar[d]^{f_{H_1}}\\
\varGamma[H^2] \ar[r]^-{f_{H_2}} & X
}\]
is also commutative.
\end{lemma}

\begin{proof}
Given a horn $\nicexy{\Lambda^d[G] \ar[r]^-{f} & X}$ of $X$, the map $\nicexy{\varGamma[H] \ar[r]^-{f_H} & X}$ corresponding to a face $\nicexy{H \ar[r]^{d'} & G}$ not equal to $d$ is the image of $\Id_H$ under the  composition
\[
\nicexy{
\varGamma^{d'}[G](H) \ar[r] & \Lambda^d[G](H) \ar[r]^-{f} & X(H).
}\]
The last diagram in the statement of the lemma is commutative because the $d$-horn $\Lambda^d[G]$ is the union of the $d'$-faces inside $\varGamma[G]$.
\end{proof}

\begin{remark}
\label{rk:hornofx}
The collection of maps $\{f_H\}$ in Lemma \ref{horndescription} is equivalent to a collection of graphices $\left\{f_H \in X(H)\right\}$ such that, if the coface square is commutative, then
\[
(a^1)^*\left(f_{H_1}\right) 
= 
(a^2)^*\left(f_{H_2}\right) \in X(J).
\]
In other words, a horn of $X$ is really a collection of graphices in $X$, one for each face not equal to the given one, that agree on common faces.
\end{remark}

\begin{example}
\label{ex:pgcorhorn}
Suppose $X \in \gupcset$, $P$ is a partially grafted corollas with top (resp., bottom) vertex $v$ (resp., $u$), and $\nicexy{C \ar[r]^-{d} & P}$ is the unique \emph{inner} coface map corresponding to the closest neighbors $u$ and $v$ in $P$.  There are only two other faces of $P$:
\begin{itemize}
\item
the \emph{outer} coface map $\nicexy{C_u \ar[r]^-{d^v} & P}$ corresponding to the almost isolated vertex $v$, and
\item
the \emph{outer} coface map $\nicexy{C_v \ar[r]^-{d^u} & P}$ corresponding to the almost isolated vertex $u$.
\end{itemize}
The common faces of the corollas $C_u$ and $C_v$ are exactly the ordinary edges in $P$, since the only coface maps into a corolla are the outer coface maps corresponding to its legs.  

Therefore, an inner horn $\nicexy{\Lambda^d[P] \ar[r]^-{f} & X}$ is equivalent to a pair of elements
\begin{equation}
\label{fufv}
(f_u, f_v) \in X(C_u) \times X(C_v),
\end{equation}
such that if
\[
\nicexy@C+10pt{
\uparrow \ar[d]_{d^e_u} \ar[r]^-{d^e_v} & C_v\\
C_u &
}\]
are the outer coface maps corresponding to an ordinary edge $e$ in $P$, then
\[
\left(d^e_u\right)^*(f_u) 
= 
\left(d^e_v\right)^*(f_v) \in X(\uparrow).
\]
In other words, such an inner horn of $X$ is a pair of $1$-dimensional elements, one for each vertex in $P$, whose profiles match along ordinary edges in $P$.  We may, therefore, visualize such an inner horn of $X$ as the following $X(\uparrow)$-colored partially grafted corollas
\begin{center}
\begin{tikzpicture}
\matrix[row sep=.2cm,column sep=1.5cm] {
\node [fatplain,label=above:$...$] (p1) {$f_v$}; \\
\node [empty]  {...}; \\
\node [fatplain,label=below:$...$] (p2) {$f_u$}; \\
};
\foreach \x in {1,2}
{
\draw [outputleg] (p\x) to +(-.6cm,.5cm);
\draw [outputleg] (p\x) to +(.6cm,.5cm);
\draw [inputleg] (p\x) to +(-.6cm,-.5cm);
\draw [inputleg] (p\x) to +(.6cm,-.5cm);
}
\draw [arrow,bend left=40] (p2) to (p1);
\draw [arrow,bend right=40] (p2) to (p1);
\end{tikzpicture}
\end{center}
decorated by $1$-dimensional elements in $X$.
\end{example}

\subsection{\texorpdfstring{$\infty$}{∞}-Properads}

We now define $\infty$-properads as graphical sets satisfying an inner horn extension property.

\begin{definition}
\label{def:infinityproperad}
Suppose $X \in \gupcset$.
\begin{enumerate}
\item
We call $X$ an \textbf{$\infty$-properad} \index{infinity properad} if for each inner horn $f$ of $X$,
\begin{equation}
\label{innerhornfiller}
\nicearrow
\xymatrix@C+12pt{
\Lambda^d[G] \ar[d]_{i} \ar[r]^-{f} & X\\
\varGamma[G] \ar@{.>}[ur]
}
\end{equation}
a dotted arrow, called an \textbf{inner horn filler},\index{inner horn filler} exists and makes the triangle commutative.
\item
A \textbf{strict $\infty$-properad} \index{strict infinity properad} is an $\infty$-properad for which each inner horn filler in \eqref{innerhornfiller} is unique.
\end{enumerate}
\end{definition}

\begin{remark}
For an $\infty$-properad, we only ask that \emph{inner} horns have fillers.  Also, uniqueness of an inner horn filler is not required, unless we are dealing with a \emph{strict} $\infty$-properad.
\end{remark}

\begin{example}
\label{ex:pgcorhornfiller}
Suppose $P$ is a partially grafted corollas.  Recall from \eqref{fufv} that an inner horn
\[
\nicexy{\Lambda^d[P] \ar[r]^-{f} & X}
\]
is exactly a pair of $1$-dimensional elements $(f_u,f_v)$, one for each vertex in $P$, with matching profiles along ordinary edges in $P$.  Then an inner horn filler for $f$ is exactly a $d$-preimage of $(f_u,f_v)$, where $d$ is the map in \eqref{fakeproduct}.
\end{example}

Recall from Remark \ref{rk:omegaingamma} that $\Gammai$ (resp., $\Gammao$) is the full subcategory of the properadic graphical category $\varGamma$ with objects in $\gupci$ (resp., $\gupco$), i.e., connected wheel-free graphs with non-empty inputs (resp., non-empty outputs).

\begin{definition}
\label{def:infpropnonempty}
A graphical set $X \in \gupciset$ (resp., $\gupcoset$) is an \textbf{$\infty$-properad with non-empty inputs} (resp., \textbf{non-empty outputs}) \index{infinity properad!with non-empty inputs} \index{infinity properad!with non-empty outputs} if it satisfies the inner horn extension property \eqref{innerhornfiller} for each $G \in \gupci$ (resp., $G \in \gupco$).
\end{definition}

\section{Properadic Segal Maps}
\label{sec:properadsegal}

The purpose of this section is to study the graphical analogs of the Segal maps.  We observe that the properadic nerve of a properad satisfies the graphical version of the Segal condition.  Moreover, a graphical set that satisfies the Segal condition must be a strict $\infty$-properad.  In section \ref{sec:fundpropsiprop} we will use these observations to characterize strict $\infty$-properads as  properadic nerves of properads up to isomorphisms, which in turn are equivalent to graphical sets that satisfy the Segal condition.

\subsection{Motivation of the Properadic Segal Maps}

To motivate the construction, first consider the simplicial Segal maps.

Fix a simplicial set $X$.  For each $n \geq 2$, there are maps $X_n \to X_1$, each one being a composition of top and bottom face maps.  These maps are coming from the maps in the finite ordinal  category $\varDelta \cong \varOmega(\ULin)$,
\[
\nicearrow
\xymatrix{
[1] = \{0,1\} \ar[d]_-{\xi_i} & 0 \ar@{|->}[d] & 1 \ar@{|->}[d]\\
[n]=\{0,1,\ldots,n\} & i & i+1
}\]
for $0 \leq i \leq n-1$.  The maps $X_n \to X_1$ fit together to form the \emph{Segal maps}\index{Segal map}
\[
X_n \to \underbrace{X_1^0 \times_{X_0} \cdots \times_{X_0} X_1^{n-1}}_{\text{$n$ factors of $X_1$}},
\]
where each $X_1^i$ is a copy of $X_1$.  The right-hand side is the limit of the diagram consisting of the top and bottom face maps
\[
\nicearrow
\xymatrix{
& X_1^{i+1} \ar[d]^{d_1}\\
X_1^i \ar[r]^-{d_0} & X_0
}\]
as $i$ runs through $\{0,\ldots,n-1\}$.  The Segal maps are well-defined because, in $\varDelta$, the diagram
\begin{equation}
\label{zeroton}
\nicearrow
\xymatrix{
[0] \ar[r]^-{d^0} \ar[d]_{d^1} & [1] \ar[d]^{\xi_i}\\
[1] \ar[r]^-{\xi_{i+1}} & [n]
}
\end{equation}
is commutative, as $0 \in [0]$ is sent to $i+1$ in both cases.

Now think about the square \eqref{zeroton} in terms of linear graphs.  The linear graph $L_n$ with $n$ vertices looks like:
\begin{center}
\begin{tikzpicture}
\matrix[row sep=.6cm,column sep=1cm] {
\node [plain] (vn) {$v_n$}; \\
\node [empty] (vdot) {$\vdots$}; \\
\node [plain] (v2) {$v_2$};\\
\node [plain] (v1) {$v_1$};\\
};
\draw [inputleg] (v1) to node[swap]{\footnotesize{$0$}} +(0,-.9cm);
\draw [arrow] (v1) to node{\footnotesize{$1$}} (v2);
\draw [arrow] (v2) to node{\footnotesize{$2$}} (vdot);
\draw [arrow] (vdot) to node{\footnotesize{$n-1$}} (vn);
\draw [outputleg] (vn) to node{\footnotesize{$n$}} +(0,.9cm);
\end{tikzpicture}
\end{center}
Of course, $L_0$ is the exceptional edge $\uparrow$.  In $\varDelta \cong \varGamma(\ULin)$ the top and bottom coface maps are exactly our outer coface maps, which delete either the top vertex or the bottom vertex together with the attached leg.  By repeatedly deleting top and bottom vertices from $L_n$ (i.e., iterating outer coface maps), one obtains the copy 
\begin{center}
\begin{tikzpicture}
\matrix[row sep=.6cm,column sep=1cm] {
\node [plain] (v) {$v_{i+1}$}; \\
};
\draw [inputleg] (v) to node[swap]{\footnotesize{$i$}} +(0,-.9cm);
\draw [outputleg] (v) to node{\footnotesize{$i+1$}} +(0,.9cm);
\end{tikzpicture}
\end{center}
of the corolla $L_1$.  There are many ways to get to this $L_1$, but their compositions are all equal.

In other words, the map $\nicexy{L_1 \ar[r]^-{\xi_i} & L_n}$ above is the composition of outer coface maps, starting at some corolla corresponding to a vertex in $L_n$, ending at $L_n$ itself.  The commutative diagram \eqref{zeroton} says that the unique ordinary edge in
\begin{center}
\begin{tikzpicture}
\matrix[row sep=.6cm,column sep=1cm] {
\node [plain] (vi2) {$v_{i+2}$};\\
\node [plain] (vi1) {$v_{i+1}$};\\
};
\draw [inputleg] (vi1) to node[swap]{\footnotesize{$i$}} +(0,-.9cm);
\draw [arrow] (vi1) to node{\footnotesize{$i+1$}} (vi2);
\draw [outputleg] (vi2) to node{\footnotesize{$i+2$}} +(0,.9cm);
\end{tikzpicture}
\end{center}
is both the \emph{output} leg of $v_{i+1}$ and the \emph{input} leg of $v_{i+2}$.

The plan for the graphical case is now as follows. The analog of $L_n \in \varDelta$ with $n \geq 2$ is a connected wheel-free graph $G$ with $n \geq 2$ vertices.  Each copy of $L_1 \in \varDelta$ within $L_n$ corresponds to a corolla $C_v$ with the same profiles as a vertex $v \in \vertex(G)$.  The map $\xi_i$ corresponds to a composition of outer coface maps $\nicexy{C_v \ar[r]^-{\xi_v} & G}$.  There are many such compositions of outer coface maps, but their compositions should all be equal.

Suppose there is an internal edge $\nicexy{u \ar[r]^-{e} & v}$ in $G$. The commutative diagram \eqref{zeroton} corresponds to the commutative square
\[
\nicearrow
\xymatrix{
\uparrow \ar[r] \ar[d]
& C_u \ar[d]^{\xi_u}\\
C_v \ar[r]^-{\xi_v} & G
}\]
in $\varGamma$.  The top map is the outer coface map corresponding to the \emph{output} leg $e \in \edge(C_u)$, while the left vertical map is the outer coface map corresponding to the \emph{input} leg $e \in \edge(C_v)$.  Observe two things:
\begin{enumerate}
\item
There may be finitely many ordinary edges from $u$ to $v$ in $G$, since $G$ may not be simply connected.  So there is one commutative square as above for each ordinary edge from $u$ to $v$.  This phenomenon does not occur for simplicial sets and dendroidal sets, since these objects correspond to linear graphs and unital trees, which are simply connected.
\item
A given vertex $u \in G$ can have ordinary edges to finitely many different vertices.  Again this phenomen does not occur for simplicial or dendroidal sets, since in a unital tree each vertex has precisely one output.
\end{enumerate}
What these two observations mean is that the targets of the properadic Segal maps are more complicated limits than in the cases of small categories and operads.

\subsection{Outer Coface Maps from Corollas}

Before we can define the properadic Segal maps, first we need to know that it is possible to connect each corolla in each graph to the graph itself using outer coface maps.  This is the graphical analog of the maps $\nicexy{L_1 \ar[r]^-{\xi} & L_n}$ above.

For a vertex $v$ in a graph $G$, recall that the notation $C_v$ denotes the corolla with the same profiles as those of $v \in \vertex(G)$.

\begin{lemma}
\label{cvtog}
Suppose $G \in \gupc$ has $n \geq 2$ vertices with $v \in \vertex(G)$.  Then there are $n-1$ composable outer coface maps
\[
\nicearrow
\xymatrix{
C_v \ar[r] & G_{n-2} \ar[r] & \cdots \ar[r] & G_1 \ar[r] & G,
}\]
whose composition sends:
\begin{itemize}
\item
the legs of $C_v$ to the corresponding edges adjacent to $v \in G$, and
\item
the unique vertex of $C_v$ to the $\ghat$-decorated corolla $C_v$.
\end{itemize}
Moreover, the composition of such outer coface maps is unique.
\end{lemma}

\begin{proof}
This follows from the second part of Corollary \ref{cor:gupcsubgraph}, since a subgraph decomposes into outer coface maps, each one increasing the number of vertices by one.
\end{proof}

\begin{definition}
\label{def:xiv}
Suppose $G \in \gupc$ with $v \in \vertex(G)$. Define
\[\label{note:corinclusion}
\nicearrow
\xymatrix{
C_v \ar[r]^-{\xi_v} & G \in \varGamma
}\]
to be:
\begin{enumerate}
\item
the unique composition of outer coface maps in Lemma \ref{cvtog} if $|\vertex(G)| \geq 2$, and
\item
the isomorphism induced by input/output relabeling $C_v \to \sigma C_v \tau$ if $G$ has $1$ vertex (so $G$ is a permuted corolla $\sigma C_v \tau$).
\end{enumerate}
\end{definition}

\begin{example}
Consider the connected wheel-free graph $G$
\begin{center}
\begin{tikzpicture}
\matrix[row sep=.4cm, column sep=.5cm]{
& \node [plain] (x) {$x$}; & \\
\node [plain] (v) {$v$}; && \node [plain] (w) {$w$};\\
& \node [plain] (u) {$u$}; &\\
};
\draw [arrow] (u) to (v);
\draw [arrow] (u) to (w);
\draw [arrow] (v) to (x);
\draw [arrow] (w) to (x);
\end{tikzpicture}
\end{center}
with four vertices and four ordinary edges.  Then there are exactly two decompositions of $\nicexy{C_v \ar[r]^-{\xi_v} & G}$ into outer coface maps, which yield the following commutative diagram.
\[
\nicearrow
\xymatrix{
& H_w \ar[r]^-{d^w} & H_x \ar[dr]^-{d^x} &\\
C_v \ar[ur]^-{d^u} \ar[dr]_-{d^x} &&& G\\
& H_w \ar[r]^-{d^w} & H_u \ar[ur]_-{d^u} &
}\]
For the top half of this diagram, the outer coface map $\nicexy{
H_x \ar[r]^-{d^x} & G}$ means $H_x$ is obtained from $G$ by deleting the almost isolated vertex $x$.  So $H_x$ is
\begin{center}
\begin{tikzpicture}
\matrix[row sep=.4cm, column sep=.5cm]{
\node [plain] (v) {$v$}; && \node [plain] (w) {$w$};\\
& \node [plain] (u) {$u$}; &\\
};
\draw [arrow] (u) to (v);
\draw [arrow] (u) to (w);
\draw [outputleg] (v) to +(0,.7cm);
\draw [outputleg] (w) to +(0,.7cm);
\end{tikzpicture}
\end{center}
which has $w$ as an almost isolated vertex.  Likewise,  $\nicexy{
H_w \ar[r]^-{d^w} & H_x}$ is the corresponding outer coface maps.  This $H_w$ is
\begin{center}
\begin{tikzpicture}
\matrix[row sep=.4cm, column sep=.5cm]{
\node [plain] (v) {$v$}; &\\
& \node [plain] (u) {$u$}; \\
};
\draw [arrow] (u) to (v);
\draw [outputleg] (u) to +(.7cm,.6cm);
\draw [outputleg] (v) to +(0,.7cm);
\end{tikzpicture}
\end{center}
which has $u$ as an almost isolated vertex.  The map $\nicexy{
C_v \ar[r]^-{d^u} & H_w}$ is the corresponding outer coface map.  The lower half of the above commutative diagram is interpreted similarly, starting with the deletion of the almost isolated vertex $u \in G$. 
\end{example}

\subsection{Compatibility over Ordinary Edges}

To define the properadic Segal maps, we need to describe how two corollas $C_u$ and $C_u$ in $G \in \gupc$ are compatible over an ordinary edge adjacent to them.  This is about outer coface maps into corollas of the form $\uparrow ~ \to C$, one for each leg of $C$.  If a leg of a corolla $C$ already has a name, say $e$, then we also write $\inp_e$ (resp., $\out_e$) for this outer coface map when $e$ is an input (resp., output) leg of $C$.  Furthermore, if $e$ is clear from the context, we will omit the subscript $e$.

\begin{lemma}
\label{uparrowtog}
Suppose $\nicexy{u \ar[r]^-{e} & v}$ is an ordinary edge in $G \in \gupc$.  Then the square
\[
\nicearrow
\xymatrix{
\uparrow \ar[r]^-{\out_e} \ar[d]_{\inp_e} & C_u \ar[d]^{\xi_u}\\
C_v \ar[r]^-{\xi_v} & G
}\]
in $\varGamma$ is commutative.
\end{lemma}

\begin{proof}
Both compositions in the square send the unique element in $\uparrow$ to the $\edge(G)$-colored exceptional edge $\uparrow_e$.
\end{proof}

\subsection{Properadic Segal Maps}

Lemma \ref{uparrowtog} ensures that the properadic Segal maps below are well-defined.  The properadic Segal maps are basically  comparison maps from each component of a graphical set to the limit determined by its $1$-dimensional components.  As we will show below, a graphical set is a strict $\infty$-properad precisely when the properadic Segal maps are all bijections.

\begin{definition}
Suppose $\sK \in \gupcset$, and $G \in \gupc$ has at least two vertices.
\begin{enumerate}
\item
Define the set
\[\label{note:corribbon}
\sK(G)_1 \defn \left(\prod_{v \in \vertex(G)} \sK(C_v)\right)_{\sK(\uparrow)}
\]
as the limit of the diagram consisting of the maps
\begin{equation}
\label{kcuv}
\nicearrow
\xymatrix@C+12pt{
& \sK(C_u) \ar[d]^{\sK(\out_e)}\\
\sK(C_v) \ar[r]^-{\sK(\inp_e)} & \sK(\uparrow)
}
\end{equation}
as $\nicearrow\xymatrix{u \ar[r]^-{e} & v}$ runs through all the ordinary edges in $G$. Call $\sK(G)_1$ the \textbf{corolla ribbon} \index{corolla ribbon!for properadic graphical set} of $\sK(G)$.
\item
Define the  \textbf{properadic Segal map} \index{properadic Segal map} \index{Segal map!properadic}
\begin{equation}
\label{segalmap}
\nicearrow
\xymatrix@C+10pt{
\sK(G) \ar[r]^-{\chi_G} & \sK(G)_1
}
\end{equation}
as the unique map induced by the maps
\[
\nicexy{\sK(G) \ar[r]^-{\sK(\xi_v)} & \sK(C_v)}
\]
as $v$ runs through all the vertices in $G$.  The properadic Segal map is well-defined by Lemma \ref{uparrowtog}.
\item
We say $\sK$ satisfies the \textbf{properadic Segal condition} \index{Segal condition!properadic} \index{properadic Segal condition} if the properadic Segal map $\chi_G$ is a bijection for every $G \in \gupc$ with at least two vertices. 
\end{enumerate}
\end{definition}

As usual we will often drop the adjective \emph{properadic} if there is no danger of confusion.

\begin{remark}
As long as we are working with connected wheel-free graphs, there is no need to consider the properadic Segal map for $G$ with exactly one vertex.  Indeed, such a $G$ must be a permuted corolla $\sigma C_v \tau$, so the corolla ribbon is $\sK(C_v)$.  The Segal map is the isomorphism
\[
\nicearrow
\xymatrix@C+12pt{
\sK(\sigma C_v \tau) \ar[r]^-{\sK(\xi_v)}_-{\cong} & \sK(C_v)
}\]
induced by input/output relabeling.  Therefore, whenever we mention the properadic Segal map, we automatically assume that $G$ has at least two vertices.
\end{remark}

\begin{remark}
\label{rk:ribbongeometric}
Let us explain the geometric meaning of the corolla ribbon $\sK(G)_1$.  An element $\theta \in \sK(G)_1$ is equivalent to the data:
\begin{itemize}
\item
a function $\nicexy{\edge(G) \ar[r]^-{\theta_0} & \sK(\uparrow)}$, and
\item
a function $\theta_1$ that assigns to each vertex $v \in \vertex(G)$ a $1$-dimensional element $\theta_1(v) \in \sK(C_v)$ with profiles corresponding to those of $v$ under $\theta_0$. 
\end{itemize}
In other words, $\theta$ is a $\sK(\uparrow)$-colored decoration of $G$, in which vertices are decorated by $1$-dimensional elements in $\sK$.
\end{remark}

\begin{example}
Consider the following connected wheel-free graph $G$.
\begin{center}
\begin{tikzpicture}
\matrix[row sep=1cm, column sep=1cm]{
\node [plain] (w) {$w$};\\
\node [plain] (v) {$v$};\\
\node [plain] (u) {$u$};\\
};
\draw [arrow] (v) to node{$c$} (w);
\draw [arrow] (u) to node[swap]{$b$} (v);
\draw [arrow, bend left=55] (u) to node{$a$} (v);
\draw [arrow, bend right=60] (u) to node[swap]{$d$} (w);
\end{tikzpicture}
\end{center}
There are three vertices and four ordinary edges.  Each vertex may have input and output legs of $G$, which are not depicted in the picture.  Then there is a commutative diagram
\[
\nicearrow
\xymatrix{
&& \uparrow_d \ar@/_3pc/[ddll]_{\out} \ar@/^3pc/[ddrr]^{\inp} &&\\
& \uparrow_a \ar[dl]_{\out} \ar[dr]^{\inp} && \uparrow_c \ar[dl]_{\out} \ar[dr]^{\inp} & \\
C_u \ar[drr]_{\xi_u} & \uparrow_b \ar[l]_{\out} \ar[r]^{\inp} & C_v \ar[d]^{\xi_v} && C_w \ar[dll]^{\xi_w}\\
&& G &&
}\]
in $\varGamma$ by Lemma \ref{uparrowtog}.  In this case, for a graphical set $\sK$, the corolla ribbon $\sK(G)_1$ is the limit of the part of this diagram above $G$, after applying $\sK$.  The properadic Segal map $\chi_G$ is induced by the three maps $\xi_*$.
\end{example}

\subsection{Alternative Description of the Properadic Segal Maps}

There is another way to define the properadic Segal maps.  Using the graphical set defined below, we will see that each properadic Segal map is represented by the corresponding properadic Segal core inclusion.

\begin{definition}
\label{def:segalcore}
For $\uparrow ~\not= G \in \gupc$, define its \textbf{properadic Segal core} \index{properadic Segal core} \index{Segal core!properadic} $\Sc[G] \in \gupcset$\label{note:segalcore} as the graphical subset
\[
\Sc[G] 
= \bigcup_{v \in \vertex(G)} \image 
\left(\nicexy{\varGamma[C_v] \ar[r]^-{\xi_v} & \varGamma[G]}\right) \subseteq \varGamma[G].
\]
Denote by $\nicexy{\Sc[G] \ar[r]^-{\epsilon_G} & \varGamma[G]}$\label{note:epsilong} the inclusion, and call it the \textbf{Segal core inclusion}.\index{Segal core!inclusion}
\end{definition}

\begin{lemma}
\label{lem:segalcoredesc}
Suppose $\sK \in \gupcset$, and $G \in \gupc$.  Then there is a commutative diagram
\[
\nicearrow
\xymatrix{
\sK(G) \ar[dd]_{\cong} \ar[rr]^-{\chi_G} && \sK(G)_1 \ar[dd]^{\cong}\\
&&\\
\gupcset\left(\varGamma[G], \sK\right) \ar[rr]^-{\epsilon_G^*} && \gupcset\left(\Sc[G], \sK\right) 
}\]
that is natural in $G$.
\end{lemma}

\begin{proof}
The two vertical bijections are by Yoneda's Lemma.  The commutativity of the diagram is by the definition of the Segal map $\chi_G$.
\end{proof}

\begin{remark}
Using Lemma \ref{lem:segalcoredesc}, one can also think of the properadic Segal map $\chi_G$ as the map $\epsilon_G^*$.
\end{remark}

\subsection{Properadic Nerves Satisfy the Segal Condition}

In the rest of this section, we connect the Segal condition with the properadic nerve   (Definition \ref{graphicalnerve})
\[
\nicexy{
\properad \ar[r]^-{N} & \gupcset
}\]
and strict $\infty$-properads.  Recall that for a properad $\sP$, the value of $N\sP$ at $G \in \varGamma$ is the set
\[
(N\sP)(G) = \properad(G,\sP).
\]
The elements in $(N\sP)(G)$ are described in Lemma \ref{npgelement}.  Briefly, each element in $(N\sP)(G)$ is a $\sP$-decoration of the edges and vertices in $G$.

The following observation says that the properadic nerve always yields a graphical set that satisfies the Segal condition.

\begin{lemma}
\label{lem:nerveissegal}
The properadic nerve of every properad satisfies the Segal condition.
\end{lemma}

\begin{proof}
Suppose $G \in \gupc$ has at least two vertices, and $\sP$ is a $\fC$-colored properad.  Since $G$ is ordinary connected wheel-free, we have
\[
\edge(G) = \bigcup_{v\in\vertex(G)} \edge(C_v) \andspace 
\vertex(G) = \coprod_{v \in \vertex(G)} \vertex(C_v).
\]
In the union on the left, we have
\[
e \in \edge(C_u) \bigcap \edge(C_v)
\]
precisely when $e$ is an ordinary edge adjacent to the vertices $u$ and $v$.  In particular, a function $\nicexy{\edge(G) \ar[r]^-{f_0} & \fC}$ of color sets is equivalent to a collection of functions of color sets
\[
\left\{
\nicexy{
\edge(C_v) \ar[r]^-{f_0^v} & \fC}
\right\}_{v \in \vertex(G)}
\]
such that
\[
f_0^u(e) = f_0^v(e)
\]
whenever $e$ is an ordinary edge in $G$ adjacent to $u$ and $v$.  Therefore, the assertion follows from the description of graphices in  $(N\sP)(G)$ in Lemma \ref{npgelement} and of the corolla ribbon in Remark \ref{rk:ribbongeometric}.
\end{proof}

\subsection{Properadic Nerve is Fully Faithful}

Before showing that the Segal condition implies a strict $\infty$-properad, we briefly digress here to observe that the properadic nerve is fully faithful.

\begin{proposition}
\label{nerve-fully-faithful}
The properadic nerve
\[
\nicexy{
\properad \ar[r]^-{N} & \gupcset
}\]
is fully faithful.
\end{proposition}

\begin{proof}
Suppose $\sP$ is a $\fC$-colored properad and $\sQ$ is a $\fD$-colored properad.  

To show that the properadic nerve $N$ is faithful, suppose $f,g : \sP \to \sQ$ are properad maps such that $Nf = Ng$ in $\gupcset$.  We must show that $f=g$.  The map $f$ is completely determined by what it does to the color set $\fC$ and the components of $\sP$.  Since $(N\sP)(\uparrow) = \fC$, the equality
\[
(Nf)(\uparrow) = (Ng)(\uparrow)
\]
implies that $f$ and $g$ are the same function on $\fC$.

To see that $f$ and $g$ are equal on the components of $\sP$, suppose $C$ is an arbitrary corolla.  Elements in $(N\sP)(C)$ are elements in $\sP$.  As $C$ runs through all the corollas, the elements in $(N\sP)(C)$ run through all the elements in all the components of $\sP$.  Since $C$ is arbitrary, the equality
\[
(Nf)(C) = (Ng)(C)
\]
implies that $f$ and $g$ are equal on all the components of $\sP$.

To show that the properadic nerve $N$ is full, suppose $\phi : N\sP \to N\sQ$ is a map of graphical sets.  We must show that $\phi = N\varphi$ for some properad map $\varphi : \sP \to \sQ$.   On color sets, $\varphi$ is defined as the map
\[
\nicexy{\fC = (N\sP)(\uparrow) \ar[r]^-{\phi}
& (N\sQ)(\uparrow) = \fD}.
\]
On elements in $\sP$, $\varphi$ is defined as the map
\[
\nicexy{(N\sP)(C) \ar[r]^-{\phi} & (N\sQ)(C)}
\]
as $C$ runs through all the corollas.  As soon as we know that $\varphi$ is a map of properads, we will know that $N\varphi = \phi$ because the properadic nerve of a properad satisfies the Segal condition (Lemma \ref{lem:nerveissegal}), hence is completely determined by what it does to corollas and $\uparrow$.

Now we show that $\varphi$ is a map of properads.  The map $\varphi$ respects colored units and the bi-equivariant structures because $\phi$ is a map of graphical sets, hence is compatible with codegeneracy maps $C_{(1;1)} \to ~\uparrow$ and input/output relabeling of corollas.  To see that $\varphi$ respects the properadic compositions, suppose $C \to D$ is the unique inner coface map of a partially grafted corollas $D$ (Example \ref{ex0:innercoface}).  Then the map
\[
\nicexy{(N\sP)(D) \ar[r] & (N\sP)(C)}
\]
is given by the properadic composition of $\sP$.  So the fact that $\phi$ is compatible with $C \to D$ implies that $\varphi$ respects the properadic compositions.
\end{proof}

\subsection{The Segal Condition Implies Strict \texorpdfstring{$\infty$}{∞}-Properad}

Next we want to show that a graphical set that satisfies the Segal condition must be a \emph{strict} $\infty$-properad.  To prove this, we first eliminate the possibility of an inner horn having multiple fillers.

\begin{lemma}
\label{lem:segalisstrict}
Suppose $\sK \in \gupcset$ satisfies the Segal condition.  Then each inner horn
\[
\nicexy@C+12pt{
\Lambda^{d^u}[G] \ar[d] \ar[r]^-{\varphi} & \sK\\
\varGamma[G] \ar@{.>}[ur]
}\]
of $\sK$ has at most one filler.
\end{lemma}

\begin{proof}
First note that $G$ has at least two vertices, and if $\nicexy{J \ar[r] & G}$ is a face of $G$, then $J$ has at least one vertex.  Moreover, $J$ has exactly one vertex if and only if it is a corolla.

The inner horn $\varphi$ is uniquely determined by the elements
\[
\left\{\varphi_J \in \sK(J) ~|~  
\nicexy{J \ar[r] & G}
\text{ face not equal to $d^u$}\right\}
\]
that agree on common faces (Remark \ref{rk:hornofx}). By the assumed Segal condition of $\sK$, each element $\varphi_J$ is uniquely determined by the elements
\[
\left\{\varphi_w^J = \xi_w^* \left(\varphi_J\right) \in \sK(C_w) ~|~ w \in \vertex(J)\right\},
\]
where $\nicexy{C_w \ar[r]^-{\xi_w} & J}$ is the composition of outer coface maps in Definition \ref{def:xiv}.  Moreover, these $\varphi_w^J$ are compatible over $\sK(\uparrow)$ for each ordinary edge in $J$ in the sense of the diagram \eqref{kcuv}.  For $v \in \vertex(G)$ that lies in two such faces $J_1$ and $J_2$, the diagram
\[
\nicearrow
\xymatrix{
C_v \ar[r]^-{\xi_v} \ar[d]_{\xi_v} & J_1 \ar[d]\\
J_2 \ar[r] & G
}\]
in $\varGamma$ is commutative. The compatibility of the elements $\varphi_J \in \sK(J)$ now implies
\[
\xi_v^* \left(\varphi_{J_1}\right) = \xi_v^* \left(\varphi_{J_2}\right) \in \sK(C_v).
\]
So we can unambiguously write $\varphi_v \in \sK(C_v)$ for such elements whenever $v \in \vertex(J_1) \bigcap \vertex(J_2)$.

Pick two distinct almost isolated vertices, say $x$ and $y$, in $G$, which exist by Corollary \ref{aiexist}.  The two corresponding outer cofaces
\[
\nicexy{J_x \ar[r]^-{d_x} & G}
\andspace
\nicexy{J_y \ar[r]^-{d_y} & G}
\]
together cover $\edge(G) \coprod \vertex(G)$, i.e.,
\[
\begin{split}
\edge(G) 
&= \edge(J_x) \bigcup \edge(J_y),\\
\vertex(G) 
&= \vertex(J_x) \bigcup \vertex(J_y).
\end{split}\]
So $\varphi_{J_x} \in \sK(J_x)$ and $\varphi_{J_y} \in \sK(J_y)$ yield a collection of $1$-dimensional elements
\begin{equation}
\label{phiphiv}
\Phi = \left\{\varphi_v \in \sK(C_v) ~|~ v \in \vertex(G)\right\}.
\end{equation}
By the assumed Segal condition of $\sK$ again, this collection $\Phi$ corresponds to a unique element in $\sK(G)$, i.e., a map $\varGamma[G] \to \sK$.

The element $\Phi \in \sK(G)$ is the only possible candidate for a filler of the inner horn $\varphi$.  Indeed, any such filler must restrict, via the respective outer coface maps, to both $\varphi_{J_x}$ and $\varphi_{J_y}$, and hence also to the $\varphi_v$ for $v \in \vertex(G)$ determined by them.  So \emph{if} $\Phi$ is a filler of $\varphi$, then it is the only one.
\end{proof}

\begin{lemma}
\label{lem2:segalisstrict}
Suppose $\sK \in \gupcset$ satisfies the Segal condition.  Then $\sK$ is a strict $\infty$-properad.
\end{lemma}

\begin{proof}
It suffices to show that, in the context of Lemma \ref{lem:segalisstrict}, the element $\Phi \in \sK(G)$ in \eqref{phiphiv} is the unique filler of the inner horn $\varphi$. We already showed in Lemma \ref{lem:segalisstrict} that $\Phi$ is the only possible candidate for a filler of $\varphi$.  Thus, it remains to show that $\Phi$ restricts to the $\varphi_J$ for faces $J$ of $G$ not equal to the given $d^u$. We consider outer and inner cofaces separately.
\begin{enumerate}
\item
For an \emph{outer} face $\nicexy{J_z \ar[r]^-{d_z} & G}$, each vertex in $J_z$ is also a vertex in $G$.  Therefore, $\Phi$ restricts to $\varphi_{J_z} \in \sK(J_z)$, i.e.,
\[
\varphi_{J_z} 
= \{\varphi_v ~|~ v \in \vertex(J_z)\} 
= d_z^{*} \Phi.
\]
\item
Suppose $\nicexy{J_v \ar[r]^{d_v} & G}$ is an \emph{inner} coface map not equal to the given $d^u$.  In particuar, there is an inner properadic factorization
\[
G = J_v(H),
\]
where $H$ is a partially grafted corollas with vertices the $v_i$ and is substituted into $v \in \vertex(J_v)$.  We need to show the equality
\[
\varphi_w 
= \xi_w^* d_v^{*} \Phi \in \sK(C_w)
\]
for every $w \in \vertex(J_v)$. First we consider vertices in $J_v$ not equal to $v$.
\begin{enumerate}
\item
A vertex $x \in \vertex(J_v) \setminus \{v\}$ is already in $G$.  So $x$ is also a vertex in some \emph{outer} face $J_z$.  The diagram
\[
\nicearrow
\xymatrix{
C_x \ar[d]_{\xi_x} \ar[r]^-{\xi_x} & J_z \ar[d]^{d_z}\\
J_v \ar[r]^-{d_v} & G
}\]
in $\varGamma$ is commutative.  The first part about outer cofaces and the commutativity of this diagram now imply:
\[
\begin{split}
\varphi_x 
&= \xi_x^* \varphi_{J_z}\\
&= \xi_x^* d_z^{*} \Phi\\
&= \xi_x^* d_v^{*} \Phi \in \sK(C_x).
\end{split}
\]
\item
It remains to prove the equality
\[
\varphi_v = \xi_v^* d_v^{*} \Phi.
\]
We first claim that there exist outer coface maps
\[
\nicearrow
\xymatrix{
H \ar`u[rrr] `[rrr]^-{\xi_{v_1,v_2}} [rrr] \ar@{<}`d[rr] `[rr]_-{\xi} [rr] \ar[r] & \cdots \ar[r] & J_w \ar[r] & G,
}\]
whose unique composition $\xi_{v_1,v_2}$ sends the edges and vertices in $H$ to the corresponding ones in $G$. The composition $\xi$ is possibly the identity map and depends on the choice of $J_w$. 

Indeed, since $G$ has two distinct inner cofaces, it has at least three vertices. There must exist an almost isolated vertex $w$ in $G$ that is not equal to either $v_i$, since otherwise the $v_i$ would not be closest neighbors in $G$ (Theorem \ref{thm:cnbdfact}).  Suppose $\nicexy{J_w \ar[r]^-{d_w} & G}$ is the corresponding outer coface obtained by deleting the almost isolated vertex $w$.  The $v_i$ are still closest neighbors in $J_w$.  If $J_w$ has at least one other vertex, then we proceed by a finite induction to obtain the desired outer coface maps.  The uniqueness of the composition $\xi_{v_1,v_2}$ follows from Lemma  \ref{lem:mapfromfree}.  This proves the claim.

We now have a commutative diagram
\[
\nicearrow
\xymatrix{
C_v \ar[d]_{\xi_v} \ar[r]^-{d_v} & H \ar[r]^-{\xi} & J_w \ar[d]^{d_w}\\
J_v \ar[rr]^-{d_v} && G
}\]
in $\varGamma$.  The top horizontal map $d_v$ is the unique inner coface of the partially grafted corollas $H$.  Graphically, this commutative diagram says that, to obtain the corolla $C_v$ from $G$, we can first smash the closest neighbors $v_i$ together, and then delete all the other vertices.  Alternatively, we may also first delete all the vertices in $G$ not equal to either $v_i$, and then smash the $v_i$ together.

Recall that the $\varphi_J$ satisfy some compatibility conditions and that $J_w \to G$ is an \emph{outer} face.  Therefore, the above commutative square now implies:
\[
\begin{split}
\varphi_v 
&= \xi_v^* \varphi_{J_v}\\
&= (\xi d_v)^* \varphi_{J_w}\\
&= d_v^{*} \xi^* d_w^* \Phi\\
&= \xi_v^* d_v^* \Phi.
\end{split}
\]
This proves that $\Phi \in \sK(G)$ restricts to $\varphi_v \in \sK(C_v)$ as well.
\end{enumerate}
\end{enumerate}
We have shown that $\Phi \in \sK(G)$ is a filler of the inner horn $\varphi$.
\end{proof}

\section{Characterization of Strict \texorpdfstring{$\infty$}{∞}-Properads}
\label{sec:fundpropsiprop}

The purpose of this section is to prove the following characterization of strict $\infty$-properads.  It says that a strict $\infty$-properad is precisely a graphical set that satisfies the Segal condition, which in turn is equivalent to being a properadic nerve up to isomorphism.

\begin{theorem}
\label{properadnerve}
Suppose $\sK \in \gupcset$.  Then the following statements are equivalent.
\begin{enumerate}
\item
There exist a properad $\sP$ and an isomorphism $\sK \cong N\sP$.
\item
$\sK$ satisfies the Segal condition.
\item
$\sK$ is a strict $\infty$-properad.\index{strict infinity properad!characterization}
\end{enumerate}
\end{theorem}

Note that the implications $(1) \Longrightarrow (2)$ and $(2) \Longrightarrow (3)$ have already been established in Lemmas \ref{lem:nerveissegal} and \ref{lem2:segalisstrict}.  Thus, it remains to prove $(3) \Longrightarrow (1)$.

\begin{remark}
The $\varDelta \cong \varGamma(\ULin)$ version of Theorem \ref{properadnerve} is of course well-known.  For example, for a simplicial set, the equivalence of being the nerve of a small category and of being a strict $\infty$-category is proved in \cite{lurie} (Proposition 1.1.2.2).  The Segal condition goes back much further \cite{segal}.
\end{remark}

\begin{remark}
The operadic version of this theorem, using the dendroidal category $\varOmega \simeq \varGamma(\uoperad)$, is also true.
The analogous equivalence of $(1)$ and $(3)$ for the dendroidal category $\varOmega$ is in \cite{mw2} (Proposition 5.3 and Theorem 6.1).  The analogous equivalence of $(1)$ and $(2)$ is in \cite{cmb} (Corollary 2.7).  In fact, the operadic version of our corolla ribbon $\sK(G)_1$ is what Cisinski and Moerdijk call $Hom(Sc[T], X)$, where $X$ is a dendroidal set, $T$ is a unital tree, and $Sc[T] \subseteq \varOmega[T]$ is their Segal core.  The last object is defined as the union of the images of the compositions of outer dendroidal coface maps $\varOmega[C_v] \to \varOmega[T]$ as $v$ runs through the vertices in $T$.
\end{remark}

\begin{remark}
The equivalence of $(1)$ and $(3)$ in Theorem \ref{properadnerve} can be interpreted as saying that our definition of an $\infty$-properad  correctly expresses the notion of an up-to-homotopy properad.  Indeed, the properadic nerve functor is fully faithful, so we may identify a properad with its properadic nerve.  In a properad, the properadic composition is uniquely defined, and the axioms (associativity, unity, and bi-equivariance) hold in the strict sense.  Just like an $\infty$-category, an $\infty$-properad should capture the notion of an up-to-homotopy properad.  So the properadic composition exists, but it is \emph{not} required to be unique.  Moreover, the axioms should only be required to hold up to homotopy, and any two such homotopies of the same axiom should be homotopic, and so forth.  Part of the theorem says that properadic nerves are exactly the \emph{strict} $\infty$-properads.  The uniqueness of the inner horn extensions are expressing the ideas that properadic compositions are uniquely defined and that axioms hold on the nose.  With the uniqueness condition removed, an $\infty$-properad is then expressing the notion of an up-to-homotopy properad.
\end{remark}

\subsection{Properad Associated to a Strict \texorpdfstring{$\infty$}{∞}-Properad}

To prove $(3) \Longrightarrow (1)$ in Theorem \ref{properadnerve}, we first need some preliminary definitions and constructions.

\begin{definition}
\label{def:properadpk}
Fix a strict $\infty$-properad $\sK$.  We define a $\sK(\uparrow)$-colored properad $\sP_{\sK}$ as follows.
\begin{enumerate}
\item
The $\Sigma_{\sS(\sK(\uparrow))}$-bimodule structure of  $\sP_{\sK}$ is that of $\sK \in \gupcset$ as in section \ref{rk:gsetproperad}.  In particular, its elements are the $1$-dimensional elements in $\sK$, with 
\begin{itemize}
\item
input/output profiles induced by the outer coface maps $\uparrow ~\to C$ into corollas, 
\item
colored units induced by the codegeneracy map  $C_{(1;1)} \to ~ \uparrow$, and 
\item
bi-equivariant structure induced by input/output relabeling.
\end{itemize}
\item
For the properadic composition, suppose
\begin{itemize}
\item
$\theta \in \sK(C_1)$ has profiles $\yxh$,
\item
$\phi \in \sK(C_2)$ has profiles $\wvh$, and
\item 
$\uw \supseteq \uw' = \ux' \subseteq \ux$ are some equal $k$-segments in $\uw$ and $\ux$ for some $k>0$.
\end{itemize}
There exists a unique $\sK(\uparrow)$-colored partially grafted corollas $D$:
\begin{itemize}
\item
whose top (resp., bottom) vertex has the input/output profiles of $C_1$ (resp., $C_2$), and
\item
whose ordinary edges correspond to the equal $k$-segments $\uw' = \ux'$.
\end{itemize}
This partially grafted corollas $D$ has:
\begin{itemize}
\item
a single inner face $\nicexy{H \ar[r]^-{d_u} & D}$, corresponding to smashing the two vertices together, and 
\item
two outer faces $\nicexy{C_i \ar[r] & D}$.
\end{itemize}
By Example \ref{ex:pgcorhorn}, the elements $\theta$ and $\phi$ determine a unique inner horn $\varphi$
\[
\nicearrow
\xymatrix@C+12pt{
\Lambda^u[D] \ar[d] \ar[r]^-{\varphi} & \sK\\
\varGamma[D] \ar@{.>}[ur]_{\exists ! \Phi}
}\]
of $\sK$.  Since $\sK$ is assumed to be a strict $\infty$-properad, there is a unique filler $\Phi \in \sK(D)$ of $\varphi$.  We now define the \textbf{properadic composition}
\begin{equation}
\label{pkcomposition}
\theta \boxtimes^{\ux'}_{\uw'} \phi \defn d_u^*(\Phi) \in \sK(H).
\end{equation}
That $d_u^*(\Phi)$ has the correct profiles, that is,
\[
\theta \boxtimes^{\ux'}_{\uw'} \phi \in \sP_{\sK}\binom{\uw \circ_{\uw'} \uy}{\ux \circ_{\ux'} \uv},
\]
is a consequence of the construction of the partially grafted corollas $D$.  We will usually abbreviate the properadic composition $\boxtimes^{\ux'}_{\uw'}$ to just $\boxtimes$.
\end{enumerate}
This finishes the definition of $\sP_{\sK}$.
\end{definition}

\begin{lemma}
\label{properadpk}
For each strict $\infty$-properad $\sK$, the above definitions give a $\sK(\uparrow)$-colored properad $\sP_{\sK}$.\index{strict infinity properad!associated properad}
\end{lemma}

\begin{proof}
We need to check the two bi-equivariance axioms, the unity axiom, and the various associativity axioms of a properad in biased form.  They are all proved in the same manner using the unity and associativity of graph substitution, which are proved with full details and generality in \cite{jy2}.  So we will prove only one properadic associativity axiom in detail to illustrate the method.

Consider a $\sK(\uparrow)$-colored connected wheel-free graph $B$ with three vertices,
\begin{center}
\begin{tikzpicture}
\matrix[row sep=.2cm, column sep=1cm]{
\node [plain] (w) {$w$}; &\\
& \node [plain] (v) {$v$};\\
\node [plain] (u) {$u$}; &\\
};
\draw [implies] (u) to (v);
\draw [implies] (u) to (w);
\draw [implies] (v) to (w);
\end{tikzpicture}
\end{center}
in which each double arrow $\Longrightarrow$ means there is at least one ordinary edge in that direction between the indicated vertices.  Each vertex may have input and output legs of $B$, which are not depicted in the picture.

Suppose $\omega$, $\theta$, and $\phi$ are elements in $\sP_{\sK}$ with the profiles of these three vertices, so the decorated graph
\begin{center}
\begin{tikzpicture}
\matrix[row sep=.2cm, column sep=1cm]{
\node [plain] (w) {$\omega$}; &\\
& \node [plain] (v) {$\theta$};\\
\node [plain] (u) {$\phi$}; &\\
};
\draw [implies] (u) to (v);
\draw [implies] (u) to (w);
\draw [implies] (v) to (w);
\end{tikzpicture}
\end{center}
makes sense.  In particular, there are properadic compositions $\omega \boxtimes \theta$, $\theta \boxtimes \phi$, $\omega \boxtimes (\theta \boxtimes \phi)$, and $(\omega \boxtimes \theta) \boxtimes \phi$.  We want to prove the properadic associativity axiom
\[
\omega \boxtimes (\theta \boxtimes \phi) = (\omega \boxtimes \theta) \boxtimes \phi
\]
in $\sP_{\sK}$. To do this, first note that $B$ has four faces:
\begin{enumerate}
\item
There is an \emph{outer} face
\[
\nicexy{C_{u,v} \ar[r]^-{d_w} & B}
\]
corresponding to deleting the almost isolated vertex $w$.  
\begin{center}
\begin{tikzpicture}
\matrix[row sep=.2cm, column sep=1cm]{
& \node [emptyvt] (o2) {}; &&& \node [plain] (w) {$w$}; &\\
\node [emptyvt] (o1) {}; & 
\node [plain] (v1) {$v$}; & 
\node [empty] (s) {}; &
\node [empty] (t) {}; && 
\node [plain] (v2) {$v$};\\
\node [plain] (u1) {$u$}; &&&&
\node [plain] (u2) {$u$}; &\\
};
\draw [arrow] (s) to node{$d_w$} (t);
\draw [implies] (u1) to (v1);
\draw [implies] (u1) to (o1);
\draw [implies] (v1) to (o2);
\draw [implies] (u2) to (v2);
\draw [implies] (u2) to (w);
\draw [implies] (v2) to (w);
\end{tikzpicture}
\end{center}
By Example \ref{ex:pgcorhorn}, the elements $\theta$ and $\phi$ determine a unique inner horn $\varphi_0$
\[
\nicearrow
\xymatrix@C+12pt{
\Lambda^t[C_{u,v}] \ar[r]^-{\varphi_0} \ar[d] & \sK\\
\varGamma[C_{u,v}] \ar@{.>}[ur]_{\exists ! \Phi_0} &
}\]
in $\sK$.  Since $\sK$ is a strict $\infty$-properad, the inner horn $\varphi_0$ has a  unique filler $\nicexy{\varGamma[C_{u,v}] \ar[r]^-{\Phi_0} & \sK}$.  This defines the properadic composition
\[
\theta \boxtimes \phi 
= d_t^* \Phi_0 \in \sK(C_{uv}).
\]
Here
\[
\nicexy{C_{uv} \ar[r]^-{d_t} & C_{u,v}}
\]
is the unique inner face of the partially grafted corollas $C_{u,v}$.
\item
There is an \emph{outer} face
\[
\nicexy{C_{v,w} \ar[r]^-{d_u} & B}
\]
corresponding to deleting the almost isolated vertex $u$.  
\begin{center}
\begin{tikzpicture}
\matrix[row sep=.2cm, column sep=1cm]{
\node [plain] (w1) {$w$}; &&&& 
\node [plain] (w2) {$w$}; &\\
\node [emptyvt] (i1) {}; 
& \node [plain] (v1) {$v$}; & 
\node [empty] (s) {}; &
\node [empty] (t) {}; && 
\node [plain] (v2) {$v$};\\
& \node [emptyvt]  (i2) {}; &&&
\node [plain] (u) {$u$}; &\\
};
\draw [arrow] (s) to node{$d_u$} (t);
\draw [implies] (v1) to (w1);
\draw [implies] (i1) to (w1);
\draw [implies] (i2) to (v1);
\draw [implies] (u) to (v2);
\draw [implies] (u) to (w2);
\draw [implies] (v2) to (w2);
\end{tikzpicture}
\end{center}
By Example \ref{ex:pgcorhorn} again, the elements $\omega$ and $\theta$ determine a unique inner horn 
\[
\nicearrow
\xymatrix@C+12pt{
\Lambda^r[C_{v,w}] \ar[r]^-{\varphi_1} \ar[d] & \sK\\
\varGamma[C_{v,w}]. \ar@{.>}[ur]_{\exists ! \Phi_1} &
}\]
Its unique filler $\nicexy{\varGamma[C_{v,w}] \ar[r]^-{\Phi_1} & \sK}$ then gives the properadic composition
\[
\omega \boxtimes \theta 
= d_r^* \Phi_1 \in \sK\left(C_{vw}\right),
\]
where
\[
\nicexy{C_{vw} \ar[r]^-{d_r} & C_{v,w}}
\]
is the unique inner face of the partially grafted corollas $C_{v,w}$.
\item
There is an \emph{inner} face
\[
\nicexy{C_{(uv),w} \ar[r]^-{d_{uv}} & B}
\]
corresponding to smashing together the closest neighbors $u$ and $v$ in $B$.
\begin{center}
\begin{tikzpicture}
\matrix[row sep=.2cm, column sep=1cm]{
\node [plain] (w1) {$w$}; &&& 
\node [plain] (w2) {$w$}; &\\
& \node [empty] (s) {}; &
\node [empty] (t) {}; && 
\node [plain] (v) {$v$};\\
\node [plain]  (uv) {$uv$}; &&&
\node [plain] (u) {$u$}; &\\
};
\draw [arrow] (s) to node{$d_{uv}$} (t);
\draw [implies] (uv) to (w1);
\draw [implies] (u) to (v);
\draw [implies] (u) to (w2);
\draw [implies] (v) to (w2);
\end{tikzpicture}
\end{center}
Once again, the elements $\omega$ and $\theta \boxtimes \phi$ determine a unique inner horn
\[
\nicearrow
\xymatrix@C+12pt{
\Lambda^q[C_{(uv),w}] \ar[r]^-{\varphi_2} \ar[d] & \sK\\
\varGamma[C_{(uv),w}]. \ar@{.>}[ur]_{\exists ! \Phi_2} &
}\]
Its unique filler $\nicexy{\varGamma[C_{(uv),w}] \ar[r]^-{\Phi_2} & \sK}$ then gives the properadic composition
\[
\omega \boxtimes (\theta \boxtimes \phi) 
= d_q^* \Phi_2 \in \sK\left(C_{(uv)w}\right),
\]
where
\[
\nicexy{C_{(uv)w} \ar[r]^-{d_q} & C_{(uv),w}}
\]
is the unique inner face of the partially grafted corollas $C_{(uv),w}$.
\item
Lastly, there is an \emph{inner} face
\[
\nicexy{C_{u,(vw)} \ar[r]^-{d_{vw}} & B}
\]
corresponding to smashing together the closest neighbors $v$ and $w$ in $B$.
\begin{center}
\begin{tikzpicture}
\matrix[row sep=.2cm, column sep=1cm]{
\node [plain] (vw) {$vw$}; &&& 
\node [plain] (w) {$w$}; &\\
& \node [empty] (s) {}; &
\node [empty] (t) {}; && 
\node [plain] (v) {$v$};\\
\node [plain]  (u1) {$u$}; &&&
\node [plain] (u2) {$u$}; &\\
};
\draw [arrow] (s) to node{$d_{vw}$} (t);
\draw [implies] (u1) to (vw);
\draw [implies] (u2) to (v);
\draw [implies] (u2) to (w);
\draw [implies] (v) to (w);
\end{tikzpicture}
\end{center}
\end{enumerate}

Consider this last inner face of $B$. The above elements $\Phi_i$ for $0 \leq i \leq 2$ in $\sK$ determine a unique inner horn
\[
\nicearrow
\xymatrix@C+12pt{
\Lambda^{d_{vw}}[B] \ar[r]^-{\psi} \ar[d] & \sK\\
\varGamma[B], \ar@{.>}[ur]_{\exists ! \Psi}
}\]
which has a unique filler $\Psi \in \sK(B)$.  We claim that
\begin{equation}
\label{omegathetaphi}
(\omega \boxtimes \theta) \boxtimes \phi 
= d_p^* d_{vw}^* \Psi \in \sK\left(C_{u(vw)}\right),
\end{equation}
where
\[
\nicexy{C_{u(vw)} \ar[r]^-{d_p} & C_{u,(vw)}}
\]
is the unique inner face of the partially grafted corollas $C_{u,(vw)}$.  To keep the flow of the proof, we will show \eqref{omegathetaphi} immediately after this proof.

On the other hand, the diagram
\[
\nicearrow
\xymatrix{
C_{(uv)w} = C_{u(vw)} \ar[r]^-{d_p} \ar[d]_{d_q} & C_{u,(vw)} \ar[d]^{d_{vw}}\\
C_{(uv),w} \ar[r]^-{d_{uv}} & B
}\]
in $\varGamma$ is commutative by the associativity of graph substitution.  Graphically, this commutative diagram says that, to compose $B$ down to a corolla, one can first smash together the closest neighbors $u$ and $v$, and then smash $w$ into it.  Alternatively, one can also first smash together the closest neighbors $v$ and $w$, and then smash $u$ into it. The commutativity of this diagram now implies:
\[
\begin{split}
\omega \boxtimes (\theta \boxtimes \phi) 
&= d_q^* \Phi_2 \\
&= d_q^* d_{uv}^* \Psi\\
&= d_p^* d_{vw}^* \Psi\\
&= (\omega \boxtimes \theta) \boxtimes \phi.
\end{split}
\]
This proves the desired associativity axiom.
\end{proof}

\begin{sublemma}
The equality \eqref{omegathetaphi} holds, i.e.,
\[
(\omega \boxtimes \theta) \boxtimes \phi 
= d_p^* d_{vw}^* \Psi \in \sK\left(C_{u(vw)}\right).
\]
\end{sublemma}

\begin{proof}
As before, the partially grafted corollas $C_{u,(vw)}$ has:
\begin{itemize}
\item
a single inner face $\nicexy{C_{u(vw)} \ar[r]^-{d_p} & C_{u,(vw)}}$ and
\item
two outer faces $\nicexy{C_{vw} \ar[r]^-{d_u} & C_{u,(vw)}}$ and $\nicexy{C_u \ar[r]^-{d} & C_{u,(vw)}}$. 
\end{itemize}
By Example \ref{ex:pgcorhorn}, the elements $\omega \boxtimes \theta \in \sK(C_{vw})$ and $\phi \in \sK(C_u)$ determine a unique inner horn
\[
\nicearrow
\xymatrix@C+12pt{
\Lambda^p[C_{u,(vw)}] \ar[d] \ar[r] & \sK\\
\varGamma[C_{u,(vw)}]. \ar@{.>}[ur]_{\exists !} &
}\]
Since $\sK$ is assumed to be a strict $\infty$-properad, there is a unique filler.  From the definition of the properadic composition in $\sP_{\sK}$, the desired equality is equivalent to $d_{vw}^* \Psi \in \sK(C_{u,(vw)})$ being the unique filler.  In other words, we need to show that the two \emph{outer} faces of $d_{vw}^* \Psi$ are exactly $\omega \boxtimes \theta$ and $\phi$.  

To prove this, first observe that there is a commutative diagram
\[
\nicearrow
\xymatrix{
C_{vw} \ar[d]_{d_u} \ar[rr]^-{d_r} && C_{v,w} \ar[d]^{d_u}\\
C_{u,(vw)} \ar[rr]^-{d_{vw}} && B\\
C_u \ar[u]^{d}  \ar[rr]^-{d_v} && C_{u,v} \ar[u]_{d_w}
}\]
in $\varGamma$. Graphically, the top commutative rectangle says that to obtain the corolla $C_{vw}$ from $B$, one can first delete the almost isolated vertex $u$, and then smash the closest neighbors $v$ and $w$ together.  Alternatively, one can first smash $v$ and $w$ together, and then delete the almost isolated vertex $u$.  The bottom commutative rectangle says that to obtain the corolla $C_u$ from $B$, one can first delete the almost isolated vertex $w$, and then delete the almost isolated vertex $v$.  Alternatively, one can first smash together the closest neighbors $v$ and $w$, and then delete the combined vertex, which is almost isolated.

The two outer faces of the partially grafted corollas $C_{u,(vw)}$ are the two left vertical maps.  Therefore, the top commutative rectangle yields:
\[
\begin{split}
\omega \boxtimes \theta 
&= d_r^* \Phi_1\\
&= d_r^* d_u^* \Psi\\
&= d_u^* d_{vw}^* \Psi.
\end{split}
\]
This shows that $\omega \boxtimes \theta$ is the expected outer face of $d_{vw}^*(\Psi)$.  Similarly, the bottom commutative rectangle yields:
\[
\begin{split}
\phi 
&= d_v^* \Phi_0\\
&= d_v^* d_w^* \Psi\\
&= d^* d_{vw}^* \Psi.
\end{split}
\]
Therefore, $\phi$ is the other expected outer face of $d_{vw}^* \Psi$.
\end{proof}

\subsection{Strict \texorpdfstring{$\infty$}{∞}-Properads are Properadic Nerves}

We now proceed to prove $(3) \Longrightarrow (1)$ in Theorem \ref{properadnerve} by showing that $\sP_{\sK}$ is the desired properad.  First we show that there is an object-wise isomorphism.

\begin{lemma}
\label{etaiso}
Suppose $\sK$ is a strict $\infty$-properad, and $G \in \gupc$.  Then there is a canonical bijection
\[
\nicearrow
\xymatrix{\sK(G) \ar[r]^-{\eta}_-{\cong} & (N\sP_{\sK})(G),
}\]
where $\sP_{\sK}$ is the properad in Lemma \ref{properadpk}.
\end{lemma}

\begin{proof}
The map $\eta$ is constructed according to the number of vertices in $G$ as follows.
\begin{enumerate}
\item
The only graph in $\gupc$ with $0$ vertex is the exceptional edge $\uparrow$.  By definition there is a bijection
\[
\nicearrow
\xymatrix{
\sK(\uparrow) \ar[d]_-{\eta^0}^-{\cong}\\
(N\sP_{\sK})(\uparrow) 
= \properad(\uparrow,\sP_{\sK})
}\]
because by Lemma \ref{npgelement} a map $\uparrow ~ \to \sP_{\sK}$ is simply an element in the color set of $\sP_{\sK}$, which is $\sK(\uparrow)$.
\item
The only graphs in $\gupc$ with $1$ vertex are the corollas and their input/output relabelings.  For a permuted corolla $\sigma C\tau$, the bijection
\[
\nicearrow
\xymatrix{
\sK(\sigma C \tau) \ar[d]_-{\eta^1}^-{\cong}\\
(N\sP_{\sK})(\sigma C\tau) 
= \properad(\sigma C \tau,\sP_{\sK})
}\]
comes from the definition of the elements in $\pk$ as the $1$-dimensional elements in $\sK$.
\item
Suppose $G \in \gupc$ has $n \geq 2$ vertices.  The map $\eta^n$ is defined as the composition
\begin{equation}
\label{etankg}
\nicearrow
\xymatrix{
\sK(G) \ar[r]^-{\eta^n} \ar[d]_-{\chi_G} 
& (N\sP_{\sK})(G)\\
\sK(G)_1 \ar@{=}[r] 
& \left[\prod \sK(C_v)\right]_{\sK(\uparrow)} .\ar[u]_-{\prod \eta^1}^-{\cong} &\\
}
\end{equation}
Here $\chi_G$ is the Segal map \eqref{segalmap}, and $\prod \eta^1$ is the bijection defined by the previous case and the fact that $N\pk$ satisfies the Segal condition (Lemma \ref{lem:nerveissegal}).
\end{enumerate}

Next we show that $\eta^n$ is a bijection by induction.  We already observed that $\eta^0$ and $\eta^1$ are bijections.

If $G \in \gupc$ has $n \geq 2$ vertices, then it must have at least one inner face, say $d_u$.  The diagram
\begin{equation}
\label{knpeta}
\nicearrow
\xymatrix@C+10pt{
\gupcset(\varGamma[G], \sK) \ar[d]_-{\cong} \ar[r]^-{\eta^n} & \gupcset(\varGamma[G], N\sP_{\sK}) \ar[d]^-{\cong}\\
\gupcset(\Lambda^u[G], \sK) \ar[r]^-{\eta^{n-1}} & \gupcset(\Lambda^u[G],N\sP_{\sK})
}
\end{equation}
is commutative by the construction of the maps $\eta$.  We want to show that $\eta^n$ is a bijection.  The left vertical map is a bijection by the strict $\infty$-properad assumption on $\sK$.  The right vertical map is a bijection because $N\sP_{\sK}$ is also a strict $\infty$-properad (Lemmas \ref{lem:nerveissegal},  \ref{lem:segalisstrict}, and \ref{lem2:segalisstrict}).  The bottom horizontal map is determined by the maps
\[
\nicearrow
\xymatrix@C+10pt{
\gupcset(\varGamma[J], \sK) \ar[r]^-{\eta^{n-1}} & \gupcset(\varGamma[J],N\sP_{\sK})
}\]
for faces $J$ of $G$ not equal to $d_u$ (Lemma \ref{horndescription}).  Each such face $J$ has $n-1$ vertices.  So by induction hypothesis, the bottom horizontal map in \eqref{knpeta} is a bijection.  Therefore, the top horizontal map $\eta^n$ is also a bijection.
\end{proof}

Next we show that the map $\eta$ has the same universal property as the unit of the adjunction object-wise.

\begin{lemma}
\label{etaknerveptwo}
Suppose $\sK$ is a strict $\infty$-properad, and $\nicearrow\xymatrix{\sK \ar[r]^-{\zeta} & N\sQ}$ is a map in $\gupcset$ for some $\fD$-colored properad $\sQ$.  Then there exists a unique map $\nicearrow\xymatrix{\sP_{\sK} \ar[r]^-{\zeta'} & \sQ}$ of properads such that the diagram
\begin{equation}
\label{etanzeta}
\nicearrow
\xymatrix{
\sK(G) \ar[d]_-{\eta}^-{\cong} \ar[r]^-{\zeta} & (N\sQ)(G)\\
(N\sP_{\sK})(G) \ar[ur]_-{N\zeta'} & 
}
\end{equation}
is commutative for each $G \in \gupc$.
\end{lemma}

\begin{proof}
First we construct the properad map $\zeta'$.
\begin{enumerate}
\item
Recall that the color set of $\sP_{\sK}$ is $\sK(\uparrow)$. 
The map $\zeta'$ on color sets is defined as the map
\[
\nicearrow
\xymatrix{
\sK(\uparrow) \ar[d]^-{\zeta}\\
(N\sQ)(\uparrow) 
= \properad(\uparrow, \sQ) 
= \fD.
}\]
\item
The component maps
\[
\nicearrow
\xymatrix{
\sP_{\sK}\yx \ar[r]^-{\zeta'} & \sQ\binom{\zeta \uy}{\zeta \ux}
}\]
are defined using Lemma \ref{npgelement} and the fact that elements in $\pk$ are $1$-dimensional elements in $\sK$.
\item
That $\zeta'$ preserves the bi-equivariant structures follows from the fact that the bi-equivariant structure in $\sP_{\sK}$ is induced by input/output relabelings on permuted corollas.  These relabelings induce maps in the graphical category $\varGamma$, which are preserved by $\zeta$.
\item
To see that $\zeta'$ preserves the colored units, recall that for $c \in \sK(\uparrow)$, the $c$-colored unit in $\sP_{\sK}$ is defined as
\[
\bone_c = s^* c, 
\]
where $s^*$ is induced by the unique codegeneracy map $\nicexy{C_{(1;1)} \ar[r]^-{s} & \uparrow}$.  That $\zeta'(\bone_c)$ is the corresponding colored unit of $\sQ$ is then a consequence of the fact that $\zeta$ is a map of graphical sets, and hence respects codegeneracy maps.
\item
To see that $\zeta'$ preserves the properadic compositions, simply note that:
\begin{enumerate}
\item
the properadic composition in $\sP_{\sK}$ is induced by coface maps in $\sK$ \eqref{pkcomposition}, which are preserved by $\zeta$, and
\item
the properadic nerve $N\sQ$ is also a strict $\infty$-properad  (Lemmas \ref{lem:nerveissegal},  \ref{lem:segalisstrict}, and \ref{lem2:segalisstrict}).
\end{enumerate}
\end{enumerate}

Next, to see that the composition $N\zeta' \circ \eta$ in \eqref{etanzeta} is equal to $\zeta$, first observe that they agree on the exceptional edge $\uparrow$ and permuted corollas. Suppose $G \in \gupc$ has $n \geq 2$ vertices.  The map $\zeta_G$ is the composition
\[
\nicearrow
\xymatrix{
\sK(G) \ar[r]^-{\zeta_G} \ar[d]_{\chi_G} & (N\sQ)(G)\\
\sK(G)_1 \ar[r]^-{\prod \zeta} & (N\sQ)(G)_1 \ar[u]_{\chi_G^{-1}}^{\cong}
}\]
because the properadic nerve $N\sQ$ satisfies the Segal condition (Lemma \ref{lem:nerveissegal}).  By the definition of the corolla ribbon, the bottom horizontal map is determined by $\zeta$ on corollas, where it agrees with the composition $N\zeta' \circ \eta$.  The agreement between $\zeta$ and $N\zeta' \circ \eta$ on $G$ now follows from the definition of $\eta^n$ \eqref{etankg}.

Finally, we observe the uniqueness of the properad map $\zeta'$ for which $N\zeta' \circ \eta = \zeta$.  Indeed, this equality already determines what the map $\zeta'$ does on color sets and elements in $\sP_{\sK}$.
\end{proof}

\begin{lemma}
\label{strictisnerve}
Suppose $\sK$ is a strict $\infty$-properad.  Then the object-wise bijections in Lemma \ref{etaiso} assemble to give an isomorphism
\[
\nicearrow
\xymatrix{
\sK \ar[r]^-{\eta}_-{\cong} & N\sP_{\sK}
}\]
in $\gupcset$.
\end{lemma}

\begin{proof}
Pick a map $G \to H \in \varGamma$.  We need to show that the square
\begin{equation}
\label{etaknervep}
\nicearrow
\xymatrix{
\sK(H) \ar[r] \ar[d]_-{\eta}^-{\cong} & \sK(G) \ar[d]_-{\cong}^-{\eta}\\
(N\sP_{\sK})(H) \ar[r] & (N\sP_{\sK})(G) 
}
\end{equation}
is commutative.  

Suppose $\nicearrow
\xymatrix{\sK \ar[r]^-{\zeta} & N\sQ}$ is an arbitrary map in $\gupcset$ for some properad $\sQ$.  Consider the following diagram.
\[
\nicearrow
\xymatrix{
\sK(H) \ar[d]_-{\eta}^-{\cong} \ar[rr]^-{\zeta} \ar[ddrr] &&  (N\sQ)(H) \ar[ddrr] &&\\
(N\sP_{\sK})(H) \ar[ddrr] \ar@{.>}[urr]_{N\zeta'} && &&\\
&& \sK(G) \ar[d]_-{\eta}^-{\cong} \ar[rr]^-{\zeta} && (N\sQ)(G)\\
&& (N\sP_{\sK})(G) \ar[urr]_-{N\zeta'} &&
}\]
Here $\nicexy{\sP_{\sK} \ar[r]^-{\zeta'} & \sQ}$ is the unique properad map determined by $\zeta$ in Lemma \ref{etaknerveptwo}.  Both $\zeta$ and $N\zeta'$ are maps in $\gupcset$.  Therefore, the two compositions in the diagram
\[
\nicearrow
\xymatrix{
\sK(H) \ar[r] \ar[d]_-{\eta}^-{\cong} & \sK(G) \ar[d]_-{\cong}^-{\eta} \ar[r]^-{\zeta} & (N\sQ)(G)\\
(N\sP_{\sK})(H) \ar[r] & (N\sP_{\sK})(G) \ar[ur]_-{N\zeta'} & 
}\]
from $\sK(H)$ to $(N\sQ)(G)$ are equal.  Since this is true for a \emph{fixed} map $G \to H$ and for \emph{all} maps from $\sK$ to an arbitrary properadic nerve $N\sQ$, the square \eqref{etaknervep} must be commutative.
\end{proof}

Theorem \ref{properadnerve} is now a consequence of Lemmas \ref{lem:nerveissegal}, \ref{lem:segalisstrict}, \ref{lem2:segalisstrict}, \ref{properadpk}, and \ref{strictisnerve}.

\subsection{Fundamental Properad}

Recall the adjoint pair
\[
\nicexy@C+10pt{
\gupcset \ar@<.5ex>[r]^-{L} 
& \properad. \ar@<.5ex>[l]^-{N}}
\]
The following definition is the properadic analog of the fundamental category of a simplicial set.

\begin{definition}
\label{def:fundproperad}
For $\sK \in \gupcset$, the image $L\sK$ is called the \textbf{fundamental properad} \index{fundamental properad} of $\sK$.
\end{definition}

\begin{corollary}
\label{cor:fundproperad}
Suppose $\sK$ is a strict $\infty$-properad.  Then the following statements hold.
\begin{enumerate}
\item
The $\sK(\uparrow)$-colored properad $\sP_{\sK}$ is canonically isomorphic to the fundamental properad of $\sK$.
\item
The map
\[
\nicearrow
\xymatrix{
\sK \ar[r]^-{\eta}_-{\cong} & N\sP_{\sK}}
\]
is the unit of the adjunction $(L,N)$.
\end{enumerate}
\end{corollary}

\begin{proof}
Lemma \ref{strictisnerve} says that $\eta$ is a map of graphical sets, while Lemma \ref{etaknerveptwo} says that it has the required universal property of the unit of the adjunction.
\end{proof}

Corollary \ref{cor:fundproperad} gives an explicit description of the fundamental properad of a strict $\infty$-properad.  The fundamental properad of a general $\infty$-properad is considered in chapter \ref{ch:fundprop}.


\chapter{Fundamental Properads of Infinity Properads}
\label{ch:fundprop}

\abstract*{We give an explicit description of the fundamental properad $L\sK$ of an $\infty$-properad $\sK$.  The fundamental properad of an $\infty$-properad consists of homotopy classes of $1$-dimensional elements.  It takes a bit of work to prove that there is a well-defined homotopy relation among $1$-dimensional elements and that a properad structure can be defined on homotopy classes.}

The purpose of this chapter is to give an explicit description of the fundamental properad (Definition \ref{def:fundproperad}) of an $\infty$-properad with non-empty inputs, with non-empty outputs, or that is reduced.  This description is directly inspired by the construction of the fundamental category of an $\infty$-category in \cite{bv} (4.11 and 4.12).  For a strict $\infty$-properad, this description reduces to the one in Corollary \ref{cor:fundproperad}.

Briefly, for an $\infty$-properad $\sK$ with some mild restriction, its fundamental properad $\qk$ is a $\sK(\uparrow)$-colored properad and has as elements \emph{homotopy classes} of $1$-dimensional elements in $\sK$.  The assignment
\[
\nicearrow
\xymatrix@C+15pt{
\left(\text{$\infty$-properad $\sK$}\right) \ar@{|->}[r] & \left(\text{fundamental properad $\qk$}\right) 
}\]
is similar in spirit to the assignment
\[
\nicearrow
\xymatrix@C+15pt{
\left(\text{path space of $X$}\right) \ar@{|->}[r] & \left(\text{fundamental groupoid $\Pi_1(X)$}\right)
}\]
for a topological space $X$.  Indeed, the left-hand side has objects that contain all higher homotopy information.  There is no unique way to compose two composable paths $g$ and $f$ in $X$.  A composition $g \circ f$ exists, and any two such compositions are homotopic.  Any two such homotopies are themselves homotopic, and so forth.  The fundamental groupoid $\Pi_1(X)$ makes all the axioms--composition, associativity, and unity--hold in the strict sense.

Likewise, in an $\infty$-properad $\sK$ with some mild restriction, there is no unique way to properadically compose two composable operations $g$ and $f$.  A properadic composition $g \boxtimes f$ exists by the existence of inner horn fillers.  Any two such properadic compositions are homotopic, as we will show below.  The fundamental properad $\qk$ makes all the properad axioms--composition, associativity, unity, and bi-equivariance--hold in the strict sense.

The homotopy relation among $1$-dimensional elements is defined in section \ref{sec:homotopygset}.  Roughly speaking, two $1$-dimensional elements $f$ and $g$ are homotopic along the $i$th input (or output) if there is an $i$th input (or output) extension of $f$ by a degenerate element, whose inner face is $g$.  We need to show that each such relation is an equivalence relation and that these equivalence relations are all equal to each other.

The properad $\qk$ associated to an $\infty$-properad $\sK$ with some mild restriction is constructed in section \ref{sec:propinfprop}.  The elements in $\qk$ are homotopy classes of $1$-dimensional elements in $\sK$.  We need to show that there is a well-defined properad structure on $\qk$ and that $\qk$ has the desired universal property.

\section{Homotopy in a Graphical Set}
\label{sec:homotopygset}

The purpose of this section is to define the homotopy relation among the $1$-dimensional elements in an $\infty$-properad.  The first task is to show that homotopy defined using each leg is actually an equivalence relation.  The second task is to show that these relations are all equal to each other.  In other words, homotopy of $1$-dimensional elements is independent of the choice of an input or an output leg in its definition.  In the next section, we will use homotopy classes of $1$-dimensional elements to define the fundamental properad of an $\infty$-properad with some mild restriction.

Everything in this section has obvious analogs for $\infty$-properads with non-empty inputs or non-empty outputs (Definition \ref{def:infpropnonempty}).

\subsection{Motivation from Fundamental Categories}

Our description and proof of the fundamental properad of an $\infty$-properad, as well as the fundamental \emph{operad} of an $\infty$-operad in \cite{mw2}, are directly inspired by the construction of the fundamental category of an $\infty$-category in \cite{bv} (4.11 and 4.12).  So let us briefly review that construction as motivation for our construction later.

Fix a simplicial set $X$.  Two $1$-simplices $f$ and $g$ in $X_1$ are said to be \textbf{homotopic},\index{homotopy!of $1$-simplices} written $f \sim g$, if there exists a $2$-simplex $H \in X_2$ such that
\[
d_2 H = f, \quad d_1 H = g, \andspace d_0 H = s_0d_0 f.
\]
In this case, write $H \colon f \sim g$, and call $H$ a \textbf{homotopy} from $f$ to $g$.  This situation can be visualized as the following diagram.
\[
\nicearrow
\xymatrix@=10pt@R-4pt{
&& 1 \ar[ddrr]^-{\text{degenerate}} &&\\
&& H &&\\
0 \ar[uurr]^-{f} \ar[rrrr]_-{g} &&&& 2
}\]
When $f$ and $g$ are homotopic, the simplicial identities ensure that $f$ and $g$ have the same faces.

When $X$ is an $\infty$-category, Boardman and Vogt \cite{bv} (4.11) showed that homotopy is an equivalence relation.  Moreover, there is a genuine category \cite{bv} (4.12), called the fundamental category, whose objects are the $0$-simplices of $X$ and whose morphisms from $x$ to $y$ are the homotopy classes of $1$-simplices $f$ with $d_1f = x$ and $d_0f = y$.

Furthermore, when $X$ is a \emph{strict} $\infty$-category (i.e., isomorphic to the nerve of a small category), the homotopy relation is the identity relation.  Indeed, given a $1$-simplex $f$, the elements $f$ and $s_0d_0f$ determine a unique inner $1$-horn in $X$.
\[
\nicearrow
\xymatrix{
& 1 \ar[dr]^-{\text{degenerate}} &\\
0 \ar[ur]^-{f} \ar@{.>}[rr]_-{\exists !} && 2
}\]
So there is a unique filler $H \in X_2$ whose $2$-face is $f$ and whose $0$-face is $s_0d_0f$.  By uniqueness its $1$-face is necessarily $f$ again.  In particular, the fundamental category of the nerve of a small category is canonically isomorphic to the given category.

Before we consider the homotopy relation for a graphical set $\sK$, let us first interpret the homotopy relation for a simplicial set in graph theoretic terms.  When we regard $\varDelta$ as the graphical subcategory $\varGamma(\ULin)$, the object $[0] \in \varDelta$ corresponds to the exceptional edge.  So $0$-simplices in a simplicial set $X$ are exactly the elements in $X(\uparrow)$.  In $\gupc$ the only graph with $0$ vertex is also the exceptional edge $\uparrow$.  So  for a graphical set $\sK$, the graphical analogue of the object set $X(\uparrow)$ is the set $\sK(\uparrow)$.

The element $[1] \in \varDelta \cong \varGamma(\ULin)$ corresponds to the linear graph $L_1$ with $1$ vertex.  The only graphs in $\gupc$ with $1$ vertex are the corollas and their input/output relabelings.  So the graphical analog of $L_1$ is a permuted corolla $\sigma C \tau$.  The graphical analog of a $1$-simplex is a $1$-dimensional element $f \in \sK(\sigma C\tau)$ for some permuted corolla, which may be depicted as follows.
\begin{center}
\begin{tikzpicture}
\matrix[row sep=1cm,column sep=1cm] {
\node [plain,label=above:$...$,label=below:$...$] (p) {$f$}; \\
};
\draw [outputleg] (p) to +(-.6cm,.5cm);
\draw [outputleg] (p) to +(.6cm,.5cm);
\draw [inputleg] (p) to +(-.6cm,-.5cm);
\draw [inputleg] (p) to +(.6cm,-.5cm);
\end{tikzpicture}
\end{center}

For the linear graph $L_2$ with $2$ vertices, $d^2$ and $d^0$ are its \emph{outer} faces, while $d^1$ is its only \emph{inner} face. The $\gupc$-analogue of $L_2$ is a connected wheel-free graph with $2$ vertices, i.e., a partially grafted corollas.  Given $f$ and $g$ in a graphical set $\sK$, a homotopy $H$ should then be an element in $\sK$ corresponding to a partially grafted corollas whose outer faces are $f$ and a degenerate element, while its unique inner face is $g$.  In this case, a degenerate element has the shape of a corolla $C_{(1;1)}$.  So the partially grafted corolla is actually a basic dioperadic graph with a corolla $C_{(1;1)}$ as one of its vertices, as in the following picture.
\begin{center}
\begin{tikzpicture}
\matrix[row sep=.2cm, column sep=1.5cm]{
&&& \node [plain] (one) {$\bone$};\\
\node [plain, label=above:$...$, label=below:$...$] (g) {$g$}; &
\node [empty] (s) {}; &
\node [empty] (t) {}; & & \\
&&& \node [plain, label=below:$...$] (f) {$f$};\\
};
\draw [arrow] (s) to node{inner} (t);
\draw [inputleg] (g) to +(-.6cm,-.5cm);
\draw [inputleg] (g) to +(.6cm,-.5cm);
\draw [outputleg] (g) to +(-.6cm,.5cm);
\draw [outputleg] (g) to +(.6cm,.5cm);
\draw [inputleg] (f) to +(-.6cm,-.5cm);
\draw [inputleg] (f) to +(.6cm,-.5cm);
\draw [outputleg] (f) to +(-.6cm,.5cm);
\draw [outputleg] (f) to +(.6cm,.5cm);
\draw [arrow] (f) to node[swap]{$i$} (one);
\draw [outputleg] (one) to +(0,.7cm);
\end{tikzpicture}
\end{center}
In this picture, the right-hand side depicts the shape of $H$ as well as its two outer faces, while the left-hand side depicts its unique inner face.

One complication here is that $f$ and $g$ in general have $m \geq 0$ inputs and $n \geq 0$ outputs.  For example, in the picture above, $i$ can be anywhere within the interval $[1,n]$.  There should also be an analogous picture where $f$ decorates the top vertex, while the degenerate element $\bone$ decorates the bottom vertex and is connected to the $j$th input of the top vertex, as in the following picture.  
\begin{center}
\begin{tikzpicture}
\matrix[row sep=.2cm, column sep=1.5cm]{
&&& \node [plain, label=above:$...$] (f) {$f$};\\
\node [plain, label=above:$...$, label=below:$...$] (g) {$g$}; &
\node [empty] (s) {}; &
\node [empty] (t) {}; && \\
&&&  \node [plain] (one) {$\bone$}; \\
};
\draw [arrow] (s) to node{inner} (t);
\draw [inputleg] (g) to +(-.6cm,-.5cm);
\draw [inputleg] (g) to +(.6cm,-.5cm);
\draw [outputleg] (g) to +(-.6cm,.5cm);
\draw [outputleg] (g) to +(.6cm,.5cm);
\draw [inputleg] (f) to +(-.6cm,-.5cm);
\draw [inputleg] (f) to +(.6cm,-.5cm);
\draw [outputleg] (f) to +(-.6cm,.5cm);
\draw [outputleg] (f) to +(.6cm,.5cm);
\draw [arrow] (one) to node[swap]{$j$} (f);
\draw [inputleg] (one) to +(0,-.7cm);
\end{tikzpicture}
\end{center}
In this case, $j$ can be anywhere within the interval $[1,m]$.  Each one of these $m+n$ scenarios can reasonably be called a homotopy from $f$ to $g$.  Therefore, part of the work is to show that the choice of an input or an output is irrelevant because they should all lead to the same homotopy relation.

\subsection{Homotopy of \texorpdfstring{$1$}{1}-Dimensional Elements}

We will use the terminology introduced in section  \ref{rk:gsetproperad} for a graphical set.  We begin by defining homotopy of $1$-dimensional elements along an input or an output.

\begin{definition}
\label{def:graphhomotopy}
Suppose $\sK \in \gupcset$, $f$ and $g$ are $1$-dimensional elements in $\sK$ with $m$ inputs and $n$ outputs, and $C = C_{(m;n)}$ is the corolla with $m$ inputs and $n$ outputs.
\begin{enumerate}
\item
For each $1 \leq i \leq m$, say \textbf{$f$ is homotopic to $g$ along the $i$th input},\index{homotopy!along an input} written $f \sim_i g$,\label{note:simi} if there exists an element
\[
H \in \sK\left(C \boxtimes^{i}_1 C_{(1;1)}\right)
\]
such that
\[
d_{\text{top}}^*(H) = f,\quad d_{\text{in}}^*(H) = g,\andspace d_{u}^*(H) = \bone_{x_i}.
\]
These elements are defined as follows.
\begin{itemize}
\item
$v$ and $u$ are, respectively, the top and bottom vertices in the dioperadic graph
\[
D_i \defn C \boxtimes^{i}_1 C_{(1;1)}
\]
with $\nicexy{C \ar[r]^-{d_v} & D_i}$ and $\nicexy{C_{(1;1)} \ar[r]^-{d_u} & D_i}$ the corresponding outer face maps.
\item
The map $d_{\text{top}}$ is the composition
\begin{equation}
\label{dtop}
\nicearrow
\xymatrix{
\sigma C \tau \ar[rr]^-{(\tau^{-1};\sigma^{-1})}_-{\cong} \ar`u[rrr] `[rrr]^-{d_{\text{top}}} [rrr] && C \ar[r]^-{d_v} & D_i},
\end{equation}
where $(\tau^{-1};\sigma^{-1})$ is the isomorphism that sends:
\begin{itemize}
\item
the $k$th input/output leg of $\sigma C \tau$ to that of $C$, and
\item
the unique vertex $v \in \sigma C \tau$ to the $\widehat{C}$-decorated graph $\sigma C \tau$.
\end{itemize}
\item
$x_{i}\in \sK(\uparrow)$ is the $i$th input profile of $f$.
\item
$\bone_{x_i} = s\left(x_{i}\right) \in \sK(C_{(1;1)})$ is the $x_i$-colored unit.
\item
The map $d_{\text{in}}$ is the composition
\[
\nicearrow\xymatrix{
\sigma' C \tau'  \ar`u[rrrr] `[rrrr]^-{d_{\text{in}}} [rrrr] \ar[rr]^-{(\tau'^{-1};\sigma'^{-1})}_-{\cong} && C \ar[rr]^-{\text{inner}}_-{\text{face}} && D_i}.
\]
\end{itemize}
In this case, we also write $H \colon f \sim_i g$, and call $H$ a \textbf{homotopy from $f$ to $g$ along the $i$th input}.
\item
For each $1 \leq j \leq n$, say  \textbf{$f$ is homotopic to $g$ along the $j$th output},\index{homotopy!along an output} written $f \sim^j g$,\label{note:simj} if there exists an element
\[
H \in \sK\left(C_{(1;1)} \boxtimes^1_{j} C  \right)
\]
such that
\[
d_u^*(H) = \bone_{y_j} ,\quad d_{\text{in}}^*(H) = g, \andspace d_{\text{bot}}^*(H) = f.
\]
These elements are defined as follows.
\begin{itemize}
\item
$u$ and $v$ are, respectively, the top and bottom vertices in the dioperadic graph
\[
D^j \defn C_{(1;1)} \boxtimes^1_{j} C
\]
with $\nicearrow\xymatrix{
C_{(1;1)} \ar[r]^-{d_u} & D^j}$ and $\nicearrow\xymatrix{C \ar[r]^-{d_v} & D^j}$ the corresponding outer face maps.
\item
The map $d_{\text{bot}}$ is the composition
\[
\nicearrow
\xymatrix{
\sigma C \tau \ar[rr]^-{(\tau^{-1};\sigma^{-1})}_-{\cong} \ar@{<}`u[rrr] `[rrr]^-{d_{\text{bot}}} [rrr] && C \ar[r]^-{d_v} & D^j}.
\]
\item
$y_j \in \sK(\uparrow)$ is the $j$th output profile of $f$.
\item
$\bone_{y_j} = s(y_j) \in \sK(C_{(1;1)})$ is the $y_j$-colored unit.
\item
The map $d_{\text{in}}$ is the composition
\[
\nicearrow\xymatrix{
\sigma' C \tau'  \ar@{<} `u[rrrr] `[rrrr]^-{d_{\text{in}}} [rrrr] \ar[rr]^-{(\tau'^{-1};\sigma'^{-1})}_-{\cong} && C \ar[rr]^-{\text{inner}}_-{\text{face}} && D^j}.
\]
\end{itemize}
In this case, we write also $H \colon f \sim^j g$, and call $H$ a \textbf{homotopy from $f$ to $g$ along the $j$th output}.
\end{enumerate}
\end{definition}

\begin{remark}
\label{rk:homotopyextension}
We may represent the dioperadic graphs $D_i$ and $D^j$ as follows.
\begin{center}
\begin{tikzpicture}
\matrix[row sep=1.2cm, column sep=4cm]{
\node [plain, label=above:$...$] (v1) {$v$}; &
\node [plain] (u2) {$u$};\\
\node [plain] (u1) {$u$}; &
\node [plain, label=below:$...$] (v2) {$v$};\\
};
\foreach \x in {1,2}
{
\draw [inputleg] (v\x) to +(-.6cm,-.5cm);
\draw [inputleg] (v\x) to +(.6cm,-.5cm);
\draw [outputleg] (v\x) to +(-.6cm,.5cm);
\draw [outputleg] (v\x) to +(.6cm,.5cm);
}
\draw [inputleg] (u1) to +(0,-.7cm);
\draw [arrow] (u1) to node{\footnotesize{$i$}} (v1);
\draw [outputleg] (u2) to +(0,.7cm);
\draw [arrow] (v2) to node{\footnotesize{$j$}} (u2);
\end{tikzpicture}
\end{center}
Therefore, one way to interpret the relation $f \sim_i g$ is that there exists an $i$th input extension $H$ of $f$ by degeneracy, whose inner face is $g$.  Likewise, one can say that the relation $f \sim^j g$ means that there exists a $j$th output extension of $f$ by degeneracy, whose inner face is $g$.  In particular, the $\sim_i$ (resp., $\sim^j$) are defined if and only if $m > 0$ (resp., $n > 0$).
\end{remark}

\begin{remark}
Let us explain more explicitly the isomorphism
\[
\nicearrow
\xymatrix{
\sigma C \tau \ar[rr]^-{(\tau^{-1};\sigma^{-1})}_-{\cong} && C}
\]
in \eqref{dtop}.  Write $i_l$ (resp., $o_j$) for the $l$th input (resp., $j$th output) leg in $C$.  This flag has the same label at the vertex $v \in C$, since $C$ is an unpermuted corolla. The edge sets $\edge(C)$ and $\edge(\sigma C\tau)$ are, by definition, equal. At the vertex $v \in \sigma C \tau$, the flag $i_l$ (resp., $o_j$) is still labeled $l$ (resp., $j$), while as an input (resp., output) flag of $\sigma C\tau$ it is labeled $\tau^{-1}(l)$ (resp., $\sigma(j)$). So sending the $k$th output leg of $\sigma C\tau$ to the $k$th output leg of $C$ for all $k$ means the assignment
\[
\nicearrow
\xymatrix@C+12pt{
\edge(\sigma C \tau) \ni o_j \ar@{|->}[r]^-{\varphi_0} & o_{\sigma(j)} \in \edge(C).
}\]
Likewise, sending the $k$th input leg of $\sigma C \tau$ to the $k$th input leg of $C$ for all $k$ means the assignment
\[
\nicearrow
\xymatrix@C+12pt{
\edge(\sigma C \tau) \ni i_l \ar@{|->}[r]^-{\varphi_0} & i_{\tau^{-1}(l)} \in \edge(C).
}\]
The profiles of $v \in \sigma C \tau$ are the same as those of $v \in C$, i.e.,
\[
\binom{\out v}{\inp v} = \binom{o_1, \ldots , o_n}{i_1, \ldots, i_m}.
\]
The $\widehat{C}$-decorated graph $\sigma C \tau$ has profiles
\[
\binom{o_{\sigma(1)}, \ldots , o_{\sigma(n)}}{i_{\tau^{-1}(1)}, \ldots , i_{\tau^{-1}(m)}} 
= \binom{\varphi_0 \out v}{\varphi_0 \inp v},
\]
which is as expected for a map in the graphical category $\varGamma$.
\end{remark}

\begin{remark}
In the above definition, $\sK$ can be any graphical set.  However, we will restrict to an $\infty$-properad below when we show that the relations are equivalence relations.
\end{remark}

First we observe that every relation in the definition above preserves profiles.  This observation will be important below when we define a properad associated to an $\infty$-properad using homotopy classes of $1$-dimensional elements.

\begin{lemma}
\label{simprofile}
Suppose $\sK \in \gupcset$, and $f$ and $g$ are $1$-dimensional elements in $\sK$ with $m$ inputs and $n$ outputs.  Suppose either
\begin{itemize}
\item
$f \sim_i g$ for some $1 \leq i \leq m$, or
\item
$f \sim^j g$ for some $1 \leq j \leq n$.
\end{itemize}
Then the profiles of $f$ and $g$ are equal.
\end{lemma}

\begin{proof}
We will use the notations in Definition \ref{def:graphhomotopy}. First suppose we are given a homotopy $H \colon f \sim_i g$ along the $i$th input for some $i$. By construction the map $(\tau^{-1};\sigma^{-1})$ preserves inputs and outputs.  Therefore, we may safely suppress the permutations associated to $f$ and $g$ in the following discussion.
\begin{enumerate}
\item
We first consider a leg of $C = C_{(m;n)}$ that is \emph{not} the $i$th input, and suppose $\nicearrow\xymatrix{\uparrow \ar[r]^-{\eta} & C}$ is the corresponding outer face map.  Then there is a commutative square
\begin{equation}
\label{etaequalize}
\nicearrow
\xymatrix{
C \ar[r]^-{d_{\inp}} & D_i\\
\uparrow \ar[r]^-{\eta} \ar[u]^-{\eta} & C \ar[u]_-{d_{\out}},
}
\end{equation}
where $d_{\inp}$ and $d_{\out}$ are the inner and outer coface maps associated to $g$ and $f$.  The commutativity of this square says that in $D_i$, each of its $n+m-1$ legs attached to the top vertex can be obtained in two different ways.  In the following picture of $D_i$, they are depicted as the dotted arrows.
\begin{center}
\begin{tikzpicture}
\matrix[row sep=1.2cm, column sep=1.5cm]{
\node [plain, label=above:$...$] (v) {$v$};\\
\node [plain] (u) {$u$}; \\
};
\draw [dottedinput] (v) to +(-.7cm,-.5cm);
\draw [dottedinput] (v) to +(.7cm,-.5cm);
\draw [dottedarrow] (v) to +(-.7cm,.5cm);
\draw [dottedarrow] (v) to +(.7cm,.5cm);
\draw [inputleg] (u) to +(0,-.7cm);
\draw [arrow] (u) to node{$i$} (v);
\end{tikzpicture}
\end{center}
Starting with $H \in \sK(D_i)$ and pulling back to $\sK(\uparrow)$ using the commutative square \eqref{etaequalize}, it follows that $f$ and $g$ have the same output profiles and $k$th input profile, provided $k \not= i$.
\item
It remains to observe that $f$ and $g$ have the same $i$th input profile.  There is a commutative square
\[
\nicearrow
\xymatrix@C+8pt{
C \ar[r]^-{d_{\inp}} & D_i\\
\uparrow \ar[u]^-{\inp_i} \ar[r]^-{\inp_1} & C_{(1;1)} \ar[u]_-{d_{\out}} \ar@<5pt>[l]^-{s},
}\]
where $d_{\out}$ is the outer coface map associated to the degenerate element. The map $s$ is the codegeneracy map, which is a section of $\inp_1$.  The commutativity of this square says that the $i$th input leg of $D_i$ is both:
\begin{itemize}
\item
the unique input leg of its lower vertex, and
\item
the $i$th input leg after the two closest neighbors in $D_i$ have been smashed together.
\end{itemize}
In the following picture of $D_i$, this input leg is depicted as the dotted arrow.
\begin{center}
\begin{tikzpicture}
\matrix[row sep=1.2cm, column sep=1.5cm]{
\node [plain, label=above:$...$] (v) {$v$};\\
\node [plain] (u) {$u$}; \\
};
\draw [inputleg] (v) to +(-.7cm,-.5cm);
\draw [inputleg] (v) to +(.7cm,-.5cm);
\draw [outputleg] (v) to +(-.7cm,.5cm);
\draw [outputleg] (v) to +(.7cm,.5cm);
\draw [dottedinput] (u) to +(0,-.9cm);
\draw [arrow] (u) to node{$i$} (v);
\end{tikzpicture}
\end{center}
Therefore, we have
\[
\begin{split}
\sK(\uparrow) \ni x_i 
&= \inp_1^* s^*(x_i)\\
&= \inp_1^*(\bone_{x_i})\\
&= \inp_1^* d_{\out}^*(H)\\
&= \inp_i^* d_{\inp}^*(H)\\
&= \inp_i^*(g).
\end{split}
\]
This means that the $i$th input profiles of $f$ and $g$ are the same.
\end{enumerate}

There is a dual proof for the case when $f \sim^j g$.
\end{proof}

\begin{remark}
Suppose $\nicearrow\xymatrix{\uparrow \ar[r]^-{\eta} & C}$ is the outer face map that identifies $\uparrow$ with the $i$th input leg of $C$.  Then the square \eqref{etaequalize} is \emph{not} commutative.  Indeed, on the one hand, the composition
\[
\nicearrow
\xymatrix{
\uparrow \ar[r]^-{\eta} & C \ar[r]^-{d_{\inp}} & D_i
}\]
identifies $\uparrow$ with the $i$th input leg of $D_i$.  On the other hand, the composition
\[
\nicearrow
\xymatrix{
\uparrow \ar[r]^-{\eta} & C \ar[r]^-{d_{\out}} & D_i
}\]
identifies $\uparrow$ with the unique ordinary edge in $D_i$.
\end{remark}

\subsection{Homotopy Relations are Equivalence Relations}

Our next objective is to show that all the $\sim_i$ and $\sim^j$ are equivalence relations.

\begin{lemma}
\label{simijrelation}
Suppose $\sK$ is an $\infty$-properad, and $m,n \geq 0$. Then the relations
\begin{itemize}
\item
$\sim_i$ for $1 \leq i \leq m$, and 
\item
$\sim^j$ for $1 \leq j \leq n$
\end{itemize}
are all equivalence relations.
\end{lemma}

\begin{proof}
The proof for $\sim^j$ is dual to that for $\sim_i$, so we only provide the proof for the latter.  We assume $m \geq 1$, since otherwise there is nothing to prove.  We check the three required conditions for $\sim_i$ to be an equivalence relation.  To improve readability, we separate the three parts into several sublemmas.
\end{proof}

\begin{sublemma}
\label{simireflexive}
The relation $\sim_i$ is reflexive.
\end{sublemma}

\begin{proof}
Pick a $1$-dimensional element $f \in \sK(\sigma C\tau)$, where $C = C_{(m;n)}$.  We want to show that $f \sim_i f$. There is a codegeneracy map
\[
\nicexy{D_i \ar[r]^-{s} & C}
\]
given by substituting the exceptional edge into the bottom vertex of $D_i$.  The composition
\[
\nicearrow
\xymatrix@C+6pt{
D_i \ar`u[rr] `[rr]^-{\zeta_i} [rr] \ar[r]^-{s} & C \ar[r]^-{(\tau;\sigma)}_-{\cong} & \sigma C \tau
}\]
yields a degeneracy element
\begin{equation}
\label{homotopyff}
\zeta_i^*(f) \in \sK(D_i)
\end{equation}
that is a homotopy from $f$ to $f$ along the $i$th input.  The reason is that the diagram
\[
\nicearrow
\xymatrix@C+6pt{
C \ar`u[rr] `[rr]^-{Id} [rr] \ar@<.6ex>[r]^-{d_v} \ar@<-.6ex>[r]_-{d_\text{in}} & D_i \ar[r]^-{s} & C\\
& C_{(1;1)} \ar[u]_-{d_u} \ar[r]^-{s} & \uparrow \ar[u]_-{\text{outer}} 
}\]
is commutative by the unity property of graph substitution.  Here $d_{\text{in}}$ is the unique inner coface map to $D_i$, while the right outer coface map identifies $\uparrow$ with the $i$th input leg of $C$.
\end{proof}

\begin{remark}
The above sublemma does not use the assumption that the graphical set $\sK$ is an $\infty$-properad.
\end{remark}

\begin{remark}
\label{rk:ignorepermute}
Given any $1$-dimensional element $f \in \sK(\sigma C\tau)$, there always exists a $1$-dimensional element $\sigma^{-1}f\tau^{-1} \in \sK(C)$ such that
\[
f \sim_i \sigma^{-1}f\tau^{-1}.
\]
In fact, using the proof of Sublemma \ref{simireflexive}, we have maps
\[
\nicearrow
\xymatrix@R-8pt{
C \ar[r]^-{d_{\text{in}}} & D_i \ar[r]^-{s} & C \ar[r]^-{(\tau;\sigma)} & \sigma C \tau,\\
d_{\text{in}}^*\zeta_i^*f & \zeta_i^*(f) \ar@{|->}[l] && f. \ar@{|->}[ll]
}\]
We can then use the element $d_{\text{in}}^*\zeta_i^*f \in \sK(C)$ as our $\sigma^{-1}f\tau^{-1}$.  This implies that in most discussion and computation that follows, we may safely suppress the permutations $\sigma$ and $\tau$ associated to a $1$-dimensional element $f \in \sK(\sigma C \tau)$.  This will greatly simplify the presentation.
\end{remark}

\begin{sublemma}
\label{simisymmetric}
The relation $\sim_i$ is symmetric.  
\end{sublemma}

\begin{proof}
Suppose there is a homotopy
\[
(H \colon f \sim_i g) \in \sK(D_i)
\]
from $f$ to $g$ along the $i$th input. We want to show $g \sim_i f$.  Consider the $3$-vertex graph
\begin{equation}
\label{cmnltwo}
\begin{split}
E_i &= D_i \boxtimes^i_1 C_{(1;1)}\\
&= \left[C_{(m;n)} \boxtimes^i_1 C_{(1;1)}\right] \boxtimes^i_1 C_{(1;1)}\\
&= C_{(m;n)} \boxtimes^i_1 \underbrace{\left[C_{(1;1)} \boxtimes^1_1 C_{(1;1)}\right]}_{L_2},
\end{split}
\end{equation}
which may be depicted as follows.
\begin{center}
\begin{tikzpicture}
\matrix[row sep=1.2cm, column sep=1.5cm]{
\node [plain, label=above:$...$] (v) {$v$};\\
\node [plain] (u) {$u$}; \\
\node [plain] (t) {$t$};\\
};
\draw [inputleg] (v) to +(-.7cm,-.5cm);
\draw [inputleg] (v) to +(.7cm,-.5cm);
\draw [outputleg] (v) to +(-.7cm,.5cm);
\draw [outputleg] (v) to +(.7cm,.5cm);
\draw [inputleg] (t) to +(0,-.7cm);
\draw [arrow] (t) to node{$e'$} (u);
\draw [arrow] (u) to node{$e$} node[swap,near end]{\footnotesize{$i$}} (v);
\end{tikzpicture}
\end{center}
In the following discussion, we suppress the permutations associated to $f$ and $g$ to simpify the presentation. We want to define an inner horn in $\sK$ corresponding to the closest neighbors $u$ and $v$, so first note that $E_i$ has four faces.
\begin{enumerate}
\item
There is an \emph{outer} face
\[
\nicexy{L_2 \ar[r]^-{d_v} & E_i}
\]
corresponding to deleting the almost isolated vertex $v$.  
\begin{center}
\begin{tikzpicture}
\matrix[row sep=1.2cm, column sep=1.5cm]{
&&& \node [plain, label=above:$...$] (v) {$v$};\\
\node [plain] (u1) {$u$}; &
\node [empty] (s) {}; & \node [empty] (t) {}; & 
\node [plain] (u2) {$u$}; \\
\node [plain] (t1) {$t$}; &&&
\node [plain] (t2) {$t$};\\
};
\draw [arrow] (s) to node{$d_v$} (t);
\draw [inputleg] (t1) to +(0,-.7cm);
\draw [arrow] (t1) to node{$e'$} (u1);
\draw [outputleg] (u1) to +(0,.7cm);
\draw [inputleg] (v) to +(-.7cm,-.5cm);
\draw [inputleg] (v) to +(.7cm,-.5cm);
\draw [outputleg] (v) to +(-.7cm,.5cm);
\draw [outputleg] (v) to +(.7cm,.5cm);
\draw [inputleg] (t2) to +(0,-.7cm);
\draw [arrow] (t2) to node{$e'$} (u2);
\draw [arrow] (u2) to node{$e$} node[swap,near end]{\footnotesize{$i$}} (v);
\end{tikzpicture}
\end{center}
Moreover, there are two codegeneracy maps
\[
\nicearrow
\xymatrix{
L_2 \ar[r]^-{s} & C_{(1;1)} \ar[r]^-{s} & \uparrow,
}\]
so there is a double degeneracy
\[
s^* \bone_{x_i} 
= s^* s^* x_i \in \sK(L_2).
\]
\item
There is another \emph{outer} face
\[
\nicexy{D_i \ar[r]^-{d_t} & E_i}
\]
corresponding to deleting the almost isolated vertex $t$.
\begin{center}
\begin{tikzpicture}
\matrix[row sep=1.2cm, column sep=1.5cm]{
\node [plain, label=above:$...$] (v1) {$v$};
&&& 
\node [plain, label=above:$...$] (v2) {$v$};\\
\node [plain] (u1) {$u$}; &
\node [empty] (s) {}; & \node [empty] (t) {}; & 
\node [plain] (u2) {$u$}; \\
&&& \node [plain] (t2) {$t$};\\
};
\draw [arrow] (s) to node{$d_t$} (t);
\draw [inputleg] (u1) to +(0,-.7cm);
\draw [arrow] (u1) to node{$e$} node[swap, near end]{\footnotesize{$i$}} (v1);
\foreach \x in {1,2}
{
\draw [inputleg] (v\x) to +(-.7cm,-.5cm);
\draw [inputleg] (v\x) to +(.7cm,-.5cm);
\draw [outputleg] (v\x) to +(-.7cm,.5cm);
\draw [outputleg] (v\x) to +(.7cm,.5cm);
}
\draw [inputleg] (t2) to +(0,-.7cm);
\draw [arrow] (t2) to node{$e'$} (u2);
\draw [arrow] (u2) to node{$e$} node[swap,near end]{\footnotesize{$i$}} (v2);
\end{tikzpicture}
\end{center}
For this outer face, we use the given homotopy
\[
(H \colon f \sim_i g) \in \sK(D_i)
\]
from $f$ to $g$ along the $i$th input.
\item
There is an \emph{inner} face
\[
\nicexy{D_i \ar[r]^-{d_{e'}} & E_i}
\]
corresponding to smashing together the closest neighbors $t$ and $u$.
\begin{center}
\begin{tikzpicture}
\matrix[row sep=1.2cm, column sep=1.5cm]{
\node [plain, label=above:$...$] (v1) {$v$};
&&& 
\node [plain, label=above:$...$] (v2) {$v$};\\
\node [plain] (u1) {$tu$}; &
\node [empty] (s) {}; & \node [empty] (t) {}; & 
\node [plain] (u2) {$u$}; \\
&&& \node [plain] (t2) {$t$};\\
};
\draw [arrow] (s) to node{$d_{e'}$} (t);
\draw [inputleg] (u1) to +(0,-.7cm);
\draw [arrow] (u1) to node{$e$} node[swap, near end]{\footnotesize{$i$}} (v1);
\foreach \x in {1,2}
{
\draw [inputleg] (v\x) to +(-.7cm,-.5cm);
\draw [inputleg] (v\x) to +(.7cm,-.5cm);
\draw [outputleg] (v\x) to +(-.7cm,.5cm);
\draw [outputleg] (v\x) to +(.7cm,.5cm);
}
\draw [inputleg] (t2) to +(0,-.7cm);
\draw [arrow] (t2) to node{$e'$} (u2);
\draw [arrow] (u2) to node{$e$} node[swap,near end]{\footnotesize{$i$}} (v2);
\end{tikzpicture}
\end{center}
For this inner face, we use the self-homotopy
\[
(\zeta^*_i(f) \colon f \sim_i f) \in \sK(D_i)
\]
from $f$ to $f$ along the $i$th input \eqref{homotopyff}.
\item
Finally, there is an \emph{inner} face
\[
\nicexy{D_i \ar[r]^-{d_{e}} & E_i}
\]
corresponding to smashing together the closest neighbors $u$ and $v$.
\begin{center}
\begin{tikzpicture}
\matrix[row sep=1.2cm, column sep=1.5cm]{
&&& 
\node [plain, label=above:$...$] (v2) {$v$};\\
\node [plain, label=above:$...$] (v1) {$uv$};
& \node [empty] (s) {}; & \node [empty] (t) {}; & 
\node [plain] (u2) {$u$}; \\
\node [plain] (t1) {$t$}; 
&&& \node [plain] (t2) {$t$};\\
};
\draw [arrow] (s) to node{$d_{e}$} (t);
\draw [inputleg] (t1) to +(0,-.7cm);
\draw [arrow] (t1) to node{$e'$} node[swap, near end]{\footnotesize{$i$}} (v1);
\foreach \x in {1,2}
{
\draw [inputleg] (v\x) to +(-.7cm,-.5cm);
\draw [inputleg] (v\x) to +(.7cm,-.5cm);
\draw [outputleg] (v\x) to +(-.7cm,.5cm);
\draw [outputleg] (v\x) to +(.7cm,.5cm);
}
\draw [inputleg] (t2) to +(0,-.7cm);
\draw [arrow] (t2) to node{$e'$} (u2);
\draw [arrow] (u2) to node{$e$} node[swap,near end]{\footnotesize{$i$}} (v2);
\end{tikzpicture}
\end{center}
\end{enumerate}

By Lemma \ref{horndescription}, the double degeneracy $s^*(\bone_{x_i})$ and the homotopies $H$ and $\zeta^*_i(f)$ define an inner horn
\[
\nicearrow
\xymatrix@C+12pt{
\Lambda^e[E_i] \ar[d] \ar[r] & \sK\\
\varGamma[E_i] \ar@{.>}[ur]_-{\exists \Phi} &
}\]
in $\sK$.  Since $\sK$ is assumed to be an $\infty$-properad, there exists a filler $\Phi \in \sK(E_i)$.  We claim that the inner face filled by it,
\[
H' \defn d_e^* \Phi \in \sK(D_i),
\]
is a homotopy from $g$ to $f$ along the $i$th input.  In other words, we need to show that
\begin{itemize}
\item
its outer faces are $g$ and a degenerate element, and
\item
its inner face is $f$.
\end{itemize}
\begin{enumerate}
\item
The reason $H'$ has $g$ as an outer face is that the square
\[
\nicearrow
\xymatrix{
D_i \ar[r]^-{d_e} & E_i\\
C \ar[u]^-{d_t} \ar[r]^-{d_e} & D_i \ar[u]_-{d_t}
}\]
is commutative, which in turn holds because deleting $t$ and smashing together the closest neighbors $u$ and $v$ are commuting operations. Therefore, we have
\[
\begin{split}
d_t^* H' 
&= d_t^* d_e^* \Phi\\
&= d_e^* d_t^* \Phi\\
&= d_e^* H\\
&= g.
\end{split}
\]
\item
The other outer face of $H'$ is degenerate because the square
\[
\nicearrow
\xymatrix@+10pt@C+6pt{
D_i \ar[r]^-{d_e} & E_i\\
C_{(1;1)} \ar[u]^-{d_{uv}} \ar[r]^-{d_u} & L_2 \ar[u]_-{d_v} \ar@<5pt>[l]^-{s}
}\]
is commutative, while the codegeneracy $s$ is a section of $d_u$.  Here $d_{uv}$ is the outer coface map corresponding to deleting the combined vertex of $u$ and $v$, which is almost isolated in $D_i$.  The commutativity of the above square is simply saying that the corolla $C_{(1;1)}$ with vertex $t$ can be obtained:
\begin{itemize}
\item
by first deleting the almost isolated vertex $v$ and then deleting the almost isolated vertex $u$, or
\item
by first smashing the closest neighbors $u$ and $v$ together and then deleting that combined vertex.
\end{itemize}
Therefore, we have
\[
\begin{split}
d_{uv}^* H' 
&= d_{uv}^* d_e^* \Phi\\
&= d_u^* d_v^*\Phi\\
&= d_u^* s^* \bone_{x_i}\\
&= \bone_{x_i}.
\end{split}
\]
\item
For the inner face of $H'$, we have a commutative diagram
\[
\nicearrow
\xymatrix{
D_i \ar[r]^-{d_e} & E_i\\
C \ar[u]^-{d_{e'}} \ar[r]^-{d_e} & D_i \ar[u]_-{d_{e'}}.
}\]
The commutativity of this square says that, to compose $E_i$ down to a corolla, the order in which $e$ and $e'$ are shrunk is irrelevant.  Therefore, we have
\[
\begin{split}
d_{e'}^* H'
&= d_{e'}^* d_e^* \Phi\\
&= d_e^* d_{e'}^* \Phi\\
&= d_e^* \zeta_i^* f\\
&= f.
\end{split}
\]
\end{enumerate}
We have shown that $H' \colon g \sim_i f$, so the relation $\sim_i$ is symmetric.
\end{proof}

\begin{sublemma}
The relation $\sim_i$ is transitive.
\end{sublemma}

\begin{proof}
Suppose given homotopies $H_1 \colon f \sim_i g$ and $H_2 \colon g \sim_i h$ along the $i$th input.  We want to show $f \sim_i h$. We will use the notations in Sublemma \ref{simisymmetric}.

Consider once again the graph $E_i$ in \eqref{cmnltwo}.  Define an inner horn
\begin{equation}
\label{innerhorneprime}
\nicearrow
\xymatrix@C+12pt{
\Lambda^{e'}[E_i] \ar[r] \ar[d] & \sK\\
\varGamma[E_i] \ar@{.>}[ur]_{\exists \Psi}
}
\end{equation}
using  Lemma \ref{horndescription} and the following elements.
\begin{enumerate}
\item
For the outer face $\nicexy{L_2 \ar[r]^-{d_v} & E_i}$, we again use the double degeneracy $s^*(\bone_{x_i}) \in \sK(L_2)$.
\item
For the outer face  $\nicexy{D_i \ar[r]^-{d_t} & E_i}$, we use the given homotopy
\[
\left(H_1 \colon f \sim_i g\right) \in \sK(D_i)
\]
from $f$ to $g$ along the $i$th input.
\item
For the inner face $\nicexy{D_i \ar[r]^-{d_{e}} & E_i}$, we use the other homotopy
\[
\left(H_2 \colon g \sim_i h\right) \in \sK(D_i)
\]
from $g$ to $h$ along the $i$th input.
\end{enumerate}

The double degeneracy $s^*(\bone_{x_i})$ and the two homotopies $H_1$ and $H_2$ define the desired inner horn in \eqref{innerhorneprime}.  Suppose $\Psi \in \sK(E_i)$ is a filler of that inner horn. Recall the inner face $\nicexy{D_i \ar[r]^-{d_{e'}} & E_i}$.  We claim that the inner face filled by $\Psi$,
\[
H'' \defn d_{e'}^* \Psi \in \sK(D_i),
\]
is a homotopy from $f$ to $h$ along the $i$th input.   In other words, we need to show that
\begin{itemize}
\item
its outer faces are $f$ and a degenerate element, and 
\item
its inner face is $h$.
\end{itemize}
The computation is almost the same as in  Sublemma \ref{simisymmetric}.
\begin{enumerate}
\item
The left commutative square and the computation in $\sK$ on the right-hand side,
\[
\nicearrow
\xymatrix{
D_i \ar[r]^-{d_{e'}} & E_i && H'' \ar@{|->}[d] & \Psi \ar@{|->}[l] \ar@{|->}[d]\\
C \ar[u]^-{d_{tu}} \ar[r]^-{d_u}_-{\text{outer}} & D_i \ar[u]_-{d_t} && f & H_1, \ar@{|->}[l]
}\]
show that the outer face $d_{tu}^* H''$ is $f$.  Here $d_{tu}$ corresponds to deleting the combined vertex of $t$ and $u$, which is almost isolated in $D_i$.
\item
The left commutative square and the computation in $\sK$ on the right-hand side,
\[
\nicearrow
\xymatrix@+8pt{
D_i \ar[r]^-{d_{e'}} & E_i && H'' \ar@{|->}[d] & \Psi \ar@{|->}[l] \ar@{|->}[d]\\
C_{(1;1)} \ar[u]^-{d_{v}} \ar[r]^-{d_{e'}} & L_2 \ar[u]_-{d_v} \ar@<5pt>[l]^-{s} && \bone & s^*(\bone), \ar@{|->}[l]
}\]
show that the other outer face $d_v^* H''$ is $\bone_{x_i}$. Here $s$ is a codegeneracy map and is a section of $d_{e'}$.
\item
The left commutative square and the computation in $\sK$ on the right-hand side,
\[
\nicearrow
\xymatrix{
D_i \ar[r]^-{d_{e'}} & E_i && H'' \ar@{|->}[d] & \Psi \ar@{|->}[l] \ar@{|->}[d]\\
C \ar[u]^-{d_{e}} \ar[r]^-{d_{e'}}_-{\text{inner}} & D_i \ar[u]_-{d_e} && h & H_2, \ar@{|->}[l]
}\]
show that the inner face $d_e^*H''$ is $h$.
\end{enumerate}
We have shown that $H'' \colon f \sim_i h$, so the relation $\sim_i$ is transitive.
\end{proof}

The proof of Lemma \ref{simijrelation} is complete.

\subsection{Same Equivalence Relation}

We now proceed to show that the above equivalence relations are all equal to each other.  We will use the same notations as above.  The following proof is quite similar to that of Lemma \ref{simijrelation}.  As before we will suppress permutations in computation to simplify the presentation.

\begin{lemma}
\label{sameirelation}
Suppose $\sK$ is an $\infty$-properad, and $m,n \geq 0$.  Then the equivalence relations $\sim_i$ for $1 \leq i \leq m$ are all equal to each other.
\end{lemma}

\begin{proof}
We assume $ m \geq 2$, since otherwise there is nothing to prove. To show that $\sim_i$ is equal to $\sim_k$ for $1 \leq i < k \leq m$, suppose given a homotopy $H \colon f \sim_i g$ along the $i$th input.  To show that $f \sim_k g$, consider the $3$-vertex graph
\begin{equation}
\label{graphdik}
\begin{split}
D_{i,k} &= D_i \boxtimes^k_1 C_{(1;1)}\\
&= \left[C \boxtimes^i_1 C_{(1;1)}\right] \boxtimes^k_1 C_{(1;1)}\\
&= \underbrace{\left[C \boxtimes^k_1 C_{(1;1)}\right]}_{D_k} \boxtimes^i_1 C_{(1;1)},
\end{split}
\end{equation}
where $C = C_{(m;n)}$.  The graph $D_{i,k}$ may be depicted as follows.
\begin{center}
\begin{tikzpicture}
\matrix[row sep=1cm, column sep=1cm]{
& \node [plain, label=above:$...$, label=below:\footnotesize{$i~k$}] (v) {$v$}; &\\
\node [plain] (t) {$t$}; 
&& \node [plain] (u) {$u$};\\
};
\draw [inputleg] (v) to +(-.8cm,-.5cm);
\draw [inputleg] (v) to +(.8cm,-.5cm);
\draw [outputleg] (v) to +(-.8cm,.5cm);
\draw [outputleg] (v) to +(.8cm,.5cm);
\draw [inputleg] (t) to +(0,-.7cm);
\draw [arrow] (t) to node[near start]{$e$} (v);
\draw [inputleg] (u) to +(0,-.7cm);
\draw [arrow] (u) to node[swap, near start]{$e'$} (v);
\end{tikzpicture}
\end{center}
First note that $D_{i,k}$ has four faces.
\begin{enumerate}
\item
There is an \emph{outer} face
\[
\nicexy{D_k \ar[r]^-{d_t} & D_{i,k}}
\]
corresponding to deleting the almost isolated vertex $t$. 
\begin{center}
\begin{tikzpicture}
\matrix[row sep=1cm, column sep=1cm]{
\node [plain, label=above:$...$] (v1) {$v$}; &
\node [empty] (s){}; & 
\node [empty] (t){}; && 
\node [plain, label=above:$...$, label=below:\footnotesize{$i~k$}] (v2) {$v$}; &\\
\node [plain] (u1) {$u$}; &&&
\node [plain] (t2) {$t$}; 
&& \node [plain] (u2) {$u$};\\
};
\draw [arrow] (s) to node{$d_t$} (t);
\foreach \x in {1,2}
{
\draw [inputleg] (v\x) to +(-.8cm,-.5cm);
\draw [inputleg] (v\x) to +(.8cm,-.5cm);
\draw [outputleg] (v\x) to +(-.8cm,.5cm);
\draw [outputleg] (v\x) to +(.8cm,.5cm);
}
\draw [inputleg] (u1) to +(0,-.7cm);
\draw [arrow] (u1) to node{$e'$} (v1);
\draw [inputleg] (t2) to +(0,-.7cm);
\draw [arrow] (t2) to node[near start]{$e$} (v2);
\draw [inputleg] (u2) to +(0,-.7cm);
\draw [arrow] (u2) to node[swap, near start]{$e'$} (v2);
\end{tikzpicture}
\end{center}
In particular, there is a self-homotopy
\[
\left(\zeta_k^*(f) \colon f \sim_k f\right) \in \sK(D_k)
\]
from $f$ to $f$ along the $k$th input \eqref{homotopyff}.
\item
There is another \emph{outer} face
\[
\nicexy{D_i \ar[r]^-{d_u} & D_{i,k}}
\]
corresponding to deleting the almost isolated vertex $u$.
\begin{center}
\begin{tikzpicture}
\matrix[row sep=1cm, column sep=1cm]{
\node [plain, label=above:$...$] (v1) {$v$}; &
\node [empty] (s){}; & 
\node [empty] (t){}; && 
\node [plain, label=above:$...$, label=below:\footnotesize{$i~k$}] (v2) {$v$}; &\\
\node [plain] (t1) {$t$}; &&&
\node [plain] (t2) {$t$}; 
&& \node [plain] (u2) {$u$};\\
};
\draw [arrow] (s) to node{$d_u$} (t);
\foreach \x in {1,2}
{
\draw [inputleg] (v\x) to +(-.8cm,-.5cm);
\draw [inputleg] (v\x) to +(.8cm,-.5cm);
\draw [outputleg] (v\x) to +(-.8cm,.5cm);
\draw [outputleg] (v\x) to +(.8cm,.5cm);
}
\draw [inputleg] (t1) to +(0,-.7cm);
\draw [arrow] (t1) to node{$e$} (v1);
\draw [inputleg] (t2) to +(0,-.7cm);
\draw [arrow] (t2) to node[near start]{$e$} (v2);
\draw [inputleg] (u2) to +(0,-.7cm);
\draw [arrow] (u2) to node[swap, near start]{$e'$} (v2);
\end{tikzpicture}
\end{center}
There is a self-homotopy  
\[
\left(\zeta_i^*(f) \colon f \sim_i f\right) \in \sK(D_i)
\]
from $f$ to $f$ along the $i$th input.
\item
There is an \emph{inner} face
\[
\nicexy{D_i \ar[r]^-{d_{e'}} & D_{i,k}}
\]
corresponding to smashing together the closest neighbors $u$ and $v$.
\begin{center}
\begin{tikzpicture}
\matrix[row sep=1cm, column sep=1cm]{
\node [plain, label=above:$...$] (v1) {$uv$}; &
\node [empty] (s){}; & 
\node [empty] (t){}; && 
\node [plain, label=above:$...$, label=below:\footnotesize{$i~k$}] (v2) {$v$}; &\\
\node [plain] (t1) {$t$}; &&&
\node [plain] (t2) {$t$}; 
&& \node [plain] (u2) {$u$};\\
};
\draw [arrow] (s) to node{$d_{e'}$} (t);
\foreach \x in {1,2}
{
\draw [inputleg] (v\x) to +(-.8cm,-.5cm);
\draw [inputleg] (v\x) to +(.8cm,-.5cm);
\draw [outputleg] (v\x) to +(-.8cm,.5cm);
\draw [outputleg] (v\x) to +(.8cm,.5cm);
}
\draw [inputleg] (t1) to +(0,-.7cm);
\draw [arrow] (t1) to node{$e$} (v1);
\draw [inputleg] (t2) to +(0,-.7cm);
\draw [arrow] (t2) to node[near start]{$e$} (v2);
\draw [inputleg] (u2) to +(0,-.7cm);
\draw [arrow] (u2) to node[swap, near start]{$e'$} (v2);
\end{tikzpicture}
\end{center}
By assumption there is a homotopy
\[
\left(H \colon f \sim_i g\right) \in \sK(D_i)
\]
from $f$ to $g$ along the $i$th input.
\item
There is another \emph{inner} face
\[
\nicexy{D_k \ar[r]^-{d_e} & D_{i,k}}
\]
corresponding to smashing together the closest neighbors $t$ and $v$.
\begin{center}
\begin{tikzpicture}
\matrix[row sep=1cm, column sep=1cm]{
\node [plain, label=above:$...$] (v1) {$tv$}; &
\node [empty] (s){}; & 
\node [empty] (t){}; && 
\node [plain, label=above:$...$, label=below:\footnotesize{$i~k$}] (v2) {$v$}; &\\
\node [plain] (u1) {$u$}; &&&
\node [plain] (t2) {$t$}; 
&& \node [plain] (u2) {$u$};\\
};
\draw [arrow] (s) to node{$d_{e}$} (t);
\foreach \x in {1,2}
{
\draw [inputleg] (v\x) to +(-.8cm,-.5cm);
\draw [inputleg] (v\x) to +(.8cm,-.5cm);
\draw [outputleg] (v\x) to +(-.8cm,.5cm);
\draw [outputleg] (v\x) to +(.8cm,.5cm);
}
\draw [inputleg] (u1) to +(0,-.7cm);
\draw [arrow] (u1) to node{$e'$} (v1);
\draw [inputleg] (t2) to +(0,-.7cm);
\draw [arrow] (t2) to node[near start]{$e$} (v2);
\draw [inputleg] (u2) to +(0,-.7cm);
\draw [arrow] (u2) to node[swap, near start]{$e'$} (v2);
\end{tikzpicture}
\end{center}
\end{enumerate}

By Lemma \ref{horndescription}, the homotopies $\zeta_k^*(f)$, $\zeta_i^*(f)$, and $H$ define an inner horn
\[
\nicearrow
\xymatrix@C+12pt{
\Lambda^e[D_{i,k}] \ar[d] \ar[r] & \sK\\
\varGamma[D_{i,k}] \ar@{.>}[ur]_-{\exists \Phi}
}\]
in $\sK$. So there exists a dotted filler $\Phi \in \sK(D_{i,k})$. We claim that the inner face filled by it,
\[
H_k \defn d_e^* \Phi \in \sK(D_k),
\]
is a homotopy from $f$ to $g$ along the $k$th input.  In other words, we need to show that
\begin{itemize}
\item
its outer faces are $f$ and a degenerate element, and
\item
its inner face is $g$.
\end{itemize}
\begin{enumerate}
\item
To see that $H_k$ has $f$ as an outer face, consider the left  commutative square and compute in $\sK$ as indicated.
\[
\nicearrow
\xymatrix{
D_k \ar[r]^-{d_e} & D_{i,k} && H_k \ar@{|->}[d] & \Phi \ar@{|->}[l]  \ar@{|->}[d] \\
C \ar[u]^-{d_u} \ar[r]^-{d_e}_-{\text{inner}} & D_i \ar[u]_-{d_u} && f & \zeta_i^*(f) \ar@{|->}[l]
}\]
\item
To see that the other outer face of $H_k$ is degenerate, consider the left commutative square and compute in $\sK$ as indicated.
\[
\nicearrow
\xymatrix{
D_k \ar[r]^-{d_e} & D_{i,k} && H_k \ar@{|->}[d] & \Phi \ar@{|->}[l]  \ar@{|->}[d] \\
C_{(1;1)} \ar[u]^-{d_{tv}} \ar[r]^-{d_v}_-{\text{outer}} & D_k \ar[u]_-{d_t} && \bone & \zeta_k^*(f) \ar@{|->}[l]
}\]
Here $d_{tv}$ is the outer face corresponding to deleting the combined vertex of $t$ and $v$, which is almost isolated in $D_k$. 
\item
Finally, to see that the inner face of $H_k$ is $g$, consider the left commutative square and compute in $\sK$ as indicated.
\[
\nicearrow
\xymatrix{
D_k \ar[r]^-{d_e} & D_{i,k} && H_k \ar@{|->}[d] & \Phi \ar@{|->}[l]  \ar@{|->}[d] \\
C \ar[u]^-{d_{e'}} \ar[r]^-{d_e}_-{\text{inner}} & D_i \ar[u]_-{d_{e'}} && g & H \ar@{|->}[l]
}\]
Here the left $d_{e'}$ is the inner face of $D_k$.  
\end{enumerate}

We have shown that $H_k \colon f \sim_k g$.  A symmetric argument shows that, if $f \sim_k g$, then $f \sim_i g$.
\end{proof}

\begin{lemma}
\label{samejrelation}
Suppose $\sK$ is an $\infty$-properad, and $m,n \geq 0$.  Then the equivalence relations $\sim^j$ for $1 \leq j \leq n$ are all equal to each other.
\end{lemma}

\begin{proof}
This is dual to the proof of Lemma \ref{sameirelation}.
\end{proof}

\begin{lemma}
\label{sameijrelation}
Suppose $\sK$ is an $\infty$-properad, and $m,n \geq 0$.  Then the equivalence relations 
\begin{itemize}
\item
$\sim_i$ for $1 \leq i \leq m$, and 
\item
$\sim^j$ for $1 \leq j \leq n$ 
\end{itemize}
are all equal to each other.
\end{lemma}

\begin{proof}
We assume that $m,n \geq 1$, since otherwise there is nothing to prove. By Lemmas \ref{sameirelation} and \ref{samejrelation}, it suffices to show that $\sim_1$ and $\sim^1$ are equal.  Suppose given a homotopy $H \colon f \sim_1 g$ along the $1$st input.  To show that $f \sim^1 g$, consider the $3$-vertex graph
\[
\begin{split}
D_1^1 &= D^1 \boxtimes^1_1 C_{(1;1)}\\
&= \left[C_{(1;1)} \boxtimes^1_1 C\right] \boxtimes^1_1 C_{(1;1)}\\
&= C_{(1;1)}  \boxtimes^1_1 
\underbrace{\left[C \boxtimes^1_1 C_{(1;1)}\right]}_{D_1},
\end{split}
\]
where $C = C_{(m;n)}$.  The graph $D^1_1$ may be depicted as follows.
\begin{center}
\begin{tikzpicture}
\matrix[row sep=.2cm, column sep=.8cm]{
\node [plain] (v) {$v$}; &\\
& \node [plain, label=above:$...$, label=below:$...$] (u) {$u$};\\
\node [plain] (t) {$t$}; &\\
};
\draw [inputleg] (t) to +(0,-.7cm);
\draw [outputleg] (v) to +(0,.7cm);
\draw [arrow] (t) to node{$e$} (u);
\draw [arrow] (u) to node[swap]{$e'$}  (v);
\draw [inputleg] (u) to +(.7cm,-.6cm);
\draw [outputleg] (u) to +(.7cm,.6cm);
\end{tikzpicture}
\end{center}
It has the following four faces.
\begin{enumerate}
\item
There is an \emph{outer} face
\[
\nicexy{D_1 \ar[r]^-{d_v} & D^1_1}
\]
corresponding to deleting the almost isolated vertex $v$.
\begin{center}
\begin{tikzpicture}
\matrix[row sep=.2cm, column sep=.8cm]{
&&&& \node [plain] (v2) {$v$}; &\\
& \node [plain, label=above:$...$, label=below:$...$] (u1) {$u$}; & 
\node [empty] (s){}; & \node [empty] (t){};
&&  \node [plain, label=above:$...$, label=below:$...$] (u2) {$u$};\\
\node [plain] (t1) {$t$}; &&&& 
\node [plain] (t2) {$t$}; &\\
};
\draw [arrow] (s) to node{$d_v$} (t);
\draw [inputleg] (t1) to +(0,-.7cm);
\draw [inputleg] (t2) to +(0,-.7cm);
\draw [arrow] (t1) to node{$e$} (u1);
\draw [arrow] (t2) to node{$e$} (u2);
\foreach \x in {1,2}
{
\draw [inputleg] (u\x) to +(.7cm,-.6cm);
\draw [outputleg] (u\x) to +(.7cm,.6cm);
}
\draw [outputleg] (u1) to +(-.7cm,.6cm);
\draw [arrow] (u2) to node[swap]{$e'$} (v2);
\draw [outputleg] (v2) to +(0,.7cm);
\end{tikzpicture}
\end{center}
There is a self-homotopy
\[
\left(\zeta_1^*(f) \colon f \sim_1 f\right) \in \sK(D_1)
\]
from $f$ to $f$ along the first input.
\item
There is an \emph{outer} face
\[
\nicexy{D^1 \ar[r]^-{d_t} & D^1_1}
\]
corresponding to deleting the almost isolated vertex $t$.
\begin{center}
\begin{tikzpicture}
\matrix[row sep=.2cm, column sep=.8cm]{
\node [plain] (v1) {$v$};
&&&& \node [plain] (v2) {$v$}; &\\
& \node [plain, label=above:$...$, label=below:$...$] (u1) {$u$}; & 
\node [empty] (s){}; & \node [empty] (t){};
&&  \node [plain, label=above:$...$, label=below:$...$] (u2) {$u$};\\
&&&& 
\node [plain] (t2) {$t$}; &\\
};
\draw [arrow] (s) to node{$d_t$} (t);
\draw [inputleg] (t2) to +(0,-.7cm);
\draw [arrow] (t2) to node{$e$} (u2);
\foreach \x in {1,2}
{
\draw [inputleg] (u\x) to +(.7cm,-.6cm);
\draw [outputleg] (u\x) to +(.7cm,.6cm);
\draw [arrow] (u\x) to node[swap]{$e'$} (v\x);
\draw [outputleg] (v\x) to +(0,.7cm);
}
\draw [inputleg] (u1) to +(-.7cm,-.6cm);
\end{tikzpicture}
\end{center}
There is a self-homotopy
\[
\left(\zeta^{1*}(f) \colon f \sim^1 f\right) \in \sK(D^1)
\]
from $f$ to $f$ along the first output.
\item
There is an \emph{inner} face
\[
\nicexy{D_1 \ar[r]^-{d_{e'}} & D^1_1}
\]
corresponding to smashing together the closest neighbors $u$ and $v$.
\begin{center}
\begin{tikzpicture}
\matrix[row sep=.2cm, column sep=.8cm]{
&&&& \node [plain] (v2) {$v$}; &\\
& \node [plain, label=above:$...$, label=below:$...$] (u1) {$uv$}; & 
\node [empty] (s){}; & \node [empty] (t){};
&&  \node [plain, label=above:$...$, label=below:$...$] (u2) {$u$};\\
\node [plain] (t1) {$t$}; &&&& 
\node [plain] (t2) {$t$}; &\\
};
\draw [arrow] (s) to node{$d_{e'}$} (t);
\draw [inputleg] (t1) to +(0,-.7cm);
\draw [inputleg] (t2) to +(0,-.7cm);
\draw [arrow] (t1) to node{$e$} (u1);
\draw [arrow] (t2) to node{$e$} (u2);
\foreach \x in {1,2}
{
\draw [inputleg] (u\x) to +(.7cm,-.6cm);
\draw [outputleg] (u\x) to +(.7cm,.6cm);
}
\draw [outputleg] (u1) to +(-.7cm,.6cm);
\draw [arrow] (u2) to node[swap]{$e'$} (v2);
\draw [outputleg] (v2) to +(0,.7cm);
\end{tikzpicture}
\end{center}
By assumption there is a homotopy
\[
\left(H \colon f \sim_1 g\right) \in \sK(D_1)
\]
from $f$ to $g$ along the first input.
\item
There is an \emph{inner} face
\[
\nicexy{D^1 \ar[r]^-{d_{e}} & D^1_1}
\]
corresponding to smashing together the closest neighbors $t$ and $u$.
\begin{center}
\begin{tikzpicture}
\matrix[row sep=.2cm, column sep=.8cm]{
\node [plain] (v1) {$v$};
&&&& \node [plain] (v2) {$v$}; &\\
& \node [plain, label=above:$...$, label=below:$...$] (u1) {$tu$}; & 
\node [empty] (s){}; & \node [empty] (t){};
&&  \node [plain, label=above:$...$, label=below:$...$] (u2) {$u$};\\
&&&& 
\node [plain] (t2) {$t$}; &\\
};
\draw [arrow] (s) to node{$d_e$} (t);
\draw [inputleg] (t2) to +(0,-.7cm);
\draw [arrow] (t2) to node{$e$} (u2);
\foreach \x in {1,2}
{
\draw [inputleg] (u\x) to +(.7cm,-.6cm);
\draw [outputleg] (u\x) to +(.7cm,.6cm);
\draw [arrow] (u\x) to node[swap]{$e'$} (v\x);
\draw [outputleg] (v\x) to +(0,.7cm);
}
\draw [inputleg] (u1) to +(-.7cm,-.6cm);
\end{tikzpicture}
\end{center}
\end{enumerate}

By Lemma \ref{horndescription}, the homotopies $\zeta_1^*(f)$, $\zeta^{1*}(f)$, and $H$ define an inner horn
\[
\nicearrow
\xymatrix@C+12pt{
\Lambda^e[D^1_1] \ar[r] \ar[d] & \sK\\
\varGamma[D^1_1] \ar@{.>}[ur]_-{\exists \Psi}
}\]
in $\sK$.  So there is a dotted filler $\Psi \in \sK(D^1_1)$.  We claim that the inner face filled by it,
\[
H^1 \defn d_e^* \Psi \in \sK(D^1),
\]
is a homotopy from $f$ to $g$ along the first \emph{output}.  In other words, we need to show that
\begin{itemize}
\item
its outer faces are $f$ and a degenerate element, and
\item
its inner face is $g$.
\end{itemize}
The argument is similar to that in the previous few lemmas.
\begin{enumerate}
\item
To see that the outer face $d_v^* H^1$ is $f$, use the left  commutative square and compute in $\sK$ as indicated.
\[
\nicearrow
\xymatrix{
D^1 \ar[r]^-{d_e} & D^1_1 && H^1 \ar@{|->}[d] & \Psi \ar@{|->}[l] \ar@{|->}[d]\\
C \ar[u]^-{d_v} \ar[r]^-{d_e}_-{\text{inner}} & D_1 \ar[u]_-{d_v} && f & \zeta_1^*(f) \ar@{|->}[l]
}\]
\item
To see that the other outer face $d_{tu}^* H^1$ is degenerate, use the left commutative square and compute in $\sK$ as indicated.
\[
\nicearrow
\xymatrix{
D^1 \ar[r]^-{d_e} & D^1_1 && H^1 \ar@{|->}[d] & \Psi \ar@{|->}[l] \ar@{|->}[d]\\
C_{(1;1)} \ar[u]^-{d_{tu}} \ar[r]^-{d_u}_-{\text{outer}} & D^1 \ar[u]_-{d_t} && \bone & \zeta^{1*}(f) \ar@{|->}[l]
}\]
\item
To see that the inner face $d_{e'}^* H^1$ is $g$, use the left commutative square and compute in $\sK$ as indicated.
\[
\nicearrow
\xymatrix{
D^1 \ar[r]^-{d_e} & D^1_1 && H^1 \ar@{|->}[d] & \Psi \ar@{|->}[l] \ar@{|->}[d]\\
C \ar[u]^-{d_{e'}} \ar[r]^-{d_e}_-{\text{inner}} & D_1 \ar[u]_-{d_{e'}} && g & H \ar@{|->}[l]
}\]
\end{enumerate}
This shows that $f \sim_1 g$ implies $f \sim^1 g$.  The converse is proved by a symmetric argument.  Therefore, the equivalence relations $\sim_1$ and $\sim^1$ are equal.
\end{proof}

\begin{definition}
\label{def:homotopyink}
Suppose $\sK$ is an $\infty$-properad. Denote by $\sim$\label{note:sim} the common equivalence relation, called \textbf{homotopy},\index{homotopy} defined by $\sim_i$ and $\sim^j$ as in Lemma \ref{sameijrelation}.  If $m=n=0$, then homotopy is defined as the equality relation. Two elements in the same homotopy class are said to be \textbf{homotopic}. The homotopy class of a $1$-dimensional element $f$ will be written as $[f]$.
\end{definition}

\begin{remark}
For a \emph{strict} $\infty$-properad $\sK$, homotopy is the identity relation.
\end{remark}

\section{Properad Associated to an \texorpdfstring{$\infty$}{∞}-Properad}
\label{sec:propinfprop}

In this section, we define a properad $\qk$ associated to an $\infty$-properad $\sK$ with non-empty inputs, with non-empty outputs, or that is reduced.  The elements in $\qk$ are homotopy classes of $1$-dimensional elements in $\sK$.  First we establish that the properadic composition is well-defined using homotopy classes.  Then we show that the properad $\qk$ has the required universal property for the fundamental properad of an $\infty$-properad.

\subsection{Properadic Composition of \texorpdfstring{$1$}{1}-Dimensional Elements}

For an $\infty$-properad $\sK$, we now define a properadic composition of $1$-dimensional elements.  Its existence is guaranteed by the inner horn extension property, but it is not unique.  We will show that, with suitable restriction on $\sK$, it is well-defined on homotopy classes and yields the fundamental properad of $\sK$.

\begin{definition}
\label{def:compqk}
Suppose:
\begin{itemize}
\item
$\sK$ is an $\infty$-properad (possibly with non-empty inputs or non-empty outputs),
\item
$f \in \sK(\sigma C_{(m;n)}\tau)$ is a $1$-dimensional element with profiles $\bah$,
\item
$g \in \sK(\sigma' C_{(p;q)}\tau')$ is a $1$-dimensional element with profiles $\dch$, and
\item
$\ub \supseteq \ub' = \uc' \subseteq \uc$ are equal $k$-segments for $k>0$.
\end{itemize}
\begin{enumerate}
\item
Define the partially grafted corollas
\[
B = C_{(p;q)} \boxtimes^{\uc'}_{\ub'} C_{(m;n)},
\]
where $m = |\ua|$, $n = |\ub|$, $p = |\uc|$, and $q = |\ud|$. Write $v$ and $u$ for its top and bottom vertices, respectively. Suppose
\[
\nicearrow
\xymatrix{
C = C_{(p+m-k; q+n-k)} \ar[r]^-{d_{\text{in}}} & B
}\]
is the inner face of $B$ corresponding to smashing together the  closest neighbors $v$ and $u$.
\item
A \textbf{properadic composition} \index{properadic composition!of $1$-dimensional elements} of $g$ and $f$ is a $1$-dimensional element $h \in \sK(C)$ such that the following four statements hold.
\begin{enumerate}
\item
There exists an element $\theta \in \sK(B)$.
\item
The composition
\[
\nicearrow
\xymatrix{
\sigma' C_{(p;q)} \tau' \ar[rr]^-{(\tau'^{-1};\sigma'^{-1})}_-{\cong} \ar`u[rrr] `[rrr]^-{d_{\text{top}}} [rrr] && C_{(p;q)} \ar[r]^-{d_v} & B}
\]
satisfies $d_{\text{top}}^*(\theta) = g$.
\item
The composition
\[
\nicearrow
\xymatrix{
\sigma C_{(m;n)} \tau \ar[rr]^-{(\tau^{-1};\sigma^{-1})}_-{\cong} \ar`u[rrr] `[rrr]^-{d_{\text{bot}}} [rrr] && C_{(m;n)} \ar[r]^-{d_u} & B}
\]
satisfies $d_{\text{bot}}^*(\theta) = f$.
\item
The inner face satisfies $d_{\text{in}}^*(\theta) = h$.
\end{enumerate}
In this case, we write
\[
\theta \colon h \simeq g \boxtimes^{\uc'}_{\ub'} f,
\]
and call $\theta$ a \textbf{witness} \index{witness!of properadic composition} of $h$ as a properadic composition of $g$ and $f$. As before, we often abbreviate $\boxtimes^{\uc'}_{\ub'}$ to just $\boxtimes$.
\end{enumerate}
\end{definition}

\begin{remark}
Properadic composition \emph{can} be defined for a general graphical set, since the above definition does not actually use the inner horn extension property of an $\infty$-properad.  However, the description of the fundamental properad in terms of homotopy classes is only valid for $\infty$-properads with suitable restriction.  So we restrict to $\infty$-properads right from the beginning.
\end{remark}

\begin{remark}
A witness $\theta$, even if it exists, is not required to be unique.  That is why we call $h$ \emph{a} properadic composition of $g$ and $f$.  When $\sK$ is an $\infty$-properad, the existence of a witness $\theta$, and hence also $h$, is guaranteed.  We want to show that properadic composition is well-defined using homotopy classes. The $\infty$-category and $\infty$-operad analogs are \cite{bv} (4.12) and \cite{mw2} (6.6).  To accomplish our task, we need the following preliminary observation.  It says that, given $f$ and $g$, any two properadic compositions of $g$ and $f$ are homotopic. In any case, one would expect this to hold if properadic composition really is well-defined using homotopy classes.  The argument is similar to what we used above when we were showing that the relations $\sim_i$ are equivalence relations.
\end{remark}

Recall from Definition \ref{def:gupcgraphicalset} that a graphical set $X$ is \emph{reduced} if the set $X(C_{\emptyprofh})$ is a singleton, where $C_{\emptyprofh}$ is the single isolated vertex.

\begin{lemma}
\label{homotopiccomp}
Suppose $\sK$ is
\begin{itemize}
\item
a reduced $\infty$-properad,
\item
an $\infty$-properad with non-empty inputs, or
\item
an $\infty$-properad with non-empty outputs.
\end{itemize}
Suppose  $f$ and $g$ are $1$-dimensional elements in $\sK$ as in Definition \ref{def:compqk}.  Suppose there are properadic compositions
\[
h \simeq g \boxtimes f \andspace h' \simeq g \boxtimes f.
\]
Then $h$ and $h'$ are homotopic.
\end{lemma}

\begin{proof}
We will use the notations in Definition \ref{def:compqk}. First we assume that $\sK$ is an $\infty$-properad with non-empty outputs. In particular, we have $q > 0$. We will construct a homotopy from $h$ to $h'$ as an inner face of some inner horn filler. Suppose $\ub' = b_{[l,l+k-1]}$, i.e., $\ub'$ begins at the $l$th entry in $\ub$.  Consider the $3$-vertex graph
\[
\begin{split}
A &= C_{(1;1)} \boxtimes^1_l \overbrace{\left[C_{(p;q)} \boxtimes^{\uc'}_{\ub'} C_{(m;n)}\right]}^{B}\\
&= \overbrace{\left[C_{(1;1)} \boxtimes^1_1 C_{(p;q)}\right]}^{D^1} \boxtimes^{\uc'}_{\ub'} C_{(m;n)},
\end{split}
\]
which may be depicted as follows.
\begin{center}
\begin{tikzpicture}
\matrix[row sep=.5cm, column sep=.8cm]{
\node [plain] (w) {$w$}; &\\
& \node [plain, label=above:$...$] (v) {$v$};\\
& \node [empty] (e) {$\be$};\\
& \node [plain, label=below:$...$] (u) {$u$};\\
};
\draw [inputleg] (u) to +(-.7cm,-.5cm);
\draw [inputleg] (u) to node[below right=.1cm]{\footnotesize{$m$}} +(.7cm,-.5cm);
\draw [outputleg] (u) to +(-.7cm,.5cm);
\draw [outputleg] (u) to node[above right=.1cm]{\footnotesize{$n$}}+(.7cm,.5cm);
\draw [inputleg] (v) to +(-.7cm,-.5cm);
\draw [inputleg] (v) to node[below right=.1cm]{\footnotesize{$p$}} +(.7cm,-.5cm);
\draw [outputleg] (v) to node[above right=.1cm]{\footnotesize{$q$}} +(.7cm,.5cm);
\draw [outputleg] (w) to +(0,.7cm);
\draw [arrow, bend left=40] (u) to (v);
\draw [arrow, bend right=40] (u) to (v);
\draw [arrow] (v) to node[swap]{$e$} (w);
\end{tikzpicture}
\end{center}
Observe that $A$ has $k+1$ ordinary edges. The $k$ internal edges from $u$ to $v$ will be collectively denoted by $\be$.

The rest of the proof follows a familiar pattern.  First note that the graph $A$ has four faces.
\begin{enumerate}
\item
There is an \emph{outer} face
\[
\nicexy{D^1 \ar[r]^-{d_u} & A}
\]
corresponding to deleting the almost isolated vertex $u$.
\medskip
\begin{center}
\begin{tikzpicture}
\matrix[row sep=.5cm, column sep=.8cm]{
\node [plain] (w1) {$w$}; &&&& 
\node [plain] (w2) {$w$}; & \\
& \node [plain, label=above:$...$, label=below:$...$] (v1) {$v$};
& \node [empty] (s){};
& \node [empty] (t){}; 
&& \node[plain, label=above:$...$] (v2) {$v$}; \\
&&&&& \node [empty] (e) {$\be$};\\
&&&&& \node [plain, label=below:$...$] (u) {$u$};\\
};
\draw [arrow] (s) to node{$d_u$} (t);
\draw [inputleg] (u) to +(-.7cm,-.5cm);
\draw [inputleg] (u) to +(.7cm,-.5cm);
\draw [outputleg] (u) to +(-.7cm,.5cm);
\draw [outputleg] (u) to +(.7cm,.5cm);
\foreach \x in {1,2}
{
\draw [inputleg] (v\x) to +(-.7cm,-.5cm);
\draw [inputleg] (v\x) to  +(.7cm,-.5cm);
\draw [outputleg] (v\x) to  +(.7cm,.5cm);
\draw [outputleg] (w\x) to +(0,.7cm);
\draw [arrow] (v\x) to node[swap]{$e$} (w\x);
}
\draw [arrow, bend left=40] (u) to (v2);
\draw [arrow, bend right=40] (u) to (v2);
\end{tikzpicture}
\end{center}
There is a self-homotopy
\[
\left(\zeta^*(g) \colon g \sim g\right) \in \sK(D^1)
\]
from $g$ to $g$ along the first output, similar to \eqref{homotopyff}.
\item
There is an \emph{outer} face
\[
\nicexy{B \ar[r]^-{d_w} & A}
\]
corresponding to deleting the almost isolated vertex $w$.
\begin{center}
\begin{tikzpicture}
\matrix[row sep=.5cm, column sep=.8cm]{
&&& 
\node [plain] (w) {$w$}; & \\
\node [plain, label=above:$...$] (v1) {$v$};
& \node [empty] (s){};
& \node [empty] (t){}; 
&& \node[plain, label=above:$...$] (v2) {$v$}; \\
\node [empty] (e) {$\be$};
&&&& \node [empty] (e) {$\be$};\\
\node [plain, label=below:$...$] (u1) {$u$};
&&&& \node [plain, label=below:$...$] (u2) {$u$};\\
};
\draw [arrow] (s) to node{$d_w$} (t);
\foreach \x in {1,2}
{
\draw [inputleg] (u\x) to +(-.7cm,-.5cm);
\draw [inputleg] (u\x) to +(.7cm,-.5cm);
\draw [outputleg] (u\x) to +(-.7cm,.5cm);
\draw [outputleg] (u\x) to +(.7cm,.5cm);
\draw [inputleg] (v\x) to +(-.7cm,-.5cm);
\draw [inputleg] (v\x) to  +(.7cm,-.5cm);
\draw [outputleg] (v\x) to  +(.7cm,.5cm);
\draw [arrow, bend left=40] (u\x) to (v\x);
\draw [arrow, bend right=40] (u\x) to (v\x);
}
\draw [outputleg] (v1) to +(-.7cm,.5cm);
\draw [outputleg] (w) to +(0,.7cm);
\draw [arrow] (v2) to node[swap]{$e$} (w);
\end{tikzpicture}
\end{center}
By assumption there is a witness
\[
\left(\theta \colon h \simeq g \boxtimes f\right) \in \sK(B).
\]
\item
There is an \emph{inner} face
\[
\nicexy{B \ar[r]^-{d_e} & A}
\]
corresponding to smashing together the closest neighbors $v$ and $w$.
\begin{center}
\begin{tikzpicture}
\matrix[row sep=.5cm, column sep=.8cm]{
&&& 
\node [plain] (w) {$w$}; & \\
\node [plain, label=above:$...$] (v1) {$vw$};
& \node [empty] (s){};
& \node [empty] (t){}; 
&& \node[plain, label=above:$...$] (v2) {$v$}; \\
\node [empty] (e) {$\be$};
&&&& \node [empty] (e) {$\be$};\\
\node [plain, label=below:$...$] (u1) {$u$};
&&&& \node [plain, label=below:$...$] (u2) {$u$};\\
};
\draw [arrow] (s) to node{$d_e$} (t);
\foreach \x in {1,2}
{
\draw [inputleg] (u\x) to +(-.7cm,-.5cm);
\draw [inputleg] (u\x) to +(.7cm,-.5cm);
\draw [outputleg] (u\x) to +(-.7cm,.5cm);
\draw [outputleg] (u\x) to +(.7cm,.5cm);
\draw [inputleg] (v\x) to +(-.7cm,-.5cm);
\draw [inputleg] (v\x) to  +(.7cm,-.5cm);
\draw [outputleg] (v\x) to  +(.7cm,.5cm);
\draw [arrow, bend left=40] (u\x) to (v\x);
\draw [arrow, bend right=40] (u\x) to (v\x);
}
\draw [outputleg] (v1) to +(-.7cm,.5cm);
\draw [outputleg] (w) to +(0,.7cm);
\draw [arrow] (v2) to node[swap]{$e$} (w);
\end{tikzpicture}
\end{center}
By assumption there is a witness
\[
\left(\theta' \colon h' \simeq g \boxtimes f\right) \in \sK(B).
\]
\item
There is an \emph{inner} face
\[
\nicexy{D^l \ar[r]^-{d_{\be}} & A}
\]
corresponding to smashing together the closest neighbors $u$ and $v$.
\bigskip
\begin{center}
\begin{tikzpicture}
\matrix[row sep=.5cm, column sep=.8cm]{
\node [plain] (w1) {$w$}; &&&&
\node [plain] (w2) {$w$}; & \\
&& \node [empty] (s){};
& \node [empty] (t){}; 
&& \node[plain, label=above:$...$] (v2) {$v$}; \\
\node [plain, label=below:$...$] (v1) {$uv$};
&&&&& \node [empty] (e) {$\be$};\\
&&&&& \node [plain, label=below:$...$] (u) {$u$};\\
};
\draw [arrow] (s) to node{$d_{\be}$} (t);
\draw [inputleg] (u) to +(-.7cm,-.5cm);
\draw [inputleg] (u) to +(.7cm,-.5cm);
\draw [outputleg] (u) to +(-.7cm,.5cm);
\draw [outputleg] (u) to +(.7cm,.5cm);
\foreach \x in {1,2}
{
\draw [outputleg] (w\x) to +(0,.7cm);
\draw [arrow] (v\x) to node{$e$} (w\x);
}
\draw [inputleg] (v1) to +(-.7cm,-.5cm);
\draw [inputleg] (v1) to node[swap, below right=.1cm]{\footnotesize{$p+m-k$}} +(.7cm,-.5cm);
\draw [outputleg] (v1) to node[swap, above right=.1cm]{\footnotesize{$q+n-k$}} +(.7cm,.5cm);
\draw [outputleg] (v1) to  +(-.7cm,.5cm);
\draw [inputleg] (v2) to +(-.7cm,-.5cm);
\draw [inputleg] (v2) to  +(.7cm,-.5cm);
\draw [outputleg] (v2) to  +(.7cm,.5cm);
\draw [arrow, bend left=40] (u) to (v2);
\draw [arrow, bend right=40] (u) to (v2);
\end{tikzpicture}
\end{center}
\end{enumerate}

By Lemma \ref{horndescription}, the elements $\zeta^*(g)$, $\theta$, and $\theta'$ define an inner horn
\[
\nicearrow
\xymatrix@C+12pt{
\Lambda^{\be}[A] \ar[d] \ar[r] & \sK\\
\varGamma[A] \ar@{.>}[ur]_-{\exists \Phi}
}\]
in $\sK$.  Since $\sK$ is an $\infty$-properad, there exists a dotted filler $\Phi \in \sK(A)$.  We claim that the inner face filled by it,
\[
H \defn d_{\be}^* \Phi \in \sK(D^l),
\]
is a homotopy from $h$ to $h'$.  In other words, we must show that
\begin{itemize}
\item
its outer faces are $h$ and a degenerate element, and
\item
its inner face is $h'$.
\end{itemize}
\begin{enumerate}
\item
To see that $h$ is an outer face of $H$, we use the left commutative square and compute in $\sK$ as indicated.
\[
\nicearrow
\xymatrix{
D^l \ar[r]^-{d_{\be}} & A && H \ar@{|->}[d] & \Phi \ar@{|->}[l] \ar@{|->}[d]\\
C \ar[u]^{d_w} \ar[r]^-{d_{\inp}} & B \ar[u]_-{d_w} && h & \theta \ar@{|->}[l]
}\]
\item
To see that the other outer face of $H$ is degenerate, we use the left commutative square and compute in $\sK$ as indicated.
\[
\nicearrow
\xymatrix{
D^l \ar[r]^-{d_{\be}} & A && H \ar@{|->}[d] & \Phi \ar@{|->}[l] \ar@{|->}[d]\\
C_{(1;1)} \ar[u]^{d_{uv}} \ar[r]^-{d_v}_-{\text{outer}} & D^1 \ar[u]_-{d_u} && \bone & \zeta^*(g) \ar@{|->}[l]
}\]
\item
To see that the inner face of $H$ is $h'$, we use the left commutative square and compute in $\sK$ as indicated.
\[
\nicearrow
\xymatrix{
D^l \ar[r]^-{d_{\be}} & A && H \ar@{|->}[d] & \Phi \ar@{|->}[l] \ar@{|->}[d]\\
C \ar[u]^{d_{\inp}} \ar[r]^-{d_{\inp}} & B \ar[u]_-{d_{e}} && h' & \theta' \ar@{|->}[l]
}\]
\end{enumerate}
We have shown that $H$ is a homotopy from $h$ to $h'$, assuming $\sK$ is an $\infty$-properad with non-empty outputs.

If $\sK$ is an $\infty$-properad with non-empty inputs, then $m>0$.  Suppose $\uc' = c_{[j,j+k-1]}$, i.e., $\uc'$ begins at the $j$th entry in $\uc$.  In this case, we use instead the $3$-vertex graph
\[
\begin{split}
A' 
&= \overbrace{\left[C_{(p;q)} \boxtimes^{\uc'}_{\ub'} C_{(m;n)}\right]}^{B} \boxtimes^j_1 C_{(1;1)} \\
&= C_{(p;q)} \boxtimes^{\uc'}_{\ub'} \overbrace{\left[C_{(m;n)} \boxtimes^1_1 C_{(1;1)}\right]}^{D_1} ,
\end{split}
\]
which may be depicted as follows.
\begin{center}
\begin{tikzpicture}
\matrix[row sep=.5cm, column sep=.8cm]{
& \node [plain, label=above:$...$] (v) {$v$};\\
& \node [empty] (e) {$\be$};\\
& \node [plain, label=below:$...$] (u) {$u$};\\
\node [plain] (w) {$w$}; &\\
};
\draw [inputleg] (u) to node[below right=.1cm]{\footnotesize{$m$}} +(.7cm,-.5cm);
\draw [outputleg] (u) to +(-.7cm,.5cm);
\draw [outputleg] (u) to node[above right=.1cm]{\footnotesize{$n$}}+(.7cm,.5cm);
\draw [inputleg] (v) to +(-.7cm,-.5cm);
\draw [inputleg] (v) to node[below right=.1cm]{\footnotesize{$p$}} +(.7cm,-.5cm);
\draw [outputleg] (v) to +(-.7cm,.5cm);
\draw [outputleg] (v) to node[above right=.1cm]{\footnotesize{$q$}} +(.7cm,.5cm);
\draw [inputleg] (w) to +(0,-.7cm);
\draw [arrow, bend left=40] (u) to (v);
\draw [arrow, bend right=40] (u) to (v);
\draw [arrow] (w) to node{$e$} (u);
\end{tikzpicture}
\end{center}
Using the four faces of $A'$, we then argue as above to obtain a homotopy from $h$ to $h'$ as an inner face of an inner horn filler.

Observe that in the above two cases, we could also have used any variation of $A$ in which the vertex $w$ is attached to any leg of the partially grafted corollas $B$.  Such a leg is guaranteed to exist because both $C_{(p;q)}$ and $C_{(m;n)}$ have non-empty outputs, or both of them have non-empty inputs.

Finally, suppose $\sK$ is a reduced $\infty$-properad.  The profiles of $h$ and $h'$ are the same, say $\yxh$, which is also the pair of profiles of the partially grafted corollas $B$.
\begin{enumerate}
\item
If $\yxh = \emptyprofh$, then $h,h' \in \sK(C_{\emptyprofh})$.  By the reduced assumption on $\sK$, this set is a singleton, so we have $h=h'$.
\item
If $\yxh \not= \emptyprofh$, then the partially grafted corollas $B$ has at least one leg.  Therefore, we may reuse the above argument, in which the vertex $w$ is attached to such a leg, to obtain a homotopy from $h$ to $h'$.
\end{enumerate}
\end{proof}

The following observation will allow us to define properadic composition using homotopy classes of $1$-dimensional elements.

\begin{lemma}
\label{propcompdefined}
Suppose $\sK$, $f$, and $g$ are as in Lemma \ref{homotopiccomp}.  Suppose there are:
\begin{itemize}
\item
homotopies $f \sim f'$ and $g \sim g'$, and
\item
properadic compositions
\[
h \simeq g \boxtimes f \andspace h' \simeq g' \boxtimes f'.
\]
\end{itemize}
Then $h$ and $h'$ are homotopic.
\end{lemma}

\begin{proof}
By Lemma \ref{homotopiccomp} it suffices to show that $h \simeq g' \boxtimes f'$.  As in the proof of that lemma, this witness will be constructed as an inner face of an inner horn filler.

Suppose $\ub' = b_{[l,l+k-1]}$ and $\uc' = \uc_{[j,j+k-1]}$, i.e., $\ub'$ and $\uc'$ begin at, respectively, the $l$th entry in $\ub$ and the $j$th entry in $\uc$.  Consider the $3$-vertex graph
\[
\begin{split}
G &= C_{(p;q)} \boxtimes^{\uc'}_{\ub'} \overbrace{\left[C_{(1;1)} \boxtimes^1_l C_{(m;n)}\right]}^{D^l}\\
&= \underbrace{\left[C_{(p;q)} \boxtimes^j_1 C_{(1;1)}\right]}_{D_j} \boxtimes^{\uc'}_{\ub'} C_{(m;n)},
\end{split}
\]
which may be depicted as follows.
\begin{center}
\begin{tikzpicture}
\matrix[row sep=1.3cm, column sep=1cm]{
\node [plain, label=above:$...$] (v) {$v$};\\
\node [plain] (t) {$t$};\\
\node [plain, label=below:$...$] (u) {$u$};\\
};
\draw [arrow] (u) to node{$e_l$} (t);
\draw [arrow] (t) to node{$e_j$} (v);
\draw [dottedarrow, bend right=45] (u) to node[swap]{$\be$} (v);
\draw [inputleg] (u) to +(-.7cm,-.4cm);
\draw [inputleg] (u) to node[below right=.1cm]{\footnotesize{$m$}} +(.7cm,-.4cm);
\draw [outputleg] (u) to +(-.7cm,.4cm);
\draw [outputleg] (u) to node[right=.1cm]{\footnotesize{$n$}} +(.7cm,.4cm);
\draw [inputleg] (v) to +(-.7cm,-.4cm);
\draw [inputleg] (v) to node[right=.1cm]{\footnotesize{$p$}} +(.7cm,-.4cm);
\draw [outputleg] (v) to +(-.7cm,.4cm);
\draw [outputleg] (v) to node[above right=.1cm]{\footnotesize{$q$}} +(.7cm,.4cm);
\end{tikzpicture}
\end{center}
The dotted arrow named $\be$ represents the collection of ordinary edges corresponding to the equal segments
\[
\left(b_{l+1}, \ldots, b_{l+k-1}\right) 
= 
\left(c_{j+1}, \ldots , c_{j+k-1}\right),
\]
which is empty if and only if $k = 1$.  In particular, if $k > 1$, then $\be$ is not empty, and $G$ is not simply-connected.

The rest of this proof follows the same pattern as in the proofs of the previous few lemmas.  The graph $G$ has the following four faces.
\begin{enumerate}
\item
There is an \emph{outer face}
\[
\nicexy{D^l \ar[r]^-{d_v} & G}
\]
corresponding to deleting the almost isolated vertex $v$.
\bigskip
\begin{center}
\begin{tikzpicture}
\matrix[row sep=1.3cm, column sep=1cm]{
&&& \node [plain, label=above:$...$] (v) {$v$};\\
\node [plain] (t1) {$t$};
& \node [empty] (s){};
& \node [empty] (t){};
& \node [plain] (t2) {$t$};\\
\node [plain, label=below:$...$] (u1) {$u$};
&&&
\node [plain, label=below:$...$] (u2) {$u$};\\
};
\draw [arrow] (s) to node{$d_v$} (t);
\draw [arrow] (u1) to node{$e_l$} (t1);
\draw [outputleg] (t1) to +(0,.7cm);
\foreach \x in {u1,u2,v}
{
\draw [inputleg] (\x) to +(-.7cm,-.4cm);
\draw [inputleg] (\x) to +(.7cm,-.4cm);
\draw [outputleg] (\x) to +(-.7cm,.4cm);
\draw [outputleg] (\x) to +(.7cm,.4cm);
}
\draw [arrow] (u2) to node{$e_l$} (t2);
\draw [arrow] (t2) to node{$e_j$} (v);
\draw [dottedarrow, bend right=45] (u2) to node[swap]{$\be$} (v);
\end{tikzpicture}
\end{center}
By assumption, there is a homotopy
\[
\left(H^l \colon f \sim f'\right) \in \sK(D^l)
\]
from $f$ to $f'$ along the $l$th output.
\item
There is an \emph{outer face}
\[
\nicexy{D_j \ar[r]^-{d_u} & G}
\]
corresponding to deleting the almost isolated vertex $u$.
\begin{center}
\begin{tikzpicture}
\matrix[row sep=1.3cm, column sep=1cm]{
\node [plain, label=above:$...$] (v1) {$v$};
&&& \node [plain, label=above:$...$] (v2) {$v$};\\
\node [plain] (t1) {$t$};
& \node [empty] (s){};
& \node [empty] (t){};
& \node [plain] (t2) {$t$};\\
&&&
\node [plain, label=below:$...$] (u2) {$u$};\\
};
\draw [arrow] (s) to node{$d_u$} (t);
\draw [arrow] (t1) to node{$e_j$} (v1);
\draw [inputleg] (t1) to +(0,-.7cm);
\draw [arrow] (u2) to node{$e_l$} (t2);
\draw [arrow] (t2) to node{$e_j$} (v2);
\draw [dottedarrow, bend right=45] (u2) to node[swap]{$\be$} (v2);
\foreach \x in {u2,v1,v2}
{
\draw [inputleg] (\x) to +(-.7cm,-.4cm);
\draw [inputleg] (\x) to +(.7cm,-.4cm);
\draw [outputleg] (\x) to +(-.7cm,.4cm);
\draw [outputleg] (\x) to +(.7cm,.4cm);
}
\end{tikzpicture}
\end{center}
By assumption, there is a homotopy
\[
\left(H_j \colon g' \sim g\right) \in \sK(D_j)
\]
from $g'$ to $g$ along the $j$th input.
\item
There is an \emph{inner face}
\[
\nicexy{B \ar[r]^-{d_{e_j}} & G}
\]
corresponding to smashing together the closest neighbors $t$ and $v$.
\bigskip
\begin{center}
\begin{tikzpicture}
\matrix[row sep=1.3cm, column sep=1cm]{
\node [plain, label=above:$...$] (v1) {$tv$};
&&&& \node [plain, label=above:$...$] (v2) {$v$};\\
&& \node [empty] (s){};
& \node [empty] (t){};
& \node [plain] (t2) {$t$};\\
\node [plain, label=below:$...$] (u1) {$u$};
&&&& \node [plain, label=below:$...$] (u2) {$u$};\\
};
\draw [arrow] (s) to node{$d_{e_j}$} (t);
\draw [arrow] (u1) to node{$e_l$} (v1);
\draw [dottedarrow, bend right=45] (u1) to node[swap]{$\be$} (v1);
\draw [arrow] (u2) to node{$e_l$} (t2);
\draw [arrow] (t2) to node{$e_j$} (v2);
\draw [dottedarrow, bend right=45] (u2) to node[swap]{$\be$} (v2);
\foreach \x in {u1,u2,v1,v2}
{
\draw [inputleg] (\x) to +(-.7cm,-.4cm);
\draw [inputleg] (\x) to +(.7cm,-.4cm);
\draw [outputleg] (\x) to +(-.7cm,.4cm);
\draw [outputleg] (\x) to +(.7cm,.4cm);
}
\end{tikzpicture}
\end{center}
By assumption, there is a witness
\[
\left(\theta \colon h \simeq g \boxtimes f\right) \in \sK(B)
\]
of $h$ as a properadic composition of $g$ and $f$.
\item
There is an \emph{inner face}
\[
\nicexy{B \ar[r]^-{d_{e_l}} & G}
\]
corresponding to smashing together the closest neighbors $t$ and $u$.  
\bigskip
\begin{center}
\begin{tikzpicture}
\matrix[row sep=1.3cm, column sep=1cm]{
\node [plain, label=above:$...$] (v1) {$v$};
&&&& \node [plain, label=above:$...$] (v2) {$v$};\\
&& \node [empty] (s){};
& \node [empty] (t){};
& \node [plain] (t2) {$t$};\\
\node [plain, label=below:$...$] (u1) {$tu$};
&&&& \node [plain, label=below:$...$] (u2) {$u$};\\
};
\draw [arrow] (s) to node{$d_{e_l}$} (t);
\draw [arrow] (u1) to node{$e_j$} (v1);
\draw [dottedarrow, bend right=45] (u1) to node[swap]{$\be$} (v1);
\draw [arrow] (u2) to node{$e_l$} (t2);
\draw [arrow] (t2) to node{$e_j$} (v2);
\draw [dottedarrow, bend right=45] (u2) to node[swap]{$\be$} (v2);
\foreach \x in {u1,u2,v1,v2}
{
\draw [inputleg] (\x) to +(-.7cm,-.4cm);
\draw [inputleg] (\x) to +(.7cm,-.4cm);
\draw [outputleg] (\x) to +(-.7cm,.4cm);
\draw [outputleg] (\x) to +(.7cm,.4cm);
}
\end{tikzpicture}
\end{center}
\end{enumerate}

By Lemma \ref{horndescription}, the elements $H^l$, $H_j$, and $\theta$ define an inner horn
\[
\nicearrow
\xymatrix@C+12pt{
\Lambda^{e_l}[G] \ar[d] \ar[r] & \sK\\
\varGamma[G] \ar@{.>}[ur]_-{\exists \Psi}
}\]
in $\sK$.  Since $\sK$ is an $\infty$-properad, there exists a dotted filler $\Psi \in \sK(G)$.  We claim that the inner face filled by it,
\[
\psi \defn d_{e_l}^* \Psi \in \sK(B),
\]
is a witness of $h \simeq g' \boxtimes f'$.  In other words, we must show that
\begin{itemize}
\item
its outer faces are $g'$ and $f'$, and
\item
its inner face is $h$.
\end{itemize}
\begin{enumerate}
\item
To see that $g'$ is an outer face of $\psi$, we use the following left commutative square and compute in $\sK$ as indicated.
\[
\nicearrow
\xymatrix{
B \ar[r]^-{d_{e_l}} & G && \psi \ar@{|->}[d] & \Psi \ar@{|->}[l] \ar@{|->}[d]\\
C_{(p;q)} \ar[u]^{d_{tu}} \ar[r]^-{d_t}_-{\text{outer}} & D_j \ar[u]_-{d_u} && g' & H_j \ar@{|->}[l]
}\]
\item
To see that $f'$ is the other outer face of $\psi$, we use the following left commutative square and compute in $\sK$ as indicated.
\[
\nicearrow
\xymatrix{
B \ar[r]^-{d_{e_l}} & G && \psi \ar@{|->}[d] & \Psi \ar@{|->}[l] \ar@{|->}[d]\\
C_{(m;n)} \ar[u]^{d_{v}} \ar[r]^-{d_{e_l}}_-{\text{inner}} & D^l \ar[u]_-{d_v} && f' & H^l \ar@{|->}[l]
}\]
\item
To see that $h$ is the inner face of $\psi$, we use the following left commutative square and compute in $\sK$ as indicated.
\[
\nicearrow
\xymatrix{
B \ar[r]^-{d_{e_l}} & G && \psi \ar@{|->}[d] & \Psi \ar@{|->}[l] \ar@{|->}[d]\\
C \ar[u]^{d_{\inp}} \ar[r]^-{d_{\inp}} & B \ar[u]_-{d_{e_j}} && h & \theta \ar@{|->}[l]
}\]
\end{enumerate}
We have shown that
\[
\psi \colon h \simeq g' \boxtimes f'.
\]
As discussed near the beginning of this proof, this suffices to finish the proof.
\end{proof}

\subsection{Properad of Homotopy Classes}

We now define the object that will be shown to be the fundamental properad of an $\infty$-properad with suitable restriction.  We will use the terminology introduced in section \ref{rk:gsetproperad} for a graphical set.

\begin{definition}
\label{def:qk}
Suppose $\sK$ is as in Lemma \ref{homotopiccomp}.  Define a $\Sigma_{\sS(\sK(\uparrow))}$-bimodule $\qk$ as follows.
\begin{enumerate}
\item
For a pair $\yxh$ of $\sK(\uparrow)$-profiles, denote by $\qk\yxh$, or $\qk\yx$, the set of homotopy classes of $1$-dimensional elements in $\sK$ with profiles $\yxh$.  This is well-defined by Lemma \ref{simprofile}.
\item
For a color $c \in \sK(\uparrow)$, define the \textbf{$c$-colored unit} of $\qk$ as the homotopy class of the degenerate element
\[
\bone_c = s^*(c) \in \sK\left(C_{(1;1)}\right),
\]
where $\nicexy{C_{(1;1)} \ar[r]^-{s} & \uparrow}$ is the codegeneracy map.
\item
Define the $\Sigma$-bimodule structure on $\qk$,
\[
\nicearrow
\xymatrix{
\qk\yx \ar[r]^-{(\pi;\lambda)} & \qk\yxlambda,
}\]
using the isomorphisms in $\sK$ induced by input/output relabelings
\[
\nicearrow
\xymatrix{
\sigma C\tau \ar[r] & \lambda(\sigma C\tau)\pi = (\lambda\sigma) C (\tau\pi)
}\]
of permuted corollas.
\item
Define a \textbf{properadic composition} \index{properadic composition!of homotopy classes} on $\qk$ using representatives of homotopy classes, i.e.,
\[
[g] \boxtimes^{\uc'}_{\ub'} [f] \defn \left[ h \right],
\]
where $h \simeq g \boxtimes^{\uc'}_{\ub'} f$ is as in Definition \ref{def:compqk}.  This is well-defined by Lemma \ref{propcompdefined}.
\end{enumerate}
\end{definition}

\begin{remark}
If, furthermore, $\sK$ is \emph{strict}, then homotopy is the identity relation.  In this case, the object $\qk$ with its structure maps is equal to the properad $\pk$ in Definition \ref{def:properadpk}.
\end{remark}

We first observe that $\qk$ forms a properad.

\begin{lemma}
\label{qkisproperad}
Suppose $\sK$ is as in Lemma \ref{homotopiccomp}.  Then $\qk$ in Definition \ref{def:qk} is a $\sK(\uparrow)$-colored properad.\index{infinity properad!associated properad}
\end{lemma}

\begin{proof}
We need to check the bi-equivariance, unity, and associativity axioms of a properad in biased form.  They are all proved in the same manner using the unity and associativity of graph substitution.  So we will prove only one properadic associativity axiom to illustrate the method.

To this end, we will reuse much of the proof of Lemma \ref{properadpk}.  Suppose $\omega$, $\theta$, and $\phi$ are representatives of homotopy classes of $1$-dimensional elements for which the decorated graph 
\begin{center}
\begin{tikzpicture}
\matrix[row sep=.2cm, column sep=1cm]{
\node [plain] (w) {$\omega$}; &\\
& \node [plain] (v) {$\theta$};\\
\node [plain] (u) {$\phi$}; &\\
};
\draw [implies] (u) to (v);
\draw [implies] (u) to (w);
\draw [implies] (v) to (w);
\end{tikzpicture}
\end{center}
makes sense.  We want to prove the properadic associativity axiom
\[
[\omega] \boxtimes \left([\theta] \boxtimes [\phi]\right) = \left([\omega] \boxtimes [\theta]\right) \boxtimes [\phi]
\]
in $\qk$. We follow the proof and notations in Lemma \ref{properadpk}.
\begin{enumerate}
\item
The elements $\theta$ and $\phi$ determine an inner horn in $\sK$, whose filler $\Phi_0$ is a witness
\[
\Phi_0 : x \defn d_t^* \Phi_0 \simeq \theta \boxtimes \phi
\]
of $x$ as a properadic composition of $\theta$ and $\phi$.
\item
The elements $\omega$ and $\theta$ determine an inner horn in $\sK$, whose filler $\Phi_1$ is a witness
\[
\Phi_1 : y \defn d_r^* \Phi_1 \simeq \omega \boxtimes \theta
\]
of $y$ as a properadic composition of $\omega$ and $\theta$.
\item
The elements $\omega$ and $x$ determine an inner horn in $\sK$, whose filler $\Phi_2$ is a witness
\begin{equation}
\label{za}
\Phi_2 : z \defn d_q^* \Phi_2 \simeq \omega \boxtimes x
\end{equation}
of $z$ as a properadic composition of $\omega$ and $x$.
\end{enumerate}
Note that, unlike the proof of Lemma \ref{properadpk}, these three witnesses are \emph{not} unique because $\sK$ here does not need to be strict.  Nonetheless, their existence is guaranteed by the assumption on $\sK$, and that is all we need.

The witnesses $\Phi_i$ for $0 \leq i \leq 2$ determine an inner horn in $\sK$, which has a filler $\Psi \in \sK(B)$.  The inner face filled by it,
\[
d_{vw}^* \Psi \in \sK(C_{u,(vw)}),
\]
has outer faces $y$ and $\phi$.  So it is a witness
\begin{equation}
\label{zb}
d_{vw}^* \Psi : d_p^*d_{vw}^* \Psi = z \simeq y \boxtimes \phi
\end{equation}
of $z$ as a properadic composition of $y$ and $\phi$.  The equations \eqref{za} and \eqref{zb} and Lemma \ref{propcompdefined} now imply the desired properadic associativity axiom.
\end{proof}

Here is the main observation of this chapter.  It says that the fundamental properad of an $\infty$-properad with suitable restriction can be described using homotopy classes of $1$-dimensional elements.

\begin{theorem}
\label{thm:qkfundamental}
Suppose $\sK$ is as in Lemma \ref{homotopiccomp}.  Then the $\sK(\uparrow)$-colored properad $\qk$ is canonically isomorphic to the fundamental properad of $\sK$.
\end{theorem}

\begin{proof}
If $\sK$ is an $\infty$-properad with non-empty inputs or non-empty outputs, then it is regarded as a properadic graphical set, which is actually an $\infty$-properad, via the left adjoint $i_!$ \eqref{gupcioleftadjoint}.

We will recycle most of the proofs of Lemmas \ref{etaiso}, \ref{etaknerveptwo}, and \ref{strictisnerve}.
\begin{enumerate}
\item
Following the proof of Lemma \ref{etaiso}, we obtain a map
\[
\nicearrow
\xymatrix{
\sK(G) \ar[r]^-{\eta} & (N\qk)(G)
}\]
for each $G \in \gupc$.  The map $\eta^0$ is still a bijection because $\qk$ is $\sK(\uparrow)$-colored.  The map $\eta^1$ sends a $1$-dimensional element in $\sK$ to its homotopy class.
\item
The object-wise map $\eta$ has the same universal property as the unit of the adjunction $(L,N)$ object-wise.  Here we follow the proof of Lemma \ref{etaknerveptwo}.  The main point is that, given a properad $\sQ$ and a map
\[
\nicearrow\xymatrix{\sK \ar[r]^-{\zeta} & N\sQ}
\]
of graphical sets, the map
\[
\nicearrow\xymatrix{\qk \ar[r]^-{\zeta'} & \sQ}
\]
is well-defined.  Indeed, homotopic $1$-dimensional elements in $\sK$ are sent to homotopic $1$-dimensional elements in $N\sQ$.  But since the nerve $N\sQ$ is a strict $\infty$-properad, homotopy is the identity relation.
\item
The object-wise map $\eta$ is actually a map
\[
\nicearrow
\xymatrix{\sK \ar[r]^-{\eta} & N\qk}
\]
of graphical sets.  Here we reuse the proof of Lemma \ref{strictisnerve} by simply replacing $\pk$ with $\qk$ and using the previous step instead of Lemma \ref{etaknerveptwo}.
\end{enumerate}
The last two steps imply that the map $\nicearrow
\xymatrix{\sK \ar[r]^-{\eta} & N\qk}$ of graphical sets has the same universal property as the unit of the adjunction.  So, up to a canonical isomorphism, $\qk$ is the fundamental properad of $\sK$.
\end{proof}



\part{Infinity Wheeled Properads}

\chapter{Wheeled Properads and Graphical Wheeled Properads}
\label{ch:wproperads}

\abstract*{We first recall from \cite{jy2} the biased and the unbiased definitions of a wheeled properad.  There is a symmetric monoidal structure on the category of wheeled properads.  Then we define graphical wheeled properads as free wheeled properads generated by connected graphs, possibly with loops and directed cycles.  With the exception of the exceptional wheel, a graphical wheeled properad has a finite set of elements precisely when the generating graph is simply connected.  So most graphical wheeled properads are infinite.  In the rest of this chapter, we discuss wheeled versions of coface maps, codegeneracy maps, and graphical maps, which are used to define the wheeled properadic graphical category $\Gammaw$.  Every wheeled properadic graphical map has a decomposition into codegeneracy maps followed by coface maps. }

In this chapter we define the graphical category $\Gammaw$ generated by connected graphs.

In section \ref{sec:wproperad} we recall the biased and the unbiased definitions of a wheeled properad.  A wheeled properad is a properad that also has a contraction operation.  Due to the presence of the contraction, in writing down the generating operations of a wheeled properad (i.e., the biased definition), one does not need a general properadic composition.  Instead, a dioperadic composition, which is a very special case of a properadic composition, together with the contraction are sufficient.  The unbiased definition of a wheeled properad describes it as an algebra over a monad defined by the set $\gwheelc$ of connected graphs.  The detailed proof of the equivalence of the two definitions of a wheeled properad is in \cite{jy2}.  We also observe that the category $\wproperad$ of wheeled properads is symmetric monoidal.

In section \ref{sec:gwproperad} we discuss graphical wheeled properads.  They are free wheeled properads generated by connected graphs.  We observe that, with the exception of the exceptional wheel $\wheel$, a graphical wheeled properad has a finite set of elements precisely when it is simply connected (Theorem \ref{thm:infinitegwproperad}).

In section \ref{sec:wcoface} we discuss coface and codegeneracy maps between graphical wheeled properads.  Since there are two generating operations in a wheeled properad besides the bi-equivariant structure and the colored units, there are two types of each of inner coface and outer coface.  Each of these four types of coface maps corresponds to a graph substitution factorization discussed in sections \ref{sec:deletable} and \ref{sec:disconnectable}.  Furthermore, due to the presence of the exceptional wheel, there is an exceptional inner coface map $\bullet \to \wheel$ from a single isolated vertex.  We will prove the wheeled analogs of the graphical identities for these coface and codegeneracy maps.  As in the case of graphical properads, most of these identities are fairly simple, except when two coface maps are composed.

In section \ref{sec:wgraphicalcat} we define the wheeled analogs of graphical maps.  They are the morphisms in the wheeled properadic graphical category $\Gammaw$, whose objects are graphical wheeled properads.  A wheeled properadic graphical map is defined as a map between graphical wheeled properads in which the map from the image to the target is a subgraph.  Therefore, before we can define a wheeled properadic graphical map, we need to discuss subgraphs and images.  The graphical category $\varGamma$ is embedded as a \emph{non-full} subcategory of the wheeled properadic graphical category $\Gammaw$ (Theorem \ref{gammanonfull}).  Also, each map in $\Gammaw$ has a factorization into codegeneracy maps followed by coface maps (Theorem \ref{thm:gwheelcepimono}).
The section ends with various lemmas approaching the uniqueness of such decompositions, which is more subtle than in the wheel-free case.

\section{Biased and Unbiased Wheeled Properads}
\label{sec:wproperad}

In this section we first recall both the biased and the unbiased definitions of a wheeled properad.  Then we construct a symmetric monoidal structure on the category of all wheeled properads.  Later this symmetric monoidal structure will induce a symmetric monoidal closed structure on the category of wheeled properadic graphical sets.

\subsection{Biased Wheeled Properads}

Let us first recall the biased definition of a wheeled properad.  In the linear setting, $1$-colored wheeled properads in biased form were introduced in \cite{mms}. The explicit biased axioms can be found in \cite{jy2}.  The following definition makes sense in any symmetric monoidal category $(\catc, \otimes, I)$ with all small colimits and initial object $\varnothing$ such that $\otimes$ commutes with colimits on both sides.

Fix a set $\fC$ of colors, so $\sS$ means $\SC$.  We will use some definitions about profiles and graphs discussed in chapter \ref{ch:graph}.

\begin{definition}
\label{def:wproperad}
A \textbf{$\fC$-colored wheeled properad}\index{wheeled properad}  $\left(\sP,\bone,\xiij,\jcompi\right)$ \index{wheeled properad} consists of:
\begin{itemize}
\item
a $\Sigma_{\sS}$-bimodule $\sP$,
\item
a $c$-colored unit
\[
\nicearrow
\xymatrix{
I \ar[r]^-{\bone_c} & \sP\ccsingle
}
\]
for each $c \in \fC$, 
\item
a \textbf{contraction} \index{contraction}
\[\label{note:xiij}
\nicearrow
\xymatrix{
\sP\dc \ar[r]^-{\xiij} & \sP\binom{\ud \setminus d_i}{\uc \setminus c_j}
}
\]
whenever $c_j = d_i$, and
\item
a \textbf{dioperadic composition} \index{dioperadic composition}
\[\label{note:jcompi}
\nicearrow
\xymatrix{
\sP\dc \otimes \sP\ba \ar[r]^-{\jcompi} & \sP\binom{\ub\circ_j \ud}{\uc \compi \ua}
}
\]
whenever $b_j = c_i$,
\end{itemize}
such that suitable axioms are satisfied.  

A \textbf{morphism} $\nicexy{\sP \ar[r]^-{f} & \sQ}$ from a $\fC$-colored wheeled properad $\sP$ to a $\fD$-colored wheeled properad $\sQ$ is a map of the underlying colored objects that respects the bi-equivariant structure, colored units, contractions, and  dioperadic compositions.  The category of all wheeled properads and morphisms is denoted by $\wproperad$.\label{note:wproperad}
\end{definition}

\begin{remark}
The colored units and the dioperadic composition can be visualized, in terms of elements, as in Remark \ref{rk:visualizeunits}.  In fact, graphically a dioperadic composition is a properadic composition with only one internal edge.  The contraction can be visualized as follows.
\begin{center}
\begin{tikzpicture}
\matrix[row sep=.8cm,column sep=1cm] {
\node [plain] (p1) {$p$}; 
& \node[empty] (s) {}; 
& \node[empty] (t) {};
& \node [plain, label=above:\footnotesize{$\ud\setminus d_i$}, label=below:\footnotesize{$\uc\setminus c_j$}] (p2) {$\xiij p$}; \\
};
\draw [mapto] (s) to node{\footnotesize{$\xiij$}} (t);
\foreach \x in {1,2}
{
\draw [outputleg] (p\x) to +(-.5cm,.3cm);
\draw [outputleg] (p\x) to +(.5cm,.3cm);
\draw [inputleg] (p\x) to +(-.5cm,-.3cm);
\draw [inputleg] (p\x) to +(.5cm,-.3cm);
}
\draw [arrow, looseness=1, in=-45, out=45, loop] (p1) to node[near start]{\footnotesize{$d_i$}} node[near end]{\footnotesize{$c_j$}} ();
\end{tikzpicture}
\end{center}
\end{remark}

\begin{remark}
\label{rk:wproperadgraph}
There are three sets of biased axioms in a wheeled properad.
\begin{enumerate}
\item
First, the $c$-colored unit and the dioperadic composition $\jcompi$ are generated by the $c$-colored exceptional edge $\uparrow_c$ (Example \ref{ex:exceptionaledge}) and the dioperadic graph $C_{\dch} \jcompi C_{\bah}$ (Example \ref{ex:dioperadic}), respectively.  Together they make $\sP$ into a dioperad.
\item
Second, the contraction $\xiij$ is generated by the contracted corolla $\xiij C_{\dch}$ (Example \ref{ex:contractedcor}).  The contractions are bi-equivariant and commute with each other with suitable shifts of indices.
\item
Finally, the dioperadic composition and the contraction are compatible with each other.  Each such compatibility axiom can be interpreted as saying that a certain connected graph with two vertices and two internal edges can be constructed in two different ways using a dioperadic graph, a contraction, and a relabeling.
\end{enumerate}
\end{remark}

\begin{remark}
Each wheeled properad has an underlying properad (Definition \ref{def:biasedproperad}).  The properadic composition is generated by the partially grafted corollas (Example \ref{ex:pgcor}).  The graph substitution decomposition of the partially grafted corollas in Example \ref{ex:pgcorgeneration} shows that each properadic composition is generated by a dioperadic composition followed by several contractions.  Moreover, the compatibility of the dioperadic composition and the contraction implies that a general properadic composition can be expressed as in the previous sentence in multiple ways.
\end{remark}

\begin{example}
As discussed in \cite{mms} (Example 2.1.1), a typical example of a $1$-colored wheeled properad is the endomorphism object (Definition \ref{def:endobject}) $\End(X)$ of a finite-dimensional vector space $X$.  The contraction on $\End(X)$ is induced by the trace, which is where finite-dimensionality is needed.  The dioperadic composition is induced by composition of maps.
\end{example}

\begin{example}
Over the category of sets, a typical example of a $\fC$-colored wheeled properad is given by the set of $\fC$-colored connected graphs $\gwheelc$ (Definition \ref{def:setsofgraphs}).
\begin{enumerate}
\item
The $\dch$-component is given by the set $\gwheelc\dch$ of connected graphs with input/output profiles $\dch$.
\item
The colored units are given by the exceptional edges $\uparrow_c$ of a single color.
\item
The dioperadic composition is given by grafting an input leg of one connected graph with an output leg with the same color of another connected graph.  So for $G_1 \in \gwheelc\bah$ and $G_2 \in \gwheelc\dch$ with $b_j = c_i$, their dioperadic composition $(G_2) \jcompi (G_1)$ can be visualized as the connected graph in the following picture.
\begin{center}
\begin{tikzpicture}
\matrix[row sep=2cm,column sep=1cm] {
\node [fatplain,label=above:\footnotesize{$\ud$}] (p1) {$G_2$}; \\
\node [fatplain,label=below:\footnotesize{$\ua$}] (p2) {$G_1$}; \\
};
\foreach \x in {1,2}
{
\draw [outputleg] (p\x) to +(-.7cm,.5cm);
\draw [outputleg] (p\x) to +(.7cm,.5cm);
\draw [inputleg] (p\x) to +(-.7cm,-.5cm);
\draw [inputleg] (p\x) to +(.7cm,-.5cm);
}
\draw [arrow] (p2) to node[swap,near end]{\footnotesize{$c_i$}} node[swap,near start]{\footnotesize{$b_j$}} (p1);
\end{tikzpicture}
\end{center}
\item
The contraction is given by connecting an output leg with an input leg of the same color of a given connected graph.  So for $G \in \gwheelc\dch$ with $d_i = c_j$, the contraction $\xiij G$ can be visualized as the connected graph in the following picture.
\begin{center}
\begin{tikzpicture}
\matrix[row sep=.8cm,column sep=.8cm] {
\node [bigplain] (v) {$G$}; \\
};
\draw [outputleg] (v) to +(-.5cm,.3cm);
\draw [outputleg] (v) to +(.5cm,.3cm);
\draw [inputleg] (v) to +(-.5cm,-.3cm);
\draw [inputleg] (v) to +(.5cm,-.3cm);
\draw [arrow, looseness=35, in=-45, out=45, loop] (v) to node[near start]{\footnotesize{$d_i$}} node[near end]{\footnotesize{$c_j$}} ();
\end{tikzpicture}
\end{center}
\end{enumerate}
All the biased axioms of a wheeled properad can be read off from this wheeled properad.
\end{example}

\subsection{Unbiased Wheeled Properads}

As for operads, properads, and other variants, there is a more conceptual, unbiased way to define a wheeled properad using a monad generated by connected graphs.

The following definition is obtained from Definition \ref{def:fgupc} by replacing $\gupc$ with $\gwheelc$.  The notation $\gwheelc\dc$ denotes the subset of $\gwheelc$ consisting of connected graphs with profiles $\dch$.

\begin{definition}
\label{def:unbiasedwproperad}
Suppose  $\sP \in \catcs$.
\begin{enumerate}
\item
Define the functor \index{free wheeled properad} \index{wheeled properad!free}
\[\label{note:fgwheelc}
F = F_{\gwheelc} \colon \catcs \longrightarrow \catcs
\]
by
\begin{equation}
\label{freewproperad}
F\sP\dc 
= \coprod_{G \in \gwheelc\dc} \sP[G] 
= \coprod_{G \in \gwheelc\dc} \bigotimes_{v \in \vertex(G)} \sP\profilev.
\end{equation}
for $\dch \in \sS$.
\item
Define the natural transformation $\nicexy{F^2 \ar[r]^-{\mu} & F}$ as the one induced by graph substitution.
\item
Define the natural transformation $\nicexy{\Id \ar[r]^-{\nu} & F}$ as the one induced by the $\dch$-corollas as $\dch$ runs through $\sS$.
\end{enumerate}
\end{definition}

As in the properad case, the associativity and unity properties of graph substitution imply the following observation.

\begin{theorem}
\label{thm:freewproperadmonad}
For each non-empty set $\fC$, there is a monad $\left(F_{\gwheelc},\mu,\nu\right)$ on $\catcs$.\index{wheeled properad!unbiased}
\end{theorem}

Using the correspondence between the generating operations and the generating graphs in Remark \ref{rk:wproperadgraph}, we have the following equivalence between the biased and the unbiased definitions of a wheeled properad.  Its detailed proof can be found in \cite{jy2}.

\begin{corollary}
\label{cor:wproperadequiv}
There is a natural bijection between $\fC$-colored wheeled properads and $F_{\gwheelc}$-algebras.
\end{corollary}

By general category theory, for a colored object $\sP$, the free $F$-algebra $F\sP$ is precisely the free wheeled properad of $\sP$.

\subsection{Symmetric Monoidal Structure}

The symmetric monoidal product on wheeled properads is defined just like the one on properads (Definition \ref{def:gpropmonoidalproduct}).

\begin{definition}
\label{def:wproperadtensor}
Suppose $\sP$ is a $\fC$-colored wheeled properad and $\sQ$ is a $\fD$-colored wheeled properad.  Define the quotient $\fC \times \fD$-colored wheeled properad\index{tensor product!of wheeled properads} \index{wheeled properad!tensor product}
\begin{equation}
\label{wproperadtensor}
\sP \otimes \sQ 
\defn 
\frac{F_{\gwheelc}(\sP \wedge \sQ)}{\text{3 types of relations}},
\end{equation}
where $F_{\gwheelc}(\sP \wedge \sQ)$ is the free wheeled properad of the colored object $\sP\wedge\sQ$  \eqref{freewproperad}.  The relations are of the following three types.
\begin{enumerate}
\item
For each color $d \in \fD$, the functions
\[
\begin{split}
\fC \ni c & \longmapsto (c,d),\\
\sP \ni p & \longmapsto p \otimes d
\end{split}
\]
are required to define a map of wheeled properads $\sP \longrightarrow \sP \otimes \sQ$.
\item
For each color $c \in \fC$, the functions
\[
\begin{split}
\fD \ni d & \longmapsto (c,d),\\
\sQ \ni q & \longmapsto c \otimes q
\end{split}
\]
are required to define a map of wheeled properads $\sQ \longrightarrow \sP \otimes \sQ$.
\item
Suppose $p \in \sP\left(\ua;\ub\right)$ and $q \in \sQ\left(\uc;\ud\right)$ with $|\ua| = k$, $|\ub| = l$, $|\uc| = m$, $|\ud| = n$, $(k,l) \not= (0,0)$, and $(m,n) \not= (0,0)$.  The relation is then the equality
\begin{equation}
\label{wdistributivity}
\begin{split}
& \{p \otimes d_j\}_{j=1}^n \times \{a_i \otimes q\}_{i=1}^k \\
&=
\sigma^n_l\left[
 \{b_j \otimes q\}_{j=1}^l \times \{p \otimes c_i\}_{i=1}^m
\right]\sigma^m_k
\end{split}
\end{equation}
in $\sP \otimes \sQ$, called \textbf{distributivity}.\index{distributivity}
\end{enumerate}
\end{definition}

\begin{remark}
The distributivity relation is discussed in section \ref{subsec:distributivity}.  Each such relation says that two decorated graphs are equal.  As in the properadic tensor product, if some of the parameters $k,l,m,n$ are $0$, then the distributivity relation requires extra interpretation to keep the decorated graphs connected.
\end{remark}

\begin{theorem}
\label{thm:wproperadmonoidal}
The category $\wproperad$ is symmetric monoidal with respect to $\otimes$.\index{wheeled properad!symmetric monoidal structure}
\end{theorem}

\begin{proof}
It is essentially the same as in the properad case (Theorem \ref{thm:propgmonoidal}) with $\gwheelc$ in place of $\gupc$.  The only difference is that the unit element here is the $1$-colored wheeled properad $\bonew$ with components
\[\label{note:bonew}
\bonew\nm = 
\begin{cases}
\{\uparrow\} & \text{ if $(m,n) = (*,*)$},\\
\{\wheel\} & \text{ if $(m,n) = \emptyprofh$},\\
\varnothing & \text{ otherwise}.
\end{cases}
\]
Its only dioperadic composition is
\[
\uparrow \onecompone \uparrow ~=~ \uparrow.
\]
Its only contraction is $\xi^1_1 \uparrow ~= \wheel$.
\end{proof}

\begin{remark}
\label{rk:freeonewproperad}
The unit $\bonew$ of the symmetric monoidal product is the free $1$-colored wheeled properad generated by the $1$-colored object $\bi$ with components
\[
\bi \nm = \varnothing
\]
for all $(m,n)$.
\end{remark}

\section{Graphical Wheeled Properads}
\label{sec:gwproperad}

In this section we define the graphical wheeled properad generated by a $1$-colored connected graph.  As for graphical properads, most graphical wheeled properads are infinite.  In fact, we observe that, with the exception of the exceptional wheel $\wheel$, the graphical wheeled properad is finite if and only if the generating connected graph is simply connected.

\subsection{Wheeled Properads Generated by Connected Graphs}

Here we define the graphical wheeled properad generated by a connected graph.  The definition is similar to the graphical properad generated by a connected wheel-free graph in section \ref{sec:graphicalproperad}.  The colors are the edges. The generating elements are the vertices.

\begin{remark}
Suppose $G \in \gwheelc$.
\begin{enumerate}
\item
Recall from Definition \ref{def:edges} that $\edge(G)$ is the set of edges in $G$.  An edge in $G$ means an exceptional edge $\uparrow$, an exceptional loop $\wheel$, an ordinary leg, an edge connecting two vertices, or a loop at a vertex.
\item
Using the set $\edge(G)$ as colors, each vertex $v \in \vertex(G)$ (if one exists) determines a corresponding pair of $\edge(G)$-profiles
\[
\profilev \in \sS(\edge(G)) 
= 
\catp(\edge(G))^{op} \times \catp(\edge(G))
\]
Since wheels are now allowed, a loop at a vertex contributes both an incoming flag and an outgoing flag at that vertex.  In this case, $\inp(v)$ and $\out(v)$ have a non-empty intersection.  Moreover, for two distinct vertices in $G$, their corresponding pairs of $\edge(G)$-profiles are different because $G$ cannot have two isolated vertices.
\end{enumerate}
\end{remark}

We now define a graphical wheeled properad.

\begin{definition}
\label{deg:ghatwproperad}
Suppose $G \in \gwheelc$ is a $1$-colored connected graph.
\begin{enumerate}
\item
Define an $\sS(\edge(G))$-colored object $\ghat$ as follows:
\[
\ghat\dc = 
\begin{cases}
\{v\} & \text{ if $\dc = \profilev$ for $\ving$},\\
\varnothing & \text{ otherwise}.
\end{cases}
\]
\item
Define the free $\edge(G)$-colored wheeled properad
\[\label{note:gammawg}
\Gammaw(G) = F_{\gwheelc}\left(\ghat\right),
\]
called the \textbf{graphical wheeled properad} generated by $G$.\index{graphical wheeled properad}
\end{enumerate}
\end{definition}

\begin{remark}
From the definition of the free wheeled properad functor $F_{\gwheelc}$, we can interpret an element in the graphical wheeled properad $\Gammaw(G)$ as an $\edge(G)$-colored $\ghat$-decorated connected graph.  So each edge is colored by an edge in $G$, and each vertex is a vertex in $G$.
\end{remark}

\begin{example}
\label{ex:unitgwproperad}
Suppose
\[
G =~ \uparrow \orspace G = \wheel,
\]
i.e., an exceptional edge (Example \ref{ex:exceptionaledge}) or an exceptional loop (Example \ref{ex:exceptionalloop})
Then
\[
\edge(G) = \{e\} \andspace 
\vertex(G) = \varnothing.
\]
So $\ghat$ is canonically isomorphic to the $1$-colored object $\bi$ with empty components discussed in Remark \ref{rk:freeonewproperad} above.   Therefore, there are canonical isomorphisms
\[
\Gammaw(\uparrow) \cong \bonew \cong \Gammaw(\wheel)
\]
of $1$-colored wheeled properads, where $\bonew$ is the free $1$-colored wheeled properad generated by $\bi$ in the proof of Theorem \ref{thm:wproperadmonoidal}.
\end{example}

\begin{example}
\label{ex:trivialgwproperad}
Consider a single isolated vertex $C_{\emptyprofh} = \bullet$ (Example \ref{ex:isolatedvt}). Then
\[
\edge(\bullet) = \varnothing \andspace 
\vertex(\bullet) = \{v\},
\]
in which the unique vertex $v$ has profiles $\emptyprofh$.  Therefore, the graphical wheeled properad $\Gammaw(\bullet)$ is the free $\varnothing$-colored wheeled properad with components
\[
\Gammaw(\bullet)\dc = 
\begin{cases}
\{\bullet\} & \text{ if $\dch = \emptyprofh$},\\
\varnothing & \text{ otherwise}.
\end{cases}
\]
Since the set $\edge(\bullet)$ of colors is empty, there are no colored units.  There are no non-identity operations in the graphical wheeled properad $\Gammaw(\bullet)$.
\end{example}

\begin{example}
\label{ex:loopgwproperad}
Consider the contracted corolla $G = \xi^1_1C_{(1;1)}$ (Example \ref{ex:contractedcor}) with one vertex $v$, one loop $e$ at $v$, and no other flags.  It can be visualized as follows.
\begin{center}
\begin{tikzpicture}
\matrix[row sep=.5cm,column sep=.3cm] {
\node [plain] (v) {$v$}; \\
};
\draw [arrow, in=-30, out=30, loop] (v) to node{\footnotesize{$e$}} ();
\end{tikzpicture}
\end{center}
Then
\[
\edge(G) = \{e\} \andspace \vertex(G) = \{v\},
\]
in which the unique vertex $v$ has profiles $(e;e)$.  So the graphical wheeled properad $\Gammaw(G)$ has the following components:
\[
\Gammaw(G)\dc = 
\begin{cases}
\left\{L_n\right\}_{n\geq 0} & \text{ if $\dch = (e;e)$},\\
\left\{\xi^1_1 L_n\right\}_{n\geq 0} & \text{ if $\dch = \emptyprofh$},\\
\varnothing & \text{ otherwise}. 
\end{cases}
\]
Here $L_n$ is the linear graph
\begin{center}
\begin{tikzpicture}
\matrix[row sep=.5cm, column sep=1cm]{
\node [empty] (v3) {};\\
\node [plain] (v2) {$v$};\\
\node [empty] (vmid) {$\vdots$}; \\
\node [plain] (v1) {$v$};\\
\node [empty] (v0) {};\\
};
\draw [arrow] (v0) to node{\footnotesize{$e$}} (v1);
\draw [arrow] (v1) to node{\footnotesize{$e$}} (vmid);
\draw [arrow] (vmid) to node{\footnotesize{$e$}} (v2);
\draw [arrow] (v2) to node{\footnotesize{$e$}} (v3);
\end{tikzpicture}
\end{center}
with $n$ copies of $v$, in which each edge is colored by $e$.  In particular, $L_0$ is the $e$-colored exceptional edge $\uparrow_e$, which serves as the $e$-colored unit.

The $\ghat$-decorated graph $\xi^1_1 L_n$ is the contraction of $L_n$, which is depicted as follows.
\begin{center}
\begin{tikzpicture}
\matrix[row sep=.5cm, column sep=1cm]{
\node [plain] (v2) {$v$}; \\
\node [empty] (vmid) {$\vdots$}; \\
\node [plain] (v1) {$v$}; \\
};
\draw [arrow] (v1) to node{\footnotesize{$e$}} (vmid);
\draw [arrow] (vmid) to node{\footnotesize{$e$}} (v2);
\draw [->,thick, looseness=1] (v2) 
to [out=30,in=-30] node{\footnotesize{$e$}} (v1);
\end{tikzpicture}
\end{center}
In particular, $\xi^1_1L_0$ is the $e$-colored exceptional loop $\wheel_e$.  The dioperadic composition on $\Gammaw(G)$ is given by grafting of the linear graphs.  The only contraction that exists in $\Gammaw(G)$ is 
\[
\nicexy{
\Gammaw(G)\binom{e}{e} \ar[r]^-{\xi^1_1} & \Gammaw(G)\emptyprof,
}\]
which takes $L_n$ to $\xi^1_1 L_n$.  Note that both non-empty components of $\Gammaw(G)$ are infinite sets.
\end{example}

\subsection{Size of Graphical Wheeled Properads}

The following observation is the $\gwheelc$ analog of Theorem \ref{omegainfinite}.  It says that most graphical wheeled properads are infinite.

\begin{theorem}
\label{thm:infinitegwproperad}
Suppose $\wheel \not= G \in \gwheelc$.  Then the graphical wheeled properad $\Gammaw(G)$ is a finite set if and only if $G$ is simply connected.\index{graphical wheeled properad!finiteness}
\end{theorem}

\begin{proof}
We reuse much of the proof of Theorem \ref{omegainfinite}.
\begin{enumerate}
\item
Using Example \ref{ex:unitgwproperad} and most of the proof of Lemma \ref{lem:omegagfinite}, we observe that if $G$ is simply connected, then the graphical wheeled properad $\Gammaw(G)$ is a finite set.  The only real change in the argument is that, even if $G$ is simply connected, $\Gammaw(G)$ can still have exceptional loops, which are \emph{not} simply connected, colored by edges of $G$.  However, since $\edge(G)$ is finite in any case, there can only be finitely many such exceptional loops in $\Gammaw(G)$.  Once these exceptional loops are taken into account, the proof of Lemma  \ref{lem:omegagfinite} works here as well.
\item
Next suppose $G$ is not simply connected.  We need to exhibit an infinite set of elements in $\Gammaw(G)$.  Since we are assuming that $G \not= \wheel$, $G$ must have either a loop at some vertex or a cycle involving at least two vertices.
\begin{enumerate}
\item
First suppose $G$ has a cycle involving at least two vertices.  In this case, we reuse the proof of Lemma \ref{lem:omegaginfinite}, which exhibits an infinite list of $\edge(G)$-colored $\ghat$-decorated connected graphs, i.e., an infinite list of elements in $\Gammaw(G)$.
\item
Next, if $G$ does not have a cycle involving at least two vertices, then $G$ must have a loop $e$ at some vertex $v$, which we depict as follows.
\begin{center}
\begin{tikzpicture}
\matrix[row sep=.5cm,column sep=.5cm] {
\node [plain] (v) {$v$};\\
};
\draw [outputleg] (v) to +(-.4cm,.3cm);
\draw [outputleg] (v) to +(.4cm,.3cm);
\draw [inputleg] (v) to +(-.4cm,-.3cm);
\draw [inputleg] (v) to +(.4cm,-.3cm);
\draw [arrow, looseness=25, in=-60, out=60, loop] (v) to node{\footnotesize{$e$}} ();
\draw [dottedarrow, looseness=20, in=200, out=160, loop] (v) to node[swap]{\footnotesize{$l$}} ();
\end{tikzpicture}
\end{center}
The dotted loop $l$ denotes the possibly-empty collection of loops at $v$ other than $e$. Since $e \in \inp(v) \cap \out(v)$, for each positive integer $n$, the following picture depicts an $\edge(G)$-colored $\ghat$-decorated connected graph $H_n$, i.e., an element in $\Gammaw(G)$.
\begin{center}
\begin{tikzpicture}
\matrix[row sep=.8cm, column sep=1cm]{
\node [plain, label=above:$...$] (v2) {$v$};\\
\node [empty] (vmid) {$\vdots$}; \\
\node [plain, label=below:$...$] (v1) {$v$};\\
};
\foreach \x in {1,2}
{
\draw [outputleg] (v\x) to +(-.5cm,.4cm);
\draw [outputleg] (v\x) to +(.5cm,.4cm);
\draw [inputleg] (v\x) to +(-.5cm,-.4cm);
\draw [inputleg] (v\x) to +(.5cm,-.4cm);
}
\draw [arrow] (v1) to node{\footnotesize{$e$}} (vmid);
\draw [arrow] (vmid) to node{\footnotesize{$e$}} (v2);
\foreach \x in {1,2}
{
\draw [dottedarrow, in=200, out=160, loop] (v\x) to node[swap]{\footnotesize{$l$}} ();
}
\end{tikzpicture}
\end{center}
In $H_n$ there are:
\begin{itemize}
\item
$n$ copies of $v$, 
\item
$n-1$ internal $e$-colored edges, each connecting two consecutive copies of $v$,
\item
loops $l$ at each copy of $v$ if $l$ is non-empty.
\end{itemize}
Therefore, $\{H_n\}$ is an infinite list of elements in $\Gammaw(G)$.
\end{enumerate}
\end{enumerate}
\end{proof}

\begin{remark}
The exclusion of the non-simply connected graph $G = \wheel$ in Theorem \ref{thm:infinitegwproperad} is due to Example \ref{ex:unitgwproperad}, where we observed that the graphical wheeled properad $\Gammaw(\wheel)$ is finite.
\end{remark}

\subsection{Maps Between Graphical Wheeled Properads}

Here we describe a wheeled properad map from a graphical wheeled properad or between two graphical wheeled properads.

\begin{lemma}
\label{lem:mapfromgwprop}
Suppose $G \in \gwheelc$, and $\sQ$ is a $\fD$-colored wheeled properad.  Then a map
\[
\nicexy{\Gammaw(G) \ar[r]^-{f} & \sQ}
\]
of wheeled properads is equivalent to a pair of functions:
\begin{enumerate}
\item
A function $\nicexy{\edge(G) \ar[r]^-{f_0} & \fD.}$
\item
A function $f_1$ that assigns to each vertex $v \in \vertex(G)$ an element $f_1(v) \in \sQ\binom{f_0 \out (v)}{f_0 \inp(v)}$.\index{graphical wheeled properad!map from}
\end{enumerate}
\end{lemma}

\begin{proof}
Since $\Gammaw(G) = F(\ghat)$ is the free wheeled properad generated by $\ghat$, a wheeled properad map $f$ is  equivalent to a map $\ghat \to \sQ$ of colored objects, which  consists of a pair of functions as stated.
\end{proof}

If we apply Lemma \ref{lem:mapfromgwprop} when $\sQ$ is a graphical wheeled properad as well, then we obtain the following observation.

\begin{lemma}
\label{lem:mapbtwgwprop}
Suppose $G,H \in \gwheelc$.  Then a map
\[
\nicexy{\Gammaw(G) \ar[r]^-{f} & \Gammaw(H)}
\]
of wheeled properads is equivalent to a pair of functions:
\begin{enumerate}
\item
A function $\nicexy{\edge(G) \ar[r]^-{f_0} & \edge(H).}$
\item
A function $f_1$ that assigns to each vertex $v \in \vertex(G)$ an $\edge(H)$-colored $\hhat$-decorated connected graph
\[
f_1(v) \in \Gammaw(H)\binom{f_0 \out (v)}{f_0 \inp(v)}.
\]
\end{enumerate} 
\end{lemma}

\section{Wheeled Properadic Coface Maps}
\label{sec:wcoface}

In this section we discuss coface and codegeneracy maps between graphical wheeled properads.  There are two types of inner (resp., outer) coface maps, one for dioperadic graphs and one for contracted corollas.  Each of these four types of coface maps corresponds to an inner or outer graph substitution factorization involving a dioperadic graph or a contracted corolla.  There is also an exceptional inner coface map $\bullet \to \wheel$ (Definition \ref{def:exceptionalinner}).

Among the wheeled analogs of the graphical identities, the most interesting one is once again about the composition of two coface maps.  We refer to this as the codimension $2$ property (Theorem \ref{gwccodimtwo}).

\subsection{Motivation for Coface Maps}

Recall from section \ref{sec:coface} that for graphical properads, coface and codegeneracy maps are induced by inner or outer properadic factorizations and degenerate reductions of connected wheel-free graphs.  Both inner and outer coface maps among graphical properads involve graph substitution decompositions involving a partially grafted corollas.  The reason partially grafted corollas play such a crucial role is that a general connected wheel-free graph can always be constructed, via iterated graph substitutions, using partially grafted corollas, permuted corollas, and exceptional edges of a single color.  This is the properadic analog of constructing a unital tree one internal edge at a time.

For graphical wheeled properads, the basic principle for coface maps is the same as before, but the actual maps look quite different.  A general connected graph can always be constructed, via iterated graph substitutions, using exceptional edges of a single color, permuted corollas, dioperadic graphs, and contracted corollas \cite{jy2}.  In other words, besides $\uparrow$ (which corresponds to codegeneracy maps) and permuted corollas (which correspond to change of listings), there are two sets of generating graphs: dioperadic graphs and contracted corollas.  Each set leads to one type of inner coface maps and one type of outer coface maps.  So for graphical wheeled properads, there are two types of inner coface maps, one for dioperadic graphs and one for contracted corollas.  Likewise, there are two types of corresponding outer coface maps.

Inner and outer factorizations corresponding to dioperadic graphs were discussed in section \ref{sec:deletable}.  Recall from Theorem \ref{thm:outerdiopfact} that, for a connected graph $G$, there is an outer dioperadic factorization $G=D(\{C_v,H\})$ (Definition \ref{def:outerdiopfact}) if and only if  $v$ is a deletable vertex in $G$ (Definition \ref{def:deletable}).  Likewise, inner dioperadic factorizations (Definition \ref{def:innerdiopfact}) correspond to internal edges connecting two distinct vertices (Theorem \ref{thm:innerdiopfact}).

Inner and outer factorizations corresponding to contracted corollas were discussed in section \ref{sec:disconnectable}.  Theorem \ref{thm:outercontfact} says that an internal edge $e$ in a connected graph $G$ is disconnectable (Definition \ref{def:disconnectable}) if and only if there is an outer contracting factorization $G=(\xi_eC)(H)$ (Definition \ref{def:outercontfact}).  Likewise, Theorem \ref{thm:innercontfact} says that inner contracting factorizations (Definition \ref{def:innercontfact}) correspond to loops at vertices.

\subsection{Coface and Codegeneracy Maps}

We first define the inner and outer coface maps corresponding to dioperadic graphs.  Then we define inner and outer coface maps corresponding to contracted corollas.  After that we define codegeneracy maps.  All the graphs under discussion are connected.  The reader may wish to review notation \ref{notation:gh} and sections \ref{sec:deletable} and \ref{sec:disconnectable} before reading the following definitions.

\begin{definition}
\label{def:innerdiopcoface}
Suppose $G,K \in \gwheelc$.  An \textbf{inner dioperadic coface map} \index{inner dioperadic coface map}
\[
\nicexy{G \ar[r]^-{d} & K}
\]
is a wheeled properad map between graphical wheeled properads corresponding to an inner dioperadic factorization $K=G(D)$ defined as follows.
\begin{enumerate}
\item
On color sets, the map
\[
\nicexy{
\edge(G) \ar[d]^-{d_0} \\
\edge(G) \coprod \{e\} = \edge(K)
}\]
is the canonical inclusion, where $e$ is the unique internal edge in the dioperadic graph $D$.
\item
For a vertex $v \in \vertex(G)$, 
\[
d_1(v) = 
\begin{cases}
D & \text{ if $v=w$},\\
C_v & \text{ otherwise},\\
\end{cases}
\]
where $w \in \vertex(G)$ is the vertex into which $D$ is substituted, and $C_v$ is the corolla with the same profiles as $v$.
\end{enumerate}
\end{definition}

\begin{remark}
For an inner dioperadic coface map as above, $K$ has one more vertex and one more internal edge than $G$.  In fact, $G$ is obtained from $K$ by shrinking away an internal edge connecting two distinct vertices, namely, the ones in the dioperadic graph $D$. The two vertices in $K$ corresponding to the ones in $D$ become a single vertex in $G$.
\end{remark}

\begin{definition}
\label{def:outerdiopcoface}
Suppose $G,K \in \gwheelc$.  An \textbf{outer dioperadic coface map} \index{outer dioperadic coface map}
\[
\nicexy{G \ar[r]^-{d} & K}
\]
is a wheeled properad map between graphical wheeled properads corresponding to an outer dioperadic factorization $K=D(\{C_v,G\})$ defined as follows.
\begin{enumerate}
\item
On color sets, the map
\[
\nicexy{
\edge(G) \ar[d]^-{d_0}\\ 
\edge(G) \coprod \left[\edge(C_v) \setminus \{e\}\right] = \edge(K)
}\]
is the canonical inclusion.  Here $e$ is the unique internal edge in the dioperadic graph $D$, which is also identified with a leg of the corolla $C_v$.
\item
For each vertex $v \in \vertex(G)$, $d_1(v) = C_v$.
\end{enumerate}
\end{definition}

\begin{remark}
\label{rk:outerdiopedge}
For an outer dioperadic coface map as above, if $G =~ \uparrow$, then $K=C_v$. The outer dioperadic coface map $\uparrow~ \to C_v$ identifies $G$ with a leg named $e$ in $K$.  So $K$ has one more vertex than $G$, but neither $G$ nor $K$ has an internal edge.

On the other hand, suppose $G \not=~ \uparrow$.  Then $G \not= \wheel$ because it must have an input or an output, which in turn is true because it is substituted into $D$.  In this case, $G$ is an ordinary connected graph, and $K$ has one more vertex and one more internal edge than $G$.  In fact, $G$ is obtained from $K$ by deleting a deletable vertex $v$.
\end{remark}

\begin{definition}
\label{def:innercontcoface}
Suppose $G,K \in \gwheelc$.  An \textbf{inner contracting coface map} \index{inner contracting coface map}
\[
\nicexy{G \ar[r]^-{d} & K}
\]
is a wheeled properad map between graphical wheeled properads corresponding to an inner contracting factorization $K=G(\xi_eC)$ defined as follows.
\begin{enumerate}
\item
On color sets, the map
\[
\nicexy{
\edge(G) \ar[d]^-{d_0}\\
\edge(G) \coprod \{e\} = \edge(K)
}\]
is the canonical inclusion, where $e$ is the unique internal edge in the contracted corolla $\xi_eC$.
\item
For a vertex $v \in \vertex(G)$, 
\[
d_1(v) = 
\begin{cases}
\xi_eC & \text{ if $v=w$},\\
C_v & \text{ otherwise},
\end{cases}
\]
where $w \in \vertex(G)$ is the vertex into which $\xi_eC$ is substituted.
\end{enumerate}
\end{definition}

\begin{remark}
For an inner contracting coface map as above, $K$ has one more internal edge (namely $e$) than $G$.  However, since the contracted corolla $\xi_eC$ has only one vertex, $G$ and $K$ have the same numbers of vertices.  In fact, $G$ is obtained from $K$ by deleting the loop $e$.
\end{remark}

For some purposes later, we also need the following exceptional inner coface map.

\begin{definition}
\label{def:exceptionalinner}
The \textbf{exceptional inner coface map} \index{exceptional inner coface map} is the wheeled properad map
\[
\Gammaw(\bullet) \to \Gammaw(\wheel)
\]
defined by sending the unique element in $\Gammaw(\bullet)\emptyprofh = \{\bullet\}$ to the unique element in $\Gammaw(\wheel)\emptyprofh = \{\wheel\}$.
\end{definition}

\begin{remark}
The graphical wheeled properads $\Gammaw(\wheel)$ and $\Gammaw(\bullet)$ were discussed in Examples \ref{ex:unitgwproperad} and \ref{ex:trivialgwproperad}, respectively.  The exceptional wheel $\wheel$ has one internal edge, while the isolated vertex has none.  Neither one of them has a non-isolated vertex.
\end{remark}

\begin{definition}
\label{def:outercontcoface}
Suppose $G,K \in \gwheelc$.  An \textbf{outer contracting coface map} \index{outer contracting coface map}
\[
\nicexy{G \ar[r]^-{d} & K}
\]
is a wheeled properad map between graphical wheeled properads corresponding to an outer contracting factorization $K=(\xi_eC)(G)$ defined as follows.
\begin{enumerate}
\item
On color sets, the map is either the bijection
\[
\nicexy{
\edge(\uparrow) = \{\uparrow\} \ar[r]^-{\cong} 
& \{\wheel\} = \edge(\wheel)
}\]
if $G=~\uparrow$ and $K=\wheel$, or the quotient map
\[
\nicexy{
\edge(G) \ar[r]^-{d_0} 
& \dfrac{\edge(G)}{(e_{-1} \sim e_1)} = \edge(K)
}\]
if $G \not=~ \uparrow$, in which the legs $e_{\pm 1}$ in $G$ are identified to form the internal edge $e$ in $K$.
\item
For each vertex $v \in \vertex(G)$, $d_1(v) = C_v$.
\end{enumerate}
\end{definition}

\begin{remark}
With the exception of the isomorphism $\uparrow~ \to \wheel$, an outer contracting coface map as above is surjective but \emph{not} injective on edges, since the two legs $e_{\pm 1}$ in $G$ are both sent to $e \in \edge(K)$.  In either case, $G$ is obtained from $K$ by disconnecting the disconnectable edge $e$, the numbers of vertices in $G$ and $K$ are the same, and $K$ has one more internal edge (namely $e$) than $G$. 
\end{remark}

\begin{definition}
Suppose $\nicexy{G \ar[r]^-{d} & K}$ is a wheeled properad map between graphical wheeled properads.
\begin{enumerate}
\item
Call $d$ a \textbf{dioperadic coface map} \index{dioperadic coface map} (resp., \textbf{contracting coface map}) \index{contracting coface map} if it is an inner or an outer dioperadic (resp., contracting) coface map.
\item
Call $d$ an \textbf{inner coface map} \index{inner coface map!for graphical wheeled properad} if it is an inner dioperadic coface map, an inner contracting coface map, or the exceptional inner coface map.
\item
Call $d$ an \textbf{outer coface map} \index{outer coface map!for graphical wheeled properad} if it is an outer dioperadic coface map or an outer contracting coface map.
\item
Call $d$ a \textbf{coface map} \index{coface map!for graphical wheeled properad} if it is an inner coface map or an outer coface map.
\end{enumerate}
\end{definition}

\begin{example}
Consider the connected graph $K$:
\begin{center}
\begin{tikzpicture}
\matrix[row sep=1cm, column sep=1cm]{
\node [plain] (v) {$v$};\\
\node [plain] (u) {$u$};\\
};
\draw [inputleg] (u) to node[swap]{\footnotesize{$i$}} +(0,-.7cm);
\draw [arrow] (u) to node{\footnotesize{$a$}} (v);
\draw [outputleg] (v) to node{\footnotesize{$o$}} +(0,.7cm);
\draw [arrow, looseness=5, in=-25, out=25, loop] (v) to node{\footnotesize{$b$}} ();
\end{tikzpicture}
\end{center}
Then there are four coface maps with target $K$.
\begin{enumerate}
\item
There is an inner dioperadic coface map $G^a \to K$ in which $G^a$ is the graph
\begin{center}
\begin{tikzpicture}
\matrix[row sep=1cm, column sep=1cm]{
\node [plain] (v) {$uv$};\\
};
\draw [inputleg] (v) to node[swap]{\footnotesize{$i$}} +(0,-.7cm);
\draw [outputleg] (v) to node{\footnotesize{$o$}} +(0,.7cm);
\draw [arrow, in=-25, out=25, loop] (v) to node{\footnotesize{$b$}} ();
\end{tikzpicture}
\end{center}
obtained from $K$ by shrinking the internal edge $a$.
\item
There is an inner contracting coface $G^b_i \to K$ in which $G^b_i$ is the graph
\begin{center}
\begin{tikzpicture}
\matrix[row sep=1cm, column sep=1cm]{
\node [plain] (v) {$v$};\\
\node [plain] (u) {$u$};\\
};
\draw [inputleg] (u) to node[swap]{\footnotesize{$i$}} +(0,-.7cm);
\draw [arrow] (u) to node{\footnotesize{$a$}} (v);
\draw [outputleg] (v) to node{\footnotesize{$o$}} +(0,.7cm);
\end{tikzpicture}
\end{center}
obtained from $K$ by deleting the loop $b$ at $v$.
\item
There is an outer dioperadic coface map $G^u \to K$ in which $G^u$ is the graph
\begin{center}
\begin{tikzpicture}
\matrix[row sep=1cm, column sep=1cm]{
\node [plain] (v) {$v$};\\
};
\draw [inputleg] (v) to node[swap]{\footnotesize{$a$}} +(0,-.7cm);
\draw [outputleg] (v) to node{\footnotesize{$o$}} +(0,.7cm);
\draw [arrow, in=-25, out=25, loop] (v) to node{\footnotesize{$b$}} ();
\end{tikzpicture}
\end{center}
obtained from $K$ by deleting the deletable vertex $u$.
\item
There is an outer contracting coface map $G^b_o \to K$ in which $G^b_o$ is the graph
\begin{center}
\begin{tikzpicture}
\matrix[row sep=1cm, column sep=1cm]{
\node [plain] (v) {$v$};\\
\node [plain] (u) {$u$};\\
};
\draw [inputleg] (u) to node[swap]{\footnotesize{$i$}} +(0,-.7cm);
\draw [arrow] (u) to node{\footnotesize{$a$}} (v);
\draw [outputleg] (v) to node{\footnotesize{$o$}} +(0,.7cm);
\draw [inputleg] (v) to node[below right=.1cm]{\footnotesize{$b_1$}} +(.6cm,-.5cm);
\draw [outputleg] (v) to node[above right]{\footnotesize{$b_{-1}$}} +(.6cm,.5cm);
\end{tikzpicture}
\end{center}
obtained from $K$ by disconnecting the disconnectable edge $b$.
\end{enumerate}
\end{example}

Next we define a codegeneracy map.

\begin{definition}
\label{def:degenreduction}
Suppose $G \in \gwheelc$ is a $1$-colored graph, and $v \in \vertex(G)$ has exactly one incoming flag and one outgoing flag.  The \textbf{degenerate reduction} \index{degenerate reduction} of $G$ at $v$ is the graph
\[
G_v = G(\uparrow)
\]
obtained from $G$ by substituting the exceptional edge into $v$ and a corolla into every other vertex of $G$.  The edge in $G_v$ corresponding to the two flags adjacent to $v$ is denoted by $e_v$.
\end{definition}

\begin{example}
\label{ex:loopreduction}
Suppose $G = \xi^1_1 C_{(1;1)}$ is the contracted corolla with one vertex $v$, one loop at $v$, and no other flags, as depicted in Example \ref{ex:loopgwproperad}.  Then there is a degenerate reduction
\[
\wheel = \left(\xi^1_1 C_{(1;1)}\right)(\uparrow)
\]
of $G$.
\end{example}

\begin{remark}
Given a degenerate reduction $G_v = G(\uparrow)$, there are canonical bijections
\begin{equation}
\label{reductionbij}
\begin{split}
\vertex(G)
&= \vertex(G_v) \coprod \{v\},\\
\edge(G_v)
&= 
\begin{cases}
\dfrac{\edge(G)}{(v_{\inp} \sim v_{\out})} & \text{ if $v_{\inp} \not= v_{\out}$},\\
\edge(G) & \text{ otherwise},
\end{cases}
\end{split}
\end{equation}
where $v_{\inp}$ and $v_{\out}$ are the incoming and outgoing flags of $v$, regarded as edges in $G$.  The case $\edge(G_v) = \edge(G)$ can only happen if the flags $v_{\inp}$ and $v_{\out}$ form a loop at $v$.  Since $v$ is assumed to have one incoming flag and one outgoing flag, these are the only flags of $v$.  By connectivity $G$ must be the contracted corolla $\xi^1_1 C_{(1;1)}$ and $G_v$ the exceptional loop, as in Example \ref{ex:loopreduction}.
\end{remark}

\begin{definition}
\label{def:wcodegen}
Suppose $G\in \gwheelc$.  A \textbf{codegeneracy map} \index{codegeneracy map}
\[
\nicexy{G \ar[r]^-{s} & G_v}
\]
is a wheeled properad map between graphical wheeled properads corresponding to a degenerate reduction $G_v = G(\uparrow)$ defined as follows.
\begin{enumerate}
\item
The map $\nicexy{\edge(G) \ar[r]^-{s_0} & \edge(G_v)}$ is either the quotient map
\[
\nicexy{
\edge(G) \ar[r] & \dfrac{\edge(G)}{(v_{\inp} \sim v_{\out})} = \edge(G_v)
}\]
if $v_{\inp} \not= v_{\out}$, or the identity map if $(G,G_v) = \left(\xi^1_1C_{(1;1)}, \wheel\right)$, using the bijections \eqref{reductionbij}.
\item
For  $u \in \vertex(G)$,
\[
s_1(u) = 
\begin{cases}
C_u & \text{ if $u\not= v$},\\
\uparrow & \text{ if $u=v$},
\end{cases}
\]
where $\uparrow$ is the $e_v$-colored exceptional edge.
\end{enumerate}
\end{definition}

\begin{example}
The exceptional inner coface map $\bullet \to \wheel$ (Definition \ref{def:exceptionalinner}) is the composition
\[
\nicexy{
\bullet \ar[r]^-{d} 
& \xi^1_1 C_{(1;1)} \ar[r]^-{s}
& \wheel,
}\]
in which $d$ is the inner contracting coface map corresponding to the unique loop in the contracted corolla $\xi^1_1 C_{(1;1)}$, and $s$ is a codegeneracy map.
\end{example}

\subsection{Graphical Identities}

Here we discuss some wheeled graphical analogs of the cosimplicial identities.  Some of the graphical identities for graphical properads (section \ref{sec:graphicalid}) can be reused here.  In particular, noting that a dioperadic graph is an example of a partially grafted corollas, we have the following $\gwheelc$ analogs of the graphical identites for graphical properads.
\begin{enumerate}
\item
The obvious $\gwheelc$ analog of Lemma \ref{graphidss}, involving two codegeneracy maps, still holds.
\item
The inner/outer dioperadic analogs of Lemmas \ref{graphidds} and \ref{graphiddstwo}, which involve an inner/outer properadic coface map and a codegeneracy map, still hold, and the proofs are essentially the same.
\end{enumerate}

The following observation is the inner contracting analog of Lemma  \ref{graphidds}.

\begin{lemma}
\label{wgraphidds}
Suppose:
\begin{enumerate}
\item
$\nicexy{G \ar[r]^-{d^e} & K}$ is an inner contracting coface map corresponding to the inner contracting factorization $K = G(\xi_eC_w)$, and 
\item
$v \in \vertex(K)$ has one incoming flag and one outgoing flag.
\end{enumerate}
Then the following statements hold.
\begin{enumerate}
\item
If $v \not= w$, then the diagram
\[
\nicexy{
G \ar[d]_{s^v} \ar[r]^-{d^e} 
& K \ar[d]^{s^v}\\
G_v \ar[r]^-{d^e} 
& K_v
}\]
is commutative.
\item
If $v=w$, then $K = \xi^1_1C_{(1;1)}$, $G = \bullet$, and the composition
\[
\nicexy{
\bullet \ar[r]^-{d^e} & \xi^1_1C_{(1;1)} \ar[r]^-{s} & \wheel
}\]
is the exceptional inner coface map.
\end{enumerate}
\end{lemma}

\begin{proof}
The first commutative square follows from the graph substitution calculation:
\[
\begin{split}
K_v
&= [G(\xi_eC)](\uparrow)\\
&= G(\{\xi_eC, \uparrow\})\\
&= [G(\uparrow)](\xi_eC)\\
&= G_v(\xi_eC).
\end{split}
\]
For the second assertion, the assumption $v=w$ implies that $C_w=C_{(1;1)}$.  Since $\xi^1_1C_{(1;1)}$ has input/output profiles $\emptyprofh$, by connectivity $G$ itself is an isolated vertex.  So we have
\[
K = (\bullet)\left(\xi^1_1C_{(1;1)}\right) = \xi^1_1C_{(1;1)}.
\]
The unique element in $\Gammaw(\bullet)\emptyprofh$ is sent by the inner contracting coface map $d^e$ to $\xi^1_1C_{(1;1)}$ in $\Gammaw(K)\emptyprofh$, which in turn is sent by the codegeneracy map $s$ to the unique element in $\Gammaw(\wheel)\emptyprofh$.  This composition is by definition the exceptional inner coface map.
\end{proof}

The following observation is the outer contracting analog of Lemma  \ref{graphiddstwo}.

\begin{lemma}
\label{wgraphiddstwo}
Suppose:
\begin{enumerate}
\item
$\nicexy{G \ar[r]^-{d^e} & K}$ is an outer contracting coface map corresponding to the outer contracting factorization $K = (\xi_eC)(G)$, and 
\item
$v \in \vertex(K)$ has one incoming flag and one outgoing flag.
\end{enumerate}
Then the diagram
\[
\nicexy{
G \ar[d]_{s^v} \ar[r]^-{d^e} 
& K \ar[d]^{s^v}\\
G_v \ar[r]^-{d^e} 
& K_v
}\]
is commutative.
\end{lemma}

\begin{proof}
This follows from the calculation:
\[
\begin{split}
K_v 
&= [(\xi_eC)(G)](\uparrow)\\
&= (\xi_eC)[G(\uparrow)]\\
&= (\xi_eC)(G_v).
\end{split}
\]
\end{proof}

\subsection{Codimensional \texorpdfstring{$2$}{2} Property}

As in the case of graphical properads, the most difficult wheeled graphical analog of the cosimplicial identities involves composition of two coface maps.

\begin{definition}
\label{def:wcodimentwo}
We say that $\gwheelc$ has the  \textbf{codimension $2$ property} \index{codimension $2$ property!for graphical wheeled properads} if the following statement holds:  Given any two composable coface maps
\[
\nicearrow
\xymatrix{
\Gammaw(G) \ar[r]^-{d^v} 
& \Gammaw(H) \ar[r]^-{d^u} 
& \Gammaw(K)
}
\]
in $\wproperad$, there exists a commutative square
\[
\nicearrow
\xymatrix{
\Gammaw(G) \ar[r]^-{d^v} \ar[d]_{d^y} 
& \Gammaw(H) \ar[d]^{d^u}\\
\Gammaw(J) \ar[r]^-{d^x} 
&\Gammaw(K)
}
\]
of coface maps such that $d^x$ is not obtainable from $d^u$ by changing the listing.
\end{definition}

\begin{theorem}
\label{gwccodimtwo}
$\gwheelc$ has the codimension $2$ property.
\end{theorem}

\begin{proof}
Suppose
\begin{equation}
\label{codimtwogwc}
\nicearrow
\xymatrix{
\Gammaw(G) \ar[r]^-{d^v} 
& \Gammaw(H) \ar[r]^-{d^u} 
& \Gammaw(K)
}
\end{equation}
are two coface maps.  We must show that the composition $d^ud^v$ has another decomposition into two coface maps.  Neither $d^u$ nor $d^v$ can be the exceptional inner coface map $\bullet \to \wheel$ because there are no coface maps with target $\bullet$ or source $\wheel$. Therefore, there are $16$ cases because each of $d^u$ and $d^v$ can be an inner/outer dioperadic/contracting coface map.  To improve readability, we will check these cases in the next four lemmas.
\end{proof}

\begin{lemma}
\label{lem1:gwccodimtwo}
If $d^u$ is an outer coface map and $d^v$ is an inner coface map in \eqref{codimtwogwc}, then $d^ud^v$ factors as an outer coface map followed by an inner coface map.
\end{lemma}

\begin{proof}
We have graph substitution factorizations
\[
K =  I(H) \andspace 
H = G(J),
\]
in which each of $I$ and $J$ can be either a dioperadic graph or a contracted corolla.  In any case, we have the following decomposition of $K$:
\[
K = I[G(J)] = [I(G)](J).
\]
Therefore, there is a commutative diagram
\[
\nicearrow
\xymatrix{
G \ar[r]^-{d^v} \ar[d]_{\footnotesize{\text{outer}}} 
& H \ar[d]^{d^u}\\
I(G) \ar[r]^-{\footnotesize{\text{inner}}} 
& K
}\]
of coface maps.
\end{proof}

\begin{lemma}
\label{lem2:gwccodimtwo}
If $d^u$ is an inner coface map and $d^v$ is an outer coface map in \eqref{codimtwogwc}, then $d^ud^v$ factors into either
\begin{enumerate}
\item
an inner coface map followed by an outer coface map, or
\item
two outer coface maps.
\end{enumerate}
\end{lemma}

\begin{proof}
We have graph substitution factorizations
\[
K = H(J) \andspace 
H = I(G)
\]
in which each of $I$ and $J$ can be either a dioperadic graph or a contracted corolla.  Note that there is a canonical inclusion $\vertex(G) \subseteq \vertex(H)$.  Suppose $J$ is substituted into $u \in \vertex(H)$.  We have
\[
K = [I(G)](J),
\]
which can be rewritten in one of two ways.
\begin{enumerate}
\item
If $u \in \vertex(G)$, then we have
\[
K = I[G(J)].
\]
Therefore, there is a commutative diagram
\[
\nicearrow
\xymatrix{
G \ar[r]^-{d^v} \ar[d]_{\footnotesize{\text{inner}}} 
& H \ar[d]^{d^u}\\
G(J) \ar[r]^-{\footnotesize{\text{outer}}} 
& K
}\]
of coface maps.
\item
If $u \not\in \vertex(G)$, then
\[
K = I\left(\{G,J\}\right),
\]
and $I$ must be a dioperdic graph because $G$ and $J$ are substituted into different vertices in $I$.
\begin{enumerate}
\item
Suppose $J$ is also a dioperadic graph.  Noting that a dioperadic graph is also a partially grafted corollas, we may therefore reuse case (2) of the proof of Lemma \ref{lem3:gupccodimtwo} here to conclude that $d^ud^v$ factors into two outer coface maps.
\item
Next suppose $J$ is a contracted corolla $\xi_fC_w$.  Then $K$ is obtained from the dioperadic graph $I$ by substituting the contracted corolla $\xi_fC_w$ into one vertex, say $x$, and $G$ into the other vertex.  For example, if $G$ is substituted into the bottom vertex of $I$, then $K$ can be visualized as follows.
\begin{center}
\begin{tikzpicture}
\matrix[row sep=1cm, column sep=1cm]{
\node [plain] (w) {$w$};\\
\node [bigplain] (g) {$G$};\\
};
\draw [arrow] (g) to (w);
\draw [outputleg] (w) to +(-.5cm,.3cm);
\draw [outputleg] (w) to +(.5cm,.3cm);
\draw [inputleg] (w) to +(-.5cm,-.3cm);
\draw [inputleg] (w) to +(.5cm,-.3cm);
\draw [arrow, looseness=5, in=-65, out=65, loop] (w) to node{\footnotesize{$f$}} ();
\end{tikzpicture}
\end{center}
At the vertex $x$ in $I$, there are no incoming/outgoing flags corresponding to $f$ because the latter is an internal edge in the contracted corolla $\xi_fC_w$.  Define the dioperadic graph $D$ as the graph obtained from $I$ by adding the legs $f_{\pm 1}$ to $x$ corresponding to the loop $f$ at $w$.  Then we have
\[
K = (\xi_fC')[D(G)],
\]
where $C'$ is the corolla with the same profiles as $D(G)$.  Therefore,  there is a commutative diagram
\[
\nicearrow
\xymatrix@C+12pt{
G \ar[r]^-{d^v} \ar[d]_{\footnotesize{\text{outer diop.}}}  
& H \ar[d]^{d^u}\\
D(G) \ar[r]^-{\footnotesize{\text{outer}}}_-{\footnotesize{\text{contract.}}} 
& K
}\]
of coface maps.
\end{enumerate}
\end{enumerate}
\end{proof}

\begin{lemma}
\label{lem3:gwccodimtwo}
If $d^u$ and $d^v$ are both outer coface maps in \eqref{codimtwogwc}, then $d^ud^v$ factors into either
\begin{enumerate}
\item
two other outer coface maps, or
\item
an outer coface map followed by an inner coface map.
\end{enumerate}
\end{lemma}

\begin{proof}
We have graph substitution factorizations
\[
K = I(H) \andspace 
H = J(G)
\]
in which each of $I$ and $J$ can be either a dioperadic graph or a contracted corolla.  In any case, we have
\[
K = I[J(G)].
\]
There are four cases.
\begin{enumerate}
\item
Suppose both $I$ and $J$ are dioperadic graphs.  Noting that dioperadic graphs are partially grafted corollas, we may therefore reuse the proof of Lemma \ref{lem2:gupccodimtwo} here to obtain an alternative factorization of $d^ud^v$ into two coface maps.
\item
Suppose both $I$ and $J$ are contracted corollas, say
\[
I = \xi_e C_1 \andspace 
J = \xi_f C_2.
\]
We can visualize $K$ as follows.
\begin{center}
\begin{tikzpicture}
\matrix[row sep=1cm, column sep=1cm]{
\node [bigplain] (g) {$G$};\\
};
\draw [arrow,  looseness=1, in=-25, out=25, loop] (g) to node{\footnotesize{$f$}} ();
\draw [arrow, looseness=5, in=205, out=155, loop] (g) to node[swap]{\footnotesize{$e$}} ();
\end{tikzpicture}
\end{center}
In a doubly-contracted corolla, the two loops can be created by contraction in either order.  So we can rewrite
\[
\begin{split}
K 
&= (\xi_e C_1)[(\xi_fC_2)(G)]\\
&= (\xi_f C_1')[(\xi_e C_2)(G)].
\end{split}
\]
Therefore,  there is a commutative diagram
\[
\nicearrow
\xymatrix@C+12pt{
G \ar[r]^-{d^v} \ar[d]
& H \ar[d]^{d^u}\\
(\xi_e C_2)(G) \ar[r]
& K
}\]
of outer contracting coface maps.
\item
Suppose $I$ is a dioperadic graph, and $J=\xi_fC$ is a contracted corolla.  For example, if $(\xi_fC)(G)$ is substituted into the bottom vertex of $I$, then $K$ can be visualized as follows.
\begin{center}
\begin{tikzpicture}
\matrix[row sep=1cm, column sep=1cm]{
\node [plain, label=above:$...$] (w) {$w$};\\
\node [bigplain] (g) {$G$};\\
};
\draw [arrow] (g) to (w);
\draw [outputleg] (w) to +(-.5cm,.3cm);
\draw [outputleg] (w) to +(.5cm,.3cm);
\draw [inputleg] (w) to +(-.5cm,-.3cm);
\draw [inputleg] (w) to +(.5cm,-.3cm);
\draw [arrow, looseness=1, in=-25, out=25, loop] (g) to node{\footnotesize{$f$}} ();
\end{tikzpicture}
\end{center}
Let $x \in \vertex(I)$ be the vertex into which $(\xi_fC)(G)$ is substituted.  There are no incoming/outgoing flags at $x$ corresponding to the internal edge $f$ of $G$.  Define $D$ as the dioperadic graph obtained from $I$ by adding the legs $f_{\pm 1}$ at $x$.  Then we have
\[
K = (\xi_fC')[D(G)],
\]
where $C'$ is the corolla with the profiles of $D(G)$. Therefore,  there is a commutative diagram
\[
\nicearrow
\xymatrix@C+12pt{
G \ar[r]^-{d^v} \ar[d]_{\footnotesize{\text{outer diop.}}}  
& H \ar[d]^{d^u}\\
D(G) \ar[r]^-{\footnotesize{\text{outer}}}_-{\footnotesize{\text{contract.}}} 
& K
}\]
of outer coface maps.
\item
Suppose $I = \xi_eC$ is a contracted corolla, and $J$ is a dioperadic graph with unique internal edge $f$.  We assume that $G$ is substituted into the bottom vertex $x$ in $J$; there is a similar argument if $G$ is substituted into the top vertex in $J$.   There are four possible shapes for $K$, depicted as follows.
\begin{center}
\begin{tikzpicture}
\matrix[row sep=1cm, column sep=2cm]{
\node [plain] (w1) {$w$};
& \node [plain] (w2) {$w$};
& \node [plain] (w3) {$w$};
& \node [plain] (w4) {$w$};\\
\node [bigplain] (g1) {$G$};
& \node [bigplain] (g2) {$G$};
& \node [bigplain] (g3) {$G$};
& \node [bigplain] (g4) {$G$};\\
};
\draw [arrow] (g1) to node{\footnotesize{$f$}}(w1);
\draw [arrow, looseness=1, in=-25, out=25, loop] (w1) to node{\footnotesize{$e$}} ();
\draw [arrow] (g2) to node{\footnotesize{$f$}}(w2);
\draw [arrow, looseness=1, in=-25, out=25, loop] (g2) to node{\footnotesize{$e$}} ();
\draw [arrow] (g3) to node{\footnotesize{$f$}}(w3);
\draw [arrow, looseness=1, in=-25, out=25] (w3) to node{\footnotesize{$e$}} (g3);
\draw [arrow, bend left=30] (g4) to node{\footnotesize{$f$}}(w4);
\draw [arrow, bend right=30] (g4) to node[swap]{\footnotesize{$e$}}(w4);
\end{tikzpicture}
\end{center}
For simplicity we did not draw the legs at the vertex $w$.  From left to right, we will refer to them as case 1 through case 4.
\begin{enumerate}
\item
For case 1, suppose $D$ is the dioperadic graph obtained from $J$ by deleting the legs $e_{\pm 1}$ at $w$.  Then we have
\[
K = [D(G)](\xi_eC_w),
\]
so there is a commutative diagram
\[
\nicearrow
\xymatrix@C+12pt{
G \ar[r]^-{d^v} \ar[d]_{\footnotesize{\text{outer diop.}}}  
& H \ar[d]^{d^u}\\
D(G) \ar[r]^-{\footnotesize{\text{inner}}}_-{\footnotesize{\text{contract.}}} 
& K
}\]
of coface maps.
\item
For case 2, suppose $D$ is the dioperadic graph obtained from $J$ by deleting the legs $e_{\pm 1}$ at $x$.  Then we have
\[
K = D[(\xi_eC_x)(G)],
\]
so there is a commutative diagram
\[
\nicearrow
\xymatrix@C+12pt{
G \ar[r]^-{d^v} \ar[d]_{\footnotesize{\text{outer contract.}}}  
& H \ar[d]^{d^u}\\
(\xi_eC_x)(G) \ar[r]^-{\footnotesize{\text{outer}}}_-{\footnotesize{\text{diop.}}} 
& K
}\]
of coface maps.
\item
For cases 3 and 4, suppose $D$ is the dioperadic graph with vertices $w$ and $x$, and unique internal edge $e$.  Then we have
\[
K = (\xi_f C')[D(G)],
\]
where $C'$ is the corolla with the profiles of $D(G)$.  Therefore, there is a commutative diagram
\[
\nicearrow
\xymatrix@C+12pt{
G \ar[r]^-{d^v} \ar[d]_{\footnotesize{\text{outer diop.}}}  
& H \ar[d]^{d^u}\\
D(G) \ar[r]^-{\footnotesize{\text{outer}}}_-{\footnotesize{\text{contract.}}} 
& K
}\]
of coface maps.
\end{enumerate}
\end{enumerate}
\end{proof}

\begin{lemma}
\label{lem4:gwccodimtwo}
If $d^u$ and $d^v$ are both inner coface maps in \eqref{codimtwogwc}, then $d^ud^v$ factors into two other inner coface maps.
\end{lemma}

\begin{proof}
We have graph substitution factorizations
\[
K = H(I) \andspace 
H = G(J)
\]
in which each of $I$ and $J$ can be either a dioperadic graph or a contracted corolla. Then we have
\[
K = [G(J)](I).
\]
So $G$ is obtained from $K$ in two steps, each involving either shrinking away an internal edge connecting two distinct vertices, or deleting a loop at some vertex.  In any case, this involves a choice of two distinct internal edges in $K$ and shrinking them away.  The order in which these two internal edges are shrunk away is irrelevant, and either order yields $G$.  Therefore, there is an alternative factorization of $d^ud^v$ into two other inner coface maps.
\end{proof}

With Lemmas \ref{lem1:gwccodimtwo}--\ref{lem4:gwccodimtwo}, the proof of Theorem \ref{gwccodimtwo} is complete.

\begin{remark}
\label{rk:uniquecodimtwoWHEELED}
A careful inspection of the proof of Theorem \ref{gwccodimtwo} reveals that a given composition $d^ud^v$ of coface maps has a \emph{unique} alternative decomposition $d^xd^y$ into coface maps.  Here uniqueness is understood to be up to listing.
\end{remark}

\section{Wheeled Properadic Graphical Category}
\label{sec:wgraphicalcat}

In this section we define the wheeled properadic graphical category $\Gammaw$, whose objects are graphical wheeled properads.  Its morphisms are wheeled properadic graphical maps, which are once again defined using the concepts of subgraphs and images.  While the graphical category $\varGamma$ can be regarded as a subcategory of $\Gammaw$, it is \emph{not} a full subcategory (Theorem \ref{gammanonfull}).  Nevertheless, each wheeled properadic graphical map has a codegeneracies-cofaces decomposition (Theorem \ref{thm:gwheelcepimono}).

\subsection{Subgraphs}

Here we define subgraphs in the wheeled graphical setting. Input/output \textbf{relabeling} of connected graphs and the induced isomorphism on graphical wheeled properads are defined exactly as in section \ref{inoutputrelabeling}. The identity map is by definition a relabeling.  Recall that there are two types of outer coface maps, namely, dioperadic and contracting.

\begin{definition}
\label{def:wsubgraph}
Suppose $G,K \in \gwheelc$.  A map $\nicexy{G \ar[r]^-{f} & K}$ of graphical wheeled properads is called a \textbf{subgraph} \index{subgraph!of graphical wheeled properads} if $f$ admits a decomposition into outer coface maps and relabelings.  In this case, we also call $G$ a \textbf{subgraph of} $K$.
\end{definition}

\begin{remark}
As in the case of graphical properads, when we discuss subgraphs, to simplify the presentation we will sometimes omit mentioning relabeling isomorphisms.
\end{remark}

The next observation is the wheeled analog of Theorem \ref{thm:gupcsubgraph}.  It gives a characterization of subgraphs in terms of graph substitution.

\begin{theorem}
\label{thm:gwheelcsubgraph}
Suppose $\nicexy{G \ar[r]^-{f} & K}$ is a map of graphical wheeled  properads.  Then the following statements are equivalent.
\begin{enumerate}
\item
$f$ is a subgraph.
\item
There exists a graph substitution decomposition
\[
K = H(G)
\]
in $\gwheelc$ such that $f$ sends the edges and vertices in $G$ to their corresponding images in $H(G)$.\index{subgraph!characterization}
\end{enumerate}
\end{theorem}

\begin{proof}
First suppose
\[
f=d_n \cdots d_1
\]
in which each $d_i$ is an outer coface map.  We show by induction on $n$ that $K$ admits a graph substitution decomposition as stated.  If $n=1$, then $f = d_1$ is an outer coface map, and $K=I(G)$ for some dioperadic graph or contracted corolla $I$.

Suppose $n > 1$.  Then there is a factorization
\[
\nicexy{
G \ar[r]^-{g} & H_{n-1} \ar[r]^-{d_n} & K
}\]
of $f$ in which
\[
g = d_{n-1} \cdots d_1.
\]
Since $d_n$ is an outer coface map, there is a graph substitution decomposition
\[
K = I(H_{n-1})
\]
for some dioperadic graph or contracted corolla $I$.  Moreover, $g$ is by definition a subgraph.  By induction hypothesis, there is a graph substitution decomposition
\[
H_{n-1} = H'(G)
\]
in $\gwheelc$.  Therefore, we have the desired decomposition
\[
K = I[H'(G)] = [I(H')](G)
\]
by associativity of graph substitution.

Conversely, suppose $K=H(G)$ as stated.  We show that $f$ is a subgraph by induction on $m=|\edgei(H)|$.  If $m=0$, then since $H$ is ordinary (because it has a vertex), it is actually a corolla $C$.  So
\[
K = C(G) = G,
\]
and $f$ is the identity map, which is a subgraph.

Suppose $m > 0$, so $H$ has at least one internal edge.  Suppose $G$ is substituted into $w \in \vertex(H)$.  There are two cases.
\begin{enumerate}
\item
If $H$ is simply connected, then since $m>0$ it has at least two vertices.  Pick a deletable vertex $v \not= w$ in $H$, which must exist by Corollary \ref{aiexist}.  The graph $H_v$ obtained from $H$ by deleting $v$ is still connected (Lemma \ref{lem:deletable}).  Since $w \in \vertex(H_v)$, the graph substitution $H_v(G)$ makes sense. Moreover, there is an outer dioperadic factorization
\[
H =D(H_v) = D(\{C_v,H_v\})
\]
for some dioperadic graph $D$.  So we have
\[
K = [D(H_v)](G) = D[H_v(G)],
\]
and $f$ factors as
\[
\nicexy{
G \ar[r] & H_v(G) \ar[r] & D[H_v(G)] = K.
}\]
Since $|\edgei(H_v)| < m$, by induction hypothesis the first map is a subgraph.  The second map is by definition an outer dioperadic coface map.  Therefore, the composition $f$ is also a subgraph.
\item
Next suppose $H$ is not simply connected, so $H$ has a cycle.  Pick an internal edge $\nicexy{u \ar[r]^-{e} & v}$ in a cycle in $H$.  It is possible that $e$ is a loop at $u$.  By Lemma \ref{disconnectablecycle} $e$ is disconnectable, so the graph $H_e$ obtained from $H$ by disconnecting $e$ is connected.  The graph substitution $H_e(G)$ still makes sense.  Furthermore, there is an outer contracting factorization
\[
H = (\xi_e C)(H_e)
\]
by Theorem \ref{thm:outercontfact}.  So we have
\[
K = [(\xi_e C)(H_e)](G) = (\xi_e C)[H_e(G)],
\]
and $f$ factors as
\[
\nicexy{
G \ar[r] & H_e(G) \ar[r] & (\xi_e C)[H_e(G)] = K.
}\]
Since $|\edgei(H_e)| < m$, by induction hypothesis the first map is a subgraph.  The second map is by definition an outer contracting coface map.  Therefore, the composition $f$ is also a subgraph.
\end{enumerate}
\end{proof}

The following observation provides some small examples of subgraphs.

\begin{corollary}
\label{cor:gwheelcsubgraph}
Suppose $K$ is a connected graph.
\begin{enumerate}
\item
For each vertex $v$ in $K$, the corolla inclusion $C_v \to K$ is a subgraph.
\item
For each edge $e$ in $K$, the edge inclusion $\uparrow_e ~\to K$ is a subgraph.
\end{enumerate}
\end{corollary}

\begin{proof}
For the first assertion, use the graph substitution decomposition
\[
K = K(C_v),
\]
in which a corolla is substituted into each vertex, and Theorem \ref{thm:gwheelcsubgraph}.  

Consider the second assertion.
\begin{enumerate}
\item
If $K$ is an exceptional edge, then $\uparrow_e ~\to K$ is the identity map.
\item
If $K$ is an exceptional wheel, then $\uparrow ~\to \wheel$ is an outer contracting coface map, and hence a subgraph.
\item
If $K \not\in \{\uparrow, \wheel\}$, then it is ordinary, so $e$ is adjacent to some vertex $v$ in $K$.  The map $\uparrow_e ~\to K$ factors into
\[
\nicexy{
\uparrow_e \ar[r] & C_v \ar[r] & K.
}\]
The first map is an outer dioperadic coface map that identifies the edge $e$ as a leg of the corolla $C_v$ (Remark \ref{rk:outerdiopedge}).  The second map is a corolla inclusion, which is a subgraph by the previous part.  Therefore, their composition is also a subgraph.
\end{enumerate}
\end{proof}

\begin{remark}
The wheeled-analogue of Lemma \ref{subgraphunique} is false: subgraphs are not uniquely determined by their input/output profiles.
As an example, consider the contracted corolla $G = \xi^1_1C_{(1;1)}$ (Example \ref{ex:contractedcor}) with one vertex $v$, one loop $e$ at $v$, and no other flags. 
\begin{center}
\begin{tikzpicture}
\matrix[row sep=.5cm,column sep=.3cm] {
\node [plain] (v) {$v$}; \\
};
\draw [arrow, in=-30, out=30, loop] (v) to node{\footnotesize{$e$}} ();
\end{tikzpicture}
\end{center}
Since $e$ is disconnectable in $G$, the corolla $C_{(e;e)}$ is a subgraph of $G$ with profiles $(e;e)$. On the other hand, the exceptional edge $\uparrow_e$ is also a subgraph of $G$ with profiles $(e;e)$, so $G$ has two distinct subgraphs with the same profiles.
\end{remark}




\subsection{Images}

Here we define image in the setting of graphical wheeled properads.

\begin{definition}
Suppose $\nicexy{G \ar[r]^{f} & K}$ is a wheeled properad map of graphical wheeled  properads.  The \textbf{image of $G$} \index{image!of graphical wheeled properad} is defined as the graph substitution
\[\label{note:fofg}
f(G) = \left[f_0 G\right] \left(\{f_1(u)\}_{u \in \vertex(G)}\right) \in \Gammaw(K),
\]
in which $f_0G$ is the graph obtained from $G$ by applying $f_0$ to its edges.
\end{definition}

\begin{example}
\label{ex:gwheelcimage}
Here we describe the images of codegeneracy maps, coface maps, subgraphs, and changes of vertex listings.
\begin{enumerate}
\item
For a codegeneracy map
\[
\nicexy{G \ar[r]^-{s} & G_v},
\]
by definition we have $G_v = G(\uparrow)$.  So the image $s(G)$ is $G_v$.
\item
For the exceptional inner coface map $\bullet \to \wheel$, the image is $\wheel$.
\item
For a non-exceptional inner coface map
\[
\nicexy{G \ar[r]^-{d_{\inp}} & K},
\]
by definition we have $K=G(I)$ for some dioperadic graph or contracted corolla $I$.  So the image $d_{\inp}(G)$ is $K$.
\item
For an outer coface map
\[
\nicexy{G \ar[r]^-{d_{\out}} & K},
\]
by definition we have $K=I(G)$ for some dioperadic graph or contracted corolla $I$.  Each vertex $v$ in $G$ is sent to a corolla $C_v$, so the image $d_{\out}(G)$ is $G \in \Gammaw(K)$, where $G$ is regarded as an $\edge(K)$-colored $\khat$-decorated graph via the maps
\[
\edge(G) \to \edge(K) \andspace
\vertex(G) \hookrightarrow \vertex(K).
\]
\item
For a subgraph 
\[
\nicexy{G \ar[r]^-{f} & K},
\]
by Theorem \ref{thm:gwheelcsubgraph} we have $K = H(G)$ for some $H \in \gwheelc$.  So the image $f(G)$ is $G \in \Gammaw(K)$.
\item
Suppose $K$ is obtained from $G$ by changing the listings at a subset of vertices, and
\[
\nicexy{G \ar[r]^-{f} & K}
\]
is the corresponding wheeled properad isomorphism.  Then the image $f(G)$ is $K$.
\end{enumerate}
\end{example}

The following observation is the wheeled analog of Lemma \ref{lem1:gupcimage}.  It says that a maps between graphical wheeled properads factors through the image.

\begin{lemma}
\label{lem1:gwheelcimage}
Suppose $\nicexy{G \ar[r]^{f} & K}$ is a wheeled properad map of graphical wheeled properads.  Then there is a canonical commutative diagram
\[
\nicexy@C+10pt{
G \ar[r]^-{g} \ar[dr]_-{f} & f(G) \ar[d]^-{h}\\
& K}
\]
of wheeled properad maps between graphical wheeled properads.
\end{lemma}

\begin{proof}
The proof is the same as that of Lemma \ref{lem1:gupcimage} with $\Gammaw(-)$ in place of $\varGamma(-)$.
\end{proof}

\begin{example}
\label{ex:wcofacesubgraph}
If $\nicexy{G \ar[r]^-{f} & K}$ is a codegeneracy map, a coface map, a subgraph, or a change of vertex listings, then the map $f(G) \to K$ is a subgraph.  Indeed, if $f$ is a codegeneracy map, an inner coface map (exceptional or not), or a change of vertex listings, then $f(G)=K$.  If $f$ is an outer coface map or, more generally, a subgraph, then $f(G)$ is $G$.
\end{example}

\subsection{Graphical Maps}

Now we define the wheeled analog of properadic graphical maps.

\begin{definition}
\label{def:wgraphicalmap}
A \textbf{wheeled properadic graphical map}, \index{wheeled properadic graphical map} \index{graphical map!wheeled properadic} or simply a \textbf{graphical map}, is defined as a wheeled properad map $\nicexy{G \ar[r]^-{f} & K}$ between graphical wheeled properads such that the map $f(G) \to K$ is a subgraph.
\end{definition}

\begin{example}
Codegeneracy maps, coface maps, subgraphs, and changes of vertex listings are all graphical maps by Example \ref{ex:wcofacesubgraph}. 
\end{example}

The following observation says that graphical maps are closed under compositions.

\begin{lemma}
\label{lem1:wgraphicalmapclosed}
Suppose $\nicexy{G \ar[r]^-{f} & K}$ and $\nicexy{K \ar[r]^-{g} & M}$ are wheeled properadic graphical maps.  Then the composition $\nicexy{G \ar[r]^-{gf} & M}$ is also a wheeled properadic graphical map.
\end{lemma}

\begin{proof}
The proof is the same as that of Lemma \ref{lem1:graphicalmapclosed}, except that the graph substitution characterization Theorem \ref{thm:gwheelcsubgraph} of subgraphs is used.
\end{proof}

Next we give another characterization of a wheeled properadic graphical map.  It says that a graphical map sends every subgraph of the source to a subgraph of the target.

\begin{theorem}
\label{thm:wgraphicalmapchar}
Suppose $\nicexy{G \ar[r]^-{f} & K}$ is a wheeled properad map between graphical wheeled properads.  Then the following statements are equivalent.
\begin{enumerate}
\item
$f$ is a wheeled properadic graphical map.
\item
For each subgraph $\nicexy{H \ar[r]^-{\varphi} & G}$, the map $\nicexy{f(H) \ar[r] & K}$ is a subgraph.\index{wheeled properadic graphical map!characterization} \index{graphical map!characterization}
\end{enumerate}
\end{theorem}

\begin{proof}
The direction $(2) \Longrightarrow (1)$ holds because the identity map on $G$ is a subgraph.

For $(1) \Longrightarrow (2)$, suppose $|\edgei(G)| = n$ and $|\edgei(H)| = m$.  Since $H$ is a subgraph of $G$, we have $m \leq n$.  We prove that $f(H) \to K$ is a subgraph by downward induction on $m$.
\begin{enumerate}
\item
If $m=n$, then $H=G$, so $f(H)=f(G)$ is a subgraph of $K$ by assumption.
\item
If $m=n-1$, then
\[
G=J(H)
\]
for some dioperadic graph or contracted corolla $J$.  Since $f(G)$ is a subgraph of $K$, by Theorem \ref{thm:gwheelcsubgraph} there is a graph substitution decomposition
\[
\begin{split}
K &= M[f(G)]\\
&= M[f(J(H))]\\
&= M[(f_0J)(f(H))]\\
&= [M(f_0J)](f(H)).
\end{split}
\]
So by Theorem \ref{thm:gwheelcsubgraph} again $f(H)$ is a subgraph of $K$.
\item
Suppose $m<n-1$.  Since $H$ is a subgraph of $G$, by Theorem \ref{thm:gwheelcsubgraph} there is a graph substitution decomposition
\[
G = J(H).
\]
The assumption  $m<n-1$ implies that $J$ has at least two internal edges.  Exactly as in the proof of Theorem \ref{thm:gwheelcsubgraph}, there is an outer dioperadic or contracting factorization
\[
J = I(M)
\]
such that the graph substitution $M(H)$ makes sense.  We have
\[
G = [I(M)](H) = I[M(H)],
\]
so $M(H) \to G$ is an outer coface map, hence a subgraph of $G$.  Note that $M$ is ordinary with at least one internal edge, so $|\edgei(M(H))| > m$.  By induction hypothesis $f(M(H))$ is a subgraph of $K$.  By Theorem \ref{thm:gwheelcsubgraph} there is a graph substitution decomposition 
\[
\begin{split}
K &= P[f(M(H))]\\
&= [P(f_0M)](f(H)).
\end{split}
\]
So by Theorem \ref{thm:gwheelcsubgraph} again $f(H)$ is a subgraph of $K$.
\end{enumerate}
\end{proof}

\begin{corollary}
\label{cor1:gwheelcgraphicalmap}
Suppose $\nicexy{G \ar[r]^-{f} & K}$ is a wheeled properadic graphical map.  Then for each vertex $v$ in $G$, $f_1(v)$ is a subgraph of $K$.
\end{corollary}

\begin{proof}
Use Theorem \ref{thm:wgraphicalmapchar} and the fact that the corolla inclusion $C_v \to G$ is a subgraph (Corollary \ref{cor:gwheelcsubgraph}).
\end{proof}

\subsection{Graphical Category}

Here we define the graphical category for connected graphs.

\begin{definition}
\label{def:wgraphicalcat}
The \textbf{wheeled properadic graphical category}, or simply the \textbf{graphical category}, \index{wheeled properadic graphical category} \index{graphical category!wheeled properadic} is the category $\Gammaw$\label{note:gammaw} with
\begin{itemize}
\item
objects the graphical wheeled properads $\Gammaw(G)$ for $G \in \gwheelc$, and
\item
morphisms $\Gammaw(G) \to \Gammaw(H) \in \wproperad$ the wheeled properadic graphical maps.
\end{itemize}
Denote by
\[
\nicexy{\Gammaw \ar[r]^-{\iota} & \wproperad}
\]
the non-full subcategory inclusion.
\end{definition}

\begin{remark}
The graphical category $\Gammaw$ is small because there is only a set of $1$-colored graphs.
\end{remark}

The properadic graphical category $\varGamma$ (Definition \ref{def:graphicalcat}) is related to the wheeled properadic graphical category $\Gammaw$ as follows.

\begin{theorem}
\label{gammanonfull}
There is a non-full embedding $\nicexy{\varGamma \ar[r]^-{\iota} & \Gammaw}$.
\end{theorem}

\begin{proof}
The functor $\iota$ sends the object $\varGamma(G)$ for $G \in \gupc \subseteq \gwheelc$ to $\Gammaw(G)$.  By Lemmas \ref{lem:mapbtwgprop} and \ref{lem:mapbtwgwprop}, we have a canonical injection
\begin{equation}
\label{propwpropinject}
\nicexy{
\properad(\varGamma(G),\varGamma(K)) \ar@{^{(}->}[d]\\
\wproperad\left(\Gammaw(G), \Gammaw(K)\right)
}
\end{equation}
for any $G,K \in \gupc$.  Suppose $f \in \varGamma(G,K)$, i.e., $f$ is a properadic graphical map.  Then the above injection says that $f$ also defines a map of graphical wheeled properads.  For the functor $\iota$ to be defined, we need to show that $f \in \Gammaw(G,K)$, i.e., $f$ is a wheeled properadic graphical map.  Since $f$ is a properadic graphical map, by Theorem \ref{thm:gupcsubgraph} there is a graph substitution decomposition
\[
K = H[f(G)]
\]
in $\gupc$, and hence also in $\gwheelc$.  By Theorem \ref{thm:gwheelcsubgraph} $f(G) \to K$ is a subgraph in $\wproperad$, so the map induced by $f$ is a wheeled properadic graphical map.  This shows that there is an embedding $\iota$.

To see that $\iota$ is \emph{not} full, consider the properad map $\nicexy{H \ar[r]^-{\psi} & G}$ in Example \ref{ex1:noninjection}.  The map $\psi$ does not belong to the subcategory $\varGamma$ because $\psi(H)$ is not a subgraph of $G$ in $\properad$.  On the other hand,
\[
\psi(H) \to G \in \wproperad
\]
is an outer contracting coface map corresponding to the disconnectable edge $e$.  So it is a subgraph, and $\psi \in \Gammaw$.
\end{proof}

\begin{remark}
The injection \eqref{propwpropinject} is strict in general.  For example, since $\varGamma(\uparrow)\emptyprofh = \varnothing$, we have
\[
\properad(\bullet, \uparrow) = \varnothing.
\]
On the other hand, there is a composition
\[
\nicexy{
\bullet \ar[r] & \wheel \ar[r]^-{\cong} & \uparrow
}\]
in $\wproperad$, in which $\bullet \to \wheel$ is the exceptional inner coface map (Definition \ref{def:exceptionalinner}), and the isomorphism is the one in Example \ref{ex:unitgwproperad}.  So
\[
\wproperad(\bullet, \uparrow) \not= \varnothing,
\]
which means \eqref{propwpropinject} is strict in this case.
\end{remark}

\subsection{Factorization of Graphical Maps}

To establish the wheeled analog of the codegeneracies-cofaces factorization, we need the following wheeled analog of Lemma \ref{innercofacefactor}.

\begin{lemma}
\label{winnercofacefactor}
Suppose $K=G(H_w)$ in $\gwheelc$ with $H_w \not\in \{\uparrow,\wheel\}$.  Then the map $\nicexy{G \ar[r]^-{\varphi} & K}$ determined by
\[
\varphi_1(u) = 
\begin{cases}
C_u & \text{ if $u\not= w$},\\
H_w & \text{ if $u=w$},
\end{cases}
\]
has a decomposition into inner coface maps and isomorphisms induced by changes of vertex listings.
\end{lemma}

\begin{proof}
Write $H$ for $H_w$.  We prove the assertion by induction on $n=|\edgei(H)|$.  If $n=0$, then $H$ is a permuted corolla, and $K$ is obtained from $G$ by changing the listing at the vertex $w$.  It induces a canonical isomorphism $\nicexy{
G \ar[r]^-{\cong} & K}$.

Suppse $n>0$.  Pick $e \in \edgei(H)$, which may connect two distinct vertices or maybe a loop.  Define $J$ as the graph obtained from $H$ by shrinking away $e$ (which means deleting $e$ in case it is a loop).  By Theorems \ref{thm:innerdiopfact} and \ref{thm:innercontfact}, there is a graph substitution factorization
\[
H = J(I)
\]
for some dioperadic graph or contracted corolla $I$.  We have
\[
K = G[J(I)] = [G(J)](I),
\]
so $\varphi$ factors as
\[
\nicexy{
G \ar[r]^-{\psi^1} & G(J) \ar[r]^-{\psi^2} & K.
}\]
Since $|\edgei(J)| = n-1$, by induction hypothesis the map $\psi^1$ has a decomposition into inner coface maps and isomorphisms induced by changes of vertex listings.  The map $\psi^2$ is by definition an inner dioperadic or contracting coface map.  This proves the lemma.
\end{proof}

The following factorization result is the wheeled analog of Theorem \ref{thm:gupcepimono}.  It uses the factorization in Lemma \ref{lem1:gwheelcimage}

\begin{theorem}
\label{thm:gwheelcepimono}
Suppose $\nicexy{G \ar[r]^-{f} & K}$ is a wheeled properadic graphical map.  Then there is a factorization
\[
\nicexy@C+10pt{
G \ar[rr]^-{f}  \ar[d]_{\sigma} \ar[drr]^-{g} && K\\
G_1 \ar[r]^-{\cong}_-{i} & G_2 \ar[r]_-{\delta} & f(G) \ar[u]_{h}
}\]
in which:
\begin{itemize}
\item
$\sigma$ is a composition of codegeneracy maps, 
\item
$i$ is an isomorphism,
\item
$\delta$ is a composition of inner coface maps, and 
\item
$h$ is a composition of outer coface maps.\index{graphical map!factorization}
\end{itemize}
\end{theorem}

\begin{proof}
If $\varphi$ is the exceptional inner coface map $\bullet \to \wheel$ (Definition \ref{def:exceptionalinner}), then there is nothing to prove.  So we assume $\varphi$ is not the exceptional inner coface map.

By assumption $\nicexy{f(G) \ar[r]^-{h} & K}$ is a subgraph, i.e., a composition of outer coface maps.  So it suffices to show that $g$ decomposes as $\delta \sigma$ up to isomorphism as stated.  If $G \in \{\uparrow,\wheel\}$, then $g$ is an isomorphism.  So we may assume that $G$ is ordinary.

We now proceed as in the proof of Theorem \ref{thm:gupcepimono}, defining
\[
\nicexy{
G \ar[r]^-{\sigma} 
& G_1 \ar[r]^-{\cong}
& G_2 \ar[r]^-{\delta}
& f(G)
}\]
in exactly the same way.  In the factorization \eqref{deltafactor} of the current $\delta$, each map $H_i \to H_{i+1}$ is a composition of inner coface maps by Lemma \ref{winnercofacefactor}.  Therefore, $\delta$ is a composition of inner coface maps.
\end{proof}

\newcommand{\Gammawimage}{\Gammaw^\blacksquare}
\newcommand{\Gammawout}{\Gammaw^{\operatorname{out}}}
\newcommand{\Gammawin}{\Gammaw^{\operatorname{in}}}

We now address the question of uniqueness of decomposition of maps in $\Gammaw$. 
We begin by defining several wide subcategories of $\Gammaw$. Let
\begin{enumerate}
	\item $\Gammawimage$\label{note:gammawimage} denote the subcategory of $\Gammaw$ consisting of all maps $f \colon G \to K$ such that $f(G) = K$, 
	\item $\Gammawout$\label{note:gammawout} denote the subcategory of $\Gammaw$ generated by outer coface maps and isomorphisms, 
	\item $\Gammawin$\label{note:gammawin} denote the subcategory of $\Gammaw$ generated by inner coface maps and isomorphisms, and
	\item $\Gammaw^\minus$\label{note:gammawminus} denote the subcategory of $\Gammaw$ generated by codegeneracies and isomorphisms.
\end{enumerate}

\begin{lemma}
\label{lem:iii.3}
	Suppose that $G \overset{f}{\to} H \overset{g}\to K$ is in $\Gammaw$. 
	If $f\in \Gammawimage$, then $(gf)(G) = g(H)$.
\end{lemma}

\begin{lemma}
\label{lem:iii.1}
The categories $\Gammaw^\minus$ and $\Gammawin$ are subcategories of $\Gammawimage$.
\end{lemma}
\begin{proof}
	It is enough to note that if $f: G\to K$ is a codegeneracy, an inner coface map, or an isomorphism, then $f(G) = K$.
\end{proof}

\begin{lemma}
\label{lem:iii.5}
Suppose that $f: G \to K$ is a map in $\Gammawimage$.
If $f$ has two decompositions $f=\delta i \sigma = \delta' i' \sigma'$,
where $\sigma, \sigma'$ are compositions of codegeneracies, $i,i'$ are isomorphisms, and $\delta, \delta'$ are compositions of inner cofaces, then there is an isomorphism $i'': G_1 \to G_1'$ making the diagram 
\[
	\begin{tikzcd}
		\, & G_1 \rar{i}  \arrow{dd}{i''} & G_2 \arrow{dr}{\delta} \\
		G \arrow{ur}{\sigma} \arrow{dr}[swap]{\sigma'} & & & K \\
		\, & G_1' \rar[swap]{i'} & G_2' \arrow{ur}[swap]{\delta'}
	\end{tikzcd}
\]
commute.
\end{lemma}
\begin{proof}
	Inner coface maps are injective on edges and 
	codegeneracy maps are surjective on edges, so by uniqueness of epi-monic factorizations for maps in $\Set$, there is a bijection $i_0''$ making the diagram
\[
	\begin{tikzcd}[background color=white]
		\, & \edge(G_1) \rar{i_0} & \edge(G_2) \arrow{dr}{\delta_0} \\
		\edge(G) \arrow{ur}{\sigma_0} \arrow{dr}[swap]{\sigma_0'} \arrow{rrr}{f_0} &  &   & \edge(K) \\
		\, & \edge(G_1') \rar[swap]{i_0'}  \arrow[crossing over, leftarrow]{uu}[near start]{i_0''} & \edge(G_2') \arrow{ur}[swap]{\delta_0'}
	\end{tikzcd}
\]	
	commute.
	Injectivity of $\delta_0 i_0$ and $\delta_0' i_0'$ imply that $\sigma$ and $\sigma'$ identify the same set of edges, so $i_0''$ gives the isomorphism $i'': G_1 \to G_1'$.
\end{proof}

The following lemma is a special case of Proposition \ref{prop:iii.4} below.

\begin{lemma}
\label{commuting_cofaces}
	Suppose that 
	\[ \begin{tikzcd}
		G \rar{g} \arrow{dr}[swap]{f} & H \dar{h} \\
		& K
	\end{tikzcd} \]
	commutes, with $f,h\in \Gammawout$ and $g\in \Gammawin$. Then $g$ is an isomorphism.
\end{lemma}
\begin{proof}
Outer coface maps and isomorphisms send each vertex to a corolla. If $g$ is not an isomorphism, then $hg$ must send at least one vertex of $G$ to a non-corolla, contradicting $hg=f \in \Gammawout$.

\end{proof}

\begin{lemma}
\label{lem:coface_codegen}
	Suppose that $f$ and $g$ are composable, and that $f$ is a composition of coface maps and isomorphisms. If $g$ is a nontrivial composition of codegeneracy maps, then $fg$ cannot be written as a composition of coface maps and isomorphisms.
\end{lemma}
\begin{proof}
	Suppose that $fg$ is a composition of coface maps and isomorphisms and $g$ is a nontrivial composition of codegeneracy maps. 
	Since no coface map has target $\bullet$, if one of the coface maps in the decomposition of $fg$ is the exceptional inner coface map $\bullet \to \wheel$ from Definition \ref{def:exceptionalinner}, then the source of $fg$ is $\bullet$. On the other hand, no codegeneracy map has source $\bullet$, so all coface maps in the decomposition of $fg$ must be non-exceptional.

	All non-exceptional coface maps take a vertex $v$ to a subgraph with at least one vertex, hence so does $fg$. 
	Since $g$ is a nontrivial composition of codegeneracy maps, then there is a vertex $v$ 
	so that 
	$fg(v) = \uparrow_e$,
	which has zero vertices.

\end{proof}

\begin{proposition}
\label{prop:iii.4}
Every map in $\Gammaw$ factors uniquely (up to isomorphism) as a map in $\Gammawimage$ followed by a map in $\Gammawout$.
\end{proposition}
\begin{proof}
	The existence of such a factorization is given in \ref{lem1:gwheelcimage}. Specifically, if $f: G\to K$ is a map, then we have that $f$ is equal to $G \overset{g}\to f(G) \overset{h}\to K$ with $h$ a composition of outer coface maps.

	Suppose that
	\[ G \overset{g'}\to g'(G) \overset{h'}\to K \]
	is a factorization of $f: G\to K$, 
	where $g' \in \Gammawimage$ and $h'\in \Gammawout$.
	Apply \ref{lem1:gwheelcimage} to the maps $h': g'(G)\to K$
	\[
	\nicexy@C+10pt{
	g'(G) \ar[r]^-{a} \ar[dr]_-{h'} & h'(g'(G)) \ar[d]^-{b}\\
	& K}
	\]
	so that $h'=ba$ with $b$ a composition of outer coface maps. 
	By Lemma \ref{lem:iii.3}, $h'(g'(G)) = f(G)$, hence $b=h$.

	Factor $a$ as in Lemma \ref{lem:iii.5}, so that $a = c e$ with $c \in \Gammawin$ and $e$ a composition of codegeneracies. Since $h' = b a = (b c) e$ is a composition of coface maps, $e$ is trivial by Lemma \ref{lem:coface_codegen}. Thus $a = c$ is in $\Gammawin$. But then $h' = h a$ and $h$ are both compositions of outer coface maps, hence $a$ is an isomorphism by Lemma \ref{commuting_cofaces}.
	The diagram
	\[
	\begin{tikzcd}
	G \rar{g} \dar[swap]{g'} & f(G) \dar{h} \\
	g'(G) \rar[swap]{h'} \arrow{ur}{a}[swap]{\cong} & K
	\end{tikzcd}
	\]
	commutes, which shows uniqueness.



\end{proof}

\chapter{Infinity Wheeled Properads}
\label{ch:infwproperad}

\abstract*{We define the adjunction
$L: \gwheelcset \rightleftarrows \wproperad : N$
between wheeled properads and wheeled properadic graphical sets.  Then we define $\infty$-wheeled properads as wheeled properadic graphical sets that satisfy an inner horn extension property.  Next we give two alternative characterizations of \emph{strict} $\infty$-wheeled properads, one in terms of the wheeled properadic Segal maps, and the other in terms of the wheeled properadic nerve.  In the last section, we give an explicit description of the fundamental wheeled properad $L\sK$ of an $\infty$-wheeled properad $\sK$ in terms of homotopy classes of $1$-dimensional elements.  }

In this chapter, we discuss (strict) $\infty$-wheeled properads.

In section \ref{sec:infwproperad} we define an adjunction
\[
\nicexy@C+10pt{
\gwheelcset \ar@<2pt>[r]^-{L} & \ar@<2pt>[l]^-{N} \wproperad}
\]
between the category $\gwheelcset$ of wheeled properadic graphical sets and the category $\wproperad$ of wheeled properads.  The right adjoint $N$ is the wheeled version of the properadic nerve.  Using the wheeled properadic nerve and representables, the symmetric monoidal product in $\wproperad$ induces a symmetric monoidal closed structure on $\gwheelcset$.  The coface maps in the graphical category $\Gammaw$ induce faces and horns of representables in $\gwheelcset$.  An $\infty$-wheeled properad is a wheeled properadic graphical set that has an inner horn extension property.  Strict $\infty$-wheeled properads are the ones with unique inner horn extensions.

In section \ref{sec:strictinfwproperad} we give two characterizations of strict $\infty$-wheeled properads.  First, a wheeled properadic graphical set is a strict $\infty$-wheeled properad if and only if it is isomorphic to the wheeled properadic nerve of some wheeled properad.  Second, a strict $\infty$-wheeled properad is equivalent to a wheeled properadic graphical set that satisfies the wheeled properadic Segal condition.  These characterizations are the wheeled versions of those in Theorem \ref{properadnerve} for strict $\infty$-properads.

In section \ref{sec:fundwproperad} we characterize the fundamental wheeled properad (i.e., the image under $L$) of a \emph{reduced} $\infty$-wheeled properad in terms of homotopy classes of $1$-dimensional elements.  The reduced condition on $\sK \in \gwheelcset$ means that $\sK(C_{\emptyprofh})$ consists of a single element, where $C_{\emptyprofh}$ is a single isolated vertex. The definition of homotopy in the properadic case (Definition \ref{def:graphhomotopy}) is used here without change.  Much of the work in this section involves showing that the dioperadic composition and the contraction are well-defined on homotopy classes, and that the object with these operations is actually a wheeled properad.

\section{\texorpdfstring{$\infty$}{∞}-Wheeled Properads}
\label{sec:infwproperad}

The main purpose of this section is to define (strict) $\infty$-wheeled properads.  Along the way, we define wheeled versions of graphical sets, the nerve, representables, faces, horns, and the symmetric monoidal closed structure on graphical sets.

\subsection{Wheeled Properadic Graphical Sets and Nerve}

Here we define the wheeled analog of the category of graphical sets, the corresponding nerve functor, and representables.

\begin{definition}
\label{def:wgraphicalset}
The diagram category $\gwheelcset$\label{note:gwcset} is called the category of \textbf{wheeled properadic graphical sets}, or simply \textbf{graphical sets}.\index{wheeled properadic graphical set} \index{graphical set!wheeled properadic}
\begin{enumerate}
\item
An object $X \in \gwheelcset$ is called a \textbf{wheeled properadic graphical set}, or simply a \textbf{graphical set}.
\item
For $G \in \gwheelc$, an element in the set $X(G)$ is called a \textbf{graphex with shape $G$}.\index{graphex}  The plural form of \emph{graphex} is \emph{graphices}.
\item
A graphical set $X$ is \textbf{reduced} \index{graphical set!reduced} if the set $X(C_{\emptyprofh})$ is a singleton, where $C_{\emptyprofh}$ is the single isolated vertex.
\end{enumerate}
\end{definition}

\begin{definition}
For a graphical set $X \in \gwheelcset$, the set $X(\uparrow)$ of \textbf{colors} of $X$, the \textbf{$1$-dimensional elements} in $X$, their \textbf{profiles}, the associated $\Sigma_{\sS(X(\uparrow))}$-bimodule, and \textbf{colored units} of $X$ are defined as in section \ref{rk:gsetproperad}.
\end{definition}

\begin{definition}
\label{wgraphicalnerve}
The \textbf{wheeled properadic nerve} \index{wheeled properadic nerve} \index{nerve!wheeled properadic} is the functor
\[\label{note:wnerve}
\nicexy{
\wproperad \ar[r]^-{N} & \gwheelcset
}\]
defined by
\[
(N\sP)(G) = \wproperad(G,\sP)
\]
for $\sP \in \wproperad$ and $G \in \Gammaw$.    In the context of the factorization of a wheeled graphical map $G \to K \in \Gammaw$ (Theorem \ref{thm:gwheelcepimono}), the map
\[
\nicexy{
(N\sP)(K) = \wproperad(K,\sP) \ar[r] &
\wproperad(G,\sP) = (N\sP)(G)}
\]
is a composition of maps of the following form:
\begin{itemize}
\item
deletion of an entry (for an outer dioperadic coface map);
\item
dioperad composition of $\sP$ (for an inner dioperadic coface map);
\item
contraction of $\sP$ (for an inner contracting coface map);
\item
isomorphism (for an outer contracting coface map);
\item
colored units of $\sP$ (for a codegeneracy map).
\end{itemize}
\end{definition}

\begin{remark}
To see that $N\sP$ is indeed in $\gwheelcset$, recall that maps in $\Gammaw$ are compositions of coface maps, codegeneracy maps, and isomorphisms.  Coface and codegeneracy maps correspond to graph substitutions.  The wheeled properad structure maps of a wheeled properad $\sP$ are associative and unital with respect to graph substitutions.
\end{remark}

\begin{lemma}
\label{wnerveadjoint}
The wheeled properadic nerve admits a left adjoint
\[
\nicexy{
\gwheelcset \ar[r]^-{L} & \wproperad}
\]
such that the diagram
\[
\nicexy{
\Gammaw \ar[r]^-{\iota} \ar[d]_{\mathrm{Yoneda}} & \wproperad\\
\gwheelcset \ar[ur]_{L} &}
\]
is commutative up to natural isomorphism.\index{wheeled properadic nerve!left adjoint}
\end{lemma}

\begin{proof}
It is the same as the proof of Lemma \ref{graphicalnerveadjoint}, with $\varGamma$ and $\properad$ replaced by $\Gammaw$ and $\wproperad$.
\end{proof}

\begin{definition}
\label{def:fundwproperad}
For $\sK \in \gwheelcset$, the image $L\sK$ is called the \textbf{fundamental wheeled properad} of $\sK$.\index{fundamental wheeled properad}
\end{definition}

The following observation says that a graphex in $(N\sP)(G)$ is really a $\sP$-decoration of $G$, which consists of a coloring of the edges in $G$ by the colors of $\sP$ and a decoration of each vertex in $G$ by an element in $\sP$ with the corresponding profiles.

\begin{lemma}
\label{wnpgelement}
Suppose $\sP$ is a $\fC$-colored wheeled properad, and $G \in \Gammaw$.  Then an element in $(N\sP)(G)$ consists of:
\begin{itemize}
\item
a function $\nicexy{\edge(G) \ar[r]^-{\varphi_0} & \fC}$, and
\item
a function $\varphi_1$ that assigns to each $v \in \vertex(G)$ an element $\varphi_1(v) \in \sP\binom{\varphi_0\out (v)}{\varphi_0 \inp(v)}$.
\end{itemize}
\end{lemma}

\begin{proof}
This is a special case of Lemma \ref{lem:mapfromgwprop}.
\end{proof}

Next we define the wheeled analog of the representable graphical set.

\begin{definition}
Suppose $G \in \Gammaw$.  The \textbf{representable wheeled graphical set} \index{representable wheeled graphical set}  $\Gammaw[G] \in \gwheelcset$\label{note:wrepresentable} is defined, for $H \in \Gammaw$, by
\[
\Gammaw[G](H) = \Gammaw(H,G),
\]
i.e., the set of wheeled properadic graphical maps $H \to G$ (Definition \ref{def:wgraphicalmap}).
\end{definition}

\begin{remark}
By Yoneda's Lemma a map $\Gammaw[G] \to X$ of wheeled properadic graphical sets is equivalent to a graphex in the set $X(G)$.
\end{remark}

\subsection{Symmetric Monoidal Closed Structure}

Here we observe that the symmetric monoidal product of wheeled properads and the wheeled properadic nerve induce a symmetric monoidal closed structure on the category of wheeled properadic graphical sets.

Note that each graphical set $X \in \gwheelcset$ can be expressed up to isomorphism as a colimit of representable graphical sets,
\[
X \cong \colim_{\Gammaw[G] \to X} \Gammaw[G],
\]
where the colimit is indexed by maps of graphical sets $\Gammaw[G] \to X$, i.e., graphices in $X$.

\begin{definition}
Suppose $X$ and $Y$ are wheeled properadic graphical sets.
\begin{enumerate}
\item
Define the wheeled properadic graphical set
\[
X \otimes Y \defn \colim_{\varGamma[G]\to X, \varGamma[G']\to Y} \left(\varGamma[G] \otimes \varGamma[G']\right),
\]
where
\[
\varGamma[G] \otimes \varGamma[G'] 
\defn
N\left(\varGamma(G) \otimes \varGamma(G')\right)
\]
with $N$ the wheeled properadic nerve and $\otimes$ the symmetric monoidal product in $\wproperad$ (Definition \ref{def:wproperadtensor}).
\item
Define the wheeled properadic graphical set $\Hom(X,Y)$ by
\[
\Hom(X,Y)(G) 
\defn 
\gwheelcset\left(X \otimes \Gammaw[G], Y\right)
\]
for $G \in \Gammaw$.
\end{enumerate}
\end{definition}

\begin{theorem}
\label{gwcsmclosed}
The category $\gwheelcset$ is symmetric monoidal closed with monoidal product $\otimes$ and internal hom $\Hom$.\index{wheeled properadic graphical set!symmetric monoidal closed structure}
\end{theorem}

\begin{proof}
It is the same as the proof of Theorem \ref{dendroidalclosed}, with $\varGamma$ and $\properad$ replaced by $\Gammaw$ and $\wproperad$.
\end{proof}

\subsection{Faces and Horns}

Here we define the wheeled analogs of the horns $\Lambda^d[G] \subset \varGamma[G]$, which we will soon use to define an $\infty$-wheeled properad.

\begin{definition}
Suppose $G \in \Gammaw$.  A \textbf{face of $G$} \index{face map!of graphical wheeled properad} is a coface map in $\Gammaw$ whose target is $G$.  An \textbf{inner/outer face of $G$} \index{inner face!of graphical wheeled properad} \index{outer face!of graphical wheeled properad} is an inner/outer coface map whose target is $G$.
\end{definition}

\begin{definition}
Suppose $X$ is a wheeled properadic graphical set.  A \textbf{graphical subset} \index{graphical subset} of $X$ is a wheeled properadic graphical set $W$ that is equipped with a map $W \to X$ of wheeled properadic graphical sets such that each component map
\[
W(G) \to X(G)
\]
for $G \in \Gammaw$ is a subset inclusion.
\end{definition}

\begin{definition}
Suppose $G \in \Gammaw$, and $\nicexy{K \ar[r]^{d} & G}$ is a face of $G$.
\begin{enumerate}
\item
The \textbf{$d$-face} of $\Gammaw[G] \in \gwheelcset$ is the graphical subset $\Gammaw^d[G]$ defined by
\[\label{note:wface}
\begin{split}
& \Gammaw^d[G](J) \\
&= 
\left\{\text{composition of }
\nicexy{J \ar[r]^-{f} & K \ar[r]^-{d} & G}
\text{ with $f \in \Gammaw[K](J)$}\right\}
\end{split}
\]
for $J \in \Gammaw$.
\item
The \textbf{$d$-horn} \index{horn!of wheeled properadic graphical set} of $\Gammaw[G]$ is the graphical subset $\Lambdaw^d[G]$ defined by
\[\label{note:whorn}
\Lambdaw^d[G](J) 
= 
\bigcup_{\text{faces $d'\not=d$}} \Gammaw^{d'}[G] (J),
\]
where the union is indexed by all the faces of $G$ \emph{except} $d$.  Write
\[
\nicexy{\Lambdaw^d[G] \ar[r]^-{i} & \Gammaw[G]}
\]
for the graphical subset inclusion.
\item
A \textbf{horn}  of $\Gammaw[G]$ is a $d$-horn for some coface map $d$.  An \textbf{inner horn} \index{inner horn!of wheeled properadic graphical set} is a $d$-horn in which $d$ is an inner coface map.
\item
Given  $X \in \gwheelcset$, a \textbf{horn of $X$} is a map
\[
\nicexy{\Lambdaw^d[G] \ar[r] & X}
\]
of graphical sets.  It is an \textbf{inner horn of $X$} if $\Lambdaw^d[G]$ is an inner horn.
\end{enumerate}
\end{definition}

\begin{remark}
\begin{enumerate}
\item
If $\nicexy{K \ar[r]^{d} & G}$ is a face of $G$, then the $d$-face $\Gammaw^d[G]$ is the image of the induced map
\[
\nicexy{\Gammaw[K] \ar[r]^{d \circ (-)} & \Gammaw[G]}
\]
of graphical sets. However, unlike the properadic analog, this induced map may \emph{not} be injective.  Indeed, when $d$ is an outer contracting coface map corresponding to a disconnectable edge $e$, it sends two legs $e_{\pm 1} \in \edge(K)$ to $e \in \edge(G)$.  When applied to $\uparrow$ and using the identification
\[
\Gammaw[K](\uparrow) = \edge(K),
\]
the induced map of $d$ is then the map $\nicexy{\edge(K) \ar[r]^-{d_0} & \edge(G)}$, which is surjective but not injective.
\item
In the definition of the $d$-horn $\Lambdaw^d[G]$, we used the same convention as before about ignoring listings (Convention \ref{graphconvention}).  In other words, when we say $d$ is excluded, we mean $d$ and every coface map $K' \to G$ obtainable from $d$ by changing the listing are all excluded.  
\item
If there is an \emph{inner} horn $\Lambdaw^d[G]$ for some inner coface map $\nicexy{K \ar[r]^-{d} & G}$, then $G$ is ordinary and has at least one internal edge because there is an inner dioperadic or contracting factorization (Definition \ref{def:innerdiopfact} and \ref{def:innercontfact}) of $G$,
\[
G = K(H_w),
\]
in which the distinguished subgraph $H_w$ is a dioperadic graph or a contracted corolla.
\end{enumerate}
\end{remark}

The following observation gives a more explicit description of a horn of a graphical set.

\begin{lemma}
\label{whorndescription}
Suppose $X \in \gwheelcset$, and $\nicexy{K \ar[r]^{d} & G}$ is a face of $G \in \Gammaw$.  Then a horn
\[
\nicexy{\Lambdaw^d[G] \ar[r]^-{f} & X}
\]
of $X$ is equivalent to a collection of maps
\[
\left\{\nicexy{\Gammaw[H] \ar[r]^-{f_H} & X} \colon 
\nicexy{H \ar[r]^{d'} & G} \text{ face $\not= d$}
\right\}
\]
such that, if
\[
\nicexy{
J \ar[r]^-{a^1} \ar[d]_-{a^2} 
& H^1 \ar[d]^{d^1}\\
H^2 \ar[r]^-{d^2} & G
}\]
is a commutative diagram of coface maps with each $d^i \not= d$, then the diagram
\[
\nicexy{
\Gammaw[J] \ar[r]^-{a^1} \ar[d]_-{a^2} 
& \Gammaw[H^1] \ar[d]^{f_{H_1}}\\
\Gammaw[H^2] \ar[r]^-{f_{H_2}} & X
}\]
is also commutative.\index{inner horn!characterization}
\end{lemma}

\begin{proof}
It is the same as the proof of Lemma \ref{horndescription}, with $\varGamma$ and $\properad$ replaced by $\Gammaw$ and $\wproperad$.
\end{proof}

\begin{remark}
\label{rk:whornofx}
The collection of maps $\{f_H\}$ in Lemma \ref{whorndescription} is equivalent to a collection of graphices $\left\{f_H \in X(H)\right\}$ such that, if the coface square is commutative, then
\[
(a^1)^*\left(f_{H_1}\right) 
= 
(a^2)^*\left(f_{H_2}\right) \in X(J).
\]
In other words, a horn of $X$ is really a collection of graphices in $X$, one for each face not equal to the given one, that agree on common faces.
\end{remark}

In the following examples, we describe explicitly an inner horn for the two types of generating graphs of $\gwheelc$.

\begin{example}
\label{ex:diophorn}
Suppose $X \in \gwheelcset$, $D$ is a dioperadic graph (Example \ref{ex:dioperadic}) with top (resp., bottom) vertex $v$ (resp., $u$), and
\[
\nicexy{C \ar[r]^-{d} & D}
\]
is the unique \emph{inner} dioperadic coface map corresponding to the only internal edge $e$ in $D$.  There are only two other faces of $D$:
\begin{itemize}
\item
the \emph{outer} dioperadic coface map
\[
\nicexy{C_u \ar[r]^-{d^v} & D}
\]
corresponding to the deletable vertex $v$, and
\item
the \emph{outer} dioperadic coface map
\[
\nicexy{C_v \ar[r]^-{d^u} & D}
\]
corresponding to the deletable vertex $u$.
\end{itemize}
The only common face of the corollas $C_u$ and $C_v$ is the internal edge $e$ in $D$.

Therefore, an inner horn
\[
\nicexy{\Lambdaw^d[D] \ar[r]^-{f} & X}
\]
is equivalent to a pair of $1$-dimensional elements
\[
(f_u, f_v) \in X(C_u) \times X(C_v),
\]
such that if
\[
\nicexy@C+10pt{
\uparrow \ar[d]_{d^e_u} \ar[r]^-{d^e_v} & C_v\\
C_u &
}\]
are the outer dioperadic coface maps corresponding to the internal edge $e$ in $D$, then
\[
\left(d^e_u\right)^*(f_u) 
= 
\left(d^e_v\right)^*(f_v) \in X(\uparrow).
\]
In other words, such an inner horn of $X$ is a pair of $1$-dimensional elements, one for each vertex in $D$, whose profiles match along the internal edge $e$ in $D$.  We may, therefore, visualize such an inner horn of $X$ as the following $X(\uparrow)$-colored dioperadic graph
\begin{center}
\begin{tikzpicture}
\matrix[row sep=1cm,column sep=1.5cm] {
\node [fatplain,label=above:$...$] (p1) {$f_v$}; \\
\node [fatplain,label=below:$...$] (p2) {$f_u$}; \\
};
\foreach \x in {1,2}
{
\draw [outputleg] (p\x) to +(-.5cm,.4cm);
\draw [outputleg] (p\x) to +(.5cm,.4cm);
\draw [inputleg] (p\x) to +(-.5cm,-.4cm);
\draw [inputleg] (p\x) to +(.5cm,-.4cm);
}
\draw [arrow] (p2) to (p1);
\end{tikzpicture}
\end{center}
decorated by $1$-dimensional elements in $X$.
\end{example}

\begin{example}
\label{ex:contracthorn}
Suppose $X \in \gwheelcset$, and $\xiij C_{(m;n)}$ is a contracted corolla (Example \ref{ex:contractedcor}) with unique vertex $v$ and loop $e$.  The contracted corolla has two faces.
\begin{enumerate}
\item
There is an \emph{inner} contracting coface
\[
\nicexy{C_{(m-1;n-1)} \ar[r]^-{d} & \xiij C_{(m;n)}}
\]
corresponding to deleting the loop $e$ at $v$.
\item
There is an \emph{outer} contracting coface
\[
\nicexy{C_{(m;n)} \ar[r]^-{d_{\out}} & \xiij C_{(m;n)}}
\]
corresponding to disconnecting the disconnectable edge $e$.
\end{enumerate}
Therefore, an inner horn
\[
\nicexy{\Lambdaw^d\left[\xiij C_{(m;n)}\right] \ar[r]^-{f} & X}
\]
is equivalent to a $1$-dimensional element
\[
x \in X(C_{(m;n)}) = X(C_v)
\]
that satisfies the following condition: If
\[
\nicexy@C+12pt{
\uparrow \ar@<.5ex>[r]^-{d^e_{\out}} \ar@<-.5ex>[r]_-{d^e_{\inp}} & C_v 
}\]
are the outer dioperadic coface maps corresponding to the two legs (one input and one output) of $C_v$ named $e$, then
\[
(d^e_{\out})^*(x) 
= (d^e_{\inp})^*(x) \in X(\uparrow).
\]
We may, therefore, visualize such an inner horn of $X$ as the following $X(\uparrow)$-colored contracted corolla
\begin{center}
\begin{tikzpicture}
\matrix[row sep=.5cm,column sep=.5cm] {
\node [plain] (v) {$x$};\\
};
\draw [outputleg] (v) to +(-.4cm,.3cm);
\draw [outputleg] (v) to +(.4cm,.3cm);
\draw [inputleg] (v) to +(-.4cm,-.3cm);
\draw [inputleg] (v) to +(.4cm,-.3cm);
\draw [arrow, looseness=25, in=-60, out=60, loop] (v) to node[near start]{\footnotesize{$i$}} node[near end]{\footnotesize{$j$}} ();
\end{tikzpicture}
\end{center}
decorated by a $1$-dimensional element in $X$.
\end{example}

\subsection{\texorpdfstring{$\infty$}{∞}-Wheeled Properads}

We now define wheeled versions of (strict) $\infty$-properads.

\begin{definition}
\label{def:infinitywproperad}
Suppose $X \in \gwheelcset$.
\begin{enumerate}
\item
We call $X$ an \textbf{$\infty$-wheeled properad} \index{infinity wheeled properad} if for each inner horn $f$ of $X$,
\begin{equation}
\label{winnerhornfiller}
\nicearrow
\xymatrix@C+12pt{
\Lambdaw^d[G] \ar[d]_{i} \ar[r]^-{f} & X\\
\Gammaw[G] \ar@{.>}[ur]
}
\end{equation}
a dotted arrow, called an \textbf{inner horn filler}, \index{inner horn filler} exists and makes the triangle commutative.
\item
A \textbf{strict $\infty$-wheeled properad} \index{strict infinity wheeled properad} is an $\infty$-wheeled properad for which each inner horn filler in \eqref{winnerhornfiller} is unique.
\end{enumerate}
\end{definition}

\begin{remark}
For an $\infty$-wheeled properad, we only ask that \emph{inner} horns have fillers.  Also, uniqueness of an inner horn filler is not required, unless we are dealing with a \emph{strict} $\infty$-wheeled properad.
\end{remark}

\begin{example}
\label{ex:diophornfiller}
Suppose $D$ is a dioperadic graph.  Recall from Example \ref{ex:diophorn} that an inner horn
\[
\nicexy{\Lambdaw^d[D] \ar[r]^-{f} & X}
\]
is exactly a pair of $1$-dimensional elements $(f_u,f_v)$, one for each vertex in $D$, with matching profiles along the unique internal edge $e$ in $D$.  Then an inner horn filler for $f$ is a graphex $y \in X(D)$ such that
\[
(d^v)^*(y) = f_u \andspace 
(d^u)^*(y) = f_v.
\]
In other words, it has $f_u$ and $f_v$ as its outer dioperadic faces.
\end{example}

\begin{example}
\label{ex:contracthornfiller}
Suppose $\xiij C_{(m;n)}$ is a contracted corolla as in Example 
\ref{ex:contracthorn}.  Recall that an inner horn 
\[
\nicexy{\Lambdaw^d\left[\xiij C_{(m;n)}\right] \ar[r]^-{f} & X}
\]
is a $1$-dimensional element $x \in X(C_{(m;n)})$ whose output profile corresponding to the internal edge $e$ in $\xiij C_{(m;n)}$ matches with its input profile corresponding to $e$.  Then an inner horn filler for $f$ is a graphex $y \in X(\xiij C_{(m;n)})$ such that
\[
d_{\out}^*(y) = x.
\]
In other words, it has $x$ as its outer contracting face.
\end{example}

\section{Characterization of Strict \texorpdfstring{$\infty$}{∞}-Wheeled Properads}
\label{sec:strictinfwproperad}

In this section we provide two alternative descriptions of strict $\infty$-wheeled properads.  This is the wheeled version of Theorem \ref{properadnerve}.  We first define wheeled analogs of the corolla ribbon, the Segal maps, and the Segal condition.  Then we prove the main Theorem \ref{wproperadnerve} using the same strategy as in the proof of Theorem \ref{properadnerve}.  The main theorem says that strict $\infty$-wheeled properads are precisely the images of the nerve up to isomorphism.  Equivalently, strict $\infty$-wheeled properads are precisely the wheeled properadic graphical sets that satisfy the wheeled properadic Segal condition.

\subsection{Outer Coface Maps from Corollas}

The definition of the wheeled properadic Segal map involves the map in the following observation, which is the wheeled analog of Lemma \ref{cvtog}.

\begin{lemma}
\label{wcvtog}
Suppose $G \in \gwheelc$ is ordinary with $n$ internal edges, and $v \in \vertex(G)$.  Then there exist outer coface maps
\[
\nicexy{
C_v \ar[r]^-{\xi_v^n} 
& G_{n-1} \ar[r]^-{\xi_v^{n-1}} 
& \cdots \ar[r]^-{\xi^2_v}
& G_1 \ar[r]^-{\xi^1_v}
& G,
}\]
whose composition is the corolla inclusion $\nicexy{C_v \ar[r]^-{\xi} & G}$.
\end{lemma}

\begin{proof}
This follows from the first part of Corollary \ref{cor:gwheelcsubgraph}, since each subgraph decomposes into outer coface maps, each increasing the number of internal edges by one.  The only exception is the outer dioperadic coface map $\uparrow ~\to C$ to a corolla, but the corolla inclusion starts at $C_v$.  So a decomposition of $\xi$ into outer coface maps does not involve the map $\uparrow ~\to C$.
\end{proof}

The following observation is the wheeled analog of Lemma \ref{uparrowtog}.  It ensures that the wheeled analog of the Segal map is well defined.  Recall that the graphical wheeled properad $\Gammaw(\uparrow)$ was discussed in Example \ref{ex:unitgwproperad}.

\begin{lemma}
\label{wuparrowtog}
Suppose $\nicexy{u \ar[r]^-{e} & v}$ is an ordinary edge in $G \in \gwheelc$, where $u=v$ is allowed.  Then the square
\[
\nicearrow
\xymatrix{
\uparrow \ar[r]^-{\out_e} \ar[d]_{\inp_e} & C_u \ar[d]^{\xi_u}\\
C_v \ar[r]^-{\xi_v} & G
}\]
in $\Gammaw$ is commutative, where $\inp_e$ (resp., $\out_e$) is the outer dioperadic coface map that identifies $e$ as an input (resp., output) leg of $C_v$ (resp., $C_u$).
\end{lemma}

\begin{proof}
Both compositions send the unique element in $\Gammaw(\uparrow)(e;e)$ (resp., $\Gammaw(\uparrow)\emptyprofh$) to the $e$-colored exceptional edge $\uparrow_e$ (resp., exceptional wheel $\wheel_e$).
\end{proof}

\subsection{Wheeled Properadic Segal Maps}

Here we define the wheeled analog of the properadic Segal maps.

\begin{definition}
Suppose $\sK \in \gwheelcset$, and $G \in \gwheelc$ is ordinary with at least one internal edge.
\begin{enumerate}
\item
Define the set
\[\label{note:wcorribbon}
\sK(G)_1 \defn \left(\prod_{v \in \vertex(G)} \sK(C_v)\right)_{\sK(\uparrow)}
\]
as the limit of the diagram consisting of the maps
\begin{equation}
\label{wkcuv}
\nicearrow
\xymatrix@C+12pt{
& \sK(C_u) \ar[d]^{\sK(\out_e)}\\
\sK(C_v) \ar[r]^-{\sK(\inp_e)} & \sK(\uparrow)
}
\end{equation}
as $\nicearrow\xymatrix{u \ar[r]^-{e} & v}$ runs through $\edgei(G)$. Call $\sK(G)_1$ the \textbf{corolla ribbon} \index{corolla ribbon!for wheeled properadic graphical set} of $\sK(G)$.
\item
Define the  \textbf{wheeled properadic Segal map} \index{wheeled properadic Segal map} \index{Segal map!wheeled properadic}
\begin{equation}
\label{wsegalmap}
\nicearrow
\xymatrix@C+10pt{
\sK(G) \ar[r]^-{\chi_G} & \sK(G)_1
}
\end{equation}
as the unique map induced by the maps
\[
\nicearrow\xymatrix{\sK(G) \ar[r]^-{\sK(\xi_v)} & \sK(C_v)}
\]
as $v$ runs through $\vertex(G)$.  The wheeled properadic Segal map is well-defined by Lemma \ref{wuparrowtog}.
\item
We say $\sK$ satisfies the \textbf{wheeled properadic Segal condition} \index{wheeled properadic Segal condition} \index{Segal condition!wheeled properadic} if the wheeled properadic Segal map $\chi_G$ is a bijection for every ordinary $G \in \gwheelc$ with at least one internal edge. 
\end{enumerate}
\end{definition}

We will often drop the phrase \emph{wheeled properadic} if there is no danger of confusion.

\begin{remark}
There is no need to consider the wheeled properadic Segal map for ordinary $G \in \gwheelc$ with $\edgei(G) = \varnothing$.  Indeed, such a $G$ must be a permuted corolla $\sigma C_v \tau$, so the corolla ribbon is $\sK(C_v)$.  The Segal map is the isomorphism
\[
\nicearrow
\xymatrix@C+12pt{
\sK(\sigma C_v \tau) \ar[r]^-{\sK(\xi_v)}_-{\cong} & \sK(C_v)
}\]
induced by input/output relabeling.  Therefore, whenever we mention the wheeled properadic Segal map, we automatically assume that $G$ has at least one internal edge.
\end{remark}

\begin{remark}
\label{rk:wribbongeometric}
For an ordinary $G \in \gwheelc$, an element $\theta \in \sK(G)_1$ is equivalent to the data:
\begin{itemize}
\item
a function $\nicexy{\edge(G) \ar[r]^-{\theta_0} & \sK(\uparrow)}$, and
\item
a function $\theta_1$ that assigns to each vertex $v \in \vertex(G)$ a $1$-dimensional element $\theta_1(v) \in \sK(C_v)$ with profiles corresponding to those of $v$ under $\theta_0$. 
\end{itemize}
In other words, $\theta$ is a $\sK(\uparrow)$-colored decoration of $G$, in which vertices are decorated by $1$-dimensional elements in $\sK$.
\end{remark}

\begin{example}
Consider the following connected graph $G$.
\begin{center}
\begin{tikzpicture}
\matrix[row sep=1cm, column sep=1cm]{
\node [plain] (v) {$v$};\\
\node [plain] (u) {$u$};\\
};
\draw [arrow, bend left=45] (u) to node{\footnotesize{$f$}} (v);
\draw [arrow, bend left=45] (v) to node{\footnotesize{$g$}} (u);
\draw [arrow, out=25, in=-25, loop] (v) to node{\footnotesize{$e$}} ();
\end{tikzpicture}
\end{center}
There are two vertices and three ordinary edges, one of which is a loop.  Each vertex may have input and output legs of $G$, which for simplicity are not depicted in the picture.  There is a commutative diagram
\[
\nicearrow
\xymatrix@C+12pt{
\uparrow_e \ar[d]_-{\inp} \ar@/_2pc/[d]_-{\out} 
& \uparrow_f \ar[dl]_-{\inp} \ar[dr]^-{\out} &\\
C_v  \ar[dr]_-{\xi_v}
& \uparrow_g \ar[l]_-{\out} \ar[r]^-{\inp}
& C_u \ar[dl]^-{\xi_u}\\
& G &
}\]
in $\Gammaw$ by Lemma \ref{wuparrowtog}.  For a graphical set $\sK$, the corolla ribbon $\sK(G)_1$ is the limit of the part of this diagram above $G$, after applying $\sK$.  The wheeled properadic Segal map $\chi_G$ is induced by the two maps $\xi_*$.
\end{example}

As in the properadic case, there is another description of the Segal map in terms of the Segal core, which we now define.

\begin{definition}
Suppose $G \in \Gammaw$ is ordinary with at least one internal edge.
\begin{enumerate}
\item
Define the \textbf{wheeled properadic Segal core} \index{wheeled properadic Segal core} \index{Segal core!wheeled properadic} $\Sc[G] \in \gwheelcset$\label{note:wsegalcore} as the colimit of the diagram consisting of the maps
\[
\nicearrow
\xymatrix{
\Gammaw[\uparrow] \ar[r]^-{\out_e} \ar[d]_{\inp_e} 
& \Gammaw[C_u] \\
\Gammaw[C_v] &
}\]
as $\nicexy{u \ar[r]^-{e} & v}$ runs through $\edgei(G)$.
\item
Denote by
\[\label{note:wepsilong}
\nicexy{\Sc[G] \ar[r]^-{\epsilon_G} & \Gammaw[G]}
\]
the map induced by the corolla inclusions
\[
\nicexy{C_v \ar[r]^-{\xi_v} & G}
\]
for $v \in \vertex(G)$, and call it the \textbf{wheeled properadic Segal core map}.\index{wheeled properadic Segal core map}
\end{enumerate}
\end{definition}

\begin{lemma}
\label{lem:wsegalcoredesc}
Suppose $\sK \in \gwheelcset$, and $G \in \gwheelc$ is ordinary with at least one internal edge.  Then there is a commutative diagram
\[
\nicearrow
\xymatrix{
\sK(G) \ar[dd]_{\cong} \ar[rr]^-{\chi_G} && \sK(G)_1 \ar[dd]^{\cong}\\
&&\\
\gwheelcset\left(\varGamma[G], \sK\right) \ar[rr]^-{\epsilon_G^*} && \gwheelcset\left(\Sc[G], \sK\right) 
}\]
that is natural in $G$.
\end{lemma}

\begin{proof}
The two vertical bijections are by Yoneda's Lemma.  The commutativity of the diagram is by the definition of the Segal map $\chi_G$.
\end{proof}

\begin{remark}
Using Lemma \ref{lem:wsegalcoredesc}, one can also think of the wheeled properadic Segal map $\chi_G$ as the map $\epsilon_G^*$.
\end{remark}

With the above definition of the Segal maps, we can now state the main result of this section.  It is the wheeled version of Theorem \ref{properadnerve}.

\begin{theorem}
\label{wproperadnerve}
Suppose $\sK \in \gwheelcset$.  Then the following statements are equivalent.
\begin{enumerate}
\item
There exist a wheeled properad $\sP$ and an isomorphism $\sK \cong N\sP$.
\item
$\sK$ satisfies the wheeled properadic Segal condition.
\item
$\sK$ is a strict $\infty$-wheeled properad.\index{strict infinity wheeled properad!characterization}
\end{enumerate}
\end{theorem}

The rest of this section is devoted to the proof of this theorem.

\subsection{Wheeled Properadic Nerves Satisfy the Segal Condition}

The following observation says that the wheeled properadic nerve always yields a graphical set that satisfies the wheeled properadic Segal condition.

\begin{lemma}
\label{lem:wnerveissegal}
The wheeled properadic nerve of every wheeled properad satisfies the Segal condition.
\end{lemma}

\begin{proof}
We reuse the proof of Lemma \ref{lem:nerveissegal}, applied to ordinary $G \in \gwheelc$ with at least one internal edge. For the descriptions of graphices in the nerve and of the corolla ribbon, instead of Lemma \ref{npgelement} and Remark \ref{rk:ribbongeometric}, here we use Lemma \ref{wnpgelement} and  Remark \ref{rk:wribbongeometric}.
\end{proof}

\subsection{Wheeled Properadic Nerve is Fully Faithful}

Before showing that the Segal condition implies a strict $\infty$-wheeled properad, we briefly digress here to observe that the wheeled properadic nerve is fully faithful.

\begin{proposition}
\label{wheeled-nerve-fully-faithful}
The wheeled properadic nerve
\[
\nicexy{
\wproperad \ar[r]^-{N} & \gwheelcset
}\]
is fully faithful.
\end{proposition}

\begin{proof}
We reuse most of the proof of Proposition \ref{nerve-fully-faithful}, which establishes the fully faithfulness of the properadic nerve.  Suppose $\sP$ is a $\fC$-colored wheeled properad and $\sQ$ is a $\fD$-colored wheeled properad.  

The proof that the wheeled properadic nerve is faithful is exactly the same as the properad case.

To show that the wheeled properadic nerve $N$ is full, suppose given a map $\phi : N\sP \to N\sQ \in \gwheelcset$.  We must show that $\phi = N\varphi$ for some wheeled properad map $\varphi : \sP \to \sQ$.   Using $\uparrow$ and corollas and proceeding as in the properad case, we obtain a map $\varphi : \sP \to \sQ$ (of color sets and components) whose image under the wheeled properadic nerve agrees with $\phi$ on $\uparrow$ and corollas.   As soon as we know that $\varphi$ is a map of wheeled properads, we will know that $N\varphi = \phi$ because the wheeled properadic nerve of a wheeled properad satisfies the Segal condition (Lemma \ref{lem:wnerveissegal}), hence is completely determined by what it does to corollas and $\uparrow$.

Now we show that $\varphi$ is a map of wheeled properads.  Again as in the properad case, we know that $\varphi$ respects the colored units and the bi-equivariant structures.  To see that $\varphi$ respects the dioperadic compositions and the contractions, we use:
\begin{itemize}
\item
the inner dioperadic coface maps $C \to D \in \Gammaw$ with $D$ a dioperadic graph and $C$ obtained from $D$ by shrinking away the internal edge;
\item
the inner contracting coface maps $C \to K \in \Gammaw$ with $K$ a contracted corolla and $C$ obtained from $K$ by deleting the loop.
\end{itemize}
Then the maps
\[
\nicexy{(N\sP)(D) \ar[r] & (N\sP)(C)} \andspace
\nicexy{(N\sP)(K) \ar[r] & (N\sP)(C)} 
\]
are given by the dioperadic composition of $\sP$ and the contraction of $\sP$.  So the fact that $\phi$ is compatible with $C \to D$ and $C \to K$ implies that $\varphi$ respects the dioperadic compositions and the contractions.
\end{proof}

\subsection{The Segal Condition Implies Strict \texorpdfstring{$\infty$}{∞}-Wheeled Properad}

Next we observe that a graphical set that satisfies the wheeled properadic Segal condition must be a strict $\infty$-wheeled properad.  As in the properad case, we first eliminiate the possibility of an inner horn having multiple fillers.

\begin{lemma}
\label{lem:wsegalisstrict}
Suppose $\sK \in \gwheelcset$ satisfies the wheeled properadic Segal condition.  Then each inner horn
\[
\nicexy@C+12pt{
\Lambdaw^{d^u}[G] \ar[d] \ar[r]^-{\varphi} & \sK\\
\Gammaw[G] \ar@{.>}[ur]
}\]
of $\sK$ has at most one filler.
\end{lemma}

\begin{proof}
We adapt the proof of Lemma \ref{lem:segalisstrict}.  Suppose $G \in \gwheelc$ is ordinary with at least one internal edge.  By Remark \ref{rk:whornofx} the inner horn $\varphi$ is uniquely determined by the graphices
\[
\left\{\varphi_J \in \sK(J) ~|~  
\nicexy{J \ar[r] & G}
\text{ face not equal to $d^u$}\right\}
\]
that agree on common faces. By the assumed Segal condition of $\sK$, each graphex $\varphi_J$ is uniquely determined by the graphices
\[
\left\{\varphi_w = \xi_w^* \left(\varphi_J\right) \in \sK(C_w) ~|~ w \in \vertex(J)\right\},
\]
where $\nicexy{C_w \ar[r]^-{\xi_w} & J}$ is the corolla inclusion, which decomposes into outer coface maps (Lemma \ref{wcvtog}). Furthermore, the compatibility condition of the $\varphi_J$ implies that $\varphi_w$ is well-defined if $w$ is a vertex in two such faces of $G$.  There are now two cases.
\begin{enumerate}
\item
If $G$ is simply connected, then since $\edgei(G) \not= \varnothing$, $G$ has at least two vertices, and hence at least two deletable vertices by Lemma \ref{deletableexists}.  Now we reuse the last two paragraphs of the proof of Lemma \ref{lem:segalisstrict} to conclude that the inner horn $\varphi$ can have at most one filler
\[
\Phi = \left\{\varphi_v \in \sK(C_v) ~|~ v \in \vertex(G)\right\} \in \sK(G).
\]
\item
Next suppose $G$ is not simply connected, so it has a cycle.  By Lemma \ref{disconnectablecycle} there exists a disconnectable edge $e$ in $G$.  By Lemma \ref{lem:geconnected} there is an outer contracting factorization
\[
G = [\xi_e C](G_e),
\]
where $G_e$ is the connected graph obtained from $G$ by disconnecting $e$, and $C$ is the corolla with the same profiles as $G_e$.  See the picture in Remark \ref{rk:gepicture}.  Since there is an outer contracting coface map
\[
\nicexy{G_e \ar[r]^-{d^e} & G},
\]
by the assumed Segal condition we have a graphex
\[
\varphi_{G_e} = 
\left\{\varphi_v = \xi_v^*\varphi_{G_e} \in \sK(C_v) ~|~ v \in \vertex(G_e) \right\} \in \sK(G_e).
\]
Since $\vertex(G_e) = \vertex(G)$, by the Segal condition again, $\varphi_{G_e} = \{\varphi_v\}$ actually determines a unique graphex
\[
\Phi = \{\varphi_v ~|~ v \in \vertex(G) \} \in \sK(G).
\]
This graphex $\Phi$ is the only possible candidate for a filler of the inner horn $\varphi$ because any such filler must restrict to the graphex $\varphi_{G_e}$, and hence also to the $\varphi_v$ for $v \in \vertex(G)$.
\end{enumerate}
\end{proof}

\begin{lemma}
\label{lem2:wsegalisstrict}
Suppose $\sK \in \gwheelcset$ satisfies the wheeled properadic Segal condition.  Then $\sK$ is a strict $\infty$-wheeled properad.
\end{lemma}

\begin{proof}
We need to show that, in the context of Lemma \ref{lem:wsegalisstrict}, the graphex
\[
\Phi = \{\varphi_v\} \in \sK(G)
\]
is a filler of the inner horn $\varphi$.  It remains to show that $\Phi$ restricts to the $\varphi_J$ for faces $J \to G$ not equal to the given inner face $d^u$. We consider outer and inner coface maps separately.
\begin{enumerate}
\item
For an \emph{outer} coface map $\nicexy{J_z \ar[r]^-{d_z} & G}$, $J_z$ is obtained from $G$ by either deleting a deletable vertex or disconnecting a disconnectable edge.  In either case, each vertex in $J_z$ is also a vertex in $G$.  Therefore, $\Phi$ restricts to $\varphi_{J_z} \in \sK(J_z)$, i.e.,
\[
\varphi_{J_z} 
= \{\varphi_v ~|~ v \in \vertex(J_z)\} 
= d_z^{*}(\Phi).
\]
\item
Suppose $\nicexy{J_v \ar[r]^{d_v} & G}$ is an \emph{inner} coface map not equal to the given $d^u$.  Since $G$ has at least two inner faces, it has at least two internal edges, and $J_v$ has at least one internal edge.  There is a graph substitution factorization
\[
G = J_v(H),
\]
where $H$ is either a dioperadic graph or a contracted corolla.  Suppose $H$ is substituted into $v \in \vertex(J_v)$.  We need to show the equality
\[
\varphi_w 
= \xi_w^* d_v^{*} \Phi \in \sK(C_w)
\]
for every $w \in \vertex(J_v)$. First we consider vertices in $J_v$ not equal to $v$.
\begin{enumerate}
\item
A vertex $x \in \vertex(J_v) \setminus \{v\}$ is already in $G$.  By Lemma \ref{wcvtog} there is an outer coface map $\nicexy{J_z \ar[r]^-{d_z} & G}$ with $x \in \vertex(J_z)$.  The diagram
\[
\nicearrow
\xymatrix{
C_x \ar[d]_{\xi_x} \ar[r]^-{\xi_x} & J_z \ar[d]^{d_z}\\
J_v \ar[r]^-{d_v} & G
}\]
in $\Gammaw$ is commutative because either composition is the corolla inclusion $C_x \to G$.  The first part about outer coface map and the commutativity of this diagram now imply:
\[
\begin{split}
\varphi_x 
&= \xi_x^* \varphi_{J_z}\\
&= \xi_x^* d_z^{*} \Phi\\
&= \xi_x^* d_v^{*} \Phi \in \sK(C_x).
\end{split}
\]
\item
It remains to prove the equality
\[
\varphi_v = \xi_v^* d_v^{*} \Phi.
\]
By the given graph substitution $G = J_v(H)$ and the characterization Theorem \ref{thm:gwheelcsubgraph} of subgraphs, the map $H \to G$ is a subgraph.  Therefore, there is a commutative diagram
\[
\nicearrow
\xymatrix{
C_v \ar[d]_{\xi_v} \ar[r]^-{d} 
& H \ar[r]^-{\xi} 
& J_w \ar[d]^{d_w}\\
J_v \ar[rr]^-{d_v} && G
}\]
in $\Gammaw$.  Here each composition is the corolla inclusion $C_v \to G$, $d$ is the unique inner dioperadic or contracting coface map into $H$, $\xi$ is a composition of outer coface maps, and $d_w$ is an outer coface map. The above commutative square now implies:
\[
\begin{split}
\varphi_v 
&= \xi_v^* \varphi_{J_v}\\
&= (\xi d)^* \varphi_{J_w}\\
&= d^{*} \xi^* d_w^* \Phi\\
&= \xi_v^* d_v^* \Phi.
\end{split}
\]
This proves that $\Phi \in \sK(G)$ restricts to $\varphi_v \in \sK(C_v)$ as well.
\end{enumerate}
\end{enumerate}
\end{proof}

\subsection{Wheeled Properad Associated to a Strict \texorpdfstring{$\infty$}{∞}-Wheeled Properad}

Here we observe that each strict $\infty$-wheeled properad has a canonically associated wheeled properad.  We will use the biased form of a wheeled properad, which was defined in Definition \ref{def:wproperad}.

\begin{definition}
\label{def:wproperadpk}
Fix a strict $\infty$-wheeled properad $\sK$.  We define a $\sK(\uparrow)$-colored wheeled properad $\sP_{\sK}$ as follows.
\begin{enumerate}
\item
The $\Sigma_{\sS(\sK(\uparrow))}$-bimodule structure of  $\sP_{\sK}$ is that of $\sK \in \gwheelcset$ as in section \ref{rk:gsetproperad}.  In particular, its elements are the $1$-dimensional elements in $\sK$, with 
\begin{itemize}
\item
input/output profiles induced by the outer dioperadic coface maps $\uparrow ~\to C$ into corollas, 
\item
colored units $\bone_c = s^*(c)$ induced by the codegeneracy map  $\nicexy{C_{(1;1)} \ar[r]^-{s} & \uparrow}$, and 
\item
bi-equivariant structure induced by input/output relabeling.
\end{itemize}
\item
For the dioperadic composition, suppose
\begin{itemize}
\item
$\theta \in \sK(C_1)$ has profiles $\yxh$,
\item
$\phi \in \sK(C_2)$ has profiles $\wvh$, and
\item
$x_i = w_j \in \sK(\uparrow)$.
\end{itemize}
Define $D$ as the dioperadic graph (Example \ref{ex:dioperadic}) whose top (resp., bottom) vertex $v_1$ (resp., $v_2$)  has the same profiles as $C_1$ (resp., $C_2$), and whose unique internal edge connects the $j$th outgoing flag of $v_2$ and the $i$th incoming flag of $v_1$.  By Example \ref{ex:diophorn} the $1$-dimensional elements $\theta$ and $\phi$ determine a unique inner horn
\[
\nicearrow
\xymatrix@C+12pt{
\Lambdaw^d[D] \ar[d] \ar[r]^-{\varphi} & \sK\\
\Gammaw[D] \ar@{.>}[ur]_-{\exists ! \Phi}
}\]
of $\sK$, where
\[
\nicexy{C \ar[r]^-{d} & D}
\]
is the unique inner dioperadic coface map into $D$.  Since $\sK$ is assumed to be a strict $\infty$-wheeled properad, there is a unique filler $\Phi \in \sK(D)$.  Now we define the \textbf{dioperadic composition} \index{dioperadic composition!of $1$-dimensional elements}
\[
\theta \jcompi \phi = d^*\Phi \in \sK(C).
\]
By construction the dioperadic composition $\theta \jcompi \phi$ has profiles $\binom{\uw \compj \uy}{\ux \compi \uv}$.
\item
For the contraction, suppose $\theta \in \sK(C)$ has profiles $\yxh$ such that $x_j = y_i \in \sK(\uparrow)$.  Define $\xi_eC$ as the contracted corollas (Example \ref{ex:contractedcor}) whose unique vertex $v$ has the same profiles as the corolla $C$, and whose loop $e$ connects the $i$th outgoing flag of $v$ with the $j$th incoming flag of $v$.  By Example \ref{ex:contracthorn} the $1$-dimensional element $\theta$ determines a unique inner horn
\[
\nicexy@C+12pt{
\Lambdaw^d\left[\xi_e C\right] \ar[r]^-{\varphi} \ar[d] & \sK\\
\Gammaw\left[\xi_eC\right] \ar@{.>}[ur]_-{\exists ! \Phi} &}
\]
of $\sK$, where
\[
\nicexy{C' \ar[r]^-{d} & \xi_eC}
\]
is the unique inner contracting coface map into $\xi_eC$.  Here $C'$ is the corolla obtained from the contracted corolla $\xi_eC$ by deleting the loop $e$.  Since $\sK$ is assumed to be a strict $\infty$-wheeled properad, there is a unique filler $\Phi \in \sK(\xi_eC)$.  Now we define the \textbf{contraction} \index{contraction!of $1$-dimensional element}
\[
\xiij(\theta) = d^*\Phi \in \sK(C').
\]
By construction the contraction $\xiij(\theta)$ has profiles $\binom{\uy\setminus y_i}{\ux\setminus x_j}$.
\end{enumerate}
This finishes the definition of $\sP_{\sK}$.
\end{definition}

\begin{lemma}
\label{wproperadpk}
For each strict $\infty$-wheeled properad $\sK$, the above definitions give a $\sK(\uparrow)$-colored wheeled properad $\sP_{\sK}$.\index{strict infinity wheeled properad!associated properad}
\end{lemma}

\begin{proof}
We need to check the various biased axioms of a wheeled properad.  They are all proved in the same manner using the unity and associativity of graph substitution, which are proved with full details and generality in \cite{jy2}.  So we will prove only one compatibility axiom between the dioperadic composition and the contraction in detail to illustrate the method.

Suppose $B$ is the $\sK(\uparrow)$-colored connected graph
\begin{center}
\begin{tikzpicture}
\matrix[row sep=1cm, column sep=1cm]{
\node [plain] (v1) {$v_1$};\\
\node [plain] (v2) {$v_2$};\\
};
\draw [arrow] (v2) to node{\footnotesize{$e$}} (v1);
\draw [arrow, looseness=2, out=45, in=-45] (v1) to node{\footnotesize{$f$}} (v2);
\foreach \x in {1,2}
{
\draw [inputleg] (v\x) to +(-.5cm,-.3cm);
\draw [inputleg] (v\x) to +(.5cm,-.3cm);
\draw [outputleg] (v\x) to +(-.5cm,.3cm);
\draw [outputleg] (v\x) to +(.5cm,.3cm);
}
\end{tikzpicture}
\end{center}
in which:
\begin{itemize}
\item
$v_1$ (resp., $v_2$) has profiles $\yxh$ (resp., $\wvh$),
\item
$e$ connects the $s$th outgoing flag of $v_2$ with the $r$th incoming flag of $v_1$, and
\item
$f$ connects the $i$th outgoing flag of $v_1$ with the $j$th incoming flag of $v_2$.
\end{itemize}
Suppose $\theta \in \sK(C_1)$ with profiles $\yxh$ and $\phi \in \sK(C_2)$ with profiles $\wvh$ are $1$-dimensional elements such that $w_s = x_r$ and $y_i = v_j$.  So the decorated graph
\begin{center}
\begin{tikzpicture}
\matrix[row sep=1cm, column sep=1cm]{
\node [plain] (v1) {$\theta$};\\
\node [plain] (v2) {$\phi$};\\
};
\draw [arrow] (v2) to node{\footnotesize{$e$}} (v1);
\draw [arrow, looseness=2, out=45, in=-45] (v1) to node{\footnotesize{$f$}} (v2);
\foreach \x in {1,2}
{
\draw [inputleg] (v\x) to +(-.5cm,-.3cm);
\draw [inputleg] (v\x) to +(.5cm,-.3cm);
\draw [outputleg] (v\x) to +(-.5cm,.3cm);
\draw [outputleg] (v\x) to +(.5cm,.3cm);
}
\end{tikzpicture}
\end{center}
makes sense.  We want to prove the compatibility axiom
\begin{equation}
\label{xicompcompat}
\xi_f(\theta \circ_e \phi)
= 
\sigma\left[ \xi_e(\phi \circ_f \theta)\right]\tau,
\end{equation}
in which
\[
\xi_e = \xi_f = \xi^{s-1+i}_{r-1+j}, \quad
\circ_e = \scompr, \quad
\circ_f = \icompj,
\]
and $\sigma$ and $\tau$ are suitable block permutations that align the input/output profiles of $\xi_e(\phi \circ_f \theta)$ with the left-hand side.

As in previous chapters, we will ignore the isomorphisms induced by relabelings to simplify the presentation.  To prove \eqref{xicompcompat} we use the same strategy as in the proof of Lemma \ref{properadpk}.  In other words, we will show that both sides of \eqref{xicompcompat} are equal to a certain codimension $2$ face of some inner horn filler for $B$.

To do this, first observe that $B$ has four faces.
\begin{enumerate}
\item
There is an outer contracting coface map
\[
\nicexy{B_f \ar[r]^-{d_f} & B}
\]
corresponding to disconnecting the disconnectable edge $f$.
\begin{center}
\begin{tikzpicture}
\matrix[row sep=1cm, column sep=1.5cm]{
\node [plain] (v11) {$v_1$};
&&& \node [plain] (v12) {$v_1$};\\
\node [plain] (v21) {$v_2$};
& \node [empty] (s) {};
& \node [empty] (t) {};
& \node [plain] (v22) {$v_2$};\\
};
\draw [arrow] (v21) to node{\footnotesize{$e$}} (v11);
\draw [arrow] (v22) to node{\footnotesize{$e$}} (v12);
\draw [arrow, looseness=2, out=45, in=-45] (v12) to node{\footnotesize{$f$}} (v22);
\foreach \x in {1,2}
\foreach \y in {1,2}
{
\draw [inputleg] (v\x\y) to +(-.5cm,-.3cm);
\draw [inputleg] (v\x\y) to +(.5cm,-.3cm);
\draw [outputleg] (v\x\y) to +(-.5cm,.3cm);
\draw [outputleg] (v\x\y) to +(.5cm,.3cm);
}
\draw [inputleg] (v21) to node[below=.1cm]{\footnotesize{$f_{1}$}} +(0,-.6cm);
\draw [outputleg] (v11) to node[above=.1cm]{\footnotesize{$f_{-1}$}} +(0,.6cm);
\draw [arrow] (s) to node{$d_f$} (t);
\end{tikzpicture}
\end{center}
The dioperadic graph $B_f$ has a unique inner face, namely, the inner dioperadic coface map
\[
\nicexy{(B_f)_{12} \ar[r]^-{d_0} & B_f},
\]
where $(B_f)_{12}$ is the corolla obtained from $B_f$ by shrinking away the internal edge $e$.  By Example \ref{ex:diophorn} the elements $\theta$ and $\phi$ determine an inner horn
\[
\nicearrow
\xymatrix@C+12pt{
\Lambdaw^{d_0}[B_f] \ar[d] \ar[r] & \sK\\
\Gammaw[B_f] \ar@{.>}[ur]_-{\exists ! \Phi_0}
}\]
of $\sK$, which has a unique filler $\Phi_0 \in \sK(B_f)$.  By construction 
\[
\theta \circ_e \phi = d_0^*\Phi_0 \in \sK\left((B_f)_{12}\right).
\]
\item
There is an outer contracting coface map
\[
\nicexy{B_e \ar[r]^-{d_e} & B}
\]
corresponding to disconnecting the disconnectable edge $e$.
\begin{center}
\begin{tikzpicture}
\matrix[row sep=1cm, column sep=1.5cm]{
\node [plain] (v11) {$v_2$};
&&& \node [plain] (v12) {$v_1$};\\
\node [plain] (v21) {$v_1$};
& \node [empty] (s) {};
& \node [empty] (t) {};
& \node [plain] (v22) {$v_2$};\\
};
\draw [arrow] (v21) to node{\footnotesize{$f$}} (v11);
\draw [arrow] (v22) to node{\footnotesize{$e$}} (v12);
\draw [arrow, looseness=2, out=45, in=-45] (v12) to node{\footnotesize{$f$}} (v22);
\foreach \x in {1,2}
\foreach \y in {1,2}
{
\draw [inputleg] (v\x\y) to +(-.5cm,-.3cm);
\draw [inputleg] (v\x\y) to +(.5cm,-.3cm);
\draw [outputleg] (v\x\y) to +(-.5cm,.3cm);
\draw [outputleg] (v\x\y) to +(.5cm,.3cm);
}
\draw [inputleg] (v21) to node[below=.1cm]{\footnotesize{$e_{1}$}} +(0,-.6cm);
\draw [outputleg] (v11) to node[above=.1cm]{\footnotesize{$e_{-1}$}} +(0,.6cm);
\draw [arrow] (s) to node{$d_e$} (t);
\end{tikzpicture}
\end{center}
The dioperadic graph $B_e$ has a unique inner face, namely, the inner dioperadic coface map
\[
\nicexy{(B_e)_{21} \ar[r]^-{d_1} & B_e},
\]
where $(B_e)_{21}$ is the corolla obtained from $B_e$ by shrinking away the internal edge $f$.  By Example \ref{ex:diophorn} the elements $\phi$ and $\theta$ determine an inner horn
\[
\nicearrow
\xymatrix@C+12pt{
\Lambdaw^{d_1}[B_e] \ar[d] \ar[r] & \sK\\
\Gammaw[B_e] \ar@{.>}[ur]_-{\exists ! \Phi_1}
}\]
of $\sK$, which has a unique filler $\Phi_1 \in \sK(B_e)$.  By construction 
\[
\phi \circ_f \theta = d_1^*\Phi_1 \in \sK\left((B_e)_{21}\right).
\]
\item
There is an inner dioperadic coface map
\[
\nicexy{B_{12} \ar[r]^-{d_{12}} & B}
\]
corresponding to shrinking away the internal edge $e$.
\begin{center}
\begin{tikzpicture}
\matrix[row sep=1cm, column sep=1.5cm]{
&&& \node [plain] (v12) {$v_1$};\\
\node [plain] (v21) {$v_1v_2$};
& \node [empty] (s) {};
& \node [empty] (t) {};
& \node [plain] (v22) {$v_2$};\\
};
\draw [arrow, out=60, in=-60, loop] (v21) to node{\footnotesize{$f$}} ();
\draw [arrow] (v22) to node{\footnotesize{$e$}} (v12);
\draw [arrow, looseness=2, out=45, in=-45] (v12) to node{\footnotesize{$f$}} (v22);
\foreach \x in {1,2}
{
\draw [inputleg] (v\x2) to +(-.5cm,-.3cm);
\draw [inputleg] (v\x2) to +(.5cm,-.3cm);
\draw [outputleg] (v\x2) to +(-.5cm,.3cm);
\draw [outputleg] (v\x2) to +(.5cm,.3cm);
}
\draw [inputleg] (v21) to +(-.5cm,-.3cm);
\draw [inputleg] (v21) to +(.5cm,-.3cm);
\draw [outputleg] (v21) to +(-.5cm,.3cm);
\draw [outputleg] (v21) to +(.5cm,.3cm);
\draw [arrow] (s) to node{$d_{12}$} (t);
\end{tikzpicture}
\end{center}
The contracted corolla $B_{12}$ has a unique inner face, namely, the inner contracting coface map
\[
\nicexy{B_{12}' \ar[r]^-{d_2} & B_{12}},
\]
where $B_{12}'$ is the corolla obtained from $B_{12}$ by deleting the loop $f$.  By Example \ref{ex:contracthorn} the graphex $\theta \circ_e \phi$ determines an inner horn
\[
\nicearrow
\xymatrix@C+12pt{
\Lambdaw^{d_2}[B_{12}] \ar[d] \ar[r] & \sK\\
\Gammaw[B_{12}] \ar@{.>}[ur]_-{\exists ! \Phi_2}
}\]
of $\sK$, which has a unique filler $\Phi_2 \in \sK(B_{12})$.  By construction 
\[
\xi_f\left(\theta \circ_e \phi\right) 
= 
d_2^*\Phi_2 \in \sK\left(B_{12}'\right).
\]
\item
Finally, there is an inner dioperadic coface map
\[
\nicexy{B_{21} \ar[r]^-{d_{21}} & B}
\]
corresponding to shrinking away the internal edge $f$.
\begin{center}
\begin{tikzpicture}
\matrix[row sep=1cm, column sep=1.5cm]{
&&& \node [plain] (v12) {$v_1$};\\
\node [plain] (v21) {$v_2v_1$};
& \node [empty] (s) {};
& \node [empty] (t) {};
& \node [plain] (v22) {$v_2$};\\
};
\draw [arrow, out=60, in=-60, loop] (v21) to node{\footnotesize{$e$}} ();
\draw [arrow] (v22) to node{\footnotesize{$e$}} (v12);
\draw [arrow, looseness=2, out=45, in=-45] (v12) to node{\footnotesize{$f$}} (v22);
\foreach \x in {1,2}
{
\draw [inputleg] (v\x2) to +(-.5cm,-.3cm);
\draw [inputleg] (v\x2) to +(.5cm,-.3cm);
\draw [outputleg] (v\x2) to +(-.5cm,.3cm);
\draw [outputleg] (v\x2) to +(.5cm,.3cm);
}
\draw [inputleg] (v21) to +(-.5cm,-.3cm);
\draw [inputleg] (v21) to +(.5cm,-.3cm);
\draw [outputleg] (v21) to +(-.5cm,.3cm);
\draw [outputleg] (v21) to +(.5cm,.3cm);
\draw [arrow] (s) to node{$d_{21}$} (t);
\end{tikzpicture}
\end{center}
The contracted corolla $B_{21}$ has a unique inner face, namely, the inner contracting coface map
\[
\nicexy{B_{21}' \ar[r]^-{d_3} & B_{21}},
\]
where $B_{21}'$ is the corolla obtained from $B_{21}$ by deleting the loop $e$.
\end{enumerate}

By Lemma \ref{whorndescription} the graphices $\Phi_i$ for $0 \leq i \leq 2$ determine an inner horn
\[
\nicearrow
\xymatrix@C+12pt{
\Lambdaw^{d_{21}}[B] \ar[d] \ar[r] & \sK\\
\Gammaw[B] \ar@{.>}[ur]_-{\exists ! \Psi}
}\]
of $\sK$, which has a unique filler $\Psi \in \sK(B)$.  Note that $B_{12}' = B_{21}'$, which is the corolla obtained from $B$ by shrinking away both internal edges.  So the diagram
\begin{equation}
\label{bonetwoprime}
\nicexy@C+10pt{
B_{12}' = B_{21}' \ar[d]_{d_2} \ar[r]^-{d_3} 
& B_{21} \ar[d]^{d_{21}}\\
B_{12} \ar[r]^-{d_{12}} & B
}
\end{equation}
is commutative.  Pulling $\Psi$ back along one composition yields
\[
\begin{split}
\xi_f(\theta \circ_e \phi) 
&= d_2^*\Phi_2\\
&= d_2^*d_{12}^*\Psi.
\end{split}
\]
So to finish the proof it suffices to show that
\begin{equation}
\label{xiphitheta}
\xi_e\left(\phi \circ_f \theta\right) = d_3^* d_{21}^* \Psi.
\end{equation}

To prove \eqref{xiphitheta}, observe that by Example \ref{ex:contracthorn} $\phi \circ_f \theta$ determines an inner horn
\[
\nicearrow
\xymatrix@C+12pt{
\Lambdaw^{d_{3}}[B_{21}] \ar[d] \ar[r] & \sK\\
\Gammaw[B_{21}] \ar@{.>}[ur]_-{\exists ! \Upsilon}
}\]
of $\sK$, which has a unique filler $\Upsilon \in \sK(B_{21})$.  By definition the contraction $\xi_e(\phi \circ_f \theta)$ is the unique inner face of $\Upsilon$.  Therefore, to prove \eqref{xiphitheta}, it suffices to show that
\[
d_{21}^*\Psi = \Upsilon.
\]
By uniqueness of $\Upsilon$, this is equivalent to
\[
d_e^* d_{21}^*\Psi = \phi \circ_f \theta.
\]
Here
\[
\nicexy{(B_{21})_e \ar[r]^-{d_e} & B_{21}}
\]
is the unique outer coface map into $B_{21}$, so $(B_{21})_e$ is obtained from $B_{21}$ by disconnecting the disconnectable edge $e$.  Note that
\[
(B_e)_{21} = (B_{21})_e,
\]
which is the corolla obtained from $B$ by disconnecting $e$ and shrinking away $f$ in either order.  Therefore, the diagram
\begin{equation}
\label{betwoone}
\nicexy@C+10pt{
(B_e)_{21} = (B_{21})_e \ar[d]_{d_1} \ar[r]^-{d_e} 
& B_{21} \ar[d]^{d_{21}}\\
B_{e} \ar[r]^-{d_{e}} & B
}
\end{equation}
is commutative.  So we have
\[
\begin{split}
d_e^* d_{21}^*\Psi
&= d_1^* d_e^* \Psi\\
&= d_1^* \Phi_1\\
&= \phi \circ_f \theta.
\end{split}
\]
As discussed above, this suffices to prove \eqref{xiphitheta}.
\end{proof}

\subsection{Strict \texorpdfstring{$\infty$}{∞}-Wheeled Properads are Nerves}

We now show that each strict $\infty$-wheeled properad $\sK$ is canonically isomorphic to the nerve of the wheeled properad $\pk$ in Lemma \ref{wproperadpk}.  The strategy is the same as in section \ref{sec:fundpropsiprop} for strict $\infty$-properads.  First we observe that there is an object-wise bijection.

\begin{lemma}
\label{wetaiso}
Suppose $\sK$ is a strict $\infty$-wheeled properad, and $G \in \gwheelc$.  Then there is a canonical bijection
\[
\nicearrow
\xymatrix{\sK(G) \ar[r]^-{\eta_G}_-{\cong} & (N\sP_{\sK})(G),
}\]
where $\sP_{\sK}$ is the wheeled properad in Lemma \ref{wproperadpk}.
\end{lemma}

\begin{proof}
We adapt the proof of Lemma \ref{etaiso} as follows.  The map $\eta$ is constructed according to the number of internal edges in $G$.
\begin{enumerate}
\item
For the exceptional edge $\uparrow$, by definition there is a bijection
\[
\nicearrow
\xymatrix{
\sK(\uparrow) \ar[d]_-{\eta_{\uparrow}}^-{\cong}\\
(N\sP_{\sK})(\uparrow) 
= \wproperad(\uparrow,\sP_{\sK})
}\]
because by Lemma \ref{wnpgelement} a map $\uparrow ~ \to \sP_{\sK}$ is simply an element in the color set of $\sP_{\sK}$, which is $\sK(\uparrow)$.
\item
For a permuted corolla $\sigma C\tau$, the bijection
\[
\nicearrow
\xymatrix{
\sK(\sigma C \tau) \ar[d]_-{\eta_{\sigma C \tau}}^-{\cong}\\
(N\sP_{\sK})(\sigma C\tau) 
= \wproperad(\sigma C \tau,\sP_{\sK})
}\]
comes from the definition of the elements in $\pk$ as the $1$-dimensional elements in $\sK$.
\item
The previous two cases take care of the case $\edgei(G) = \varnothing$.  The only \emph{exceptional} connected graph with $\edgei(G) \not= \varnothing$ is the exceptional wheel.  The bijection $\eta_{\wheel}$ is defined as the composition of isomorphisms
\[
\nicexy@C+10pt{
\sK(\wheel) \ar[d]_{\cong} \ar[r]^-{\eta_{\wheel}} 
& N\pk(\wheel) \\
\sK(\uparrow) \ar[r]^-{\eta_{\uparrow}}_-{\cong}
& N\pk(\uparrow) \ar[u]_{\cong}.
}\]
The two vertical isomorphisms are induced by the isomorphism $\Gammaw(\wheel) \cong \Gammaw(\uparrow)$ in $\Gammaw$ (Example \ref{ex:unitgwproperad}).
\item
Suppose $G \in \gwheelc$ is ordinary with $\edgei(G) \not= \varnothing$.  The map $\eta_G$ is defined as the composition
\begin{equation}
\label{wetankg}
\nicearrow
\xymatrix@C+12pt{
\sK(G) \ar[r]^-{\eta_G} \ar[d]_-{\chi_G} 
& (N\sP_{\sK})(G)\\
\left[\prod \sK(C_v)\right]_{\sK(\uparrow)} \ar[r]^-{\prod \eta_C}_-{\cong} 
& \left[\prod (N\pk)(C_v)\right]_{N\pk(\uparrow)} \ar[u]^-{\cong}_-{\chi_G^{-1}}\\
}
\end{equation}
Here the $\chi_G$ are the Segal maps \eqref{wsegalmap}, and $\prod \eta_C$ is induced by the bijections defined previously for corollas and the exceptional edge.  The Segal map for the nerve $N\pk$ is a bijection by Lemma \ref{lem:wnerveissegal}.
\end{enumerate}

Next we show that $\eta_G$ is a bijection by induction on $n=|\edgei(G)|$.  We already observed that $\eta_G$ are bijections for $G \in \{\uparrow, \wheel, \sigma C\tau\}$.  So we now assume $G$ is ordinary with $n \geq 1$.  Then there must be at least one inner dioperadic or contracting coface map $\nicexy{K \ar[r]^-{d_u} & G}$.  Using the inner horn $\Lambdaw^u[G] \to \Gammaw[G]$ defined by $d_u$, the diagram
\begin{equation}
\label{wknpeta}
\nicearrow
\xymatrix@C+10pt{
\gwheelcset\left(\Gammaw[G], \sK\right) \ar[d]_-{\cong} \ar[r]^-{\eta_G} 
& \gwheelcset\left(\Gammaw[G], N\sP_{\sK}\right) \ar[d]^-{\cong}\\
\gwheelcset\left(\Lambdaw^u[G], \sK\right) \ar[r]^-{\eta} & \gwheelcset\left(\Lambdaw^u[G],N\sP_{\sK}\right)
}
\end{equation}
is commutative by the construction of the maps $\eta$.  We want to show that $\eta_G$ is a bijection.  The left vertical map is a bijection by the strict $\infty$-wheeled properad assumption on $\sK$.  The right vertical map is a bijection because $N\sP_{\sK}$ is also a strict $\infty$-wheeled properad (Lemmas \ref{lem:wnerveissegal},  \ref{lem:wsegalisstrict}, and \ref{lem2:wsegalisstrict}).  The bottom horizontal map is determined by the maps
\[
\nicearrow
\xymatrix@C+10pt{
\gwheelcset\left(\Gammaw[J], \sK\right) \ar[r]^-{\eta_J} & \gwheelcset\left(\Gammaw[J],N\sP_{\sK}\right)
}\]
for coface maps $J \to G$ not equal to $d_u$ (Lemma \ref{whorndescription}).  Each such $J$ is ordinary and has $n-1$ internal edges.  So by induction hypothesis, the bottom horizontal map in \eqref{wknpeta} is a bijection.  Therefore, the top horizontal map $\eta_G$ is also a bijection.
\end{proof}

Next we observe that the maps $\eta_G$ in the previous lemma have the same universal property as the unit of the adjunction object-wise.

\begin{lemma}
\label{wetaknerveptwo}
Suppose $\sK$ is a strict $\infty$-wheeled properad, and $\nicearrow\xymatrix{\sK \ar[r]^-{\zeta} & N\sQ}$ is a map in $\gwheelcset$ for some $\fD$-colored wheeled properad $\sQ$.  Then there exists a unique map $\nicearrow\xymatrix{\sP_{\sK} \ar[r]^-{\zeta'} & \sQ}$ of wheeled properads such that the diagram
\[
\nicearrow
\xymatrix{
\sK(G) \ar[d]_-{\eta_G}^-{\cong} \ar[r]^-{\zeta} & (N\sQ)(G)\\
(N\sP_{\sK})(G) \ar[ur]_-{N\zeta'} & 
}\]
is commutative for each $G \in \gwheelc$.
\end{lemma}

\begin{proof}
We reuse the proof of Lemma \ref{etaknerveptwo} as follows.  The map $\nicexy{\pk \ar[r]^-{\zeta'} & \sQ}$ on color sets is the map
\[
\nicexy{
\sK(\uparrow) \ar[r]^-{\zeta} & (N\sQ)(\uparrow) = \fD.
}
\]
The map $\zeta'$ on elements of $\pk$ is given by the maps
\[
\nicexy{
\sK(\sigma C\tau) \ar[r]^-{\zeta} & (N\sQ)(\sigma C \tau)
}\]
on permuted corollas, using Lemma \ref{wnpgelement} to interpret graphices in $(N\sQ)(\sigma C \tau)$ as elements in $\sQ$.  To see that $\zeta'$ preserves the wheeled properad structure--namely, the bi-equivariant structure, colored units, dioperadic compositions, and contractions--note that the structure maps in $\pk$ are all induced by maps in $\Gammaw$ and that $N\sQ$ is a strict $\infty$-wheeled properad (Lemmas \ref{lem:wnerveissegal},  \ref{lem:wsegalisstrict}, and \ref{lem2:wsegalisstrict}).

To see that the composition $N\zeta' \circ \eta$ is equal to $\zeta$, first observe that they agree on $G \in \{\uparrow, \wheel, \sigma C\tau\}$.  Suppose $G \in \gwheelc$ is ordinary with $\edgei(G) \not= \varnothing$.  The map $\zeta_G$ is the composition
\[
\nicearrow
\xymatrix@C+10pt{
\sK(G) \ar[r]^-{\zeta_G} \ar[d]_{\chi_G} & (N\sQ)(G)\\
\sK(G)_1 \ar[r]^-{\prod \zeta} & (N\sQ)(G)_1 \ar[u]_{\chi_G^{-1}}^{\cong}
}\]
because the nerve $N\sQ$ satisfies the Segal condition (Lemma \ref{lem:wnerveissegal}).  By the definition of the corolla ribbon, the bottom horizontal map is determined by $\zeta$ on corollas and the exceptional edge, where it agrees with the composition $N\zeta' \circ \eta$.  The agreement between $\zeta$ and $N\zeta' \circ \eta$ on $G$ now follows from the definition of $\eta_G$ \eqref{wetankg}.

Finally, we observe the uniqueness of the wheeled properad map $\zeta'$ for which $N\zeta' \circ \eta = \zeta$.  Indeed, this equality already determines what the map $\zeta'$ does on color sets and elements in $\sP_{\sK}$.
\end{proof}

\begin{lemma}
\label{wstrictisnerve}
Suppose $\sK$ is a strict $\infty$-wheeled properad.  Then the object-wise bijections in Lemma \ref{wetaiso} assemble to give an isomorphism
\[
\nicearrow
\xymatrix{
\sK \ar[r]^-{\eta}_-{\cong} & N\sP_{\sK}
}\]
in $\gwheelcset$.
\end{lemma}

\begin{proof}
We reuse the proof of Lemma \ref{strictisnerve} essentially verbatim, with $\varGamma$, $\properad$, and Lemma \ref{etaknerveptwo} replaced by $\Gammaw$, $\wproperad$, and Lemma \ref{wetaknerveptwo}, respectively.
\end{proof}

With Lemmas \ref{lem:wnerveissegal}, \ref{lem:wsegalisstrict}, \ref{lem2:wsegalisstrict}, \ref{wproperadpk}, and \ref{wstrictisnerve},  the proof of Theorem \ref{wproperadnerve} is complete.

\begin{corollary}
\label{cor:fundwproperad}
Suppose $\sK$ is a strict $\infty$-wheeled properad.  Then the following statements hold.
\begin{enumerate}
\item
The $\sK(\uparrow)$-colored wheeled properad $\sP_{\sK}$ is canonically isomorphic to the fundamental wheeled properad of $\sK$.
\item
The map
\[
\nicearrow
\xymatrix{
\sK \ar[r]^-{\eta}_-{\cong} & N\sP_{\sK}}
\]
is the unit of the adjunction $(L,N)$.
\end{enumerate}
\end{corollary}

\begin{proof}
Lemma \ref{wstrictisnerve} says that $\eta$ is a map of wheeled properadic graphical sets, while Lemma \ref{wetaknerveptwo} says that it has the required universal property of the unit of the adjunction.
\end{proof}

\section{Fundamental Wheeled Properads of \texorpdfstring{$\infty$}{∞}-Wheeled Properads}
\label{sec:fundwproperad}

In this section, we observe that for a \emph{reduced} $\infty$-wheeled properad, its fundamental wheeled properad can be described using homotopy classes of $1$-dimensional elements.  This is the wheeled analog of Theorem \ref{thm:qkfundamental}.  First we define the homotopy relation of $1$-dimensional elements.  Then we show that a dioperadic composition and a contraction can be defined on homotopy classes of $1$-dimensional elements.  Next we show that the object with these operations is a wheeled properad.  Finally, we put the pieces together and prove the main Theorem \ref{thm:wqkfundamental}.

\subsection{Homotopy of \texorpdfstring{$1$-Dimensional}{1-Dimensional} Elements}

Recall the definitions in section \ref{rk:gsetproperad}, which also apply in $\gwheelcset$.  So for $\sK \in \gwheelcset$, a $1$-dimensional element  with $m$ inputs and $n$ outputs means a graphex in $\sK\left(\sigma C\tau\right)$ for some permuted corolla $\sigma C \tau$, where $C=C_{(m;n)}$ is the corolla with $m$ inputs and $n$ outputs.  The input/output profiles of a $1$-dimensional element is the pair of $\sK(\uparrow)$-profiles given by the images under the maps
\[
\sK\left(\sigma C\tau\right) \to \sK(\uparrow)
\]
induced by the outer dioperadic coface maps $\uparrow ~\to \sigma C \tau$ corresponding to the legs of $\sigma C\tau$.

\begin{definition}
Suppose $\sK \in \gwheelcset$, $f$ and $g$ are $1$-dimensional elements in $\sK$ with $m$ inputs and $n$ outputs.  For $1 \leq i \leq m$ (resp., $1 \leq j \leq n$), define 
\[\label{note:wsimi}
H \colon f \sim_i g \quad (\text{resp.}, H \colon f \sim^j g),
\]
exactly as in Definition \ref{def:graphhomotopy}.  If $H \colon f \sim_i g$ (resp., $H \colon f \sim^j g$), then we say \textbf{$f$ is homotopic to $g$ along the $i$th input} (resp., \textbf{$j$th output}).
\end{definition}

\begin{remark}
As explained in Remark \ref{rk:homotopyextension}, graphically the relation $f \sim_i g$ (resp., $f \sim^j g$) means that there is an $i$th input (resp., $j$th output) extension $H$ of $f$ by a degenerate element $\bone$, as in the pictures 
\begin{center}
\begin{tikzpicture}
\matrix[row sep=1.2cm, column sep=4cm]{
\node [plain, label=above:$...$] (v1) {$f$}; &
\node [plain] (u2) {$\bone$};\\
\node [plain] (u1) {$\bone$}; &
\node [plain, label=below:$...$] (v2) {$f$};\\
};
\foreach \x in {1,2}
{
\draw [inputleg] (v\x) to +(-.6cm,-.5cm);
\draw [inputleg] (v\x) to node[below right=.1cm]{\footnotesize{$m$}} +(.6cm,-.5cm);
\draw [outputleg] (v\x) to +(-.6cm,.5cm);
\draw [outputleg] (v\x) to node[above right=.1cm]{\footnotesize{$n$}} +(.6cm,.5cm);
}
\draw [inputleg] (u1) to +(0,-.7cm);
\draw [arrow] (u1) to node{\footnotesize{$i$}} (v1);
\draw [outputleg] (u2) to +(0,.7cm);
\draw [arrow] (v2) to node{\footnotesize{$j$}} (u2);
\end{tikzpicture}
\end{center}
whose inner dioperadic face is $g$.  To say that $H$ is an input or output extension of $f$ by $\bone$ means that the two outer dioperadic faces of $H$ are $f$ and $\bone$.
\end{remark}

\begin{lemma}
\label{wproperadhtpy}
Suppose $\sK \in \gwheelcset$, and $f$ and $g$ are $1$-dimensional elements in $\sK$ with $m$ inputs and $n$ outputs. 
\begin{enumerate}
\item
Suppose either
\begin{itemize}
\item
$f \sim_i g$ for some $1 \leq i \leq m$, or
\item
$f \sim^j g$ for some $1 \leq j \leq n$.
\end{itemize}
Then the profiles of $f$ and $g$ are equal.
\item
Suppose $\sK$ is an $\infty$-wheeled properad.  Then:
\begin{enumerate}
\item
The relations $\sim_i$ for $1\leq i \leq m$ and $\sim^j$ for $1 \leq j \leq n$ are all equivalence relations.
\item
The equivalence relations $\sim_i$ and $\sim^j$ are all equal to each other.
\end{enumerate}
\end{enumerate}
\end{lemma}

\begin{proof}
Simply reuse the proofs of Lemmas \ref{simprofile}, \ref{simijrelation}, and \ref{sameijrelation}.
\end{proof}

\begin{definition}
\label{def:whomotopy}
Suppose $\sK$ is an $\infty$-wheeled properad. Denote by $\sim$\label{note:wsim} the common equivalence relation, called \textbf{homotopy}, defined by $\sim_i$ and $\sim^j$ as in Lemma \ref{wproperadhtpy}.  If $m=n=0$, then homotopy is defined as the equality relation. Two elements in the same homotopy class are said to be \textbf{homotopic}. The homotopy class of a $1$-dimensional element $f$ will be written as $[f]$.
\end{definition}

As in previous chapters, in what follows we often omit isomorphisms induced by change of listings to simplify the presentation.  It is understood that such isomorphisms are applied wherever necessary.

\subsection{Dioperadic Compositions of \texorpdfstring{$1$-Dimensional}{1-Dimensional} Elements}

The following definition is the dioperadic version of Definition \ref{def:compqk}.

\begin{definition}
\label{def:diopqk}
Suppose:
\begin{itemize}
\item
$\sK$ is an $\infty$-wheeled properad,
\item
$f \in \sK(C_{(m;n)})$ is a $1$-dimensional element with profiles $\bah$, 
\item
$g \in \sK(C_{(p;q)})$ is a $1$-dimensional element with profiles$\dch$, and 
\item
$b_j=c_i$ for some $1 \leq j \leq n$ and $1 \leq i \leq p$.
\end{itemize}
\begin{enumerate}
\item
Define the dioperadic graph
\[
D = C_{\dch} \jcompi C_{\bah}
\]
with top vertex $v$, bottom vertex $u$, and unique internal edge $e$ (Example \ref{ex:dioperadic}). Suppose
\[
\nicearrow
\xymatrix{
C = C_{(\uc \compi \ua; \ub \compj \ud)} \ar[r]^-{d_{\inp}} & D
}\]
is the inner dioperadic coface map corresponding to shrinking away the unique internal edge $e$.
\item
By Example \ref{ex:diophorn} $g$ and $f$ define an inner horn
\[
\nicexy@C+12pt{
\Lambdaw^{d_{\inp}}[D] \ar[d] \ar[r] & \sK\\
\Gammaw[D] \ar@{.>}[ur]_-{\exists \theta} &
}\]
of $\sK$.  Since $\sK$ is an $\infty$-wheeled properad, there exists a filler $\theta \in \sK(D)$.  If
\[
h = d_{\inp}^*\theta \in \sK(C),
\]
then we write
\[
\theta \colon h \simeq g \jcompi f.
\]
Call $h$ a \textbf{dioperadic composition of $g$ and $f$}, and call $\theta$ a \textbf{witness}.\index{witness!of dioperadic composition}
\end{enumerate}
\end{definition}

\begin{remark}
In the context of Definition \ref{def:diopqk}, a witness is guaranteed to exist, but it is not necessarily unique because $\sK$ may not be strict.  The two outer dioperadic faces of a witness are $g$ and $f$, and its inner dioperadic face is a dioperadic composition of $g$ and $f$.
\end{remark}

Since a witness is not unique, neither is a dioperadic composition.  Next we will observe that the homotopy class of a dioperadic composition is well-defined.  First we observe that, given $g$ and $f$, any two dioperadic compositions are homotopic.  Recall that $\sK \in \gwheelcset$ is \emph{reduced} if $\sK(C_{\emptyprofh}) = \{*\}$.

\begin{lemma}
\label{lem1:diophtpy}
Suppose:
\begin{itemize}
\item
$\sK$ is a reduced $\infty$-wheeled properad,
\item
$g \in \sK(C_{(p;q)})$ and $f \in \sK(C_{(m;n)})$ are $1$-dimensional elements in $\sK$ as in Definition \ref{def:diopqk}, and
\item
there exist witnesses
\[
\theta \colon h \simeq g \jcompi f \andspace
\theta' \colon h' \simeq g \jcompi f.
\]
\end{itemize}
Then $h$ and $h'$ are homotopic.
\end{lemma}

\begin{proof}
The graphices $h$ and $h'$ have the same profiles, say $\yxh$.  If $\yxh = \emptyprofh$, then there is nothing to prove because $\sK(C_{\emptyprofh})$ is a singleton.

So suppose $\yxh \not= \emptyprofh$.  Then at least one of four things must happen: $q \geq 1$, $m \geq 1$, $p \geq 2$, or $n \geq 2$.  Assume $q \geq 1$; similar proofs exist for the other cases.  Since $q \geq 1$, the $3$-vertex graph
$A$, depicted as
\begin{center}
\begin{tikzpicture}
\matrix[row sep=1.5cm, column sep=1cm]{
\node [plain] (w) {$w$}; \\
\node [plain, label=above:$...$] (v) {$v$};\\
\node [plain, label=below:$...$] (u) {$u$};\\
};
\draw [inputleg] (u) to +(-.6cm,-.4cm);
\draw [inputleg] (u) to node[below right=.1cm]{\footnotesize{$m$}} +(.6cm,-.4cm);
\draw [outputleg] (u) to +(-.6cm,.4cm);
\draw [outputleg] (u) to node[above right=.1cm]{\footnotesize{$n$}}+(.6cm,.4cm);
\draw [inputleg] (v) to +(-.6cm,-.4cm);
\draw [inputleg] (v) to node[below right=.1cm]{\footnotesize{$p$}} +(.6cm,-.4cm);
\draw [outputleg] (v) to node[above right=.1cm]{\footnotesize{$q$}} +(.6cm,.4cm);
\draw [outputleg] (w) to +(0,.7cm);
\draw [arrow] (u) to node[near end]{\footnotesize{$i$}} node[near start]{\footnotesize{$j$}} node[swap]{\footnotesize{$\alpha$}} (v);
\draw [arrow, bend left=45] (v) to (w);
\end{tikzpicture}
\end{center}
is defined.  More precisely, $A$ is the simply connected graph
\[
\begin{split}
A 
&=C_{(1;1)} \jcompone D\\
&= C_{(1;1)} \jcompone \left[C_{(p;q)} \jcompi C_{(m;n)}\right].
\end{split}
\]
We now reuse the proof of Lemma \ref{homotopiccomp}, with $\left(\boxtimes, \be\right)$ replaced by $\left(\jcompi, \alpha\right)$, to see that there exist an inner horn filler $\Phi \in \sK(A)$ and an inner dioperadic face $H$ of $\Phi$ that is a homotopy from $h$ to $h'$.

If $m \geq 1$ (resp., $p \geq 2$, or $n \geq 2$), then we modify $A$ above by grafting the corolla $C_{(1;1)}$ with vertex $w$ to the dioperadic graph $D$ via an input leg of $u$ (resp., input leg of $v$, or output leg of $u$).
\end{proof}

The following observation will ensure that, for a reduced $\infty$-wheeled properad, dioperadic compositions can be defined on homotopy classes.

\begin{lemma}
\label{lem2:diophtpy}
Suppose:
\begin{itemize}
\item
$\sK$, $g \in \sK(C_{(p;q)})$, and $f \in \sK(C_{(m;n)})$ are as in Lemma \ref{lem1:diophtpy},
\item
there exist homotopies $f \sim f'$ and $g \sim g'$, and
\item
there exist witnesses
\[
\theta \colon h \simeq g \jcompi f \andspace
\theta' \colon h' \simeq g' \jcompi f'.
\]
\end{itemize}
Then $h$ and $h'$ are homotopic.
\end{lemma}

\begin{proof}
By Lemma \ref{lem1:diophtpy} it suffices to show that $h \simeq g' \jcompi f'$.  Consider the $3$-vertex graph
\[
\begin{split}
G 
&= C_{(p;q)} \jcompi \left[C_{(1;1)} \jcompone C_{(m;n)}\right]\\
&= \left[C_{(p;q)} \onecompi C_{(1;1)}\right] \jcompi C_{(m;n)},
\end{split}
\]
which may be depicted as follows.
\begin{center}
\begin{tikzpicture}
\matrix[row sep=1.3cm, column sep=1cm]{
\node [plain, label=above:$...$] (v) {$v$};\\
\node [plain] (t) {$t$};\\
\node [plain, label=below:$...$] (u) {$u$};\\
};
\draw [arrow] (u) to node[near start]{\footnotesize{$j$}} (t);
\draw [arrow] (t) to node[near end]{\footnotesize{$i$}} (v);
\draw [inputleg] (u) to +(-.6cm,-.4cm);
\draw [inputleg] (u) to node[below right=.1cm]{\footnotesize{$m$}} +(.6cm,-.4cm);
\draw [outputleg] (u) to +(-.6cm,.4cm);
\draw [outputleg] (u) to node[above right=.1cm]{\footnotesize{$n$}} +(.6cm,.4cm);
\draw [inputleg] (v) to +(-.6cm,-.4cm);
\draw [inputleg] (v) to node[below right=.1cm]{\footnotesize{$p$}} +(.6cm,-.4cm);
\draw [outputleg] (v) to +(-.6cm,.4cm);
\draw [outputleg] (v) to node[above right=.1cm]{\footnotesize{$q$}} +(.6cm,.4cm);
\end{tikzpicture}
\end{center}
We now reuse the proof of Lemma \ref{propcompdefined}, ignoring $\be$ and replacing $\boxtimes$ with $\jcompi$, to see that there exist an inner horn filler $\Psi \in \sK(G)$ and an inner dioperadic face $\psi$ of $\Psi$ that is a witness of $h \simeq g' \jcompi f'$.
\end{proof}

\subsection{Contractions of \texorpdfstring{$1$-Dimensional}{1-Dimensional} Elements}

The following definition is the contraction analog of Definition \ref{def:diopqk}.

\begin{definition}
\label{def:contractqk}
Suppose:
\begin{itemize}
\item
$\sK$ is an $\infty$-wheeled properad,
\item
$g \in \sK(C_{(p;q)})$ is a $1$-dimensional element with profiles $\dch$, and 
\item
$c_j = d_i$ for some $1 \leq j \leq p$ and $1 \leq i \leq q$.
\end{itemize}
\begin{enumerate}
\item
Define the contracted corolla
\[
\xi_eC_v = \xiij C_{\dch}
\]
with vertex $v$ and loop $e$ (Example \ref{ex:contractedcor}). Suppose
\[
\nicearrow
\xymatrix{
C = C_{(\uc \setminus c_j; \ud \setminus d_i)} \ar[r]^-{d_{\inp}} & \xi_eC_v
}\]
is the inner contracting coface map corresponding to deleting the loop $e$.
\item
By Example \ref{ex:contracthorn} $g$ defines an inner horn
\[
\nicexy@C+12pt{
\Lambdaw^{d_{\inp}}[\xi_eC_v] \ar[d] \ar[r] & \sK\\
\Gammaw[\xi_eC_v] \ar@{.>}[ur]_-{\exists \theta} &
}\]
of $\sK$.  Since $\sK$ is an $\infty$-wheeled properad, there exists a filler $\theta \in \sK(\xi_eC_v)$.  If
\[
h = d_{\inp}^*\theta \in \sK(C),
\]
then we write
\[
\theta \colon h \simeq \xiij g.
\]
Call $h$ a \textbf{contraction of $g$}, and call $\theta$ a \textbf{witness}.\index{witness!of contraction}
\end{enumerate}
\end{definition}

\begin{remark}
In the context of Definition \ref{def:contractqk}, the contracted corolla $\xi_eC_v$ has a unique inner face $d_{\inp}$ and a unique outer face
\[
\nicexy{
C_{\dch} \ar[r]^-{d_{\out}} & \xi_eC_v
}\]
corresponding to disconnecting the disconnectable edge $e$.  For a witness $\theta$, we have 
\[
d_{\out}^*\theta = g.
\]
\end{remark}

Now we observe that, for a reduced $\infty$-wheeled properad,  contraction can be defined on homotopy classes.

\begin{lemma}
\label{lem:contracthtpy}
Suppose:
\begin{itemize}
\item
$\sK$ is a reduced $\infty$-wheeled properad,
\item
$g \in \sK(C_{(p;q)})$ is a $1$-dimensional element in $\sK$ as in Definition \ref{def:contractqk},
\item
there exists a homotopy $g \sim g'$, and
\item
there exist witnesses
\[
\theta \colon h \simeq \xiij g \andspace
\theta' \colon h' \simeq \xiij g'.
\]
\end{itemize}
Then $h$ and $h'$ are homotopic.
\end{lemma}

\begin{proof}
The contractions $h$ and $h'$ have the same profiles $\yxh$.  If $\yxh = \emptyprofh$, then the reduced assumption on $\sK$ implies $h=h'$.  So we assume $\yxh \not= \emptyprofh$.  This means that $q \geq 2$, $p \geq 2$, or both.  We assume $q \geq 2$ and $i \not= 1$; similar proofs exist in the other cases.

Consider the graph $B$
\begin{center}
\begin{tikzpicture}
\matrix[row sep=1.5cm, column sep=1.5cm]{
\node [plain] (w) {$w$};\\
\node [plain] (v) {$v$};\\
};
\draw [arrow, out=55, in=-55, loop] (v) to 
node[near start]{\tiny{$i$}}
node{\footnotesize{$e$}} 
node[near end]{\tiny{$j$}} ();
\draw [arrow, bend left=45] (v) to node{\footnotesize{$f$}} (w);
\draw [inputleg] (v) to +(-.5cm,-.2cm);
\draw [inputleg] (v) to node[below right=-.1cm]{\tiny{$p$}} +(.5cm,-.2cm);
\draw [outputleg] (v) to node[above right=-.1cm]{\tiny{$q$}} +(.5cm,.2cm);
\draw [outputleg] (w) to +(0,.6cm);
\end{tikzpicture}
\end{center}
with two vertices and two internal edges, in which $e$ is a loop at $v$.  More precisely,
\[
B = \left(\xiij C_{(p;q)}\right)\left[ C_{(1;1)} \onecompone C_{(p;q)}\right].
\]
We will construct the desired homotopy $h \sim h'$ as an inner face of an inner horn filler of $\sK$ with respect to $B$.  First observe that $B$ has $4$ faces.
\begin{enumerate}
\item
There is an outer contracting coface map
\[
\nicexy{B_e \ar[r]^-{d_e} & B}
\]
corresponding to the disconnectable edge $e$.
\begin{center}
\begin{tikzpicture}
\matrix[row sep=1.5cm, column sep=1.5cm]{
\node [plain] (w1) {$w$};
&&& \node [plain] (w2) {$w$};\\
\node [plain, label=above:$...$, label=below:$...$] (v1) {$v$};
& \node [empty] (s) {}; 
& \node [empty] (t) {};
& \node [plain] (v2) {$v$};\\
};
\draw [arrow] (s) to node{$d_e$} (t);
\draw [arrow, out=55, in=-55, loop] (v2) to 
node[near start]{\tiny{$i$}}
node{\footnotesize{$e$}} 
node[near end]{\tiny{$j$}} ();
\foreach \x in {1,2}
{
\draw [inputleg] (v\x) to +(-.5cm,-.2cm);
\draw [inputleg] (v\x) to node[below right=-.1cm]{\tiny{$p$}} +(.5cm,-.2cm);
\draw [outputleg] (v\x) to node[above right=-.1cm]{\tiny{$q$}} +(.5cm,.2cm);
\draw [outputleg] (w\x) to +(0,.6cm);
\draw [arrow, bend left=45] (v\x) to node{\footnotesize{$f$}} (w\x);
}
\end{tikzpicture}
\end{center}
So $B_e$ is obtained from $B$ by disconnecting $e$.  By hypothesis there is a homotopy
\[
\left(H_0 \colon g \sim g'\right) \in \sK(B_e)
\]
along the first output.  This means that the outer dioperadic faces (resp., inner dioperadic face) of $H_0$ are $g$ and a degenerate element $\bone$ (resp., is $g'$).
\item
There is an outer dioperadic coface map
\[
\nicexy{B_w \ar[r]^-{d_w} & B}
\]
corresponding to the deletable vertex $w$.
\begin{center}
\begin{tikzpicture}
\matrix[row sep=1.5cm, column sep=1.5cm]{
&&& \node [plain] (w) {$w$};\\
\node [plain] (v1) {$v$};
& \node [empty] (s) {}; 
& \node [empty] (t) {};
& \node [plain] (v2) {$v$};\\
};
\draw [arrow] (s) to node{$d_w$} (t);
\foreach \x in {1,2}
{
\draw [arrow, out=55, in=-55, loop] (v\x) to 
node[near start]{\tiny{$i$}}
node{\footnotesize{$e$}} 
node[near end]{\tiny{$j$}} ();
\draw [inputleg] (v\x) to +(-.5cm,-.2cm);
\draw [inputleg] (v\x) to node[below right=-.1cm]{\tiny{$p$}} +(.5cm,-.2cm);
\draw [outputleg] (v\x) to node[above right=-.1cm]{\tiny{$q$}} +(.5cm,.2cm);
}
\draw [outputleg] (v1) to +(-.5cm,.2cm);
\draw [outputleg] (w) to +(0,.6cm);
\draw [arrow, bend left=45] (v2) to node{\footnotesize{$f$}} (w);
\end{tikzpicture}
\end{center}
So $B_w$ is obtained from $B$ by deleting the deletable vertex $w$.  There is a witness
\[
\left(\theta \colon h \simeq \xiij g\right) \in \sK(B_w).
\]
This means that $\theta$ has outer (resp., inner) contracting face $g$ (resp., $h$).
\item
There is an inner dioperadic coface map
\[
\nicexy{B^f \ar[r]^-{d^f} & B}
\]
corresponding to the internal edge $f$.
\begin{center}
\begin{tikzpicture}
\matrix[row sep=1.5cm, column sep=1.5cm]{
&&& \node [plain] (w) {$w$};\\
\node [plain] (v1) {$vw$};
& \node [empty] (s) {}; 
& \node [empty] (t) {};
& \node [plain] (v2) {$v$};\\
};
\draw [arrow] (s) to node{$d^f$} (t);
\foreach \x in {1,2}
{
\draw [arrow, out=55, in=-55, loop] (v\x) to 
node[near start]{\tiny{$i$}}
node{\footnotesize{$e$}} 
node[near end]{\tiny{$j$}} ();
\draw [inputleg] (v\x) to +(-.5cm,-.2cm);
\draw [inputleg] (v\x) to node[below right=-.1cm]{\tiny{$p$}} +(.5cm,-.2cm);
\draw [outputleg] (v\x) to node[above right=-.1cm]{\tiny{$q$}} +(.5cm,.2cm);
}
\draw [outputleg] (v1) to +(-.5cm,.2cm);
\draw [outputleg] (w) to +(0,.6cm);
\draw [arrow, bend left=45] (v2) to node{\footnotesize{$f$}} (w);
\end{tikzpicture}
\end{center}
So $B^f$ is obtained from $B$ by shrinking away the internal edge $f$.  There is a witness
\[
\left(\theta' \colon h' \simeq \xiij g'\right) \in \sK(B^f).
\]
This means that $\theta'$ has outer (resp., inner) contracting face $g'$ (resp., $h'$).
\item
There is an inner contracting coface map
\[
\nicexy{B^e \ar[r]^-{d} & B}
\]
corresponding to the loop $e$.
\begin{center}
\begin{tikzpicture}
\matrix[row sep=1.5cm, column sep=1.5cm]{
\node [plain] (w1) {$w$};
&&& \node [plain] (w2) {$w$};\\
\node [plain, label=above:$...$, label=below:$...$] (v1) {$v'$};
& \node [empty] (s) {}; 
& \node [empty] (t) {};
& \node [plain] (v2) {$v$};\\
};
\draw [arrow] (s) to node{$d$} (t);
\draw [arrow, out=55, in=-55, loop] (v2) to 
node[near start]{\tiny{$i$}}
node{\footnotesize{$e$}} 
node[near end]{\tiny{$j$}} ();
\draw [inputleg] (v1) to +(-.5cm,-.2cm);
\draw [inputleg] (v1) to node[below right=-.1cm]{\tiny{$p-1$}} +(.5cm,-.2cm);
\draw [outputleg] (v1) to node[above right=-.1cm]{\tiny{$q-1$}} +(.5cm,.2cm);
\draw [inputleg] (v2) to +(-.5cm,-.2cm);
\draw [inputleg] (v2) to node[below right=-.1cm]{\tiny{$p$}} +(.5cm,-.2cm);
\draw [outputleg] (v2) to node[above right=-.1cm]{\tiny{$q$}} +(.5cm,.2cm);
\foreach \x in {1,2}
{
\draw [outputleg] (w\x) to +(0,.6cm);
\draw [arrow, bend left=45] (v\x) to node{\footnotesize{$f$}} (w\x);
}
\end{tikzpicture}
\end{center}
So $B^e$ is obtained from $B$ by deleting the loop $e$.
\end{enumerate}

By Lemma \ref{whorndescription} the graphices $H_0$, $\theta$, and $\theta'$ define an inner horn
\[
\nicearrow
\xymatrix@C+12pt{
\Lambdaw^{d_e}[B] \ar[d] \ar[r] & \sK\\
\Gammaw[B] \ar@{.>}[ur]_-{\exists \Psi}
}\]
of $\sK$.  So there is an inner horn filler $\Psi \in \sK(B)$.  Consider the inner contracting face
\[
H \defn d^*\Psi \in \sK(B^e).
\]
We claim that $H$ is a homotopy $h \sim h'$ along the first output.  In other words, we need to show that $H$ has
\begin{itemize}
\item
outer dioperadic faces $h$ and a degenerate element, and 
\item
inner dioperadic face $h'$.
\end{itemize}
\begin{enumerate}
\item
There is a commutative diagram
\[
\nicexy{
C_w=C_{(1;1)} \ar[d]_{d_v} \ar[r]^-{d_{v'}} & B^e \ar[d]^{d}\\
B_e \ar[r]^-{d_e} & B
}\]
in $\Gammaw$, in which the left vertical (resp., top horizontal) map is the outer dioperadic coface map corresponding to deleting the deletable vertex $v$ in $B_e$ (resp., $v'$ in $B^e$).  Therefore, we have
\[
\begin{split}
\bone
&= d_v^*H_0\\
&= d_v^* d_e^*\Psi\\
&= d_{v'}^* d^*\Psi\\
&= d_{v'}^* H.
\end{split}
\]
This shows that an outer dioperadic face of $H$ is a degenerate element.
\item
There is a commutative diagram
\[
\nicexy{
C_{(p-1;q-1)} \ar[d]_{d^e} \ar[r]^-{d_{w}} 
& B^e \ar[d]^{d}\\
B_w \ar[r]^-{d_w} & B
}\]
in $\Gammaw$, in which the left vertical (resp., top horizontal) map is the inner contracting (resp., outer dioperadic) coface map corresponding to deleting the loop $e$ in $B_w$ (resp., deleting the deletable vertex $w$ in $B^e$).  Therefore, we have
\[
\begin{split}
h
&= (d^e)^* \theta\\
&= (d^e)^* d_w^*\Psi\\
&= d_{w}^* d^*\Psi\\
&= d_{w}^* H.
\end{split}
\]
This shows that the other outer dioperadic face of $H$ is $h$.
\item
There is a commutative diagram
\[
\nicexy{
C_{(p-1;q-1)} \ar[d]_{d} \ar[r]^-{d^f} 
& B^e \ar[d]^{d}\\
B^f \ar[r]^-{d^f} & B
}\]
in $\Gammaw$, in which the left vertical (resp., top horizontal) map is the inner contracting (resp., inner dioperadic) coface map corresponding to deleting the loop $e$ in $B^f$ (resp., shrinking away the internal edge $f$ in $B^e$).  Therefore, we have
\[
\begin{split}
h'
&= d^* \theta'\\
&= d^* (d^f)^* \Psi\\
&= (d^f)^* d^*\Psi\\
&= (d^f)^* H.
\end{split}
\]
This shows that the inner dioperadic face of $H$ is $h'$.
\end{enumerate}
We have shown that $H$ is a homotopy $H \colon h \sim h'$.

Still assuming $q \geq 2$, if $i=1$, then we use the modification
\[
\left(\xi^1_j C_{(p;q)}\right)\left[ C_{(1;1)} \kcompone C_{(p;q)}\right]
\]
of $B$ above in which the edge $f$ is connected to the $k$th outgoing flag of $v$ for some $1<k\leq q$.  

Finally, if $q=1$, then $i=1$ and $p \geq 2$.  Then we use a modification of $B$ in which the corolla $C_w = C_{(1;1)}$ is grafted to the contracted corolla $\xi^1_j C_v$ via any input leg.  In other words, in any case, since $\yxh \not= \emptyprofh$, the vertex $v$ must have an incoming flag or an outgoing flag $a$ that is not part of the loop $e$.  Form a variation of $B$ by grafting  $C_w=C_{(1;1)}$ to $\xi_eC_v$ via this flag $a$, and then argue as above.
\end{proof}

\subsection{Wheeled Properad of Homotopy Classes}

We now define the object that will be shown to be the fundamental wheeled properad of a reduced $\infty$-wheeled properad.

\begin{definition}
\label{def:wqk}
Suppose $\sK$ is a reduced $\infty$-wheeled properad.  Define a $\Sigma_{\sS(\sK(\uparrow))}$-bimodule $\qk$ as follows.
\begin{enumerate}
\item
For a pair $\yxh$ of $\sK(\uparrow)$-profiles, denote by $\qk\yxh$, or $\qk\yx$, the set of homotopy classes of $1$-dimensional elements in $\sK$ with profiles $\yxh$.  This is well-defined by Lemma \ref{wproperadhtpy}.
\item
For a color $c \in \sK(\uparrow)$, define the \textbf{$c$-colored unit} of $\qk$ as the homotopy class of the degenerate element
\[
\bone_c = s^*(c) \in \sK\left(C_{(1;1)}\right),
\]
where $\nicexy{C_{(1;1)} \ar[r]^-{s} & \uparrow}$ is the codegeneracy map.
\item
Define the $\Sigma$-bimodule structure on $\qk$,
\[
\nicearrow
\xymatrix{
\qk\yx \ar[r]^-{(\pi;\lambda)} & \qk\yxlambda,
}\]
using the isomorphisms in $\sK$ induced by input/output relabelings
\[
\nicearrow
\xymatrix{
\sigma C\tau \ar[r] & \lambda(\sigma C\tau)\pi = (\lambda\sigma) C (\tau\pi)
}\]
of permuted corollas.
\item
Define a \textbf{dioperadic composition} on $\qk$ using representatives of homotopy classes, i.e.,\index{dioperadic composition!of homotopy classes}
\[
[g] \jcompi [f] \defn \left[ h \right],
\]
where $h \simeq g \jcompi f$ is as in Definition \ref{def:diopqk}.  This is well-defined by Lemma \ref{lem2:diophtpy}.
\item
Define a \textbf{contraction} on $\qk$ using representatives of homotopy classes, i.e.,\index{contraction!of homotopy class}
\[
\xiij [g] \defn \left[ h \right],
\]
where $h \simeq \xiij g$ is as in Definition \ref{def:contractqk}.  This is well-defined by Lemma \ref{lem:contracthtpy}.
\end{enumerate}
\end{definition}

We first observe that $\qk$ forms a wheeled properad.  The following observation is the wheeled analog of Lemma \ref{qkisproperad}.

\begin{lemma}
\label{qkwproperad}
Suppose $\sK$ is a reduced $\infty$-wheeled properad.  Then $\qk$ in Definition \ref{def:wqk} is a $\sK(\uparrow)$-colored wheeled properad.\index{infinity wheeled properad!associated properad}
\end{lemma}

\begin{proof}
We need to check the various biased axioms of a wheeled properad.  They are all proved in the same manner using the unity and associativity of graph substitution, which are proved with full details and generality in \cite{jy2}.  So we will prove only one compatibility axiom between the dioperadic composition and the contraction in detail to illustrate the method.

For this purpose, we will reuse much of the proof of Lemma \ref{wproperadpk} as well as the notations therein.  Suppose $B$, $\theta \in \sK(C_1)$, and $\phi \in \sK(C_2)$ are as in that proof.  So the decorated graph
\begin{center}
\begin{tikzpicture}
\matrix[row sep=1cm, column sep=1cm]{
\node [plain] (v1) {$\theta$};\\
\node [plain] (v2) {$\phi$};\\
};
\draw [arrow] (v2) to node{\footnotesize{$e$}} (v1);
\draw [arrow, looseness=2, out=45, in=-45] (v1) to node{\footnotesize{$f$}} (v2);
\foreach \x in {1,2}
{
\draw [inputleg] (v\x) to +(-.5cm,-.3cm);
\draw [inputleg] (v\x) to +(.5cm,-.3cm);
\draw [outputleg] (v\x) to +(-.5cm,.3cm);
\draw [outputleg] (v\x) to +(.5cm,.3cm);
}
\end{tikzpicture}
\end{center}
makes sense.  We want to prove the compatibility axiom
\[
\xi_f\left([\theta] \circ_e [\phi]\right)
= 
\sigma\left[ \xi_e([\phi] \circ_f [\theta])\right]\tau.
\]
We proceed as in the proof of Lemma \ref{wproperadpk}, in particular suppressing the permutations $\sigma$ and $\tau$.
\begin{enumerate}
\item
The elements $\theta$ and $\phi$ determine an inner horn
\[
\nicexy{
\Lambdaw^{d_0}[B_f] \ar[r] & \sK,\\
}\]
which has a filler $\Phi_0 \in \sK(B_f)$.  By definition it is a witness
\[
\Phi_0 : d_0^*\Phi_0 \simeq \theta \circ_e \phi,
\]
where
\[
\nicexy{(B_f)_{12} \ar[r]^-{d_0} & B_f}
\]
is the inner dioperadic coface map corresponding to shrinking away the internal edge $e$ in $B_f$.
\item
The elements $\phi$ and $\theta$ determine an inner horn
\[
\nicexy{
\Lambdaw^{d_1}[B_e] \ar[r] & \sK,\\
}\]
which has a filler $\Phi_1 \in \sK(B_e)$.  By  definition it is a witness
\[
\Phi_1 : d_1^*\Phi_1 \simeq \phi \circ_f \theta,
\]
where
\[
\nicexy{(B_e)_{21} \ar[r]^-{d_1} & B_e}
\]
is the inner dioperadic coface map corresponding to shrinking away the internal edge $f$ in $B_e$.
\item
The element $d_0^* \Phi_0$ determines an inner horn
\[
\nicexy{
\Lambdaw^{d_2}[B_{12}] \ar[r] & \sK,\\
}\]
which has a filler $\Phi_2 \in \sK(B_{12})$.  By  definition it is a witness
\[
\Phi_2 : d_2^*\Phi_2 \simeq \xi_f\left(d_0^* \Phi_0\right),
\]
where
\[
\nicexy{B_{12}' \ar[r]^-{d_2} & B_{12}}
\]
is the inner contracting coface map corresponding to deleting the loop $f$ in $B_{12}$.
\end{enumerate}
These three witnesses are \emph{not} unique, but their existence is all that we need.

By Lemma \ref{whorndescription} the graphices $\Phi_i$ for $0 \leq i \leq 2$ above determine an inner horn
\[
\nicearrow
\xymatrix{
\Lambdaw^{d_{21}}[B] \ar[r] & \sK\\
}\]
of $\sK$, which has a filler $\Psi \in \sK(B)$.  Consider the two commutative diagrams \eqref{bonetwoprime} and \eqref{betwoone}.
Pulling back $\Psi$ along the first commutative diagram yields:
\[
\begin{split}
d_3^* d_{21}^* \Psi
&= d_2^* d_{12}^* \Psi\\
&= d_2^*\Phi_2\\ 
&\simeq \xi_f(d_0^* \Phi_0).
\end{split}
\]
Pulling back $\Psi$ along the the second commutative diagram yields:
\[
\begin{split}
d_e^* d_{21}^* \Psi
&= d_1^* d_{e}^* \Psi\\
&= d_1^* \Phi_1\\ 
&\simeq \phi \circ_f \theta.
\end{split}
\]
Note that $d_3^* d_{21}^* \Psi$ is the inner contracting face of $d_{21}^* \Psi$, while $d_e^* d_{21}^* \Psi$ is the outer contracting face of $d_{21}^* \Psi$.  Therefore, combining the above calculation, we now have:
\[
\begin{split}
\xi_f\left([\theta] \circ_e [\phi]\right)
&= \xi_f [d_0^* \Phi_0]\\
&= [d_2^* \Phi_2]\\
&= \xi_e [d_1^* \Phi_1]\\
&= \xi_e\left([\phi] \circ_f [\theta]\right).
\end{split}
\]
This proves the desired compatibility axiom in $\qk$.
\end{proof}

\subsection{Fundamental Wheeled Properad}

Below is the main observation of this section.  It says that the fundamental wheeled properad of a reduced $\infty$-wheeled properad can be described using homotopy classes of $1$-dimensional elements.  It is the wheeled analog of Theorem \ref{thm:qkfundamental}.

\begin{theorem}
\label{thm:wqkfundamental}
Suppose $\sK$ is a reduced $\infty$-wheeled  properad.  Then the  wheeled properad $\qk$ is canonically isomorphic to the fundamental wheeled properad of $\sK$.
\end{theorem}

\begin{proof}
We will reuse most of the proofs of Lemmas \ref{wetaiso}, \ref{wetaknerveptwo}, and \ref{wstrictisnerve}.
\begin{enumerate}
\item
Following the proof of Lemma \ref{wetaiso}, we obtain a map
\[
\nicearrow
\xymatrix{
\sK(G) \ar[r]^-{\eta_G} & (N\qk)(G)
}\]
for each $G \in \gwheelc$.  The map $\eta_{\uparrow}$ is a bijection because $\qk$ is $\sK(\uparrow)$-colored, and so $\eta_{\wheel}$ is also a bijection.  The map $\eta_{\sigma C\tau}$ for a permuted corolla $\sigma C \tau$ is a surjection because it sends a $1$-dimensional element in $\sK$ to its homotopy class.
\item
The object-wise map $\{\eta_G\}$ has the same universal property as the unit of the adjunction $(L,N)$ object-wise.  Here we follow the proof of Lemma \ref{wetaknerveptwo}.  The main point is that, given a wheeled properad $\sQ$ and a map
\[
\nicearrow\xymatrix{\sK \ar[r]^-{\zeta} & N\sQ}
\]
of wheeled properadic graphical sets, the map
\[
\nicearrow\xymatrix{\qk \ar[r]^-{\zeta'} & \sQ}
\]
is well-defined.  Indeed, homotopic $1$-dimensional elements in $\sK$ are sent to homotopic $1$-dimensional elements in $N\sQ$.  But since the nerve $N\sQ$ is a strict $\infty$-wheeled properad, homotopy is the identity relation.
\item
The object-wise map $\{\eta_G\}$ yields a map
\[
\nicearrow
\xymatrix{\sK \ar[r]^-{\eta} & N\qk}
\]
of wheeled properadic graphical sets.  Here we reuse the proof of Lemma \ref{wstrictisnerve} by simply replacing $\pk$ with $\qk$ and using the previous step instead of Lemma \ref{wetaknerveptwo}.
\end{enumerate}
The last two steps imply that the map $\nicearrow
\xymatrix{\sK \ar[r]^-{\eta} & N\qk}$ of graphical sets has the same universal property as the unit of the adjunction.  So, up to a canonical isomorphism, $\qk$ is the fundamental wheeled properad of $\sK$.
\end{proof}

\chapter{What's Next?}\label{C:applications}

In this brief chapter, we mention several problems related to infinity properads that we find interesting. 
These cover applications as well as conceptual understanding of infinity properads.

\section{Homotopy theory of infinity properads}  

In a future paper, we will show that the category of graphical sets admits a Quillen model structure \cite{quillen} so that the fibrant objects are precisely the infinity-properads. 
This is an extension of the Joyal model structure on simplicial sets \cite{joyal_theory,lurie} and the Cisinski-Moerdijk model structure on dendroidal sets \cite{cmb}.

The first two authors constructed \emph{simplicial} models of $\infty$-props in \cite{hr2},  which one can show induces a simplicial model for $\infty$-properads.  More generally, one can show that for any wheel-free pasting scheme $\cg$, the category of all simplicial $\cg$-props admits a cofibrantly generated  model category structure. This model structure is right proper and, we suspect, satisfies a type of relative left properness such as that in our paper \cite{hry}. 

It remains to show that these two models of $\infty$-properads are closely related.
In future work, we will construct a homotopy coherent graphical nerve as part of an adjoint pair between graphical sets and simplicial properads. 
The homotopy coherent graphical nerve takes an entrywise fibrant simplicial properad to an $\infty$-properad, and we expect that the adjunction is a Quillen equivalence.

\section{String topology infinity properad}\label{section string topology}

In \cite{poirier-rounds}, Poirier and Rounds initiate a study of chain-level string topology operations, which uses short geodesic arcs to define the operations, rather than the usual transversality arguments.
A \emph{string diagram} of type $(g,k,\ell)$ is a `fat graph' which determines a Riemann surface of genus $g$ with $k+\ell$ boundary components.
There is a compactified moduli space of string diagrams $\overline{SD}(g,k,\ell)$ so that the following hold.
\begin{itemize}
	\item The chains on $\overline{SD}(g,k,\ell)$ act on the chains of a free loop space of a manifold $M^d$, i.e.\ there is a chain map
	\begin{equation}\label{P-R Chain map}
		C_*(\overline{SD}(g,k,\ell)) \otimes C_*(LM) \to C_{*+(2-2g-k-\ell)d}(LM).
	\end{equation}
	\item There is an equivalence relation on the cells of $\overline{SD}(g,k,\ell)$ so that
	\[
		\left( \overline{SD}(g,k,\ell) / \sim \right) \simeq Sull(g,k,\ell)
	\]
	where the latter is the moduli space of Sullivan chord diagrams.
	\item In homology, the diagram
	\[ \begin{tikzcd}
		H_*(\overline{SD}(g,k,\ell)) \otimes H_*(LM) \dar \rar & H_{*+(2-2g-k-\ell)d}(LM) \\
		H_*(Sull(g,k,\ell)) \otimes H_*(LM), \arrow{ur}[swap]{CGC}
	\end{tikzcd} \]
	commutes, where the map $CGC$ gives the string topology operations defined by Cohen and Godin \cite{cohen-godin} and Chataur \cite{chataur}.
\end{itemize}
Poirier and Rounds point out that $\{ C_*(\overline{SD}(g,k,\ell) / \sim)\}$ is not a properad, as the composition (given by gluing of string diagrams) is only associative up to homotopy.
A question, then, is if this object is an infinity properad.
One possible approach is to show that the collection $C_*(\overline{SD}(g,k,\ell) / \sim)$ constitutes a strong homotopy properad in the sense of Gran{\aa}ker \cite{granaker} and complete the problem described in \ref{section sh properads}.
Another possibility is that the moduli spaces $\overline{SD}(g,k,\ell) / \sim$ themselves form an infinity properad.

\section{Operadic approach}\label{operadic approach}

When working with a fixed set of colors, properads are algebras over an operad.
More precisely, suppose that $\fC$ is a set of colors. 
Then there is an $\SC$-colored operad $\sO_{\fC}$ whose category of algebras is equivalent to the category of properads with color set $\fC$ and color-preserving maps.
The operad $\sO_{\fC}$ is generated in arity 2, and the generators are in bijection with the set of partially-grafted corollas in $\fC$.

Look at algebras (in the sense of \cite{heuts1}) over the infinity operad $N(\sO_{\fC})$. 
Each such algebra is then an up-to-homotopy $\fC$-colored properad, and it should be possible to produce a graphical set from such an algebra.
We expect that this graphical set is an infinity properad.

A related approach is as follows. Resolve the operad $\sO_{\fC}$ using the Boardman-Vogt $W$-construction to get an operad $W\sO_{\fC}$.
An algebra over this operad is an up-to-homotopy $\fC$-colored properad, which in analogy with the operad case in \cite{brinkmeier} we can call a  \emph{lax properad}.
Lax properads should give examples of infinity properads.
For simplicity, suppose that $\fC = *$ is a single point.
There is an operad $\sO_*^{\bullet \rightarrow \bullet}$ whose algebras are morphisms of monochrome properads.
Taking the Boardman-Vogt $W$ construction, we obtain an operad $W(\sO_*^{\bullet \rightarrow \bullet})$, whose algebras are homotopy homomorphisms of monochrome lax properads (compare to \cite[\S 19]{brinkmeier}).
A graph $G$ determines a free monochrome operad $F(G)$ with operations generated by the vertices, and this functor $F: \varGamma \to \properad$ induces a nerve functor $N_{lax}$ from lax properads to graphical sets.
For a lax properad $\sP$, define $N_{lax}(\sP)_G$ to be the set of homotopy homomorphisms $F(G) \rightsquigarrow \sP$.
We expect that this graphical set $N_{lax}(\sP)$ is an infinity properad.

\section{Strong homotopy properads}\label{section sh properads}

Lurie defined a nerve functor \cite[\S 1.3.1]{lurie2} which takes a differential graded category $\mathcal C$ and produces a quasi-category.
This was extended by Faonte \cite{faonte} in a conceptual way, so that the input $\mathcal C$ can be an $A_\infty$-category instead.
In brief, this nerve is described by defining the simplicial set whose $n$-simplices are the $A_\infty$-category functors from $\mathbb K [n]$ to $\mathcal C$.

On the other hand, Pepijn van der Laan \cite{vanderlaan} defined an operadic analogue of $A_\infty$-algebras, dubbed \emph{strong homotopy operads} or \emph{sh operads}. 
A natural question is whether these sh operads provide examples of infinity operads.
An affirmative answer was recently given by \cite{legrignou}.
In this paper, a colored version of Van der Laan's sh operads is constructed, called `strict unital homotopy colored operads', which are to sh operads as $A_\infty$-categories are to $A_\infty$-algebras.
Le Grignou defines the nerve functor on this expanded category, and shows that it produces infinity operads. 

Gran{\aa}ker \cite{granaker} defines a \emph{strong homotopy} (sh) {properad} to be a $\varSigma$-bimodule $\mathcal E$ together with a codifferential $\partial_{\mathcal E}$ on the cofree coproperad $\bar {\mathcal F}^c (\mathcal E [1])$.
A morphism $(\mathcal E, \partial_{\mathcal E}) \to (\mathcal E', \partial_{\mathcal E'})$ is a map of dg coproperads $\bar {\mathcal F}^c (\mathcal E [1]) \to \bar {\mathcal F}^c (\mathcal E' [1])$.
The present authors suggest that there should be a nerve functor from the category of sh properads to infinity properads, perhaps by pursuing a properadic version of \cite{legrignou}.

\section{Deformation theory}

At least over a field of characteristic zero, every deformation problem \cite{ger2} is governed by a deformation complex, i.e. a differential graded (dg) Lie algebra via solutions (modulo gauge action) of the Maurer-Cartan equation. By governed we mean that one can use the dg Lie algebra to functorially construct a deformation functor.  Although the interpretation of deformation problems in terms of solutions of the Maurer-Cartan equation is very useful, in many situations this is too rigid for a complete theory. The appropriate way of extending this category is to associate to each deformation problem a homotopy Lie algebra or $L_{\infty}$-algebras. 


Operads and, more recently properads, have proven extremely useful in deformation theory, since, among other things, operadic methods consolidate the construction of the deformation complexes mentioned above. The basic idea is that morphism complexes of (pr)operads have canonical $L_{\infty}$-structures.  In their papers \cite{mv1,mv2}, Merkulov and Vallette define the deformation theory of a morphism of prop(erad)s by showing that the vector space of morphisms between dg properads has a canonical filtered $L_{\infty}$-structure, where the Maurer-Cartan elements are morphisms of dg properads. They show that if one fixes a Maurer-Cartan element, $\gamma$, and let $Q^{\gamma}$ be the associated twisting of the canonical $L_{\infty}$-algebra by $\gamma$ then the deformation complex of $\gamma$ is defined by this twisted $L_\infty$-algebra.  

Several authors \cite{gk,gan,vallette} have explicitly computed deformation complexes by exploiting Koszul duality. Traditionally,  the deformation complex is given by the canonical cofibrant resolution of a properad. Since a general cofibrant resolution is often quite large, Merkulov and Vallette construct minimal models using ``homotopy'' Koszul duality of properads \cite{mv1,mv2}. A homotopy Koszul properad has a space of generators equal to the Koszul dual of a quadratic properad associated to it. Merkulov and Vallette compute the differentials of homotopy Koszul properads by appealing to the (dual) of Gran{\aa}ker \cite{granaker} strong homotopy properad construction. Given that we expect Gran{\aa}ker's construction to give rise to examples of infinity properads, we think it would be worth considering how this homotopy Koszul duality compares with Koszul duality of infinity properads. 

\section{Weber theory}

There is a well known and often exploited relationship between categories and simplicial sets via the nerve functor. 
In the classical example, the functor $N: \category \rightarrow \Set^{\varDelta^{op}}$ takes a small category $\mathcal{C}$ 
to the simplicial set $ N(\mathcal{C})$ with $N(\mathcal{C})_{n}$ the set of all composable chains of morphisms in $\mathcal{C}$ 
\[c_{0}\rightarrow c_{1}\rightarrow \ldots\rightarrow c_{n}\] 
of length $n$.  
The nerve functor is full and faithful, so it determines $\mathcal{C}$ up to isomorphism. Conversely,  we may view the notion of a simplicial set as a generalization of the notion of a category. 
We can tell if a simplicial set arises as the nerve of a category via very explicit conditions (e.g. Segal conditions or unique filling of inner horns). This implies that a category is a simplicial set that falls into this essential image and a functor is a map between simplicial sets of this form. 

Recent work of M. Weber \cite{weber,leinsterNlab} generalizes this relationship to construct nerves of objects more general than categories.  One can think of the category of all small categories as the category of algebras of the free category monad on the category of directed graphs. 
The generalized nerve construction says that in favorable cases, a monad $T$ induces a category $\varTheta(T)$, consisting of certain `linear' free $T$-algebras. This category $\varTheta(T)$ is a small full subcategory of $\alg(T)$.
Furthermore, there is an induced generalized nerve functor $N_T:\alg(T)\rightarrow \Set^{\varTheta(T)^{op}}$ which is full and faithful. The essential image of $N_T$ consists of the presheaves preserving certain limits.

Weber's machinery works on many known extensions of $\category$, including operads.  
In particular, Weber's nerve recovers the dendroidal nerve. 
Joyal and Kock \cite{joyalkock} have also applied this machinery to colored modular operads. 
Weber's construction applied to properads leads to a graphical category that is different than $\varGamma$ in this text. 
Our graphical category $\varGamma$ is Reedy, and we do not expect this to be true of Weber's $\varTheta(T)$. Further exploration of this idea could be quite interesting.

\backmatter


\newcommand{\where}[1]{\> \> \pageref{#1} \> \>}
\newcommand{\blob}{\> \> \> \> \hspace{1em}}

\Extrachap{Notation}

\begin{tabbing}
\textbf{Notation} \= \hspace{1.3cm}\= \textbf{Page}\= \hspace{.8cm}\=\textbf{Description} \\
$\frakC$ \where{note:colors} set of colors \\
$\uc$ or $c_{[1,m]}$ \where{note:profile} profile of colors \\
$\uc \circ_{\uc'} \ud$ \where{note:subprofile} substitution of profiles \\
$\pofc$ \where{note:pc} the groupoid of $\fC$-profiles \\
$\SC$ \where{note:sc} the category of pairs of $\fC$-profiles \\
$\dch$ or $\dc$ \where{note:dc} a pair of profiles \\
$\Flag$ \where{note:flag} the set of flags of a graph \\
$\vertex$ \where{note:vt} the set of vertices of a graph \\
$\Leg$ \where{note:leg} the set of legs of a graph \\
$\edge$ \where{note:edge} the set of edges of a graph \\
$\edgei$ \where{note:edge} the set of internal edges of a graph \\
$\inp(u)$ \where{note:in} the set of inputs \\
$\out(u)$ \where{note:in} the set of outputs \\
$\bullet$ \where{note:isovt} an isolated vertex \\
$\uparrow_c$ \where{note:uparrow} an exceptional edge \\
$\wheel_c$ \where{note:wheel} an exceptional wheel \\
$C_{\dch}$ \where{note:corolla} the $\dch$-corolla \\
$\sigma C\tau$ \where{note:permcor} a permuted corolla\\
$\xiij C$ or $\xi_eC$ \where{note:contcor} a contracted corolla \\
$C_{\dch} \boxtimes^{\uc'}_{\ub'} C_{\bah}$ \where{note:pgcor} a partially grafted corollas \\
$C_{\dch} \icompj C_{\bah}$ \where{note:dioperadic} a dioperadic graph \\
$\gwheelc$ \where{note:gwheelc} the set of connected graphs \\
$\gupc$ \where{note:gupc} the set of connected wheel-free graphs \\
$\gupci$ \where{note:gupci} the set of connected wheel-free graphs\\
\blob with non-empty inputs\\
$\gupco$ \where{note:gupci} the set of connected wheel-free graphs\\
\blob with non-empty outputs\\
$\gupcs$ \where{note:gupcs} the set of special connected wheel-free graphs\\
$\gupd$ \where{note:gupd} the set of simply connected graphs\\
$\uoperad$ \where{note:uoperad} the set of unital trees\\
$\ULin$ \where{note:ulin} the set of linear graphs\\
$G(\{H_v\})$ \where{note:graphsub} graph substitution \\
$G(H_w)$ \where{note:ghw} a special kind of graph substitution\\
$\dis(\SC)$ \where{note:dissc} the discrete category of $\SC$ \\
$\bone_c$ \where{note:coloredunit} the $c$-colored unit\\
$\boxtimes$ \where{note:propcomp} a properadic composition\\
$\properad$ \where{note:properad} the category of properads\\
$\properadi$ \where{note:properadi} the category of properads with\\
\blob non-empty inputs\\
$\properado$ \where{note:properadi} the category of properads with\\
\blob non-empty outputs\\
$\properads$ \where{note:properadi} the category of special properads\\
$\sP[G]$ \where{wheeldecoration} $\sP$-decorated graph\\
$F_{\gupc}$ \where{fgx} the free properad monad\\
$\sP \wedge \sQ$ \where{note:smash} smash product of colored objects\\
$\sP \otimes \sQ$ \where{gproptensor} tensor product of properads\\
$\Hom(\sP,\sQ)$ \where{note:hom} internal Hom\\
$\profilev$ \where{note:profv} the profiles of a vertex\\
$\ghat$ \where{note:ghat} the colored object of vertices in $G$\\
$\varGamma(G)$ \where{note:gammag} the graphical properad generated by $G$\\
$\sigma G \tau$ \where{note:relabeling} relabeling of a graph\\
$f(G)$ \where{note:image} the image of $G$ under $f$\\
$\varGamma$ \where{note:gamma} the properadic graphical category\\
$\varDelta$ \where{note:delta} the finite ordinal category\\
$\varOmega$ \where{note:omega} the dendroidal category\\
$\varTheta$ \where{note:theta} the graphical subcategory of\\
\blob simpy connected graphs\\
$\Gammai$ \where{note:gammai} the graphical subcategory of connected\\
\blob wheel-free graphs with non-empty inputs\\
$\Gammao$ \where{note:gammao} the graphical subcategory of connected\\
\blob wheel-free graphs with non-empty outputs\\
$\Iso(\mathcal{C})$ \where{note:isoc} maximal sub-groupoid\\
$\deg(G)$ \where{note:deg} degree of a graph\\
$\varGamma^+$ \where{note:gammaplus} a wide subcategory of $\varGamma$\\
$\varGamma^-$ \where{note:gammaplus} a wide subcategory of $\varGamma$\\
$\gupcset$ \where{note:gupcset} the category of properadic graphical sets\\
$N$ \where{note:nerve} the properadic nerve\\
$\varGamma[G]$ \where{note:representable} representable graphical set\\
$\varGamma^d[G]$ \where{note:dface} a face of $\varGamma[G]$\\
$\Lambda^d[G]$ \where{note:dhorn} a horn of $\varGamma[G]$\\
$\xi_v$ \where{note:corinclusion} corolla inclusion\\
$\sK(G)_1$ \where{note:corribbon} corolla ribbon\\
$\chi_G$ \where{segalmap} the properadic Segal map\\
$\Sc[G]$ \where{note:segalcore} the properadic Segal core\\
$\epsilon_G$ \where{note:epsilong} the Segal core inclusion\\
$\sim_i$ \where{note:simi} homotopy along the $i$th input\\
$\sim^j$ \where{note:simj} homotopy along the $j$th output\\
$\sim$ \where{note:sim} homotopy\\
$\xiij$ \where{note:xiij} contraction\\
$\jcompi$ \where{note:jcompi} dioperadic composition\\
$\wproperad$ \where{note:wproperad} the category of wheeled properads\\
$F_{\gwheelc}$ \where{note:fgwheelc} the free wheeled properad monad\\
$\sP \otimes \sQ$ \where{wproperadtensor} tensor product of wheeled properads\\
$\bonew$ \where{note:bonew} the unit of $\otimes$\\
$\Gammaw(G)$ \where{note:gammawg} graphical wheeled properad\\
$f(G)$ \where{note:fofg} the image of $G$ under $f$\\
$\Gammaw$ \where{note:gammaw} the wheeled properadic graphical category\\
$\Gammawimage$ \where{note:gammawimage} maps whose target and image coincide\\
$\Gammawout$ \where{note:gammawout} subcategory generated by \\
\blob outer cofaces and isomorphisms\\
$\Gammawin$ \where{note:gammawin} subcategory generated by \\
\blob inner cofaces and isomorphisms\\
$\Gammaw^\minus$ \where{note:gammawminus} a wide subcategory of $\Gammaw$ \\
$\gwheelcset$ \where{note:gwcset} the category of wheeled properadic\\
\blob graphical sets\\
$N$ \where{note:wnerve} the wheeled properadic graphical nerve\\
$\Gammaw[G]$ \where{note:wrepresentable} representable wheeled graphical set\\
$\Gammaw^d[G]$ \where{note:wface} a face of $\Gammaw[G]$\\
$\Lambdaw^d[G]$ \where{note:whorn} a horn of $\Gammaw[G]$\\
$\sK(G)_1$ \where{note:wcorribbon} corolla ribbon\\
$\chi_G$ \where{wsegalmap} the wheeled properadic Segal map\\
$\Sc[G]$ \where{note:wsegalcore} the wheeled properadic Segal core\\
$\epsilon_G$ \where{note:wepsilong} the wheeled properadic Segal core map\\
$\sim_i$ \where{note:wsimi} homotopy along the $i$th input\\
$\sim^j$ \where{note:wsimi} homotopy along the $j$th output\\
$\sim$ \where{note:wsim} homotopy\\

\end{tabbing}

%

\bibliographystyle{amsalpha}

\printindex

\end{document}